\documentclass[11pt,a4paper]{article}
\usepackage[utf8]{inputenc}
\usepackage[british,UKenglish,USenglish,english,american]{babel}
\usepackage[all]{xy}
\usepackage{enumitem}
\usepackage[dvipsnames]{xcolor}
\usepackage{graphicx}
\usepackage{stmaryrd}
\usepackage{amsfonts,amstext,amssymb,amsmath,amsthm,pspicture,bigints}
\usepackage{fullpage}
\usepackage{moreverb}
\usepackage{pifont}
\usepackage{pst-all}
\usepackage[left=2cm,right=2cm,top=2cm,bottom=2cm]{geometry}
\usepackage{esint}
\usepackage{fourier-orns}
\usepackage{mathrsfs}
\usepackage{array}
\newcommand{\T}{\mathbb{T}}
\usepackage{hyperref}
\newcommand{\ii }{{\rm i} }
\usepackage{authblk}
\newcolumntype{C}[1]{>{\centering\arraybackslash}b{#1}}
\newcolumntype{R}[1]{>{\raggedleft\arraybackslash}b{#1}}
\newcolumntype{L}[1]{>{\raggedright\arraybackslash}b{#1}}
\newcolumntype{M}[1]{>{\centering}m{#1}}
\newtheorem{theo}{Theorem}[section]
\newtheorem{defin}{Definition}[section]
\newtheorem{lem}{Lemma}[section]
\newtheorem{prop}{Proposition}[section]

\newtheorem{remark}{Remark}[section]

\newtheorem{cor}{Corollary}[section]
\usepackage{bbm}
\numberwithin{equation}{section}

\date{}
\title{Time quasi-periodic vortex patches for\\
	quasi-geostrophic shallow-water equations}
\author{Taoufik Hmidi\thanks{NYUAD Research Institute, New York University Abu Dhabi, PO Box 129188, Abu Dhabi, United Arab Emirates. Univ Rennes, CNRS, IRMAR – UMR 6625, F-35000 Rennes, France\\
Email address: thmidi@univ-rennes1.fr.} \hspace{0.5cm} and \hspace{0.5cm}Emeric Roulley\thanks{Univ Rennes, CNRS, IRMAR – UMR 6625, F-35000 Rennes, France\\
E-mail address : emeric.roulley@univ-rennes1.fr}}
\begin{document}
	\maketitle
	\begin{abstract}
		In this paper, we shall implement KAM theory in order to construct a large class of  time quasi-periodic solutions for an active scalar model arising in fluid dynamics. More precisely, the construction of invariant tori is performed for quasi-geostrophic shallow-water equations when the {\it Rossby deformation length} belongs to a massive Cantor set.  As a consequence, we  construct   pulsating  vortex patches   whose boundary is localized  in a thin annulus for any time.
	\end{abstract}
	\tableofcontents
	%\newpage
	\section{Introduction}
	We shall discuss in this introduction the quasi-geostrophic shallow-water equation which is a nonlinear and nonlocal transport equation generalizing 2D Euler equations and used to describe large scale motion for  the atmosphere and the ocean circulation.  A particular concern will be addressed to  the long time dynamics behavior  and especially to  the emergence of long-lived structures in the vortex patch setting. The  main  novelty in this work is to explore quasi-periodic solutions around the Rankine vortices  using KAM theory and Nash-Moser scheme in the spirit recent work of  Baldi, Berti, Haus and  Montalto  \cite{BBMH18}. 
	\subsection{Model, relative equilibria from periodic to  quasi-periodic solutions}
	\noindent The current paper  deals with the  quasi-geostrophic shallow-water equations  $(\textnormal{QGSW})_{\lambda},$   which is considered as  one of the most common  asymptotic models used to describe the large scale motion of the atmospheric and oceanic circulation and can be derived   asymptotically from the rotating shallow-water equations when Rossby and Froude numbers are small enough, for more details we refer to \cite{DP11,V17} and the references therein. This model is planar and the evolution of the potential vorticity $\boldsymbol{q}$  takes the form of a nonlinear and nonlocal transport equation,
	$$(\textnormal{QGSW})_{\lambda}\left\lbrace\begin{array}{ll}
		\partial_{t}\boldsymbol{q}+\mathbf{v}\cdot\nabla\boldsymbol{q}=0, & \mbox{in }\mathbb{R}_{+}\times\mathbb{R}^{2}\\
		\mathbf{v}=\nabla^{\perp}(\Delta-\lambda^{2})^{-1}\boldsymbol{q}, & \\
		\boldsymbol{q}(0,\cdot)=\boldsymbol{q}_{0},&
	\end{array}\right.\quad \mbox{ where }\quad \nabla^{\perp}=\left(\begin{array}{c}
		-\partial_{2}\\
		\partial_{1}
	\end{array}\right).$$
	Here $\mathbf{v}$ denotes  the velocity field which is solenoidal   and $\boldsymbol{q}$ is a scalar   function. 
	Physically, the parameter $\lambda$ is defined by 
	$$\lambda=\frac{\omega_{c}}{\sqrt{gH}},$$
	where $g$ is the gravity constant, $H$ is the mean active layer depth and $\omega_{c}$ is the Coriolis frequency, assumed to be  constant. In the literature, the number $\frac{1}{\lambda}$ is called the  \textit{Rossby  deformation length} or \textit{Rossby radius} and measures the length scale at which the rotation effects are balanced by  the  stratification. Notice that small values of $\lambda$ corresponds to a free surface which is nearly rigid and when $\lambda=0$ we get Euler equations written in the formulation velocity-vorticity.
	
	The main purpose of this paper is to explore the emergence of time quasi-periodic solutions  for $(\textnormal{QGSW})_{\lambda}$ when  $\lambda$ belongs to a massive  Cantor like set. This goal will be accomplished in the special class of vortex patches. Before coming to the details we shall first discuss the existence of relative equilibria which are solutions that do not change their shape during the motion. This subject is well-explored for Euler equations and new results have been established recently for $(\textnormal{QGSW})_{\lambda}$ in \cite{DHR19}.
	Next we shall discuss the vortex patch problem and explore some of its specific feature.\\
	
	\noindent $\blacktriangleright$ {\it Vortex patches/relative equilibria.} We consider an initial condition taking the form of a vortex patch, that is,  $\boldsymbol{q}_{0}=\mathbf{1}_{D_{0}}$ where $D_{0}$ is a bounded domain in $\mathbb{R}^{2}$.  Since $\mathbf{1}_{D_{0}}\in L^{1}(\mathbb{R}^{2})\cap L^{\infty}(\mathbb{R}^{2})$, then in a similar way to Euler equations  Yudovich's theory applies and implies that this system admits a unique global in time weak solution. Then from the transport equation governing the potential vorticity we find that  the patch structure is preserved during the motion and one gets
	$$\boldsymbol{q}(t,\cdot)=\mathbf{1}_{D_{t}}\quad \mbox{ where }\quad D_{t}:=\mathbf{\Phi}_{t}(D_{0}),$$
	where $\displaystyle \mathbf{\Phi}_{t}:\mathbb{C}\rightarrow\mathbb{C}$ stands for the flow map associated to the velocity field $\mathbf{v}$ and defined by
	\begin{equation}\label{definition of the flow map}
		\partial_{t}\mathbf{\Phi}_{t}(z)=\mathbf{v}(t,\mathbf{\Phi}_{t}(z))\quad \mbox{ and }\quad \mathbf{\Phi}_{0}=\textnormal{Id}_{\mathbb{C}}.
	\end{equation}
	The boundary motion in  the smooth case reduces to tackle the evolution of  a curve in the complex plane surrounding a constant area domain and subject to the deformation induced by its own effect. Local/global in time persistence of the boundary regularity  is a  relevant subject in fluid dynamics and has attracted a lot of attention during the past decades, not only for  Euler equations but also for similar active scalar equations  such as  generalized surface quasi-geostrophic equations,  the aggregation equation.% For more details about this subject we refer for instance to \cite{BC93,BGLV16,CMOV,C95,DZ16,DZ16-2,HH15-1,H05,MZ20} and the references therein. 
	The literature is huge and we shall restrict the discussion to some suitable contributions that fit with our main task. Let us now briefly see how to write down the contour dynamics equations, more details can be found in \cite{HMV13,HMV15}. 
	Given a smooth parametrization $z(t,\cdot):\mathbb{T}\rightarrow\partial D_{t}$ of the boundary of the patch, then as particles located  at the boundary move with the boundary then we get the evolution equation 
	\begin{equation}\label{Lagran-Fo}
		\left[\partial_{t}z(t,\theta)-\mathbf{v}(t,z(t,\theta))\right]\cdot\mathbf{n}(t,z(t,\theta))=0,
	\end{equation}
	where $\mathbf{n}(t,\cdot)$ is the outward normal vector to the boundary $\partial D_{t}$ of $D_{t}.$ This equation reflects the fact that the particle velocity and the boundary velocity admit the same normal components.  As we shall see later in Section \ref{sec contour}, one may deduce from the preceding equation that  the boundary equation will evolve through a nonlinear integro-differential equation. Looking for particular solutions where the domain moves without any shape  deformation  is a traditional subject in fluid dynamics and important developments have been performed   for long time ago. In the literature,  these  structures appear under different names: relative equilibria, V-states, long-lived structures, vortex crystals, etc\ldots A particular case  is given by  rotating patches which are solutions rotating uniformly with constant angular velocity about their center of mass that can be fixed at the origin, namely,
	$$\boldsymbol{q}(t,z)={\bf{1}}_{D_{t}}\quad \mbox{ with }\quad D_{t}=e^{it\Omega}D_{0}.$$
	These solutions are periodic in time with period $\frac{2\pi}{\Omega}$ or equivalently with frequency $\Omega.$
	For Euler equations, we know two explicit examples of rotating patches. The first one is the so-called Rankine vortices given by the discs which are stationary solutions not only for Euler but also for $(\textnormal{QGSW})_{\lambda}$. The second example is called Kirchhoff ellipses which rotate with the angular velocity $\Omega=\frac{ab}{(a+b)^{2}}$ where $a$ and $b$ are the semi-axes of the ellipse, we refer to \cite{K74} and \cite[p. 304]{BM02} for more explainations. Numerical experiments achieved  by Deem and Zabusky in \cite{DZ78} show the existence  of rotating solutions with $\mathbf{m}$-fold symmetry for $\mathbf{m}=3,4,5.$ An analytical proof based on the bifurcation theory and complex analysis tools was devised by Burbea  in \cite{B82} who proved the existence of $\mathbf{m}$-fold (for any $\mathbf{m}\in\mathbb{N}^{*}$) symmetric V-states bifurcating from Rankine vortices with angular velocity $\Omega_{\mathbf{m}}:=\frac{\mathbf{m}-1}{2\mathbf{m}}$. More investigations on the  V-states in different settings like  the doubly-connected patches, vortex pairs, boundary effects or for different models has been implemented during the past decade by several authors, for more details we refer to \cite{CCG16,CCG16-3,G20,G21,G19,HH15,HH21,HHH18, HHHM15,HMW20,HW21,HHMV16,HM16,HM16-2,HM17}. Concerning $(\textnormal{QGSW})_{\lambda}$ there are a few results dealing with relative equilibria.  Interesting numerical simulations showing the complexity and the richness  of the  bifurcation diagram with respect to the parameter $\lambda$ was studied  in  \cite{DJ20,DHR19}. In \cite{DHR19}, using bifurcation tools  the authors proved analogous results to those of Burbea. They show  in particular  the existence of branches  of $\mathbf{m}$-fold symmetric V-states ($\mathbf{m}\geq 2$) bifurcating from Rankine vortex ${\bf{1}}_{\mathbb{D}}$  with the angular velocity
	$$\Omega_{\mathbf{m}}(\lambda)=I_{1}(\lambda)K_{1}(\lambda)-I_{\mathbf{m}}(\lambda)K_{\mathbf{m}}(\lambda),$$
	where $I_m$ and $K_m$ are the modified Bessel functions of first and second kind. For more details about these functions, we refer to the \mbox{Appendix \ref{appendix Bessel}}. Notice that in the same paper the authors explored the two-fold branch when $\lambda$ is small and proved  first that it is located close to the ellipse branch of Euler equations and second it is not connected and from numerical simulations they put in evidence the fragmentation of this branch in multiple connected pieces. The second bifurcation from this branch was also analyzed leading to similar results as for Euler equations.\\
	\indent It is worthy to point out that the model under consideration is typically  a reversible Hamiltonian equation with  one degree of freedom given by the external the parameter $\lambda$.  Therefore, it is a natural issue to explore   whether or not quasi-periodic solutions constructed from the  linearized operator can survive at the  the nonlinear level  when $\lambda$ is selected in a suitable Cantor set. Studying  the  persistence of invariant tori is a relevant subject from KAM theory where a lot of important developments has been done in different directions. %Our main goal in this paper is to explore the existence of quasi-periodic solutions for $(\textnormal{QGSW})_{\lambda}$  in the setting of vortex patches.
	In the next paragraph we shall review and recall some basic notions and results in this subject.\\
	
	\noindent $\blacktriangleright$ {\it Quasi-periodic solutions.} A  function $f:\mathbb{R}\rightarrow\mathbb{R}$ is called quasi-periodic if there exists $F:\mathbb{T}^{d}\rightarrow\mathbb{R}$ such that  
	$$\forall\, t\in\mathbb{R},\quad f(t)=F(\omega t)$$
	for some frequency vector $\omega\in\mathbb{R}^{d}$ ($d\in\mathbb{N}^{*}$) which is non-resonant, that is
	\begin{equation}\label{nonresonnace omega}
		\forall\, l\in\mathbb{Z}^{d}\setminus\{ 0\},\quad \omega\cdot l\neq 0.
	\end{equation}
	Here and in the sequel, we denote $\mathbb{T}^{d}=\mathbb{R}^{d}/(2\pi\mathbb{Z})^{d}$ the flat torus of dimension $d$.  In the case $d=1$, we recover from this definition periodic functions with frequency $\omega\in\mathbb{R}^{*}.$\\
	The study of quasi-periodic solutions to Hamiltonian systems goes back to the pioneering  works of Kolmogorov \cite{K54}, Arnold \cite{A63} and Moser \cite{M62} where they proved, in finite dimension and under suitable non degeneracy and smoothness conditions, the persistence of invariant tori for  small perturbation of integrable Hamiltonian systems. Namely, using action-angle variables $(I,\vartheta)$ (integrability in the Liouville sense) we may write
	$$H(I,\vartheta)=h(I)+\varepsilon P(I,\vartheta).$$	
	The phase space of the integrable Hamiltonian system associated to $h$ is foliated by Lagrangian invariant tori carrying a resonant or non-resonant quasi-periodic dynamics. Roughly speaking, Kolmogorov's theorem asserts that, for  small  $\varepsilon$ and when the perturbation is analytic, many non-resonant Lagrangian invariant tori persist. Kolmogorov's proof is based on the reduction of the Hamiltonian to an integrable one using symplectic change of coordinates. This is done using a Newton scheme where at any level we may replace the remainder by an integrable contribution (depending only on the action) up to a new remainder which is smaller than the initial one. To do that we should solve a functional equation called the {\it homological equation} where we should avoid resonances and deal with small divisors problem. For a long time people like Poincar\'e thought that it was not possible to make the scheme  convergent due to the uncontrolled  loss of regularity. The key idea of Kolmogorov was to introduce Diophantine conditions to control the small divisors problem and get only an algebraic loss of regularity. Arnold made rigorous the idea of Kolmogorov and Moser extended it to a differentiable case using what is now called "Nash-Moser procedure". This is a modification  of  the standard Newton scheme making appeal to  regularizing operators  $(\Pi_{N})_N$ in order to solve an equation $F(U)=0$ in a Banach scale  allowing some fixed loss of derivatives at each step. 
	This  strategy was first introduced by Nash in \cite{N54} to prove the isometric embedding theorem and improved later by Moser. The theory born from these works was named afterwards KAM theory in their honor.\\
	\indent Later on, this  theory was explored and developed in partial differential equations  by several authors leading to important contributions and opening new perspectives. The complexity of the problem depends on the space dimension and on the structure of the equations. For example in the semi-linear case the nonlinearity can be seen as a bounded perturbation of the linear problem and this simplifies  a lot the problem of finding a right inverse of the linearized operator around a state close to the equilibrium. However in the quasi-linear case where the nonlinearity is unbounded and has the same order as the linear part the situation turns to be much more tricky. This is the case for instance in  the water-wave equations where several results has been obtained in the past few years on the periodic or quasi-periodic, standing or traveling settings, see \cite{AB11,BBMH18,BFM21,BFM21-1,BM18,IPT05,PT01}.\\
	Next we shall give some insights on  the general scheme performed  to construct quasi-periodic solutions in the quasi-linear setting  that was developed in particular by Berti and Bolle in \cite{BB15}. This approach is  robust and flexible and will  be adapted to  our framework with the suitable  modifications. The first step is to write in a standard way  the equations  using the action-angle  variables for the tangential part. When we linearize the nonlinear functional around a state near the equilibrium  we end with an operator with variable coefficients that we should invert approximately up to small errors provided  the external parameters belong to a suitable Cantor set defined through various  Diophantine conditions. To do that we first look for an approximate inverse using  an isotropic torus  built around the  initial one. It has the advantage to transform the linearized operator via symplectic change of coordinates into a triangular system up to errors that vanish when  tested against an invariant torus.  Notice that the outcome is that  the Hamiltonian has a good normal form structure such that we can almost decouple the dynamics in the phase space in tangential and  normal modes. On the tangential part the system can be solved in a triangular way provided we can invert the linearized operator on the normal part up to a small coupling error term. This is more or less a finite dimensional KAM theory appearing here. Then, the analysis reduces to invert the linearized operator on the normal part with is a small perturbation of a diagonal infinite dimensional matrix. This is done by conjugating  the  linearized operator to  a diagonal one with constant coefficients. This step is long and technical and most of the non-resonance conditions in the Cantor set arise  during this process.   This allows the construction of  an approximate inverse for the  linearized operator with adequate  tame estimates required along    Nash-Moser scheme. We point out that this approach has been successfully implemented    to generate   quasi-periodic solutions to several quasi-linear and fully nonlinear autonomous PDE's, see for instance  \cite{BBMH18,BBM14,BBM16,BM18}. Recent progress towards applying KAM theory for the vortex dynamics has been performed for gSQG equation \cite{HHM21} and Euler equation \cite{BHM}. Finally, we point out that the use of suitable isotropic tori is a commodity  but it is not essential to get the triangular structure up to small errors. This point will be discussed in Section \ref{sec Ham-tool}, see also \cite{HHM21}.
	
	\subsection{Main result and sketch of the proof}
	The contour dynamics equation stated in  \eqref{Lagran-Fo} can be written in a more tractable way using polar coordinates. This is meaningful at least for short time when the initial patch is close to the equilibrium state given by the Rankine vortex $\mathbf{1}_{\mathbb{D}}$ where $\mathbb{D}$ is the unit disc of $\mathbb{R}^{2}.$  Thus the boundary $\partial D_t$ will be parametrized as follows
	\begin{align}\label{BSSS}
		z(t,\theta)=R(t,\theta)e^{i\theta}\quad\mbox{ with }\quad R(t,\theta)=\left(1+2r(t,\theta)\right)^{\frac{1}{2}}.
	\end{align}
	We shall prove in Section \ref{sec contour} that the function $r$ satisfies a nonlinear and non-local transport equation taking the form   
	\begin{equation}\label{equation introduction}
		\partial_{t}r+F_{\lambda}[r]=0,
	\end{equation}
	with
	\begin{equation*}
		F_{\lambda}[r](t,\theta)=\bigintsss_{\mathbb{T}}K_{0}\left(\lambda A_{r}(t,\theta,\eta)\right)\partial_{\theta\eta}^{2}\Big(R(t,\eta)R(t,\theta)\sin(\eta-\theta)\Big)d\eta
	\end{equation*}
	and
	\begin{equation*}
		A_{r}(t,\theta,\eta)=\left|R(t,\theta)e^{i\theta}-R(t,\eta)e^{i\eta}\right|.
	\end{equation*}
	The function $K_0$ is a Bessel function of imaginary parts and it is defined in Appendix \ref{appendix Bessel}. 
	Next, we take a parameter $\Omega\neq0$ and look for the solutions in the form
	\begin{equation}\label{ansatz1}
		r(t, \theta)=\tilde r(t,\theta+\Omega t),
	\end{equation} 
	then the equation \eqref{equation introduction} is equivalent to (to alleviate the notation $\tilde r$ will be denoted by $r$)
	\begin{equation}\label{equation-Master}
		\partial_{t} r+\Omega\partial_{\theta}  r
		+F_{\lambda}[ r]=0.
	\end{equation}
	We point out  that the introduction of  the parameter $\Omega$ seems  at this level  artificial  but it will be used later  to fix  the degeneracy of the first eigenvalue associated  with  the linearized operator at the equilibrium state.
	% In the quasi-periodic setting, we should find a non-resonant  vector frequency  $\omega\in\mathbb{R}^{d}$  such that the equation \eqref{equation-Master} admits a solution  in the form  $\tilde{r}(t,\theta)=\widehat{r}(\omega t,\theta)$ with $\widehat{r}:\mathbb{T}^d\times\mathbb{T}\to\mathbb{R}$ being a smooth $(2\pi)^{d+1}$-periodic function. Then we can easily check that  $\widehat{r}$ , still denoted in what follows by $r$, satisfies
	%\begin{equation}\label{QP-E}
	%\omega\cdot\partial_{\varphi}r+\Omega\partial_{\theta}r+F_{\alpha}[r]=0.
	%\end{equation}
	As we shall see in Proposition \ref{proposition Hamiltonian formulation of the equation}, the equation \eqref{equation-Master} has an Hamiltonian structure
	\begin{equation}\label{hamiltionian formulation of the equation introduction}
		\partial_{t}r=\partial_{\theta}\nabla H(r),
	\end{equation}
	where the  Hamiltonian $H$ is related to the kinetic energy and the angular momentum which are prime integrals of the system. In the quasi-periodic setting, we should find a frequency vector  \mbox{$\omega\in\mathbb{R}^{d}$,} such that  the equation \eqref{equation-Master} admits solutions in the form  $r(t,\theta)=\widehat{r}(\omega t,\theta)$ with $\widehat{r}$ being a smooth $(2\pi)^{d+1}-$periodic function. Then $\widehat{r}$ satisfies (to alleviate the notation we keep the  notation $r$ for $\widehat{r}$)
	\begin{equation}\label{quasi-periodic equation introduction}
		\omega\cdot\partial_{\varphi}r+\Omega\partial_{\theta}r+F_{\lambda}[r]=0.
	\end{equation}
	To explore quasi-periodic solutions we should first check their existence at the linear level. Then according to Lemma \ref{lemma general form of the linearized operator} the linearized operator to \eqref{quasi-periodic equation introduction} around a given  small state $r$ is given by the  linear  Hamiltonian equation,
	\begin{align}\label{Linearized-op}
		\mathcal{L}_{r}\rho=0\quad\textnormal{with}\quad\mathcal{L}_{r}=\partial_{t}+\partial_{\theta}\left[V_{r}\cdot-\mathbf{L}_{r}\right],
	\end{align}
	where $V_{r}$ is a scalar function  defined by
	\begin{equation}\label{V-r-Hm1}
		V_{r}(\lambda,t,\theta)=\Omega+\frac{1}{R(t,\theta)}\bigintsss_{\mathbb{T}}K_{0}\left(\lambda A_{r}(t,\theta,\eta)\right)\partial_{\eta}\left(R(t,\eta)\sin(\eta-\theta)\right)d\eta
	\end{equation} 
	and $\mathbf{L}_{r}$ is a non-local operator in the form
	\begin{equation}\label{L-r-Hm1}
		\mathbf{L}_{r}(\rho)(\lambda,t,\theta)=\int_{\mathbb{T}}K_{0}\left(\lambda A_{r}(t,\theta,\eta)\right)\rho(t,\eta)d\eta.
	\end{equation}
	At the equilibrium state $r\equiv 0$, we find that the linearized operator is a Fourier multiplier, see Lemma \ref{lemma linearized operator at equilibrium},
	\begin{equation}\label{Linearized-op-eq}
		\mathcal{L}_{0}\rho=\partial_{t}\rho+V_{0}(\lambda)\partial_{\theta}\rho-\partial_{\theta}\mathcal{K}_{\lambda}\ast\rho\cdot
	\end{equation}
	where $\ast$ denotes the convolution product in the variable $\theta$ and 
	$$V_{0}(\lambda)=\Omega+I_{1}(\lambda)K_{1}(\lambda)\quad\textnormal{and}\quad\mathcal{K}_{\lambda}(\theta)=K_{0}\left(2\lambda\left|\sin\left(\tfrac{\theta}{2}\right)\right|\right).$$
	Expanding into Fourier series 
	$$\rho(t,\theta)=\sum_{j\in\mathbb{Z}}\rho_{j}(t)e^{\ii j\theta},$$
	yields to 
	\begin{align}\label{FME}
		\rho\in\ker(\mathcal{L}_{0})\quad\Longleftrightarrow\quad\rho(t,\theta)=\sum_{j\in\mathbb{Z}}\rho_{j}(0)e^{\ii(j\theta-\Omega_{j}(\lambda)t)},
	\end{align}
	where the eigenvalues $\Omega_{j}$ are defined by
	\begin{align}\label{Freq-equilibrium}
		\Omega_{j}(\lambda)=j\Big(\Omega+I_{1}(\lambda)K_{1}(\lambda)-I_{j}(\lambda)K_{j}(\lambda)\Big)
	\end{align}
	and  the Bessel functions of imaginary argument $I_n$ and $K_n$ are given by   \eqref{definition of modified Bessel function of first kind}. It is worthy to point out  that the frequency associated to the mode $j=0$ is vanishing and therefore it creates trivial resonance. This can be fixed by imposing  a zero space average which  can be maintained at the nonlinear level by virtue of the structure of  \eqref{hamiltionian formulation of the equation introduction}. Hence we shall work with the phase space of real functions enjoying this property, namely, 
	$$
	H_{0}^{s}:=H_{0}^{s}(\mathbb{T},\mathbb{R})=\Big\{ r(\theta)=\sum_{j\in\mathbb{Z}^{*}}r_{j}e^{\ii j\theta}\quad \textnormal{s.t.}\quad r_{-j}=\overline{r_{j}}\quad\textnormal{and}\quad\| r\|_{s}^{2}=\sum_{j\in\mathbb{Z}^{*}}|r_{j}|^{2}|j|^{2s}<\infty\Big\}.
	$$
	Another similar comment concerns the mode $j=1$ which vanishes  for any $\lambda$ when $\Omega=0$. This is why we have introduced $\Omega$ which should be strictly positive  to remedy to this defect and avoid any resonance at higher frequencies. The reversibility of the  system \eqref{hamiltionian formulation of the equation introduction} can be also exploited to find the requested parity property of the solutions. Actually, we can check that if $(t,\theta)\mapsto r(t,\theta)$ is a solution then  $(t,\theta)\mapsto r(-t,-\theta)$ is a solution too.  
	Then the solutions to the linear problem with this symmetry are in the form
	\begin{equation}\label{reversible solutions for the linearized equation at the equilibrium state}
		\rho(t,\theta)=\sum_{j\in\mathbb{Z}^{*}}\rho_{j}\cos\big(j\theta-{\Omega_{j}(\lambda)} t\big).
	\end{equation}
	Now, in order to generate quasi-periodic solutions to the linear problem it suffices to excite a finite number of frequencies from the linear spectrum. We shall then consider the following frequency vector.  
	$$\omega_{\textnormal{Eq}}(\lambda)=(\Omega_{j}(\lambda))_{j\in\mathbb{S}}\quad\hbox{with}\quad \mathbb{S}=\{j_{1},\ldots,j_{d}\}\subset \mathbb{N}^*.
	$$
	Note that throughout the paper, we use the notation
	$$\mathbb{N}=\{0,1,2,\ldots\}\quad\textnormal{and}\quad\mathbb{N}^*=\{1,2,\ldots\}.$$
	Notice that the vector $\omega_{\textnormal{Eq}}(\lambda)$ gives periodic solutions provided that it satisfies the non-resonant condition \eqref{nonresonnace omega}. This property holds true for almost all the values of $\lambda$ as it is proved in Proposition \ref{lemma sol Eq}. Our main result concerns the persistence of quasi-periodic solutions for the nonlinear model \eqref{hamiltionian formulation of the equation introduction} when the perturbation is small enough and the parameter $\lambda$ is subject to be in a massive Cantor set. 
	%%
	%
	%{main theorem}
	%
	%
	\begin{theo}\label{main theorem}
		Let $\lambda_{1}>\lambda_{0}>0$,  $d\in\mathbb{N}^{*}$ and  $\mathbb{S}\subset\mathbb{N}^*$ with $|\mathbb{S}|=d.$
		There exist  $\varepsilon_{0}\in(0,1)$ small enough with the following properties :  For every amplitudes ${\mathtt{a}}=(\mathtt{a}_{j})_{j\in\mathbb{S}}\in(\mathbb{R}_{+}^{*})^{d}$ satisfying
		$$|{\mathtt{a}}|\leqslant\varepsilon_{0},$$ 
		there exists a Cantor-like set $\mathcal{C}_{\infty}\subset(\lambda_{0},\lambda_{1})$ with asymptotically full Lebesgue  measure as ${\mathtt{a}}\rightarrow 0,$ i.e.
		$$\lim_{{\mathtt{a}}\rightarrow 0}|\mathcal{C}_{\infty}|=\lambda_{1}-\lambda_{0},$$
		such that for any $\lambda\in\mathcal{C}_{\infty}$, the  equation \eqref{hamiltionian formulation of the equation introduction} admits a time quasi-periodic solution with diophantine frequency vector ${\omega}_{\tiny{\textnormal{pe}}}(\lambda,{\mathtt{a}}):=(\omega_{j}(\lambda,{\mathtt{a}}))_{j\in\mathbb{S}}\in\mathbb{R}^{d}$ and taking  the form
		$$r(t,\theta)=\sum_{j\in\mathbb{S}}{\mathtt{a}_{j}}\cos\big(j\theta+\omega_{j}(\lambda,{\mathtt{a}})t\big)+\mathtt{p}({\omega}_{\tiny{\textnormal{pe}}}t,\theta),$$
		with
		$${\mathtt{\omega}}_{\tiny{\textnormal{pe}}}(\lambda,{\mathtt{a}})\underset{{\mathtt{a}}\rightarrow 0}{\longrightarrow}(-\Omega_{j}(\lambda))_{j\in\mathbb{S}},$$
		where $\Omega_{j}(\lambda)$ are the equilibrium frequencies defined in \eqref{Freq-equilibrium} and  the perturbation $\mathtt{p}:\T^{d+1}\to\mathbb{R}$ is an even function satisfying
		$$\| \mathtt{p}\|_{H^{{s}}(\mathbb{T}^{d+1},\mathbb{R})}\underset{{\mathtt{a}}\rightarrow 0}{=}o(|{\mathtt{a}}|)$$
		for some large   index of regularity $s.$ 
	\end{theo} 
	
	Before discussing the main steps of the proof, some remarks are in order.
	\begin{remark}
		\begin{enumerate}[label=\textbullet]
			\item  Combining this theorem with   \eqref{ansatz1} and \eqref{BSSS} we find that   the boundary shape of the patch can be parametrized  in polar coordinates as follows
			\begin{equation*}
				z(t,\theta)=R(t,\theta+\Omega t)e^{\ii \theta}\quad\mbox{ with }\quad R(t,\theta)=\left(1+2r(t,\theta)\right)^{\frac{1}{2}}
			\end{equation*}
			and $r$ is described  as in the theorem. 
			The time evolution of the shape is given  by small pulsation around the unit disc and the boundary is  localized in an annulus around the unit circle. 
			
			%			\item The quasi-periodic solutions constructed are perturbations, both for the frequencies and for the structure, of the quasi-periodic solutions for the linearized equation at the equilibrium state \eqref{reversible solutions for the linearized equation at the equilibrium state}.
			\item It is not at clear what could happen when $\lambda$ is in the complement of the Cantor set. We expect the linear invariant tori to be destroyed by the nonlinearity and filamentation may take place  generating fast increase of the curvature and the boundary length. 		
			%						\item The proof is based on a Nash-Moser scheme implemented in the Banach scale of Sobolev spaces. One may expect to have a similar result in the Banach scale of Hölder spaces. \textcolor{red}{Avec les injections de Sobolev, cette remarque est-elle vraiment pertinente ?}
			%			\item The case $d=1$ gives periodic solutions.
		\end{enumerate}
	\end{remark}
	
	We shall now outline the  main steps of the proof which will be developed following standard scheme as in   the preceding works \cite{BBMH18, BBM16, BB15, BM18} with different variations connected to the models structure.  We mainly use techniques from  KAM theory combined with Nash Moser scheme. This will be implemented along  several steps which are detailed below.\\
	
	\noindent $\blacktriangleright$ {{\bf{Step 1.} {\it Action-angle  reformulation.}} We first notice that the equation \eqref{hamiltionian formulation of the equation introduction} can be seen as a perturbation of the integrable system given by the linear dynamics at the equilibrium state. Indeed, by combining  \eqref{Linearized-op-eq},\eqref{Freq-equilibrium}  and  \eqref{hamiltionian formulation of the equation introduction} we may write
		$$\partial_{t}r=\partial_{\theta}\mathrm{L}(\lambda)(r)+X_{P}(r),$$
		where $\mathrm{L}(\lambda)$ and  the perturbed Hamiltonian vector field  $X_{P}$ are defined by
		$$\mathrm{L}(\lambda)(r)=-\big(\Omega+(I_1 K_1)(\lambda)\big) r+\mathcal{K}_\lambda\ast r\quad\hbox{and}\quad X_{P}(r)=I_{1}(\lambda)K_{1}(\lambda)\partial_{\theta}r-\partial_{\theta}\mathcal{K}_{\lambda}\ast r-F_{\lambda}[r].$$
		Since we are looking for small solutions then we find it  convenient to  rescale the solution $r\leadsto \varepsilon r$ with $\varepsilon$ a small positive number and consequently the new unknown still denoted by $r$ satisfies		\begin{equation*}
			\partial_{t}r=\partial_{\theta}\mathrm{L}(\lambda)(r)+\varepsilon X_{P_{\varepsilon}}(r),
		\end{equation*}
		where $X_{P_{\varepsilon}}$ is the Hamiltonian vector field defined by
		$X_{P_{\varepsilon}}(r)=\varepsilon^{-2}X_{P}(\varepsilon r).$ Then finding  quasi-periodic solutions with frequency $\omega\in \mathbb{R}^d$ amounts to solve the  equation
		\begin{equation*}
			\omega\cdot\partial_\varphi r=\partial_{\theta}\mathrm{L}(\lambda)(r)+\varepsilon X_{P_{\varepsilon}}(r).
		\end{equation*}
		Here we still use the same notation $r$ for the new profile which depends in the variables $(\varphi,\theta)\in\mathbb{T}^{d+1}$.
		The next step consists in  splitting the phase space $H_{0}^{s}$ into an orthogonal sum of  tangential and normal subspaces  as follows
		$$H_{0}^{s}=H_{\overline{\mathbb{S}}}\overset{\perp}{\oplus }H_{\perp}^{s},
		$$
		where $H_{\overline{\mathbb{S}}}$ is the finite dimensional subspace of real functions generated by $\{e^{\ii j\theta}, j\in\overline{\mathbb{S}}\}$ with $\overline{\mathbb{S}}=\mathbb{S}\cup(-\mathbb{S})$.
		For more details on this description we refer to Section \ref{subsec act-angl}. To track the dynamics it seems to be  more suitable to use the  action-angle variables $(I,\vartheta)$ seen  as symplectic polar variables for the Fourier coefficients of the tangential part  in $H_{\overline{\mathbb{S}}}.$ This leads to reformulate  the problem in terms of the embedded torus,
		$$i:\begin{array}[t]{rcl}
			\mathbb{T}^{d} & \rightarrow & \mathbb{T}^{d}\times\mathbb{R}^{d}\times H_{\perp}^{s}\\
			\varphi & \mapsto & (\vartheta(\varphi),I(\varphi),z(\varphi)),
		\end{array}$$
		with
		\begin{align*}
			r(\varphi,\theta)		&=\underbrace{v(\vartheta(\varphi),I(\varphi))(\theta)}_{\in H_{\overline{\mathbb{S}}}}+\underbrace{z(\varphi,\theta)}_{\in H_{\perp}^{s}}:=A(i(\varphi))(\theta)
		\end{align*}
		and
		$$
		v (\vartheta, I) = \sum_{j \in \overline{\mathbb{S}}}    
		\sqrt{\mathtt{a}_{j}^2+\tfrac{|j|}{2}I_j}\,  e^{\ii \vartheta_j}e_j,\quad e_j(\theta)=e^{\ii j\theta}\, . 
		$$
		Notice that the action and angle variables should satisfy the symmetry properties
		\begin{equation*}
			\forall j\in\overline{\mathbb{S}},\quad I_{-j}=I_j\in\mathbb{R}\quad \textnormal{and}\quad \vartheta_{-j}=-\vartheta_j\in \mathbb{R}.
		\end{equation*}
		Therefore we reduce the problem in the new variables  to construct  invariant tori with non-resonant frequency vector $\omega$ to the system
		\begin{equation}\label{Modi-Eq}
			\omega\cdot\partial_{\varphi}i(\varphi)=X_{H_{\varepsilon}}(i(\varphi)),
		\end{equation}
		where $X_{H_{\varepsilon}}$ is  the Hamiltonian vector field  associated to the Hamiltonian $H_{\varepsilon}$ given by 
		$$
		H_{\varepsilon}=-\omega_{\textnormal{Eq}}(\lambda)\cdot I+\tfrac{1}{2}\langle\mathrm{L}(\lambda)z,z\rangle_{L^{2}(\mathbb{T})}+\varepsilon\mathcal{P}_{\varepsilon},
		$$
		with $\mathcal{P}_{\varepsilon}$ defined by $\mathcal{P}_{\varepsilon}=P_{\varepsilon}\circ A.$ A useful  trick used by Berti and Bolle in  \cite{BB15} consists to solve first  the relaxed problem
		$$
		\omega\cdot\partial_{\varphi}i(\varphi)=X_{H_{\varepsilon}^{\alpha}}(i(\varphi)),
		$$
		where the vector field $X_{H_{\varepsilon}^{\alpha}}$ is associated to the modified Hamiltonian $H_{\varepsilon}^{\alpha}$ given by 
		$$H_{\varepsilon}^{\alpha}=\alpha\cdot I+\tfrac{1}{2}\langle\mathrm{L}(\lambda)z,z\rangle_{L^{2}(\mathbb{T})}+\varepsilon\mathcal{P}_{\varepsilon}.
		$$
		The advantage of this procedure  is to get one degree of freedom with the
		vector $\alpha\in\mathbb{R}^d$ that will be used  to ensure some compatibility assumptions  during the construction of an approximate
		inverse of the linearized operator. At the end of Nash-Moser scheme we shall  adjust implicitly the frequency $\omega$
		so that  $\alpha$ coincides with the equilibrium frequency $-\omega_{\textnormal{Eq}}(\lambda)$, which enables    to finally get  solutions to the original Hamiltonian equation. The relaxed problem can be written in the following form
		$$\mathcal{F}(i,\alpha,\lambda
		,\omega,\varepsilon)=0,$$
		with
		\begin{align}\label{PG1-C}
			\nonumber\mathcal{F}(i,\alpha,\lambda
			,\omega,\varepsilon)&:=\omega\cdot\partial_{\varphi}i(\varphi)-X_{H_{\varepsilon}^{\alpha}}(i(\varphi))\\
			&=\left(\begin{array}{c}
				\omega\cdot\partial_{\varphi}\vartheta(\varphi)-\alpha-\varepsilon\partial_{I}\mathcal{P}_{\varepsilon}(i(\varphi))\\
				\omega\cdot\partial_{\varphi}I(\varphi)+\varepsilon\partial_{\vartheta}\mathcal{P}_{\varepsilon}(i(\varphi))\\
				\omega\cdot\partial_{\varphi}z(\varphi)-\partial_{\theta}\left(\mathrm{L}(\lambda)z(\varphi)+\varepsilon\nabla_{z}\mathcal{P}_{\varepsilon}(i(\varphi))\right)
			\end{array}\right).
		\end{align}
		We point out that the linear   torus corresponding to the  linear  solution
		$$
		r(\varphi,\theta)= \sum_{j \in \overline{\mathbb{S}}}    
		{\mathtt{a}_{j}}\,  e^{\ii \varphi_j}e^{\ii j \theta}
		$$
		is given in the new coordinates system  by 
		$i_{\textnormal{\tiny{flat}}}(\varphi)=(\varphi,0,0)$
		and it is obvious that 
		$$
		\mathcal{F}\big(i_{\textnormal{\tiny{flat}}},-\omega_{\textnormal{Eq}}(\lambda),\lambda,-\omega_{\textnormal{Eq}}(\lambda),0\big)=0.
		$$
		We emphasize that  at this stage the classical implicit function theorem  does not work  because  the  linearized operator at the equilibrium state is not invertible due to  resonances.  One can  avoid resonances   by restricting the parameter $\lambda$ to a suitable Cantor set  according to some Diophantine conditions on  the linear  frequency $\omega_{\textnormal{Eq}}(\lambda)$ allowing in particular  to control the small divisors problem. By this way we get an  inverse at the equilibrium state but with algebraic loss of regularity. Unfortunately, this is not enough   to apply Nash-Moser scheme which requires   to construct a right inverse with tame estimates  in a small neighborhood of the equilibrium and this is the challenging deal in this topic. Indeed, the linearized operator is no longer with constant coefficients as for the integrable case and its  main part is not a Fourier multiplier.    At this level we are dealing with a quasilinear problem where  the  perturbation is unbounded.  \\
		\\
		$\blacktriangleright$ {{\bf{Step 2.} {\it  Approximate inverse.}}} Let $\alpha_0\in\mathbb{R}^d$ (actually $\alpha_0$ is a function of the parameters $\omega$ and $\lambda$) and  consider an embedded torus $i_0=(\vartheta_{0},I_0,z_0)$ near the flat one  with the  reversible structure, 
		$$\vartheta_{0}(-\varphi)=-\vartheta_{0}(\varphi),\quad I_0(-\varphi)=I_0(\varphi)\quad\textnormal{and}\quad z_{0}(-\varphi,-\theta)=z_{0}(\varphi,\theta).$$
		To deal with the    linearized operator $d_{i,\alpha}\mathcal{F}(i_0,\alpha_0)$, which exhibits    complicated structure,  and see whether we can construct an approximate inverse we should fix two important issues. One is related to the coupling structure in  the new coordinates system and the second is that the linearized operator is with variable coefficients. For this aim we shall follow the  approach conceived  by Berti and Bolle in \cite{BB15} with making suitable modifications. This approach consists  in  linearizing   around an isotropic torus  close enough to the original one  and then use   a symplectic change of coordinates leading   to a triangular   system  up to   small  errors, essentially of  "type $Z$" or highly decaying in frequency, that can be incorporated  in Nash-Moser scheme. Therefore to invert this  triangular system it suffices to  get an approximate  right  inverse  for the linearized operator in the normal direction, denoted in what follows by $\widehat{\mathcal{L}}_{\omega}$.  We notice that in Section \ref{sec approx-inv}, and similarly to \cite{HHM21}, we can bypass the use of isotropic torus by a slight modification of  Berti-Bolle approach.  Actually, according to Proposition \ref{prop conj DGO}, we  can conjugate the linearized operator with the  transformation  described by \eqref{trasform-sympl}  computed at the torus $i_0$  and get a triangular system with small errors mainly of  "type $Z$". The computations are performed  in a straightforward way using in a crucial way the  Hamiltonian structure of the original system. The main  advantage that simplifies some arguments   is  to require the invertibility  for   the linearized operator only at the torus itself and not necessary at a closer  isotropic one. By this way, we can avoid the accumulation of   different errors induced by the isotropic torus that one  encounters for example  in the estimates of the approximate inverse or in the multiple  Cantor sets generated along the different reduction steps where the coefficients should be computed at the isotropic torus. The final outcome of this first step is to reduce the invertibility to finding an approximate inverse of  $\widehat{\mathcal{L}}_{\omega}$ which takes, according to  Proposition \ref{lemma setting for Lomega}, the form
		$$\widehat{\mathcal{L}}_{\omega}=\Pi_{\mathbb{S}_0}^{\perp}\left(\mathcal{L}_{\varepsilon r}-\varepsilon\partial_{\theta}\mathcal{R}\right)\Pi_{\mathbb{S}_0}^{\perp}\quad\hbox{with}\quad\mathcal{L}_{\varepsilon r}=\omega\cdot\partial_\varphi+\partial_{\theta}\left[V_{\varepsilon r}\cdot-\mathbf{L}_{\varepsilon r}\right],
		$$
		where $\varepsilon\partial_{\theta}\mathcal{R}$ is a perturbation of finite rank, the function $V_{\varepsilon r}$ and  the nonlocal  operator $\mathbf{L}_{\varepsilon r}$  are  defined in \eqref{V-r-Hm1} and \eqref{L-r-Hm1}, respectively.
		At the equilibrium state (corresponding to $\varepsilon=0$)  $\widehat{\mathcal{L}}_{\omega}$ is diagonal and we shall see that the set of parameters  $(\lambda,\omega)$ leading to the existence of a right inverse is almost full. Now remark that even for $\varepsilon$ small, the perturbation affects the main part of the operator in a similar way to water waves \cite{BBMH18,BM18} or generalized SQG equation \cite{HHM21} and then we should construct the suitable change of coordinates in order to reduce the positive part of the operator to a diagonal operator. Later  we should implement KAM scheme to diagonalize the zero-order part.  This will be done in three steps.\\
		\foreach \x in {\bf a} {%
			\textcircled{\x}
		}  {\bf Reduction of the transport part.}}  This procedure  will be discussed  in  Proposition \ref{reduction of the transport part} and \mbox{Proposition \ref{Action on the non local part}.} We basically  use KAM techniques as in \cite{BFM21,FGMP19} in order to conjugate the operator $\mathcal{L}_{\varepsilon r}$, through  a suitable quasi-periodic symplectic change of coordinates $\mathscr{B}$, to a transport operator with constant coefficients. Indeed, we may construct an invertible transformation 
	\begin{equation*}
		\mathscr{B}\rho(\varphi,\theta)=\big(1+\partial_{\theta}\beta(\varphi,\theta)\big)\rho\big(\varphi,\theta+\beta(\varphi,\theta)\big)
	\end{equation*}
	and  a constant $c_{i_{0}}(\lambda,\omega)$ such that for any  given number $n\in\mathbb{N}$, if the parameter $(\lambda,\omega)$ belongs to  the  truncated set defined through the first order Melnikov condition 
	$$
	\mathcal{O}_{\infty,n}^{\gamma,\tau_{1}}(i_0)=\bigcap_{(l,j)\in\mathbb{Z}^{d}\times\mathbb{Z}\setminus\{(0,0)\}\atop|l|\leqslant N_{n}}\left\lbrace(\lambda,\omega)\in\mathcal{O}\quad\textnormal{s.t.}\quad\left|\omega\cdot l+jc_{i_{0}}(\lambda,\omega)\right|>\tfrac{4\gamma^{\upsilon}\langle j\rangle}{\langle l\rangle^{\tau_{1}}}\right\rbrace,
	$$
	then we have
	\begin{equation}\label{Asym-HHLL}
		\mathfrak{L}_{\varepsilon r}:=\mathscr{B}^{-1}\mathcal{L}_{\varepsilon r}\mathscr{B}=\omega\cdot\partial_{\varphi}+c_{i_{0}}(\lambda,\omega)\partial_{\theta}-\partial_{\theta}\mathcal{K}_{\lambda}\ast\cdot+\partial_{\theta}\mathfrak{R}_{\varepsilon r}+\mathtt{E}_{n}^{0},
	\end{equation}
	with $N_n= N_0^{(\frac{3}{2})^n}, N_0\geqslant 2,\,\upsilon\in(0,1),\, \mathcal{O}=(\lambda_{0},\lambda_{1})\times \mathscr{U},\,0<\lambda_{0}<\lambda_{1}$ and  $\mathscr{U}=B(0,R_{0})$ being  an open ball of $\mathbb{R}^{d}$ containing the curve of the linear vector frequency $\lambda\in(\lambda_0,\lambda_1)\mapsto -{\omega}_{\textnormal{Eq}}(\lambda)$. The \mbox{operator $\mathfrak{R}_{\varepsilon r}$}  is a self-adjoint Toeplitz integral operator satisfying the estimates	
	$$\forall s\in[s_{0},S],\quad \max_{k\in\{0,1,2\}}\| \partial_{\theta}^k\mathfrak{R}_{\varepsilon r}\|_{\textnormal{\tiny{O-d}},q,s}^{\gamma,\mathcal{O}}\lesssim \varepsilon\gamma^{-1}\left(1+\|\mathfrak{I}_{0}\|_{q,s+\sigma}^{\gamma,\mathcal{O}}\right),$$
	where the 	off-diagonal norm $\| \cdot\|_{\textnormal{\tiny{O-d}},q,s}^{\gamma,\mathcal{O}}$ is defined in \eqref{Top-NormX} and the loss of regularity $\sigma$ is connected \mbox{to $\tau_1$} and $d$ but it is independent of the index regularity $s$. Concerning the operator $\mathtt{E}_n^{0}$, we can show that it  is a small fast decaying remainder with the following  estimate for low regularity		
	\begin{equation}\label{decay-EE}
		\|\mathtt{E}_{n}^0\rho\|_{q,s_0}^{\gamma,\mathcal{O}}\lesssim \varepsilon\,N_0^{\mu_{2}}N_n^{-\mu_{2}}\|\rho\|_{q,s_0+2}^{\gamma,\mathcal{O}},
	\end{equation}
	where the weighted norms $\|\cdot\|_{q,s_0}^{\gamma,\mathcal{O}}$ are defined in \eqref{Norm-def}. For the  number $\mu_2$, it is connected to the regularity of the torus  $i_0$ and can be taken  large enough allowing to  identify the contributions of $\mathtt{E}_{n}^0$  as small errors  in the construction of the approximate inverse.  The next step will be discussed in  Proposition \ref{projection in the normal directions} where  we explore   the effect of the  transport  reduction  on the original operator $\widehat{\mathcal{L}}_{\omega}$ which is localized to the normal direction.  We prove that with   the localized transformation  $\mathscr{B}_{\perp}$ defined by 
	$$\mathscr{B}_{\perp}=\Pi_{\mathbb{S}_0}^\perp\mathscr{B}\Pi_{\mathbb{S}_0}^\perp,
	$$
	one obtains in  the Cantor set $\mathcal{O}_{\infty,n}^{\gamma,\tau_1}(i_0),$ 
	\begin{equation}\label{Remaind-KLP1}
		\mathscr{B}_{\perp}^{-1}\widehat{\mathcal{L}}_{\omega}\mathscr{B}_{\perp}=\omega\cdot\partial_{\varphi}\Pi_{\mathbb{S}_{0}}^{\perp}+\mathscr{D}_{0}+\mathscr{R}_{0}+\mathtt{E}_{n}^{1}, 
	\end{equation}
	where $\mathscr{D}_{0}$ is a diagonal operator whose spectrum $\{\ii \mu_{j}^{0}, j\in\mathbb{S}_0^c\}$ satisfies
	$$\mu_{j}^{0}(\lambda,\omega,i_{0})=\Omega_{j}(\lambda)+jr^{1}(\lambda,\omega,i_{0})\quad\mbox{ with }\quad \|r^1\|_{q}^{\gamma,\mathcal{O}}\lesssim {\varepsilon }$$
	and $\mathscr{R}_{0}$ is a remainder term taking the form of  an integral operator with  Toeplitz and reversibility structures with the estimates the asymptotic 
	\begin{equation}\label{estim partial thetamathfrakRr}
		\forall s\in[s_{0},S],\quad \max_{k\in\{0,1\}}\| \partial_{\theta}^k\mathscr{R}_{0}\|_{\textnormal{\tiny{O-d}},q,s}^{\gamma,\mathcal{O}}\lesssim \varepsilon\gamma^{-1}\left(1+\|\mathfrak{I}_{0}\|_{q,s+\sigma}^{\gamma,\mathcal{O}}\right),
	\end{equation}
	We remark  that the operator $\mathtt{E}_n^1$   satisfies similar estimates as for $\mathtt{E}_n^0$ seen before in  \eqref{decay-EE}. Finally, we want to emphasize that the  derivation of  the asymptotic structure of the operator $\mathfrak{L}_{\varepsilon r}$ seen before in \eqref{Asym-HHLL} requires some refined analysis. The delicate point concerns the expansion of the operator $\mathbf{L}_{r}$ defined in  \eqref{L-r-Hm1} and for this part we use the kernel structure detailed in  \eqref{explicit form for K0}  
	
	\begin{equation*}
		K_{0}(z)=-\log\left(\tfrac{z}{2}\right)I_{0}(z)+\sum_{m=0}^{\infty}\tfrac{\psi(m+1)}{(m!)^{2}}\left(\tfrac{z}{2}\right)^{2m}.
	\end{equation*}
	with $I_0$ being analytic.	This is different from the cases discussed before as for the   water waves in \cite{BBMH18, BM18}  where the kernel is given by that of Euler equations (corresponding to $\lambda=0$), that is, $K(z)=-\log\left(\frac{z}{2}\right)$. In this latter case  the deformed kernel enjoys the structure
	$$-\log\left(2A_{r}(t,\theta,\eta)\right)=-\log\left|\sin\left(\tfrac{\theta-\eta}{2}\right)\right|+\hbox{smooth nonhomogeneous kernel}.
	$$
	This means that the associated operator is given by a diagonal operator  of order $-1$ up to a smoothing non diagonal  pseudo-differential operator in $OPS^{-\infty}.$ In our context, this decomposition fails for $\lambda>0$ and we   get a similar one but with less smoothing operator. Actually we obtain from \eqref{fundamental decomposition of K0(Ar)} the splitting
	\begin{equation}\label{Split-Kernel-Sep}
		K_{0}\left(2\lambda A_{r}(t,\theta,\eta)\right)=K_{0}\left(2\lambda\sin\left(\tfrac{\theta-\eta}{2}\right)\right)+\mathscr{K}(\eta-\theta)\mathscr{K}_{r,1}^{1}(\lambda,\varphi,\theta,\eta)+\mathscr{K}_{r,1}^{2}(\lambda,\varphi,\theta,\eta),
	\end{equation} 			
	where the kernels  $\mathscr{K}_{r,1}^{1}$ and $\mathscr{K}_{r,1}^{2}$ are smooth whereas  $\mathscr{K}$ is slightly singular taking the form
	\begin{equation*}
		\mathscr{K}(\theta)=\sin^{2}\left(\tfrac{\theta}{2}\right)\log\left(\left|\sin\left(\tfrac{\theta}{2}\right)\right|\right).
	\end{equation*}
	%\ding{226} \textbf{Sub-step2: 
	%
	\foreach \x in {\bf b} {%
		\textcircled{\x}
	}    \textit{\bf KAM reduction of the remainder term.} 
	This is the main target  of  Section \ref{KAM-red-Feb} and the result is  stated  in Proposition \ref{reduction of the remainder term}. The goal is to conjugate  the remainder $\mathscr{R}_{0}$ of \eqref{Remaind-KLP1} and transform  it   into a diagonal operator. This will be developed in a standard way by constructing successive transformations through the KAM reduction  allowing to replace at each step the old remainder by a new one which is much smaller provided that we make the suitable parameters extraction. This scheme can be achieved unless  we solve the associated  {\it{homological equation}}. To avoid resonances we should at each step  make an extraction from the parameters set through the second order Melnikov conditions and the final outcome is as follows,
	$$
	\mathscr{L}_{\infty}:=\Phi_{\infty}^{-1}\mathscr{L}_{0}\Phi_{\infty}=\omega\cdot\partial_{\varphi}\Pi_{\mathbb{S}_0}^{\perp}+\mathscr{D}_{\infty},
	$$
	where $\mathscr{D}_{\infty}=\left(\ii\mu_{j}^{\infty}(\lambda,\omega,i_{0})\right)_{(l,j)\in\mathbb{Z}^{d}\times\mathbb{S}_{0}^{c}}$ is a diagonal operator  with pure imaginary spectrum and $\Phi_{\infty}$ is a reversible invertible operator. This reduction is possible when the parameters $(\lambda,\omega)$ belong to  the following  Cantor-like set,
	$$\mathcal{O}_{\infty,n}^{\gamma,\tau_{1},\tau_{2}}(i_{0})=\bigcap_{(l,j,j_{0})\in\mathbb{Z}^{d}\times(\mathbb{S}_{0}^{c})^{2}\atop|l|\leqslant N_{n}}\left\lbrace(\lambda,\omega)\in\mathcal{O}_{\infty,n}^{\gamma,\tau_{1}}(i_{0})\quad\textnormal{s.t.}\quad\left|\omega\cdot l+\mu_{j}^{\infty}(\lambda,\omega,i_{0})-\mu_{j_{0}}^{\infty}(\lambda,\omega,i_{0})\right|>\tfrac{2\gamma\langle j-j_{0}\rangle}{\langle l\rangle^{\tau_{2}}}\right\rbrace.$$
	The eigenvalues admit the following asymptotic
	$$
	\mu_{j}^{\infty}(\lambda,\omega,i_{0})=\Omega_{j}(\lambda)+jr^{1}(\lambda,\omega,i_{0})+r_j^\infty(\lambda,\omega,i_{0}),
	$$
	where  $r^{1}$ and $r_j^\infty$ are real small coefficients with  Lipschitz dependence with respect to the torus. Indeed, we have  
	\begin{align*}
		\|r^1\|_{q}^{\gamma,\mathcal{O}}+\sup_{j\in\mathbb{S}_0^c}|j|\|r_j^\infty\|_{q}^{\gamma,\mathcal{O}}&\lesssim\varepsilon\gamma^{-1},\nonumber\\
		\|\Delta_{12}r^1\|_{q}^{\gamma,\mathcal{O}}+\sup_{j\in\mathbb{S}_0^c}\|\Delta_{12}r_j^{\infty}\|_{q}^{\gamma,\mathcal{O}}&\lesssim\varepsilon\gamma^{-1}\|\Delta_{12}i\|_{q,\overline{s}_h+\sigma}^{\gamma,\mathcal{O}},
	\end{align*}
	for some index regularity $\overline{s}_h+\sigma$ and $\Delta_{12}r^1=r^{1}(\lambda,\omega,i_{1})-r^{1}(\lambda,\omega,i_{2}).$\\
	\foreach \x in {\bf c} {%
		\textcircled{\x}
	}    \textit{\bf Construction of the approximate inverse.}
	The next step is to invert approximately  the operator $\widehat{\mathcal{L}}_{\omega}$ detailed in Proposition \ref{inversion of the linearized operator in the normal directions}.  First we establish an approximate inverse for  $\mathscr{L}_{\infty},$ on the Cantor set 
	$$
	\Lambda_{\infty,n}^{\gamma,\tau_{1}}(i_0)=\bigcap_{(l,j)\in\mathbb{Z}^{d}\times\mathbb{S}_0^c\atop|l|\leqslant N_{n}}\left\lbrace(\lambda,\omega)\in\mathcal{O}\quad\textnormal{s.t.}\quad|\omega\cdot l+\mu_{j}^{\infty}(\lambda,\omega)|>\tfrac{\gamma\langle j\rangle}{\langle l\rangle^{\tau_{1}}}\right\rbrace.
	$$
	Then, introducing  the  Cantor set 
	$$\mathcal{G}_{n}(\gamma,\tau_1,\tau_2,i_0)=\mathcal{O}_{\infty,n}^{\gamma,\tau_{1}}(i_{0})\cap\mathcal{O}_{\infty,n}^{\gamma,\tau_{1},\tau_{2}}(i_{0})\cap\Lambda_{\infty,n}^{\gamma,\tau_{1}}(i_0),
	$$
	we are able to construct an approximate inverse of $\widehat{\mathcal{L}}_{\omega}$ in the following sense,
	$$ \widehat{\mathcal{L}}_{\omega}\mathtt{T}_{\omega,n}=\textnormal{Id}+\mathtt{E}_n\quad\textnormal{in }\mathcal{G}_n,$$
	where $\mathtt{E}_n$ is a fast frequency decaying operator as in \eqref{decay-EE} and $\mathtt{T}_{\omega,n}$ satisfies tame estimates uniformly in $n$.  Therefore coming back to Section \ref{sec approx-inv}, more precisely to Theorem \ref{proposition alomst approximate inverse}, this enables to construct an approximate right inverse $\mathrm{T}_0$  for the full differential  $d_{i,\alpha}\mathcal{F}(i_0,\alpha_0)$  enjoying suitable  tame estimates.\\	
	In what follows we want to make  some comments. The first one concerns the Lipschitz dependence of the eigenvalues with respsect to the torus. This is required in studying the stability of Cantor sets in Nash-Moser scheme and allows to construct a final massive Cantor set. As to the second one, it concerns  the  KAM reduction which  allows to diagonalize the operator when the parameters belong to a Cantor set like, even though   all the involved  transformations and  operators can be extended   in the whole set of parameters using standard cut-off functions for the Fourier coefficients. This extension with adequate estimates will be needed later during the implementation of Nash-Moser scheme. This is not  the only way to produce suitable extensions and one  expects    Whitney extension to be also well adapted  as  in   \cite{BBMH18, BM18}.  In our case we privilege the first procedure  which  can be easily set up and  manipulated using classical \mbox{ functional tools.}
	The last comment  is related to a technical point in KAM reduction, Contrary to the preceding  papers such as  \cite{BBMH18,BM18},  we do not need to use pseudo-differential operators techniques  in the description of the aforementioned  asymptotic structures of $\mathfrak{L}_{\varepsilon r}$ and $\mathscr{B}_{\perp}^{-1}\widehat{\mathcal{L}}_{\omega}\mathscr{B}_{\perp}$. In fact, they can be avoided since  all the involved operators can be described through their kernels  and therefore instead of splitting the symbols we simply expand the kernels  as in \eqref{Split-Kernel-Sep} which sounds to be  more appropriate in our context.  \\
	%		
	%		The last comment concerns a technical point in KAM reduction.  Contrary to  the preceding  papers like  \cite{BBMH18,BM18} we do not use  the formalism of  pseudo-differential operators which can be avoided here. Actually, the involved operators are of integral type with slightly singular  kernels and one may perform the estimates directly.  For the reduction of the transport part we need just to conjugate the remainder with the transformation $\mathscr{B}$ which has the advantage to be explicit and  stabilize the kernel structure. Concerning the smoothing effect of the remainder which allows to get the asymptotic of the eigenvalues, it can also be checked using the kernel structures. 
	%		.\\
	
	\noindent $\blacktriangleright$ {{\bf{Step 3.} }{\it Nash-Moser scheme.}
		This is the main  purpose of Section \ref{Sec Nash-Moser} where we 
		construct   zeros for the nonlinear function $\mathcal{F}$ defined in \eqref{PG1-C} for small $\varepsilon$  following  Nash-Moser scheme in the spirit of the papers \cite{BBMH18,BM18}. Let us quickly sketch this scheme.  We build by induction a sequence of approximate solutions $U_{n}$
		$$U_{n+1}=U_{n}+H_{n+1}\quad\textnormal{with}\quad H_{n+1}=-\Pi_{N_{n}}\mathrm{T}_{n}\Pi_{N_{n}}\mathcal{F}(U_{n}).$$
		with $\mathrm{T}_n$ an approximate inverse of  $d_{i,\alpha}\mathcal{F}(U_n)$ constructed in {\bf Step 2}. Thus using Taylor Formula we may write 
		$$\mathcal{F}(U_{n+1})=\Pi_{N_{n}}^{\perp}\mathcal{F}(U_{n})-\Pi_{N_{n}}\big(L_{n}\mathrm{T}_{n}-\hbox{Id}\big)\Pi_{N_{n}}\mathcal{F}(U_{n})+\big(L_{n}\Pi_{n}^{\perp}-\Pi_{N_{n}}^{\perp}L_{n}\big)\mathtt{T}_{n}\Pi_{N_{n}}\mathcal{F}(U_{n})+Q_{n}, $$
		where $Q_n$ is a quadratic functional. Consider the  Cantor set $$\mathcal{A}_{n}^{\gamma}=\displaystyle{\bigcap_{k=0}^{n-1}}\mathcal{G}_{k}(\gamma_{k+1},\tau_{1},\tau_{2},i_{k}),$$ with $\gamma_{n}=\gamma(1+2^{-n})$,  then we  show 
		by induction that 
		\begin{equation}\label{decayFUn}
			\| U_{n}\|_{q,s_{0}+\overline{\sigma}}^{\gamma,\mathcal{O}}\lesssim \varepsilon\gamma^{-1}N_{0}^{q\overline{a}},\quad  \| U_{n}\|_{q,b_{1}+\overline{\sigma}}^{\gamma,\mathcal{O}}\lesssim \varepsilon\gamma^{-1}N_{n-1}^{\mu},\quad \|\mathcal{F}(U_n)\|_{q,s_0}^{\gamma,\mathrm{O}_n^\gamma}\lesssim\varepsilon N_{n-1}^{-a_1}
		\end{equation}
		for a suitable choice of the parameters  $a_1, b_{1},\overline{a}, \mu,\overline\sigma>0$ and $\mathrm{O}_n^\gamma$ is an open enlargement of $\mathcal{A}_n^\gamma$ needed to construct classical extensions to the whole set of parameters $\mathcal{O}$.  Actually, we get a precise statement  in Proposition \ref{Nash-Moser} allowing to  deduce that the sequence  $(U_n)_n$ converges in a strong topology towards a sufficient smooth  profile $(\lambda,\omega)\in\mathcal{O}\mapsto U_\infty(\lambda,\omega)=\big(i_{\infty}(\lambda,\omega),\alpha_{\infty}(\lambda,\omega),(\lambda,\omega)\big)
		$ with 
		$$\forall(\lambda,\omega)\in\mathcal{G}_{\infty}^{\gamma},\quad\mathcal{F}(U_{\infty}(\lambda,\omega))=0,\quad \mathcal{G}_{\infty}^{\gamma}=\bigcap_{n\in\mathbb{N}}\mathcal{A}_{n}^{\gamma}.$$
		Moreover, we get in view of Corollary \ref{Corollary NM} a smooth   function $\lambda\in(\lambda_{0},\lambda_{1})\mapsto(\lambda,\omega(\lambda,\varepsilon))$ with
		\begin{equation}\label{alpha infty-1-intro}
			\omega(\lambda,\varepsilon)=-\omega_{\textnormal{Eq}}(\lambda)+\bar{r}_{\varepsilon}(\lambda),\quad \|\bar{r}_{\varepsilon}\|_{q}^{\gamma,\mathcal{O}}\lesssim\varepsilon\gamma^{-1}N_{0}^{q\overline{a}}
		\end{equation}
		and 
		$$\forall\lambda\in \mathcal{C}_{\infty}^{\varepsilon},\quad \mathcal{F}\big(U_{\infty}(\lambda,\omega(\lambda,\varepsilon))\big)=0\quad\textnormal{with}\quad\alpha_{\infty}\big(\lambda,\omega(\lambda,\varepsilon)\big)=-\omega_{\textnormal{Eq}}(\lambda),
		$$
		where the  Cantor set $\mathcal{C}_{\infty}^{\varepsilon}$ is defined by 
		\begin{equation}\label{def final Cantor set intro}
			\mathcal{C}_{\infty}^{\varepsilon}=\Big\{\lambda\in(\lambda_{0},\lambda_{1})\quad\textnormal{s.t.}\quad\big(\lambda,\omega(\lambda,\varepsilon)\big)\in\mathcal{G}_{\infty}^{\gamma}\Big\}.
		\end{equation}
		This gives solutions to the original equation \eqref{Modi-Eq} provided that $\lambda$ belongs to the final Cantor set $\mathcal{C}_{\infty}^{\varepsilon}$ and the last point to deal with  aims to measure this set.\\
		
		\noindent $\blacktriangleright$ {{\bf{Step 4.} }{\it  Measure estimates.}}
		The measure of the final Cantor set $\mathcal{C}_{\infty}^{\varepsilon}$ will be explored in \mbox{Section \ref{Section measure of the final Cantor set}.} We show in  Proposition \ref{lem-meas-es1} that by fixing $\gamma=\varepsilon^a$ for some small $a$ we get 
		$$\left|\mathcal{C}_{\infty}^{\varepsilon}\right|\geqslant(\lambda_{1}-\lambda_{0})-C\varepsilon^{\delta},$$
		with small $\delta$ connected to the geometry of the Cantor set and the non degeneracy of the equilibrium spectrum.
		There are two main  ingredients to get this result. The first one is   the stability of the intermediate Cantor sets $(\mathcal{A}_{n}^{\gamma})_n$ following from the fast convergence  of  Nash-Moser scheme. However the second one is the transversality property stated in Lemma \ref{lemma Russeman condition for the perturbed frequencies}  used in the spirit of \cite{BBM11} and \cite{R01}. This property will  be first established for the  linear frequencies in Proposition \ref{lemma transversality}, using the analyticity of the eigenvalues and their asymptotic behavior. Then the  
		extension of the transversality assumption to the perturbed frequencies is done using perturbative arguments together with the asymptotic description of the approximate  eigenvalues detailed  in \eqref{asy-z1}, \eqref{uniform estimate rjinfty} and \eqref{uniform estimate r1}. \\We emphasize  that the transversality is strongly related to the non-degeneracy of the eigenvalues in the sense of 
		the Definition \ref{def-degenerate} . For instance, we show that  the curve $\lambda\in[\lambda_0,\lambda_1] \mapsto\big(\Omega_{j_1} (\lambda),...,\Omega_{j_d} (\lambda)\big)$ is  not contained in any vectorial plane, that is, if there exists a constant vector $c=(c_1,..,c_d)$ such that 
		$$
		\forall \lambda\in[\lambda_0,\lambda_1] ,\quad \sum_{j=1}^d c_k \Omega_{j_k}(\lambda)=0,
		$$ 
		then $c=0.$ This is proved in Lemma \ref{Non-degenracy1} and  follows from the asymptotic of the eigenvalues for large values of $\lambda$ according to the law \eqref{asymptotic expansion of large argument} combined with the invertibility of Vandermonde matrices.		
				
			%We refer for instance to \cite{BBM16,BBMH18} for examples. This is not the case here.
			\section{Hamiltonian formulation of the patch motion}
			In this section we shall set up  the contour dynamics equation governing the patch motion. A particular attention will be focused on the vortex patch equation in the polar coordinates system. We shall see that the Hamiltonian structure still survives at the level of the patch dynamics, which is the  starting point towards    the construction  of  quasi-periodic solutions. 
			\subsection{Contour dynamics equation in polar coordinates}\label{sec contour}
			 Here and in the sequel, we identify $\mathbb{C}$ with $\mathbb{R}^{2}$ equipped with the canonical Euclidean structure  through the standard inner  product defined for all $z_{1}=a_{1}+ib_{1},\,z_{2}=a_{2}+ib_{2}\in\mathbb{C}$ by
			\begin{equation}\label{definition of the scalar product on R2}
				z_{1}\cdot z_{2}:=\langle z_{1},z_{2}\rangle_{\mathbb{R}^{2}}=\mbox{Re}\left(z_{1}\overline{z_{2}}\right)=a_{1}a_{2}+b_{1}b_{2}.%=\mbox{Im}\left(iz_{1}\overline{z_{2}}\right)=a_{1}a_{2}+b_{1}b_{2}.
			\end{equation}
			The planar set  $\mathbb{D}$ stands for the open   unit disc of $\mathbb{R}^{2}$ and  the Rankine vortex   $\mathbf{1}_{\mathbb{D}}$ (actually any radial function) is a stationary solution to $(\textnormal{QGSW})_{\lambda}$. To look for    ordered structure like periodic or quasi-periodic vortex patches $t\mapsto\mathbf{1}_{D_{t}}$ around this equilibrium state, we find it  convenient to consider a polar parametrization of the boundary  
			\begin{align}\label{Param-curv}
				z:\begin{array}[t]{rcl}
					\mathbb{R}_{+}\times[0,2\pi] & \mapsto & \mathbb{C}\\
					(t,\theta) & \mapsto & z(t,\theta)=\left(1+2r(t,\theta)\right)^{\frac{1}{2}}e^{i\theta}.
				\end{array}
			\end{align}
			Here $r$ is the radial deformation of the patch which is small, namely $|r(t,\theta)|\ll1.$ Taking $r=0$ gives a parametrization of the  unit circle $\mathbb{T}.$  We shall introduce the new symplectic unknown \begin{equation}\label{definition of R}
				R(t,\theta)=\left(1+2r(t,\theta)\right)^{\frac{1}{2}}.
			\end{equation}
			which will be useful to write down the equations into the Hamiltonian form. In what follows we want  to explicit  the contour dynamics equation with the polar coordinates.  It is a classical fact, see for instance \cite{HMV13,HMV15}, that the particles on the boundary move with the flow and remain at the boundary  and therefore in the smooth case one has   
			$$
			\left[\partial_{t}z(t,\theta)-\mathbf{v}(t,z(t,\theta))\right]\cdot\mathbf{n}(t,z(t,\theta))=0,
			$$
			where $\mathbf{n}(t,z(t,\theta))$ is the outward normal vector to the boundary $\partial D_{t}$ of $D_{t}$ at the point $z(t,\theta)$. Since  one has, up to a real constant of renormalization, $\mathbf{n}(t,z(t,\theta))=-\ii\partial_{\theta}z(t,\theta)$, then  we find the complex form of the contour dynamics motion, 
			\begin{equation}\label{complex vortex patch equation}
				\mbox{Im}\left(\left[\partial_{t}z(t,\theta)-\mathbf{v}(t,z(t,\theta))\right]\overline{\partial_{\theta}z(t,\theta)}\right)=0.
			\end{equation}
			In order to transform it into   a nonlinear  PDE, we  need  to recover  the velocity field $\mathbf{v}$ from the patch parametrization. To do so, recall that  $\mathbf{v}=\nabla^{\perp}\mathbf{\Psi}$ where $\mathbf{\Psi}$ is the stream function associated to the vorticity and governed by  Helmholtz equation,
			$$(\Delta-\lambda^{2})\mathbf{\Psi}(t,\cdot)=\mathbf{1}_{D_{t}}.$$
			To invert this operator we shall make appeal to the Green function $T_{\lambda}$ solution of the equation
			$$(-\Delta+\lambda^{2})T_{\lambda}=\delta_{0}\quad\mbox{in }\mathcal{S}'(\mathbb{R}^{2}).$$
			Using the Fourier transform yields
			$$\forall\xi\in\mathbb{R}^{2},\quad\widehat{T_{\lambda}}(\xi)=\frac{1}{|\xi|^{2}+\lambda^{2}}\cdot
			$$
			Thus by Fourier inversion  theorem and using a scaling argument, we find 
			$$T_{\lambda}(z)=T_{1}(\lambda z)\quad\mbox{with}\quad T_{1}(z)=\frac{1}{4\pi^{2}}\int_{\mathbb{R}^{2}}\frac{e^{iz\cdot\xi}}{1+|\xi|^{2}}d\xi.$$
			Applying  a polar change of variables gives
			$$T_{1}(z)= \displaystyle\frac{1}{4\pi^{2}}\int_{0}^{\infty}\frac{r}{1+r^{2}}\int_{0}^{2\pi}\cos(|z|r\cos(\theta))d\theta dr.$$
			Simple arguments based on the symmetry of trigonometric functions allow to get the identity
						$$
						\int_{0}^{2\pi}\cos(|z|r\cos(\theta))d\theta=2\int_{0}^{\pi}\cos(|z|r\sin(\theta))d\theta.
						$$
			Consequently, we  get in view of \eqref{Besse-repr}
			$$T_{1}(z)= \frac{1}{2\pi}\int_{0}^{\infty}\frac{rJ_{0}(|z|r)}{1+r^{2}}dr,
			$$
			where $J_n$ denotes the Bessel function.
			 Applying \eqref{integral representation for Knu} with $\nu=\mu=0,$ $a=1$ and $b=|z|$, we finally deduce the representation
			$$T_{1}(z)=\frac{1}{2\pi}K_{0}(|z|).$$
			Therefore one gets the formula,
			\begin{equation}\label{velocity potential}
				\mathbf{\Psi}(t,z)=-\frac{1}{2\pi}\int_{\mathbb{R}^{2}}K_{0}(\lambda|z-\xi|)\mathbf{1}_{D_{t}}(\xi)dA(\xi),
			\end{equation}
			where $dA$ denotes the planar Lebesgue measure. To get get explicit form of the velocity in terms of the patch boundary  we shall use  the complex version of Stokes theorem 
			\begin{equation}\label{complex version of Stokes theorem}
				2\ii\int_{D}\partial_{\overline{\xi}}f(\xi,\overline{\xi})dA(\xi)=\int_{\partial D}f(\xi,\overline{\xi})d\xi.
			\end{equation}
			In view of   $\mathbf{v}(t,z)=2\ii\partial_{\overline{z}}\mathbf{\Psi}(t,z)$,  one deduces that
			\begin{equation}\label{velocity}
				\mathbf{v}(t,z)=\frac{1}{2\pi}\int_{\partial D_{t}}K_{0}(\lambda|z-\xi|)d\xi.
			\end{equation}
			Next we intend to write down the boundary motion in terms of the contour dynamics.	 First, from   the polar parametrization, it is easy to check from \eqref{Param-curv} that
			$$
			\mbox{Im}\left(\partial_{t}z(t,\theta)\overline{\partial_{\theta}z(t,\theta)}\right)=-\partial_{t}r(t,\theta).
			$$
			On the other hand, using \eqref{velocity} and \eqref{symmetry Bessel}, we infer
			$$\mbox{Im}\left(\mathbf{v}(t,z(t,\theta))\overline{\partial_{\theta}z(t,\theta)}\right)=\displaystyle\int_{\mathbb{T}}K_{0}\left(\lambda|z(t,\theta)-z(t,\eta)|\right)\mbox{Im}\left(\partial_{\eta}z(t,\eta)\overline{\partial_{\theta}z(t,\theta)}\right)d\eta.$$
			Here and throughout this paper, we shall work with the following notation 
			\begin{equation}\label{Convention1}
			\displaystyle\int_{\mathbb{T}}f(\eta)d\eta:=\tfrac{1}{2\pi}\int_{0}^{2\pi} f(\eta)d\eta.
			\end{equation}
			Next we observe that,
			\begin{align*}\mbox{Im}\left(\partial_{\eta}z(t,\eta)\overline{\partial_{\theta}z(t,\theta)}\right)&=\partial^2_{\theta\eta}\mbox{Im}\left(z(t,\eta)\overline{z(t,\theta)}\right)\\
				&=\partial_{\theta\eta}^{2}\Big(R(t,\eta)R(t,\theta)\sin(\eta-\theta)\Big).
			\end{align*}
			Thus, by setting
			\begin{equation}\label{definition of Ar}
				A_{r}(t,\theta,\eta):=|z(t,\theta)-z(t,\eta)|=|R(t,\theta)e^{i\theta}-R(t,\eta)e^{i\eta}|
			\end{equation}
			and
			\begin{equation}\label{definition of Flambda}
				F_{\lambda}[r](t,\theta):=\int_{\mathbb{T}}K_{0}\left(\lambda A_{r}(t,\theta,\eta)\right)\partial_{\theta\eta}^{2}\Big(R(t,\eta)R(t,\theta)\sin(\eta-\theta)\Big)d\eta,
			\end{equation}
			we get  the vortex patch equation in the polar coordinates
			\begin{equation}\label{equationL1}
				\partial_{t}r(t,\theta)+F_{\lambda}[r](t,\theta)=0.
			\end{equation}
			Now, we fix a parameter $\Omega$ that will be used later to get rid of trivial resonances, and  we shall  look for solutions % in the rotating frame, namely 
  in   the form
\begin{equation}
r(t, \theta)=\tilde r(t,\theta+\Omega t).
\end{equation}  
%Differentiating the last expression with respect to $t$ gives
% $$
%\partial_t w(t, \theta) = 
%- \ii \Omega e^{ -\ii \Omega t} z(t,\theta) +
%e^{ -\ii \Omega t} \partial_t z(t,\theta)   \, , 
%$$
Then elementary change of variables applied with   \eqref{definition of Flambda} show that  
\begin{equation}\label{u-rot}
F_{\lambda}[\tilde r](t,\theta{+\Omega t})= 
F_{\lambda}[ r](t,\theta) \, . 
\end{equation}
Thus, the equation  \eqref{equationL1} becomes (to alleviate the notation we simply use $r$ instead of $\tilde{r}$)
\begin{equation}\label{equation}
\partial_{t} r(t,\theta)+\Omega\partial_{\theta}  r(t,\theta)
 +F_{\lambda}[ r](t,\theta)=0,
\end{equation}
	which is a nonlinear and nonlocal transport PDE. 
			To fix the terminology, we mean by a  time quasi-periodic solution of \eqref{equation}, a solution in  the form  
			$$r(t,\theta)=\widehat{r}(\omega t,\theta),$$
			where  $\widehat{r}:\,(\varphi,\theta)\in\mathbb{T}^{d+1}\mapsto \,\widehat{r}(\varphi,\theta)\in \mathbb{R},$ $\omega\in\mathbb{R}^{d}$ and $d\in\mathbb{N}^{*}$. Hence in this setting,  the equation \eqref{equation} becomes  
			\begin{equation*}
				\omega\cdot\partial_{\varphi}\widehat{r}(\varphi,\theta)+\Omega\partial_{\theta}\widehat{r}(\varphi,\theta)+F_{\lambda}[\widehat{r}](\varphi,\theta)=0.
			\end{equation*}
			In the sequel, we shall  alleviate the notation and denote $\widehat{r}$ simply by $r$ and the foregoing equation writes
			\begin{equation}\label{equation with quasi-periodic ansatz}
				\forall (\varphi,\theta)\in\mathbb{T}^{d+1},\quad\omega\cdot\partial_{\varphi}{r}(\varphi,\theta)+\Omega\partial_{\theta}{r}(\varphi,\theta)+F_{\lambda}[{r}](\varphi,\theta)=0.
			\end{equation}
			
			\subsection{Hamiltonian structure}
			We  now move  to a new consideration related to the analysis of the  Hamiltonian structure behind the transport equation  \eqref{equation}. This structure sounds to be  essential if one wants to explore quasi-periodic solutions  near Rankine vortices. Notice that it is a classical fact that  incompressible active scalar equations such as 2D Euler equations are  Hamiltonian and as we shall see in this section,  we can find a suitable  interpretation of  this property   at the level of the contour dynamics equations which is a stronger reformulation.  
			
			\subsubsection{Hamiltonian reformulation}
			We consider the  kinetic energy and the angular impulse associated to the patch $\omega(t)={\bf{1}}_{D_t}$ and defined by 
			\begin{equation}\label{definition of the pseudo-energy and angular momentum}
				E(t)=\frac{1}{2\pi}\int_{D_{t}}\mathbf{\Psi}(t,z)dA(z)\quad\mbox{ and }\quad J(t)=\frac{1}{2\pi}\int_{D_{t}}|z|^{2}dA(z),
			\end{equation}
			where the stream function $\mathbf{\Psi}$ is defined according to  \eqref{velocity potential}. The following result dealing with the time  conservation of the preceding quantities  is classical and can be proved in a similar way to Euler equations.
			\begin{lem}
				The kinetic energy $E$ and the angular impulse $J$ are conserved during the motion,
					$$\frac{d E(t)}{dt}=0=\frac{dJ(t)}{dt}\cdot$$
			\end{lem}
			In what follows we shall state the main result of this section on the Hamiltonian structure governing   the equation \eqref{equation}.	
					\begin{prop}\label{proposition Hamiltonian formulation of the equation}
				{The equation \eqref{equation} is  Hamiltonian and takes  the form 
					\begin{equation}\label{Hamiltonian formulation of the equation}
						\partial_{t}r=\partial_{\theta}\nabla H(r),
					\end{equation}
					where $\nabla$ is the $L^{2}(\mathbb{T}_{\theta})$-gradient with respect to the $L^{2}(\mathbb{T}_{\theta})$-normalized scalar product defined by 
					$$\langle \rho_{1},\rho_{2}\rangle_{L^{2}(\mathbb{T})}=\int_{\mathbb{T}}\rho_{1}(\theta)\rho_{2}(\theta)d\theta$$
					and  the hamiltonian $H$ is defined by 
					$$H(r)=\tfrac{1}{2}\left(E(r)-\Omega J(r)\right).$$}
					In particular, we get the conservation of the average, that is, 
					\begin{equation}\label{pres avrg}
						\frac{d}{dt}\displaystyle\int_{\mathbb{T}}r(t,\theta)d\theta=0.
					\end{equation}
			\end{prop}
			\begin{proof}
				$\blacktriangleright$ Using  Stokes formula \eqref{complex version of Stokes theorem}, we may write
				$$J(r)(t)=\tfrac{1}{8\ii \pi}\int_{\partial D_{t}}|\xi|^{2}\overline{\xi}d\xi.
				$$
				Then from   the parametrization detailed in \eqref{Param-curv} one gets easily 				\begin{align*}
					J(r)(t) & =  \displaystyle\tfrac{1}{4\ii}\int_{\mathbb{T}}|z(t,\theta)|^{2}\overline{z(t,\theta)}\partial_{\theta}z(t,\theta)d\theta\\
					& =  \displaystyle\tfrac{1}{16\ii}\int_{\mathbb{T}}\partial_{\theta}\left(R^{4}(t,\theta)\right)d\theta+\tfrac{1}{4}\int_{\mathbb{T}}R^{4}(t,\theta)d\theta\\
					& =  \displaystyle\tfrac{1}{4}\int_{\mathbb{T}}R^{4}(t,\theta)d\theta.
				\end{align*}
				Consequently,
				\begin{equation}\label{expression for the angular momentum}
					J(r)(t)=\tfrac{1}{4}\int_{\mathbb{T}}\left(1+2r(t,\theta)\right)^{2}d\theta.
				\end{equation}
				Differentiating in $r$ one gets for  $\rho\in L^{2}(\mathbb{T})$
				$$\langle\nabla J(r),\rho\rangle_{L^{2}(\mathbb{T})}(t)=\int_{\mathbb{T}}(1+2r(t,\theta))\rho(\theta)d\theta,\quad \mbox{ i.e. }\quad \nabla J(r)=1+2r.$$
				It follows that 
				\begin{equation}\label{link J and r}
					\tfrac12\Omega\partial_{\theta}\nabla J(r)=\Omega\partial_{\theta}r.
				\end{equation}
				$\blacktriangleright$ Next, we shall   compute the G\^ateaux derivative of $E$ in a given direction $\rho\in L^{2}(\mathbb{T}).$ We point out that the computations done below are formal but they can be justified rigorously in a classical way. The first step is to  express the energy
				$$E(t)=\tfrac{1}{2\pi}\int_{D_{t}}\mathbf{\Psi}(t,z)dA(z)$$
				in terms of the boundary parametrization of  $\partial D_{t}.$ According to 
				Stokes theorem \eqref{complex version of Stokes theorem} we have 
				\begin{align*}
					\mathbf{\Psi}(t,z) & =  \displaystyle-\tfrac{1}{2\pi}\int_{D_{t}}K_{0}(\lambda|\xi-z|)dA(\xi)\\
					& = \displaystyle\tfrac{1}{\pi\lambda^{2}}\bigintsss_{D_{t}}\partial_{\overline{\xi}}\left(\tfrac{(\overline{\xi}-\overline{z})\left[\lambda|\xi-z|K_{1}(\lambda|\xi-z|)-1\right]}{|\xi-z|^{2}}\right)dA(\xi)\\
					& =  \displaystyle\tfrac{1}{2\ii\pi\lambda^{2}}\bigintsss_{\partial D_{t}}\tfrac{(\overline{\xi}-\overline{z})\left[\lambda|\xi-z|K_{1}(\lambda|\xi-z|)-1\right]}{|\xi-z|^{2}}d\xi.
				\end{align*}
				To prove the second equality above, it suffices to find  an anti-derivative of $K_{0}(\lambda|\xi-z|)$ with respect to $\overline{\xi}.$ We shall search it in  the form 
				$$(\overline{\xi}-\overline{z})f(\lambda|\xi-z|).$$
				Then we should get
				$$K_{0}(\lambda|\xi-z|)=\partial_{\overline{\xi}}\left((\overline{\xi}-\overline{z})f(\lambda|\xi-z|)\right)=f(\lambda|\xi-z|)+\tfrac{\lambda|\xi-z|}{2}f'(\lambda|\xi-z|).$$
				Hence $f$ is a  solution on $\mathbb{R}_{+}^{*}$ of the ordinary differential equation 
				$$\tfrac12{x}f^\prime(x)+f(x)=K_{0}(x),\quad\mbox{ i.e. }\quad(x^{2}f(x))'=2xK_{0}(x).$$
				Using \eqref{Bessel and anti-derivatives}, we obtain
					$$f(x)=-\tfrac{2xK_{1}(x)+C}{x^{2}},$$
					where $C$ is a constant to be fixed  so that the integral converges. Using \eqref{power series Kj}, one has on the real line
					$$K_{1}(x)\underset{x\rightarrow 0}{=}\tfrac{1}{x}+\tfrac{x}{2}\log\left(\tfrac{x}{2}\right)+o\left(x\log\left(\tfrac{x}{2}\right)\right),$$
					so that
					$$xK_{1}(x)\underset{x\rightarrow 0}{=}1+\tfrac{x^{2}}{2}\log\left(\tfrac{x}{2}\right)+o\left(x^{2}\log\left(\tfrac{x}{2}\right)\right).$$
					Making the choice  $C=-2$ we get 
					\begin{equation}
						f(x)=-\tfrac{2\left(xK_{1}(x)-1\right)}{x^{2}},
					\end{equation}
					which behaves like a logarithm near $0$ and thus it  is integrable. Therefore using the parametrization \eqref{Param-curv} we find $$\mathbf{\Psi}(t,z(t,\theta))=\tfrac{1}{\ii\lambda^{2}}\bigintsss_{\mathbb{T}}\frac{(\overline{z}(t,\eta)-\overline{z}(t,\theta))\left[\lambda|z(t,\theta)-z(t,\eta)|K_{1}\left(\lambda|z(t,\eta)-z(t,\theta)|\right)-1\right]}{|z(t,\eta)-z(t,\theta)|^{2}}\partial_{\eta}z(t,\eta)d\eta.$$
				Making appeal to $f$ and removing  the time  dependence, we get
				\begin{equation}\label{expression of psi(z(theta))}
					\mathbf{\Psi}(z(\theta))=\tfrac{\ii}{2}\int_{\mathbb{T}}(\overline{z}(\theta)-\overline{z}(\eta))f(\lambda|z(\theta)-z(\eta)|)\partial_{\eta}z(\eta)d\eta.
				\end{equation}
				At this stage  we need to look for   an anti-derivative with respect to $\overline{z}$ of $\frac{-1}{2}(\overline{\xi}-\overline{z})f(\lambda|\xi-z|)$ in the form 
				$$(\overline{\xi}-\overline{z})^{2}g(\lambda|\xi-z|).$$
				Therefore  we deduce the constraint
				$$\tfrac{-1}{2}(\overline{\xi}-\overline{z})f(\lambda|\xi-z|)=\partial_{\overline{z}}\left((\overline{\xi}-\overline{z})^{2}g(\lambda|\xi-z|)\right)=-(\overline{\xi}-\overline{z})\left(2g(\lambda|\xi-z|)+\tfrac{\lambda|\xi-z|}{2}g'(\lambda|\xi-z|)\right).$$
				Hence, $g$ should be a solution on $\mathbb{R}_{+}^{*}$ of the ordinary differential equation
				\begin{equation}\label{ode for g in hamiltonian formulation}
					\tfrac{x}{2}g'(x)+2g(x)=\tfrac{f(x)}{2},\quad \mbox{ i.e. }\quad (x^{4}g(x))'=x^{3}f(x)=2x-2x^{2}K_{1}(x).
				\end{equation}
				Using once again \eqref{Bessel and anti-derivatives} yields
				$$g(x)=\tfrac{x^{2}+2x^{2}K_{2}(x)+C}{x^{4}},
				$$
				where $C$ is again a constant used to cancel the violent singularity. From \eqref{power series Kj} and \eqref{definition of modified Bessel function of first kind}, one obtains the asymptotic
				$$K_{2}(x)\underset{x\rightarrow 0}{=}\tfrac{2}{x^{2}}-\tfrac{1}{2}+O\left(x^{2}\log(x)\right).$$
				Thus
				$$x^{2}K_{2}(x)\underset{x\rightarrow 0}{=}2-\tfrac{x^{2}}{2}+O(x^{4}\log(x)).$$
				Then by  choosing $C=-4$ we deduce  that the function below 
				$$g(x)=\tfrac{x^{2}+2x^{2}K_{2}(x)-4}{x^{4}}$$
				is integrable. Hence, applying once again  Stokes theorem \eqref{complex version of Stokes theorem}, we infer
				\begin{align*}
					E(r)(t)&=-\tfrac{1}{4\pi^{2}\lambda^{4}}\bigintsss_{\partial D_{t}}\bigintsss_{\partial D_{t}}\frac{(\overline{\xi}-\overline{z})^{2}\left[\lambda^{2}|\xi-z|^{2}\left(1+2K_{2}(\lambda|\xi-z|)\right)-4\right]}{|\xi-z|^{4}}dz d\xi\\
					&=\tfrac{-1}{\lambda^{4}}\bigintsss_{\mathbb{T}}\bigintsss_{\mathbb{T}}\tfrac{(\overline{z}(t,\eta)-\overline{z}(t,\theta))^{2}\left[\lambda|z(t,\eta)-z(t,\theta)|\left(1+2K_{2}(\lambda|z(t,\eta)-z(t,\theta)|)\right)-4\right]}{|z(t,\eta)-z(t,\theta)|^{2}}\partial_{\theta}z(t,\theta)\partial_{\eta}z(t,\eta)d\eta d\theta.
				\end{align*}
				Hence using  $g$ and removing the dependence in time, we find
				\begin{equation}\label{expression for the energy}
					E(r)=-\tfrac{1}{2}\int_{\mathbb{T}}\int_{\mathbb{T}}(\overline{z}(\theta)-\overline{z}(\eta))^{2}g(\lambda|z(\theta)-z(\eta)|)\partial_{\theta}z(\theta)\partial_{\eta}z(\eta)d\theta d\eta.
				\end{equation}
				The next goal is to  compute the derivative of $E$ with respect to $r$ in the direction $\rho,$ which is straightforward 
				\begin{align*}
					\langle\nabla E(r),&\rho\rangle_{L^{2}(\mathbb{T})}  =\displaystyle -\bigintsss_{\mathbb{T}}\bigintsss_{\mathbb{T}}\left(\overline{z}(\theta)-\overline{z}(\eta)\right)g\left(\lambda|z(\theta)-z(\eta)|\right)\left(\tfrac{\rho(\theta)e^{-\ii\theta}}{R(\theta)}-\tfrac{\rho(\eta)e^{-\ii\eta}}{R(\eta)}\right)\partial_{\theta}z(\theta)\partial_{\eta}z(\eta)d\theta d\eta\\
					&-  \displaystyle\tfrac{\lambda}{2}\bigintsss_{\mathbb{T}}\bigintsss_{\mathbb{T}}\tfrac{\left(\overline{z}(\theta)-\overline{z}(\eta)\right)^{2}}{|z(\theta)-z(\eta)|}g'\left(\lambda|z(\theta)-z(\eta)|\right)\tfrac{\rho(\theta)}{R(\theta)}\mbox{Re}\left((z(\theta)-z(\eta))e^{-\ii\theta}\right)\partial_{\theta}z(\theta)\partial_{\eta}z(\eta)d\theta d\eta\\
					& - \displaystyle \tfrac{\lambda}{2}\bigintsss_{\mathbb{T}}\bigintsss_{\mathbb{T}}\tfrac{\left(\overline{z}(\theta)-\overline{z}(\eta)\right)^{2}}{|z(\theta)-z(\eta)|}g'\left(\lambda|z(\theta)-z(\eta)|\right)\tfrac{\rho(\eta)}{R(\eta)}\mbox{Re}\left((z(\eta)-z(\theta))e^{-\ii\eta}\right)\partial_{\theta}z(\theta)\partial_{\eta}z(\eta)d\theta d\eta\\
					&-\displaystyle \tfrac{1}{2}\bigintsss_{\mathbb{T}}\bigintsss_{\mathbb{T}}\left(\overline{z}(\theta)-\overline{z}(\eta)\right)^{2}g\left(\lambda|z(\theta)-z(\eta)|\right)\partial_{\theta}\left(\tfrac{\rho(\theta)e^{\ii\theta}}{R(\theta)}\right)\partial_{\eta}z(\eta)d\theta d\eta\\
					&- \displaystyle \tfrac{1}{2}\bigintsss_{\mathbb{T}}\bigintsss_{\mathbb{T}}\left(\overline{z}(\theta)-\overline{z}(\eta)\right)^{2}g\left(\lambda|z(\theta)-z(\eta)|\right)\partial_{\theta}z(\theta)\partial_{\eta}\left(\tfrac{\rho(\eta)e^{\ii\eta}}{R(\eta)}\right) d\theta d\eta.
				\end{align*}
				By exchanging in the double integral  $\theta$ and $\eta$, we deduce
				\begin{align*}
					\langle\nabla E(r),&\rho\rangle_{L^{2}(\mathbb{T})} = \displaystyle -2\bigintsss_{\mathbb{T}}\bigintsss_{\mathbb{T}}\left(\overline{z}(\theta)-\overline{z}(\eta)\right)g\left(\lambda|z(\theta)-z(\eta)|\right)\tfrac{\rho(\theta)e^{-\ii\theta}}{R(\theta)}\partial_{\theta}z(\theta)\partial_{\eta}z(\eta)d\theta d\eta\\
					&  \displaystyle -\lambda\bigintsss_{\mathbb{T}}\bigintsss_{\mathbb{T}}\tfrac{\left(\overline{z}(\theta)-\overline{z}(\eta)\right)^{2}}{|z(\theta)-z(\eta)|}g'\left(\lambda|z(\theta)-z(\eta)|\right)\tfrac{\rho(\theta)}{R(\theta)}\mbox{Re}\left((z(\theta)-z(\eta))e^{-\ii\theta}\right)\partial_{\theta}z(\theta)\partial_{\eta}z(\eta)d\theta d\eta\\
					& \displaystyle -\bigintsss_{\mathbb{T}}\bigintsss_{\mathbb{T}}\left(\overline{z}(\theta)-\overline{z}(\eta)\right)^{2}g\left(\lambda|z(\theta)-z(\eta)|\right)\partial_{\theta}\left(\tfrac{\rho(\theta)e^{\ii\theta}}{R(\theta)}\right)\partial_{\eta}z(\eta)d\theta d\eta.
				\end{align*}
				An integration by parts in the last integral leads to
				$$\begin{array}{l}
					\displaystyle -\bigintsss_{\mathbb{T}}\bigintsss_{\mathbb{T}}\left(\overline{z}(\theta)-\overline{z}(\eta)\right)^{2}g\left(\lambda|z(\theta)-z(\eta)|\right)\partial_{\theta}\left(\tfrac{\rho(\theta)e^{\ii\theta}}{R(\theta)}\right)\partial_{\eta}z(\eta)d\theta d\eta\\
					=\displaystyle 2\bigintsss_{\mathbb{T}}\bigintsss_{\mathbb{T}}(\overline{z}(\theta)-\overline{z}(\eta))g(\lambda|z(\theta)-z(\eta)|)\tfrac{\rho(\theta)e^{\ii\theta}}{R(\theta)}\partial_{\theta}\overline{z}(\theta)\partial_{\eta}z(\eta)d\theta d\eta\\
					\mbox{\hspace{0.5cm}}\displaystyle+\lambda\bigintsss_{\mathbb{T}}\bigintsss_{\mathbb{T}}\tfrac{(\overline{z}(\theta)-\overline{z}(\eta))^{2}}{|z(\theta)-z(\eta)|}g'(\lambda|z(\theta)-z(\eta)|)\tfrac{\rho(\theta)e^{\ii\theta}}{R(\theta)}\mbox{Re}\left[(z(\theta)-z(\eta))\partial_{\theta}\overline{z}(\theta)\right]\partial_{\eta}z(\eta)d\theta d\eta.
				\end{array}$$
				Using  the identities
				$$e^{\ii\theta}\partial_{\theta}\overline{z}(\theta)-e^{-\ii\theta}\partial_{\theta}z(\theta)=-2\ii R(\theta)$$
				and 
				$$\mbox{Re}\left[(z(\theta)-z(\eta))\partial_{\theta}\overline{z}(\theta)\right]e^{\ii\theta}-\partial_{\theta}z(\theta)\mbox{Re}\left[(z(\theta)-z(\eta))e^{-\ii\theta}\right]=-\ii R(\theta)(z(\theta)-z(\eta)),$$
				we infer
				$$\begin{array}{rcl}
					\langle\nabla E(r),\rho\rangle_{L^{2}(\mathbb{T})} & = & \displaystyle\tfrac{4}{\ii}\int_{\mathbb{T}}\int_{\mathbb{T}}(\overline{z}(\theta)-\overline{z}(\eta))g(\lambda|z(\theta)-z(\eta)|)\partial_{\eta}z(\eta)\rho(\theta)d\theta d\eta\\
					& & \displaystyle+\tfrac{\lambda}{\ii}\int_{\mathbb{T}}\int_{\mathbb{T}}(\overline{z}(\theta)-\overline{z}(\eta))|z(\theta)-z(\eta)|g'(\lambda|z(\theta)-z(\eta)|)\partial_{\eta}z(\eta)\rho(\theta)d\theta d\eta.
				\end{array}$$
				Applying  \eqref{ode for g in hamiltonian formulation}, we find
				$$\langle\nabla E(r),\rho\rangle_{L^{2}(\mathbb{T})}=\displaystyle\tfrac{1}{\ii}\int_{\mathbb{T}}\int_{\mathbb{T}}(\overline{z}(\theta)-\overline{z}(\eta))f(\lambda|z(\theta)-z(\eta)|)\partial_{\eta}z(\eta)\rho(\theta)d\theta d\eta,$$
				which implies by virtue of   \eqref{expression of psi(z(theta))} 
				$$\nabla E(r)=\tfrac{1}{\ii}\int_{\mathbb{T}}(\overline{z}(\theta)-\overline{z}(\eta))f(\lambda|z(\theta)-z(\eta)|)\partial_{\eta}z(\eta)d\eta=-2\mathbf{\Psi}(z(\theta)).$$
				Now, using the complex notation we deduce that
				\begin{align*}
					\partial_{\theta}\mathbf{\Psi}(z(\theta)) & = \nabla\mathbf{\Psi}(z(\theta))\cdot \partial_{\theta}z(\theta)\\
					%2\mbox{Re}\left(\partial_{z}\mathbf{\Psi}(z(\theta))\partial_{\theta}z(\theta)\right)\\
					%& = & 2\mbox{Im}\left(i\overline{\partial_{z}\mathbf{\Psi}(z(\theta))}\overline{\partial_{\theta}z(\theta)}\right)\\
					%& = & \mbox{Im}\left(2i\partial_{\overline{z}}\mathbf{\Psi}(z(\theta))\overline{\partial_{\theta}z(\theta)}\right)\\
					& =  \mbox{Im}\left(\mathbf{v}(z(\theta))\overline{\partial_{\theta}z(\theta)}\right)\\
					& =  F_{\lambda}[r](\theta),
				\end{align*}
				where we used \eqref{definition of the scalar product on R2} and the facts that $\nabla^{\perp}\mathbf{\Psi}=\mathbf{v}$ and $\mathbf{\Psi}$ is real-valued. Recall  that  the functional $F_{\lambda}[r]$ was introduced in \eqref{definition of Flambda}.
				Hence 
				$$\partial_{\theta}\nabla E(r)=-2\partial_{\theta}\mathbf{\Psi}(z(\theta))=-2F_{\lambda}[r](\theta).$$
				Finally we get 
				\begin{equation}\label{link energy and Flambda}
					\tfrac{1}{2}\partial_{\theta}\nabla E(r)=-F_{\lambda}[r](\theta).
				\end{equation}
				The conservation of the average is easy to check from the Hamiltonian equation. Therefore the proof of Proposition \ref{proposition Hamiltonian formulation of the equation} is achieved.
			\end{proof}
			\subsubsection{Reversibility}
			The main concern is to investigate the reversibility of the Hamiltonian \mbox{equation \eqref{Hamiltonian formulation of the equation}.} This property will be used in a crucial way to fix the symmetry in the function spaces and by this way remove from the phase space the trivial resonances. For more details we refer to Section \ref{Functionspsaces} and Section \ref{sec Ham-tool}. To define the reversibility, we shall introduce the involution $\mathscr{S}$  
			\begin{equation}\label{definition of the involution mathcal S}
				(\mathscr{S}r)(\theta):=r(-\theta),
			\end{equation}
			which satisfies the obvious properties
			\begin{equation}\label{properties of the involution mathcal S}
				\mathscr{S}^{2}=\hbox{Id}\qquad%,\quad \mathscr{S}^{T}=\mathscr{S}\quad\mbox{ and }\quad\partial_{\theta}
				\hbox{and}\qquad \partial_{\theta}\circ\mathscr{S}=-\mathscr{S}\circ\partial_{\theta}.
			\end{equation}
			The following elementary result is useful and can be easily checked from the structure of the Hamiltonian. Actually, it suffices to make changes of variables. 
			\begin{lem}\label{lem rev H-X}
				{The Hamiltonian $H$ and its associated vector field $X:=\partial_{\theta}\nabla H$ satisfy
					$$H\circ\mathscr{S}=H\quad \mbox{ and }\quad X\circ\mathscr{S}=-\mathscr{S}\circ X.$$}
			\end{lem}
			\section{Linearization and frequencies structure}
			This section is devoted to  some aspects of  the linearized operator associated to the evolution equation \eqref{equation} or its Hamiltonian version  \eqref{Hamiltonian formulation of the equation}.  We shall in particular compute it at any state close to the equilibrium and reveal some of its main general  feature. As we shall see,  the radial shape is very special  and gives rise to  a  Fourier multiplier and thus the spectral properties follow immediately.  This latter case   serves as a toy model to check the emergence of  quasi-periodic solutions at the linear level provided that the Rossby radius $\lambda$ belongs to a Cantor set, see Proposition \ref{lemma sol Eq} . However, around this ideal state the situation is roughly uncontrolled and the operator is no longer diagonal and    its spectral study is extremely  delicate due to resonances that prevent to diagonalize the operator.  To deal with this problem  we will implement some important tools borrowed from KAM theory as we shall see in  Section \ref{sec red lin op}.     
			\subsection{Linearized operator}
			The main goal of this section to compute the differential of the nonlinear operator in \eqref{equation} for any small state $r$. The computations will be conducted at a  formal level by simply computing Gateaux derivatives which are related  to  Frechet derivatives. This formal part can be justified rigorously in a classical way  for  the  suitable functional setting fixed in Section \ref{Functionspsaces}.		
			%
			%\subsubsection{General form}
				\subsubsection{The general form}
			In what follows we shall derive a formula for the linearized operator associated to the equation \eqref{Hamiltonian formulation of the equation}. We shall see that it can be split into a transport part with variable coefficients  and a nonlocal operator of order zero. More precisiely, we shall establish the following lemma.
			\begin{lem}\label{lemma general form of the linearized operator}
				{The linearized  equation of \eqref{Hamiltonian formulation of the equation} at a given small state $r$ is given by the time-dependent linear  Hamiltonian equation,
					$$\partial_{t}\rho(t,\theta)=\partial_{\theta}\Big(-V_{r}(\lambda,t,\theta)\rho(t,\theta)+\mathbf{L}_{r}\rho(\lambda,t,\theta)\Big),$$
					where $V_{r}$ is a scalar function  defined by
					\begin{equation}\label{definition of Vr}
						V_{r}(\lambda,t,\theta)=\Omega+\tfrac{1}{R(t,\theta)}\int_{\mathbb{T}}K_{0}\left(\lambda A_{r}(t,\theta,\eta)\right)\partial_{\eta}\left(R(t,\eta)\sin(\eta-\theta)\right)d\eta
					\end{equation} and $\mathbf{L}_{r}$ is a non-local operator given by
					\begin{equation}\label{definition of mathbfLr}
						\mathbf{L}_{r}(\rho)(\lambda,t,\theta)=\int_{\mathbb{T}}K_{0}\big(\lambda A_{r}(t,\theta,\eta)\big)\rho(t,\eta)d\eta.
					\end{equation}
					We recall that $K_{0}$, $A_{r}$ and $R$ are  defined by \eqref{explicit form for K0}, \eqref{definition of Ar} and \eqref{definition of R}, respectively. \\
					Moreover, if $r(-t,-\theta)=r(t,\theta)$, then
					\begin{equation}\label{symmetry for Vr}
						V_{r}(\lambda,-t,-\theta)=V_{r}(\lambda,t,\theta).
				\end{equation}}
			\end{lem}
			\begin{proof} Throughout  the proof, we shall remove the  time  dependency  of the involved  quantities  except  when it is relevant to keep it. The computations of  the G\^ateaux derivative of $F_{\lambda}$ defined by \eqref{definition of Flambda} at a point $r$ in the direction $\rho$ are straightforward and standard and we shall only sketch the main lines. Notice that the functional $F_{\lambda}$ is smooth in a suitable functional setting and therefore its  differential  should be recovered from   its G\^ateaux derivative.  First, we observe that the   function $A_{r}$ defined in \eqref{definition of Ar} can be written in the form
				\begin{align}\label{formula}
					\nonumber A_{r}(\theta,\eta)&=\left(R^{2}(\theta)+R^{2}(\eta)-2R(\theta)R(\eta)\cos(\eta-\theta)\right)^{\frac{1}{2}}\\
					&=\left((R(\theta)-R(\eta))^{2}+4R(\theta)R(\eta)\sin^{2}\left(\tfrac{\eta-\theta}{2}\right)\right)^{\frac{1}{2}}.
				\end{align}
				This identity \eqref{formula} will be of constant use in the sequel. Second, after straightforward computations, we obtain \mbox{from \eqref{definition of Flambda},}
				\begin{align*}
					d_{r}F_{\lambda}[r](\rho)&=\partial_\tau F_{\lambda}[r+\tau\rho]_{|_{\tau=0}}\\
					&=\mathcal{I}_{1}+\mathcal{I}_{2}+\mathcal{I}_{3}+\mathcal{I}_{4},
				\end{align*}
				where
				\begin{align*}
					\mathcal{I}_{1}&:=\lambda\rho(\theta)\displaystyle\bigintssss_{\mathbb{T}}B_{r}(\theta,\eta)K_{0}'\left(\lambda A_{r}(\theta,\eta)\right)\partial^{2}_{\theta\eta}\left(R(\theta)R(\eta)\sin(\eta-\theta)\right)d\eta,\\
					\mathcal{I}_{2}&:=\lambda\displaystyle\bigintssss_{\mathbb{T}}\rho(\eta)B_{r}(\eta,\theta)K_{0}'\left(\lambda A_{r}(\theta,\eta)\right)\partial^{2}_{\theta\eta}\left(R(\theta)R(\eta)\sin(\eta-\theta)\right)d\eta,\\
					\mathcal{I}_{3}&:=\displaystyle\bigintssss_{\mathbb{T}}K_{0}\left(\lambda A_{r}(\theta,\eta)\right)\partial^{2}_{\theta\eta}\left(\rho(\theta)\tfrac{R(\eta)\sin(\eta-\theta)}{R(\theta)}\right)d\eta,\\
					\mathcal{I}_{4}&:=\displaystyle\bigintssss_{\mathbb{T}}K_{0}\left(\lambda A_{r}(\theta,\eta)\right)\partial^{2}_{\theta\eta}\left(\rho(\eta)\tfrac{R(\theta)\sin(\eta-\theta)}{R(\eta)}\right)d\eta,
				\end{align*}
			with 
				\begin{equation}\label{definition of Br}
					B_{r}(\theta,\eta):=\frac{R(\theta)-R(\eta)\cos(\eta-\theta)}{R(\theta)A_{r}(\theta,\eta)}\cdot
				\end{equation}
				Next, we shall compute  $\mathcal{I}_{1}+\mathcal{I}_{3}.$ To do that, we split $\mathcal{I}_{3}$ into two terms as follows,
				\begin{align*}
					\mathcal{I}_{3}&={\partial_{\theta}\rho(\theta)}\bigintssss_{\mathbb{T}}K_{0}\left(\lambda A_{r}(\theta,\eta)\right)\partial_{\eta}\left(\tfrac{R(\eta)\sin(\eta-\theta)}{R(\theta)}\right)d\eta\\
					&\quad+\rho(\theta)\bigintssss_{\mathbb{T}}K_{0}\left(\lambda A_{r}(\theta,\eta)\right)\partial_{\theta\eta}^{2}\left(\tfrac{R(\eta)\sin(\eta-\theta)}{R(\theta)}\right)d\eta\\
					&:=\partial_{\theta}\rho(\theta)\overline{V}_{r}(\lambda,\theta)+\rho(\theta)\overline{\mathcal{I}}_{3}.
				\end{align*}
				An integration by parts in $\overline{\mathcal{I}}_{3}$ allows to get,
				$$\overline{\mathcal{I}}_{3}=-\lambda\displaystyle\bigintssss_{\mathbb{T}}\partial_{\eta}A_{r}(\theta,\eta)K_{0}'\left(\lambda A_{r}(\theta,\eta)\right)R(\eta)\partial_{\theta}\left(\tfrac{\sin(\eta-\theta)}{R(\theta)}\right)d\eta.$$
				Putting together the preceding identities yields to
				\begin{align}\label{papapa1}
				\mathcal{I}_{1}+\mathcal{I}_{3}=\partial_{\theta}\rho(\theta)\overline{V_{r}}(\lambda,\theta)+ \rho(\theta)V_{1}(\lambda,\theta)
				\end{align}
		with 
%				\begin{align}\label{Vr-tilde}\overline{V_{r}}(\lambda,\theta)=\bigintsss_{\mathbb{T}}K_{0}\left(\lambda A_{r}(\theta,\eta)\right)\partial_{\eta}\left(\frac{R(\eta)\sin(\eta-\theta)}{R(\theta)}\right)d\eta
%				\end{align}
%				and 
				\begin{align}\label{V1}
					V_{1}(\lambda,\theta)&:=\lambda\int_{\mathbb{T}}B_{r}(\theta,\eta)\partial_{\theta\eta}^{2}\big(R(\theta)R(\eta)\sin(\eta-\theta)\big)K_{0}'\left(\lambda A_{r}(\theta,\eta)\right)d\eta\nonumber\\
					&\quad-\lambda\int_{\mathbb{T}}\partial_{\eta}A_{r}(\theta,\eta)\partial_{\theta}\left(\tfrac{R(\eta)\sin(\eta-\theta)}{R(\theta)}\right)K_{0}'\left(\lambda A_{r}(\theta,\eta)\right)d\eta.
				\end{align}
				Differentiating term by term $\overline{V_{r}}$ with respect to $\theta$ gives
					\begin{align*}
						\partial_{\theta}\overline{V_{r}}(\lambda,\theta)&=\lambda\bigintssss_{\mathbb{T}}\partial_{\theta}A_{r}(\theta,\eta)K_{0}^\prime\left(\lambda A_{r}(\theta,\eta)\right)\partial_{\eta}\left(\tfrac{R(\eta)\sin(\eta-\theta)}{R(\theta)}\right)d\eta\\
						&\quad+\bigintssss_{\mathbb{T}}K_{0}\left(\lambda A_{r}(\theta,\eta)\right)\partial^2_{\theta\eta}\left(\tfrac{R(\eta)\sin(\eta-\theta)}{R(\theta)}\right)d\eta\\
						&:=\mathcal{J}_{1}+\mathcal{J}_{2}.
					\end{align*}
					Integrating  by parts in $\mathcal{J}_{2}$ yields
					$$\mathcal{J}_{2}=-\lambda\int_{\mathbb{T}}R(\eta)\partial_{\eta}A_{r}(\theta,\eta)K_{0}'\left(\lambda A_{r}(\theta,\eta)\right)\partial_{\theta}\left(\tfrac{\sin(\eta-\theta)}{R(\theta)}\right)d\eta.$$
					Combining the preceding identities allows to  we deduce that 
					\begin{align*}
						\partial_{\theta}\overline{V_{r}}(\lambda,\theta)&=\lambda\bigintssss_{\mathbb{T}}K_{0}^\prime\left(\lambda A_{r}(\theta,\eta)\right)\left[\partial_{\theta}A_{r}(\theta,\eta)\partial_{\eta}\left(\tfrac{R(\eta)\sin(\eta-\theta)}{R(\theta)}\right)-\partial_{\eta}A_{r}(\theta,\eta)\partial_{\theta}\left(\tfrac{R(\eta)\sin(\eta-\theta)}{R(\theta)}\right) \right]d\eta.
					\end{align*}
					Next we shall check the following identity
					\begin{equation}\label{magical simplification}
						\partial_{\theta}A_{r}(\theta,\eta)\partial_{\eta}\left(\tfrac{R(\eta)\sin(\eta-\theta)}{R(\theta)}\right)=B_{r}(\theta,\eta)\partial_{\theta\eta}^{2}\left(R(\theta)R(\eta)\sin(\eta-\theta)\right)-\partial_{\eta}A_{r}(\theta,\eta).
					\end{equation}
					Indeed, by \eqref{formula} and \eqref{definition of Br}, one finds
					$$\partial_{\theta}A_{r}(\theta,\eta)\partial_{\eta}\left(\tfrac{R(\eta)\sin(\eta-\theta)}{R(\theta)}\right)=\partial_{\theta}R(\theta)B_{r}(\theta,\eta)\partial_{\eta}(R(\eta)\sin(\eta-\theta))-\tfrac{R(\eta)\sin(\eta-\theta)\partial_{\eta}(R(\eta)\sin(\eta-\theta))}{A_{r}(\theta,\eta)}$$
					and 
					\begin{align*}B_{r}(\theta,\eta)\partial_{\theta\eta}^{2}\left(R(\theta)R(\eta)\sin(\eta-\theta)\right)&=\partial_{\theta}R(\theta)B_{r}(\theta,\eta)\partial_{\eta}(R(\eta)\sin(\eta-\theta))\\
						&\quad-\tfrac{(R(\theta)-R(\eta)\cos(\eta-\theta))\partial_{\eta}(R(\eta)\cos(\eta-\theta))}{A_{r}(\theta,\eta)}\cdot
					\end{align*}
					Putting together the foregoing identities leads to
					$$\frac{\partial_{\theta}A_{r}(\theta,\eta)\partial_{\eta}(R(\eta)\sin(\eta-\theta))}{R(\theta)}=B_{r}(\theta,\eta)\partial_{\theta\eta}^{2}\left(R(\theta)R(\eta)\sin(\eta-\theta)\right)+g(\theta,\eta),$$
					where
					\begin{align*}
						g(\theta,\eta)&:=\tfrac{1}{A_{r}(\theta,\eta)}\left[(R(\theta)-R(\eta)\cos(\eta-\theta))\partial_{\eta}(R(\eta)\cos(\eta-\theta))-R(\eta)\sin(\eta-\theta)\partial_{\eta}(R(\eta)\sin(\eta-\theta))\right]\\
						&=\tfrac{R(\theta)\partial_{\eta}(R(\eta)\cos(\eta-\theta))-R(\eta)\partial_{\eta}R(\eta)}{A_{r}(\theta,\eta)}\\
						&=-\partial_{\eta}A_{r}(\theta,\eta).
					\end{align*}
					This achieves the proof of \eqref{magical simplification}. From the periodicity we get
					$$\int_{\mathbb{T}}\lambda\partial_{\eta}A_{r}(\theta,\eta)K_{0}'\left(\lambda A_{r}(\theta,\eta)\right)d\eta=\int_{\mathbb{T}}\partial_{\eta}\left[K_{0}\left(\lambda A_{r}(\theta,\eta)\right)\right]d\eta=0$$
					and thus we get the following important identity
					$$\partial_{\theta}\overline{V_{r}}(\lambda,\theta)=V_{1}(\lambda,\theta).
					$$
			Plugging this into \eqref{papapa1} allows to get
				$$
				\mathcal{I}_{1}+\mathcal{I}_{3}=\partial_{\theta}\left(\overline{V_{r}}(\lambda,\theta)\rho(\theta)\right).
				$$
				Notice that it is easy to check that  if $r(-t,-\theta)=r(t,\theta)$, then 
				\begin{equation}\label{symmetry for V_{0}}
					\overline{V_{r}}(\lambda,-t,-\theta)=\overline{V_{r}}(\lambda,t,\theta).
				\end{equation}
				The next task is to  compute $\mathcal{I}_{2}+\mathcal{I}_{4}.$ Using  integration by parts in $\mathcal{I}_{4}$ gives,
				$$\mathcal{I}_{4}=-\lambda\displaystyle\int_{\mathbb{T}}\rho(\eta)\partial_{\eta}A_{r}(\theta,\eta)K_{0}'\left(\lambda A_{r}(\theta,\eta)\right)\partial_{\theta}\left(\tfrac{R(\theta)\sin(\eta-\theta)}{R(\eta)}\right)d\eta.$$
				From the symmetry property  $A_{r}(\theta,\eta)=A_{r}(\eta,\theta)$ and by  exchanging the roles of $\theta$ and $\eta$ in \eqref{magical simplification}, one deduces  
				$$\begin{array}{l}
					B_{r}(\eta,\theta)\partial_{\theta\eta}^{2}\left(R(\theta)R(\eta)\sin(\eta-\theta)\right)-\partial_{\eta}A_{r}(\theta,\eta)\partial_{\theta}\left(\tfrac{R(\theta)\sin(\eta-\theta)}{R(\eta)}\right)=-\partial_{\theta}A_{r}(\theta,\eta).
				\end{array}$$
				Therefore we obtain 
				$$\mathcal{I}_{2}+\mathcal{I}_{4}=-\partial_{\theta}\left(\int_{\mathbb{T}}\rho(\eta)K_{0}\left(\lambda A_{r}(\theta,\eta)\right)d\eta\right):=-\partial_{\theta}\mathbf{L}_{r}(\rho)(\lambda,\theta).$$
				Finally, by setting  
				$$V_{r}(\lambda,t,\theta)=\Omega+\overline{V_{r}}(\lambda,t,\theta)$$
				and combining the preceding identities, we end the proof of Lemma \ref{lemma general form of the linearized operator}.
			\end{proof}
			\subsubsection{The integrable case}
			The main purpose here is to explore the structure of the linearized operator at the equilibrium state. We shall see that the radial shape is reflected on the structure   the linearized operator  which  is a Fourier multiplier (of a convolution type). More precisely, we have the following result.
			\begin{lem}\label{lemma linearized operator at equilibrium}
				\begin{enumerate}
					\item 
					The linearized  equation of \eqref{Hamiltonian formulation of the equation} at the equilibrium state  $(r=0)$ writes,
					\begin{equation}\label{hamiltonian equation at equilibrium}
						\partial_{t}\rho=\partial_{\theta}\mathrm{L}(\lambda)\rho=\partial_{\theta}\nabla H_{\mathrm{L}}(\rho),
					\end{equation}
					where $\mathrm{L}(\lambda)$ is the self-adjoint operator defined by $\mathrm{L}(\lambda)=-V_{0}(\lambda)+\mathcal{K}_{\lambda}\ast_{\theta}$ with
					\begin{equation}\label{definition of V0}
						V_{0}(\lambda)=\Omega+I_{1}(\lambda)K_{1}(\lambda)
					\end{equation}
					and
					\begin{equation}\label{definition of mathcalKlambda}
						\mathcal{K}_{\lambda}(\theta)=K_{0}\left(2\lambda\left|\sin\left(\tfrac{\theta}{2}\right)\right|\right).
					\end{equation}
					We refer to the Appendix A for the definitions of the modified Bessel functions $I_{1},$ $K_{1}$ and $K_{0}.$\\
					Moreover, the  Hamiltonian $H_{\mathrm{L}}$ is quadratic and takes the form
					$$H_{\mathrm{L}}(\rho)=\tfrac{1}{2}\langle\mathrm{L}(\lambda)\rho,\rho\rangle_{L^{2}(\mathbb{T})}.$$
					\item  The solutions to  \eqref{hamiltonian equation at equilibrium} with zero space  average are given by 
					%\forall j\in\mathbb{Z},\dot{\rho_{j}}=-i\Omega_{j}(\lambda)\rho_{j}\mbox{ (Fourier multiplier)}\mbox{ that is }
					\begin{equation}\label{sol eq}
						\rho(t,\theta)=\sum_{j\in\mathbb{Z}^{*}}\rho_{j}(0)e^{\ii(j\theta-\Omega_{j}(\lambda)t)},
					\end{equation}
					with
					\begin{equation}\label{def eigenvalues at equilibrium}
						\Omega_{j}(\lambda)=j\big[\Omega+(I_{1}K_{1})(\lambda)-(I_{j}K_{j})(\lambda)\big].
					\end{equation}
					and for $\displaystyle{\rho(\theta)=\sum_{j\in\mathbb{Z}^{*}}\rho_{j}e^{\ii j\theta}}$ we have 
					\begin{equation}\label{defLHL}
						\mathrm{L}(\lambda)\rho(\theta)=-\sum_{j\in\mathbb{Z}^{*}}\tfrac{\Omega_{j}(\lambda)}{j}\rho_{j}e^{\ii j\theta}\quad \mbox{ and }\quad H_{\mathrm{L}}\rho=-\sum_{j\in\mathbb{Z}^{*}}\tfrac{\Omega_{j}(\lambda)}{2j}|\rho_{j}|^{2},
					\end{equation}

				\end{enumerate}
			\end{lem}
			Before proceeding with the proof we want to give some remarks.
			\begin{remark}
				\begin{enumerate}[label=\textbullet]
					\item When $\Omega=0$ the eigenvalue $\Omega_{1}(\lambda)$   vanishes for any $\lambda$ due to the rotation invariance of the equation   and the use of the free parameter $\Omega$ is   to avoid this degeneracy. However the trivial resonance $\Omega_0(\lambda)=0$ can be removed   by imposing the zero space average  which is preserved by
	the nonlinear dynamics from the Hamiltonian structure as we have seen before in \eqref{pres avrg}.
	\item  The solutions to the linear equation at the equilibrium are aperiodic and if we excite only a finite number of frequencies with non-resonances assumption we get quasi-periodic solutions. We will make a precise comment later on Proposition $\ref{lemma sol Eq}.$
%					\item The frequencies expression agree   with the formula obtained in \cite[Prop. 5.8.]{DHR19} by taking $\Omega=0.$ In addition when $\lambda=\Omega=0$ we obtain in view of \eqref{asymptotic expansion of small argument} Euler frequencies $\Omega_j(0)=\frac{j-1}{2{j}}\cdot$
				\end{enumerate}
			\end{remark}
			\begin{proof}
				${\bf 1.}$ First observe that from  \eqref{formula}, one deduces for $r=0$ that  $A_{0}(\theta,\eta)=2\left|\sin\left(\frac{\eta-\theta}{2}\right)\right|$. Then we obtain from \eqref{definition of Vr} and  \eqref{definition of mathbfLr},
				\begin{align*}
					\mathbf{L}_{0}\rho(\lambda,\theta)&=\bigintssss_{\mathbb{T}}\rho(\eta)K_{0}\left(2\lambda\left|\sin\left(\tfrac{\eta-\theta}{2}\right)\right|\right)d\eta\\
					&=\mathcal{K}_{\lambda}\ast\rho(\theta),
				\end{align*}
				with $\mathcal{K}_{\lambda}$ defined in \eqref{definition of mathcalKlambda} and using the change of variables $\eta\mapsto\eta+\theta$ we obtain
				\begin{align*}
					V_{0}(\lambda,\theta) & =  \displaystyle\Omega+\bigintssss_{\mathbb{T}}K_{0}\left(2\lambda\left|\sin\left(\tfrac{\eta-\theta}{2}\right)\right|\right)\cos(\eta-\theta)d\eta\\
					& = \displaystyle\Omega+\bigintssss_{\mathbb{T}}K_{0}\left(2\lambda\left|\sin\left(\tfrac{\eta}{2}\right)\right|\right)\cos(\eta)d\eta\\
					&:=V_{0}(\lambda).
				\end{align*}
				We remark that if we write $e_{j}(\theta)=e^{\ii j\theta},$ then direct computations yield \begin{align*}
					(\mathcal{K}_{\lambda}\ast e_{j})(\theta) & =  \displaystyle\bigintssss_{\mathbb{T}}K_{0}\left(2\lambda\left|\sin\left(\tfrac{\eta}{2}\right)\right|\right)e^{\ii j(\theta-\eta)}d\eta\\
					& = \displaystyle e_{j}(\theta)\bigintssss_{\mathbb{T}}K_{0}\left(2\lambda\left|\sin\left(\tfrac{\eta}{2}\right)\right|\right)e^{-\ii j\eta}d\eta.
				\end{align*}
				Since the function $\eta\mapsto K_{0}\left(2\lambda\left|\sin\left(\tfrac{\eta}{2}\right)\right|\right)$ is even, we deduce using the change of variables $\eta=2\tau+\pi$ and the formula \eqref{Bessel} that
				\begin{align*}
					\displaystyle\bigintssss_{\mathbb{T}}K_{0}\left(2\lambda\left|\sin\left(\tfrac{\eta}{2}\right)\right|\right)e^{-\ii j\eta}d\eta & =  \displaystyle\bigintssss_{\mathbb{T}}K_{0}\left(2\lambda\left|\sin\left(\tfrac{\eta}{2}\right)\right|\right)\cos\left(j\eta\right)d\eta\\
					& = \displaystyle\tfrac{(-1)^{j}}{\pi}\bigintssss_{-\frac{\pi}{2}}^{\frac{\pi}{2}}K_{0}(2\lambda\cos(\tau))\cos(2j\tau)d\tau\\
					%& =  \displaystyle\tfrac{2(-1)^{j}}{\pi}\bigintssss_{0}^{\frac{\pi}{2}}K_{0}(2\lambda\cos(\tau))\cos(2j\tau)d\tau\\
					& = (I_{j}K_{j})(\lambda).
				\end{align*}
				Hence, the Fourier coefficients of $\mathcal{K}_{\lambda}$ are 
				\begin{equation}\label{Fourier coefficients of mathcalKlambda}
					\left(\mathcal{K}_{\lambda}\right)_{j}=(I_{j}K_{j})(\lambda).
				\end{equation}
				Similar arguments as  before with $j=1$ allow to get
				$$V_{0}(\lambda)=\Omega+(I_{1}K_{1})(\lambda).$$
				Recall that  $\mathcal{K}_{\lambda}$ is even and then we find that  $\mathrm{L}(\lambda)$ is self-adjoint in $L^2(\mathbb{T})$.\\
				${\bf 2.}$
				Starting from the Fourier  expansion  $\rho(t,\theta)=\displaystyle\sum_{j\in\mathbb{Z}^{*}}\rho_{j}(t)e^{\ii j\theta}$, then we can easily ensure from  direct computations using the previous results,  that $\rho$ solves the equation \eqref{hamiltonian equation at equilibrium} if and only if
				$$\dot{\rho_{j}}=-\ii\Omega_{j}(\lambda)\rho_{j}\quad\mbox{ with }\quad\Omega_{j}(\lambda)=j\big[\Omega+(I_{1}K_{1})(\lambda)-(I_{j}K_{j})(\lambda)\big],$$
				and therefore 
				$$\rho(t,\theta)=\sum_{j\in\mathbb{Z}^{*}}\rho_{j}(0)e^{\ii(j\theta-\Omega_{j}(\lambda)t)}.
				$$
				Concerning the identities \eqref{defLHL} they can be obtained from straightforward computations. This ends the proof of Lemma \ref{lemma linearized operator at equilibrium}.
			\end{proof}
			\subsection{Structure of the linear frequencies}
			The main target  in this section is to explore some
			interesting  structures of the equilibrium frequencies. We shall in particular focus on  their
			monotonicity and detail some asymptotic behavior for large modes. Another important discussion will be devoted  to  the non-degeneracy of these frequencies through the so-called R\"usseman conditions.
			This is the cornerstone step  in measuring   the final Cantor set giving rise to quasi-periodic solutions for the
			linear/nonlinear problems. Actually, in the nonlinear case the final Cantor appears as a perturbation of the Cantor set constructed from
			the equilibrium eigenvalues and therefore perturbative arguments based on their non-degeneracy are
			very useful and will be performed in Section
			\ref{Section measure of the final Cantor set}.
			
			\subsubsection{Monotonicity and asymptotic behaviour}
			Our purpose is   to establish some useful
			properties related to the monotonicity and the asymptotic behavior for large modes of the eigenvalues
			of the linearized operator at the equilibrium state. Notice that their explicit values are detailed
			in \eqref{def eigenvalues at equilibrium}. Our result reads as follows.
			\begin{lem}\label{lemma properties linear frequencies}
				Let $\Omega>0$ and $ \lambda\in\mathbb{R}$, then the  frequencies $(\Omega_{j}(\lambda))_{j\in\mathbb{Z}^*}$ satisfy the following properties.
				\begin{enumerate}[label=(\roman*)]
					\item For any $ j\in\mathbb{Z}^*,\,\lambda>0$ we have $\Omega_{-j}(\lambda)=-\Omega_{j}(\lambda).$
					\item For any $\lambda>0,$ the sequences $(\Omega_{j}(\lambda))_{j\in\mathbb{N}^*}$ and $\left(\tfrac{\Omega_{j}(\lambda)}{j}\right)_{j\in\mathbb{N}^{*}}$ are strictly increasing.
					\item For any $\lambda>0,$ the following expansion holds
					\begin{equation}\label{asymptotic expansion of the eigenvalues at the equilibrium}
						\Omega_{j}(\lambda)\underset{j\rightarrow\infty}{=}V_{0}(\lambda)j-\tfrac{1}{2}+\tfrac{\lambda^{2}}{4j^{2}}+O_{\lambda}\left(\tfrac{1}{j^{4}}\right),
					\end{equation} where $V_{0}(\lambda)$ is defined in \eqref{definition of V0}.
					\item For any $j\in\mathbb{Z}^*,\,\lambda>0$ we have 
					$$|\Omega_{j}(\lambda)|\geqslant \Omega|j|.$$
					\item Given $0<\lambda_{0}<\lambda_{1}$, there exists $C_0>0$ such that	
									$$
					\forall \lambda\in[\lambda_{0},\lambda_{1}],\forall j,j_0\in\mathbb{Z},\quad   |\Omega_{j}(\lambda)\pm\Omega_{j_{0}}(\lambda)|\geqslant C_{0}|j\pm j_{0}|.
					$$
				
					\item Given $0<\lambda_{0}<\lambda_{1}$ and  $q_{0}\in\mathbb{N},$ there exists $C_0>0$ such that   
					$$
					\forall\, j,j_0\in\mathbb{Z}^*,\quad \max_{q\in\llbracket 0,q_{0}\rrbracket}\sup_{\lambda\in[\lambda_{0},\lambda_{1}]}\left|\partial_{\lambda}^{q}\left(\Omega_{j}(\lambda)-\Omega_{j_{0}}(\lambda)\right)\right|\leqslant C_{0}|j-j_{0}|.
					$$
					
				\end{enumerate}
			\end{lem}
			\begin{proof}
				\textbf{(i)} It is an immediate consequence of \eqref{def eigenvalues at equilibrium} and \eqref{symmetry Bessel}.\\
				\textbf{(ii)} The monotonicity of the sequence   $\left(\tfrac{\Omega_{j}(\lambda)}{j}\right)_{j\in\mathbb{N}^{*}}$ is proved in \cite[Prop. 5.9. (1)]{DHR19}, see also the Appendix \ref{appendix Bessel}. It follows that 
				the sequence $(\Omega_{j}(\lambda))_{j\in\mathbb{N}^{*}}$ is strictly increasing as the product of two strictly increasing  positive sequences.\\
				\textbf{(iii)} It is an immediate consequence of \eqref{def eigenvalues at equilibrium} and the asymptotic expansion \eqref{asymptotic expansion of high order}\\
				\textbf{(iv)}
				Recall that  $j\mapsto \Omega_j(\lambda)$ is odd and vanishes at $j=0$. Then it suffices to check  the result for $j\in\mathbb{N}^*.$ According to the Appendix \ref{appendix Bessel}, the sequence $j\mapsto (I_{j}K_{j})(\lambda)$ is decreasing and therefore
				\begin{align}\label{zita1}
					\forall \lambda> 0,\quad  (I_{1}K_{1})(\lambda)- (I_{j}K_{j})(\lambda)\geqslant 0.
				\end{align}
				It follows that 
				$$\forall\, \lambda>0,\quad |\Omega_{j}(\lambda)|\geqslant \Omega j.
				$$  
				\textbf{(v)} By the oddness of $j\mapsto \Omega_j(\lambda)$ it is enough  to establish the estimate for $j,j_0\in\mathbb{N}.$ We shall first focus on the estimate of the difference $\Omega_{j}(\lambda)-\Omega_{j_{0}}(\lambda)$. Without  loss of generality we can assume that   $j> j_{0}\geqslant 1,$ (The case $j=j_0$ is obvious and the case $j_{0}=0$ brings us back to the previous point). One may write by \eqref{def eigenvalues at equilibrium} that  for  $\lambda>0$,
				\begin{align}\label{diff eig eq}
					\nonumber\Omega_{j}(\lambda)-\Omega_{j_{0}}(\lambda)&=(j-j_{0})\Big(\Omega+I_{1}(\lambda)K_{1}(\lambda)-I_{j}(\lambda)K_{j}(\lambda)\Big)\\
					&\quad+j_{0}\Big(I_{j_{0}}(\lambda)K_{j_{0}}(\lambda)-I_{j}(\lambda)K_{j}(\lambda)\Big).
				\end{align}
				Combining this estimate with \eqref{zita1} yields
				\begin{align}\label{Mazh1}
					\Omega_{j}(\lambda)-\Omega_{j_{0}}(\lambda)\geqslant(j-j_{0})\Omega+j_{0}\Big(I_{j_{0}}(\lambda)K_{j_{0}}(\lambda)-I_{j}(\lambda)K_{j}(\lambda)\Big).
				\end{align}
				We need to get refined estimate for the last term of the right hand side. For this goal we use the formulae  \eqref{form-Lebe0} to write
				\begin{align}\label{form-Lebe1}
					(I_{n}K_{n})(\lambda)=\tfrac12\int_{0}^{\infty} J_0\big(2\lambda\sinh(\tfrac{t}{2})\big)e^{-nt}dt.
				\end{align}
				This allows to construct for a fixed $\lambda$ a smooth extension $n\in(0,\infty)\mapsto (I_{n}K_{n})(\lambda)$. Thus differentiating term by term using change of variable we get   for any $m\in\mathbb{N}$
				\begin{align}
					\label{decayk}
					\nonumber \sup_{\lambda\in\mathbb{R}}\left|\partial_n^m(I_{n}K_{n})(\lambda)\right|&\leqslant \tfrac12\int_{0}^{\infty} t^me^{-nt}dt\\
					&\leqslant \tfrac{m!}{2 n^{m+1}},
				\end{align}
				where we have used the classical estimates for  Bessel functions (applied with $n=q=0$)
				\begin{align}\label{Bessel-est1}
					\sup_{n,q\in\mathbb{N}\atop x\in\mathbb{R}}|J_n^{(q)}(x)|\leqslant 1,
				\end{align}
				which follows easily from the integral representation \eqref{Besse-repr}. In particular, for $m=1$ we find that for  for any $n\geqslant 1$
				$$
				\sup_{\lambda\in\mathbb{R}}\Big|\partial_n(I_{n}K_{n})(\lambda)\Big|\leqslant  \tfrac{1}{2n^2}\cdot
				$$
				Therefore applying   Taylor Formula we infer for $j>j_0\geqslant 1$
				\begin{align}\label{HUP1}
					\nonumber \sup_{\lambda\in\mathbb{R}}\Big|(I_{j}K_{j})(\lambda)-(I_{j_0}K_{j_0})(\lambda)\Big|&\leqslant \tfrac12\int_{j_0}^j\frac{dn}{n^2}\\
					&\leqslant \tfrac{|j-j_0|}{2\,j\,j_0}\cdot
				\end{align}
				Inserting this estimate into \eqref{Mazh1} gives
				$$\Omega_{j}(\lambda)-\Omega_{j_{0}}(\lambda)\geqslant (j-j_{0})\Big(\Omega-\tfrac{1}{2j}\Big).
				$$
				Therefore for $j> N=\big[\Omega^{-1}\big]$ and $j> j_0\geqslant 1$ we get 
				\begin{align}\label{Tyy1}
					\Omega_{j}(\lambda)-\Omega_{j_{0}}(\lambda)\geqslant\tfrac12\Omega\, (j-j_{0}).
				\end{align}
				Now for  $j\neq j_{0}\in\llbracket 1,N\rrbracket$ we get from the  point (ii) that the map $\lambda\in[\lambda_0,\lambda_1]\mapsto\Omega_{j}(\lambda)-\Omega_{j_{0}}(\lambda)$ does not vanish and therefore we can find by a compactness argument a constant $C>0$ such that 
				$$\forall\lambda\in[\lambda_{0},\lambda_{
					1}],\quad|\Omega_{j}(\lambda)-\Omega_{j_{0}}(\lambda)|\geqslant C |j-j_{0}|.$$
				Taking $C_0=\min(C,\frac12\Omega)$ and combining  the preceding inequality with \eqref{Tyy1} we obtain
				$$ 
				\forall\lambda\in[\lambda_{0},\lambda_{1}],\forall j\geqslant j_{0}\geqslant 1,\quad |\Omega_{j}(\lambda)-\Omega_{j_{0}}(\lambda)|\geqslant C_{0}|j-j_{0}|.
				$$
				Finally we get 
				$$
				\forall\lambda\in[\lambda_{0},\lambda_{1}],\forall j,j_{0}\in\mathbb{N},\quad |\Omega_{j}(\lambda)-\Omega_{j_{0}}(\lambda)|\geqslant C_{0}|j-j_{0}|.
				$$
				Let us now move to the estimate  $\Omega_{j}(\lambda)+\Omega_{j_{0}}(\lambda)$ for $j,j_0\in\mathbb{N}.$ Since both quantities are positive then using the point (iv) yields
				$$
				\forall\lambda\in[\lambda_{0},\lambda_{1}],\quad|\Omega_{j}(\lambda)+\Omega_{j_{0}}(\lambda)|=\Omega_{j}(\lambda)+\Omega_{j_{0}}(\lambda)\geqslant \Omega(j+j_{0})\geqslant C_0(j+j_{0}).
				$$
				This completes the proof of the desired estimate.\\
				\textbf{(vi)} Let $q_{0}\in\mathbb{N}^*.$ let $q\in\llbracket 0,q_{0}\rrbracket.$ Differentiating $q$ times \eqref{diff eig eq} in $\lambda$, one obtains
				\begin{align}\label{deriv diff eig}
					\partial_{\lambda}^{q}\big(\Omega_{j}(\lambda)-\Omega_{j_{0}}(\lambda)\big)&=(j-j_{0})\Big(\partial_{\lambda}^{q}\Omega+\partial_{\lambda}^{q}\big(I_{1}(\lambda)K_{1}(\lambda)\big)-\partial_{\lambda}^{\alpha}\big(I_{j}(\lambda)K_{j}(\lambda)\big)\Big)\nonumber\\
					&\quad+j_{0}\partial_{\lambda}^{q}\Big(I_{j_{0}}(\lambda)K_{j_{0}}(\lambda)-I_{j}(\lambda)K_{j}(\lambda)\Big).
				\end{align}
				Similarly, we get by differentiating $q$ times in $\lambda$ the identity  \eqref{form-Lebe1} 
				\begin{align}\label{Diff-YU1}
					\partial_\lambda^q(I_{n}K_{n})(\lambda)=2^{q-1}\int_{0}^\infty J_0^{(q)}\big(2\lambda\sinh(\tfrac{t}{2})\big) \sinh^q(\tfrac{t}{2})e^{-nt}dt.
				\end{align}
				From \eqref{Bessel-est1} we deduce for any $\lambda\in[\lambda_0,\lambda_1],$
				\begin{align*}
					\left|\partial_\lambda^q(I_{n}K_{n})(\lambda)\right|&\leqslant  2^{q-1}\int_{0}^\infty  \sinh^q(\tfrac{t}{2})e^{-nt}dt.
				\end{align*}
				Then using the inequality $\sinh x\leqslant \frac{e^x}{2}$ for $x\geqslant 0$  we get for $n>\frac{q}{2}$
				\begin{align}\label{decay deriv InKn}
					\left|\partial_\lambda^q(I_{n}K_{n})(\lambda)\right|&\leqslant  \tfrac{1}{2}\int_{0}^\infty  e^{\left(\tfrac{q}{2}-n\right)t}dt\nonumber\\
					&\leqslant  \tfrac{1}{2n-q}\cdot
				\end{align}
				By compactness argument, we deduce that
				\begin{equation}\label{est deriv prod IjKj}
					\sup_{j\in\mathbb{N}}\max_{q\in\llbracket 0,q_{0}\rrbracket}\sup_{\lambda\in[\lambda_{0},\lambda_{1}]}\left|\partial_{\lambda}^{q}\big(I_{j}(\lambda)K_{j}(\lambda)\big)\right|\leqslant C.
				\end{equation}
				Differentiating in $n$ \eqref{Diff-YU1} yields
				\begin{align*}
					\partial_\lambda^q\partial_n(I_{n}K_{n})(\lambda)=-2^{q-1}\int_{0}^\infty J_0^{(q)}\big(2\lambda\sinh(t/2)\big) \sinh^q(t/2)te^{-nt}dt.
				\end{align*}
				Therefore  applying similar arguments used to show \eqref{decay deriv InKn} gives for $2n>q$
				\begin{align}\label{Es-W-L}
					\nonumber \left|\partial_\lambda^q\partial_n(I_{n}K_{n})(\lambda)\right|&\leqslant \tfrac{1}{2}\int_{0}^\infty te^{-\left(n-\tfrac{q}{2}\right)t}dt\\
					&\leqslant \frac{2}{(2n-q)^2}\cdot
				\end{align}
				Then Taylor Formula allows to get for $j, j_{0}>\frac{q}{2}$
				\begin{equation}\label{est deriv diff IjKj-1}
					\sup_{\lambda\in[\lambda_0,\lambda_1]}\left|\partial_\lambda^q(I_{j}K_{j}-I_{j_{0}}K_{j_{0}})(\lambda)\right|\leqslant C\tfrac{|j-j_{0}|}{j j_{0}}\cdot
				\end{equation}
				Setting $N=\lfloor\frac{q_{0}}{2}\rfloor+1$, one obtains for any $j,j_{0}\geqslant N$
				\begin{align*}
					\nonumber \max_{q\in\llbracket 0,q_{0}\rrbracket}\sup_{\lambda\in[\lambda_0,\lambda_1]}\left|\partial_\lambda^q(I_{j}K_{j}-I_{j_{0}}K_{j_{0}})(\lambda)\right|&\leqslant C|j-j_{0}|\cdot
				\end{align*}
				By compactness argument, one obtains for any $j,j_{0}\in\llbracket 1,N\rrbracket$
				\begin{align*}
					\nonumber \max_{q\in\llbracket 0,q_{0}\rrbracket}\sup_{\lambda\in[\lambda_0,\lambda_1]}\left|\partial_\lambda^q(I_{j}K_{j}-I_{j_{0}}K_{j_{0}})(\lambda)\right|&\leqslant C|j-j_{0}|\cdot
				\end{align*}
				Now for the remaining case  $j_0\in\llbracket 1,N\rrbracket$ and $j\geqslant N$ one has gathering the previous two estimates 
				\begin{align*}
					\nonumber \max_{q\in\llbracket 0,q_{0}\rrbracket}\sup_{\lambda\in[\lambda_0,\lambda_1]}\left|\partial_\lambda^q(I_{j}K_{j}-I_{j_{0}}K_{j_{0}})(\lambda)\right|&\leqslant\max_{q\in\llbracket 0,q_{0}\rrbracket}\sup_{\lambda\in[\lambda_0,\lambda_1]}\left|\partial_\lambda^q(I_{j}K_{j}-I_{N}K_{N})(\lambda)\right|\\
					&\quad+\max_{q\in\llbracket 0,q_{0}\rrbracket}\sup_{\lambda\in[\lambda_0,\lambda_1]}\left|\partial_\lambda^q(I_{N}K_{N}-I_{j_{0}}K_{j_{0}})(\lambda)\right|\\
					&\leqslant C|j-N|+C|N-j_{0}|
					\leqslant C|j-j_{0}|\cdot
				\end{align*}
				Thus we can find $C>0$ such that for any  $j,j_{0}\in\mathbb{N}^*$
				\begin{equation}\label{est deriv diff IjKj}
					\nonumber \max_{q\in\llbracket 0,q_{0}\rrbracket}\sup_{\lambda\in[\lambda_0,\lambda_1]}\left|\partial_\lambda^q(I_{j}K_{j}-I_{j_{0}}K_{j_{0}})(\lambda)\right|\leqslant C|j-j_{0}|\cdot
				\end{equation}
				Putting together \eqref{deriv diff eig}, \eqref{est deriv prod IjKj} and \eqref{def eigenvalues at equilibrium} yields				$$\max_{q\in\llbracket 0,q_{0}\rrbracket}\sup_{\lambda\in[\lambda_{0},\lambda_{1}]}\left|\partial_{\lambda}^{q}\left(\Omega_{j}(\lambda)-\Omega_{j_{0}}(\lambda)\right)\right|\leqslant C|j-j_{0}|.$$
				This ends the proof of Lemma \ref{lemma properties linear frequencies}.
			\end{proof}
			\subsubsection{Non-degeneracy and transversality}
			Fix finitely many tangential sites
			$$\mathbb{S}:=\{j_{1},\dots,j_{d}\}\subset\mathbb{N}^*\quad\textnormal{with}\quad d\geqslant 1\quad\textnormal{and}\quad 1\leqslant  j_{1}<\dots<j_{d}.$$
			We consider the linear vector frequency  at the equilibrium state
			\begin{equation}\label{def linear frequency vector}
				\omega_{\textnormal{Eq}}(\lambda):=(\Omega_{j}(\lambda))_{j\in\mathbb{S}},
			\end{equation}
			where $\Omega_{j}(\lambda)$ is defined by \eqref{def eigenvalues at equilibrium}. The main purpose  is to  study  some  Diophantine structure of the  curve $\lambda\in(\lambda_0,\lambda_1)\mapsto{\omega}_{\textnormal{\tiny{Eq}}}(\lambda)$ for fixed $0<\lambda_{0}<\lambda_{1}.$ In particular, we shall focus on  the non-degeneracy and the transversality conditions of  these eigenvalues which are essential in getting non trivial Cantor set from which  quasi-periodic solutions emerge  at the linear and nonlinear levels.  Notice that the approach that we shall implement  here has been developed  before in several papers such as \cite{BBMH18, BBM11, R01}. Before exploring these properties we need to fix  some definitions. 
			\begin{defin}\label{def-degenerate} 
				Given two numbers $\lambda_0<\lambda_1$ and $d\in\mathbb{N}^*$, a vector-valued function $f = (f_1, ..., f_d ) : [\lambda_0,\lambda_1] \to \mathbb{R}^d$ is called non-degenerate if, for any vector $c = (c_1,...,c_d) \in  \mathbb{R}^d \setminus \{0\}$, the function $f \cdot c = f_1c_1 + ... + f_dc_d$ is not identically zero on the whole interval $[\lambda_0,\lambda_1]$. This means that the curve of $f$ is not contained in an hyperplane.
			\end{defin}
			Now we shall prove the  following  result on the non-degeneracy of the linear frequencies which   is related to the asymptotic behavior of Bessel functions  $(I_j K_j)(\lambda)$ for large values of $\lambda$.  This property will be crucial  to check a suitable transversality assumption. 
			\begin{lem}\label{Non-degenracy1}%[non-degeneracy of unperturbed linear frequencies]
				Let $\Omega\in\mathbb{R}^*$ and $0<\lambda_0<\lambda_1$, then the  frequency curve ${\omega}_{\textnormal{\tiny{Eq}}}$  defined by \eqref{def linear frequency vector} and the vector-valued function  $\lambda\mapsto(\Omega+I_{1}K_{1},{\omega}_{\textnormal{\tiny{Eq}}})\in\mathbb{R}^{d+1}$ are non degenerate on $[\lambda_0,\lambda_1]$ in the sense of the Definition $\ref{def-degenerate}.$
			\end{lem} 
			
			\begin{proof}
				$\blacktriangleright$ Let us start with checking the non-degeneracy of ${\omega}_{\textnormal{\tiny{Eq}}}$. For this aim, we shall argue by contradiction and assume    the existence of a fixed vector $c=(c_{k})_{0\leqslant k\leqslant d }\in\mathbb{R}^{d}$ such that 
				\begin{align}\label{Syst1-1}
					\forall\lambda\in[\lambda_0,\lambda_1],\quad \sum_{k=1}^{d}c_{k}\Omega_{j_{k}}(\lambda)=0.
				\end{align}
				Since for all $j\in\mathbb{N}^{*}$, the application $\lambda\mapsto (I_{j}K_{j})(\lambda)$ admits  a holomorphic extension in the open connected set $\big\{ \lambda\in\mathbb{C}, \textnormal{Re}(\lambda)>0\big\}$ (see Appendix \ref{appendix Bessel}) then by the continuation principle we obtain
				\begin{align}\label{iden-YG1}
					\forall\lambda>0,\quad  \sum_{k=1}^{d }c_{k}j_{k}(I_{j_{k}}K_{j_{k}})(\lambda)=\left(\sum_{k=1}^{d }c_{k}j_{k}\right)\big((I_{1}K_{1})(\lambda)+\Omega\big).
				\end{align}
				Using the asymptotic expansion \eqref{asymptotic expansion of large argument} obtained  for   $I_{j}K_{j}$ with large $\lambda,$  we first get  
				$$\forall j\in\mathbb{N}^{*},\quad \lim_{\lambda\to \infty}(I_{j}K_{j})(\lambda)=0.$$
				Then taking the limit in \eqref{iden-YG1} as $\lambda \to \infty$ implies 
				$$
				\Omega \sum_{k=1}^{d }c_{k}j_{k}=0.
				$$
				Since we assumed that  $\Omega\neq 0,$ then necessary  we find that $\displaystyle\sum_{k=1}^{d }c_{k}j_{k}=0$ which implies in turn according to \eqref{iden-YG1}  
				$$
				\forall\lambda>0,\quad  \sum_{k=1}^{d }c_{k}j_{k}(I_{j_{k}}K_{j_{k}})(\lambda)=0.
				$$
				Applying once again the expansion \eqref{asymptotic expansion of large argument}  yields 
				\begin{equation}\label{iden-YG2}
					\forall m\in\llbracket 1,d \rrbracket,\quad  \sum_{k=1}^{d }c_{k}j_{k}\alpha_{j_{k},m}=0.
				\end{equation}
				We consider the matrix $A_d=(A_{m,k})_{1\leqslant m,k\leqslant d }\in M_{d }(\mathbb{R})$ defined by 
				$$\forall(m,k)\in\llbracket 1,d \rrbracket^{2},\quad A_{m,k}=j_{k}\alpha_{j_{k},m}.
				$$
				Then the system \eqref{iden-YG2} is equivalent to $A_dc=0$ with $c=\left(\begin{array}{ccc}
					c_1 \\
					\vdots \\
					c_d
				\end{array}\right).$ To get the desired result, $c=0$, it suffices to check that $\textnormal{det } A_d\neq 0$. Using the expression of the coefficients $\alpha_{j_{k},m}$ in \eqref{Form-Pol}  one deduces that
				\begin{align}\label{ZRT1}
					\alpha_{j_{k},m}= a_m(\mu_{j_k}-1)Q_{m}\big(\mu_{j_k}\big),\quad a_m=(-1)^{m}\tfrac{(2m)!}{4^m\big(m!\big)^2},\quad \mu_j=4{j^2},
				\end{align}
				with $Q_1(X)=1 $ and for $m\geq 2$
				\begin{align*}
					Q_m(X)=\prod_{\ell=2}^{m}\big(X-(2\ell-1)^2\big).
				\end{align*}
				Remark that $Q_m$ is a unitary polynomial of degree $m-1$. Using the homogeneity of the determinant with respect to each column and row we find
				$$
				\textnormal{det } A_d=\prod_{m,k=1}^{d}a_m(\mu_{j_k}-1) \,\textnormal{det }B_d,
				$$ 
				with $B_d$ the matrix given by 
				$$B_d=\left(\begin{array}{ccc}
					Q_{1}(\mu_{j_{1}}) & \cdots & Q_{1 }(\mu_{j_{d}})\\
					\vdots &  & \vdots\\
					Q_{d}(\mu_{j_{1}}) & \cdots & Q_{d }(\mu_{j_{d}})
				\end{array}\right).$$
				Therefore we infer  that $A_d$ is nonsingular if  $\textnormal{det}B_d\neq0$. On the other hand, 
				the computation of  $\textnormal{det}B_d$ can be done   in a similar way to  Vandermonde determinant. Indeed, define the polynomial given by the determinant 
				$$P(X)=\left|\begin{array}{cccc}
					Q_{1}(\mu_{j_{1}}) & \cdots&Q_{1 }(\mu_{j_{d-1}}) & Q_{1 }(X)\\
					\vdots &  & \vdots&\vdots\\
					Q_{d}(\mu_{j_{1}}) & \cdots &Q_{d }(\mu_{j_{d-1}})& Q_{d }(X)
				\end{array}\right|.$$
				Then $P$ is a polynomial of degree $d-1$ and vanishes at all the points $X=\mu_{j_k}$ for $k\in\llbracket 1,d-1\rrbracket.$ Consequently, we get
				$$
				\textnormal{det}B_d=P(\mu_{j_d})=\textnormal{det}B_{d-1}\prod_{k=1}^{d-1}\big( \mu_{j_d}-\mu_{j_k}\big).
				$$
				Therefore, iterating this identity yields
				$$
				\textnormal{det}B_d=\prod_{1\leqslant k<\ell\leqslant d-1}\big( \mu_{j_\ell}-\mu_{j_k}\big).
				$$
				Since  $\mu_{j_\ell}\neq\mu_{j_k}$ for $\ell\neq k$ we get $\textnormal{det}B_d\neq0$ which achieves the proof of the first point.\\
				$\blacktriangleright$ Next we move to the second point of the lemma and show that if 
				$$
				\forall\lambda\in[\lambda_0,\lambda_1],\quad c_{0}\Big(\Omega+(I_{1}K_{1})(\lambda)\Big)+\sum_{k=1}^{d}c_{k}j_k\Big(\Omega+(I_{1}K_{1})(\lambda)-(I_{j_k}K_{j_k})(\lambda)\Big)=0,
				$$
				then necessary $c_0=...=c_d=0.$ As before we can extend by analyticity the preceding identity to $(0,\infty)$
				By checking the terms in $\frac1\lambda$  in the preceding identity using  \eqref{asymptotic expansion of large argument} we find immediately that  $c_0=0$. Therefore the system reduces to \eqref{Syst1-1} and then we may apply the result of the first point in order to  get $c_1=...=c_d=0.$ This completes the proof of Lemma  \ref{Non-degenracy1}.
			\end{proof}
			
			The next goal is to check that  R\"{u}ssemann transversality conditions  are satisfied for the linear frequencies of the equilibrium state.  Namely, we shall prove the following result in  the spirit of the papers  \cite{BBMH18,BBM11,R01}.
			\begin{lem}{\textnormal{[Transversality]}}\label{lemma transversality}
				{Given $0<\lambda_0<\lambda_1$, there exist $q_{0}\in\mathbb{N}$ and $\rho_{0}>0$ such that the following results hold true. Recall that $\omega_{\textnormal{Eq}}$ and $\Omega_j$ are defined in \eqref{def linear frequency vector} and \eqref{def eigenvalues at equilibrium} respectively.
					\begin{enumerate}[label=(\roman*)]
						\item For any $l\in\mathbb{Z}^{d }\setminus\{0\},$ we have
						$$
						 \inf_{\lambda\in[\lambda_{0},\lambda_{1}]}\max_{q\in\llbracket 0, q_{0}\rrbracket}|\partial_{\lambda}^{q}\omega_{\textnormal{Eq}}(\lambda)\cdot l|\geqslant\rho_{0}\langle l\rangle.
						$$
						\item For any $ (l,j)\in(\mathbb{Z}^{d }\times\mathbb{N})\setminus\{(0,0)\}$ 
						$$
						\quad\inf_{\lambda\in[\lambda_{0},\lambda_{1}]}\max_{q\in\llbracket 0, q_{0}\rrbracket}\big|\partial_{\lambda}^{q}\big(\omega_{\textnormal{Eq}}(\lambda)\cdot l\pm j(I_{1}K_{1})(\lambda)\big)\big|\geqslant\rho_{0}\langle l\rangle.
						$$
						\item  For any $ (l,j)\in\mathbb{Z}^{d }\times (\mathbb{N}^*\setminus\mathbb{S})$ 
						$$
						\quad\inf_{\lambda\in[\lambda_{0},\lambda_{1}]}\max_{q\in\llbracket 0, q_{0}\rrbracket}\big|\partial_{\lambda}^{q}\left(\omega_{\textnormal{Eq}}(\lambda)\cdot l\pm\Omega_{j}(\lambda)\right)\big|\geqslant\rho_{0}\langle l\rangle.
						$$
						\item For any $ l\in\mathbb{Z}^{d }, j,j^\prime\in\mathbb{N}^*\setminus\mathbb{S}$  with $(l,j)\neq(0,j^\prime),$ we have
						$$\,\quad\inf_{\lambda\in[\lambda_{0},\lambda_{1}]}\max_{q\in\llbracket 0, q_{0}\rrbracket}\big|\partial_{\lambda}^{q}\big(\omega_{\textnormal{Eq}}(\lambda)\cdot l+\Omega_{j}(\lambda)\pm\Omega_{j^\prime}(\lambda)\big)\big|\geqslant\rho_{0}\langle l\rangle.$$	
				\end{enumerate}}
			\end{lem}
			\begin{proof} 
				\textbf{(i)} We argue  by contradiction by assuming that  for any $q_{0}\in\mathbb{N}$ and  $\rho_{0}>0$, there exist $l\in\mathbb{Z}^{d }\setminus\{0\}$ and $\lambda\in[\lambda_{0},\lambda_{1}]$ such that 
				$$
				\max_{q\in\llbracket 0, q_{0}\rrbracket}|\partial_{\lambda}^{q}(\omega_{\textnormal{Eq}}(\lambda)\cdot l)|<\rho_{0}\langle l\rangle.
				$$
				It follows that for any $m\in\mathbb{N}$, and  by  taking $q_{0}=m$ and $\rho_{0}=\frac{1}{m+1}$, there exist $l_{m}\in\mathbb{Z}^{d}\setminus\{0\}$ and $\lambda_{m}\in[\lambda_{0},\lambda_{1}]$ such that 
				$$
				\max_{q\in\llbracket 0, m\rrbracket}|\partial_{\lambda}^{q}\omega_{\textnormal{Eq}}(\lambda_{m})\cdot l_{m}|<\tfrac{\langle l_m\rangle}{m+1}
				$$
				and therefore 
				\begin{equation}\label{Rossemann 0}
					\forall q\in\mathbb{N},\quad \forall m\geqslant q, \quad \left|\partial_{\lambda}^{q}\omega_{\textnormal{Eq}}(\lambda_{m})\cdot \tfrac{l_{m}}{\langle l_m\rangle}\right|<\tfrac{1}{m+1}\cdot
				\end{equation}
				Since the sequences $\left(\frac{l_{m}}{\langle l_m\rangle}\right)_{m}$ and $(\lambda_{m})_{m}$ are bounded, then  by compactness and up to an extraction we can assume that
				$$\lim_{m\to\infty}\tfrac{l_{m}}{\langle l_m\rangle}=\bar{c}\neq 0\quad\hbox{and}\quad \lim_{m\to\infty}\lambda_{m}=\bar{\lambda}.
				$$
				Hence, passing to the limit  in \eqref{Rossemann 0}  as  $m\rightarrow\infty$ leads to
				$$\forall q\in\mathbb{N},\quad \partial_{\lambda}^{q}\omega_{\textnormal{Eq}}(\bar{\lambda})\cdot\bar{c}=0.$$
				Thus, we conclude that  the real analytic function $\lambda\mapsto\omega_{\textnormal{Eq}}(\lambda)\cdot\bar{c}$ is identically zero which contradicts the non-degeneracy condition stated in  Lemma \ref{Non-degenracy1}.\\
				%\textbf{(ii)}
				\\
				\textbf{(ii)} We shall first check the result for  the case $l=0$ and $j\in\mathbb{N}^*$. Obviously, one has from the monotonicity of $\lambda\mapsto I_{1}(\lambda)K_{1}(\lambda)$ stated in Appendix \ref{appendix Bessel},
				\begin{align*}
						\quad\inf_{\lambda\in[\lambda_{0},\lambda_{1}]}\max_{q\in\llbracket 0, q_{0}\rrbracket}\big|\partial_{\lambda}^{q}\big( j(I_{1}K_{1})(\lambda)\big)\big|&\geqslant (I_{1}K_{1})(\lambda_1)\\
						&\geqslant\rho_{0}\langle l\rangle,
						\end{align*}
for some $\rho_0>0.$ Now let us consider $l\in\mathbb{Z}^d\backslash\{0\}$ and $j\in\mathbb{N}.$	Then 			
				 we may  write according to  the triangle and Cauchy-Schwarz inequalities combined with the boundedness of $\omega_{\textnormal{Eq}}$ and the monotonicity of $\lambda\mapsto I_{1}(\lambda)K_{1}(\lambda)$ stated in Appendix \ref{appendix Bessel},
				$$|\omega_{\textnormal{Eq}}(\lambda)\cdot l\pm jI_{1}(\lambda)K_{1}(\lambda)|\geqslant jI_{1}(\lambda_{1})K_{1}(\lambda_{1})-|\omega_{\textnormal{Eq}}(\lambda)\cdot l|\geqslant c_{0}j-C\langle l\rangle\geqslant \langle l\rangle$$
				provided that $j\geqslant C_{0}\langle l\rangle$ for some $C_{0}>0.$ Therefore we reduce  the proof to indices $j$ and $l$ with 
				\begin{equation}\label{parameter condition I1K1}
					0\leqslant j< C_{0}\langle l\rangle,\quad j\in\mathbb{N}\quad  \hbox{and}\quad l\in\mathbb{Z}^d\backslash\{0\}.
				\end{equation}
				Arguing  by contradiction as in the previous case, we may assume the existence of sequences  $l_{m}\in\mathbb{Z}^{d }\setminus\{0\}$, $j_{m}\in\mathbb{N}$ satisfying \eqref{parameter condition I1K1} and $\lambda_{m}\in[\lambda_{0},\lambda_{1}]$ such that 
				$$\max_{q\in\llbracket 0, m\rrbracket}\left|\partial_{\lambda}^{q}\left(\omega_{\textnormal{Eq}}(\lambda)\cdot\tfrac{l_{m}}{|l_{m}|}\pm j_m\tfrac{(I_{1}K_{1})(\lambda)}{|l_{m}|}\right)_{|\lambda=\lambda_m}\right|<\tfrac{1}{m+1}$$ 
				and therefore
				\begin{equation}\label{Rossemann I1K1}
					\forall q\in\mathbb{N},\quad \forall m\geqslant q,\quad \left|\partial_{\lambda}^{q}\left(\omega_{\textnormal{Eq}}(\lambda)\cdot\tfrac{l_{m}}{|l_{m}|}\pm\tfrac{j_{m}}{|l_{m}|}(I_{1}K_{1})(\lambda)\right)_{|\lambda=\lambda_m}\right|<\tfrac{1}{m+1}\cdot
				\end{equation}
				Since the sequences $\left(\frac{l_{m}}{|l_{m}|}\right)_{m}$, $\left(\frac{j_{m}}{|l_{m}|}\right)_{m}$ and $(\lambda_{m})_{m}$ are bounded, then up to an extraction we can assume  that  
				$$
				\lim_{m\to\infty}\tfrac{l_{m}}{| l_{m}|}=\bar{c}\neq 0,\quad \lim_{m\to\infty}\tfrac{j_{m}}{| l_{m}|}=\bar{d}\quad\hbox{and}\quad \lim_{m\to\infty}\lambda_{m}=\bar{\lambda}.
				$$
				Hence, by letting $m\rightarrow\infty$ in \eqref{Rossemann I1K1}, using that $\lambda\mapsto (I_{1}K_{1})(\lambda)$ is smooth, we find
				$$\forall q\in\mathbb{N},\quad\partial_{\lambda}^{q}\big(\omega_{\textnormal{Eq}}({\lambda})\cdot\bar{c}\pm\bar{d}\,(I_{1}K_{1})({\lambda})\big)_{|\lambda=\overline\lambda}=0.$$
				Thus, the real analytic function $\lambda\mapsto\omega_{\textnormal{Eq}}(\lambda)\cdot\bar{c}\pm \bar{d}\,I_{1}(\lambda)K_{1}(\lambda)$ with $(\bar{c},\bar{d})\neq(0,0)$ is identically zero and this contradicts Lemma \ref{Non-degenracy1}.\\
				\\
				\textbf{(iii)} Consider $(l,j)\in\mathbb{Z}^{d }\times (\mathbb{N}^*\setminus\mathbb{S})$. Then applying  the  triangle inequality and Lemma \ref{lemma properties linear frequencies}-(iv), yields
				\begin{align*}
					|\omega_{\textnormal{Eq}}(\lambda)\cdot l\pm\Omega_{j}(\lambda)|&\geqslant|\Omega_{j}(\lambda)|-|\omega_{\textnormal{Eq}}(\lambda)\cdot l|\\
					&\geqslant \Omega j-C|l|\geqslant \langle l\rangle
				\end{align*}
				provided $j\geqslant C_{0} \langle l\rangle$ for some $C_{0}>0.$ Thus  as before we shall restrict the proof to indices $j$ and $l$ with 
				\begin{equation}\label{parameter condition 1}
					0\leqslant j< C_{0} \langle l\rangle,\quad j\in\mathbb{N}^*\setminus\mathbb{S}\quad\hbox{and}\quad l\in\mathbb{Z}^d\backslash\{0\}.
				\end{equation}
				Proceeding   by contradiction as in the previous case, we may assume the existence of sequences  $l_{m}\in\mathbb{Z}^{d }\setminus\{0\}$, $j_{m}\in\mathbb{N}\setminus\mathbb{S}$ satisfying \eqref{parameter condition 1} and $\lambda_{m}\in[\lambda_{0},\lambda_{1}]$ such that 
				$$
				\max_{q\in\llbracket 0, m\rrbracket}\ \left|\partial_{\lambda}^{q}\left(\omega_{\textnormal{Eq}}(\lambda)\cdot\tfrac{l_{m}}{| l_{m}|}\pm\tfrac{\Omega_{j_{m}}(\lambda)}{| l_{m}|}\right)_{|\lambda=\lambda_m}\right|<\tfrac{1}{m+1}
				$$ 
				and therefore
				\begin{equation}\label{Rossemann 1}
					\forall q\in\mathbb{N},\quad \forall m\geqslant q,\quad \left|\partial_{\lambda}^{q}\left(\omega_{\textnormal{Eq}}(\lambda)\cdot\tfrac{l_{m}}{| l_{m}|}\pm\tfrac{\Omega_{j_{m}}(\lambda)}{| l_{m}|}\right)_{|\lambda=\lambda_m}\right|<\tfrac{1}{m+1}\cdot
				\end{equation}
				Since  the sequences $\left(\tfrac{l_{m}}{|l_{m}|}\right)_{m}$ and $(\lambda_{m})_{m}$ are bounded, then up to an extraction we can assume  that 
				$$\lim_{m\to\infty}\tfrac{l_{m}}{| l_{m}|}=\bar{c}\neq 0\quad\hbox{and}\quad \lim_{m\to\infty}\lambda_{m}=\bar{\lambda}.
				$$
				Now we shall distinguish  two cases.\\
				$\blacktriangleright$ Case \ding{182} : $(l_{m})_{m}$ is bounded. In this case, by \eqref{parameter condition 1} we  find  that $(j_{m})_{m}$ is bounded too and thus up to to an extraction we may assume  $ \displaystyle \lim_{m\to\infty}l_{m}=\bar{l}$ and $  \displaystyle \lim_{m\to\infty}j_{m}=\bar{j}.$
				Since $(j_{m})_{m}$ and $(|l_{m}|)_{m}$ are sequences of integers, then they are necessary stationary. In particular, the condition \eqref{parameter condition 1} implies $\bar{l}\neq 0.$ Hence, taking the limit $n\rightarrow\infty$ in \eqref{Rossemann 1}, yields
				$$\forall q\in\mathbb{N},\quad \partial_{\lambda}^{q}\left(\omega_{\textnormal{Eq}}({\lambda})\cdot\bar{l}\pm \Omega_{\bar{j}}({\lambda})\right)_{|\lambda=\overline\lambda}=0.$$
				Thus, the analytic function $\lambda\mapsto\omega_{\textnormal{Eq}}(\lambda)\cdot\bar{l}\pm\Omega_{\bar{j}}(\lambda)$ with $(\bar{l},1)\neq (0,0)$ is identically zero which contradicts Lemma \ref{Non-degenracy1}.\\
				$\blacktriangleright$ Case \ding{183} : $(l_{m})_{m}$ is unbounded. Up to an extraction  we can assume that $\displaystyle \lim_{m\to\infty}|l_{m}|=\infty.$
				We have two sub-cases.\\
				$\bullet$ Sub-case \ding{172} : $(j_{m})_{m}$ is bounded. In this case and up to an extraction   we can assume that it converges. Then, taking the limit $m\rightarrow\infty$ in \eqref{Rossemann 1}, we find
				$$\forall q\in\mathbb{N},\quad \partial_{\lambda}^{q}\omega_{\textnormal{Eq}}(\bar{\lambda})\cdot\bar{c}=0.$$
				As before we conclude that  function $\lambda\mapsto\omega_{\textnormal{Eq}}(\lambda)\cdot\bar{c}$ with $\bar{c}\neq 0$ is identically zero which contradicts Lemma \ref{Non-degenracy1}.\\
				$\bullet$ Sub-case \ding{173} : $(j_{m})_{m}$ is unbounded. Then  up to an extraction we can assume that $  \displaystyle \lim_{m\to\infty}j_{m}=\infty$. We write  according to \eqref{def eigenvalues at equilibrium}
				\begin{align}\label{FHK1}
					\tfrac{\Omega_{j_{m}}(\lambda)}{| l_{m}|}=\tfrac{j_{m}}{| l_{m}|}\left(\Omega+(I_{1}K_{1})(\lambda)-(I_{j_{m}}K_{j_{m}})(\lambda)\right).
				\end{align}
				By \eqref{parameter condition 1}, the sequence $\left(\frac{j_{m}}{| l_{m}|}\right)_{n}$ is bounded, thus  up to an extraction  we can assume that it converges to $\bar{d}.$ Using the first inequality of \eqref{decayk} we deduce that
				$$
				\forall \, m\in\mathbb{N},\quad \sup_{\lambda\in\mathbb{R}}\big|(I_{j_{m}}K_{j_{m}})(\lambda)\big|\leq \tfrac{1}{2j_m},
				$$
				which implies that
				$$
				\lim_{m\to\infty} \sup_{\lambda\in\mathbb{R}}(I_{j_{m}}K_{j_{m}})(\lambda)=0.
				$$
				Moreover by \eqref{decay deriv InKn}, we have
				\begin{align}\label{KGA2}
					\lim_{m\to\infty}\sup_{\lambda\in[\lambda_0,\lambda_1]}\left|\partial_\lambda^q(I_{j_{m}}K_{j_{m}})(\lambda)\right|&=0.
				\end{align}
				Taking the limit in \eqref{FHK1} and using \eqref{KGA2} yields
				\begin{align*}
					\lim_{m\to\infty}\tfrac{\partial_\lambda^q\Omega_{j_{m}}(\lambda_{m})}{| l_{m}|}=\partial_\lambda^q\Big(\overline{d}\Big(\Omega+\big(I_{1}K_{1}\big)(\overline{\lambda})\Big)\Big).
				\end{align*}
				Consequently, taking the limit $m\rightarrow\infty$ in \eqref{Rossemann 1}, we have 
				$$\forall q\in\mathbb{N},\quad \partial_{\lambda}^{q}\left(\omega_{\textnormal{Eq}}({\lambda})\cdot\bar{c}\pm \bar{d}\big(\Omega+\big(I_{1}K_{1}\big)({\lambda})\big)\right)_{|\lambda=\lambda_m}=0.$$
				By continuation the  analytic function $\lambda\mapsto\omega_{\textnormal{Eq}}(\lambda)\cdot\bar{c}\pm \bar{d}(\Omega+I_{1}(\lambda)K_{1}(\lambda))$ with $(\bar{c},\bar{d})\neq 0$ is identically zero which contradicts  Lemma \ref{Non-degenracy1}.\\
				\\
				\textbf{(iv)} Consider $ l\in\mathbb{Z}^{d }, j,j^\prime\in\mathbb{N}^*\setminus\mathbb{S}$  with $(l,j)\neq(0,j^\prime).$ Then applying the  triangle inequality combined with Lemma \ref{lemma properties linear frequencies}- (v), we infer
				$$|\omega_{\textnormal{Eq}}(\lambda)\cdot l+\Omega_{j}(\lambda)\pm\Omega_{j'}(\lambda)|\geqslant|\Omega_{j}(\lambda)\pm\Omega_{j'}(\lambda)|-|\omega_{\textnormal{Eq}}(\lambda)\cdot l|\geqslant C_{0}|j\pm j'|-C|l|\geqslant \langle l\rangle$$
				provided that $|j\pm j'|\geqslant c_{0}\langle l\rangle$ for some $c_{0}>0.$ Then it remains to check the proof for  indices  satisfying 
				\begin{equation}\label{parameter condition 2}
					|j\pm j'|< c_{0}\langle l\rangle,\quad  l\in\mathbb{Z}^{d }\backslash\{0\}\quad\textnormal{and}\quad j,j^\prime\in\mathbb{N}^*\setminus\mathbb{S}.
				\end{equation}
				Reasoning by contradiction as in the previous cases, we get   for all $m\in\mathbb{N}$, real numbers  $l_{m}\in\mathbb{Z}^{d }\setminus\{0\}$, $j_{m},j'_{m}\in\mathbb{N}^*\setminus\mathbb{S}$ satisfying \eqref{parameter condition 2} and $\lambda_{m}\in[\lambda_{0},\lambda_{1}]$ such that 
				$$
				\max_{q\in\llbracket 0, m\rrbracket}\left|\partial_{\lambda}^{q}\left(\omega_{\textnormal{Eq}}(\lambda)\cdot\tfrac{l_{m}}{| l_{m}|}+\tfrac{\Omega_{j_{m}}(\lambda)\pm \Omega_{j'_{m}}(\lambda)}{| l_{m}|}\right)_{|\lambda=\lambda_m}\right|<\tfrac{1}{m+1}
				$$ 
				implying in turn that
				\begin{equation}\label{Rossemann 2}
					\forall q\in\mathbb{N},\quad \forall m\geqslant q,\quad \left|\partial_{\lambda}^{q}\left(\omega_{\textnormal{Eq}}(\lambda)\cdot\tfrac{l_{m}}{| l_{m}|}+\tfrac{\Omega_{j_{m}}(\lambda)\pm \Omega_{j'_{m}}(\lambda)}{| l_{m}|}\right)_{|\lambda=\lambda_m}\right|<\tfrac{1}{m+1}\cdot
				\end{equation}
				Up to an extraction we can assume that $\displaystyle \lim_{m\to\infty}\tfrac{l_{m}}{| l_{m}|}=\bar{c}\neq 0$ and $\displaystyle  \lim_{m\to\infty}\lambda_{m}=\bar{\lambda}.$ \\
				As before we shall distinguish two cases.\\
				$\blacktriangleright$ Case \ding{182} : $(l_{m})_{m}$ is bounded. We shall only focus on the most delicate case associated to the difference $\Omega_{j_m}-\Omega_{j^\prime_m}$. Up to an extraction we may assume  that  $\displaystyle \lim_{m\to\infty}l_{m}=\bar{l}\neq 0.$ Now according  to  \eqref{parameter condition 2} we have two sub-cases to discuss depending whether the sequences $(j_{m})_{m}$ and $(j'_{m})_{m}$ are simultaneously bounded  or unbounded.\\
				$\bullet$ Sub-case \ding{172} : $(j_{m})_{m}$ and $(j'_{m})_{m}$ are bounded. In this case, up to an extraction we may assume that these sequences are stationary  $j_{m}=\bar{j}$ and $j'_{m}=\bar{j'}$ with $ \bar{j},\bar{j^\prime}\in\mathbb{N}^*\setminus\mathbb{S}.$
				Hence taking the limit as $m\rightarrow\infty$ in \eqref{Rossemann 2}, we infer
				$$
				\forall q\in\mathbb{N},\quad \partial_{\lambda}^{q}\left(\omega_{\textnormal{Eq}}({\lambda})\cdot\bar{l}+\Omega_{\bar{j}}(\bar{\lambda})-\Omega_{\bar{j}^\prime}({\lambda})\right)_{\lambda=\overline\lambda}=0.
				$$
				Thus, the analytic function $\lambda\mapsto\omega_{\textnormal{Eq}}(\lambda)\cdot\bar{l}+\Omega_{\bar{j}}(\lambda)-\Omega_{\bar{j'}}(\lambda)$ is identically zero. If $\bar{j}=\bar{j^\prime}$ then this  contradicts Lemma \ref{Non-degenracy1} since $\bar{l}\neq 0.$ However in the case  $\bar{j}\neq\bar{j^\prime}\in\mathbb{N}^*\setminus\mathbb{S}$ this still contradicts this lemma applied with the vector frequency $(\omega_{\textnormal{Eq}},\Omega_{\bar{j}},\Omega_{\bar{j'}})$ instead of $\omega_{\textnormal{Eq}}.$\\
				$\bullet$ Sub-case \ding{173} : $(j_{m})_{m}$ and $(j'_{m})_{m}$ are both unbounded and without loss of generality we can assume that $\displaystyle \lim_{m\to\infty}j_{m}= \lim_{m\to\infty}j'_{m}=\infty$. From \eqref{est deriv diff IjKj-1} combined with  \eqref{parameter condition 2} and the boundedness of $(l_{m})_{m}$ we deduce that
				\begin{align*}
					\left|\partial_\lambda^q(I_{j_m}K_{j_m}-I_{j_m^\prime}K_{j_m^\prime})(\lambda_m)\right|&\leqslant \tfrac{C}{j_m j_m^\prime},
				\end{align*}
				which implies in turn
				\begin{align}\label{rdd11}
					\lim_{m\to\infty} j_m^\prime\,\partial_\lambda^q(I_{j_m}K_{j_m}-I_{j_m^\prime}K_{j_m^\prime})(\lambda_m)
					=0.
				\end{align}
				Coming back to  \eqref{def eigenvalues at equilibrium} we get the splitting 
				\begin{align}\label{Appla1}
					\nonumber\Omega_{j_m}(\lambda)-\Omega_{j_{m}^\prime}(\lambda)=&(j_m-j^\prime_{m})\big(\Omega+(I_{1}K_{1})(\lambda)\big)-(j_m-j^\prime_{m})(I_{j_m}K_{j_m})(\lambda)\\
					&+j^\prime_{m}\Big((I_{j^\prime_{m}}K_{j^\prime_{m}})(\lambda)-(I_{j_m}K_{j_m})(\lambda)\Big).
				\end{align}
				Therefore by applying \eqref{KGA2} and \eqref{rdd11} we get for any $q\in\mathbb{N},$
				\begin{align*}
					\lim_{m\to\infty}
					\partial_\lambda^q\Big(\Omega_{j_m}(\lambda)-\Omega_{j_{m}^\prime}(\lambda)-&(j_m-j^\prime_{m})\big(\Omega+(I_{1}K_{1})(\lambda)\big)\Big)_{\lambda=\lambda_m}=0.
				\end{align*}
				Using once again \eqref{parameter condition 2} and up to an extraction we have   $\displaystyle \lim_{m\to\infty}\tfrac{j_{m}-j'_{m}}{|l_m|}=\bar{d}.$
				Thus
				$$
				\lim_{m\to\infty}|l_m|^{-1}\partial_{\lambda}^{q}\left({\Omega_{j_{m}}(\lambda)-\Omega_{j'_{m}}(\lambda)}\right)_{|\lambda=\lambda_m}=\bar{d}\,\partial_{\lambda}^{q}\,\left(\Omega+(I_{1}K_{1})({\lambda})\right)_{|\lambda=\overline\lambda}.$$
				By taking the limit as $m\rightarrow\infty$ in \eqref{Rossemann 2}, we find				$$
				\forall q\in\mathbb{N},\quad \partial_{\lambda}^{q}\left({\omega}_{\textnormal{Eq}}({\lambda})\cdot\bar{c}+\bar{d}\big(\Omega+(I_{1}K_{1})({\lambda})\big)\right)_{|\lambda=\overline\lambda}=0.
				$$
				Thus, the analytic function $\lambda\mapsto{\omega}_{\textnormal{Eq}}(\lambda)\cdot\bar{c}+\bar{d}(\Omega+I_{1}(\lambda)K_{1}(\lambda))$ with $(\bar{c},\bar{d})\neq 0$ is vanishing  which contradicts  Lemma \ref{Non-degenracy1}. Now we shall move to the second case.\\
				%Case \ding{183}
				$\blacktriangleright$  Case \ding{183} : $(l_{m})_{m}$ is unbounded. Up to an extraction  we can assume that $\displaystyle \lim_{m\to\infty}|l_{m}|=\infty.$\\
				We shall distinguish three sub-cases.\\
				$\bullet$ Sub-case \ding{172}. The sequences  $(j_{m})_{m}$ and $(j'_{m})_{m}$ are bounded. In this case and  up to an extraction  they will  converge and then taking the limit in \eqref{Rossemann 2} yields,
				$$\forall q\in\mathbb{N},\quad\partial_{\lambda}^{q}{\omega}_{\textnormal{Eq}}(\bar{\lambda})\cdot\bar{c}=0.
				$$
				which leads to a contradiction as before. \\
				$\bullet$ Sub-case \ding{173}. The sequences  $(j_{m})_{m}$ and $(j'_{m})_{m}$ are both unbounded. This is similar to  the sub-case \ding{173} of the case \ding{182}.\\
				$\bullet$ Sub-case \ding{174}. The sequence $(j_{m})_{m}$ is unbounded and $(j'_{m})_{m}$ is bounded (the symmetric case is similar).  Without loss of generality we can assume that $\displaystyle \lim_{m\to\infty}j_m=\infty$ and $  j_{m}^\prime=\overline{j}.$ By \eqref{parameter condition 2} and  up to an extraction one gets  $\displaystyle \lim_{m\to\infty}\tfrac{j_{m}\pm j'_{m}}{| l_{m}|}=\bar{d}.$ One may use \eqref{def eigenvalues at equilibrium}
				combined with  \eqref{KGA2} and \eqref{rdd11} in order to get for any $q\in\mathbb{N},$
				\begin{align*}
					\lim_{m\to\infty}|l_m|^{-1}
					\partial_\lambda^q\Big(\Omega_{j_m}(\lambda)\pm\Omega_{j_{m}^\prime}(\lambda)-(j_m\pm j^\prime_{m})\big(\Omega+(I_{1}K_{1})(\lambda)\big)\Big)_{|\lambda=\lambda_m}&=\\
					\lim_{m\to\infty}
					\partial_\lambda^q\left(\tfrac{(j_m\pm j^\prime_m)}{|l_m|}(I_{j_{m}} K_{j_{m}})(\lambda)\pm \tfrac{j^\prime_{m}}{|l_m|}\Big((I_{j_{m}} K_{j_{m}})(\lambda)-(I_{j_{m}^\prime}  K_{j_{m}^\prime})(\lambda)\right)_{|\lambda=\lambda_m}&=0.
				\end{align*}
				Hence, taking the limit  in \eqref{Rossemann 2} implies 
				$$\forall q\in\mathbb{N},\partial_{\lambda}^{q}\left({\omega}_{\textnormal{Eq}}({\lambda})\cdot\bar{c}+\bar{d}\big(\Omega+(I_{1}K_{1})({\lambda})\big)\right)_{\lambda=\overline\lambda}=0.$$
				Thus, the analytic function $\lambda\mapsto{\omega}_{\textnormal{Eq}}(\lambda)\cdot\bar{c}+\bar{d}(\Omega+I_{1}(\lambda)K_{1}(\lambda))$ is identically zero with $(\bar{c},\bar{d})\neq0$ which contradicts Lemma \ref{Non-degenracy1}.This completes the proof of Lemma \ref{lemma transversality}.
			\end{proof}
		\subsubsection{Linear quasi-periodic solutions}
		Notice that all the solutions of
		\eqref{hamiltonian equation at equilibrium} taking  the form \eqref{sol eq} are either periodic, quasi-periodic or almost periodic in time, with linear frequencies of oscillations
		$\Omega_{j}(\lambda)$  defined by \eqref{def eigenvalues at equilibrium} These different notions differs on the irrationality properties of the frequencies $\{\Omega_{j}(\lambda)\}$ and on the cardinality of the Fourier-space support (finite for quasi-periodic functions and possibly infinite for almost periodic ones). Remark that we have the implications
		$$\textnormal{Periodic}\quad\Rightarrow\quad\textnormal{Quasi-periodic}\quad\Rightarrow\quad\textnormal{Almost periodic}.$$ 
		We shall prove here the existence of quasi-periodic solutions for the linear equation \eqref{hamiltonian equation at equilibrium} when $\lambda$ belongs to a massive Cantor set. 	\begin{prop}\label{lemma sol Eq}
			Let $\lambda_1>\lambda_0> 0,$ $d\in\mathbb{N}^{*}$ and $\mathbb{S}\subset\mathbb{N}^*$ with $|\mathbb{S}|=d.$ Then, there exists a Cantor-like set $\mathcal{C}\subset[\lambda_0,\lambda_1]$ satisfying $|\mathcal{C}|=\lambda_1-\lambda_0$ and such that for all $\lambda\in\mathcal{C}$, every function in the form
			\begin{equation}\label{rev sol eq}
				\rho(t,\theta)=\sum_{j\in\mathbb{S}}\rho_{j}\cos(j\theta-\Omega_{j}(\lambda)t),\quad\rho_{j}\in\mathbb{R}^*
			\end{equation} 
			is a time quasi-periodic reversible solution to the equation \eqref{hamiltonian equation at equilibrium} with the vector frequency 
			$$\omega_{\textnormal{Eq}}(\lambda)=(\Omega_{j}(\lambda))_{j\in\mathbb{S}}.$$
					\end{prop}
	\begin{proof}
		It is easy to check that any function in the form \eqref{rev sol eq} is  a reversible solution to \eqref{hamiltonian equation at equilibrium}, that is a solution satisfying the property  $$r(-t,-\theta)=r(t,\theta).$$
%		
%		Given a real valued solution $\rho$ of \eqref{hamiltonian equation at equilibrium}, one can write in view of Lemma \ref{lemma linearized operator at equilibrium} 
%			$$\rho(t,\theta)=\sum_{j\in\mathbb{Z}^{*}}\rho_{j}(0)e^{i(j\theta-\Omega_{j}(\lambda)t)}\quad\textnormal{with}\quad \rho_{-j}(0)=\overline{\rho_{j}(0)}.$$
%			Using Lemma \ref{lemma properties linear frequencies}-(i), we deduce that
%		\begin{align*}
%			\rho(-t,-\theta)&=\sum_{j\in\mathbb{Z}^{*}}\rho_{j}(0)e^{i(-j\theta+\Omega_{j}(\lambda)t)}\\
%			&=\sum_{j\in\mathbb{Z}^{*}}\rho_{j}(0)e^{i(-j\theta-\Omega_{-j}(\lambda)t)}\\
%			&=\sum_{j\in\mathbb{Z}^{*}}\rho_{-j}(0)e^{i(j\theta-\Omega_{j}(\lambda)t)}.
%	\end{align*}
%We then conclude, by uniqueness of Fourier coefficients, that
%$$\rho(-t,-\theta)=\rho(t,\theta)\quad \Leftrightarrow\quad \forall j\in\mathbb{Z}^{*},\,\,\rho_{j}(0)\in\mathbb{R}.$$
%Thus, any real valued solution of \eqref{hamiltonian equation at equilibrium} satisfying the reveribility condition $\rho(-t,-\theta)=\rho(t,\theta)$ writes
%\begin{align*}
%	\rho(t,\theta)&=\frac{1}{2}\Big(\rho(t,\theta)+\rho(-t,-\theta)\Big)\\
%	&=\sum_{j\in\mathbb{Z}^{*}}\rho_{j}(0)\frac{e^{i(j\theta-\Omega_{j}(\lambda)t)}+e^{-i(j\theta-\Omega_{j}(\lambda)t)}}{2}\\
%	&=\sum_{j\in\mathbb{Z}^{*}}\rho_{j}(0)\cos(j\theta-\Omega_{j}(\lambda)t).
%\end{align*} 
%		As a consequence, any function $\rho$ in the form \eqref{rev sol eq} is a reversible solution of \eqref{hamiltonian equation at equilibrium}. 
		Then, it remains to check the non-resonance condition \eqref{nonresonnace omega} for the frequency vector $\omega_{\textnormal{Eq}}$ for almost every  $\lambda\in[\lambda_0,\lambda_1]$.
%		 We already know from Lemma \ref{Non-degenracy1} that there exists at least one value of $\lambda$ such that the condition \eqref{nonresonnace omega} is satisfied. But we can obtain a better result and show that it happens for values of $\lambda$ in a good Cantor set with almost full Lebesgue measure in $(\lambda_0,\lambda_1)$.
		  For that purpose, we consider $\tau_1>0, \gamma\in(0,1)$ and define the set $\mathcal{C}_{\gamma}$ by
		$$\mathcal{C}_{\gamma}=\bigcap_{l\in\mathbb{Z}^{d}\setminus\{0\}}\left\lbrace\lambda\in[\lambda_0,\lambda_1]\quad\textnormal{s.t.}\quad \left|\omega_{\textnormal{Eq}}(\lambda)\cdot l\right|>\tfrac{\gamma}{\langle l\rangle^{\tau_{1}}}\right\rbrace.$$
		Therefore  its complement  set takes the form 
		$$[\lambda_0,\lambda_1]\setminus\mathcal{C}_{\gamma}=\bigcup_{l\in\mathbb{Z}^{d}\backslash\{0\}}\mathcal{R}_{l}\quad \textnormal{where}\quad \mathcal{R}_{l}=\left\lbrace\lambda\in[\lambda_0,\lambda_1]\quad\textnormal{s.t.}\quad \left|\omega_{\textnormal{Eq}}(\lambda)\cdot l\right|\leqslant\tfrac{\gamma}{\langle l\rangle^{\tau_{1}}}\right\rbrace.
		$$	
	It follows that 
		$$
		\Big|[\lambda_0,\lambda_1]\setminus\mathcal{C}_\gamma\Big|\leqslant\sum_{l\in\mathbb{Z}^{d}\backslash\{0\}}\left|\mathcal{R}_{l}\right|.$$
		Now applying Lemma \ref{lemma useful for measure estimates} together with Lemma \ref{lemma transversality}-(i), one obtains
		$$\left|\mathcal{R}_{l}\right|\lesssim\gamma^{\frac{1}{q_0}}\langle l\rangle^{-1-\frac{\tau_{1}+1}{q_0}}.$$
		Then by imposing  
		$$\tau_{1}>(d-1)q_0-1,$$
		one  gets a convergent   series  with 
		$$\big|[\lambda_0,\lambda_1]\setminus\mathcal{C}_\gamma\big|\leqslant C\gamma^{\frac{1}{q_0}}.$$
		Now, we define the Cantor set 
		$${\mathcal{C}=\displaystyle\bigcup_{\gamma>0}\mathcal{C}_\gamma.}$$
		Then one gets easily for any $\gamma>0$
		$$
		\lambda_1-\lambda_0-C\gamma^{\frac{1}{q_0}}\leqslant |\mathcal{C}_\gamma|\leqslant |\mathcal{C}|\leqslant \lambda_1-\lambda_0 .
		$$
	Passing to the limit as $\gamma\to0$ yields
	$$
	 |\mathcal{C}|= \lambda_1-\lambda_0,
	$$	
	which achieves the proof of Proposition \ref{lemma sol Eq}.	
	\end{proof}

		In the previous proof, we used the following Lemma whose proof can be found in \cite[Thm. 17.1]{R01}. Notice that in all the paper, we use the notation $|A|$ as the Lebesgue measure of a given measurable set $A$.
		\begin{lem}\label{lemma useful for measure estimates}
			Let $q_{0}\in\mathbb{N}^{*}$ and $\alpha,\beta\in\mathbb{R}_{+}.$ Let $f\in C^{q_{0}}([a,b],\mathbb{R})$ such that
			$$\inf_{x\in[a,b]}\max_{k\in\llbracket 0,q_{0}\rrbracket}|f^{(k)}(x)|\geqslant\beta.$$
			Then, there exists $C=C(a,b,q_{0},\| f\|_{C^{q_{0}}([a,b],\mathbb{R})})>0$ such that 
			$$\Big|\left\lbrace x\in[a,b]\quad\textnormal{s.t.}\quad |f(x)|\leqslant\alpha\right\rbrace\Big|\leqslant C\tfrac{\alpha^{\frac{1}{q_{0}}}}{\beta^{1+\frac{1}{q_{0}}}}\cdot$$
		\end{lem}
	
			\section{Functional setting and technical Lemmas}\label{sec funct set}
			In this section, we set up  the general topological framework for both the functions and the  operators classes. We also provide  some   classical results on the law products, composition rule, Toeplitz operators, etc\ldots \\
			Next, we intend to introduce some parameters with some restrictions that will be used later.
			\begin{equation}\label{initial parameter condition}
				\gamma\in(0,1),\quad q,d\in\mathbb{N}^{*},\quad S\geqslant s\geqslant s_{0}>\tfrac{d+1}{2}+q+2,
			\end{equation}
			where $S$ is a fixed large number. 
			\begin{equation}\label{setting tau1 and tau2}
				\tau_{2}>\tau_{1}>d.
			\end{equation}
			Let 
			$$0<\lambda_{0}<\lambda_{1}.$$
			Since the mapping  $\omega_{\textnormal{Eq}}$, defined by \eqref{def linear frequency vector}, is  continuous, then we can find  a radius $R_{0}>0$ such that
			$$\omega_{\textnormal{Eq}}\left((\lambda_{0},\lambda_{1})\right)\subset\mathscr{U}:=B(0,R_{0}).$$ We consider $\mathcal{O}$ the open bounded subset of $\mathbb{R}^{d+1}$ defined by
			\begin{equation}\label{def initial parameters set}
				\mathcal{O}=(\lambda_{0},\lambda_{1})\times \mathscr{U}.
			\end{equation}
			Now, we explain the role played by the foregoing parameters throughout this paper. 
			\begin{remark}
				\begin{enumerate}[label=\textbullet]
					\item The parameter $\lambda$ comes from the model   $(QGSW)_{\lambda}$ and it is free in a fixed  interval $(\lambda_{0},\lambda_{1}).$ However at the end it  will belong to a Cantor set  for which invariant  torus can be constructed.
					\item The integer $d$ is the number of excited frequencies that will generate the  quasi-periodic solutions. This is the dimension of the space where lies the  frequency vector  $\omega\in \mathscr{U}\subset\mathbb{R}^{d},$ that will be a perturbation of the equilibrium frequency ${\omega}_{\textnormal{\tiny{Eq}}}(\lambda)$.
					\item The real number $s$ is the Sobolev index  regularity of the functions in the variables $\varphi$ and $\theta.$ The index $s$ will vary between $s_{0}$ and a large enough  parameter $S$ and at this end of Nash-Moser scheme it will be fixed as a large number related to the geometry of the intermediate  Cantor sets.
					\item The integer  $q$ is the index of regularity of our functions/operators with respect to the parameters $\lambda$ and $\omega.$ We have to consider such regularity in order to perform measure estimates in Section $\ref{Section measure of the final Cantor set}$ by checking the  R\"ussemann conditions. Its value will be fixed equal to $q_{0}+1$, where $q_0$ is the non degeneracy index of the tangential  frequencies given in Lemma $\ref{lemma transversality}.	$
									\item  All the remaining parameters $\gamma$, $\tau_{1}$ and $\tau_{2}$  are linked to  different Diophantine conditions, see for instance  Lemma $\ref{L-Invert0}$ and Propositions $\ref{reduction of the transport part}$ and $\ref{reduction of the remainder term}.$ The choice of $\tau_{1}$ and $\tau_{2}$ will be finally fixed in \eqref{choice tau 1 tau2 upsilon}. We point out that the parameter  $\gamma$ appears in the weighted Sobolev spaces and will be fixed in Proposition $\ref{Nash-Moser}$ with respect to  the rescaling parameter $\varepsilon$ giving  the smallness  condition of the solutions around the equilibrium.
				\end{enumerate}
			\end{remark}
			\subsection{Function spaces}\label{Functionspsaces}
			We shall  introduce the function spaces  that will be frequently used  along the paper. They are given by   weighted Sobolev spaces with respect to a parameter   $\gamma\in(0,1)$ used in defining the Cantor sets  to track the regularity with respect to the external  parameters $\lambda$ and $\omega$ of the solutions to  the nonlinear equation. 
%			 and . They take into acount a $C^{q}$ regularity with respect to the parameters $(\lambda,\omega)\in\mathcal{O}$ and a Sobolev regularity with respect to the time-space variables $(\varphi,\theta)\in\mathbb{T}^{d+1}.$
			We denote by  $(\mathbf{e}_{l,j})_{(l,j)\in\mathbb{Z}^{d}\times\mathbb{Z}}$ the Hilbert basis of the complex Hilbert space  $L^{2}(\mathbb{T}^{d+1},\mathbb{C})$ defined by
			$$\mathbf{e}_{l,j}(\varphi,\theta)=e^{\ii(l\cdot\varphi+j\theta)}.$$
			We endow this space   with the Hermitian inner  product
			\begin{align*}
				\big\langle\rho_{1},\rho_{2}\big\rangle_{L^{2}(\mathbb{T}^{d+1},\mathbb{C})}&=\bigintssss_{\mathbb{T}^{d+1}}\rho_{1}(\varphi,\theta)\overline{\rho_{2}(\varphi,\theta)}d\varphi d\theta\\
				&=\frac{1}{(2\pi)^{d+1}}\bigintssss_{[0,2\pi]^{d+1}}\rho_{1}(\varphi,\theta)\overline{\rho_{2}(\varphi,\theta)}d\varphi d\theta.
			\end{align*}
			To get  the last line of the preceding identity we use the notation \eqref{Convention1}. Given 
			 $\rho\in L^{2}(\mathbb{T}^{d+1},\mathbb{C}),$ we may decompose it  in Fourier expansion  as 
			$$\rho=\sum_{(l,j)\in\mathbb{Z}^{d+1 }}\rho_{l,j}\,\mathbf{e}_{l,j}\quad \mbox{ where }\quad \rho_{l,j}=\big\langle\rho,\mathbf{e}_{l,j}\big\rangle_{L^{2}(\mathbb{T}^{d+1},\mathbb{C})}.%=\int_{\mathbb{T}^{d+1}}\rho(\varphi,\theta)\mathbf{e}_{-l,-j}(\varphi,\theta)d\varphi d\theta\in\mathbb{C}.
			$$
			Next, we introduce for $s\in\mathbb{R}$ the complex  Sobolev space $H^{s}(\mathbb{T}^{d +1},\mathbb{C})$   by 
			$$
			H^{s}(\mathbb{T}^{d +1},\mathbb{C})=\Big\{\rho\in L^{2}(\mathbb{T}^{d +1},\mathbb{C})\quad\textnormal{s.t.}\quad\| \rho\|_{H^{s}}^{2}:=\sum_{(l,j)\in\mathbb{Z}^{d+1 }}\langle l,j\rangle^{2s}|\rho_{l,j}|^{2}<\infty\Big\},
			$$
			where $\langle l,j\rangle:=\max(1,|l|,|j|)$ with $|\cdot|$ denoting either the $\ell^{1}$ norm in $\mathbb{R}^{d }$ or the absolute value in $\mathbb{R}.$\\
			The real Sobolev spaces can be viewed as closed sub-spaces of the preceding one, 
			\begin{align*}
				\nonumber H^{s}=H^{s}(\mathbb{T}^{d+1},\mathbb{R})&=\left\lbrace\rho\in H^{s}(\mathbb{T}^{d +1},\mathbb{C})\quad\textnormal{s.t.}\quad\forall\, (\varphi,\theta)\in\mathbb{T}^{d+1 },\,\rho(\varphi,\theta)=\overline{\rho(\varphi,\theta)}\right\rbrace\\
				&=\Big\{\rho\in H^{s}(\mathbb{T}^{d +1},\mathbb{C})\quad\textnormal{s.t.}\quad\forall \,(l,j)\in\mathbb{Z}^{d+1},\,\rho_{-l,-j}=\overline{\rho_{l,j}}\Big\}.
			\end{align*}
			We shall also make use of the following subspaces of $H^{s}$ taking into account of some particular symmetries on odd and even functions,
			\begin{align*}H_{\mbox{\tiny{even}}}^{s}&=\Big\{\rho\in H^{s}\quad\textnormal{s.t.}\quad\forall\, (\varphi,\theta)\in\mathbb{T}^{d+1 },\,\rho(-\varphi,-\theta)=\rho(\varphi,\theta)\Big\}\\
				&=\Big\{\rho\in H^{s}\quad\textnormal{s.t.}\quad \forall\, (l,j)\in\mathbb{Z}^{d+1 },\,\rho_{-l,-j}=\rho_{l,j}\Big\}
			\end{align*}
			and
			\begin{align*}H_{\mbox{\tiny{odd}}}^{s}&=\Big\{\rho\in H^{s}\quad\textnormal{s.t.}\quad\forall\, (\varphi,\theta)\in\mathbb{T}^{d+1 },\,\rho(-\varphi,-\theta)=-\rho(\varphi,\theta)\Big\}\\
				&=\Big\{\rho\in H^{s}\quad\textnormal{s.t.}\quad\forall \,(l,j)\in\mathbb{Z}^{d+1 },\, \rho_{-l,-j}=-\rho_{l,j}\Big\}.
			\end{align*}
			For $N\in\mathbb{N}^{*},$ we define the cut-off frequency projectors on  $H^{s}(\mathbb{T}^{d +1},\mathbb{C})$ as follows 
			\begin{equation}\label{definition of projections for functions}
				\Pi_{N}\rho=\sum_{\underset{\langle l,j\rangle\leqslant N}{(l,j)\in\mathbb{Z}^{d+1 }}}\rho_{l,j}\mathbf{e}_{l,j}\quad \mbox{ and }\quad \Pi^{\perp}_{N}=\textnormal{Id}-\Pi_{N}.
			\end{equation}
			We shall also make use of the following mixed weighted Sobolev spaces. %Let $\mathcal{O}$ be an open bounded set of $\mathbb{R}^{d+1}$ and define   
			\begin{align*}
			W^{q,\infty,\gamma}(\mathcal{O},H^{s})&=\Big\lbrace \rho:\mathcal{O}\rightarrow H^{s}\quad\textnormal{s.t.}\quad\|\rho\|_{q,s}^{\gamma,\mathcal{O}}<\infty\Big\rbrace,\\
				W^{q,\infty,\gamma}(\mathcal{O},\mathbb{C})&=\Big\lbrace\rho:\mathcal{O}\rightarrow\mathbb{C}\quad\textnormal{s.t.}\quad\|\rho\|_{q}^{\gamma,\mathcal{O}}<\infty\Big\rbrace,
			\end{align*}
			where $\mu\in\mathcal{O}\mapsto \rho(\mu)\in H^s$ and 
			\begin{align}\label{Norm-def}
				\|\rho\|_{q,s}^{\gamma,\mathcal{O}}&=\sum_{\underset{|\alpha|\leqslant q}{\alpha\in\mathbb{N}^{d+1}}}\gamma^{|\alpha|}\sup_{\mu\in{\mathcal{O}}}\|\partial_{\mu}^{\alpha}\rho(\mu,\cdot)\|_{H^{s-|\alpha|}},\nonumber\\
				\|\rho\|_{q}^{\gamma,\mathcal{O}}&=\sum_{\underset{|\alpha|\leqslant q}{\alpha\in\mathbb{N}^{d+1}}}\gamma^{|\alpha|}\sup_{\mu\in{\mathcal{O}}}|\partial_{\mu}^{\alpha}\rho(\mu)|.
			\end{align}
			Note that a function $\rho\in W^{q,\infty,\gamma}(\mathcal{O},H^{s})$ can be written in the form 
			$$\rho(\mu,\varphi,\theta)=\sum_{(l,j)\in\mathbb{Z}^{d+1}}\rho_{l,j}(\mu)\mathbf{e}_{l,j}(\varphi,\theta).$$
			\begin{remark}
				\begin{enumerate}[label=\textbullet]
					\item From  Sobolev embeddings, we obtain
					$$W^{q,\infty,\gamma}(\mathcal{O},H^{s})\hookrightarrow C^{q-1}(\mathcal{O},H^{s})\quad \mbox{ and }\quad W^{q,\infty,\gamma}(\mathcal{O},\mathbb{C})\hookrightarrow C^{q-1}(\mathcal{O},\mathbb{C}).$$
					\item The spaces  $\left(W^{q,\infty,\gamma}(\mathcal{O},H^{s}),\|\cdot\|_{q,s}^{\gamma,\mathcal{O}}\right)$ and $\left(W^{q,\infty,\gamma}(\mathcal{O},\mathbb{C}),\|\cdot\|_{q}^{\gamma,\mathcal{O}}\right)$ are complete.
				\end{enumerate}
			\end{remark}
			In the next lemma we  collect some useful classical results dealing with  various operations in  weighted  Sobolev spaces. The proofs are very close to those in \cite{BFM21,BFM21-1,BM18}, so we omit them.
			\begin{lem}\label{Lem-lawprod}%[classical properties of the norm $\|\cdot\|_{q,s}^{\gamma,\mathcal{O}}$]
				{Let $(\gamma,q,d,s_{0},s)$ satisfying \eqref{initial parameter condition}, then the following assertions hold true.
					\begin{enumerate}[label=(\roman*)]
						\item Space translation invariance:  Let $\rho\in W^{q,\infty,\gamma}(\mathcal{O},H^{s}),$ then for all $\eta\in\mathbb{T},$ the function $(\varphi,\theta)\mapsto\rho(\varphi,\eta+\theta)$ belongs to $W^{q,\infty,\gamma}(\mathcal{O},H^{s})$, and satisfies
						$$\|\rho(\cdot,\eta+\cdot)\|_{q,s}^{\gamma,\mathcal{O}}=\|\rho\|_{q,s}^{\gamma,\mathcal{O}}.$$
						\item Projectors properties: Let $\rho\in W^{q,\infty,\gamma}(\mathcal{O},H^{s}),$ then for all $N\in\mathbb{N}^{*}$ and for all $t\in\mathbb{R}_{+}^{*},$
						$$\|\Pi_{N}\rho\|_{q,s+t}^{\gamma,\mathcal{O}}\leqslant N^{t}\|\rho\|_{q,s}^{\gamma,\mathcal{O}}\quad\mbox{ and }\quad\|\Pi_{N}^{\perp}\rho\|_{q,s}^{\gamma,\mathcal{O}}\leqslant N^{-t}\|\rho\|_{q,s+t}^{\gamma,\mathcal{O}},
						$$
						where the projectors are defined in \eqref{definition of projections for functions}.
						\item Interpolation inequality: 
						%\begin{enumerate}[label=(\alph*)]
						%\item 
						Let $q<s_{1}\leqslant s_{3}\leqslant s_{2}$ and $\overline{\theta}\in[0,1],$ with  $s_{3}=\overline{\theta} s_{1}+(1-\overline{\theta})s_{2}.$\\
						If $\rho\in W^{q,\infty,\gamma}(\mathcal{O},H^{s_{2}})$, then  $\rho\in W^{q,\infty,\gamma}(\mathcal{O},H^{s_{3}})$ and
						$$\|\rho\|_{q,s_{3}}^{\gamma,\mathcal{O}}\lesssim\left(\|\rho\|_{q,s_{1}}^{\gamma,\mathcal{O}}\right)^{\overline{\theta}}\left(\|\rho\|_{q,s_{2}}^{\gamma,\mathcal{O}}\right)^{1-\overline{\theta}}.$$
						%\textcolor{blue}{\item Let $0\leqslant q_{1}\leqslant q_{3}\leqslant q_{2}<s$ and $t\in[0,1]$ with $q_{3}=tq_{1}+(1-t)q_{2}.$ \\
						%If $\rho\in W^{q_{2},\infty,\gamma}(\mathcal{O},H^{s})$, then $\rho\in W^{q_{3},\infty,\gamma}(\mathcal{O},H^{s})$ and 
						%$$\|\rho\|_{q_{3},s}^{\gamma,\mathcal{O}}\lesssim\left(\|\rho\|_{q_{1},s}^{\gamma,\mathcal{O}}\right)^{t}\left(\|\rho\|_{q_{2},s}^{\gamma,\mathcal{O}}\right)^{1-t}.$$}
						%\end{enumerate}
						\item Law products:
						\begin{enumerate}[label=(\alph*)]
							\item Let $\rho_{1},\rho_{2}\in W^{q,\infty,\gamma}(\mathcal{O},H^{s}).$ Then $\rho_{1}\rho_{2}\in W^{q,\infty,\gamma}(\mathcal{O},H^{s})$ and 
							$$\| \rho_{1}\rho_{2}\|_{q,s}^{\gamma,\mathcal{O}}\lesssim\| \rho_{1}\|_{q,s_{0}}^{\gamma,\mathcal{O}}\| \rho_{2}\|_{q,s}^{\gamma,\mathcal{O}}+\| \rho_{1}\|_{q,s}^{\gamma,\mathcal{O}}\| \rho_{2}\|_{q,s_{0}}^{\gamma,\mathcal{O}}.$$
							\item Let $\rho_{1},\rho_{2}\in W^{q,\infty,\gamma}(\mathcal{O},\mathbb{C}).$ Then $\rho_{1}\rho_{2}\in W^{q,\infty,\gamma}(\mathcal{O},\mathbb{C})$ and $$\| \rho_{1}\rho_{2}\|_{q}^{\gamma,\mathcal{O}}\lesssim\| \rho_{1}\|_{q}^{\gamma,\mathcal{O}}\| \rho_{2}\|_{q}^{\gamma,\mathcal{O}}.$$
							\item Let $(\rho_{1},\rho_{2})\in W^{q,\infty,\gamma}(\mathcal{O},\mathbb{C})\times W^{q,\infty,\gamma}(\mathcal{O},H^{s}).$ Then $\rho_{1}\rho_{2}\in W^{q,\infty,\gamma}(\mathcal{O},H^{s})$ and 
							$$\| \rho_{1}\rho_{2}\|_{q,s}^{\gamma,\mathcal{O}}\lesssim\| \rho_{1}\|_{q}^{\gamma,\mathcal{O}}\| \rho_{2}\|_{q,s}^{\gamma,\mathcal{O}}.$$
						\end{enumerate}
						\item Composition law: Let $f\in C^{\infty}(\mathcal{O}\times\mathbb{R},\mathbb{R})$ and  $\rho_{1},\rho_{2}\in W^{q,\infty,\gamma}(\mathcal{O},H^{s})$  such that $$\| \rho_{1}\|_{q,s}^{\gamma,\mathcal{O}},\|\rho_{2}\|_{q,s}^{\gamma,\mathcal{O}}\leqslant C_{0}$$ for an  arbitrary  constant  $C_{0}>0$ and define the pointwise composition $$\forall (\mu,\varphi,\theta)\in \mathcal{O}\times\mathbb{T}^{d+1},\quad f(\rho)(\mu,\varphi,\theta):= f(\mu,\rho(\mu,\varphi,\theta)).$$
						Then $f(\rho_{1})-f(\rho_{2})\in W^{q,\infty,\gamma}(\mathcal{O},H^{s})$ with 
						$$\| f(\rho_{1})-f(\rho_{2})\|_{q,s}^{\gamma,\mathcal{O}}\leqslant C(s,d,q,f,C_{0})\| \rho_{1}-\rho_{2}\|_{q,s}^{\gamma,\mathcal{O}}.$$
						\item Composition law 2: Let $f\in C^{\infty}(\mathbb{R},\mathbb{R})$ with bounded derivatives. Let $\rho\in W^{q,\infty,\gamma}(\mathcal{O},\mathbb{C}).$ Then 
						$$\|f(\rho)-f(0)\|_{q}^{\gamma,\mathcal{O}}\leqslant C(q,d,f)\|\rho\|_{q}^{\gamma,\mathcal{O}}\left(1+\|\rho\|_{L^{\infty}(\mathcal{O})}^{q-1}\right).
						$$
				\end{enumerate}}
			\end{lem}
			The following technical  lemma turns out to be very useful in the study of the linearized operator.
			\begin{lem}\label{cheater lemma}
				Let $(\gamma,q,d,s_{0},s)$ satisfy  \eqref{initial parameter condition} and  $f\in W^{q,\infty,\gamma}(\mathcal{O},H^{s}).$\\
				We consider the function $g:\mathcal{O}\times\mathbb{T}_{\varphi}^{d}\times\mathbb{T}_{\theta}\times\mathbb{T}_{\eta}\rightarrow\mathbb{C}$ defined by
				$$g(\mu,\varphi,\theta,\eta)=\left\lbrace\begin{array}{ll}
					\frac{f(\mu,\varphi,\eta)-f(\mu,\varphi,\theta)}{\sin\left(\tfrac{\eta-\theta}{2}\right)} & \textnormal{if }\theta\neq \eta\\
					2\partial_{\theta}f(\mu,\varphi,\theta)& \textnormal{if }\theta=\eta.
				\end{array}\right.$$
				Then 
				$$\forall k\in\mathbb{N},\quad\|(\partial_{\theta}^kg)(\ast,\cdot,\centerdot,\eta+\centerdot)\|_{q,s}^{\gamma,\mathcal{O}}\lesssim\|\partial_{\theta}f\|_{q,s+k}^{\gamma,\mathcal{O}}\lesssim\|f\|_{q,s+k+1}^{\gamma,\mathcal{O}}.$$
			\end{lem}
			
			\begin{proof}
				Since the differentiation with respect to $\mu$ can be transported from $g$ to $f$, then it is enough to
				check the result for $q = 0$ and therefore we shall remove the dependence in $\mu$. We start with expanding
				$f$ into its Fourier series,
				$$f(\varphi, \theta)= \sum_{(l,j)\in\mathbb{Z}^{d+1}}f_{l,j}\mathbf{e}_{l,j }(\varphi, \theta).$$
				Thus, one can write
				\begin{align*}
					g(\varphi,\theta,\eta)&=\sum_{(l,j)\in\mathbb{Z}^{d+1}}f_{l,j}\tfrac{e^{\ii j\eta}-e^{\ii j\theta}}{\sin\left(\tfrac{\eta-\theta}{2}\right)}e^{il\cdot\varphi}\\
					%&=\sum_{(l,j)\in\mathbb{Z}^{d+1}}f_{l,j}e^{\ii j\tfrac{\theta+\eta}{2}}\tfrac{e^{\ii j\tfrac{\eta-\theta}{2}}-e^{-\ii j\tfrac{\eta-\theta}{2}}}{\sin\left(\tfrac{\eta-\theta}{2}\right)}e^{il\cdot\varphi}\\
					&=2\ii\sum_{(l,j)\in\mathbb{Z}^{d+1}\atop j\neq 0}f_{l,j}e^{\ii j\tfrac{\theta+\eta}{2}}\tfrac{\sin\left(j\tfrac{\eta-\theta}{2}\right)}{\sin\left(\tfrac{\eta-\theta}{2}\right)}e^{il\cdot\varphi}.
				\end{align*}
				We shall introduce the Chebychev polynomials of second kind $(U_n)_{n\in\mathbb{N}}.$ They are defined for all $n\in\mathbb{N}$ by the following relation
				$$\forall \theta\in\mathbb{R},\quad \sin(\theta)U_{n}(\cos(\theta))=\sin((n+1)\theta).$$
				Using these polynomials, we obtain a new formulation for $g$, namely
				$$g(\varphi,\theta,\eta)=2\ii\sum_{(l,j)\in\mathbb{Z}^{d+1}\atop j\neq 0}\tfrac{jf_{l,j}}{|j|}e^{\ii j\tfrac{\theta+\eta}{2}} U_{|j|-1}\left(\cos\left(\tfrac{\theta-\eta}{2}\right)\right)e^{il\cdot\varphi}.$$
				Differentiating in $\theta$ yields by Leibniz rule 
				\begin{align}\label{dtk g}
					\partial_{\theta}^{k}g(\varphi,\theta,\eta)=2\ii \sum_{(l,j)\in\mathbb{Z}^{d+1}\atop j\neq 0}\sum_{m=0}^{k}\binom{k}{m}\tfrac{j^{k+1-m}f_{l,j}\ii^{k-m}}{|j|2^{k-m}}e^{\ii j\tfrac{\theta+\eta}{2}}\partial_{\theta}^{m}\Big(U_{|j|-1}\left(\cos\left(\tfrac{\theta-\eta}{2}\right)\right)\Big).
				\end{align}
				For all $j\in\mathbb{N}^{*}$, we consider the function $f_{j}$ defined by
				$$f_{j}(\theta)=U_{j-1}(\cos(\theta))=\frac{\sin(j\theta)}{\sin(\theta)}.$$
				Notice that $f_{j}$ is even and $2\pi$-periodic. Thus, we restrict its study to the interval $[0,\pi].$ Also remark that
				$$f_{j}(\pi-\theta)=(-1)^{j}f_{j}(\theta).$$
				Hence, we restrict the study to the interval $[0,\tfrac{\pi}{2}].$ We first consider the function $f_j$ on the interval $[\tfrac{\pi}{6},\tfrac{\pi}{2}].$ There, the function $f_{j}$ writes as the quotient of two smooth functions with non vanishing denominator. Therefore, differentiating in $\theta$ leads to
				$$\forall k\in\mathbb{N},\quad\sup_{\theta\in[\frac{\pi}{6},\tfrac{\pi}{2}]}\left|\partial_{\theta}^{k}f_{j}(\theta)\right|\lesssim|j|^{k}.$$
				Now we look at the behaviour close to $0$ by looking at the function $f_{j}$ restricted to $[0,\tfrac{\pi}{4}].$ Using Taylor Formula, we can write
				\begin{align*}
					f_{j}(\theta)&=\tfrac{\sin(j\theta)}{\theta}\times\tfrac{\theta}{\sin(\theta)}\\
					&=j\int_{0}^{1}\cos(tj\theta)dt\times\tfrac{\theta}{\sin(\theta)}.
				\end{align*}
				The function $\theta\mapsto\tfrac{\theta}{\sin(\theta)}$ being smooth on $[0,\tfrac{\pi}{4}]$, then differentiating in $\theta$ leads to
				$$\forall k\in\mathbb{N},\quad\sup_{\theta\in[0,\tfrac{\pi}{4}]}\left|\partial_{\theta}^{k}f_{j}(\theta)\right|\lesssim|j|^{k+1}.$$
				Combining the previous estimates, one gets
				\begin{equation}\label{e-ptk Uj}
					\forall j\in\mathbb{N}^*,\quad\forall k\in\mathbb{N},\quad\sup_{\theta\in\mathbb{R}}\left|\partial_{\theta}^{k}\Big(U_{j-1}\left(\cos\left(\theta\right)\right)\Big)\right|\lesssim |j|^{k+1}.
				\end{equation}
				Gathering \eqref{dtk g} and \eqref{e-ptk Uj}, we deduce that
				\begin{align*}
					(\partial_{\theta}^{k}g)(\varphi,\theta,\eta+\theta)=\sum_{(l,j)\in\mathbb{Z}^{d+1}\atop j\neq 0}c_{l,j,k}(\eta)\mathbf{e}_{l,j}(\varphi,\theta),
				\end{align*}
				with
				$$\sup_{\eta\in\mathbb{T}}\left|c_{l,j,k}(\eta)\right|\lesssim|j|^{k+1}|f_{l,j}|.$$
				Therefore,
				\begin{align*}
					\sup_{\eta\in\mathbb{T}}\|(\partial_{\theta}^{k}g)(\cdot,\centerdot,\eta+\centerdot)\|_{H_{\varphi,\theta}^{s}}^{2}&=\sum_{(l,j)\in\mathbb{Z}^{d+1}\atop j\neq 0}\langle l,j\rangle^{2s}\sup_{\eta\in\mathbb{T}}|c_{l,j,k}(\eta)|^{2}\\
					&\lesssim\sum_{(l,j)\in\mathbb{Z}^{d+1}}\langle l,j\rangle^{2s}|j|^{2k+2}|f_{l,j}|^{2}\\
					&\lesssim\|\partial_{\theta}f\|_{H^{s+k}}^{2}.
				\end{align*}
				This concludes the proof of Lemma \ref{cheater lemma}.
			\end{proof}
			\subsection{Operators}\label{section-ope}
			We shall focus in this section on some useful  norms related to suitable  operators class. These notions were used before in \cite{BBMH18,BFM21,BFM21-1,BM18}.
			We consider a smooth family of bounded operators on Sobolev spaces $H^s(\T^{d+1},\mathbb{C})$, that is a smooth map $T: \mu=(\lambda,\omega)\in \mathcal{O}\mapsto T(\mu)\in  \mathcal{L}(H^s(\T^{d+1},\mathbb{C})) $ of linear continuous operators on Sobolev space  $H^s(\T^{d+1},\mathbb{C})$, with  $\mathcal{O}$ being  an open bounded set of $\mathbb{R}^{d+1}.$ % $T(\lambda,\omega): H^s(\T^{d+1})\to H^s(\T^{d+1}) $.
			 Then we find it convenient to encode   $T(\mu)$ in terms of  the infinite dimensional matrix $\left(T_{l_{0},j_{0}}^{l,j}(\mu)\right)_{\underset{(j,j_{0})\in\mathbb{Z}^{2}}{(l,l_{0})\in(\mathbb{Z}^{d })^{2}}}$ with 
			 $$T(\mu)\mathbf{e}_{l_{0},j_{0}}=\sum_{(l,j)\in\mathbb{Z}^{d+1}}T_{l_{0},j_{0}}^{l,j}(\mu)\mathbf{e}_{l,j}$$
			and
			\begin{align}\label{F-C-O}
				T_{l_{0},j_{0}}^{l,j}(\mu)&=\big\langle T(\mu)\mathbf{e}_{l_{0},j_{0}},\mathbf{e}_{l,j}\big\rangle_{L^{2}(\mathbb{T}^{d+1})}.
				%\\
				%&=\int_{\mathbb{T}^{d+1}}\big(T(\mu)\mathbf{e}_{l_{0},j_{0}}\big)(\varphi,\theta)\,\mathbf{e}_{-l,-j}(\varphi,\theta)d\varphi d\theta.
			\end{align}
			Next, we need to fix a notation that we are implicitly using along the paper. For a given   family of  multi-parameter operators  $T(\mu)$, it acts on $W^{q,\infty,\gamma}(\mathcal{O},H^{s}(\mathbb{T}^{d+1},\mathbb{C}))$ in the following sense, 
			$$
			\rho\in W^{q,\infty,\gamma}(\mathcal{O},H^{s}(\mathbb{T}^{d+1},\mathbb{C})),\quad \,\quad (T\rho)(\mu,\varphi,\theta):=T(\mu)\rho(\mu,\varphi,\theta).
			$$
			\subsubsection{Toeplitz in time operators}\label{Top-Sec-11}
			In this short section we shall introduce a suitable class of Toeplitz  operators.
			\begin{defin}
				We say that an operator $T(\mu)$ is Toeplitz in time (actually in the variable $\varphi$) if its Fourier coefficients defined by \eqref{F-C-O}, satisfy
				$$\forall\, l_o,l,j,j_0\in\mathbb{Z},\quad T_{l_{0},j_{0}}^{l,j}(\mu)=T_{0,j_{0}}^{l-l_0,j}(\mu).$$
				Or equivalently 
				$$T_{l_{0},j_{0}}^{l,j}(\mu)=T_{j_{0}}^{j}(\mu,l-l_{0}),$$
				with $T_{j_{0}}^{j}(\mu,l):=T_{0,j_{0}}^{l,j}(\mu).$
			\end{defin}
			The action of a Toeplitz operator  $T(\mu)$ on a function $\rho=\displaystyle\sum_{(l_{0},j_{0})\in\mathbb{Z}^{d+1}}\rho_{l_{0},j_{0}}\mathbf{e}_{l_{0},j_{0}}$ is then given by
			\begin{equation}\label{action of Toeplitz in time operators}
				T(\mu)\rho=\sum_{(l,l_{0})\in(\mathbb{Z}^{d})^{2}\\\atop (j,j_{0})\in\mathbb{Z}^{2}}T_{j_{0}}^{j}(\mu,l-l_{0})\rho_{l_{0},j_{0}}\mathbf{e}_{l,j}.
			\end{equation}
			In this paper, we will encounter several operators acting only on the variable $\theta$ and that can be considered as  $\varphi$-dependent operators $T(\mu,\varphi)$ taking the form
			\begin{align*}
				T(\mu,\varphi)\rho(\varphi,\theta)=\int_{\T} K(\mu,\varphi,\theta,\eta)\rho(\varphi,\eta)d\eta.
			\end{align*}
			%\textcolor{blue}{Pas uniquement sous forme intégrale. L'opérateur $\mathscr{B}$ par exemple.}\\
			One can easily check that those operators are Toeplitz and therefore they  satisfy \eqref{action of Toeplitz in time operators}.\\
			\noindent For $q\in\mathbb{N}$ and $s\in\mathbb{R},$ we can equip Toeplitz operators with the off-diagonal norm  given by,
			\begin{align}\label{Top-NormX}
				\| T\|_{\textnormal{\tiny{O-d}},q,s}^{\gamma,\mathcal{O}}=\sum_{\underset{|\alpha|\leqslant q}{\alpha\in\mathbb{N}^{d+1}}}\gamma^{|\alpha|}\sup_{\mu \in{\mathcal{O}}}\|\partial_{\mu}^{\alpha}(T)(\mu)\|_{\textnormal{\tiny{O-d}},s-|\alpha|},
			\end{align}
			where 
			$$\| T\|_{\textnormal{\tiny{O-d}},s}^{2}=\sum_{(l,m)\in\mathbb{Z}^{d+1}}\langle l,m\rangle^{2s}\sup_{j-k=m}|T_{j}^{k}(l)|^{2}.$$
			This norm will be of important use later during the KAM reduction of the remainder.
			%\textcolor{red}{On peut fusionner dans une seule proposition les lemmes 3.5, 3.6 et Prop3.6. On peut aussi enlever les demos si on a une refernce}\\
			%As in the case of functions, 
			The cut-off projectors $(P_N)_ {N\in\mathbb{N}^{*}}$ are defined as follows: 
			\begin{equation}\label{definition of projections for operators}
				\left(P_{N}T(\mu)\right)\mathbf{e}_{l_{0},j_{0}}=\sum_{\underset{|l-l_{0}|,|j-j_{0}|\leqslant N}{(l,j)\in\mathbb{Z}^{d+1}}}T_{l_{0},j_{0}}^{l,j}(\mu)\mathbf{e}_{l,j}\quad\mbox{and}\quad P_{N}^{\perp}T=T-P_{N}T.
			\end{equation}
			In the next lemma we shall gather classical results whose  proofs are very close to those in \cite{BM18} concerning pseudo-differential operators. We recall that the weighted norms on functions  that will be used below are defined in \eqref{Norm-def}.
				\begin{lem}\label{properties of Toeplitz in time operators}
					Let $(\gamma,q,d,s_{0},s)$ satisfying \eqref{initial parameter condition}.  Let $T,$ $T_{1}$ and $T_{2}$ be Toeplitz in time operators.
					\begin{enumerate}[label=(\roman*)]
						\item Projectors properties : Let $N\in\mathbb{N}^{*}.$ Let $t\in\mathbb{R}_{+}.$ Then
						$$\| P_{N}T\rho\|_{\textnormal{\tiny{O-d}},q,s+t}^{\gamma,\mathcal{O}}\leqslant N^{t}\| T\rho\|_{\textnormal{\tiny{O-d}},q,s}^{\gamma,\mathcal{O}}\quad\mbox{and}\quad\| P_{N}^{\perp}T\rho\|_{\textnormal{\tiny{O-d}},q,s}^{\gamma,\mathcal{O}}\leqslant N^{-t}\| T\rho\|_{\textnormal{\tiny{O-d}},q,s+t}^{\gamma,\mathcal{O}}.$$
						\item Interpolation inequality : Let $q<s_{1}\leqslant s_{3}\leqslant s_{2}, \,\overline{\theta}\in[0,1]$ with $s_{3}=\overline{\theta}s_{1}+(1-\overline{\theta})s_{2}.$  Then 
						$$\| T\|_{\textnormal{\tiny{O-d}},q,s_{3}}^{\gamma,\mathcal{O}}\lesssim\left(\| T\|_{\textnormal{\tiny{O-d}},q,s_{1}}^{\gamma,\mathcal{O}}\right)^{\overline{\theta}}\left(\| T\|_{\textnormal{\tiny{O-d}},q,s_{2}}^{\gamma,\mathcal{O}}\right)^{1-\overline{\theta}}.$$
						\item Composition law :
						$$\| T_{1}T_{2}\|_{\textnormal{\tiny{O-d}},q,s}^{\gamma,\mathcal{O}}\lesssim\| T_{1}\|_{\textnormal{\tiny{O-d}},q,s}^{\gamma,\mathcal{O}}\| T_{2}\|_{\textnormal{\tiny{O-d}},q,s_{0}}^{\gamma,\mathcal{O}}+\| T_{1}\|_{\textnormal{\tiny{O-d}},q,s_{0}}^{\gamma,\mathcal{O}}\| T_{2}\|_{\textnormal{\tiny{O-d}},q,s}^{\gamma,\mathcal{O}}.$$ 
						\item Link between operators and off-diagonal norms :
						$$\| T\rho\|_{q,s}^{\gamma,\mathcal{O}}\lesssim\| T\|_{\textnormal{\tiny{O-d}},q,s_{0}}^{\gamma,\mathcal{O}}\|\rho\|_{q,s}^{\gamma,\mathcal{O}}+\| T\|_{\textnormal{\tiny{O-d}},q,s}^{\gamma,\mathcal{O}}\|\rho\|_{q,s_{0}}^{\gamma,\mathcal{O}}.$$
						In particular
						$$\| T \rho\|_{q,s}^{\gamma,\mathcal{O}}\lesssim\| T\|_{\textnormal{\tiny{O-d}},q,s}^{\gamma,\mathcal{O}}\|\rho\|_{q,s}^{\gamma,\mathcal{O}}.$$
					\end{enumerate}
			\end{lem}
			\subsubsection{Reversible and reversibility preserving operators}
			In this section we intend to  collect some definitions and properties related to different reversibility notions for operators and give practical characterizations. We shall also come back to Toeplitz operators defined before in  Section \ref{Top-Sec-11} and discuss two important examples frequently encountered during  this paper and given by multiplications and integral operators.\\
			First, we  give  the following definitions following \cite[Def. 2.2]{BBM14}.
			\begin{defin}\label{Def-Rev}
					Introduce the following involution
					\begin{equation}\label{definition involution mathcalS2}
						(\mathscr{S}_{2}\rho)(\varphi,\theta)=\rho(-\varphi,-\theta).
					\end{equation}
					We say that an operator $T(\mu)$ is 
					\begin{enumerate}[label=\textbullet]
						\item real if for all $\rho\in L^{2}(\mathbb{T}^{d+1},\mathbb{C}),$ we have 
						$$\overline{\rho}=\rho\quad\Longrightarrow\quad\overline{T\rho}=T\rho.$$
						\item reversible if
						$$T(\mu)\circ\mathscr{S}_{2}=-\mathscr{S}_{2}\circ T(\mu).$$
						\item reversibility preserving if
						$$T(\mu)\circ\mathscr{S}_{2}=\mathscr{S}_{2}\circ T(\mu).$$
					\end{enumerate}
				\end{defin}
				We now detail the following   characterizations needed at several places in this paper and the proofs are  quite easy and follow from Fourier expansion. One can find a similar result in \cite[Lem. 2.6]{BBM14}.
								\begin{prop}\label{characterization of real operator by its Fourier coefficients}
					Let $T$ be an operator. Then $T$ is 
					\begin{enumerate}[label=\textbullet]
						\item real if and only if
						$$\forall(l,l_{0},j,j_{0})\in(\mathbb{Z}^{d})^{2}\times\mathbb{Z}^{2},\quad T_{-l_{0},-j_{0}}^{-l,-j}=\overline{T_{l_{0},j_{0}}^{l,j}}.$$
						\item reversible if and only if
						$$\forall(l,l_{0},j,j_{0})\in(\mathbb{Z}^{d})^{2}\times\mathbb{Z}^{2},\quad T_{-l_{0},-j_{0}}^{-l,-j}=-T_{l_{0},j_{0}}^{l,j}.$$
						\item reversibility-preserving if and only if
						$$\forall(l,l_{0},j,j_{0})\in(\mathbb{Z}^{d})^{2}\times\mathbb{Z}^{2},\quad T_{-l_{0},-j_{0}}^{-l,-j}=T_{l_{0},j_{0}}^{l,j}.$$
					\end{enumerate}
				\end{prop}
			In what follows,  we shall focus on  two particular cases of operators which will be of constant use throughout this paper. Namely, multiplication and integral operators.
			\begin{defin}\label{Defin-Rever-1}
				Let $T$ be an operator as in Section $\ref{section-ope}.$
 We say that 
				\begin{enumerate}[label=\textbullet]
					\item $T$ is a multiplication operator if there exists a function $M:(\mu,\varphi,\theta)\mapsto M(\mu,\varphi,\theta)$ such that
					$$(T\rho)(\mu,\varphi,\theta)=M(\mu,\varphi,\theta)\rho(\mu,\varphi,\theta).$$
					\item $T$ is an integral operator if there exists a function (called the kernel) $K:(\mu,\varphi,\theta,\eta)\mapsto K(\mu,\varphi,\theta,\eta)$ such that
					$$(T\rho)(\mu,\varphi,\theta)=\int_{\mathbb{T}}\rho(\mu,\varphi,\eta)K(\mu,\varphi,\theta,\eta)d\eta.$$
				
				\end{enumerate}
			\end{defin}
			We intend to prove the following lemma.
			\begin{lem}\label{lemma symmetry and reversibility}
				{Let $(\gamma,q,d,s_{0},s)$ satisfy \eqref{initial parameter condition}, then the following assertions hold true.
					\begin{enumerate}[label=(\roman*)]
						\item Let $T$ be a multiplication operator by a real-valued function $M$, then 						\begin{enumerate}[label=\textbullet]
							\item If $M(\mu,-\varphi,-\theta)=M(\mu,\varphi,\theta)$, then $T$ is  real and reversibility preserving Toeplitz in time and space operator.
							\item If $M(\mu,-\varphi,-\theta)=-M(\mu,\varphi,\theta)$, then $T$ is  real and reversible Toeplitz in time and space operator.
						\end{enumerate}
						Moreover,
						$$\| T\|_{\textnormal{\tiny{O-d}},q,s}^{\gamma,\mathcal{O}}\lesssim\| M\|_{q,s+s_{0}}^{\gamma,\mathcal{O}}.$$
						\item Let $T$ be an integral operator with a real-valued kernel $K$.
						\begin{enumerate}[label=\textbullet]
							\item If $K(\mu,-\varphi,-\theta,-\eta)=K(\mu,\varphi,\theta,\eta)$, then $T$ is  a real and reversibility preserving Toeplitz in time operator.
							\item If $K(\mu,-\varphi,-\theta,-\eta)=-K(\mu,\varphi,\theta,\eta)$, then $T$ is  a real and reversible Toeplitz in time operator.
						\end{enumerate}
						In addition,
						$$\| T\|_{\textnormal{\tiny{O-d}},q,s}^{\gamma,\mathcal{O}}\lesssim \int_{\T}\|K(\ast,\cdot,\centerdot,\eta+\centerdot)\|_{q,s+s_{0}}^{\gamma,\mathcal{O}} d\eta$$
						and
						$$\| T\rho\|_{q,s}^{\gamma,\mathcal{O}}\lesssim \|\rho\|_{q,s_0}^{\gamma,\mathcal{O}} \int_{\T}\|K(\ast,\cdot,\centerdot,\eta+\centerdot)\|_{q,s}^{\gamma,\mathcal{O}} d\eta+\|\rho\|_{q,s}^{\gamma,\mathcal{O}} \int_{\T}\|K(\ast,\cdot,\centerdot,\eta+\centerdot)\|_{q,s_0}^{\gamma,\mathcal{O}} d\eta,
						$$
						where the notation $\ast,\cdot,\centerdot$ denote $\mu,\varphi,\theta$, respectively.
				\end{enumerate}}
			\end{lem}
			\begin{proof} We point out that the proofs will be implemented for the particular case $q=0$ and the general case can  be  done similarly by  differentiating with respect to $\mu$ and using Leibniz rule. 
				\\ 
				\textbf{(i)}  Since $M$ is a real-valued function, then we get by the definition
				\begin{align*}\overline{T_{-l_{0},-j_{0}}^{-l,-j}}&=\displaystyle\int_{\mathbb{T}^{d+1}}M(\varphi,\theta)\overline{\mathbf{e}_{-l_{0},-j_{0}}(\varphi,\theta)}\overline{\mathbf{e}_{l,j}(\varphi,\theta)}d\varphi d\theta\\
					&=\displaystyle\int_{\mathbb{T}^{d+1}}M(\varphi,\theta)\mathbf{e}_{l_{0},j_{0}}(\varphi,\theta)\mathbf{e}_{-l,-j}(\varphi,\theta)d\varphi d\theta=T_{l_{0},j_{0}}^{l,j}.
				\end{align*}
				This shows in view of  Proposition \ref{characterization of real operator by its Fourier coefficients} that the operator   $T$ is a real. It remains to check the reversibility preserving property.
				We write from the definition
				\begin{align*}
					T(\mathscr{S}_{2}\rho)(\varphi,\theta)&=M(\varphi,\theta)\rho(-\varphi,-\theta)\\
					&= M(-\varphi,-\theta)\rho(-\varphi,-\theta)\\
					&=\mathscr{S}_{2}\left(T\rho\right)(\varphi,\theta).
				\end{align*}
				This gives the desired result. As  to the reversible Toepliz structure, it can be checked in a similar way.  To achieve the proof of the first point it remains to establish the suitable estimate.
				Using a duality argument $H^{s+s_{0}}-H^{-s-s_{0}}$, we may write,
				$$|T_{j}^{j'}(l)|=\left|\int_{\mathbb{T}^{d+1}}M(\varphi,\theta)\mathbf{e}_{l,j-j'}(\varphi,\theta)d\varphi d\theta\right|\lesssim\langle l,j-j'\rangle^{-s-s_{0}}\| M\|_{H^{s+s_{0}}}.$$
				It follows that
				\begin{align*}\| T\|_{\textnormal{\tiny{O-d}},s}^{2}=&\sum_{(l,m)\in\mathbb{Z}^{d+1}}\langle l,m\rangle^{2s}\sup_{j-j^\prime=m}|T_{j}^{j^\prime}(l)|^{2}\\
					\lesssim&\| M\|_{H^{s+s_{0}}}^2\sum_{(l,m)\in\mathbb{Z}^{d+1 }}\langle l,m\rangle^{2s}\langle l,m\rangle^{-2s-2s_{0}}\\
					\lesssim&\| M\|_{H^{s+s_{0}}}^2.
				\end{align*}
				Therefore we find
				$$\| T\|_{\textnormal{\tiny{O-d}},s}\lesssim\| M\|_{H^{s+s_{0}}}.$$
				\textbf{(ii)}  By assumption,  $K$ is real  and thus 
				\begin{align*}
					\overline{T_{-l_{0},-j_{0}}^{-l,-j}}&=\displaystyle\int_{\mathbb{T}^{d+2}}K(\varphi,\theta,\eta)\overline{\mathbf{e}_{-l_{0},-j_{0}}(\varphi,\eta)}\overline{\mathbf{e}_{l,j}(\varphi,\theta)}d\varphi d\theta\\
					&=\displaystyle\int_{\mathbb{T}^{d+2}}K(\varphi,\theta,\eta)\mathbf{e}_{l_{0},j_{0}}(\varphi,\eta)\mathbf{e}_{-l,-j}(\varphi,\theta)d\varphi d\theta d\eta=T_{l_{0},j_{0}}^{l,j}.
				\end{align*}
				{This implies, according to Proposition \ref{characterization of real operator by its Fourier coefficients},  that $T$ is a real operator. Now we shall check the reversibility preserving. The reversibility can be checked in a similar way. By the change of variables $\eta\mapsto-\eta,$ we may write,
					\begin{align*}
						T(\mathscr{S}_{2}\rho)(\varphi,\theta)&=\int_{\mathbb{T}}K(\varphi,\theta,\eta)\rho(-\varphi,-\eta)d\eta\\
						&=\int_{\mathbb{T}}K(-\varphi,-\theta,-\eta)\rho(-\varphi,-\eta)d\eta\\
						&=\int_{\mathbb{T}}K(-\varphi,-\theta,\eta)\rho(-\varphi,\eta)d\eta=\mathscr{S}_{2}\left(T\rho\right)(\varphi,\theta).
				\end{align*}}
				From Fubini's theorem  and the  duality $H^{s+s_{0}}_{\varphi,\theta}-H_{\varphi,\theta}^{-s-s_{0}}$, we infer,
				$$\begin{array}{rcl}
					|T_{j}^{j'}(l)| & = & \displaystyle\left|\int_{\mathbb{T}^{d+2}}K(\varphi,\theta,\eta)e^{\ii(l\cdot\varphi+j\theta-j'\eta)}d\varphi d\theta d\eta\right|\\
					& = & \displaystyle\left|\int_{\mathbb{T}^{d+1}} e^{\ii(l\cdot\varphi+(j-j')\theta)}\left(\int_{\T}K(\varphi,\theta,\eta+\theta)e^{-\ii j'\eta} d\eta\right) d\varphi d\theta \right|\\
					& \lesssim & \displaystyle\langle l,j-j'\rangle^{-s-s_{0}} \int_{\T}\|K(\ast,\cdot,\centerdot,\eta+\centerdot)\|_{H^{s+s_0}_{\varphi,\theta}} d\eta.
				\end{array}$$
				Hence, we deduce that
				$$\| T\|_{\textnormal{\tiny{O-d}},s}\lesssim\int_{\T}\|K(\ast,\cdot,\centerdot,\eta+\centerdot)\|_{H^{s+s_0}_{\varphi,\theta}} d\eta
				.$$
				The last estimate   in  Lemma \ref{lemma symmetry and reversibility} can be  obtained from the expression
				$$(T\rho)(\varphi,\theta)=\int_{\mathbb{T}}\rho(\varphi,\theta+\eta)K(\varphi,\theta,\theta+\eta)d\eta,
				$$
				combined with the  law products and the translation invariance in  Lemma \ref{Lem-lawprod}-(i)-(iv).\\
				This concludes the proof of Lemma \ref{lemma symmetry and reversibility}.
			\end{proof}
			\section{Hamiltonian toolkit and  approximate inverse}\label{sec Ham-tool}
			In this section, we shall  reformulate the problem into the form of searching for zeros of a functional $\mathcal{F}.$ We first rescale the equation by introducing a small parameter $\varepsilon$. This allows us to see the Hamiltonian equation \eqref{Hamiltonian formulation of the equation} as a perturbation of the equilibirum one \eqref{hamiltonian equation at equilibrium}. The latter being integrable and admitting quasi-periodic solutions in view of Lemma \ref{lemma linearized operator at equilibrium}-2 and Lemma \ref{lemma sol Eq}, we can hope using KAM technics to find quasi-periodic solutions to the first one. This approach has been intensively used before in \cite{BBMH18,BBM16,BFM21,BFM21-1,BM18}. According to Proposition \ref{proposition Hamiltonian formulation of the equation}, it seems more convenient to work  with the phase space
			$$H_{0}^{s}=H_{0}^{s}(\mathbb{T},\mathbb{R})=\left\lbrace r(\theta)=\sum_{j\in\mathbb{Z}^{*}}r_{j}e^{ij\theta}\quad\textnormal{s.t.}\quad r_{-j}=\overline{r_{j}}\quad\mbox{ and }\quad\| r\|_{s}^{2}=\sum_{j\in\mathbb{Z}^{*}}|r_{j}|^{2}|j|^{2s}<\infty\right\rbrace.$$
			Therefore, we select finitely-many tangential sites $\mathbb{S}$ and decompose the phase space into tangential and normal subspaces described by the selection of Fourier modes belonging to $\mathbb{S}$ or not. On the tangential part,  containing the main part of the  quasi-periodic solutions, we introduce action-angle variables allowing to reformulate the problem in terms of  embedded invariant tori.  We shall also be concerned with some  regularity aspects for the perturbed Hamiltonian vector field appearing in $\mathcal{F}$ and needed during the Nash-Moser scheme. Finally, we construct an approximate right inverse for the linearized operator associated to $\mathcal{F}.$\\
The symplectic structure on $L^2_{0}(\T)$ (corresponding to the subspace of $L^2(\T)$ of real functions with zero average) induced by \eqref{Hamiltonian formulation of the equation} is given by the symplectic $2$-form
			\begin{equation}\label{def symp-form}
				\mathcal{W}(r,h)=\int_{\mathbb{T}}\partial_{\theta}^{-1}r(\theta)h(\theta)d\theta\quad\mbox{ with }\quad\partial_{\theta}^{-1}r(\theta)=\sum_{j\in\mathbb{Z}^{*}}\frac{r_{j}}{\ii j}e^{\ii j\theta}.
			\end{equation}
			%The corresponding Hamiltonian vector field is $X_{H}(r)=\partial_{\theta}\nabla H(r)$ (where $\nabla$ is the $L^{2}$-gradient).
			Then for a given function  $H$,  its  symplectic gradient $X_{H}$  is defined through the identity
			$$dH(r)[\cdot]=\mathcal{W}(X_{H}(r),\cdot).$$
			Using the  Fourier expansion
			$$r(\theta)=\sum_{j\in\mathbb{Z}^{*}}r_{j}e^{\ii j\theta}\quad\textnormal{with}\quad r_{-j}=\overline{r_{j}},$$
			we easily find that the symplectic form $\mathcal{W}$ writes
			$$\mathcal{W}(r,h)=\sum_{j\in\mathbb{Z}^{*}}\frac{1}{\ii j}r_{j}h_{-j}=\sum_{j\in\mathbb{Z}^{*}}\frac{1}{\ii j}r_{j}\overline{h_{j}},$$
			that is
			\begin{equation}\label{sympl ref}
				\mathcal{W}=\sum_{j\in\mathbb{Z}^{*}}\frac{1}{\ii j}dr_{j}\wedge dr_{-j}=2\sum_{j\in\mathbb{N}^{*}}\frac{1}{\ii j}dr_{j}\wedge dr_{-j}.
			\end{equation}
			Next,   with the result of Lemma \ref{lemma linearized operator at equilibrium} we can easily check that the equation \eqref{Hamiltonian formulation of the equation}  can be written in the form 
			$$\partial_{t}r=\partial_{\theta}\mathrm{L}(\lambda)(r)+X_{P}(r),$$
			where $X_{P}$ is the  Hamiltonian vector field defined by
			\begin{equation}\label{def XP}
				X_{P}(r)=I_{1}(\lambda)K_{1}(\lambda)\partial_{\theta}r-\partial_{\theta}\mathcal{K}_{\lambda}\ast r-F_{\lambda}[r].
			\end{equation}
			Remind that $F_{\lambda}[r]$ is introduced in  \eqref{definition of Flambda} and the convolution kernel is stated in \eqref{definition of mathcalKlambda}. To measure the smallness condition it  seems to be more convenient  to  introduce a small parameter $\varepsilon$ and rescale the Hamiltonian as done  for instance in the papers \cite{BBMH18,BM18}. To do that  we rescale the solution as follows  $r\mapsto\varepsilon r$ with $r$ bounded. Therefore the Hamiltonian equation  takes the form
			\begin{equation}\label{perturbed hamiltonian}
				\partial_{t}r=\partial_{\theta}\mathrm{L}(\lambda)(r)+\varepsilon X_{P_{\varepsilon}}(r),
			\end{equation}
			where $X_{P_{\varepsilon}}$ is the rescaled Hamiltonian vector field defined by
			$X_{P_{\varepsilon}}(r):=\varepsilon^{-2}X_{P}(\varepsilon r).$
			Notice  that \eqref{perturbed hamiltonian} can be recast in the Hamiltonian form
			\begin{equation}\label{MM-E}
				\partial_{t}r=\partial_{\theta}\nabla \mathcal{H}_{\varepsilon}(r),
			\end{equation}
			where the rescaled Hamiltonian  $\mathcal{H}_{\varepsilon}(r)$ is given by 
			\begin{align}\label{HEE}
				\nonumber \mathcal{H}_{\varepsilon}(r)&=\varepsilon^{-2}H(\varepsilon r)\\
				&:=H_{\mathrm{L}}(r)+\varepsilon P_{\varepsilon}(r),
			\end{align}
			with  $H_{\mathrm{L}}$ being  the quadratic Hamiltonian defined in Lemma \ref{lemma linearized operator at equilibrium} and 
			$\varepsilon P_{\varepsilon}(r)$ is composed with terms  of  higher order more than  cubic. 
			\subsection{Action-angle reformulation}\label{subsec act-angl}
			Let us consider finitely many Fourier-frenquencies, called tangential sites, gathered in the tangantial set $\mathbb{S}$ defined by 
			$$\quad{\mathbb S} = \{ j_1, \ldots, j_d \}\subset\mathbb{N}^* \quad\textnormal{with}\quad 
			1 \leqslant j_1 < j_2 < \ldots < j_d.$$
			We now define the symmetrized tangential sets $\overline{\mathbb{S}}$ and $\mathbb{S}_{0}$ by
			\begin{align}\label{tangent-set}
				\overline{{\mathbb S}}={\mathbb S}\cup (-{\mathbb S})= \{ \pm j,\,\,j\in {\mathbb S}\}\quad\textnormal{and}\quad\mathbb{S}_{0}=\overline{\mathbb{S}}\cup\{0\}.
			\end{align}
			% $|{\mathbb S}|$ is required to be even because the solutions $r$ of the nonlinear  equation \eqref{eq} have to be real valued. 
			Recall from \eqref{def linear frequency vector} that we denote
			the unperturbed tangential frequency vector by 
			\begin{equation}\label{definition linear frequency vector}
				{\omega}_{\textnormal{\tiny{Eq}}}(\lambda)=(\Omega_{j}(\lambda))_{j\in\mathbb{S}},
			\end{equation}
			where $\Omega_j (\lambda) $ are given by \eqref{def eigenvalues at equilibrium}. For $s\in\mathbb{R}$, we decompose the phase space of $H^s_0 (\mathbb{T})$ as the direct sum 
			\begin{align}\label{decoacca}
				H^s_0 (\mathbb{T}) =& 
				{H}_{\overline{\mathbb{S}}}
				\overset{\perp}{\oplus} {H}^{s}_{\bot},\\
				\nonumber {H}_{\overline{\mathbb{S}}} = \Big\{v=\sum_{ j\in \overline{ \mathbb{S}}} r_j e^{\ii j\theta},\;   \overline{r_j}=r_{-j} \Big\}  &, \quad
				{H}^{s}_{\bot}= \Big\{ z = 
				\sum_{j\in \mathbb{Z}\setminus\mathbb{S}_0} z_j e_j\in H^s,\, \overline{z_j}=z_{-j}  \Big\} \, ,
			\end{align}
			where $e_{j}(\theta)=e^{\ii j\theta}.$ We denote by $\Pi_{\overline{\mathbb S}}, \Pi^\bot_{\mathbb S_0}$ the corresponding orthogonal projectors defined by
			\begin{equation}\label{def proj}
				r=v+z,\quad v:=\Pi_{\overline{\mathbb S}}r:=\sum_{ j\in \overline{\mathbb{S}}} r_j e_j,\quad z:=\Pi^\bot_{\mathbb S_0}r:=\sum_{ j\in \mathbb{Z}\setminus\mathbb{S}_0} r_j e_j\,,
			\end{equation}
			where $v$ and $z$ are  called the tangential and normal variables, respectively.
			Fix some small amplitudes $(\mathtt{a}_{j})_{j\in\mathbb{S}}\in(\mathbb{R}_{+}^{*})^{d}$ and set $\mathtt{a}_{-j}=\mathtt{a}_{j}$.  We shall now introduce the action-angle variables on the tangential set $H_{\overline{\mathbb{S}}}$ by making the following symplectic polar change of coordinates
%			\begin{equation}\label{ham-syst-k}
%				\left\lbrace\begin{array}{ll}
%					\displaystyle r_j  =
%					\sqrt{\mathtt{a}_{j}+|j|I_j}\,  e^{i \vartheta_j} \quad&\textnormal{if} \quad  j\in \mathbb{S}_0,
%					\\
%					r_j = z_j \quad&\textnormal{if} \quad  j\in \mathbb{Z}^{*}\setminus\mathbb{S}_0
%				\end{array}\right.
%			\end{equation}
\begin{equation}\label{ham-syst-k}
	\forall\, j\in\overline{\mathbb{S}},\quad r_j  =\sqrt{\mathtt{a}_{j}^{2}+\tfrac{|j|}{2}I_j}\, e^{\ii \vartheta_j},
\end{equation}
			where
			\begin{equation}\label{sym I-vartheta}
				\forall\, j\in\overline{\mathbb{S}},\quad I_{-j}=I_j\in\mathbb{R}\quad \textnormal{and}\quad \vartheta_{-j}=-\vartheta_j\in \mathbb{R}.
			\end{equation}
			Thus, any function
			of the phase space $H_0^s$ decomposes as
			\begin{equation}\label{definition of A action-angle-normal}
				r= A( \vartheta,I,z) :=  
				v (\vartheta,I)+ z  \quad
				\textnormal{where} \quad  
				v (\vartheta, I) := \sum_{j \in \overline{\mathbb{S}}}    
				\sqrt{\mathtt{a}_{j}^{2}+\tfrac{|j|}{2}I_j}\,  e^{\ii \vartheta_j}e_j \, . 
			\end{equation}
			In these coordinates the solutions \eqref{rev sol eq} of the linear system \eqref{hamiltonian equation at equilibrium} simply read as $v (-{\omega}_{\textnormal{Eq}} (\lambda)t,I)$ where ${\omega}_{\textnormal{Eq}}$ is defined in \eqref{definition linear frequency vector} and $I\in \mathbb{R}^d$ such that the quantity under the square root is positive. The involution $ {\mathscr S}$ defined in \eqref{definition of the involution mathcal S} now reads in the new variables
			\begin{equation}\label{rev_aa}
				{\mathfrak S} : (   \vartheta,I, z)\mapsto ( -\vartheta,I, {\mathscr S} z )
			\end{equation}
			and the symplectic  $2$-form in \eqref{sympl ref} becomes after straightforward computations using \eqref{ham-syst-k} and \eqref{sym I-vartheta}
			\begin{equation}\label{sympl_form}
				{\mathcal W} =  
				\sum_{j \in \mathbb{S}} d\vartheta_j  \wedge  d I_j   +\sum_{j \in \mathbb{Z} \setminus \mathbb{S}_0}\frac{1}{\ii j}dr_j\wedge dr_{-j}=
				\Big(\sum_{j \in \mathbb{S}}d\vartheta_j  \wedge  d I_j   \Big)  \oplus {\mathcal W}_{| H^s_\bot} , 
			\end{equation}
			where ${\mathcal W}_{|{ H}_{\bot}^s}$ denotes the restriction of $\mathcal{W}$ to ${ H}_{\bot}^s$. %  and $\Lambda$ is the contact 1-form on $\mathbb{T}^\nu\times\mathbb{R}^\nu\times{\mathbb H}_{\mathbb{S}^\bot}$ defined by $\Lambda_{( \vartheta, I, z)}:\mathbb{R}^\nu\times\mathbb{R}^\nu\times{\mathbb H}_{\mathbb{S}^\bot}\to \mathbb{R}$
			Note that ${\mathcal W}  $ is an exact $ 2 $-form as 
			$$
			\mathcal{W} = d \Lambda ,
			$$
			where $ \Lambda $ is the Liouville $ 1$-form defined by
			\begin{equation}\label{Lambda 1 form}
				\Lambda_{( \vartheta, I, z)}[
				\widehat \vartheta, \widehat I, \widehat z] =-\sum_{j \in \mathbb{S}} I_j   \widehat\vartheta_j 
				+ \tfrac12 \langle  \partial_\theta^{-1}  z, \widehat z \rangle_{L^2 (\mathbb{T})} \, .    
			\end{equation}
%			The Poisson structure $\mathcal{J}$ corresponding to $\mathcal{W}$, defined by the identity $\{F,G\}=\mathcal{W}(X_F,X_G)=\langle\nabla F, \mathcal{J}\nabla G\rangle$ where $\langle \cdot,\cdot\rangle$ is  the inner product, defined by
%			$$
%			\langle ( \vartheta_1,  I_1,  z_1),(\vartheta_2,I_2, z_2)\rangle=\vartheta_1\cdot\vartheta_2+  I_1\cdot I_2+(z_1,  z_2 )_{L^2 (\mathbb{T})},
%			$$
%			is the unbounded operator 
%			$$
%			\mathcal{J}: (\vartheta,  I,  z)\mapsto (  I,- \vartheta, \partial_\theta  z)
%			$$ 
%			
%			\smallskip
%			
%			
			The next goal is to   study the Hamiltonian system  generated by the Hamiltonian $ {\mathcal H}_\varepsilon   $ in \eqref{HEE}, 
			in  the action-angle and normal  coordinates 
			$ 
			(\vartheta, I, z) \in  \mathbb{T}^\nu \times \mathbb{R}^\nu \times {H}^{s}_{\bot} \, . 
			$ 
			We consider the Hamiltonian $ H_{\varepsilon} $ defined by 
			\begin{equation}\label{Hepsilon}
				H_{\varepsilon} =\mathcal{H}_{\varepsilon} \circ  A,
			\end{equation}
			where  
			$ A $ is the map described before  in \eqref{definition of A action-angle-normal}.  Since $ \mathrm{L}(\lambda) $  in \eqref{defLHL} is a Fourier multiplier keeping invariant  the subspaces $H_{\mathbb{S}_0}$ and ${H}^{s}_{\bot}$, then the quadratic Hamiltonian $ H_{\rm L}$ in \eqref{defLHL}
			in the variables $ (\vartheta, I, z) $ reads, up to an additive constant which can be removed since it does not change the dynamics in view of \eqref{Hamiltonian formulation of the equation},   
			\begin{align}\label{QHAM}
				H_\mathrm{L} \circ A =  -\sum_{j\in\mathbb{S}} \, \Omega_j(\lambda)I_j+ \tfrac12  \langle  \mathrm{L}(\lambda)\, z, z \rangle_{L^2(\mathbb{T})} =   -{\omega}_{\textnormal{Eq}}(\lambda)\cdot I
				+ \tfrac12  \langle  \mathrm{L}(\lambda)\, z, z \rangle_{L^2(\mathbb{T})}  ,
			\end{align}
			where $ {\omega}_{\textnormal{Eq}} \in \mathbb{R}^d $ is the unperturbed 
			tangential frequency vector defined by \eqref{def eigenvalues at equilibrium}.
			According to  \eqref{HEE} and \eqref{QHAM}, one deduces that  
			the Hamiltonian $H_{\varepsilon} $ in \eqref{Hepsilon} has the form
			\begin{equation}\label{cNP}
				H_{\varepsilon} = 
					{\mathcal N} + \varepsilon \mathcal{ P}_{\varepsilon}  \quad {\rm with} \quad{\mathcal N} :=   -{\omega}_{\textnormal{Eq}}(\lambda)\cdot I  + \tfrac12  \langle  \mathrm{L}(\lambda)\, z, z \rangle_{L^2(\mathbb{T})}  
					\quad \textnormal{and}
					\quad \mathcal{ P}_{\varepsilon} :=   P_\varepsilon \circ A .  
			\end{equation}
			We look for an embedded invariant torus
			\begin{equation}\label{rev-torus}
				i:\begin{array}[t]{rcl} \mathbb{T}^d &\rightarrow&
				\mathbb{R}^d \times \mathbb{R}^d \times {H}^{s}_{\bot} \\
				\varphi &\mapsto& i(\varphi):= (  \vartheta(\varphi), I(\varphi), z(\varphi))
				\end{array}  
			\end{equation}
			of the Hamiltonian vector field 
			\begin{equation}\label{hamiltonian vector filed associated with Hepsilon}
			X_{H_{\varepsilon}}:= 
			(\partial_I H_{\varepsilon} , -\partial_\vartheta H_{\varepsilon} , \Pi_{\mathbb{S}_0}^\bot
			\partial_\theta \nabla_{z} H_{\varepsilon} ) 
			\end{equation}
			filled by quasi-periodic solutions with Diophantine frequency 
			vector $\omega$. 
			Remark that  for the value   $\varepsilon=0$, the Hamiltonian system reduces to the linear equation 
			$$\omega\cdot\partial_\varphi i (\varphi) = X_{H_0} ( i (\varphi))$$   which admits the trivial solution given by  the flat  torus  $i_{\textnormal{\tiny{flat}}}(\varphi)=(\varphi,0,0)$ provided that $\omega=-{\omega}_{\textnormal{Eq}}(\lambda).$
			In what follows we shall consider the  modified Hamiltonian equation indexed with a   parameter $\alpha\in\mathbb{R}^d$,
			\begin{equation}\label{H alpha}
				\begin{aligned}
					H_\varepsilon^\alpha := {\mathcal N}_\alpha +\varepsilon  {\mathcal P}_{\varepsilon} \quad  \textnormal{where}\quad{\mathcal N}_\alpha :=  
					\alpha \cdot I 
					+ \frac12 \langle \mathrm{L}(\lambda)\, z, z\rangle_{L^2(\mathbb{T})}.
				\end{aligned}
			\end{equation}
			For the value $\alpha=-{\omega}_{\textnormal{Eq}}(\lambda)$ we have $H_\varepsilon^\alpha= H_\varepsilon$.
			The parameter $\alpha$ will play the role of a Lagrangian multiplier in order to satisfy a compatibility condition during  the approximate inverse process. Notice that the initial problem reduces after this multiple transformations to find 
			 zeros of the nonlinear operator
			\begin{equation}\label{main function}
				\begin{array}{l}
					\mathcal{F}(i,\alpha,\mu,\varepsilon):=\omega\cdot\partial_{\varphi}i(\varphi)-X_{H_{\varepsilon}^{\alpha}}(i(\varphi))\\
					\mbox{\hspace{2cm}}=\left(\begin{array}{c}
						\omega\cdot\partial_{\varphi}\vartheta(\varphi)-\alpha-\varepsilon\partial_{I}\mathcal{P}_{\varepsilon}(i(\varphi))\\
						\omega\cdot\partial_{\varphi}I(\varphi)+\varepsilon\partial_{\theta}\mathcal{P}_{\varepsilon}(i(\varphi))\\
						\omega\cdot\partial_{\varphi}z(\varphi)-\partial_{\theta}\big[\mathrm{L}(\lambda)z(\varphi)+\varepsilon\nabla_{z}\mathcal{P}_{\varepsilon}\big(i(\varphi)\big)\big]
					\end{array}\right),\,\mu=(\lambda,\omega),
				\end{array}
			\end{equation}
			where $\mathcal{P}_{\varepsilon}$ is defined in \eqref{HEE}. We point out that we can easily check that 
			the Hamiltonian $H_{\varepsilon}^{\alpha}$ is reversible in the sense of the Definition \ref{Def-Rev},  that is, 
			\begin{equation}\label{rev_Halpha}
				H_{\varepsilon}^{\alpha}\circ\mathfrak{S}=H_{\varepsilon}^{\alpha},
			\end{equation}
			where the involution $\mathfrak{S}$ is defined in \eqref{rev_aa}. Thus, we shall look for reversible solutions of 
			$$\mathcal{F}(i,\alpha,\mu,\varepsilon)=0,$$
			that is, solutions  satisfying
			$$\mathfrak{S}i(\varphi)=i(-\varphi),$$ 
			or equivalently,
			\begin{equation}\label{reversibility condition in the variables theta I and z}
				\begin{array}{ccc}
					\vartheta(-\varphi)=-\vartheta(\varphi), & I(-\varphi)=I(\varphi), & z(-\varphi)=(\mathscr{S}z)(\varphi).
				\end{array}
			\end{equation}
			We define the periodic component $\mathfrak{I}$ of the torus $i$ by
			$$\mathfrak{I}(\varphi):=i(\varphi)-(\varphi,0,0)=(\Theta(\varphi),I(\varphi),z(\varphi))\quad \mbox{ with }\quad \Theta(\varphi)=\vartheta(\varphi)-\varphi.$$
			We define the weighted Sobolev norm of $\mathfrak{I}$ as 
			$$\|\mathfrak{I}\|_{q,s}^{\gamma,\mathcal{O}}:=\|\Theta\|_{q,s}^{\gamma,\mathcal{O}}+\| I\|_{q,s}^{\gamma,\mathcal{O}}+\| z\|_{q,s}^{\gamma,\mathcal{O}}.$$
			\subsection{Hamiltonian regularity}
			This section is devoted to some regularity aspects of the Hamiltonian vector field introduced in \eqref{def XP}, together with the rescaled one associated to the Hamiltonian described in \eqref{cNP}. The first main result reads as follows.	\begin{lem}\label{lemma estimates vector field XP}
				{Let $(\gamma,q,s_{0},s)$ satisfying \eqref{initial parameter condition}. There exists $\varepsilon_{0}\in(0,1)$ such that if
					$$\| r\|_{q,s_{0}+2}^{\gamma,\mathcal{O}}\leqslant\varepsilon_{0},$$
					then the vector field $X_{P}$ defined in \eqref{def XP} satisfies the following estimates
					\begin{enumerate}[label=(\roman*)]
						\item $\| X_{P}(r)\|_{q,s}^{\gamma,\mathcal{O}}\lesssim \| r\|_{q,s+1}^{\gamma,\mathcal{O}}.$
						\item $\| d_{r}X_{P}(r)[\rho]\|_{q,s}^{\gamma,\mathcal{O}}\lesssim\|\rho\|_{q,s+1}^{\gamma,\mathcal{O}}+\| r\|_{q,s+2}^{\gamma,\mathcal{O}}\|\rho\|_{q,s_{0}+1}^{\gamma,\mathcal{O}}.$
						\item 
						$\| d_r^{2}X_{P}(r)[\rho_{1},\rho_{2}]\|_{q,s}^{\gamma,\mathcal{O}}\lesssim\|\rho_{1}\|_{q,s_{0}+1}^{\gamma,\mathcal{O}}\|\rho_{2}\|_{q,s+2}^{\gamma,\mathcal{O}}+\|\rho_{1}\|_{q,s+2}^{\gamma,\mathcal{O}}\|\rho_{2}\|_{q,s_{0}+1}^{\gamma,\mathcal{O}}+\| r\|_{q,s+2}^{\gamma,\mathcal{O}}\|\rho_{1}\|_{q,s_{0}+1}^{\gamma,\mathcal{O}}\|\rho_{2}\|_{q,s_{0}+1}^{\gamma,\mathcal{O}}$.
				\end{enumerate}}
			\end{lem}
			\begin{proof}
				\textbf{(i)} According to \eqref{Fourier coefficients of mathcalKlambda}, the Fourier coefficients of $\partial_{\theta}\mathcal{K}_{\lambda}$ are $\left(\ii jI_{j}(\lambda)K_{j}(\lambda)\right)_{j\in\mathbb{Z}}$. Hence
				$$\|\partial_{\theta}\mathcal{K}_{\lambda}\ast r\|_{H^{s}}^{2}=\sum_{(l,j)\in\mathbb{Z}^{d }\times\mathbb{Z}}\langle l,j\rangle^{2s}j^{2}I_{|j|}^{2}(\lambda)K_{|j|}^{2}(\lambda)|r_{l,j}|^{2}\leqslant\tfrac{1}{4}\| r\|_{H^{s}}^{2}.$$
				Notice that the last  inequality is obtained by the decay property of the product $I_{j}K_{j}$ on $\mathbb{R}_{+}^{*}$, \eqref{symmetry Bessel} and \eqref{asymptotic expansion of small argument}. Thus we deduce that
				$$\| \partial_{\theta}\mathcal{K}_{\lambda}\ast r\|_{H^{s}}\leqslant\tfrac{1}{2}\| r\|_{H^{s}}\leqslant\| r\|_{H^{s}}.$$
				Now we claim that
				\begin{equation}\label{estimate kernel equilibrium}
					\|\partial_{\theta}\mathcal{K}_{\lambda}\ast r\|_{q,s}^{\gamma ,\mathcal{O}}\lesssim\| r\|_{q,s}^{\gamma ,\mathcal{O}}.
				\end{equation}
				Indeed, from  \eqref{Fourier coefficients of mathcalKlambda}, we infer that  
					\begin{align}\label{Split-K}\partial_{\theta}\mathcal{K}_{\lambda}\ast r=\sum_{(l,j)\in\mathbb{Z}^{d}\times\mathbb{Z}}\ii jI_{j}(\lambda)K_{j}(\lambda)r_{l,j}(\lambda,\omega)\mathbf{e}_{l,j}.
					\end{align}
					At this stage we need to explore the regularity of the multiplier with respect to $\lambda$. By using \eqref{Bessel}, we write
					$$I_{j}(\lambda)K_{j}(\lambda)=\tfrac{2(-1)^{j}}{\pi}\int_{0}^{\frac{\pi}{2}}K_{0}\left(2\lambda\cos(\tau)\right)\cos(2j\tau)d\tau.$$
					From  \eqref{explicit form for K0}, we have the decomposition
					\begin{align}\label{Split-Kern0}
						K_{0}(z)=-\log(z/2)I_{0}(z)+f(z),
					\end{align}
					with $I_{0}$ being the modified Bessel function of the  first kind and $f$  an analytic function. By the morphism property of the logarithm, we get
					\begin{align*}
						I_{j}(\lambda)K_{j}(\lambda)&=-\log(\lambda)\tfrac{2(-1)^{j}}{\pi}\int_{0}^{\frac{\pi}{2}}I_{0}(2\lambda\cos(\tau))\cos(2j\tau)d\tau\\
						&\quad-\tfrac{2(-1)^{j}}{\pi}\int_{0}^{\frac{\pi}{2}}\log(\cos(\tau))\cos(2j\tau)d\tau\\
						&\quad-\tfrac{2(-1)^{j}}{\pi}\int_{0}^{\frac{\pi}{2}}\log(\cos(\tau))\left(I_{0}(2\lambda\cos(\tau))-1\right)\cos(2j\tau)d\tau\\
						&\quad+\tfrac{2(-1)^{j}}{\pi}\int_{0}^{\frac{\pi}{2}}f (2\lambda\cos(\tau))\cos(2j\tau)d\tau\\
						&:=\mathcal{I}_{1,j}(\lambda)+\mathcal{I}_{2,j}+\mathcal{I}_{3,j}(\lambda)+\mathcal{I}_{4,j}(\lambda).
					\end{align*}
					Since $I_{0}$ and $f$ are analytic, then the above expressions are smooth with respect to the parameter $\lambda\in(\lambda_{0},\lambda_{1})\subset\mathbb{R}_{+}^{*}.$ An integration by parts in $\mathcal{I}_{1,j}(\lambda)$ and $\mathcal{I}_{4,j}(\lambda)$ yields
					$$\forall\, i\in\{1,4\},\quad\sup_{j\in\mathbb{Z}}\left(|j|\max_{n\in\llbracket 0,q\rrbracket}\|\partial_{\lambda}^{(n)}\mathcal{I}_{i,j}\|_{L^{\infty}([\lambda_{0},\lambda_{1}])}\right)\lesssim 1.$$
					Looking at the definition of $I_{0}$ in \eqref{definition of modified Bessel function of first kind}, we see that we have uniformly in $\lambda\in [\lambda_{0},\lambda_{1}]$, 
					$$\forall\, n\in\llbracket 0,q\rrbracket,\quad\partial_{\lambda}^{(n)}\left(I_{0}\left(2\lambda\cos(\tau)\right)-1\right)=O\left(\cos(\tau)\right).$$
					Hence, an integration by parts in $\mathcal{I}_{3,j}(\lambda)$ yields 
					$$\sup_{j\in\mathbb{Z}}\left(|j|\max_{n\in\llbracket 0,q\rrbracket}\|\partial_{\lambda}^{(n)}\mathcal{I}_{3,j}\|_{L^{\infty}([\lambda_{0},\lambda_{1}])}\right)\lesssim 1.$$
					It remains to study the integral $\mathcal{I}_{2,j}.$ One can easily check from the above decomposition that 
					$$
					\mathcal{I}_{2,j}=\lim_{\lambda\to 0^+} I_{j}(\lambda)K_{j}(\lambda).
					$$
					Using \eqref{asymptotic expansion of small argument}, we then find
					$$
					\mathcal{I}_{2,j}=\frac{1}{2j}\cdot
					$$
					Putting together the preceding estimates, we obtain
					$$\sup_{j\in\mathbb{Z}}\left(|j|\max_{n\in\llbracket 0,q\rrbracket}\| \partial_{\lambda}^{(n)}\left(I_{j}K_{j}\right)\|_{L^{\infty}([\lambda_{0},\lambda_{1}])}\right)\lesssim 1.$$
					Then coming back to \eqref{Split-K} and using Leibniz formula, we obtain \eqref{estimate kernel equilibrium}. On the other hand, applying Lemma \ref{Lem-lawprod}-(iv)-(c) we get
				\begin{equation}\label{estimate-T1}
					\|({I}_1{K}_1)\,\partial_\theta r\|_{q,s}^{\gamma ,\mathcal{O}}\lesssim \| r\|_{q,s+1}^{\gamma ,\mathcal{O}}.
				\end{equation}
				Next we shall move to the estimate of $F_\lambda[r]$ defined in \eqref{definition of Flambda}. According to \eqref{formula}
				we may write
				\begin{align}\label{expression of Ar with vr1}
					\nonumber A_{r}(\varphi,\theta,\eta)&=2\left|\sin\left(\tfrac{\eta-\theta}{2}\right)\right|\left(\left(\tfrac{R(\varphi,\eta)-R(\varphi,\theta)}{2\sin\left(\frac{\eta-\theta}{2}\right)}\right)^{2}+R(\varphi,\eta)R(\varphi,\theta)\right)^{\frac{1}{2}}\\
					&:=2\left|\sin\left(\tfrac{\eta-\theta}{2}\right)\right|v_{r,1}(\varphi,\theta,\eta).
				\end{align}
				Notice that $v_{r,1}$ is smooth when $r$ is smooth and small enough, and $v_{0,1}=1.$ More precisely, Lemma \ref{Lem-lawprod}-(v) combined with  \mbox{Lemma \ref{cheater lemma}} allow  to get 
				\begin{align}\label{estimate vr1}
					\sup_{\eta\in\mathbb{T}}\| v_{r,1}(\ast,\cdot,\centerdot,\eta+\centerdot)-1\|_{q,s}^{\gamma,\mathcal{O}}&\lesssim \| r\|_{q,s+1}^{\gamma,\mathcal{O}},\nonumber\\
					\forall k\in\mathbb{N}^{*},\,\sup_{\eta\in\mathbb{T}}\| (\partial_{\theta}^{k}v_{r,1})(\ast,\cdot,\centerdot,\eta+\centerdot)\|_{q,s}^{\gamma,\mathcal{O}}&\lesssim \| r\|_{q,s+1+k}^{\gamma,\mathcal{O}}.
				\end{align}
			Here and in the sequel, the symbols $\ast,\cdot,\centerdot$ denote the variables $\mu=(\lambda,\omega),\varphi,\theta$, respectively. 				Then from the identity \eqref{Split-Kern0} we infer
				\begin{align}\label{ff}
						K_{0}(\lambda A_{r}(\varphi,\theta,\eta)) & =  K_{0}\left(2\lambda\left|\sin\left(\tfrac{\eta-\theta}{2}\right)\right|\right)+\log\left(\lambda\left|\sin\left(\tfrac{\eta-\theta}{2}\right)\right|\right)\left[I_{0}\left(2\lambda\left|\sin\left(\tfrac{\eta-\theta}{2}\right)\right|\right)-I_{0}\big(\lambda A_{r}(\varphi,\theta,\eta)\big)\right]\nonumber\\
						&  -\log\big(v_{r,1}(\varphi,\theta,\eta)\big)I_{0}\big(\lambda A_{r}(\varphi,\theta,\eta)\big)+f(\lambda A_{r}\big(\varphi,\theta,\eta)\big)-f\left(2\lambda\left|\sin\left(\tfrac{\eta-\theta}{2}\right)\right|\right).
				\end{align}
				By virtue of the expansion   \eqref{definition of modified Bessel function of first kind},  we can write 
				$$I_{0}\left(2\lambda\left|\sin\left(\tfrac{\eta-\theta}{2}\right)\right|\right)-I_{0}\big(\lambda A_{r}(\varphi,\theta,\eta)\big)=\sin^{2}\left(\tfrac{\eta-\theta}{2}\right)\mathscr{K}_{r,1}^{1}(\lambda,\varphi,\theta,\eta),$$
				with $\mathscr{K}_{r,1}^{1}$ being  smooth and vanishing at $r=0$. More precisely, we have the expansion
				\begin{equation}\label{definition of mathscrK11}
					\mathscr{K}_{r,1}^{1}(\lambda,\varphi,\theta,\eta)=\sum_{m=1}^{\infty}\tfrac{(2\lambda)^{2m}}{(m!)^{2}}\sin^{2m-2}\left(\tfrac{\eta-\theta}{2}\right)\left(1-v_{r,1}^{2m}(\varphi,\theta,\eta)\right).
				\end{equation}
				Now our aim is to establish the following estimate.
				\begin{equation}\label{estimate kernel Kr11}
					\forall\, k\in\mathbb{N},\quad\sup_{\eta\in\mathbb{T}}\|(\partial_{\theta}^{k}\mathscr{K}_{r,1}^{1})(\ast,\cdot,\centerdot,\eta+\centerdot)\|_{q,s}^{\gamma,\mathcal{O}}\lesssim\| r\|_{q,s+1+k}^{\gamma,\mathcal{O}}.
				\end{equation}
				For this goal we apply  Taylor Formula at the  order $2$,
					\begin{align*}
						&I_{0}\big(\lambda A_{r}(\varphi,\theta,\eta)\big)-I_{0}\left(2\lambda\left|\sin\left(\tfrac{\eta-\theta}{2}\right)\right|\right)
						=2\lambda\left|\sin\left(\tfrac{\eta-\theta}{2}\right)\right|\Big(v_{r,1}(\varphi,\theta,\eta)-1\Big)I_{0}'\left(2\lambda\left|\sin\left(\tfrac{\eta-\theta}{2}\right)\right|\right)\\
						&+4\lambda^{2}\sin^{2}\left(\tfrac{\eta-\theta}{2}\right)\Big(v_{r,1}(\varphi,\theta,\eta)-1\Big)^{2}\bigintssss_{0}^{1}(1-t)I_{0}''\left(2\lambda\left|\sin\left(\tfrac{\eta-\theta}{2}\right)\right|\big(1-t+tv_{r,1}(\varphi,\theta,\eta)\big)\right)dt.
					\end{align*}
					Consequently,  the kernel $\mathscr{K}_{r,1}^{1}$ can be rewritten into the form
					\begin{align}\label{ExpressionK}
					 \mathscr{K}_{r,1}^{1}(\lambda,\varphi,\theta,\eta)&=2\lambda\Big(1-v_{r,1}(\varphi,\theta,\eta)\Big)\frac{I_{0}'\left(2\lambda\left|\sin\left(\tfrac{\eta-\theta}{2}\right)\right|\right)}{\left|\sin\left(\tfrac{\eta-\theta}{2}\right)\right|}\\
							\nonumber&-4\lambda^{2}\Big(v_{r,1}(\varphi,\theta,\eta)-1\Big)^{2}\int_{0}^{1}(1-t)I_{0}''\left(2\lambda\left|\sin\left(\tfrac{\eta-\theta}{2}\right)\right|(1-t+tv_{r,1}(\varphi,\theta,\eta))\right)dt.
				\end{align}
				Using the structure \eqref{definition of modified Bessel function of first kind} and Lemma \ref{Lem-lawprod}-(iv)-(v) combined with \eqref{estimate vr1} we deduce  the estimate \eqref{estimate kernel Kr11}.
				%\begin{equation}\label{estimate kernel Kr11}
				%\forall k\in\{0,1,2\},\sup_{\eta\in\mathbb{T}}\|(\partial_{\theta}^{k}\mathscr{K}_{r,1}^{1})(\lambda,\varphi,\theta,\eta+\theta)\|_{q,s}^{\gamma,\mathcal{O}}\lesssim\| r\|_{q,s+2+k}^{\gamma,\mathcal{O}}.
				%\end{equation}
				Coming back to \eqref{ff} and set  
					\begin{align}\label{definition of mathscrK12}
						\mathscr{K}_{r,1}^{2}(\lambda,\varphi,\theta,\eta) & =  \log(\lambda)\sin^{2}\left(\tfrac{\eta-\theta}{2}\right)\mathscr{K}_{r,1}^{1}(\lambda,\varphi,\theta,\eta)-\log(v_{r,1}(\varphi,\theta,\eta))I_{0}(\lambda A_{r}(\varphi,\theta,\eta))\nonumber\\
						& \quad +f(\lambda A_{r}(\varphi,\theta,\eta))-f\left(2\lambda\left|\sin\left(\tfrac{\eta-\theta}{2}\right)\right|\right).
					\end{align}
				Then, by virtue of the law products and the composition laws of Lemma \ref{Lem-lawprod} combined with  \eqref{estimate vr1},  \eqref{estimate kernel Kr11} and the fact that $f$ is analytic and even, we get
				\begin{equation}\label{estimate kernel Kr12}
					\forall k\in\mathbb{N},\quad\sup_{\eta\in\mathbb{T}}\|(\partial_{\theta}^{k}\mathscr{K}_{r,1}^{2})(\ast,\cdot,\centerdot,\eta+\centerdot)\|_{q,s}^{\gamma,\mathcal{O}}\lesssim\| r\|_{q,s+1+k}^{\gamma,\mathcal{O}}.
				\end{equation}
				Consequently we obtain the decomposition 
				\begin{equation}\label{fundamental decomposition of K0(Ar)}
					K_{0}(\lambda A_{r}(\varphi,\theta,\eta)) =K_{0}\left(2\lambda\left|\sin\left(\tfrac{\eta-\theta}{2}\right)\right|\right)+\mathscr{K}(\eta-\theta)\mathscr{K}_{r,1}^{1}(\lambda,\varphi,\theta,\eta)+\mathscr{K}_{r,1}^{2}(\lambda,\varphi,\theta,\eta),
				\end{equation}
				where $\mathscr{K}$ is defined by
				\begin{equation}\label{definition of mathscrK}
					\mathscr{K}(\theta)=\sin^{2}\left(\tfrac{\theta}{2}\right)\log\left(\left|\sin\left(\tfrac{\theta}{2}\right)\right|\right)
				\end{equation}
				and the functions $\mathscr{K}_{r,1}^{1}$ and $\mathscr{K}_{r,1}^{2}$ satisfy the estimates \eqref{estimate kernel Kr11} and \eqref{estimate kernel Kr12}. We can obviously check that  $\mathscr{K}$ is an even function satisfying
				\begin{equation}\label{remark on the integrability of mathscrK}
					\mathscr{K},\,\partial_{\theta}\mathscr{K}\in L^{\infty}(\mathbb{T},\mathbb{R})\subset L^{1}(\mathbb{T},\mathbb{R})\quad\mbox{ and }\quad\partial_{\theta}^{2}\mathscr{K}\in L^{1}(\mathbb{T},\mathbb{R})\setminus L^{\infty}(\mathbb{T},\mathbb{R}).
				\end{equation}
				Introduce 
				\begin{equation}\label{definition of mathbbK1}
					\mathbb{K}_{r,1}(\lambda,\varphi,\theta,\eta)=\mathscr{K}(\eta-\theta)\mathscr{K}_{r,1}^{1}(\lambda,\varphi,\theta,\eta)+\mathscr{K}_{r,1}^{2}(\lambda,\varphi,\theta,\eta).
				\end{equation}
			Hence, putting together \eqref{estimate kernel Kr11}, \eqref{estimate kernel Kr12} and \eqref{remark on the integrability of mathscrK}, we obtain 
				\begin{equation}\label{estimate kernel mathbbK1}
					\forall k\in\{0,1\},\quad\sup_{\eta\in\mathbb{T}}\|(\partial_{\theta}^{k}\mathbb{K}_{r,1})(\ast,\cdot,\centerdot,\eta+\centerdot)\|_{q,s}^{\gamma,\mathcal{O}}\lesssim\| r\|_{q,s+1+k}^{\gamma,\mathcal{O}}.
				\end{equation}
				In addition, if $r(-\varphi,-\theta)=r(\varphi,\theta)$, then the kernel $\mathbb{K}_{r,1}$ satisfies the following symmetry property
				\begin{equation}\label{symmetry kernel K1}
					\mathbb{K}_{r,1}(\lambda,-\varphi,-\theta,-\eta)=\mathbb{K}_{r,1}(\lambda,\varphi,\theta,\eta).
				\end{equation}
				Plugging \eqref{definition of mathbbK1} and \eqref{fundamental decomposition of K0(Ar)} into $F_\lambda[r]$ defined in \eqref{definition of Flambda} yields
				$$F_\lambda[r](\varphi,\theta)= \displaystyle\bigintssss_{\mathbb{T}}\left(K_{0}\left(2\lambda\left|\sin\left(\tfrac{\eta-\theta}{2}\right)\right|\right)+\mathbb{K}_{r,1}(\lambda,\varphi,\theta,\eta)\right)\partial_{\theta\eta}^{2}\big[R(\varphi,\theta)R(\varphi,\eta)\sin(\eta-\theta)\big]d\eta.$$
				We denote
				$$f(\varphi,\theta,\eta)=\partial_{\theta\eta}^{2}\Big(R(\varphi,\theta)R(\varphi,\eta)\sin(\eta-\theta)\Big),
				$$
				then straightforward computations yield
				$$f(\varphi,\theta,\eta)=\Big(\partial_{\theta}R(\theta)\partial_{\eta}R(\eta)+R(\theta)R(\eta)\Big)\sin(\eta-\theta)+\Big(\partial_{\theta}R(\theta)R(\eta)-\partial_{\eta}R(\eta)R(\theta)\Big)\cos(\eta-\theta).$$
				We immediately deduce by law products, translation invariance property and composition laws in Lemma \ref{Lem-lawprod} that
				\begin{equation}\label{estimate partialthetaeta}
					\| f(\ast,\cdot,\centerdot,\eta+\centerdot)-\sin(\eta)\|_{q,s}^{\gamma,\mathcal{O}}\lesssim\| r\|_{q,s+1}^{\gamma,\mathcal{O}}.
				\end{equation}
				Using a change of variables, we obtain
				\begin{align*}F_\lambda[r](\varphi,\theta)&=\bigintssss_{\mathbb{T}}\Big(K_{0}\left(2\lambda\left|\sin\left(\tfrac{\eta}{2}\right)\right|\Big)+\mathbb{K}_{r,1}(\lambda,\varphi,\theta,\theta+\eta)\right)\Big(f(\varphi,\theta,\eta+\theta)-\sin\eta\Big)d\eta\\
				&\quad+\bigintssss_{\mathbb{T}}\Big(K_{0}\left(2\lambda\left|\sin\left(\tfrac{\eta}{2}\right)\right|\Big)+\mathbb{K}_{r,1}(\lambda,\varphi,\theta,\theta+\eta)\right)\sin\eta\, d\eta.				
				\end{align*}
			By symmetry, we  find
			$$
			\bigintssss_{\mathbb{T}}K_{0}\left(2\lambda\left|\sin\left(\tfrac{\eta}{2}\right)\right|\right)\sin\eta\, d\eta=0,
			$$
			which allows to get
			\begin{align*}F_\lambda[r](\varphi,\theta)&=\bigintssss_{\mathbb{T}}\Big(K_{0}\left(2\lambda\left|\sin\left(\tfrac{\eta}{2}\right)\right|\Big)+\mathbb{K}_{r,1}(\lambda,\varphi,\theta,\theta+\eta)\right)\Big(f(\varphi,\theta,\eta+\theta)-\sin\eta\Big)d\eta\\
				&\quad+\bigintssss_{\mathbb{T}}\mathbb{K}_{r,1}\big(\lambda,\varphi,\theta,\theta+\eta\big)\sin\eta\, d\eta.	
				\end{align*}			
				Recall that $K_{0}$ admits a logarithmic  behavior around  $0$, hence it is integrable at $0.$ Therefore, using the  law products of Lemma \ref{Lem-lawprod}, \eqref{estimate kernel mathbbK1}, \eqref{estimate partialthetaeta} and the smallness property on $r$, we infer 
				\begin{equation}\label{estimateFlambda}
					\|F_{\lambda}[r]\|_{q,s}^{\gamma,\mathcal{O}}\lesssim\|r\|_{q,s+1}^{\gamma,\mathcal{O}}.
				\end{equation}
				Combining \eqref{estimateFlambda} with \eqref{estimate kernel equilibrium} and \eqref{estimate-T1} achieves the proof of the first point.\\
				\\
				\textbf{(ii)} From \eqref{definition of mathbfLr}, \eqref{definition of mathcalKlambda}, \eqref{fundamental decomposition of K0(Ar)} and \eqref{definition of mathbbK1} we deduce that the operator $\mathbf{L}_{r}$ writes 
				\begin{equation}\label{link non local operators at state r and at equilibrium}
					\mathbf{L}_{r}=\mathcal{K}_{\lambda}\ast\cdot+\mathbf{L}_{r,1},
				\end{equation}
				where $\mathbf{L}_{r,1}$ is the integral operator of kernel $\mathbb{K}_{r,1}$ introduced in \eqref{definition of mathbbK1}. From Lemma \ref{lemma general form of the linearized operator} and its proof we find 
				$$d_{r}F_{\lambda}[r]\rho=\partial_{\theta}\left({({V}_{r}-\Omega)}\rho\right)-\partial_{\theta}\mathcal{K}_{\lambda}\ast\rho-\partial_{\theta}\mathbf{L}_{r,1}\rho.$$
				Thus, we get according to the definition \eqref{def XP}
				\begin{equation}\label{expression of dXP(r)}
					d_{r}X_{P}(r)\rho=\partial_{\theta}\mathbf{L}_{r,1}\rho-\partial_{\theta}\left((V_{r}-V_{0})\rho\right).
				\end{equation}
				Coming back to \eqref{definition of Vr} and using the kernel  decomposition \eqref{fundamental decomposition of K0(Ar)} together with the law products, the  composition laws in Lemma \ref{Lem-lawprod} and  the smallness condition, we deduce for any $s\geqslant s_0$,
				\begin{align}\label{Estim-Vr-October}
				\|V_{r}-V_{0}\|_{q,s}^{\gamma,\mathcal{O}}\lesssim \|r\|_{q,s+1}^{\gamma,\mathcal{O}}.
				\end{align}
				Therefore, we obtain from  the law products, \eqref{Estim-Vr-October} and   the smallness property on $r$, 
				\begin{align*}
					\|\partial_{\theta}\left(\left(V_{r}-V_{0})\right)\rho\right)\|_{q,s}^{\gamma,\mathcal{O}}&\lesssim\|V_{r}-V_{0}\|_{q,s+1}^{\gamma,\mathcal{O}}\|\rho\|_{q,s_{0}+1}^{\gamma,\mathcal{O}}+\|V_{r}-V_{0}\|_{q,s_{0}+1}^{\gamma,\mathcal{O}}\|\rho\|_{q,s+1}^{\gamma,\mathcal{O}}\\
					&\lesssim\|\rho\|_{q,s+1}^{\gamma,\mathcal{O}}+\|r\|_{q,s+2}^{\gamma,\mathcal{O}}\|\rho\|_{q,s_{0}+1}^{\gamma,\mathcal{O}}.
				\end{align*}
				Now, by using the last point in  Lemma \ref{lemma symmetry and reversibility} with \eqref{estimate kernel mathbbK1} and the smallness property on $r$, we obtain
				$$\|\partial_{\theta}\mathbf{L}_{r,1}\rho\|_{q,s}^{\gamma,\mathcal{O}}\lesssim\|\rho\|_{q,s}^{\gamma,\mathcal{O}}+\|r\|_{q,s+2}^{\gamma,\mathcal{O}}\|\rho\|_{q,s_{0}}^{\gamma,\mathcal{O}}.$$
				Putting together \eqref{expression of dXP(r)} and the last two estimates gives
				$$\|d_{r}X_{P}(r)\rho\|_{q,s}^{\gamma,\mathcal{O}}\lesssim\|\rho\|_{q,s+1}^{\gamma,\mathcal{O}}+\|r\|_{q,s+2}^{\gamma,\mathcal{O}}\|\rho\|_{q,s_{0}+1}^{\gamma,\mathcal{O}}.$$\\
				\textbf{(iii)} Differentiating in $r$ the identity  \eqref{expression of dXP(r)} yields,
				\begin{equation}\label{Second-diff}
					d_{r}^{2}X_{P}(r)[\rho_{1},\rho_{2}]=\partial_{\theta}\left(d_{r}\mathbf{L}_{r,1}(r)[\rho_{2}]\rho_{1}\right)-\partial_{\theta}\left(\big(d_{r}{V}_{r}(r)[\rho_{2}]\big)\rho_{1}\right).
				\end{equation}
				For the first member of the right-hand side we first recall from \eqref{definition of mathbbK1} that 
				\begin{equation*}%\label{definition of mathbbK1D}
					\mathbf{L}_{r,1}\rho(\varphi,\theta)=\int_{\mathbb{T}}\rho(\varphi,\eta)\big[\mathscr{K}(\eta-\theta)\mathscr{K}_{r,1}^{1}(\lambda,\varphi,\theta,\eta)+\mathscr{K}_{r,1}^{2}(\lambda,\varphi,\theta,\eta)\big] d\eta.
				\end{equation*} 
				Hence, by differentiation and change of variables, we obtain
				\begin{align}\label{definition of mathbbK1M}
					 d_{r}\mathbf{L}_{r,1}(r)[\rho_{2}]\rho_{1}(\varphi,\theta)&=\int_{\mathbb{T}}\rho_1(\varphi,\eta)\big[\mathscr{K}(\eta-\theta)\big(d_{r}\mathscr{K}_{r,1}^{1}\big)[\rho_2](\varphi,\theta,\eta)+d_{r}\mathscr{K}_{r,1}^{2}\big)[\rho_2](\varphi,\theta,\eta)\big] d\eta\\
					\nonumber&=\int_{\mathbb{T}}\rho_1(\varphi,\theta+\eta)\big[\mathscr{K}(\eta)\big(d_{r}\mathscr{K}_{r,1}^{1}\big)[\rho_2](\varphi,\theta,\theta+\eta)+d_{r}\mathscr{K}_{r,1}^{2}\big)[\rho_2](\varphi,\theta,\theta+\eta)\big] d\eta.
				\end{align} 
				Coming back to \eqref{ExpressionK}, we emphasize   that the dependence in $r$ of the functional $ \mathscr{K}_{r,1}^{1}$ is smooth  since the function $v_{r,1}$, introduced in \eqref{expression of Ar with vr1}, depends smoothly in $r$.  In addition   $d_{r}\mathscr{K}_{r,1}^{1}$ can be easily related  to  $d_{r}v_{r,1}$.
				From straightforward calculus we see that, for the sake of simple notation  we remove the dependence in the parameters and $\varphi$,  
				\begin{align}\label{First-diffvr}
					d_{r}v_{r,1}(r)[\rho](\theta,\eta)&=\frac{1}{v_{r,1}(\theta,\eta)}\left(\frac{R(\theta)-R(\eta)}{\sin^2\left(\frac{\eta-\theta}{2}\right)}\left(\frac{\rho(\theta)}{R(\theta)}-\frac{\rho(\eta)}{R(\eta)}\right)+\frac{\rho(\theta)R^{2}(\eta)+\rho(\eta)R^{2}(\theta)}{2R(\theta)R(\eta)}\right).
				\end{align}
				Therefore using \eqref{estimate vr1} combined with the  law products stated in Lemma \ref{Lem-lawprod}, Lemma \ref{cheater lemma} and the smallness condition of Lemma \ref{lemma estimates vector field XP}, we find that
				\begin{equation}\label{estimate vr1PP}
					\sup_{\eta\in\mathbb{T}}\| d_{r}v_{r,1}(r)[\rho](\ast,\cdot,\centerdot,\eta+\centerdot)\|_{q,s}^{\gamma,\mathcal{O}}\lesssim \| \rho\|_{q,s+1}^{\gamma,\mathcal{O}}+\| \rho\|_{q,s_0+1}^{\gamma,\mathcal{O}}\| r\|_{q,s+1}^{\gamma,\mathcal{O}}.
				\end{equation}
				Similarly to \eqref{estimate vr1PP}, one gets from \eqref{ExpressionK} and \eqref{definition of mathscrK12},
				\begin{equation}\label{estimate vr1PM}
					\forall \, i\in\{1,2\},\quad \sup_{\eta\in\mathbb{T}}\| d_{r}\mathscr{K}_{r,1}^{i}[\rho](\ast,\cdot,\centerdot,\eta+\centerdot)\|_{q,s}^{\gamma,\mathcal{O}}\lesssim \| \rho\|_{q,s+1}^{\gamma,\mathcal{O}}+\| \rho\|_{q,s_0+1}^{\gamma,\mathcal{O}}\| r\|_{q,s+1}^{\gamma,\mathcal{O}}.
				\end{equation}
				Inserting \eqref{estimate vr1PM} into \eqref{definition of mathbbK1M} and using once again  the  law products and the smallness condition we obtain,
				\begin{align}\label{First-member}
					\nonumber \|\partial_{\theta}d_{r}\mathbf{L}_{r,1}(r)[\rho_{2}]\rho_{1}\|_{q,s}^{\gamma,\mathcal{O}}\lesssim&\|d_{r}\mathbf{L}_{r,1}(r)[\rho_{2}]\rho_{1}\|_{q,s+1}^{\gamma,\mathcal{O}}\\
					\lesssim&\|\rho_{1}\|_{q,s+1}^{\gamma,\mathcal{O}}\|\rho_{2}\|_{q,s_{0}+1}^{\gamma,\mathcal{O}}+\|\rho_{1}\|_{q,s_{0}}^{\gamma,\mathcal{O}}\|\rho_{2}\|_{q,s+2}^{\gamma,\mathcal{O}}+\| r\|_{q,s+2}^{\gamma,\mathcal{O}}\|\rho_{1}\|_{q,s_{0}}^{\gamma,\mathcal{O}}\|\rho_{2}\|_{q,s_{0}+1}^{\gamma,\mathcal{O}}.
				\end{align}
				Next we shall move to the estimate of the last member of \eqref{Second-diff}. Differentiating the definition of ${V}_{r}$ in the proof of Lemma \ref{lemma general form of the linearized operator}, we infer
				\begin{align*}
					\nonumber d_{r}{V}_{r}(r)[\rho_{2}](\theta) & =  \displaystyle\bigintssss_{\mathbb{T}}K_{0}\left(\lambda A_{r}(\theta,\eta)\right)\partial_{\eta}\left(\tfrac{\rho_{2}(\eta)R^2(\theta)-\rho_2(\theta)R^2(\eta)}{R^3(\theta)R(\eta)}\sin(\eta-\theta)\right)d\eta\\
					\nonumber&\quad+  \displaystyle\tfrac{\lambda}{R(\theta)}\bigintssss_{\mathbb{T}}\tfrac{(R(\theta)-R(\eta))\left(\tfrac{\rho_{2}(\theta)}{R(\theta)}-\tfrac{\rho_{2}(\eta)}{R(\eta)}\right)}{A_{r}(\theta,\eta)}K_{0}'\left(\lambda A_{r}(\theta,\eta)\right)\partial_{\eta}(R(\eta)\sin(\eta-\theta))d\eta\\
					\nonumber& \quad+ \displaystyle{2\lambda}\bigintssss_{\mathbb{T}}\tfrac{\rho_{2}(\theta)R^{2}(\eta)+\rho_{2}(\eta)R^{2}(\theta)}{R^2(\theta)R(\eta)A_{r}(\theta,\eta)}\sin^{2}\left(\tfrac{\eta-\theta}{2}\right)K_{0}'\left(\lambda A_{r}(\theta,\eta)\right)\partial_{\eta}(R(\eta)\sin(\eta-\theta))d\eta\\
					&:=\mathcal{I}_1(\theta)+\mathcal{I}_2(\theta)+\mathcal{I}_3(\theta).
				\end{align*}
				The estimate of $\mathcal{I}_1$ is similar to \eqref{estimateFlambda}, we use the same tools and one finds
				\begin{equation}\label{estimateFlambdaP-0}
					\| \mathcal{I}_1\|_{q,s}^{\gamma,\mathcal{O}}\lesssim  \| \rho_2\|_{q,s+1}^{\gamma,\mathcal{O}}+\| \rho_2\|_{q,s_0+1}^{\gamma,\mathcal{O}}\| r\|_{q,s+1}^{\gamma,\mathcal{O}}.
				\end{equation}
				For the terms $\mathcal{I}_2$ and $\mathcal{I}_3$ the computations are straightforward and we shall only extract their main parts and give the suitable estimates. For this goal  we differentiate    \eqref{explicit form for K0}, leading to
				$$K_{0}'(z)=\frac{-1}{z}+\log(z)F(z)+G(z),$$
				with $F$ and $G$ being entire functions. Hence, applying \eqref{expression of Ar with vr1}, we deduce that  $\mathcal{I}_2$ takes the form 
				$$\mathcal{I}_2(\theta)=-\displaystyle\frac{1}{4}\bigintsss_{\mathbb{T}}\tfrac{(R(\theta)-R(\eta))\left(\frac{\rho_{2}(\theta)}{R(\theta)}-\frac{\rho_{2}(\eta)}{R(\eta)}\right)}{R^2(\theta)v_{r,1}^{2}(\theta,\eta)\sin^2\left(\frac{\eta-\theta}{2}\right)}\partial_{\eta}(R(\eta)\sin(\eta-\theta))d\eta+\textnormal{l.o.t}.$$
				Hence we proceed as for \eqref{estimate vr1PM} and one finds
				\begin{equation}\label{estimateFlambdaP-1}
					\| \mathcal{I}_2\|_{q,s}^{\gamma,\mathcal{O}}\lesssim  \| \rho_2\|_{q,s+1}^{\gamma,\mathcal{O}}+\| \rho_2\|_{q,s_0+1}^{\gamma,\mathcal{O}}\| r\|_{q,s+1}^{\gamma,\mathcal{O}}.
				\end{equation}
				As to the last term $ \mathcal{I}_3$, we write 
				$$\mathcal{I}_3(\theta)= -\frac12\displaystyle\bigintssss_{\mathbb{T}}\tfrac{\rho_{2}(\theta)R^{2}(\eta)+\rho_{2}(\eta)R^{2}(\theta)}{R^2(\theta)R(\eta)v_{r,1}^{2}(\theta,\eta)}\partial_{\eta}(R(\eta)\sin(\eta-\theta))d\eta+\textnormal{l.o.t}.$$
				Then, we get
				\begin{equation}\label{estimateFlambdaP-2}
					\| \mathcal{I}_3\|_{q,s}^{\gamma,\mathcal{O}}\lesssim  \| \rho_2\|_{q,s}^{\gamma,\mathcal{O}}+\| \rho_2\|_{q,s_0}^{\gamma,\mathcal{O}}\| r\|_{q,s+1}^{\gamma,\mathcal{O}}.
				\end{equation}
				Putting together \eqref{estimateFlambdaP-0}, \eqref{estimateFlambdaP-1} and \eqref{estimateFlambdaP-2} yields
				\begin{equation}\label{secondpart}
					\|d_{r}{V}_{r}(r)[\rho_{2}]\|_{q,s}^{\gamma,\mathcal{O}}\lesssim  \| \rho_2\|_{q,s+1}^{\gamma,\mathcal{O}}+\| \rho_2\|_{q,s_0+1}^{\gamma,\mathcal{O}}\| r\|_{q,s+1}^{\gamma,\mathcal{O}}.
				\end{equation}
				Therefore we obtain according to the law products in  Lemma \ref{Lem-lawprod}, \eqref{secondpart} and the smallness condition,
				\begin{align*}\big\| \partial_{\theta}\big(d_{r}{V}_{r}(r)[\rho_{2}]\rho_{1}\big)\big\|_{q,s}^{\gamma,\mathcal{O}}&\lesssim\big\| d_{r}{V}_{r}(r)[\rho_{2}]\big\|_{q,s+1}^{\gamma,\mathcal{O}}\big\| \rho_{1}\big\|_{q,s_0}^{\gamma,\mathcal{O}}+\big\| d_{r}{V}_{r}(r)[\rho_{2}]\big\|_{q,s_0}^{\gamma,\mathcal{O}}\big\| \rho_{1}\big\|_{q,s+1}^{\gamma,\mathcal{O}}\\
					&\lesssim\|\rho_{1}\|_{q,s_{0}}^{\gamma,\mathcal{O}}\|\rho_{2}\|_{q,s+2}^{\gamma,\mathcal{O}}+\| r\|_{q,s+2}^{\gamma,\mathcal{O}}\|\rho_{1}\|_{q,s_{0}}^{\gamma,\mathcal{O}}\|\rho_{2}\|_{q,s_{0}+1}^{\gamma,\mathcal{O}}+\|\rho_{1}\|_{q,s+1}^{\gamma,\mathcal{O}}\|\rho_{2}\|_{q,s_{0}+1}^{\gamma,\mathcal{O}}.
				\end{align*}
				Combining the latter estimate with \eqref{Second-diff} and  \eqref{First-member} allows to get
				$$\| d_r^{2}X_{P}(r)[\rho_{1},\rho_{2}]\|_{q,s}^{\gamma,\mathcal{O}}\lesssim\|\rho_{1}\|_{q,s_{0}}^{\gamma,\mathcal{O}}\|\rho_{2}\|_{q,s+2}^{\gamma,\mathcal{O}}+\| r\|_{q,s+2}^{\gamma,\mathcal{O}}\|\rho_{1}\|_{q,s_{0}}^{\gamma,\mathcal{O}}\|\rho_{2}\|_{q,s_{0}+1}^{\gamma,\mathcal{O}}+\|\rho_{1}\|_{q,s+1}^{\gamma,\mathcal{O}}\|\rho_{2}\|_{q,s_{0}+1}^{\gamma,\mathcal{O}}.$$
				Using Sobolev embeddings we get the desired result. This achieves the proof of Lemma \ref{lemma estimates vector field XP}.
			\end{proof}
As an application of Lemma \ref{lemma estimates vector field XP}, we shall  establish tame estimates for the Hamiltonian vector field 
		$$X_{\mathcal{P}_{\varepsilon}}=(\partial_{I}\mathcal{P}_{\varepsilon},-\partial_{\vartheta}\mathcal{P}_{\varepsilon},\Pi_{\mathbb{S}}^{\perp}\partial_{\theta}\nabla_{z}\mathcal{P}_{\varepsilon})$$
		defined through  \eqref{cNP} and \eqref{hamiltonian vector filed associated with Hepsilon}.
		\begin{lem}\label{tame estimates for the vector field XmathcalPvarepsilon}
			{Let $(\gamma,q,s_{0},s)$ satisfy \eqref{initial parameter condition}.
				There exists $\varepsilon_0\in(0,1)$ such that if 
				$$\varepsilon\leqslant\varepsilon_0\quad\textnormal{and}\quad\|\mathfrak{I}\|_{q,s_{0}+2}^{\gamma,\mathcal{O}}\leqslant 1,$$ 
				then  the  perturbed Hamiltonian vector field $X_{\mathcal{P}_{\varepsilon}}$ satisfies the following  estimates,
				\begin{enumerate}[label=(\roman*)]
					\item $\| X_{\mathcal{P}_{\varepsilon}}(i)\|_{q,s}^{\gamma,\mathcal{O}}\lesssim 1+\|\mathfrak{I}\|_{q,s+2}^{\gamma,\mathcal{O}}.$
					\item $\big\| d_{i}X_{\mathcal{P}_{\varepsilon}}(i)[\,\widehat{i}\,]\big\|_{q,s}^{\gamma,\mathcal{O}}\lesssim \|\,\widehat{i}\,\|_{q,s+1}^{\gamma,\mathcal{O}}+\|\mathfrak{I}\|_{q,s+2}^{\gamma,\mathcal{O}}\|\,\widehat{i}\,\|_{q,s_{0}+1}^{\gamma,\mathcal{O}}.$
					\item $\big\| d_{i}^{2}X_{\mathcal{P}_{\varepsilon}}(i)[\,\widehat{i},\widehat{i}\,]\big\|_{q,s}^{\gamma,\mathcal{O}}\lesssim \|\,\widehat{i}\,\|_{q,s+2}^{\gamma,\mathcal{O}}\|\,\widehat{i}\,\|_{q,s_{0}+2}^{\gamma,\mathcal{O}}+\|\mathfrak{I}\|_{q,s+2}^{\gamma,\mathcal{O}}\left(\|\,\widehat{i}\,\|_{q,s_{0}+2}^{\gamma,\mathcal{O}}\right)^{2}.$
			\end{enumerate}}
		\end{lem}
		\begin{proof}These estimates can be recovered  from  Lemma \ref{lemma estimates vector field XP} combined with  the following estimate on the action-angle change of variables introduced in \eqref{definition of A action-angle-normal}
			\begin{align}\label{v-theta1}
				\forall\alpha,\beta\in\mathbb{N}^{d},\quad \|\partial_{\vartheta}^{\alpha}\partial_{I}^{\beta}v(\vartheta,I)\|_{q,s}^{\gamma,\mathcal{O}}&\lesssim  1+\|\mathfrak{I}\|_{q,s}^{\gamma,\mathcal{O}}.
			\end{align}
		This estimate follows from  Lemma \ref{Lem-lawprod}-(iv)-(v) provided that $\|\vartheta\|_{q,s_0}^{\gamma,\mathcal{O}},\|I\|_{q,s_0}^{\gamma,\mathcal{O}}\leqslant 1.$ This latter condition is satisfied due to the smallness condition in the Lemma. For more details, we refer to \cite[Lem. 5.1]{BM18}.

		\end{proof}
\subsection{Berti-Bolle approach for the approximate inverse}\label{sec approx-inv}
In this section, we shall follow the remarkable procedure developed by Berti and Bolle in \cite{BB15} to construct an  approximate right inverse for  the linearized operator
\begin{equation}\label{lin func}
	d_{i,\alpha}\mathcal{F}(i_{0},\alpha_{0})[\widehat{i},\widehat{\alpha}]=\omega\cdot\partial_{\varphi}\widehat{i}-d_{i}X_{H_{\varepsilon}^{\alpha_{0}}}(i_{0}(\varphi))[\widehat{i}]-(\widehat{\alpha},0,0),
\end{equation}
where $\mathcal{F}$ is the nonlinear functional  defined in \eqref{main function}. This construction is crucial  for the Nash-Moser scheme that we shall perform later in Section \ref{Nash-Moser scheme and measure estimates}. From \eqref{rev-torus}, we denote by $i_0$ an embedded torus with
$$i_{0}(\varphi)=(\vartheta_{0}(\varphi),I_{0}(\varphi),z_{0}(\varphi))\quad\textnormal{and}\quad\mathfrak{I}_{0}(\varphi)=i_{0}(\varphi)-(\varphi,0,0).$$
Throughout  this section, we shall assume the following smallness condition : the application $(\lambda,\omega)\mapsto\mathfrak{I}_{0}(\lambda,\omega)$ is $q$-times differentiable on $\mathcal{O}$ and there exists $\varepsilon_0\in(0,1)$ (small enough) such that
\begin{equation}\label{small-C3}
	\|\mathfrak{I}_{0}\|_{q,s_{0}+2}^{\gamma,\mathcal{O}}+\|\alpha_0-\omega\|_{q}^{\gamma,\mathcal{O}}\leqslant \varepsilon_0.
\end{equation}
We mainly follow the same approach  as in  \cite{BB15}	which reduces the search of an approximate right inverse of \eqref{lin func} to the search of an approximate right inverse in the normal directions. The main difference with \cite{BB15} is to be able to bypass the use  of the  isotropic torus in a similar way to the recent paper \cite{HHM21}.  
\subsubsection{Triangularization up to error terms}
Given a linear operator $A\in\mathcal{L}(\mathbb{R}^{d},H_{\perp}^{s}),$ we define the transposed operator $A^{\top}:H_{\perp}^{s}\rightarrow\mathbb{R}^{d}$ by the duality relation
\begin{equation}\label{dual rel}
	\forall (u,v)\in H_{\perp}^{s}\times\mathbb{R}^{d},\quad\langle A^{\top}u,v\rangle_{\mathbb{R}^{d}}=\langle u,Av\rangle_{L^{2}(\mathbb{T})}.
\end{equation}
We introduce  the following change of coordinates
$ G_0 : (\phi, y, w) \to (\vartheta, I, z)$ of the phase space $\mathbb{T}^d \times \mathbb{R}^d \times H_{\perp}^{s}$ defined by
\begin{equation}\label{trasform-sympl}
	\begin{pmatrix}
		\vartheta \\
		I \\
		z
	\end{pmatrix} := G_0 \begin{pmatrix}
		\phi \\
		y \\
		w
	\end{pmatrix} := 
	\begin{pmatrix}
		\vartheta_0(\phi) \\
		I_0 (\phi) + L_1(\varphi)y+L_{2}(\varphi) w \\
		z_0(\phi) + w
	\end{pmatrix} ,
\end{equation}
where 
\begin{align}
	L_1(\varphi)&:=  [\partial_\varphi \vartheta_0(\varphi)]^{-\top},\label{L1}\\
	L_2(\varphi)&:= - [(\partial_\vartheta \widetilde{z}_0)(\vartheta_0(\varphi))]^\top \partial_\theta ^{-1},\label{L2}\\
	\widetilde{z}_0 (\vartheta) &:= z_0 (\vartheta_0^{-1} (\vartheta)),\label{z0t}
\end{align}
provided that $\vartheta_0:\mathbb{R}^d\to\mathbb{R}^d$ is a diffeomorphism.
Notice that one recovers the torus $i_0$ by taking in the new coordinates, the flat torus $i_{\textnormal{\tiny{flat}}}(\varphi)=(\varphi,0,0)$ namely
$$G_{0}(i_{\textnormal{\tiny{flat}}}(\varphi))=i_{0}(\varphi).$$
Next, we shall adopt the notation  ${\mathtt u}=(\phi, y,w)$ to denote  the new coordinates induced by $G_0$ in \eqref{trasform-sympl} and  we  simply set  ${\mathtt u}_0(\varphi)=i_{\textnormal{\tiny{flat}}}(\varphi)$. Now, to  measure to which extent an  embedded torus $i_0(\mathbb{T})$  is close to be   invariant for the Hamiltonian vector field $X_{H_\varepsilon^{\alpha_{0}}}$, we shall make appeal to the error function     
\begin{equation} \label{def Z}
	Z(\varphi) :=  (Z_1, Z_2, Z_3) (\varphi) := {\mathcal F}(i_0, \alpha_{0}) (\varphi) =
	\omega \cdot \partial_\varphi i_0(\varphi) - X_{H^{\alpha_{0}}_\varepsilon}(i_0(\varphi)) \, .
\end{equation}
We say that a quantity is of "type Z"  is it is  $O(Z)$, and particular it is vanishing  at an exact solution.
In the next Proposition, we study the conjugation of the linear operator $d_{i,\alpha}\mathcal{F}(i_{0},\alpha_{0})$ by the linear change of variables induced by $G_0$ defined in \eqref{trasform-sympl},
\begin{equation}\label{DG0}
	D G_0(\varphi, 0, 0) 
	\begin{pmatrix}
		\widehat \phi \, \\
		\widehat y \\
		\widehat w
	\end{pmatrix} 
	:= 
	\begin{pmatrix}
		\partial_\varphi \vartheta_0(\varphi) & 0 & 0 \\
		\partial_\varphi I_0(\varphi) &L_1(\varphi) & 
		L_2(\varphi) \\
		\partial_\varphi z_0(\varphi) & 0 & I
	\end{pmatrix}
	\begin{pmatrix}
		\widehat \phi \, \\
		\widehat y \\
		\widehat w
	\end{pmatrix} .
\end{equation}
The following result is proved in \cite{HHM21}.
\begin{prop}\label{prop conj DGO}
	The conjugation of the linearized operator $d_{i,\alpha} {\mathcal F} (i_0,\alpha_{0})$ by the linear change of variables $D G_0(\mathtt{u}_0)$ writes as follows
	\begin{align}\label{bbL}
		[D G_0({\mathtt u}_0)]^{-1}  d_{i,\alpha} {\mathcal F} (i_0,\alpha_0) D \widetilde{G}_0({\mathtt u}_0)
		[\widehat \phi, \widehat y, \widehat w, \widehat \alpha ]&=\mathbb{D}[\widehat \phi, \widehat y, \widehat w, \widehat \alpha ]+\mathbb{E}[\widehat \phi, \widehat y, \widehat w],
	\end{align}
where $\widetilde{G}_0$ is defined by
$$\widehat{G}_0(\mathtt{u},\alpha):=(G_{0}(\mathtt{u}),\alpha)$$
and where
\begin{enumerate}[label=(\roman*)]
	\item the operator $\mathbb{D}$ admits a triangular structure in the variables $(\widehat{\phi},\widehat{y},\widehat{w})$ in the form
		\begin{align*}\mathbb{D}[\widehat \phi, \widehat y, \widehat w, \widehat \alpha ]:=\left(
		\begin{array}{c}
			\omega\cdot\partial_\varphi  \widehat\phi-\big[\mathsf{K}_{20}(\varphi) \widehat y+\mathsf{K}_{11}^\top(\varphi) \widehat w+L_1^\top (\varphi)\widehat \alpha\big]
			\\
			\omega\cdot\partial_\varphi  \widehat y+\mathtt{B}(\varphi)\widehat \alpha \\
			\omega\cdot\partial_\varphi \widehat w-\partial_\theta\big[\mathsf{K}_{11}(\varphi) \widehat y+\mathsf{K}_{02}(\varphi)\widehat w -L_2^{\top}(\varphi) \widehat\alpha   \big] 
		\end{array}
		\right)\,,
	\end{align*}
$\mathtt{B}(\varphi)$ and $\mathsf{K}_{20}(\varphi)$ are $d\times d$ real matrices given by
\begin{align}
	\mathtt{B}(\varphi)&:=[\partial_{\varphi}\vartheta_{0}(\varphi)]^{\top}\partial_{\varphi}I_{0}(\varphi)L_{1}^{\top}(\varphi)+[\partial_{\varphi}z_0(\varphi)]^{\top}L_{2}^{\top}(\varphi),\label{ttb}\\
	\mathsf{K}_{20}(\varphi)&:=\varepsilon L_{1}^{\top}(\varphi)(\partial_{II}\mathcal{P}_{\varepsilon})(i_0(\varphi))L_{1}(\varphi),\label{K20}
\end{align}
$\mathsf{K}_{02}(\varphi)$ is a linear self-adjoint operator of $H_{\perp}^{s}$ in the form
\begin{align}
	\mathsf{K}_{02}(\varphi)&:=( \partial_{z}\nabla_z H_\varepsilon^{\alpha_{0}}) (i_0(\varphi))  +\varepsilon L_2^\top(\varphi) ( \partial_{II} \mathcal{P}_\varepsilon)(i_0(\varphi)) L_2(\varphi)\nonumber\\ 
	&\quad+\varepsilon L_2^\top(\varphi)( \partial_{zI} \mathcal{P}_\varepsilon) (i_0(\varphi)) + \varepsilon  (\partial_I\nabla_z \mathcal{P}_\varepsilon) (i_0(\varphi))  L_2(\varphi) \, \label{K02}
\end{align}
and $\mathsf{K}_{11}(\varphi)\in\mathcal{L}(\mathbb{R}^{d},H_{\perp}^s)$ is given by
\begin{align}
	\mathsf{K}_{11}(\varphi):=\varepsilon L_{2}^\top(\varphi)(\partial_{II}\mathcal{P}_{\varepsilon})(i_0(\varphi))L_{1}(\varphi)+\varepsilon(\partial_{I}\nabla_{z}\mathcal{P}_{\varepsilon})(i_0(\varphi))L_{1}(\varphi).\label{K11}
\end{align}
\item the operator $\mathbb{E}$ is an error term in the form
\begin{align*}
	\mathbb{E}[\widehat \phi, \widehat y, \widehat w]&:=
	[D G_0({\mathtt u}_0)]^{-1}    \partial_\varphi Z(\varphi) \widehat\phi\\ 
	&\quad + \left(
	\begin{array}{c}
		0 
		\\
		\mathtt{A}(\varphi)\big[\mathsf{K}_{20}(\varphi) \widehat y+\mathsf{K}_{11}^\top(\varphi) \widehat w\big]-R_{10}(\varphi) \widehat y -R_{01}(\varphi) \widehat w 
		\\
		0   
	\end{array}
	\right) ,
\end{align*}
where $\mathtt{A}(\varphi)$ and $R_{10}(\varphi)$ are $d\times d$ matrices defined by
\begin{align}\label{devcalA}
	\mathtt{A}(\varphi)&:=[\partial_\varphi \vartheta_0(\varphi)]^\top \partial_\varphi I_0(\varphi)-[\partial_\varphi I_0(\varphi)]^\top \partial_\varphi \vartheta_0(\varphi)  +[\partial_\varphi z_0(\varphi)]^{\top} \partial_\theta ^{-1} \partial_\varphi z_0(\varphi),\\
	R_{10}(\varphi)&:=[\partial_{\varphi}Z_{1}(\varphi)]^{\top}L_{1}(\varphi)\label{R10}
\end{align}
and $R_{01}(\varphi)\in\mathcal{L}(H_{\perp}^s,\mathbb{R}^{d})$ with
\begin{equation}\label{R01}
	R_{01}(\varphi):=-[\partial_{\varphi}Z_{1}(\varphi)]^{\top}L_{2}(\varphi)+[\partial_{\varphi}Z_{3}(\varphi)]^{\top}\partial_{\theta}^{-1}.
\end{equation}
\end{enumerate}
\end{prop}
Now  we recall the following result, for the proof we refer to \cite[Lemma 6.1]{HHM21} and  Lemmata 5.6-5.7  in \cite{BM18},\begin{lem}\label{lem:est-G} 
	The following assertions hold true.
	\begin{enumerate}[label=(\roman*)]
		\item The operator  $DG_0({\mathtt u}_0)$ and $[DG_0({\mathtt u}_0)]^{- 1} $ satisfy for all $\widehat{{\mathtt u}}=(\widehat{\phi},\widehat{y},\widehat{w})$,
		\begin{align*}
			\forall s\in[s_0,S],\quad\| [DG_0({\mathtt u}_0)]^{\pm 1}[\widehat{{\mathtt u}}]\|_{q,s}^{\gamma,\mathcal{O}}
			\lesssim\|\widehat{{\mathtt u}}\|_{q,s}^{\gamma,\mathcal{O}}+\|\mathfrak{I}_{0}\|_{q,s+1}^{\gamma,\mathcal{O}}\|\widehat{{\mathtt u}}\|_{q,s_0}^{\gamma,\mathcal{O}}\, . 
		\end{align*}
		\item  The operators $R_{10}$ and  $ R_{01}$,  defined in \eqref{R10} and \eqref{R01}, satisfy the estimates
		\begin{align*}
			\forall s\in[s_0,S],\quad& \| R_{10} \widehat y \|_{q,s}^{\gamma,\mathcal{O}} \lesssim \| Z \|_{q,s+ 1 }^{\gamma,\mathcal{O}}  \| \widehat y\|_{q,s_0 + 1 }^{\gamma,\mathcal{O}}
			+ \| Z \|_{q,s_0 + 1 }^{\gamma,\mathcal{O}} \|\mathfrak{I}_{0}\|_{q,s+1}^{\gamma,\mathcal{O}}\| \widehat y \|_{q,s_0 +1}^{\gamma,\mathcal{O}},
			\\
			\forall s\in[s_0,S],\quad& \| R_{01} \widehat w  \|_{q,s}^{\gamma,\mathcal{O}} \lesssim  \| Z \|_{q,s+ 1 }^{\gamma,\mathcal{O}}  \| \widehat w \|_{q,s_0 + 1 }^{\gamma,\mathcal{O}}
			+ \| Z \|_{q,s_0 + 1 }^{\gamma, \mathcal{O}} \|\mathfrak{I}_{0}\|_{q,s+1}^{\gamma,\mathcal{O}}\| \widehat w \|_{q,s_0 +1}^{\gamma,\mathcal{O}}
			\, . 
		\end{align*}
		\item The operators $\mathsf{K}_{20}$ and  $ \mathsf{K}_{11}$,  defined in \eqref{K20} and \eqref{K11}, satisfy the estimates
		\begin{align*}
			\forall s\in[s_0,S],\quad& \| \mathsf{K}_{20}  \|_{q,s}^{\gamma,\mathcal{O}} \lesssim \varepsilon \big( 1 + \| \mathfrak{I}_{0}\|_{q,s +3
			}^{\gamma,\mathcal{O}} \big) \, ,  \\ 
			\forall s\in[s_0,S],\quad& \| \mathsf{K}_{11} \widehat y \|_{q,s}^{\gamma,\mathcal{O}} 
			\lesssim \varepsilon \big(\| \widehat y \|_{q,s+3}^{\gamma,\mathcal{O}}
			+ \| \mathfrak{I}_{0} \|_{q,s +3}^{\gamma,\mathcal{O}}  
			\| \widehat y \|_{q,s_0+3}^{\gamma,\mathcal{O}} \big) \, , \\ 
			\forall s\in[s_0,S],\quad&  \| \mathsf{K}_{11}^\top \widehat w \|_{q,s}^{\gamma,\mathcal{O}}
			\lesssim \varepsilon \big(\| \widehat w \|_{q,s + 3}^{\gamma,\mathcal{O}}
			+  \| \mathfrak{I}_{0} \|_{q,s + 3}^{\gamma,\mathcal{O}}
			\| \widehat w \|_{q,s_0 + 3}^{\gamma,\mathcal{O}} \big)\, . 
		\end{align*}
	\item The matrices $\mathtt{A}$ and $\mathtt{B}$ defined in \eqref{devcalA} and \eqref{ttb} satisfy
	$$\forall s\in[s_0,S],\quad\|\mathtt{A}\|_{q,s}^{\gamma,\mathcal{O}}+\|\mathtt{B}\|_{q,s}^{\gamma,\mathcal{O}}\lesssim\|\mathfrak{I}_0\|_{q,s+1}^{\gamma,\mathcal{O}}.$$
	\end{enumerate}
\end{lem}
Notice that the matrix $\mathtt{A}(\varphi)$ measures the defect of the symplectic structure. In the following, we shall see  that it is of order  $O(Z)$. 
Notice that according to \eqref{devcalA} and \cite[Lem. 5]{BB15}, the coefficients $\mathtt{A}_{jk}$ of the matrix $\mathtt{A}$ can be written
\begin{align}\label{exp-Ajk}
	\mathtt{A}_{jk}(\varphi)&=\partial_{\varphi_k}I_0(\varphi)\cdot\partial_{\varphi_j}\vartheta_{0}(\varphi)-\partial_{\varphi_k}\vartheta_{0}(\varphi)\cdot\partial_{\varphi_j}I_0(\varphi)+\langle\partial_{\theta}^{-1}\partial_{\varphi_k}z_0(\varphi),\partial_{\varphi_j}z_{0}(\varphi)\rangle_{L^{2}(\mathbb{T})},
\end{align} 
and satisfy 
	\begin{align}\label{esp-ompartAjk}
		\omega\cdot\partial_\varphi \mathtt{A}_{jk}(\varphi)&=  {\mathcal W}\big( \partial_\varphi Z(\varphi) \underline{e}_k ,  \partial_\varphi i_0(\varphi)  \underline{e}_j \big) 
		+ 
		{\mathcal W} \big(\partial_\varphi i_0(\varphi) \underline{e}_k , \partial_\varphi Z(\varphi) \underline{e}_j \big) ,
	\end{align}
	where $\mathcal{W}$ is the symplectic form defined in \eqref{def symp-form} and  $(\underline{e}_1,\ldots,\underline{e}_d)$ denotes the canonical basis of $\mathbb{R}^d$.
In order to estimate $\mathtt{A}_{jk}(\varphi),$ we shall discuss the invertibility of  the operator $\omega\cdot\partial_{\varphi}.$  This task was accomplished in several paper \cite{BBMH18,BB15,BM18,HHM21}. and we shall outline here the main lines.\\
Let  $\gamma \in(0,1]$ and $\tau_{1}>0$ be defined as in \eqref{initial parameter condition}.   % $\lambda=(\alpha,\omega)\in\mathcal{O}\mapsto c_{\lambda} $ be   a smooth  function. 
We introduce  the Diophantine Cantor set  
$$
\mathtt {DC}(\gamma, \tau_1) := \bigcap_{l\in\mathbb{Z}^{d}\setminus\{0\}}\Big\{ \omega \in \mathbb{R}^d \quad \textnormal{s.t.} \quad 
|\omega \cdot l | >\frac{\gamma}{\langle l \rangle^{\tau_1}}\Big\}$$
and for $N\in\mathbb{N}^{*}$ we define the truncated Diophantine Cantor set
\begin{equation}\label{DC tau gamma}
	\mathtt {DC}_{N} (\gamma, \tau_1) :=\bigcap_{l\in\mathbb{Z}^{d}\setminus\{0\}\atop|l|\leqslant N} \Big\{ \omega \in \mathbb{R}^d \quad \textnormal{s.t.} \quad  
	|\omega \cdot l | > \frac{\gamma}{\langle l \rangle^{\tau_1}} \Big\} \, .
\end{equation}
Given $f:\mathcal{O}\times\T^{d+1}\to \mathbb{R}$ a smooth function with zero $\varphi$-average, that can be expanded in Fourier series as follows
$$
f=\sum_{(l,j)\in\mathbb{Z}^{d+1}\atop l\neq 0}f_{l,j}(\lambda,\omega) {\bf{e}}_{l,j},\quad {\bf{e}}_{l,j}(\varphi,\theta):=e^{i(l\cdot\varphi+j\theta)}.
$$
If  $\omega\in\mathtt {DC}(\gamma, \tau_1)$ then the equation $ \omega\cdot\partial_\varphi u=f$
has a periodic solution $u:\mathbb{T}^{d+1}\to\mathbb{R}$ given by
$$
u(\lambda,\varphi,\theta)=-\ii\sum_{(l,j)\in\mathbb{Z}^{d+1}\atop l\neq 0}\frac{f_{l,j}(\lambda)}{\omega\cdot l}{\bf{e}}_{l,j}(\varphi,\theta).
$$
For all $\omega\in\mathcal{O},$  we define the smooth extension of $u$ by
\begin{align}\label{Extend00}
	(\omega\cdot\partial_\varphi)_{\textnormal{ext}}^{-1}f:=-\ii\sum_{(l,j)\in\mathbb{Z}^{d+1}\atop l\neq 0}\frac{\chi\big(\gamma^{-1}\langle l\rangle^{\tau_1}\,\omega\cdot l \big)f_{l,j}(\lambda)}{\omega\cdot l}{\bf{e}}_{l,j}.
\end{align}
where $\chi\in\mathscr{C}^\infty(\mathbb{R},[0,1])$ is an even positive cut-off function  such that 
\begin{equation}\label{properties cut-off function first reduction} 
	\chi(\xi)=\left\{ \begin{array}{ll}
		0\quad \hbox{if}\quad |\xi|\leqslant\frac13&\\
		1\quad \hbox{if}\quad |\xi|\geqslant\frac12.
	\end{array}\right.
\end{equation}
Notice that this operator is well-defned in the whole set of parameters $\mathcal{O}$ and coincides with
the formal inverse of $(\omega\cdot\partial_\varphi)^{-1}$ when the frequency $\omega$ belongs to $ \mathtt {DC}(\gamma, \tau_1)$. The next result is the fundamental theorem of calculus in the quasi-periodic setting. It is proved in \cite[Lem. 2.5]{BBMH18} and \cite[Lem. 5.4]{HHM21}.
\begin{lem}\label{L-Invert0}
	Let $\gamma\in(0,1], q\in\mathbb{N}^{*}$. Then for any $ s\geqslant q$ we have
	$$
	\big\|(\omega\cdot\partial_\varphi)_{\textnormal{ext}}^{-1}f\big\|_{q,s}^{\gamma,\mathcal{O}}\lesssim\gamma^{-1}  \|f\|_{q,s+\tau_{1}q+\tau_{1}}^{\gamma,\mathcal{O}}.
	$$
	In addition, for any $N\in\mathbb{N}^{*}$ and for any $\omega\in \mathtt {DC}_{N} (\gamma, \tau_1)$ we have
	$$
	(\omega\cdot\partial_\varphi)(\omega\cdot\partial_\varphi)_{\textnormal{ext}}^{-1}\Pi_N=\Pi_N,
	$$
	where $\Pi_N$ is the orthogonal projection defined by
	$$
	\Pi_N\sum_{(l,j)\in\mathbb{Z}^{d+1}}f_{l,j} {\bf{e}}_{l,j}=\sum_{(l,j)\in\mathbb{Z}^{d+1}\atop |l|\leqslant N}f_{l,j} {\bf{e}}_{l,j}.
	$$
\end{lem}
For later purposes we need to fix some notation that will be adopted in the sequel. Take $N_0\geq 2$ and define  the sequence
\begin{equation}\label{definition of Nm}
				N_{-1}=1,\quad \forall n\in\mathbb{N},\quad N_{n}=N_{0}^{\left(\frac{3}{2}\right)^{n}}.
			\end{equation}
%\textcolor{red}{Recall that the Nash-Moser scheme finds a zero by constructing a sequence of approximate solutions living in finite dimensional subspaces obtained in our case by troncation up to order $N_n$ for the Fourier modes at the current step $n$ of the scheme. To make the Nash-Moser scheme convergent, we need to consider at each step, only a finite number of Diophantine conditions. We then need to split our operators in such a way they respect this troncation. The remainder terms being treated as prturbations.}
%Due to some  technical reasons related to the stabilization of Cantor sets during Nash-Moser scheme and the measure of the final Cantor set we need to work with finite number of non-resonant conditions imposed by a simple truncation in the time mode. %Therefore,  to perform the approximate inverse we need first to  decompose the operator $\omega\cdot\partial_{\varphi}$ into higher and lower frequencies and establish suitable estimates. The first result deals with  the following elementary lemma on the construction of a suitable splitting that respects the  reversibility.
Next, we shall  split the coefficients of the matrix $\mathtt{A}=\mathtt{A}(\varphi)$  defined in  \eqref{devcalA} as
\begin{equation}\label{decomp-Ajk}
	\mathtt{A}_{kj}=\mathtt{A}_{kj}^{(n)}+\mathtt{A}_{kj}^{(n),\perp}, \qquad \mathtt{A}_{kj}^{(n)}:= \Pi_{N_n}\mathtt{A}_{kj},\qquad \mathtt{A}_{kj}^{(n),\perp}:= \Pi_{N_n}^\perp\mathtt{A}_{kj}.
\end{equation}
%Then, by the
%smoothing properties \eqref{} and by \eqref{def:akj} we get
%\begin{equation}\label{stima A ij}
%%\| A_{k j}^{(n),\perp} \|_s^{q,\overline\gamma} \lesssim_s  \|  {\mathfrak I}_0 \|_{s+2}^{q,\overline\gamma},\qquad \,
%  \| \mathtt{A}_{k j}^{(n),\perp} \|_{s}^{q,\overline\gamma} \lesssim_{s, b} N_n^{-b}  \|  {\mathfrak I}_0 \|_{s+2+b}^{q,\overline\gamma}, \quad \forall b\geq 0.
%\end{equation}
The proof of the following lemma is quite similar to  Lemma 5.3. in \cite{BM20} with the a minor difference in the weighted norms. See also \cite[Lem. 6.3]{HHM21}.
\begin{lem}\label{lem:est-akj}
	Let $n\in\mathbb{N},$ then the following results  hold true.
	\begin{enumerate}[label=(\roman*)]
		\item The function  $\mathtt{A}_{kj}^{(n),\perp}$ satisfies  
		\begin{equation*}%\label{stima A ij}
			\forall b\geqslant 0,\quad\forall s\in[s_0,S],\quad \| \mathtt{A}_{k j}^{(n),\perp} \|_{q,s}^{\gamma,\mathcal{O}} \lesssim N_n^{-b}  \|  {\mathfrak I}_0 \|_{q,s+2+b}^{\gamma,\mathcal{O}}.
		\end{equation*}
		\item There exist functions $\mathtt{A}_{kj}^{(n),\textnormal{ext}}$ defined for any $(\lambda,\omega)\in \mathcal{O}$, $q$-times differentiable with respect to $\lambda$ and satisfying the estimate
		\begin{equation*}%\label{estim A ij}
			\forall s\in[s_0,S],\quad\|  \mathtt{A}_{k j}^{(n),\textnormal{ext}} \|_{q,s}^{\gamma,\mathcal{O}} \lesssim \gamma^{-1}
			\big(\| Z \|_{s+\tau_{1}q + \tau_{1}+1}^{\gamma,\mathcal{O}} + \| Z \|_{q,s_0+1}^{\gamma,\mathcal{O}} \|  {\mathfrak I}_0 \|_{q,s+\tau_{1}q + \tau_{1} +1}^{\gamma,\mathcal{O}} \big)\,.
		\end{equation*}
		Moreover,  $\mathtt{A}_{k j}^{(n),\textnormal{ext}}$ coincides with $ \mathtt{A}_{k j}^{(n)}$ in the Cantor set $\mathtt {DC}_{N_n} (\gamma, \tau_1) $.
	\end{enumerate}
\end{lem}
%\begin{proof}
%	\textbf{(i)} Follows from \eqref{decomp-Ajk}, \eqref{exp-Ajk} and Lemma \ref{Lem-lawprod}-(ii).\\
%	\textbf{(ii)} Applying $\Pi_{N_n}$ to \eqref{esp-ompartAjk}, one obtains
%	$$\omega\cdot\partial_{\varphi}\mathtt{A}_{jk}^{(n)}(\varphi)=\Pi_{N_n}\Big[{\mathcal W}\big( \partial_\varphi Z(\varphi) \underline{e}_k ,  \partial_\varphi i_0(\varphi)  \underline{e}_j \big) 
%	+ 
%	{\mathcal W} \big(\partial_\varphi i_0(\varphi) \underline{e}_k , \partial_\varphi Z(\varphi) \underline{e}_j \big)\Big].$$
%	Using Lemma \ref{Lem-lawprod}-(ii) and the continuity of the bilinear application $\mathcal{W}$, one gets
%	$$\|\omega\cdot\partial_{\varphi}\mathtt{A}_{jk}^{(n)}\|_{q,s}^{\gamma,\mathcal{O}}\lesssim\|Z\|_{q,s+1}^{\gamma,\mathcal{O}}+\|\mathfrak{I}_0\|_{q,s+1}^{\gamma,\mathcal{O}}\|Z\|_{q,s_{0}+1}^{\gamma,\mathcal{O}}.$$
%	We define the function $\mathtt{A}_{jk}^{(n),\textnormal{ext}}$ as 
%	$$\mathtt{A}_{jk}^{(n),\textnormal{ext}}(\varphi):=(\omega\cdot\partial_{\varphi})_{\textnormal{ext}}^{-1}\Pi_{N_n}\Big[{\mathcal W}\big( \partial_\varphi Z(\varphi) \underline{e}_k ,  \partial_\varphi i_0(\varphi)  \underline{e}_j \big) 
%	+ 
%	{\mathcal W} \big(\partial_\varphi i_0(\varphi) \underline{e}_k , \partial_\varphi Z(\varphi) \underline{e}_j \big)\Big].$$
%	Gathering Lemma \ref{L-Invert0} with the previous estimate gives the desired result.
%\end{proof}
\subsubsection{Construction of the approximate inverse}
This section is devoted to the  construction of  an approximate  right inverse of the operator $d_{i,\alpha}\mathcal{F}(i_{0},\alpha_{0})$ that will be discussed  in  Theorem \ref{proposition alomst approximate inverse}. One first may observe according to Proposition \ref{prop conj DGO}-(ii) and Lemmas \ref{lem:est-G} and \ref{lem:est-akj},  that the operator $\mathbb{E}$ vanishes at an exact solution up to fast decaying remainder terms. As a consequence, getting an approximate inverse for the full operator $d_{i,\alpha}\mathcal{F}(i_0,\alpha_0)$  amounts simply   to  invert the operator $\mathbb{D}$ up to small errors of type "Z" mixed with fast frequency decaying error. Let us  consider the triangular system given by
\begin{equation}\label{triangle-eq}
	\mathbb{D}[\widehat{\phi},\widehat{y},\widehat{w},\widehat{\alpha}]=\begin{pmatrix}
	g_{1}\\
	g_{2}\\
	g_{3}
\end{pmatrix},
\end{equation}
where $\mathbb{D}$ is defined in Proposition \ref{prop conj DGO}-(i). The system \eqref{triangle-eq} writes more explicitly in the following way
\begin{equation}\label{triangle-eq-1}
	\left\lbrace\begin{array}{l}
	\omega\cdot\partial_{\varphi}\widehat{\phi}=g_1+[\mathsf{K}_{20}(\varphi)\widehat{y}+\mathsf{K}_{11}^{\top}(\varphi)\widehat{w}+L_{1}^{\top}(\varphi)\widehat{\alpha}]\\
	\omega\cdot\partial_{\varphi}\widehat{y}=g_{2}-\mathtt{B}(\varphi)\widehat{\alpha}\\
	\big(\omega\cdot\partial_{\varphi}-\partial_{\theta}\mathsf{K}_{02}(\varphi)\big)\widehat{w}=g_3+\partial_{\theta}[\mathsf{K}_{11}(\varphi)\widehat{y}-L_{2}^{\top}(\varphi)\widehat{\alpha}].
\end{array}\right.
\end{equation}
The strategy to solve the above system in the variables $(\widehat{\phi},\widehat{y},\widehat{w})$ is first to solve the second action-component equation, then to solve the third normal-component equation and finally to solve the first angle-component equation. \\ Due to the fact that the Cantor set should be truncated then we need to solve approximately  the system \eqref{triangle-eq-1} and for this aim  we need the following statement proved in \cite[Lem. 6.4]{HHM21} and \cite{BBMH18}.

% almost inversion of the operator $\omega\cdot\partial_{\varphi}$ is given by the following lemma  which provides a suitable decomposition for $\omega\cdot\partial_{\varphi}$ respecting the finitely-many non-resonant conditions needed for the Nash-Moser scheme together
%
% with the preservation of the reversibility property.
\begin{lem}\label{lem:est-Dn}
	The following results hold true.
	\begin{enumerate}[label=(\roman*)]
		\item There exists a function $\mathtt{g}:\mathbb{Z}^d\setminus \{0\}\to \{-1,1\}$ such that  
		$$\forall l\in \mathbb{Z}^d\setminus \{0\},\quad \mathtt{g}(-l)=-\mathtt{g}(l).$$
		\item For all $(\lambda,\omega)\in\mathcal{O}$ the operator  $\omega\cdot \partial_\varphi$ can be split  as follows 
		$$\omega\cdot \partial_\varphi=\mathcal{D}_{(n)} +\mathcal{D}_{(n)}^{\perp},$$
		with 
		\begin{align*}
			\mathcal{D}_{(n)}&:= \omega\cdot \partial_\varphi\, \Pi_{N_n}+ \Pi_{N_n,\mathtt{g}}^\perp\\
			\mathcal{D}_{(n)}^{\perp}&:=   \omega\cdot \partial_\varphi\, \Pi_{N_n}^\perp- \Pi_{N_n, \mathtt{g}}^\perp\,, 
		\end{align*}
		where 
		$$
		\Pi_{N_n, \mathtt{g}}^\perp\sum_{(l,j)\in\mathbb{Z}^{d+1}}f_{l,j} {\bf{e}}_{l,j}=\sum_{(l,j)\in\mathbb{Z}^{d+1}\atop |l|> N_n} \mathtt{g}(l) f_{l,j} {\bf{e}}_{l,j}.
		$$
		\item The operator $\mathcal{D}_{(n)}^{\perp}$ satisfies
		$$
		\forall b\geqslant 0,\quad\forall s\in[s_0,S],\quad \|\mathcal{D}_{(n)}^{\perp}h\|_{q,s}^{\gamma,\mathcal{O}}\leqslant N_n^{-b}\|h\|_{q,s+b+1}^{\gamma,\mathcal{O}}.
		$$
		\item There exists a family of linear operators $\big([\mathcal{D}_{(n)}]_{\textnormal{ext}}^{-1} \big)_n$ satisfying, for any $h\in W^{q,\infty,\gamma}(\mathcal{O},H_0^s(\mathbb{T}^{d+1}))$,  
		$$
		\forall s\in[s_0,S],\quad \sup_{n\in\mathbb{N}} \|[\mathcal{D}_{(n)}]_{\textnormal{ext}}^{-1} h\|_{q,s}^{\gamma,\mathcal{O}}\lesssim \gamma^{-1}\|h\|_{q,s+\tau_{1}q+\tau_{1}}^{\gamma,\mathcal{O}}.
		$$
		Moreover, for all $\omega\in  \mathtt {DC}_{N_n} (\gamma, \tau_1) $ one has the identity
		\begin{equation}\label{DnDnext-1_Id}
		\mathcal{D}_{(n)}[\mathcal{D}_{(n)}]_{\textnormal{ext}}^{-1} =\textnormal{Id}.
	\end{equation}
	\end{enumerate}
\end{lem}
%
%\begin{proof}
%	\textbf{(i)} For all $l=(l_1,...,l_d)\in \mathbb{Z}^d\setminus \{0\}$ we define $\mathtt{g}(l)$ as the sign of the first non-zero component in the vector $l$. This defines a function enjoying the expected properties. \\
%	\textbf{(ii)} Obvious.\\
%	\textbf{(iii)} Follows immediately from Lemma \ref{Lem-lawprod}-(ii) and the fact that $|\mathtt{g}|=1$.\\
%	\textbf{(iv)} Follows from ${\rm (ii)}$, Lemma \ref{Lem-lawprod}-(ii) and Lemma \ref{L-Invert0}.
%\end{proof}
Consider    the linearized operator restricted to  the normal directions $\widehat{\mathcal{L}}_{\omega}$ and defined by
\begin{equation}\label{Lomega def}
	\widehat{\mathcal{L}}_{\omega} := \Pi_{\mathbb{S}_0}^\bot \big(\omega\cdot \partial_\varphi   - 
	\partial_\theta  \mathsf{K}_{02}(\varphi) \big)\Pi_{\mathbb{S}_0}^{\perp},
\end{equation}
which appears in the last equation of \eqref{triangle-eq-1}.
 The construction of an approximate right inverse of this operator is the heart part of this paper and will be discussed in  Proposition
\ref{inversion of the linearized operator in the normal directions}.  Here  we give only a partial  statement.\begin{prop}\label{thm:inversion of the linearized operator in the normal directions}
	Let $(\gamma,q,d,\tau_{1},\tau_2,s_{0},s_{h},\mu_{2},S)$ satisfy  \eqref{initial parameter condition} \eqref{setting tau1 and tau2} and \eqref{Conv-T2}. There exist $\varepsilon_0>0$ and  $\sigma=\sigma(\tau_1,\tau_2,q,d)>0$ such that if
	\begin{equation}\label{small-ilo}
		\varepsilon\gamma^{-2-q}N_0^{\mu_2}\leqslant\varepsilon_0\quad\textnormal{and}\quad\|\mathfrak{I}_0\|_{q,s_h+\sigma}^{\gamma,\mathcal{O}}\leqslant 1,
	\end{equation} 
then there exist a family of linear  operator $(\mathtt{T}_{\omega,n})_n$  satisfying 	\begin{equation}\label{estimate mathcalTomega0}
		\forall \, s\in\,[ s_0, S] ,\quad\sup_{n\in\mathbb{N}}\|\mathtt{T}_{\omega,n}h\|_{q,s}^{\gamma,\mathcal{O} }\lesssim\gamma^{-1}\left(\|h\|_{q,s+{\sigma}}^{\gamma,\mathcal{O}}+\| \mathfrak{I}_{0}\|_{q,s+{\sigma}}^{\gamma,\mathcal{O} }\|h\|_{q,s_{0}+{\sigma}}^{\gamma,\mathcal{O}}\right)
	\end{equation}
	and   a  family of  Cantor sets $\{\mathcal{G}_{n}=\mathcal{G}_{n}(\gamma,\tau_{1},\tau_{2},i_{0})\}_n$, satisfying the inclusion 
	$$\mathcal{G}_{n}\subset(\lambda_0,\lambda_1)\times\mathtt{DC}_{N_n}(\gamma,\tau_{1})$$ 
	such that in each  set $\mathcal{G}_{n}$
	we have the splitting
	$$
	\widehat{\mathcal{L}}_{\omega}=\widehat{\mathtt{L}}_{\omega,n}+\widehat{\mathtt{R}}_n, %\quad\hbox{with}\quad \widehat{\mathtt{L}}_{\omega,n}\mathtt{T}_{\omega,n}=\textnormal{Id} %\quad\hbox{and}\quad \widehat{\mathtt{R}}_n={\mathtt{E}_n}\widehat{\mathtt{L}}_{\omega}
	$$
	with
	\begin{equation}\label{Lnhn_Id}
	\widehat{\mathtt{L}}_{\omega,n}\mathtt{T}_{\omega,n}=\textnormal{Id},
\end{equation}
	where the operators $\widehat{\mathtt{L}}_{\omega,n}$ and $\widehat{\mathtt{R}}_n$ are defined in the whole set $\mathcal{O}$ with the estimates
	%and
	\begin{align*}
		\forall\, s\in [s_0,S],\quad &\|\widehat{\mathtt{L}}_{\omega,n} \rho\|_{q,s}^{\gamma,\mathcal{O}}\lesssim \|\rho\|_{q,s+1}^{\gamma,\mathcal{O}}+\varepsilon\gamma^{-2}\|\mathfrak{I}_{0}\|_{q,s+\sigma}^{\gamma,\mathcal{O}}\|\rho\|_{q,s_{0}+1}^{\gamma,\mathcal{O}},\\
		%\forall\, s\in [s_0,S],\quad &\sup_{n\in\mathbb{N}}\|\widehat{\mathtt{R}}_n\rho\|_{q,s}^{\gamma,\mathcal{O} } \lesssim \overline\gamma^{-1}\|\rho\|_{q,s+\sigma}^{\gamma,\mathcal{O} }+{\varepsilon\gamma^{-3}}\| \mathfrak{I}_{0}\|_{q,s+\sigma}^{\gamma,\mathcal{O} }\|\rho\|_{q,s_{0}+\sigma}^{\gamma,\mathcal{O}},\\
		\forall\, b\in [0,S],\quad & \|\widehat{\mathtt{R}}_n\rho\|_{q,s_0}^{\gamma,\mathcal{O}}
		\lesssim N_n^{-b}\gamma^{-1}\Big( \|\rho\|_{q,s_0+b+\sigma}^{\gamma,\mathcal{O}}+{\varepsilon\gamma^{-2}}\| \mathfrak{I}_{0}\|_{q,s_0+b+\sigma}^{\gamma,\mathcal{O}}\|\rho\|_{q,s_0+\sigma}^{\gamma,\mathcal{O}} \Big)\\
		&\qquad\qquad\qquad+ \varepsilon\overline\gamma^{-3}N_{0}^{\mu_{2}}{N_{n+1}^{-\mu_{2}}} \|\rho\|_{q,s_0+\sigma}^{\gamma,\mathcal{O}}.
	\end{align*}
\end{prop}
For the splitting below which follows from  the foregoing results we refer to  (6.45) in \cite{HHM21}. Consider the linear operator $\mathbb{L}_{\textnormal{ext}}$ defined by 
\begin{equation}\label{def:lext}
	\mathbb{L}_{\textnormal{ext}}= \mathbb{D}_n+{\mathbb E}_{n}^{\textnormal{ext}}+{\mathscr P}_{n}+{\mathscr Q}_{n},
\end{equation}
where,  for any $(\widehat \phi, \widehat y, \widehat w, \widehat \alpha)\in  \mathbb{T}^d \times \mathbb{R}^d \times H_{\perp}^s \times \mathbb{R}^d $
\begin{align}
	{\mathbb D}_n [\widehat \phi, \widehat y, \widehat w, \widehat \alpha ] & :=
	\left(
	\begin{array}{c}
		\mathcal{D}_{(n)}  \widehat\phi-\mathsf{K}_{20}(\varphi)\widehat y-\mathsf{K}_{11}^\top(\varphi) \widehat w -L_1^{\top}(\varphi)\widehat \alpha
		\\
		\mathcal{D}_{(n)} \widehat y+\mathtt{B}(\varphi) \widehat \alpha \\
		\widehat{\mathtt{L}}_{\omega,n} \widehat w -\partial_\theta\big[\mathsf{K}_{11}(\varphi) \widehat y -L_2^{\top}(\varphi) \widehat\alpha   \big] 
	\end{array}
	\right), \label{def:Dn}
	\\
	\nonumber{\mathbb E}_{n}^{\textnormal{ext}} [\widehat \phi, \widehat y, \widehat w, \widehat\alpha ] & :=
	[D G_0({\mathtt u}_0(\varphi))]^{-1}    \partial_\varphi Z(\varphi) \widehat\phi - \left(
	\begin{array}{c}
		0 
		\\
		R_{10}(\varphi)\widehat y+R_{01}(\varphi)\widehat w
		\\
		0   
	\end{array}
	\right) \\ &\quad +\quad \left(
	\begin{array}{c}
		0 
		\\
		\mathtt{A}^{(n),\textnormal{ext}}(\varphi)\big[\mathsf{K}_{20}(\varphi) \widehat y+\mathsf{K}_{11}^\top(\varphi) \widehat w\big]
		\\
		0   
	\end{array}
	\right) , \label{def:Rn}
	\\
	{\mathscr P}_{n}  [\widehat \phi, \widehat y, \widehat w, \widehat\alpha ] & :=
	\left(
	\begin{array}{c}
		\mathcal{D}_{(n)}^{\perp}  \widehat\phi
		\\
		\mathcal{D}_{(n)}^{\perp} \widehat y + \mathtt{A}^{(n),\perp}(\varphi)\big[\mathsf{K}_{20}(\varphi) \widehat y+\mathsf{K}_{11}^\top(\varphi) \widehat w\big]
		\\
		0
	\end{array}
	\right), \label{def:Pn}
	\\
	{\mathscr Q}_{n}  [\widehat \phi, \widehat y, \widehat w, \widehat\alpha ] & :=
	\left(
	\begin{array}{c}
		0
		\\
		0
		\\
		\widehat{\mathtt{R}}_n[ \widehat  w]
	\end{array}
	\right). \label{def:Qn}
\end{align}
Then,  the operator $\mathbb{L}_{\textnormal{ext}}$ is defined on the whole set $ \mathcal{O} $ and when it is  restricted to the Cantor set $\mathcal{G}_{n}$ it coincides with the conjugated linearized operator obtained in \eqref{bbL}, that is,
\begin{equation}\label{lext-f}
	\mathbb{L}_{\textnormal{ext}}=[D G_0({\mathtt u}_0)]^{-1}  d_{i,\alpha} {\mathcal F} (i_0,\alpha_0) D \widetilde{G}_0({\mathtt u}_0) \quad \textnormal{in }\mathcal{G}_n.
\end{equation}
In the next result, we give some useful  estimates for the different terms appearing in $\mathbb{L}_{\textnormal{ext}}$ needed to obtain good tame estimates for the approximate inverse. 
\begin{prop}\label{prop:decomp-lin}
	Let $(\gamma,q,d,\tau_{1},s_{0},\mu_{2})$ satisfy \eqref{initial parameter condition} and \eqref{Conv-T2} and assume the  conditions \eqref{small-C3} and  \eqref{small-ilo}. Then, denoting $\widehat{\mathtt v}=(\widehat \phi, \widehat y, \widehat w, \widehat\alpha)$, the following assertions hold true.
	\begin{enumerate}[label=(\roman*)]
		
		\item The operator ${\mathbb E}_n^{\textnormal{ext}}$ satisfies the estimate 
		\begin{align*} 
			\| {\mathbb E}_n^{\textnormal{ext}} [\widehat{\mathtt v}] \|_{q,s_{0}}^{\gamma,\mathcal{O}}
			&\lesssim \| Z \|_{q,s_0 + \sigma }^{\gamma,\mathcal{O}}  \| \widehat{\mathtt v} \|_{q,s_0 + \sigma }^{\gamma,\mathcal{O}}.
		\end{align*}
		
		\item The operator $ {\mathscr P}^{(n)} $ satisfies the estimate
		\begin{align*} 
			\forall b\geqslant 0,\quad  
			\|  {\mathscr P}_{n} [\widehat{\mathtt v}] \|_{q,s_0}^{\gamma,\mathcal{O}} &\lesssim N_n^{-b}  \big( \| \widehat{\mathtt v} \|_{q,s_0 + \sigma + b }^{\gamma,\mathcal{O}}+\varepsilon
			\| \mathfrak{I}_{0} \|_{q,s_0+  \sigma  +b    }^{\gamma,\mathcal{O}} \big \| \widehat{\mathtt v} \|_{q,s_0 +\sigma}^{\gamma,\mathcal{O}} \big).
		\end{align*}
		\item The operator $ {\mathscr Q}_{n} $ satisfies the estimate
	\begin{align*}
		\forall\, b\in [0,S],\quad & \|\mathscr{Q}_{n}\widehat{\mathtt{v}}\|_{q,s_0}^{\gamma,\mathcal{O}}
		\lesssim N_n^{-b}\gamma^{-1}\Big( \|\widehat{w}\|_{q,s_0+b+\sigma}^{\gamma,\mathcal{O}}+{\varepsilon\gamma^{-2}}\| \mathfrak{I}_{0}\|_{q,s_0+b+\sigma}^{\gamma,\mathcal{O}}\|\widehat{w}\|_{q,s_0+\sigma}^{\gamma,\mathcal{O}} \Big)\\
		&\qquad\qquad\qquad+ \varepsilon\overline\gamma^{-3}N_{0}^{\mu_{2}}{N_{n+1}^{-\mu_{2}}} \|\widehat{w}\|_{q,s_0+\sigma}^{\gamma,\mathcal{O}}.
	\end{align*}
		
		\item There exists a family of operators $\big([{\mathbb D}_n]_{\textnormal{ext}}^{-1}\big)_n$  such that
		for all $ g := (g_1, g_2, g_3) $ 
		satisfying  the reversibility property 
		\begin{equation*}
			g_1(\varphi) = g_1(- \varphi)\,,\quad g_2(\varphi) = - g_2(- \varphi)\,,\quad g_3(\varphi) = - ({\mathscr S} g_3)(- \varphi)\,,
		\end{equation*}
		the function 
		$ [{\mathbb D}_n]_{\textnormal{ext}}^{-1} g $ 
		satisfies the estimate
		\begin{equation*} 
			\forall\, s\in [s_0,S],\quad \| [{\mathbb D}_n]_{\textnormal{ext}}^{-1}g \|_{q,s}^{\gamma,\mathcal{O}}
			\lesssim \gamma^{-1} \big( \| g \|_{q,s + \sigma }^{\gamma,\mathcal{O}}
			+  \| {\mathfrak I}_0  \|_{q,s + \sigma}^{\gamma,\mathcal{O}}
			\| g \|_{q,s_0 + \sigma}^{\gamma,\mathcal{O}}  \big)
		\end{equation*}
		and for all $(\lambda,\omega)\in \mathcal{G}_{n}$ one has
		$$
		{\mathbb D}_n \,[{\mathbb D}_n]_{\textnormal{ext}}^{-1} =\textnormal{Id}.
		$$
	\end{enumerate}
\end{prop}
\begin{proof}
	{\bf (i)} The estimate of ${\mathbb E}_n^{\textnormal{ext}}$ is obtained from \eqref{def:Rn}, Lemma \ref{lem:est-G}, Lemma \ref{Lem-lawprod}-(iv) and  \mbox{Lemma \ref{lem:est-akj}-{\rm (ii)}.}\\
	{\bf (ii)} From  \eqref{def:Pn}, Lemma \ref{lem:est-Dn}-{\rm (iii)},   Lemma \ref{Lem-lawprod}-(iv),  Lemma \ref{lem:est-akj}-{\rm (i)}, Lemma \ref{lem:est-G}-{\rm (ii)} we obtain the estimate on  ${\mathscr P}_{n}$.\\
	{\bf (iii)} It is a consequence of  \eqref{def:Qn} and Proposition \ref{thm:inversion of the linearized operator in the normal directions}.\\
	{\bf (iv)} The proof can be found in \cite[Prop. 6.3]{HHM21} and for the sake of completeness we shall sketch the main ideas. We intend to look for  an exact inverse of ${\mathbb D}_n$ by solving  the system 
	\begin{equation}\label{operatore inverso approssimato proiettato}
		{\mathbb D}_n [\widehat \phi, \widehat y, \widehat w, \widehat \alpha]  
		= \begin{pmatrix}
			g_1  \\
			g_2  \\
			g_3 
		\end{pmatrix},
	\end{equation}
	where $(g_1, g_2, g_3)$ satisfy the reversibility property 
	\begin{equation}\label{parita g1 g2 g3}
		g_1(\varphi) = g_1(- \varphi)\,,\quad g_2(\varphi) = - g_2(- \varphi)\,,\quad g_3(\varphi) = - ({\mathscr S} g_3)(- \varphi),
	\end{equation}
with   $\mathscr{S}$ being the involution defined in \eqref{definition of the involution mathcal S}. Note that in view of \eqref{def:Dn}, the system \eqref{operatore inverso approssimato proiettato} writes
\begin{equation}\label{triangle-eq-n}
	\left\lbrace\begin{array}{l}
		\mathcal{D}_{(n)}\widehat{\phi}=g_1+[\mathsf{K}_{20}(\varphi)\widehat{y}+\mathsf{K}_{11}^{\top}(\varphi)\widehat{w}+L_{1}^{\top}(\varphi)\widehat{\alpha}]\\
		\mathcal{D}_{(n)}\widehat{y}=g_{2}-\mathtt{B}(\varphi)\widehat{\alpha}\\
		\widehat{\mathtt{L}}_{\omega,n}\widehat{w}=g_3+\partial_{\theta}[\mathsf{K}_{11}(\varphi)\widehat{y}-L_{2}^{\top}(\varphi)\widehat{\alpha}].
	\end{array}\right.
\end{equation}
	We first consider the second action-component equation in \eqref{triangle-eq-n}, namely 
	$$ \mathcal{D}_{(n)} \widehat y  =
	g_2  -  \mathtt{B}(\varphi)\widehat \alpha. $$ 
	In view of \eqref{parita g1 g2 g3}, \eqref{ttb} and \eqref{exp-Ajk}, $g_2$ and $\mathtt{B}$ are odd in the variable $\varphi.$ Thus, the $\varphi$-average of the right hand side of this equation is zero. Then, by Lemma \ref{lem:est-Dn}-(iv) its solution in the Cantor \mbox{set $\mathtt{DC}_{N_n}(\gamma,\tau_1)$ is given by} 
	\begin{equation}\label{soleta}
		\widehat y := [\mathcal{D}_{(n)}]_{\textnormal{ext}}^{-1} \big(
		g_2  -  \mathtt{B}(\varphi)\widehat \alpha \big)\,.  
	\end{equation}
	Then we turn to the third normal-component equation in \eqref{triangle-eq-n}, namely
	$$\widehat{\mathtt{L}}_{\omega,n} \widehat w = g_3 + \partial_\theta[ \mathsf{K}_{11}(\varphi) \widehat y -L_2^\top(\varphi)  \widehat \alpha].$$
	By Proposition \ref{thm:inversion of the linearized operator in the normal directions}, this equation admits as a solution 
	\begin{equation}\label{normalw}
		\widehat w := \mathtt{T}_{\omega,n} \big( g_3 + \partial_\theta [\mathsf{K}_{11}(\varphi) \widehat y - L_2^\top(\varphi)  \widehat \alpha]\big) \, .  
	\end{equation}
	Finally, we solve the first angle-equation in \eqref{triangle-eq-n}, 
	which, substituting \eqref{soleta}, \eqref{normalw}, becomes
	\begin{equation}\label{equazione psi hat}
		\mathcal{D}_{(n)} \widehat \phi  = 
		g_1 +  M_1(\varphi)\widehat\alpha + M_2(\varphi) g_2 + M_3(\varphi) g_3\,,
	\end{equation}
	where
	\begin{align}\label{M1}
		M_1(\varphi) &:=L_1^\top(\varphi)- M_2(\varphi) \mathtt{B}(\varphi) - M_3(\varphi)  \partial_{\theta}L_{2}^{\top}(\varphi) \,, \\
		\label{cal M2}
		M_2(\varphi) &:=  \mathsf{K}_{20}(\varphi) [\mathcal{D}_{(n)}]_{\textnormal{ext}}^{-1}  + \mathsf{K}_{11}^\top(\varphi)\mathtt{T}_{\omega,n} \partial_\vartheta \mathsf{K}_{11}(\varphi)[\mathcal{D}_{(n)}]_{\textnormal{ext}}^{-1}  \, , \\
		M_3(\varphi) &:= \mathsf{K}_{11}^\top (\varphi) \mathtt{T}_{\omega,n}  \, .    \label{cal M3}
	\end{align}
	%%%%%%%%%%%%%%%%%%%%%%%
	To solve the equation \eqref{equazione psi hat} we choose $ \widehat \alpha $ such that the right hand side  has zero $\varphi$-average. Notice that Lemma \ref{lem:est-G}, \eqref{small-C3}, \eqref{estimate mathcalTomega0} and Lemma \ref{lem:est-Dn}-(ii) imply
	\begin{equation}\label{est:M2-M3}
		\forall s\in[s_0,S],\quad \|M_{2}g_{2}\|_{q,s}^{\gamma,\mathcal{O}}+\|M_{3}g_{3}\|_{q,s}^{\gamma,\mathcal{O}}\lesssim \varepsilon\left(\|g\|_{q,s+\sigma}^{\gamma,\mathcal{O}}+\|\mathfrak{I}_0\|_{q,s+\sigma}^{\gamma,\mathcal{O}}\|g\|_{q,s_0+\sigma}^{\gamma,\mathcal{O}}\right).
		\end{equation}
	By Lemma \ref{lem:est-G}-{\rm (iii)}, \eqref{small-C3}, the $\phi$-averaged matrix is
	$ \langle M_1 \rangle = {\rm Id} + O( \varepsilon \gamma^{-1 }) $.  
	Therefore, for $ \varepsilon \gamma^{- 1} $ small enough,  
	$ \langle M_1 \rangle$ is invertible and $\langle M_1 \rangle^{-1} = {\rm Id} 
	+ O(\varepsilon\gamma^{- 1})$. We thus define 
	\begin{equation}\label{sol alpha}
		\widehat\alpha  := - \langle M_1 \rangle^{-1} 
		( \langle g_1 \rangle + \langle M_2 g_2 \rangle + \langle M_3 g_3 \rangle ) \, .
	\end{equation}
Remark that $\widehat{\alpha}$ satisfies
\begin{equation}\label{est:alphah}
	\|\widehat{\alpha}\|_{q}^{\gamma,\mathcal{O}}\lesssim\|g\|_{q,s_0+\sigma}^{\gamma,\mathcal{O}}.
\end{equation}
Coming back to \eqref{soleta} and using \eqref{est:alphah}, \eqref{small-C3} together with Lemma \ref{lem:est-Dn}-(iv) and Lemma \ref{lem:est-G}-(iv), we obtain
\begin{equation}\label{est:yh}
	\forall s\in[s_0,S],\quad \|\widehat{y}\|_{q,s}^{\gamma,\mathcal{O}}\lesssim\gamma^{-1}\left(\|g\|_{q,s+\sigma}^{\gamma,\mathcal{O}}+\|\mathfrak{I}_0\|_{q,s+\sigma}^{\gamma,\mathcal{O}}\|g\|_{q,s_0+\sigma}^{\gamma,\mathcal{O}}\right).
\end{equation}
Putting together \eqref{normalw}, \eqref{estimate mathcalTomega0}, Lemma \ref{lem:est-G}-(iii), \eqref{est:alphah}, \eqref{est:yh} and \eqref{small-C3}, one should get, up to redefine the value of $\sigma$,
\begin{equation}\label{est:wh}
	\forall s\in[s_0,S],\quad \|\widehat{w}\|_{q,s}^{\gamma,\mathcal{O}}\lesssim\gamma^{-1}\left(\|g\|_{q,s+\sigma}^{\gamma,\mathcal{O}}+\|\mathfrak{I}_0\|_{q,s+\sigma}^{\gamma,\mathcal{O}}\|g\|_{q,s_0+\sigma}^{\gamma,\mathcal{O}}\right).
\end{equation}
	With the choice \eqref{sol alpha} of $ \widehat \alpha $, the
	equation \eqref{equazione psi hat} admits as a solution
	\begin{equation}\label{sol psi}
		\widehat \phi :=
		[\mathcal{D}_{(n)}]_{\textnormal{ext}}^{-1} \big( g_1 + M_1(\varphi)\widehat\alpha+ M_2(\varphi) g_2 + M_3(\varphi) g_3 \big) \, . 
	\end{equation}
Putting together \eqref{sol psi}, Lemma \ref{lem:est-Dn}-(ii), \eqref{est:alphah} and \eqref{est:M2-M3}, one obtains
\begin{equation}\label{est:phih}
	\forall s\in[s_0,S],\quad \|\widehat{\phi}\|_{q,s}^{\gamma,\mathcal{O}}\lesssim\gamma^{-1}\left(\|g\|_{q,s+\sigma}^{\gamma,\mathcal{O}}+\|\mathfrak{I}_0\|_{q,s+\sigma}^{\gamma,\mathcal{O}}\|g\|_{q,s_0+\sigma}^{\gamma,\mathcal{O}}\right).
\end{equation}
	In conclusion, we have obtained
	a solution  $(\widehat \phi, \widehat y, \widehat w, \widehat\alpha):=[\mathbb{D}_n]_{\textnormal{ext}}^{-1}g$ of the linear system \eqref{operatore inverso approssimato proiettato} satisfying in virtue of \eqref{est:alphah}, \eqref{est:phih}, \eqref{est:wh} and \eqref{est:yh},
	$$\forall s\in[s_0,S],\quad \|[\mathbb{D}_n]_{\textnormal{ext}}^{-1}g\|_{q,s}^{\gamma,\mathcal{O}}\lesssim\gamma^{-1}\left(\|g\|_{q,s+\sigma}^{\gamma,\mathcal{O}}+\|\mathfrak{I}_0\|_{q,s+\sigma}^{\gamma,\mathcal{O}}\|g\|_{q,s_0+\sigma}^{\gamma,\mathcal{O}}\right).$$ 
	Notice that the relation 
	$$\mathbb{D}_n[\mathbb{D}_n]_{\textnormal{ext}}^{-1}=\textnormal{Id}\quad\textnormal{in }\mathcal{G}_n$$
	is a direct consequence of \eqref{DnDnext-1_Id} and \eqref{Lnhn_Id}.
\end{proof}
The last point is to  prove that the operator 
\begin{equation}\label{definizione T} 
	{\rm T}_0 := {\rm T}_0(i_0) := (D { \widetilde G}_0)(\mathtt{u}_0) \circ [{\mathbb D}_n]_{\textnormal{ext}}^{-1}\circ (D G_0) (\mathtt{u}_0)^{-1}
\end{equation}
is an approximate right  inverse for $d_{i,\alpha} {\mathcal F}(i_0 ,\alpha_{0}).$  
\begin{theo}[Approximate inverse]  \label{proposition alomst approximate inverse}
	Let $(\gamma,q,d,\tau_{1},\tau_2,s_{0},s_h,\mu_{2},S)$ satisfy  \eqref{initial parameter condition}, \eqref{setting tau1 and tau2}, \eqref{Conv-T2} and \eqref{param}. 
	Then there exists $ \overline{\sigma}= \overline{\sigma}(\tau_1,\tau_2,d,q)>0$ such that if the smallness conditions \eqref{small-C3} and \eqref{small-ilo} hold, then  
	the operator $ {\rm T}_0  $ defined in \eqref{definizione T} is reversible and satisfies for all $ g = (g_1, g_2, g_3) $, \mbox{with  \eqref{parita g1 g2 g3},}
	\begin{equation}\label{estimate on T0}
		\forall s\in [s_0,S],\quad \| {\rm T}_0 g\|_{q,s}^{\gamma,\mathcal{O}}\lesssim\gamma^{-1}\left(\|g\|_{q,s+\overline{\sigma}}^{\gamma,\mathcal{O}}+\|\mathfrak{I}_{0}\|_{q,s+\overline{\sigma}}^{\gamma,\mathcal{O}}\|g\|_{q,s_{0}+\overline{\sigma}}^{\gamma,\mathcal{O}}\right).
	\end{equation}
	Moreover  ${\rm T}_0$ is an almost-approximate  
	right 
	inverse of $d_{i, \alpha} 
	\mathcal{ F}(i_0,\alpha_0)$ in the Cantor set $ \mathcal{G}_{n}$. More precisely,   for all $ (\lambda,\omega)\in \mathcal{G}_{n}$ one has
	\begin{equation}\label{splitting per approximate inverse}
		d_{i , \alpha} \mathcal{ F} (i_0) \circ {\rm T}_0
		- {\rm Id} = \mathcal{E}^{(n)}_1+\mathcal{E}^{(n)}_2+\mathcal{E}^{(n)}_3,  
	\end{equation}
	where the operators $\mathcal{E}^{(n)}_1$, $\mathcal{E}^{(n)}_2$ and $\mathcal{E}^{(n)}_3$ are defined in the  set $\mathcal{O}$ with the estimates
	\begin{align}
		\|\mathcal{E}^{(n)}_1 g \|_{q,s_0}^{\gamma,\mathcal{O}} & \lesssim \gamma^{-1 } \| \mathcal{ F}(i_0, \alpha_0) \|_{q,s_0 +\overline{\sigma}}^{\gamma,\mathcal{O}} \| g \|_{q,s_0 + \overline{\sigma}}^{\gamma,\mathcal{O}},\label{eaai-E1}\\
		\forall b \geqslant 0,\quad \| \mathcal{E}^{(n)}_2 g \|_{q,s_0}^{\gamma,\mathcal{O}}& \lesssim 
		\gamma^{- 1} N_n^{- b } \big( \| g \|_{q,s_0+b+\overline{\sigma} }^{\gamma,\mathcal{O}}+\varepsilon
		\| \mathfrak{I}_{0} \|_{q,s_0
			+b+\overline{\sigma}    }^{\gamma,\mathcal{O}} \big \| g \|_{q,s_0 +\overline{\sigma}}^{\gamma,\mathcal{O}} \big)\,,
		\label{eaai-E2}\\
		\forall b\in [0,S], \quad \| \mathcal{E}^{(n)}_3 g \|_{q,s_0}^{\gamma,\mathcal{O}}& \lesssim N_n^{-b}\gamma^{-2}\Big( \|g\|_{q,s_0+b+\overline{\sigma}}^{\gamma,\mathcal{O}}+{\varepsilon\gamma^{-2}}\| \mathfrak{I}_{0}\|_{q,s_0+b+\overline{\sigma}}^{\gamma,\mathcal{O}}\|g\|_{q,s_0+\overline{\sigma}}^{\gamma,\mathcal{O}} \Big)\nonumber\\
		&\quad+ \varepsilon\gamma^{-4}N_{0}^{{\mu}_{2}}{N_{n}^{-\mu_{2}}} \|g\|_{q,s_0+\overline{\sigma}}^{\gamma,\mathcal{O}}. \label{eaai-E3} 
	\end{align}
\end{theo}

\begin{proof}
	The estimate \eqref{estimate on T0} is a consequence of \eqref{definizione T}, Proposition \ref{prop:decomp-lin}-{\rm (iv)} and Lemma \ref{lem:est-G}-{\rm (i)}. Then, according to  \eqref{def:lext} and 
	\eqref{lext-f}, in the Cantor set $ \mathcal{G}_{n} $ we have the decomposition
	\begin{equation*}
		%\label{dimuf} 
		\begin{aligned}
			d_{i,\alpha}\mathcal{F}(i_{0},\alpha_{0})&=DG_{0}({\mathtt u}_0)\circ \mathbb{L}_{\textnormal{ext}} \circ D[\widetilde{G}_{0}({\mathtt u}_0)]^{-1}
			\\
			&=DG_{0}({\mathtt u}_0)\circ {\mathbb{D}}_n\circ D[\widetilde{G}_{0}({\mathtt u}_0)]^{-1}+ DG_{0}({\mathtt u}_0) \circ \mathbb{E}_{n}^{\textnormal{ext}} \circ[\widetilde{G}_{0}({\mathtt u}_0)]^{-1}\\ &\quad+DG_{0}({\mathtt u}_0) \circ {\mathscr P}_{n}\circ[\widetilde{G}_{0}({\mathtt u}_0)]^{-1}+DG_{0}({\mathtt u}_0) \circ {\mathscr Q}_{n}\circ[\widetilde{G}_{0}({\mathtt u}_0)]^{-1}.
		\end{aligned}
	\end{equation*}
	By applying ${\rm T}_0$, defined in \eqref{definizione T},  to the last identity   we get for all $(\lambda,\omega)\in  \mathcal{G}_{n}$
	$$d_{i,\alpha}\mathcal{F}(i_{0},\alpha_{0})\circ {\rm T}_{0}-\textnormal{Id}=\mathcal{E}^{(n)}_1+\mathcal{E}^{(n)}_2+\mathcal{E}^{(n)}_3,
	$$
	with
	\begin{align*}
		&\mathcal{E}^{(n)}_1:= DG_{0}({\mathtt u}_0) \circ \mathbb{E}_{n}^{\textnormal{ext}} \circ[\widetilde{G}_{0}({\mathtt u}_0)]^{-1}\circ{\rm T}_{0},
		\\
		&\mathcal{E}^{(n)}_2:= DG_{0}({\mathtt u}_0) \circ {\mathscr P}_{n} \circ[\widetilde{G}_{0}({\mathtt u}_0)]^{-1}\circ{\rm T}_{0},
		\\
		&\mathcal{E}^{(n)}_3:= DG_{0}({\mathtt u}_0) \circ {\mathscr Q}_{n} \circ[\widetilde{G}_{0}({\mathtt u}_0)]^{-1}\circ{\rm T}_{0}.
	\end{align*}
	The estimates on $\mathcal{E}^{(n)}_1$, $\mathcal{E}^{(n)}_2$ and $\mathcal{E}^{(n)}_3$  come  from \eqref{estimate on T0}, Proposition \ref{prop:decomp-lin} and Lemma \ref{lem:est-G}-{\rm (i)}. 
	
\end{proof}
%\textcolor{red}{stop here}
	\section{Reduction of the linearized operator in the normal directions}\label{sec red lin op}
			\indent In this section, we fix a torus $i_{0}=(\vartheta_{0},I_{0},z_{0})$ close to the flat one and  satisfying the reversibility condition \eqref{reversibility condition in the variables theta I and z}, that is
			\begin{equation}\label{reversibility condition in the variables theta I and z for the torus i0}
				\begin{array}{ccc}
					\vartheta_{0}(-\varphi)=-\vartheta_{0}(\varphi), & I_{0}(-\varphi)=I_{0}(\varphi), & z_{0}(-\varphi)=(\mathscr{S}z_{0})(\varphi).
				\end{array}
			\end{equation}
			As in the previous section, we denote
			$\displaystyle\mathfrak{I}_{0}(\varphi)=i_{0}(\varphi)-(\varphi,0,0).$ 
			Our main goal here is to explore the  invertibility of  the operator 
			\begin{equation}\label{Norm-proj}
				\widehat{\mathcal{L}}_{\omega}=\widehat{\mathcal{L}}_{\omega}(i_{0})=\Pi_{\mathbb{S}_0}^{\perp}\left(\omega\cdot\partial_{\varphi}-\partial_{\theta}\mathsf{K}_{02}(\varphi)\right)\Pi_{\mathbb{S}_0}^{\perp}
			\end{equation}
			defined through \eqref{Lomega def} and \eqref{K02} with  the suitable tame estimates for the inverse. For a precise statement we refer to  Proposition  \ref{inversion of the linearized operator in the normal directions}.  Notice that this operator will be described as  a quasilinear perturbation of the  diagonal operator stated in  Lemma \ref{lemma general form of the linearized operator} and we expect that suitable standard  reductions can be performed to  conjugate it  to a diagonal one provided that the exterior parameters are subject to live in  a Cantor set allowing  to prevent resonances.  For this aim, we shall implement with suitable adaptions the strategy   developed in  the works \cite{BBMH18,BM18}. We distinguish  two long  reduction steps. First, we  perform  a quasi-periodic change of variables such that in the new coordinates  system the transport part is straightened to  a constant coefficient operator. The construction of this transformation is based on a KAM reducibility procedure as in \cite{FGMP19}. The outcome of this first step is a new operator whose positive part is diagonal with a small  nonlocal perturbation of order $-1$. Then the second step consists in applying KAM scheme in order to reduce the remainder and  conjugate the resulting operator from step 1 into a diagonal one up to small errors. The proof follows  basically   a common  procedure that can be found for instance in \cite{B19}. We point out that our results differ slightly from the preceding ones in \cite{BBMH18,BM18}, especially at the level of Cantor sets which are constructed over the  final targets.  \\
	We shall use throughout the proofs  some  frequency cut-offs with respect to  the  sequence defined in \eqref{definition of Nm}, 
		with $N_{0}$ a constant  needed to be  large enough.
					In the current section, the numbers  $N_{0}\geqslant 2$ and $\gamma\in(0,1)$ are a priori free parameters, but during  the Nash-Moser scheme, see  Proposition \ref{Nash-Moser}, they will be adjusted with respect to   $\varepsilon$ according to the relations		
			$$
			N_{0}=\gamma^{-1}\quad\mbox{ and }\quad \gamma=\varepsilon^{a}\quad\mbox{ for some }a>0.
			$$
			We shall set the following parameters required along  the different reductions that we intend to perform, 
			\begin{equation}\label{param}
				\begin{array}{ll} 
					s_{l}:=s_0+\tau_1q+\tau_1+2,& \overline{\mu}_2:=4\tau_1q+6\tau_1+3, \\
					\overline{s}_{l}:=s_l+\tau_2q+\tau_2, & \overline{s}_{h}:=\frac{3}{2}\overline{\mu}_{2}+s_{l}+1,  
				\end{array}
			\end{equation}
			supplemented with the assumptions \eqref{initial parameter condition} and \eqref{setting tau1 and tau2}.
%		$$s_{h}\geqslant\max(\frac{3}{2}\mu_{2}+\overline{s_{l}}+1,\overline{s}_{h}+\mathtt{p})\quad \mu_{2}\geqslant\overline{\mu}_{2}+2\tau_2q+2\tau_2.$$
%		Notice that the parmeter $\mathtt{p}\geqslant 0,$ seems a priori free and has to be well-chosen according to the context in the various proofs. The important fact is that it remains bounded, namely
%		$$0\leqslant\mathtt{p}\leqslant\mathtt{p}_{\textnormal{\tiny{max}}}:=s_0+4\tau_{2}q+4\tau_{2}.$$
%			We shall also make  use of  the following smallness assumption : there exists $\varepsilon_{0}>0$ small enough and a constant $\sigma:=\sigma(\tau_{1},\tau_{2},q,d,s_{0})>0$ such that
%			\begin{equation}\label{small}
%				\varepsilon\gamma^{-2-q}N_{0}^{\mu_{2}}\leqslant\varepsilon_{0}\quad\textnormal{and}\quad \|\mathfrak{I}_{0}\|_{q,s_{h}+\sigma}^{\gamma,\mathcal{O}}\leqslant 1.
%			\end{equation}
%			 The constant $\sigma$ appears in Proposition \ref{inversion of the linearized operator in the normal directions} and represents the final loss of derivatives accumulated along the different reductions  of this section. What is important to mention  is the fact that $\sigma$   is independent of the Sobolev index $s$ and depends only on the geometry of the Cantor sets built in this section.
			 %
			 %\subsection{Structure of the linearized operator in the normal directions}
			 %
			\subsection{Localization on the normal directions}
			According to   Theorem \ref{proposition alomst approximate inverse}, the construction  of an approximate inverse for $d_{i,\alpha}\mathcal{F}(i_{0},\alpha_{0})$ is based on  Proposition \ref{thm:inversion of the linearized operator in the normal directions} dealing with finding an approximate right inverse for  the operator $\widehat{\mathcal{L}}_{\omega}$. This program will be achieved along  several steps and in  the first one we shall  describe  its asymptotic structure around the linearized operator at  the equilibrium state described in Lemma \ref{lemma general form of the linearized operator}.  More precisely, we shall prove the following result.			
			\begin{prop}\label{lemma setting for Lomega}
				{Let $(\gamma,q,d,s_{0})$ satisfy \eqref{initial parameter condition}. 
					Then the operator $\widehat{\mathcal{L}}_{\omega}$ defined in \eqref{Norm-proj} takes the form 
					$$\widehat{\mathcal{L}}_{\omega}=\Pi_{\mathbb{S}_0}^{\perp}\left(\mathcal{L}_{\varepsilon r}-\varepsilon\partial_{\theta}\mathcal{R}\right)\Pi_{\mathbb{S}_0}^{\perp},\quad \mathcal{L}_{\varepsilon r}=\omega\cdot\partial_{\varphi}+\partial_{\theta}\left(V_{\varepsilon r}\cdot\right)-\partial_{\theta}\mathbf{L}_{\varepsilon r},
					$$
					where $V_{\varepsilon r}$ and $\mathbf{L}_{\varepsilon r}$ 		
					 are defined in Lemma $\ref{lemma general form of the linearized operator},$
					and from \eqref{definition of A action-angle-normal} we have 
					\begin{align*}r(\varphi)
					&=A\big(\vartheta_{0}(\varphi),\, I_{0}(\varphi),z_0(\varphi)\big)\\
					&=v\big(\vartheta_{0}(\varphi),I_{0}(\varphi)\big)+z_{0}(\varphi),
					\end{align*}
					supplemented with the reversibility assumption
					\begin{equation}\label{symmetry for r}
						r(\lambda,\omega,-\varphi,-\theta)=r(\lambda,\omega,\varphi,\theta).
					\end{equation}
					Moreover,   $\mathcal{R}$ is an integral operator in the sense of the Definition $ \ref{Defin-Rever-1},$ whose   kernel $J$ satisfies the  symmetry property
					\begin{equation}\label{symmetry for the kernel J}
						J(\lambda,\omega,-\varphi,-\theta,-\eta)=J(\lambda,\omega,\varphi,\theta,\eta).
					\end{equation}
					and  under the assumption
					\begin{equation}\label{frakI0 bnd}
					\|\mathfrak{I}_0\|_{q,s_0}^{\gamma,\mathcal{O}}\leqslant  1,	
					\end{equation} we have for all $s\geqslant s_{0}$, 
\begin{enumerate}[label=(\roman*)]
\item  The function $r$ satisfies the estimates,	\begin{equation}\label{estimate r and mathfrakI0}
	\| r\|_{q,s}^{\gamma,\mathcal{O}}\lesssim 1+\|\mathfrak{I}_{0}\|_{q,s}^{\gamma,\mathcal{O}}
\end{equation}
and 
\begin{equation}\label{control difference r}
	\|\Delta_{12}r\|_{q,s}^{\gamma,\mathcal{O}}\lesssim\|\Delta_{12}i\|_{q,s}^{\gamma,\mathcal{O}}+\| \Delta_{12}i\|_{q,s_0}^{\gamma,\mathcal{O}}\max_{j=1,2}\|\mathfrak{I}_{j}\|_{q,s}^{\gamma,\mathcal{O}}.
\end{equation}
\item The kernel $J$ satisfies the following estimates for all $\ell\in\mathbb{N}$,
						\begin{equation}\label{estimate J}
							\sup_{\eta\in\mathbb{T}}\|(\partial_{\theta}^{\ell}J)(\ast,\cdot,\centerdot,\eta+\centerdot)\|_{q,s}^{\gamma,\mathcal{O}}\lesssim 1+\|\mathfrak{I}_{0}\|_{q,s+3+\ell}^{\gamma,\mathcal{O}}
						\end{equation}
						and
						\begin{equation}\label{differences J} 
							\sup_{\eta\in\mathbb{T}}\|\Delta_{12}(\partial_{\theta}^{\ell}J)(\ast,\cdot,\centerdot,\eta+\centerdot)\|_{q,s}^{\gamma,\mathcal{O}}\lesssim\|\Delta_{12}i\|_{q,s+3+\ell}^{\gamma,\mathcal{O}}+\|\Delta_{12}i\|_{q,s_0+3}^{\gamma,\mathcal{O}}\max_{j=1,2}\|\mathfrak{I}_{j}\|_{q,s+3+\ell}^{\gamma,\mathcal{O}}.
						\end{equation}
						Here  $\ast,\cdot,\centerdot$ stand for $(\lambda,\omega),\varphi,\theta$, respectively and $\displaystyle\mathfrak{I}_{\ell}(\varphi)=i_{\ell}(\varphi)-(\varphi,0,0).$ In addition,  for any function $f$,   $\Delta_{12} f:=f(i_1)-f(i_2)$ refers   for the difference of $f$ taken  at two different states $i_1$ and $i_2$ satisfying \eqref{frakI0 bnd}. 
				\end{enumerate}} 
			\end{prop}
	
	\begin{proof} To alleviate the notation we shall at  several stages of the proof remove the dependence  of the involved  functions/operators with respect to $(\lambda,\omega)$ and keep it  when we deem it  relevant. Recall that the operator $\widehat{\mathcal{L}}_{\omega}$ is defined in  \eqref{Norm-proj}. To describe $\mathsf{K}_{02}(\varphi)$ we  follow \cite{BBMH18,BM18}. First, we observe  from   \eqref{K02} and \eqref{H alpha} that 
				\begin{align}\label{Y-1}
					\nonumber\mathsf{K}_{02}(\varphi)=\mathrm{L}(\lambda)+\varepsilon\partial_{w}\nabla_{w}(\mathcal{P}_{\varepsilon}(i_{0}(\varphi))+\varepsilon\mathcal{R}(\varphi),
				\end{align}
				with
				$$\mathcal{R}(\varphi)=\mathcal{R}_{1}(\varphi)+\mathcal{R}_{2}(\varphi)+\mathcal{R}_{3}(\varphi),$$
				where 
				\begin{align*}
					\mathcal{R}_{1}(\varphi)&:=L_{2}^{\top}(\varphi)\partial_{I}\nabla_I\mathcal{P}_{\varepsilon}(i_{0}(\varphi))L_{2}(\varphi),\\ \mathcal{R}_{2}(\varphi)&:=L_{2}^{\top}(\varphi)\partial_{z}\nabla_{I}\mathcal{P}_{\varepsilon}(i_{0}(\varphi)),\\
					\mathcal{R}_{3}(\varphi)&:=\partial_{I}\nabla_{z}\mathcal{P}_{\varepsilon}(i_{0}(\varphi))L_{2}(\varphi).
				\end{align*}
				As we shall see, all the operators  $\mathcal{R}_{1}(\varphi),$ $\mathcal{R}_{2}(\varphi)$ and  $\mathcal{R}_{3}(\varphi)$   have a  finite-dimensional rank. This property is obvious  for  the operator $L_{2}(\varphi)$  defined in \eqref{L2}, which sends in view of \eqref{dual rel} the space  $H_{\perp}^{s}$ to $\mathbb{R}^{d}$ and therefore for any $\rho\in H_{\perp}^{s}$ we write 
				$$L_{2}(\varphi)[\rho]=\sum_{k=1}^{d}\big\langle L_{2}(\varphi)[\rho],\underline{e}_{k}\big\rangle_{\mathbb{R}^{d}}\,\underline{e}_{k}=\sum_{k=1}^{d}\big\langle\rho,L_{2}^{\top}(\varphi)[\underline{e}_{k}]\big\rangle_{L^{2}(\mathbb{T})}\,\underline{e}_{k},$$
				with $\displaystyle(\underline{e}_{k})_{k=1}^d$ being the canonical basis of $\mathbb{R}^d.$
				Hence
				$$\begin{array}{ll}
					\displaystyle \mathcal{R}_{1}(\varphi)[\rho]=\sum_{k=1}^{d}\big\langle\rho,L_{2}^{\top}(\varphi)[\underline{e}_{k}]\big\rangle_{L^{2}(\mathbb{T})}A_{1}(\varphi)[\underline{e}_{k}] & \quad\mbox{with }\quad A_{1}(\varphi)=L_{2}^{\top}(\varphi)\partial_{I}\nabla_I\mathcal{P}_{\varepsilon}(i_{0}(\varphi)),\\
					\displaystyle \mathcal{R}_{3}(\varphi)[\rho]=\sum_{k=1}^{d}\big\langle\rho,L_{2}^{\top}(\varphi)[\underline{e}_{k}]\big\rangle_{L^{2}(\mathbb{T})}A_{3}(\varphi)[\underline{e}_{k}] &\quad \mbox{with }\quad A_{3}(\varphi)=\partial_{I}\nabla_{z}\mathcal{P}_{\varepsilon}(i_{0}(\varphi)).
				\end{array}$$
				In a similar way, by setting  $A_{2}(\varphi):=\partial_{z}\nabla_{I}\mathcal{P}_{\varepsilon}(i_{0}(\varphi)):H_{\perp}^{s}\rightarrow\mathbb{R}^{d},$ then we may write
				$$
				\mathcal{R}_{2}(\varphi)[\rho]=\sum_{k=1}^{d}\big\langle\rho,A_{2}^{\top}(\varphi)[\underline{e}_{k}]\big\rangle_{L^{2}\,(\mathbb{T})}L_{2}^{\top}(\varphi)[\underline{e}_{k}].
				$$
				Define 
				$$
				g_{k,1}(\varphi,\theta)=g_{k,3}(\varphi,\theta)=\chi_{k,2}(\varphi,\theta):=L_{2}^{\top}(\varphi)[\underline{e}_{k}](\theta),\quad g_{k,2}(\varphi,\theta):=A_{2}^{\top}(\varphi)[\underline{e}_{k}](\theta)
				$$
				and
				$$
				\mbox{ }\chi_{k,1}(\varphi,\theta):=A_{1}(\varphi)[\underline{e}_{k}](\theta),\quad \chi_{k,3}(\varphi,\theta):=A_{3}(\varphi)[\underline{e}_{k}](\theta),
				$$
				then we can see that   the operator $\mathcal{R}$ takes the integral form
				\begin{align*}\mathcal{R}\rho(\varphi,\theta)&=\sum_{k'=1}^{3}\sum_{k=1}^{d}\langle\rho(\varphi,\cdot),g_{k,k'}(\varphi,\cdot)\rangle_{L^{2}(\mathbb{T})}\chi_{k,k'}(\varphi,\theta)\\
					&=\int_{\mathbb{T}}\rho(\varphi,\eta)J(\varphi,\theta,\eta)d\eta,
				\end{align*}
				with 
				$$J(\varphi,\theta,\eta):=\sum_{k'=1}^{3}\sum_{k=1}^{d}g_{k,k'}(\varphi,\eta)\chi_{k,k'}(\varphi,\theta).$$
				Now we remark that  by construction $g_{k,k'},\chi_{k,k'}\in H_{\perp}^{s}$  with
				\begin{equation}\label{estimate gkk' and chikk'}
					\| g_{k,k'}\|_{q,s}^{\gamma,\mathcal{O}}+\|\chi_{k,k'}\|_{q,s}^{\gamma,\mathcal{O}}\lesssim 1+\|\mathfrak{I}_{0}\|_{q,s+3}^{\gamma,\mathcal{O}}
				\end{equation}
				and straightforward computations yield
				\begin{equation}\label{estimate differential gkk' and chikk'}
					\| d_{i}g_{k,k'}[\widehat{i}]\|_{q,s}^{\gamma,\mathcal{O}}+\| d_{i}\chi_{k,k'}[\widehat{i}]\|_{q,s}^{\gamma,\mathcal{O}}\lesssim\|\widehat{i}\|_{q,s+2}^{\gamma,\mathcal{O}}+\|\mathfrak{I}_{0}\|_{q,s+4}^{\gamma,\mathcal{O}}\|\widehat{i}\|_{q,s_{0}+2}^{\gamma,\mathcal{O}}.
				\end{equation}
			On the other hand, one has from direct computations that
			$$\forall \ell\in\mathbb{N},\quad (\partial_{\theta}^{\ell}J)(\varphi,\theta,\eta+\theta)=\sum_{k'=1}^{3}\sum_{k=1}^{d}g_{k,k'}(\varphi,\eta+\theta)(\partial_{\theta}^{\ell}\chi_{k,k'})(\varphi,\theta).
			$$
			Hence,  we may combine  \eqref{estimate gkk' and chikk'} with Lemma \ref{Lem-lawprod}-(iv) and \eqref{frakI0 bnd} allowing to get  
				\begin{align*}
					\sup_{\eta\in\mathbb{T}}\| (\partial_{\theta}^{\ell}J)(\ast,\cdot,\centerdot,\eta+\centerdot)\|_{q,s}^{\gamma,\mathcal{O}}&\lesssim \sum_{k'=1}^{3}\sum_{k=1}^{d}\|g_{k,k'}(\ast,\cdot,\eta+\centerdot)\|_{q,s}^{\gamma,\mathcal{O}}\|\chi_{k,k'}(\ast,\cdot,\centerdot)\|_{q,s_{0}+\ell}^{\gamma,\mathcal{O}}\\
					&\quad+\sum_{k'=1}^{3}\sum_{k=1}^{d}\|g_{k,k'}(\ast,\cdot,\eta+\centerdot)\|_{q,s_{0}}^{\gamma,\mathcal{O}}\|\chi_{k,k'}(\ast,\cdot,\centerdot)\|_{q,s+\ell}^{\gamma,\mathcal{O}}\\
					&\lesssim1+\|\mathfrak{I}_{0}\|_{q,s+3+\ell}^{\gamma,\mathcal{O}},
				\end{align*}
				where we have used the  interpolation inequality: for $s\geqslant s_0$
				\begin{align*}
				\|g_{k,k'}(\ast,\cdot,\eta+\centerdot)\|_{q,s}^{\gamma,\mathcal{O}}\|\chi_{k,k'}(\ast,\cdot,\centerdot)\|_{q,s_{0}+\ell}^{\gamma,\mathcal{O}}&\lesssim  \|g_{k,k'}(\ast,\cdot,\eta+\centerdot)\|_{q,s+\ell}^{\gamma,\mathcal{O}}\|\chi_{k,k'}(\ast,\cdot,\centerdot)\|_{q,s_{0}}^{\gamma,\mathcal{O}}\\
				&\quad+\|g_{k,k'}(\ast,\cdot,\eta+\centerdot)\|_{q,s_0}^{\gamma,\mathcal{O}}\|\chi_{k,k'}(\ast,\cdot,\centerdot)\|_{q,s+\ell}^{\gamma,\mathcal{O}}
				\end{align*}
			In addition, to estimate the difference we simply write
			\begin{align*}
				\forall \ell\in\mathbb{N},\quad \Delta_{12}(\partial_{\theta}^{\ell}J)(\varphi,\theta,\eta+\theta)&=\sum_{k'=1}^{3}\sum_{k=1}^{d}\Delta_{12}g_{k,k'}(\varphi,\eta+\theta)(\partial_{\theta}^{\ell}(\chi_{k,k'})_{r_{1}})(\varphi,\theta)\\
				&\quad+\sum_{k'=1}^{3}\sum_{k=1}^{d}(g_{k,k'})_{r_{2}}(\varphi,\eta+\theta)(\Delta_{12}\partial_{\theta}^{\ell}\chi_{k,k'})(\varphi,\theta).
			\end{align*}
				By applying the mean value theorem combined  with \eqref{estimate differential gkk' and chikk'} and \eqref{frakI0 bnd} combined with interpolation inequalities
				\begin{align*}
				\sup_{\eta\in\mathbb{T}}\| \Delta_{12}(\partial_{\theta}^{\ell}J)(\ast,\cdot,\centerdot,\eta+\centerdot)\|_{q,s}^{\gamma,\mathcal{O}}&\lesssim \sum_{k'=1}^{3}\sum_{k=1}^{d}\|\Delta_{12}g_{k,k'}(\ast,\cdot,\eta+\centerdot)\|_{q,s}^{\gamma,\mathcal{O}}\|\chi_{k,k'}(\ast,\cdot,\centerdot)\|_{q,s_{0}+\ell}^{\gamma,\mathcal{O}}\\
				&\quad+\sum_{k'=1}^{3}\sum_{k=1}^{d}\|\Delta_{12}g_{k,k'}(\ast,\cdot,\eta+\centerdot)\|_{q,s_0}^{\gamma,\mathcal{O}}\|\chi_{k,k'}(\ast,\cdot,\centerdot)\|_{q,s+\ell}^{\gamma,\mathcal{O}}\\
					&\quad+\sum_{k'=1}^{3}\sum_{k=1}^{d}\|g_{k,k'}(\ast,\cdot,\eta+\centerdot)\|_{q,s_{0}}^{\gamma,\mathcal{O}}\|\Delta_{12}\chi_{k,k'}(\ast,\cdot,\centerdot)\|_{q,s+\ell}^{\gamma,\mathcal{O}}\\
					&\quad+\sum_{k'=1}^{3}\sum_{k=1}^{d}\|g_{k,k'}(\ast,\cdot,\eta+\centerdot)\|_{q,s}^{\gamma,\mathcal{O}}\|\Delta_{12}\chi_{k,k'}(\ast,\cdot,\centerdot)\|_{q,s_0+\ell}^{\gamma,\mathcal{O}}\\
					&\lesssim\|\Delta_{12}i\|_{q,s+3+\ell}^{\gamma,\mathcal{O}}+\|\Delta_{12}i\|_{q,s_0+3}^{\gamma,\mathcal{O}}\max_{j=1,2}\|\mathfrak{I}_{j}\|_{q,s+3+\ell}^{\gamma,\mathcal{O}}.
				\end{align*}
				The symmetry property detailed in  \eqref{symmetry for the kernel J} is a consequence of the definition of $r$ and the reversibility condition \eqref{reversibility condition in the variables theta I and z for the torus i0} imposed on the torus $i_{0}.$ Consequently, putting together \eqref{HEE} and \eqref{definition of A action-angle-normal} gives
				\begin{align*}
					\mathsf{K}_{02}(\varphi) & =  \mathrm{L}(\lambda)\Pi_{\mathbb{S}_0}^{\perp}+\varepsilon\partial_{z}\nabla_{z}\mathcal{P}_{\varepsilon}(i_{0}(\varphi))+\varepsilon\mathcal{R}(\varphi)\\
					& =  \mathrm{L}(\lambda)\Pi_{\mathbb{S}_0}^{\perp}+\varepsilon\Pi_{\mathbb{S}_0}^{\perp}\partial_{r}\nabla_{r}P_{\varepsilon}(A(i_{0}(\varphi)))\Pi_{\mathbb{S}_0}^{\perp}+\varepsilon\mathcal{R}(\varphi)\\
					& =  \Pi_{\mathbb{S}_0}^{\perp}\partial_{r}\nabla_{r}\mathcal{H}_{\varepsilon}(A(i_{0}(\varphi)))\Pi_{\mathbb{S}_0}^{\perp}+\varepsilon\mathcal{R}(\varphi)\\
					& =  \Pi_{\mathbb{S}_0}^{\perp}\partial_{r}\nabla_{r}H(\varepsilon A(i_{0}(\varphi)))\Pi_{\mathbb{S}_0}^{\perp}+\varepsilon\mathcal{R}(\varphi).
				\end{align*}
			Recall from \eqref{definition of A action-angle-normal} that 
			\begin{equation}\label{defr}
				r(\varphi,\cdot)=A(i_{0}(\varphi)),
			\end{equation}
				then according to the general form of the linearized operator stated  in Lemma \ref{lemma general form of the linearized operator} one has
				$$-\partial_{\theta}\partial_{r}\nabla_{r}H(\varepsilon r(\varphi,\cdot))=\partial_{\theta}\left(V_{\varepsilon r}\cdot\right)-\partial_{\theta}\mathbf{L}_{\varepsilon r},$$
				which implies in turn
				$$\begin{array}{rcl}
					-\mathsf{K}_{02}(\varphi) 
					= \Pi_{\mathbb{S}_0}^{\perp}\big(\partial_{\theta}\left(V_{\varepsilon r}\cdot\right)-\partial_{\theta}\mathbf{L}_{\varepsilon r})-\varepsilon\mathcal{R}(\varphi)\big)\Pi_{\mathbb{S}_0}^{\perp}.
				\end{array}$$
				Plugging this identity into \eqref{Norm-proj} gives the desired result. Next, using \eqref{defr},  \eqref{definition of A action-angle-normal} and \eqref{v-theta1}, we obtain
					\begin{align*}
						\nonumber
						\| r\|_{q,s}^{\gamma,\mathcal{O}}&\lesssim \| v(\vartheta_0,I_0)\|_{q,s}^{\gamma,\mathcal{O}}+\| z_0\|_{q,s}^{\gamma,\mathcal{O}}\\
						&\lesssim 1+\|\mathfrak{I}_{0}\|_{q,s}^{\gamma,\mathcal{O}}.
					\end{align*}
					We shall now move  to the proof of the bound \eqref{control difference r}. First, we observe from \eqref{definition of A action-angle-normal} that
					\begin{align*}
						\nonumber
						\| \Delta_{12}r\|_{q,s}^{\gamma,\mathcal{O}}&\lesssim \| \Delta_{12}v(\vartheta,I)\|_{q,s}^{\gamma,\mathcal{O}}+\| \Delta_{12}z\|_{q,s}^{\gamma,\mathcal{O}}.
					\end{align*}
					Therefore,  Taylor Formula with \eqref{v-theta1} and law products allow to get 
					\begin{align*}
						\| \Delta_{12}v(\vartheta,I)\|_{q,s}^{\gamma,\mathcal{O}}&\lesssim\| \Delta_{12}(I,\vartheta)\|_{q,s}^{\gamma,\mathcal{O}}+\| \Delta_{12}(I,\vartheta)\|_{q,s_0}^{\gamma,\mathcal{O}}\max_{j=1,2}\|\mathfrak{I}_{j}\|_{q,s}^{\gamma,\mathcal{O}},
					\end{align*}
					which implies  that
					\begin{align*}
						\nonumber
						\| \Delta_{12}r\|_{q,s}^{\gamma,\mathcal{O}}
						&\lesssim \| \Delta_{12}i\|_{q,s}^{\gamma,\mathcal{O}}+\| \Delta_{12}i\|_{q,s_0}^{\gamma,\mathcal{O}}\max_{j=1,2}\|\mathfrak{I}_{j}\|_{q,s}^{\gamma,\mathcal{O}}.
					\end{align*}
		This achieves the proof of Proposition \ref{lemma setting for Lomega}.		
				
			\end{proof}
			\subsection{Reduction of order 1}\label{Reduction of order 1-1}
			In this section, we perform the reduction of the transport part of the linearized operator $\mathcal{L}_{\varepsilon r}$ described in Proposition  \ref{lemma setting for Lomega}. More precisely, we conjugate the operator $\mathcal{L}_{\varepsilon r}$ by a quasi-periodic symplectic change of variables $\mathscr{B}$ leading to  a  transport part with constant coefficients depending only on the torus $i_{0}$ and the parameters  $\varepsilon,$ $\lambda$ and $\omega.$ To get a precise information on the remainder, which is of order $-1$ in $\theta$, we  need to describe  the action of this conjugation on the nonlocal term using the kernel structure rather than pseudo-differential theory. The reduction to a constant coefficient operator is based on KAM scheme through the construction of successive  quasi-periodic symplectic change of coordinates. This will be implemented in the same spirit of  \cite{BFM21-1,FGMP19}. Here we need to extend their construction to the framework of  of symplectic change of coordinates with  $C^{q}$ regularity. We point out that  similar results with slight variations  have been established in \cite{BM20,BFM21-1} in a non-symplectic framework.\\
			\subsubsection{Reduction of the transport part}
			Before stating our result, we need to introduce some transformations. 
			Let $\beta: \mathcal{O}\times \T^{d+1}\to \mathbb{T}$ be a smooth function such that $\displaystyle\sup_{\mu\in \mathcal{O}}\|\beta(\mu,\cdot,\centerdot)\|_{\textnormal{Lip}}<1$ 
			then  the map
			$$(\varphi,\theta)\in\T^{d+1}\mapsto (\varphi, \theta+\beta(\mu,\varphi,\theta))\in\T^{d+1}$$
			is a diffeomorphism and its  inverse takes the form  
			$$(\varphi,\theta)\in\T^{d+1}\mapsto (\varphi, \theta+\widehat{\beta}(\mu,\varphi,\theta))\in\T^{d+1}.$$
			The relation between $\beta$ and $\widehat\beta$ is described through,
			\begin{equation}\label{def betahat}
				y=\theta+\beta(\mu,\varphi,\theta)\Longleftrightarrow \theta=y+\widehat\beta(\mu,\varphi,y).
			\end{equation}
			Now we define the operators
			\begin{equation}\label{definition symplectic change of variables}
				\mathscr{B}=(1+\partial_{\theta}\beta)\mathcal{B},
			\end{equation}
			with $$\mathcal{B}\rho(\mu,\varphi,\theta)=\rho\big(\mu,\varphi,\theta+\beta(\mu,\varphi,\theta)\big).$$
			Direct computations show that the inverse $\mathscr{B}^{-1}$ keeps the same form, that is, 
			\begin{equation}\label{mathscrB1}
				\mathscr{B}^{-1}\rho(\mu,\varphi,y)=\Big(1+\partial_y\widehat{\beta}(\mu,\varphi,y)\Big) \rho\big(\mu,\varphi,y+\widehat{\beta}(\mu,\varphi,y)\big)
			\end{equation}
			and
			$$
			\mathcal{B}^{-1} \rho(\mu,\varphi,y)=\rho\big(\mu,\varphi,y+\widehat{\beta}(\mu,\varphi,y)\big).
			$$
			We shall now give some elementary algebraic properties for $\mathcal{B}^{\pm 1}$ and $\mathscr{B}^{\pm 1}$ which can be checked by straightforward computations.
			\begin{lem}\label{algeb1}
				The following assertions hold true.
				\begin{enumerate}[label=(\roman*)]
					%\item The operator $\mathscr{B}$ is invertible  and its inverse $\mathscr{B}^{-1}$ keeps the same form, that is, 
					%\begin{equation*}
					% \mathscr{B}^{-1}h(\varphi,y)=(1+\partial_y\widehat{\beta}(\varphi,y)) h\big(\varphi,y+\widehat{\beta}(\varphi,y)\big)
					%\end{equation*}
					\item The action of $\mathscr{B}^{-1}$ on the derivative is given by
					\begin{equation*}
						\mathscr{B}^{-1}\partial_{\theta}=\partial_{\theta}\mathcal{B}^{-1}.
					\end{equation*}
					\item The conjugation of the transport operator by $\mathscr{B}$  keeps the same structure
					$$
					\mathscr{B}^{-1}\Big(\omega\cdot\partial_\varphi+\partial_\theta\big(V(\varphi,\theta)\cdot\big)\Big)\mathscr{B}=\omega\cdot\partial_\varphi+\partial_y(\mathscr{V}(\varphi,y)\cdot\big),
					$$
					with
					$$
					\mathscr{V}(\varphi,y)=\mathcal{B}^{-1}\Big(\omega\cdot\partial_{\varphi} \beta(\varphi,\theta)+V(\varphi,\theta)\big(1+\partial_\theta \beta(\varphi,\theta)\big)\Big).
					$$
					\item Denote by $\mathscr{B}^\star$ the $L^2_\theta(\T)$-adjoint of $\mathscr{B}$, then
					$$
					\mathscr{B}^\star=\mathcal{B}^{-1}\quad\hbox{and}\quad \mathcal{B}^\star=\mathscr{B}^{-1}.
					$$
					%\item We have the estimates
					%$$
					%\|\mathscr{B} h\|_{s}^{\gamma,q}\lesssim (1+\|\beta\|_{q,s_0+1}^{\gamma,\mathcal{O}}\big)^2\|h\|_{s}^{\gamma,q}+(1+\|\beta\|_{q,s_0+1}^{\gamma,\mathcal{O}}\big)(1+\|\beta\|_{q,s+1}^{\gamma,\mathcal{O}}\big)\|h\|_{s_0}^{\gamma,q}
					%$$
					%and
					%$$
					%\|\mathscr{B}^{-1} h\|_{s}^{\gamma,q}\lesssim (1+\|\widehat\beta\|_{q,s_0+1}^{\gamma,\mathcal{O}}\big)^2\|h\|_{s}^{\gamma,q}+(1+\|\widehat\beta\|_{q,s_0+1}^{\gamma,\mathcal{O}}\big)(1+\|\widehat\beta\|_{q,s+1}^{\gamma,\mathcal{O}}\big)\|h\|_{s_0}^{\gamma,q}
					%$$
					%In addition 
					%$$
					%\|\mathcal{B} h\|_{s}^{\gamma,q}\lesssim (1+\|\beta\|_{q,s_0+1}^{\gamma,\mathcal{O}}\big)\|h\|_{s}^{\gamma,q}+(1+\|\beta\|_{q,s+1}^{\gamma,\mathcal{O}}\big)\|h\|_{s_0}^{\gamma,q}
					%$$
					
				\end{enumerate}
			\end{lem}
			Now we shall state the following result proved in \cite{FGMP19} for $q=1$ and which can be obtained by induction for a general $q\in\mathbb{N}^{*}$ up to slight modifications. We also refer to \cite[(A.2)]{BFM21-1}.
			\begin{lem}\label{Compos1-lemm}
				Let $(q,d\gamma,s_0)$ as in \eqref{initial parameter condition}.
				Let $\beta\in W^{q,\infty,\gamma}\big(\mathcal{O},H^{\infty}(\T^{d+1})\big) $ such that 
				\begin{equation}\label{small beta lem}
					\|\beta \|_{q,2s_0}^{\gamma,\mathcal{O}}\leqslant \varepsilon_0,
				\end{equation}
				with $\varepsilon_{0}$ small enough. Then the following assertions hold true.
				\begin{enumerate}[label=(\roman*)]
					\item The linear operators $\mathcal{B},\mathscr{B}:W^{q,\infty,\gamma}\big(\mathcal{O},H^{s}(\T^{d+1})\big)\to W^{q,\infty,\gamma}\big(\mathcal{O},H^{s}(\T^{d+1})\big)$ are continuous and invertible, with 
					\begin{equation}\label{tame comp}
						\forall s\geqslant s_0,\quad \|\mathcal{B}^{\pm1}\rho\|_{q,s}^{\gamma,\mathcal{O}}\leqslant \|\rho\|_{q,s}^{\gamma,\mathcal{O}}\left(1+C\|\beta\|_{q,s_{0}}^{\gamma,\mathcal{O}}\right)+C\|\beta\|_{q,s}^{\gamma,\mathcal{O}}\|\rho\|_{q,s_{0}}^{\gamma,\mathcal{O}}
					\end{equation}
					and 
					\begin{equation}\label{tame comp symp}
						\forall s\geqslant s_0,\quad\|\mathscr{B}^{\pm1}\rho\|_{q,s}^{\gamma,\mathcal{O}}\leqslant \|\rho\|_{q,s}^{\gamma,\mathcal{O}}\left(1+C\|\beta\|_{q,s_{0}}^{\gamma,\mathcal{O}}\right)+C\|\beta\|_{q,s+1}^{\gamma,\mathcal{O}}\|\rho\|_{q,s_{0}}^{\gamma,\mathcal{O}}.
					\end{equation}
					%\begin{equation}
					%\|\mathcal{B}^{\pm1}\rho-\rho\|_{q,s}^{\gamma,\mathcal{O}}\lesssim_{s,q} \|\rho\|_{q,s_{0}}^{\gamma,\mathcal{O}}\|\beta\|_{q,s+d}^{\gamma,\mathcal{O}}+\|\beta\|_{q,s_{0}}^{\gamma,\mathcal{O}}\|\rho\|_{q,s+1}^{\gamma,\mathcal{O}}.
					%\end{equation}
					\item The functions $\beta$ and $\widehat{\beta}$ are linked through
					\begin{equation}\label{link betah and beta}
						\forall s\geqslant s_0,\quad\|\widehat{\beta}\|_{q,s}^{\gamma,\mathcal{O}}\leqslant C\|\beta\|_{q,s}^{\gamma,\mathcal{O}}.
					\end{equation}
					\item Let $\beta_{1},\beta_{2}\in W^{q,\infty,\gamma}(\mathcal{O},H^{\infty}(\mathbb{T}^{d+1}))$ satisfying \eqref{small beta lem}. If we denote 
					$$\Delta_{12}\beta=\beta_{1}-\beta_{2}\quad\textnormal{ and }\quad\Delta_{12}\widehat{\beta}=\widehat{\beta}_{1}-\widehat{\beta}_{2},$$
					then they are linked through
					\begin{equation}\label{link diff beta hat and diff beta}
						\forall s\geqslant s_0,\quad\|\Delta_{12}\widehat{\beta}\|_{q,s}^{\gamma,\mathcal{O}}\leqslant C\left(\|\Delta_{12}\beta\|_{q,s}^{\gamma,\mathcal{O}}+\|\Delta_{12}\beta\|_{q,s_{0}}^{\gamma,\mathcal{O}}\max_{j\in\{1,2\}}\|\beta_{j}\|_{q,s+1}^{\gamma,\mathcal{O}}\right).
					\end{equation}
				\end{enumerate}
			\end{lem}
			\begin{proof}
				\textbf{(i)}-\textbf{(ii)} For \eqref{tame comp} and \eqref{link betah and beta}, we refer to \cite[(A.2)]{BFM21-1} and \cite[Lem. A.3.]{FGMP19}. The estimate \eqref{tame comp symp} is obtained from \eqref{tame comp} and law product in Lemma \ref{Lem-lawprod}.\\
				\textbf{(iii)} One has by Taylor Formula
				\begin{align*}
					\Delta_{12}\widehat{\beta}(y)&=\widehat{\beta}_{1}(y)-\widehat{\beta}_{2}(y)\\
					&=\beta_{2}(y+\widehat{\beta}_{2}(y))-\beta_{1}(y+\widehat{\beta}_{1}(y))\\
					&=-\Delta_{12}\beta(y+\widehat{\beta}_{2}(y))-\Delta_{12}\widehat{\beta}(y)\int_{0}^{1}\partial_{\theta}\beta_{1}(y+\widehat{\beta}_{1}(y)-t\Delta_{12}\widehat{\beta}(y))dt.
				\end{align*}
				Hence
				$$\Delta_{12}\widehat{\beta}(y)=\tfrac{-\mathcal{B}_{2}^{-1}\Delta_{12}\beta(y)}{1+\mathscr{I}(y)}\quad\textnormal{with}\quad\mathscr{I}(y):=\int_{0}^{1}\partial_{\theta}\beta_{1}(y+\widehat{\beta}_{1}(y)-t\Delta_{12}\widehat{\beta}(y))dt.$$
				By composition estimate in Lemma \ref{Lem-lawprod}, one has
				$$\Big\|\tfrac{1}{1+\mathscr{I}}\Big\|_{q,s}^{\gamma,\mathcal{O}}\lesssim 1+\|\mathscr{I}\|_{q,s}^{\gamma,\mathcal{O}}.$$
				Thus, applying the  law product in Lemma \ref{Lem-lawprod} implies
				$$\|\Delta_{12}\widehat{\beta}\|_{q,s}^{\gamma,\mathcal{O}}\lesssim\left(1+\|\mathscr{I}\|_{q,s}^{\gamma,\mathcal{O}}\right)\|\mathcal{B}_{2}^{-1}\Delta_{12}\beta\|_{q,s_{0}}^{\gamma,\mathcal{O}}+\left(1+\|\mathscr{I}\|_{q,s_{0}}^{\gamma,\mathcal{O}}\right)\|\mathcal{B}_{2}^{-1}\Delta_{12}\beta\|_{q,s}^{\gamma,\mathcal{O}}.$$
				Using \eqref{tame comp}, \eqref{link betah and beta} and \eqref{small beta lem} yields
				\begin{align*}
					\|\mathcal{B}_{2}^{-1}\Delta_{12}\beta\|_{q,s}^{\gamma,\mathcal{O}}&\lesssim\|\Delta_{12}\beta\|_{q,s}^{\gamma,\mathcal{O}}\left(1+\|\widehat{\beta}_{2}\|_{q,s_{0}}^{\gamma,\mathcal{O}}\right)+\|\widehat{\beta}_{2}\|_{q,s}^{\gamma,\mathcal{O}}\|\Delta_{12}\beta\|_{q,s_{0}}^{\gamma,\mathcal{O}}\\
					&\lesssim \|\Delta_{12}\beta\|_{q,s}^{\gamma,\mathcal{O}}+\|\widehat{\beta}_{2}\|_{q,s}^{\gamma,\mathcal{O}}\|\Delta_{12}\beta\|_{q,s_{0}}^{\gamma,\mathcal{O}}\\
					&\lesssim \|\Delta_{12}\beta\|_{q,s}^{\gamma,\mathcal{O}}+\|\beta_{2}\|_{q,s}^{\gamma,\mathcal{O}}\|\Delta_{12}\beta\|_{q,s_{0}}^{\gamma,\mathcal{O}}
				\end{align*}
				and
				\begin{align*}
					\|\mathscr{I}\|_{q,s}^{\gamma,\mathcal{O}}&\lesssim\|\beta_{1}\|_{q,s+1}^{\gamma,\mathcal{O}}\left(1+\|\widehat{\beta}_{1}\|_{q,s_{0}}^{\gamma,\mathcal{O}}+\|\Delta_{12}\widehat{\beta}\|_{q,s_{0}}^{\gamma,\mathcal{O}}\right)+\left(\|\widehat{\beta}_{1}\|_{q,s}^{\gamma,\mathcal{O}}+\|\Delta_{12}\widehat{
						\beta}\|_{q,s}^{\gamma,\mathcal{O}}\right)\|\beta_{1}\|_{q,s_{0}+1}^{\gamma,\mathcal{O}}\\
					&\lesssim \|\beta_{1}\|_{q,s+1}^{\gamma,\mathcal{O}}+\|\Delta_{12}\widehat{\beta}\|_{q,s}^{\gamma,\mathcal{O}}\|\beta_{1}\|_{q,s_{0}+1}^{\gamma,\mathcal{O}}.
				\end{align*}
				Putting together the foregoing estimates gives
				\begin{align}\label{est diff betah proof}
					\|\Delta_{12}\widehat{\beta}\|_{q,s}^{\gamma,\mathcal{O}}&\leqslant C\left(1+\|\beta_{1}\|_{q,s+1}^{\gamma,\mathcal{O}}+\|\Delta_{12}\widehat{\beta}\|_{q,s}^{\gamma,\mathcal{O}}\|\beta_{1}\|_{q,s_{0}+1}^{\gamma,\mathcal{O}}\right)\left(1+\|\beta_{2}\|_{q,s_{0}}^{\gamma,\mathcal{O}}\right)\|\Delta_{12}\beta\|_{q,s_{0}}^{\gamma,\mathcal{O}}\nonumber\\
					&+C\left(1+\|\beta_{1}\|_{q,s_{0}+1}^{\gamma,\mathcal{O}}+\|\Delta_{12}\widehat{\beta}\|_{q,s_{0}}^{\gamma,\mathcal{O}}\|\beta_{1}\|_{q,s_{0}+1}^{\gamma,\mathcal{O}}\right)\left(\|\Delta_{12}\beta\|_{q,s}^{\gamma,\mathcal{O}}+\|\beta_{2}\|_{q,s}^{\gamma,\mathcal{O}}\|\Delta_{12}\beta\|_{q,s_{0}}^{\gamma,\mathcal{O}}\right).
				\end{align}
				From the triangle inequality, \eqref{link betah and beta} and \eqref{small beta lem}, one has
				\begin{align*}
					\|\Delta_{12}\widehat{\beta}\|_{q,s_{0}}^{\gamma,\mathcal{O}}&\leqslant\|\widehat{\beta}_{1}\|_{q,s_{0}}^{\gamma,\mathcal{O}}+\|\widehat{\beta}_{2}\|_{q,s_{0}}^{\gamma,\mathcal{O}}\\
					&\leqslant \|\beta_{1}\|_{q,s_{0}}^{\gamma,\mathcal{O}}+\|\beta_{2}\|_{q,s_{0}}^{\gamma,\mathcal{O}}\\
					&\leqslant 2\varepsilon_{0}.
				\end{align*}
				From Sobolev embeddings we infer that
				$$
				\max_{j\in\{1,2\}}\|\beta_{j}\|_{q,s_{0}+1}^{\gamma,\mathcal{O}}\leqslant \max_{j\in\{1,2\}}\|\beta_{j}\|_{q,2s_{0}}^{\gamma,\mathcal{O}}\leqslant\varepsilon_0.
				$$
								Thus, by choosing $\varepsilon_{0}$ small enough, we can ensure
				$$C\|\Delta_{12}\widehat{\beta}\|_{q,s}^{\gamma,\mathcal{O}}\|\beta_{1}\|_{q,s_{0}+1}^{\gamma,\mathcal{O}}\left(1+\|\beta_{2}\|_{q,s_{0}}^{\gamma,\mathcal{O}}\right)\|\Delta_{12}\beta\|_{q,s_{0}}^{\gamma,\mathcal{O}}\leqslant\frac{1}{2}\|\Delta_{12}\widehat{\beta}\|_{q,s}^{\gamma,\mathcal{O}}.$$
				Inserting  this term into  the left hand side in \eqref{est diff betah proof} and using Sobolev embeddings, we find
				$$\|\Delta_{12}\widehat{\beta}\|_{q,s}^{\gamma,\mathcal{O}}\leqslant C\left(\|\Delta_{12}\beta\|_{q,s}^{\gamma,\mathcal{O}}+\|\Delta_{12}\beta\|_{q,s_{0}}^{\gamma,\mathcal{O}}\max_{j\in\{1,2\}}\|\beta_{j}\|_{q,s+1}^{\gamma,\mathcal{O}}\right).$$
				This ends  the proof of Lemma \ref{Compos1-lemm}.
			\end{proof}
		
			Now we shall state the main result of this section concerning the reduction of the transport part of the linearized operator $\mathcal{L}_{\varepsilon r}.$
			\begin{prop}\label{reduction of the transport part}
				Let $(\gamma,q,d,\tau_{1},s_{0},S,s_l,\overline{s}_h,\overline{\mu}_{2})$ satisfy \eqref{initial parameter condition}, \eqref{setting tau1 and tau2} and \eqref{param}. Let  $\upsilon\in\left(0,\frac{1}{q+2}\right].$ We set
				\begin{equation}\label{sigma1}
					\sigma_{1}=s_0+\tau_{1}q+2\tau_{1}+4.
				\end{equation}
				For any $(\mu_2,\mathtt{p},s_h)$ satisfying
					\begin{equation}\label{param-trans}
						\mu_{2}\geqslant \overline{\mu}_{2}:=4\tau_{1}q+6\tau_{1}+3,\quad \mathtt{p}\geqslant 0, \quad s_{h}\geqslant\max\left(\frac{3}{2}\mu_{2}+s_{l}+1,\overline{s}_{h}+\mathtt{p}\right),
					\end{equation}
					there exists $\varepsilon_{0}>0$ such that if
					\begin{equation}\label{smallness condition transport}
						\varepsilon\gamma^{-1}N_{0}^{\mu_{2}}\leqslant\varepsilon_{0}\quad\textnormal{and}\quad \|\mathfrak{I}_{0}\|_{q,s_{h}+\sigma_{1}}^{\gamma,\mathcal{O}}\leqslant 1,
					\end{equation}
				 there exist
					$$c_{i_{0}}\in W^{q,\infty,\gamma }(\mathcal{O},\mathbb{R})\quad\mbox{ and }\quad\beta\in \bigcap_{s\in[s_{0},S]}W^{q,\infty,\gamma }(\mathcal{O},H_{\mbox{\tiny{\textnormal{odd}}}}^{s})$$
					such that with $\mathscr{B}$ defined in \eqref{definition symplectic change of variables} one gets the following results.
					\begin{enumerate}[label=(\roman*)]
						\item The function $c_{i_{0}}$ satisfies the following estimate,
						\begin{equation}\label{estimate r1}
							\| c_{i_{0}}-V_{0}\|_{q}^{\gamma ,\mathcal{O}}\lesssim \varepsilon,
						\end{equation}
						where $V_{0}$ is defined in Lemma $\ref{lemma linearized operator at equilibrium}.$
						\item The transformations $\mathscr{B}^{\pm 1},\mathcal{B}^{\pm 1}, {\beta}$ and $\widehat{\beta}$ satisfy the following estimates for all $s\in[s_{0},S]$ 
						\begin{equation}\label{estimate on the first reduction operator and its inverse}
							\|\mathscr{B}^{\pm 1}\rho\|_{q,s}^{\gamma ,\mathcal{O}}+\|\mathcal{B}^{\pm 1}\rho\|_{q,s}^{\gamma ,\mathcal{O}}\lesssim\|\rho\|_{q,s}^{\gamma ,\mathcal{O}}+\varepsilon\gamma ^{-1}\| \mathfrak{I}_{0}\|_{q,s+\sigma_{1}}^{\gamma ,\mathcal{O}}\|\rho\|_{q,s_{0}}^{\gamma ,\mathcal{O}}
						\end{equation}
						and 
						\begin{equation}\label{estimate beta and r}
							\|\widehat{\beta}\|_{q,s}^{\gamma,\mathcal{O}}\lesssim\|\beta\|_{q,s}^{\gamma ,\mathcal{O}}\lesssim \varepsilon\gamma ^{-1}\left(1+\| \mathfrak{I}_{0}\|_{q,s+\sigma_{1}}^{\gamma ,\mathcal{O}}\right).
						\end{equation}
						\item Let $n\in\mathbb{N}$, then in the truncated Cantor set
						$$\mathcal{O}_{\infty,n}^{\gamma,\tau_{1}}(i_{0})=\bigcap_{(l,j)\in\mathbb{Z}^{d}\times\mathbb{Z}\setminus\{(0,0)\}\atop|l|\leqslant N_{n}}\left\lbrace(\lambda,\omega)\in \mathcal{O}\quad\textnormal{s.t.}\quad\big|\omega\cdot l+jc_{i_{0}}(\lambda,\omega)\big|>\tfrac{4\gamma^{\upsilon}\langle j\rangle}{\langle l\rangle^{\tau_{1}}}\right\rbrace,$$
						we have 
						$$\mathscr{B}^{-1}\big(\omega\cdot\partial_{\varphi}+\partial_{\theta}\big(V_{\varepsilon r}\cdot\big)\big)\mathscr{B}=\omega\cdot\partial_{\varphi}+c_{i_{0}}\partial_{\theta}+\mathtt{E}_{n}^{0},$$
						with $\mathtt{E}_{n}^{0}=\mathtt{E}_{n}^{0}(\lambda,\omega,i_{0})$ a linear operator satisfying
%						\begin{equation}\label{estim En0 s}
%							\forall s\in[s_{0},S],\quad\|\mathtt{E}_{n}^{0}\rho\|_{q,s}^{\gamma,\mathcal{O}}\lesssim\|\rho\|_{q,s+2}^{\gamma,\mathcal{O}}+\varepsilon\gamma^{-1}\|\mathfrak{I}_{0}\|_{q,s+\sigma_{1}}^{\gamma,\mathcal{O}}\|\rho\|_{q,s_{0}+2}^{\gamma,\mathcal{O}}
%						\end{equation}
%						and
						\begin{equation}\label{estim En0 s0}
							\|\mathtt{E}_{n}^{0}\rho\|_{q,s_{0}}^{\gamma,\mathcal{O}}\lesssim\varepsilon N_{0}^{\mu_{2}}N_{n+1}^{-\mu_{2}}\|\rho\|_{q,s_{0}+2}^{\gamma,\mathcal{O}}.
						\end{equation}
						\item Given two tori $i_{1}$ and $i_{2}$ both satisfying \eqref{smallness condition transport}, we have 
						\begin{equation}\label{difference ci}
							\|\Delta_{12}c_{i}\|_{q}^{\gamma,\mathcal{O}}\lesssim\varepsilon\| \Delta_{12}i\|_{q,\overline{s}_{h}+2}^{\gamma,\mathcal{O}}
						\end{equation}
						and
						\begin{equation}\label{difference beta}
							\|\Delta_{12}\beta\|_{q,\overline{s}_{h}+\mathtt{p}}^{\gamma,\mathcal{O}}+\|\Delta_{12}\widehat\beta\|_{q,\overline{s}_{h}+\mathtt{p}}^{\gamma,\mathcal{O}}\lesssim\varepsilon\gamma^{-1}\|\Delta_{12}i\|_{q,\overline{s}_{h}+\mathtt{p}+\sigma_{1}}^{\gamma,\mathcal{O}}.
						\end{equation}
					\end{enumerate}
			\end{prop}
			%\begin{proof}
			%The proof of this result, based on KAM technics, is now rather classical (see \cite{FGMP19} for the Lipschitz case and \cite{BM20} \textcolor{blue}{article pas encore publié, juste sur arxiv} for a more similar result in a non-symplectic case). Here the proposition is given in a symplectic context, which is a small modification of the previous results mentioned above. For completeness and convinience of the reader, we put the whole proof of this result in appendix B.
			%\end{proof}
			Before giving the proof, some remarks are in order.
			\begin{remark}
				\begin{enumerate}[label=\textbullet]
					\item The final Cantor set $\mathcal{O}_{\infty,n}^{\gamma,\tau_{1}}(i_{0})$ is constructed over the limit coefficient $c_{i_0}$ but it is still truncated in the time frequency, that is  $|l|\leqslant N_n$,  leading to a residual remainder with enough decay. This induces a suitable stability property that is crucial during the Nash-Moser scheme achieved with  the nonlinear functional.
					\item  Notice that, since $4\gamma^{\upsilon}\geqslant\gamma$, then looking at $j=0$ we find that the Cantor set $\mathcal{O}_{\infty,n}^{\gamma,\tau_{1}}(i_0)$ is contained in the Diophantine Cantor set $(\lambda_0,\lambda_1)\times\mathtt{DC}_{N_n}(\gamma,\tau_{1})$ introduced in \eqref{DC tau gamma}.
					\item The parameter $\upsilon$ is introduced for technical reasons appearing later in the measure estimates of the final Cantor set and it will be fixed in \eqref{choice tau 1 tau2 upsilon}.
					\item The constant $4$ used in the definition of the Cantor set $\mathcal{O}_{\infty,n}^{\gamma,\tau_1}(i_{0})$ is useful to ensure the inclusion of this set in all the Cantor sets built in the KAM procedure (see \eqref{inclusion Oinftyn in Om+1} in the proof below) and also to establish some inclusions related to  the final Cantor set (see the proof of Lemma $\ref{some cantor set are empty}$).
					\item We emphasize here that the functions $\beta$ and $\widehat{\beta}$ are odd in the sense
					\begin{equation}\label{symmetry for beta}
						\beta(\lambda,\omega,-\varphi,-\theta)=-\beta(\lambda,\omega,\varphi,\theta)\quad\textnormal{and}\quad\widehat{\beta}(\lambda,\omega,-\varphi,-\theta)=-\widehat{\beta}(\lambda,\omega,\varphi,\theta)
					\end{equation}
					which will be crucial later to get the Toeplitz structure  of the new  remainder term emerging  after this reduction.
					%\item We point out  that the conditions stated in  \eqref{param-trans} and \eqref{smallness condition transport} are much weaker than those of  \eqref{param} and \eqref{small} where we should take into account of several restrictions arising during the different reduction steps. 
				\end{enumerate}
			\end{remark}
			\begin{proof}
				Since we are looking at a state near the disc, we can split $V_{\varepsilon r}$ defined  by \eqref{definition of Vr} according to 

				\begin{equation}\label{V=V0+f0}
					V_{\varepsilon r}(\lambda,\varphi,\theta)=V_{0}(\lambda)+f_{0}(\lambda,\varphi,\theta),
				\end{equation}
				with $f_{0}$ being a  perturbation term of small size. We refer to \eqref{initial smallness condition in sh norm} for a more precise quantification of this smallness. The proof is an iteration process introducing at each step a linear quasi-periodic symplectic change of coordinates. This transformation is linked to the remainder term of the previous step. Roughly speaking, if the latter is of size $\varepsilon$, then we choose the change of coordinates in such a way that we extract the main diagonal part of the previous remainder and keep a new perturbation term of size $\varepsilon^2.$ The choice of the transformation is done through the resolution of an homological equation requiring non-resonance conditions capted by a suitable selection of the parameters of the system. Thus, by iteration, we can construct a final Cantor set gathering all the parameters restrictions of all steps in which we completely reduced the transport operator into a constant coefficient one. We shall now explain a typical step of the procedure Later, we shall implement the scheme.\\
				\textbf{(i)-(ii)} $\blacktriangleright$ \textbf{KAM step}.
				Let us consider a transport operator in the form,	$$\omega\cdot\partial_{\varphi}+\partial_{\theta}\Big(V+f\Big)$$
				for suitable parameters $(\lambda,\omega)$ that belong to  a subset $\mathcal{O}_{-}^{\gamma}\subset\mathcal{O}$, where $\mathcal{O}$ is the ambient set and 
				$$
				V=V(\lambda,\omega)\quad\textnormal{and}\quad f=f(\lambda,\omega,\varphi,\theta),
				$$
				where $f$ enjoys the following symmetry condition
				\begin{equation}\label{symmetry for f}
					f(\lambda,\omega,-\varphi,-\theta)=f(\lambda,\omega,\varphi,\theta).
				\end{equation}
				To alleviate the notations we shall use during the proof the variable $\mu:=(\lambda,\omega)$. We consider a symplectic quasi-periodic change of coordinates close to the identity taking the form
				\begin{equation}\label{change of variables}
					\begin{array}{rcl}
						\mathscr{G}\rho(\mu,\varphi,\theta) &:=& \big(1+\partial_{\theta}g(\mu,\varphi,\theta)\big)\mathcal{G}\rho(\mu,\varphi,\theta)\\
						& := & \big(1+\partial_{\theta}g(\mu,\varphi,\theta)\big)\rho\big(\mu,\varphi,\theta+g(\mu,\varphi,\theta)\big),
					\end{array}
				\end{equation}
				where $g:\mathcal{O}\times\mathbb{T}^{d+1}\rightarrow\mathbb{R}$ is a small which will be later linked to $f.$ Then, by using Lemma \ref{algeb1}, we can write for any  $N\geqslant 2$
				\begin{equation}\label{transformation KAM step transport}
					\mathscr{G}^{-1}\Big(\omega\cdot\partial_{\varphi}+\partial_{\theta}\left(V+f\right)\Big)\mathscr{G}=\omega\cdot\partial_{\varphi}+\partial_{\theta}\mathcal{G}^{-1}\left(V+\omega\cdot\partial_{\varphi}g+V\partial_{\theta}g+\Pi_{N}f+\Pi_{N}^{\perp}f+f\partial_{\theta}g\right).
				\end{equation}
				Recall that the projections $\Pi_{N}$ are defined in \eqref{definition of projections for functions}. The basic idea is to obtain after this transformation a new transport operator in the form
				\begin{equation}\label{link V,f and V+,f+}
					\mathscr{G}^{-1}\Big(\omega\cdot\partial_{\varphi}+\partial_{\theta}\left(V+f\right)\Big)\mathscr{G}=\omega\cdot\partial_{\varphi}+\partial_{\theta}\big(V_{+}+f_{+}\big),
				\end{equation}
				where
				$$
				V_{+}=V_{+}(\mu)\quad\textnormal{and}\quad f_{+}=f_{+}(\mu,\varphi,\theta),
				$$
				with $f_{+}$ quadratically smaller than $f.$ In order to get rid of the terms wich are not small of quadratic in $f$, then, in view of \eqref{transformation KAM step transport}, we shall select $g$ solving the following \textit{homological equation}
				\begin{equation}\label{equation satisfied by g}
					\omega\cdot\partial_{\varphi}g+V\partial_{\theta}g+\Pi_{N}f=\langle f\rangle_{\varphi,\theta},
				\end{equation}
				where $$\langle f\rangle_{\varphi,\theta}(\mu):=\int_{\mathbb{T}^{d+1}}f(\mu,\varphi,\theta)d\varphi d\theta.$$
				To find a solution to the \textit{homological equation} \eqref{equation satisfied by g}, we use Fourier decomposition and look for $g$ in the form
				\begin{equation}\label{definition of g}
					g(\mu,\varphi,\theta):=\ii\sum_{(l,j)\in\mathbb{Z}^{d+1}\setminus\{0\}\atop\langle l,j\rangle\leqslant N}\tfrac{f_{l,j}(\mu)}{\omega\cdot l+jV(\mu)}e^{\ii(l\cdot\varphi+j\theta)}.
				\end{equation}
				The denominators appearing in the Fourier decomposition of $g$ may be small and generate problems in the convergence of the series in \eqref{definition of g} for large $N.$ This is a well-knonw phenomenon in KAM theory called "small divisors problem". To overcome this difficulty, one has to avoid the resonances and, following the ideas of Kolmogorov, we introduce Diophantine conditions gathered in the following Cantor set	\begin{equation}\label{definition of O+}
					\mathcal{O}_{+}^{\gamma}:=\bigcap_{(l,j)\in\mathbb{Z}^{d+1 }\setminus\{0\}\atop\langle l,j\rangle\leqslant N}\left\lbrace\mu:=(\lambda,\omega)\in \mathcal{O}_{-}^{\gamma}\quad\textnormal{s.t.}\quad\big|\omega\cdot l+jV(\mu)\big|>\tfrac{\gamma ^{\upsilon}\langle j\rangle}{\langle l\rangle^{\tau_{1}}}\right\rbrace.
				\end{equation}
				Such a selection of the external parameters allows us to control the size of the denominators in \eqref{definition of g}. As we shall see in \eqref{control of g by f}, the quantification of this control, linked to the parameters $\gamma$ and $\tau_1,$ allows to get suitable  estimates for $g$ with some loss of regularity  uniform with respect to  $N$. Before performing  this estmate, we shall first  construct  an extension of $g$ to the whole set $\mathcal{O}.$ In what follows, we still denote $g$ this extension. This is done by extending the Fourier coefficients of $g$ using the cut-off function $\chi$ defined in \eqref{properties cut-off function first reduction}. More precisely, we define
				\begin{align}\label{def extension glj}
					g_{l,j}(\mu)&:=\ii\tfrac{\chi\big((\omega\cdot l+jV(\mu))(\gamma ^{\upsilon}\langle j\rangle)^{-1}\langle l\rangle^{\tau_{1}}\big)}{\omega\cdot l+jV(\mu)}f_{l,j}(\mu)\\
					&:=\widetilde{g_{l,j}}(\mu)f_{l,j}(\mu)\nonumber.
				\end{align}
				Notice that the extension $g$ is a solution to \eqref{equation satisfied by g} only when the parameters are restricted to the Cantor set $\mathcal{O}_{+}^{\gamma}.$  Then, we define
				$$V_{+}=V+\langle f\rangle_{\varphi,\theta}\quad\textnormal{and}\quad f_{+}=\mathcal{G}^{-1}\big(\Pi_{N}^{\perp}f+f\partial_{\theta}g\big),$$
so that in restriction to the Cantor set $\mathcal{O}_{+}^{\gamma},$ the identity \eqref{link V,f and V+,f+} holds. Remark that $V_{+}$ and $f_{+}$ are well-defined in the whole set of parameters $\mathcal{O}$ and the function $g$ is smooth since it is generated by a finite number of frequencies. According to \eqref{symmetry for f}, we obtain that $g$ is odd. As a consequence, 
\begin{equation}\label{symmetry for g}
					g\in\bigcap_{s\geqslant 0}W^{q,\infty,\gamma}(\mathcal{O},H_{\mbox{\tiny{odd}}}^{s}).
				\end{equation}
				Our next task is to estimate the Fourier coefficients $\widetilde{g_{l,j}}$ defined by \eqref{def extension glj}. Notice that we can write them in the following form 
\begin{align}\label{set gtlj}
					\widetilde{g_{l,j}}(\mu)&=\ii\,a_{l,j}\widehat{\chi}(a_{l,j}A_{l,j}(\mu)),\quad\widehat{\chi}(x):=\tfrac{\chi(x)}{x}\\ A_{l,j}(\mu)&:=\omega\cdot l+jV(\mu),\quad a_{l,j}:=(\gamma^{\upsilon}\langle j\rangle)^{-1}\langle l\rangle^{\tau_{1}}.\nonumber
				\end{align}
				Since $\widehat{\chi}$ is $C^{\infty}$ with bounded derivatives and $\widehat{\chi}(0)=0,$ then applying Lemma \ref{Lem-lawprod}-(vi), we obtain
				$$\forall q'\in\llbracket0,q\rrbracket,\quad \|\widetilde{g_{l,j}}\|_{q'}^{\gamma,\mathcal{O}}\lesssim a_{l,j}^{2}\|A_{l,j}\|_{q'}^{\gamma,\mathcal{O}}\Big(1+a_{l,j}^{q'-1}\|A_{l,j}\|_{L^{\infty}(\mathcal{O})}^{q'-1}\Big).$$
Direct computations lead to
\begin{align*}
	\forall(l,j)\in\mathbb{Z}^{d+1},\,\forall\alpha\in\mathbb{N}^{d+1},\quad|\alpha|\leqslant q,\quad \sup_{\mu\in\mathcal{O}}\left|\partial_{\mu}^{\alpha}A(\mu)\right|&\lesssim\langle l,j\rangle\max\left(1,\sup_{\mu\in\mathcal{O}}\left|\partial_{\mu}^{\alpha}V(\mu)\right|\right)\\
					&\lesssim\gamma^{-|\alpha|}\langle l,j\rangle\max\left(1,\|V\|_{q}^{\gamma,\mathcal{O}}\right).
				\end{align*}
				Assuming
				\begin{equation}\label{boundedness assumption on V}
					\|V\|_{q}^{\gamma,\mathcal{O}}\leqslant C,
				\end{equation}
				we then obtain
				\begin{equation}\label{est Alj}
					\forall q'\in\llbracket0,q\rrbracket,\quad\forall(l,j)\in\mathbb{Z}^{d+1},\quad \|A_{l,j}\|_{q'}^{\gamma,\mathcal{O}}\lesssim\langle l,j\rangle.
				\end{equation}
				Added to the fact that $0\leqslant a_{l,j}\leqslant\gamma^{-\upsilon}\langle l\rangle^{\tau_{1}}$, we then find that
				\begin{equation}\label{gljtqprm}
					\forall q'\in\llbracket0,q\rrbracket,\quad\|\widetilde{g_{l,j}}\|_{q'}^{\gamma,\mathcal{O}}\lesssim\gamma^{-\upsilon(q'+1)}\langle l,j\rangle^{\tau_{1}q'+\tau_{1}+q'}.
				\end{equation}
%				$$\forall(l,j)\in\mathbb{Z}^{d+1},\,\forall\alpha\in\mathbb{N}^{d},\quad|\alpha|\leqslant q\Rightarrow\sup_{\mu\in\mathcal{O}}\left|\partial_{\mu}^{\alpha}A(\mu)\right|\lesssim\gamma^{-|\alpha|}\langle l,j\rangle.$$
%				By using Fa\`a di Bruno formula and the fact that $\upsilon\in(0,1)$, we find 
%				\begin{align*}
%					\forall(l,j)\in\mathbb{Z}^{d+1},\forall\alpha\in\mathbb{N}^{d+1},\quad|\alpha|\leqslant q,\quad \sup_{(\lambda,\omega)\in\mathcal{O}}\left|\partial_{\mu}^{\alpha}g_{l,j}(\mu)\right|&\lesssim a^{|\alpha|+1}\gamma^{-|\alpha|}\sum_{\sum km_{k}=|\alpha|}\prod_{p=1}^{|\alpha|}\langle l,j\rangle^{m_{p}}\\
%					&\lesssim\gamma^{-(q+1)\upsilon}\gamma^{-|\alpha|}\langle l,j\rangle^{\tau_{1}(|\alpha|+1)+|\alpha|}.
%				\end{align*}
				Our choice of $\upsilon$ in Proposition \ref{reduction of the transport part} implies in particular that
\begin{equation}\label{first choice of upsilon}
	\upsilon\leqslant\tfrac{1}{q+1}.
\end{equation}
Therefore, we deduce from \eqref{def extension glj} and Leibniz rule that for all $\alpha\in\mathbb{N}^{d+1}$ with $|\alpha|\leqslant q$
				\begin{align*}
					\gamma^{2|\alpha|}\|\partial_{\mu}^{\alpha}g(\mu,\cdot,\centerdot)\|_{H^{s-|\alpha|}}^{2}&\lesssim\sum_{(l,j)\in\mathbb{Z}^{d+1}\setminus\{(0,0)\}\atop\langle l,j\rangle\leqslant N}\sum_{\beta\in\mathbb{N}^{d+1}\atop\beta\leqslant\alpha}\gamma^{2|\alpha|-2|\beta|}\big|\partial_{\mu}^{\alpha-\beta}\widetilde{g_{l,j}}(\mu)\big|^{2}\gamma^{2|\beta|}\big|\partial_{\mu}^{\beta}f_{l,j}(\mu)\big|^{2}\langle l,j\rangle^{2s-2|\alpha|}\\
					&\lesssim\sum_{(l,j)\in\mathbb{Z}^{d+1}\setminus\{(0,0)\}\atop\langle l,j\rangle\leqslant N}\sum_{\beta\in\mathbb{N}^{d+1}\atop\beta\leqslant\alpha}\left(\|\widetilde{g_{l,j}}\|_{|\alpha|-|\beta|}^{\gamma,\mathcal{O}}\right)^{2}\gamma^{2|\beta|}\big|\partial_{\mu}^{\beta}f_{l,j}(\mu)\big|^{2}\langle l,j\rangle^{2s-2|\alpha|}\\
					&\lesssim\sum_{(l,j)\in\mathbb{Z}^{d+1}\setminus\{(0,0)\}\atop\langle l,j\rangle\leqslant N}\sum_{\beta\in\mathbb{N}^{d+1}\atop\beta\leqslant\alpha}\gamma^{-2}\gamma^{2|\beta|}\big|\partial_{\mu}^{\beta}f_{l,j}(\mu)\big|^{2}\langle l,j\rangle^{2(s+\tau_{1}q+\tau_{1}-|\beta|)}.
				\end{align*}
			As a consequence, by interverting the summation symbols, we find
				\begin{equation}\label{control of g by f}
					\|g\|_{q,s}^{\gamma ,\mathcal{O}}\lesssim\gamma ^{-1}\|\Pi_{N}f\|_{q,s+\tau_{1} q+\tau_{1}}^{\gamma ,\mathcal{O}}.
				\end{equation}
				Assume now that
				\begin{equation}\label{smallness assumption f}
					\gamma^{-1}N^{\tau_{1}q+\tau_{1}+1}\|f\|_{q,s_{0}}^{\gamma,\mathcal{O}}\leqslant\varepsilon_{0}.
				\end{equation}
				Then  added to \eqref{control of g by f} and Lemma \ref{Lem-lawprod}-(ii), we get
				$$\|g\|_{q,s_{0}}^{\gamma,\mathcal{O}}\leqslant C\gamma^{-1}N^{\tau_{1}q+\tau_{1}}\|f\|_{q,s_{0}}^{\gamma,\mathcal{O}}\leqslant C\varepsilon_{0}.$$
				On the other hand if we assume
				$$\|f\|_{q,s}^{\gamma,\mathcal{O}}\lesssim\varepsilon\left(1+\|\mathfrak{I}_0\|_{q,s+1}^{\gamma,\mathcal{O}}\right),$$
				then \eqref{control of g by f} gives
				\begin{align*}
					\|g\|_{q,2s_0+1}^{\gamma,\mathcal{O}}&\lesssim\gamma^{-1}\|f\|_{q,2s_0+\tau_{1}q+\tau_{1}+1}^{\gamma,\mathcal{O}}\\
					&\lesssim\varepsilon\gamma^{-1}\left(1+\|\mathfrak{I}_0\|_{q,2s_0+\tau_{1}q+\tau_{1}+2}^{\gamma,\mathcal{O}}\right)\\
					&\lesssim\varepsilon\gamma^{-1}\left(1+\|\mathfrak{I}_0\|_{q,s_h+\sigma_{1}}^{\gamma,\mathcal{O}}\right).
				\end{align*}
			Notice that to obtain the last inequality we used the fact that \eqref{param-trans} and \eqref{sigma1} imply
			$$2s_0+\tau_{1}q+\tau_{1}+2\leqslant s_h+\sigma_{1}.$$
				Using interpolation inequality and \eqref{smallness condition transport}, one gets for some $\overline{\theta}\in(0,1).$
				\begin{align}
					\|g\|_{q,2s_0}^{\gamma,\mathcal{O}}&\lesssim\left(\|g\|_{q,s_0}^{\gamma,\mathcal{O}}\right)^{\overline{\theta}}\left(\|g\|_{q,2s_0+1}^{\gamma,\mathcal{O}}\right)^{1-\overline{\theta}}\nonumber\\
					&\lesssim\varepsilon_0.\label{e-g-2s_0}
				\end{align}
				Thus, taking $\varepsilon_{0}$ small enough, we can ensure  the smallness condition in Lemma \ref{Compos1-lemm} and get  that the linear operator $\mathscr{G}$ is  invertible. Now, we introduce
				$$u=\Pi_{N}^{\perp}f+f\partial_{\theta}g.$$
				By the triangle  inequality,  Lemma \ref{Lem-lawprod}-(ii) and \eqref{control of g by f}, we obtain for all $s\in[s_{0},S]$ 
				\begin{align*}
					\|u\|_{q,s}^{\gamma ,\mathcal{O}} 
					& \leqslant  \|\Pi_{N}^{\perp}f\|_{q,s}^{\gamma ,\mathcal{O}}+C\Big(\|f\|_{q,s_{0}}^{\gamma ,\mathcal{O}}\|\partial_{\theta}g\|_{q,s}^{\gamma ,\mathcal{O}}+\|f\|_{q,s}^{\gamma ,\mathcal{O}}\|\partial_{\theta}g\|_{q,s_{0}}^{\gamma ,\mathcal{O}}\Big)\\
					& \leqslant  \|\Pi_{N}^{\perp}f\|_{q,s}^{\gamma ,\mathcal{O}}+C\Big(\|f\|_{q,s_{0}}^{\gamma ,\mathcal{O}}\|g\|_{q,s+1}^{\gamma ,\mathcal{O}}+\|f\|_{q,s}^{\gamma ,\mathcal{O}}\|g\|_{q,s_{0}+1}^{\gamma ,\mathcal{O}}\Big)\\
					& \leqslant  \|\Pi_{N}^{\perp}f\|_{q,s}^{\gamma ,\mathcal{O}}+C\gamma ^{-1}N^{\tau_{1}q+\tau_{1}+1}\|f\|_{q,s_{0}}^{\gamma ,\mathcal{O}}\|f\|_{q,s}^{\gamma ,\mathcal{O}}.
				\end{align*}
				Combined with  Lemma \ref{Compos1-lemm}, Lemma \ref{Lem-lawprod}-(ii) and \eqref{smallness assumption f}, we get for all $s\in[s_{0},S]$
				\begin{align*}
					\|f_{+}\|_{q,s}^{\gamma ,\mathcal{O}} & =  \|\mathcal{G}^{-1}(u)\|_{q,s}^{\gamma ,\mathcal{O}}\\
					& \leqslant  \|u\|_{q,s}^{\gamma ,\mathcal{O}}+C\Big(\| u\|_{q,s}^{\gamma ,\mathcal{O}}\|\widehat{g}\|_{q,s_{0}}^{\gamma ,\mathcal{O}}+\|\widehat{g}\|_{q,s}^{\gamma ,\mathcal{O}}\|u\|_{q,s_{0}}^{\gamma ,\mathcal{O}}\Big)\\
					& \leqslant  \|u\|_{q,s}^{\gamma ,\mathcal{O}}+C\Big(\|u\|_{q,s}^{\gamma ,\mathcal{O}}\|g\|_{q,s_{0}}^{\gamma ,\mathcal{O}}+\|g\|_{q,s}^{\gamma ,\mathcal{O}}\|u\|_{q,s_{0}}^{\gamma ,\mathcal{O}}\Big)\\
					& \leqslant  \|\Pi_{N}^{\perp}f\|_{q,s}^{\gamma ,\mathcal{O}}+C\gamma ^{-1}N^{\tau_{1} q+\tau_{1}+1}\|f\|_{q,s_{0}}^{\gamma ,\mathcal{O}}\|f\|_{q,s}^{\gamma ,\mathcal{O}}.
				\end{align*}
				Using once again Lemma \ref{Lem-lawprod}-(ii), we find for $S\geqslant\overline{s}\geqslant s\geqslant s_{0}$
				\begin{equation}\label{estimate KAM step transport}
					\|f_{+}\|_{q,s}^{\gamma ,\mathcal{O}}\leqslant N^{s-\overline{s}}\|f\|_{q,\overline{s}}^{\gamma ,\mathcal{O}}+C\gamma ^{-1}N^{\tau_{1} q+\tau_{1}+1}\|f\|_{q,s_{0}}^{\gamma ,\mathcal{O}}\|f\|_{q,s}^{\gamma ,\mathcal{O}}.
				\end{equation}
				$\blacktriangleright$ \textbf{KAM scheme.} Let us now assume that we have constructed $V_{m}$ and $f_{m}$, well-defined in the whole set of parameters $\mathcal{O}$ and satisfying the assumptions \eqref{boundedness assumption on V} and \eqref{smallness assumption f}. We shall now construct the corresponding quantity at the next order, namely $V_{m+1}$ and $f_{m+1}$, still satisfying  \eqref{boundedness assumption on V} and \eqref{smallness assumption f}. For this aim, we shall implement the KAM step with $(V,f,V_{+},f_{+},N)$ replaced by $(V_{m},f_{m},V_{m+1},f_{m+1},N_{m}).$ More  precisely, we will shall prove by induction the existence of a sequence $\{V_m,f_m\}_{m\in\mathbb{N}}$ such that 
				\begin{equation}\label{hypothesis of induction deltam}
					\delta_{m}(s_{l})\leqslant\delta_{0}(s_{h})N_{0}^{\mu_{2}}N_{m}^{-\mu_{2}}\quad \mbox{ and }\quad\delta_{m}(s_{h})\leqslant \left(2-\tfrac{1}{m+1}\right)\delta_{0}(s_{h})
				\end{equation}
				and
				\begin{equation}\label{assumptions KAM iterations}
					\|V_{m}\|_{q}^{\gamma,\mathcal{O}}\leqslant C\quad\textnormal{and}\quad N_{m}^{\tau_{1}q+\tau_{1}+1}\delta_{m}(s_{0})\leqslant\varepsilon_{0},
				\end{equation}
				with $f_m$ satisfying the following symmetry condition
				\begin{equation}\label{symmetry for f_{m}}
					f_{m}(\mu,-\varphi,-\theta)=f_{m}(\mu,\varphi,\theta)			\end{equation}
				and where we denote
				$$\delta_{m}(s):=\gamma ^{-1}\| f_{m}\|_{q,s}^{\gamma ,\mathcal{O}}.$$
 Recall that the parameters $s_l$ and $s_h$ were intruduced in \eqref{param} and \eqref{param-trans}.\\				\ding{226} \textit{Initialization.} We shall first check that the estimates \eqref{hypothesis of induction deltam} and  \eqref{assumptions KAM iterations} are satisfied for $m=0.$ In which case the functions $V_0$ and $f_0$ are defined by \eqref{definition of V0} and \eqref{V=V0+f0}. By \eqref{Estim-Vr-October} and \eqref{estimate r and mathfrakI0} we infer 
				\begin{align}\label{estimate delta0 and I0}
					\delta_{0}(s)&=\gamma^{-1}\|V_{\varepsilon r}-V_{0}\|_{q,s}^{\gamma,\mathcal{O}}\nonumber\\
					&\lesssim\varepsilon\gamma^{-1}\|r\|_{q,s+1}^{\gamma,\mathcal{O}}\nonumber\\
					&\lesssim\varepsilon\gamma^{-1}\left(1+\|\mathfrak{I}_{0}\|_{q,s+1}^{\gamma,\mathcal{O}}\right).
				\end{align}
				Thus, the notation \eqref{param-trans} and the smallness condition \eqref{smallness condition transport} imply that
				\begin{equation}\label{initial smallness condition in sh norm}
					N_{0}^{\mu_{2}}\delta_{0}(s_{h})\leqslant C\varepsilon_{0}.
				\end{equation}
				In addition, by \eqref{symmetry for r} and \eqref{symmetry for Vr}, we deduce that $f_{0}$ satisfies the following symmetry condition 
				\begin{equation}\label{symmetry for f0}
					f_{0}(\lambda,-\varphi,-\theta)=f_{0}(\lambda,\varphi,\theta).
				\end{equation}
				We set $\mathcal{O}_{0}^{\gamma}=\mathcal{O}$ and consider $N_{0}\geqslant 2.$ Our next task is to check that the assumptions \eqref{boundedness assumption on V} and \eqref{smallness assumption f} are satisfied by $V_{0}$ and $f_{0}.$ First recall that $V_{0}$ is defined by
				$$V_{0}(\lambda)=\Omega+I_{1}(\lambda)K_{1}(\lambda).$$
				Using the smooth regularity of \eqref{form-Lebe0}, we obtain
				\begin{equation}\label{boundedness of V0}
					\|V_{0}\|_{q}^{\gamma,\mathcal{O}}\leqslant C.
				\end{equation}
				Therefore, the required boundedness property \eqref{boundedness assumption on V} is satisfied with $V=V_{0}.$ Now by \eqref{param-trans}, we have 
				\begin{equation}\label{mu2-1}
					\mu_{2}\geqslant\tau_{1}q+\tau_1+2.
				\end{equation}
			Hence, using \eqref{initial smallness condition in sh norm}, we obtain
				\begin{align*}
					\gamma^{-1}N_{0}^{\tau_{1}q+\tau_{1}+1}\|f_{0}\|_{q,s_{0}}^{\gamma,\mathcal{O}}&= N_{0}^{\tau_{1}q+\tau_{1}+1}\delta_{0}(s_{0})\\
					&\leqslant N_{0}^{\tau_{1}q+\tau_{1}+1-\mu_{2}}N_{0}^{\mu_{2}}\delta_{0}(s_{h})\\
					&\leqslant C\varepsilon_{0}N_{0}^{-1}.
				\end{align*}
				By taking  $N_{0}$ large enough we get
				\begin{equation}\label{Est-NO}
					CN_{0}^{-1}\leqslant 1,
				\end{equation}
				so that
				$$\gamma^{-1}N_{0}^{\tau_{1}q+\tau_{1}+1}\|f_{0}\|_{q,s_{0}}^{\gamma,\mathcal{O}}\leqslant\varepsilon_{0}.$$
				Hence, the assumption \eqref{smallness assumption f} is satisfied for $f=f_{0}.$ This ends the initialization step. \\
				\ding{226} {\textit{Iteration.}}  let us now assume  that we have constructed $V_m$ and $f_m$ enjoying the  properties \eqref{hypothesis of induction deltam}, \eqref{assumptions KAM iterations} and \eqref{symmetry for f_{m}}. We shall see how to construct $V_{m+1}$ and $f_{m+1}$. According to the KAM step, we consider a symplectic quasi-periodic change of variables $\mathscr{G}_{m}$ taking the form 
				\begin{align*}
					\mathscr{G}_{m}\rho(\mu,\varphi,\theta)&:=\big(1+\partial_{\theta}g_{m}(\mu,\varphi,\theta)\big)\mathcal{G}_{m}\rho(\mu,\varphi,\theta)\\
					&=\big(1+\partial_{\theta}g_{m}(\mu,\varphi,\theta)\big)\rho\big(\mu,\varphi,\theta+g_{m}(\mu,\varphi,\theta)\big),
				\end{align*}
				with \begin{equation}\label{def gm}
					g_{m}(\mu,\varphi,\theta):=\ii\sum_{(l,j)\in\mathbb{Z}^{d+1}\setminus\{0\}\atop\langle l,j\rangle\leqslant N_{m}}\tfrac{\chi\big((\omega\cdot l+jV_{m}(\mu))(\gamma^{\upsilon}\langle j\rangle)^{-1}\langle l\rangle^{\tau_{1}}\big)}{\omega\cdot l+jV_{m}(\mu)}(f_{m})_{l,j}(\mu)\,e^{i(l\cdot\varphi+j\theta)},
				\end{equation}
				where $\chi$ is the cut-off function introduced in \eqref{properties cut-off function first reduction} and $N_{m}$ is defined in \eqref{definition of Nm}. As explained in the KAM step, $g_{m}$ is well-defined on the whole set of parameters $\mathcal{O}$ and solves the \textit{homological equation}
				$$\omega\cdot\partial_{\varphi}g_{m}+V_{m}\partial_{\theta}g_{m}+\Pi_{N_{m}}f_{m}=\langle f_{m}\rangle_{\varphi,\theta}
				$$				when   restricted to the Cantor set 
				\begin{equation}\label{definition mathcal Om+1}
					\mathcal{O}_{m+1}^{\gamma}:=\bigcap_{\underset{\langle l,j\rangle\leqslant N_{m}}{(l,j)\in\mathbb{Z}^{d+1}\setminus\{0\}}}\left\lbrace\mu=(\lambda,\omega)\in \mathcal{O}_{m}^{\gamma }\quad\textnormal{s.t.}\quad\big|\omega\cdot l+jV_{m}(\mu)\big|>\tfrac{\gamma^{\upsilon}\langle j\rangle}{\langle l\rangle^{\tau_{1}}}\right\rbrace.
				\end{equation}
Hence, in the Cantor set $\mathcal{O}_{m+1}^{\gamma },$  the following reduction holds
				$$\mathscr{G}_{m}^{-1}\Big(\omega\cdot\partial_{\varphi}+\partial_{\theta}(V_{m}+f_{m})\Big)\mathscr{G}_{m}=\omega\cdot\partial_{\varphi}+\partial_{\theta}(V_{m+1}+f_{m+1}),$$
				with $V_{m+1}$ and $f_{m+1}$ defined by
				\begin{equation}\label{definition Vm+1 and fm+1}
					\left\lbrace\begin{array}{l}
						V_{m+1}=V_{m}+\langle f_{m}\rangle_{\varphi,\theta}\\
						f_{m+1}=\mathcal{G}_{m}^{-1}\Big(\Pi_{N_{m}}^{\perp}f_{m}+f_{m}\partial_{\theta}g_{m}\Big).
					\end{array}\right.
				\end{equation}
				In view of \eqref{symmetry for f_{m}}, the function $f_m$ is even and therefore $g_m$ is odd. Consequently, we deduce through elementary manipulations that $f_{m+1}$ is also even. This allows us to follow the symmetry persistence along the scheme.
				Besides, in a similar way to \eqref{symmetry for g}, one obtains
				\begin{equation}\label{symmetry for g_{m}}
					g_{m}\in\bigcap_{s\geqslant 0}W^{q,\infty,\gamma}(\mathcal{O},H_{\mbox{\tiny{odd}}}^{s}).
				\end{equation}
				Now, we set 
				$$\mathscr{B}_{-1}=\mathscr{G}_{-1}=\textnormal{Id}\quad\textnormal{ and }\quad \forall m\in\mathbb{N},\,\,\mathscr{B}_{m}=\mathscr{G}_{0}\circ\mathscr{G}_{1}\circ...\circ\mathscr{G}_{m}.$$ 
				One easily finds that
				\begin{align*}
					\mathscr{B}_{m}\rho(\mu,\varphi,\theta)&=\big(1+\partial_{\theta}\beta_{m}(\mu,\varphi,\theta)\big)\mathcal{B}_{m}\rho(\mu,\varphi,\theta)\\
					&=\big(1+\partial_{\theta}\beta_{m}(\mu,\varphi,\theta)\big)\rho\big(\mu,\varphi,\theta+\beta_{m}(\mu,\varphi,\theta)\big),
				\end{align*}
				where the sequence $(\beta_{m})_{m\in\mathbb{N}}$ is defined by $\beta_{-1}=g_{-1}=0$ and  
				\begin{equation}\label{definition betam}
					\beta_{0}=g_{0}\quad\mbox{ and }\quad\beta_{m}(\mu,\varphi,\theta)=\beta_{m-1}(\mu,\varphi,\theta)+g_m\big(\mu,\varphi,\theta+\beta_{m-1}(\mu,\varphi,\theta)\big).
				\end{equation}
				A trivial induction based on \eqref{symmetry for g_{m}} yields
				\begin{equation}\label{symmetry for h_{m}}
					\beta_{m}\in\bigcap_{s\geqslant 0}W^{q,\infty,\gamma}(\mathcal{O},H_{\mbox{\tiny{odd}}}^{s}).
				\end{equation}
				According to  Sobolev embeddings, \eqref{definition Vm+1 and fm+1} and the induction assumption \eqref{hypothesis of induction deltam}, we infer
				\begin{align}\label{Cauchy Vm}
					\| V_{m}-V_{m-1}\|_{q}^{\gamma ,\mathcal{O}}&=\|\langle f_{m-1}\rangle_{\varphi,\theta}\|_{q}^{\gamma ,\mathcal{O}}\nonumber\\
					&\leqslant\| f_{m-1}\|_{q,s_{0}}^{\gamma ,\mathcal{O}}\nonumber\\
					&=\gamma\delta_{m-1}(s_{0})\nonumber\\
					&\leqslant\gamma\delta_{0}(s_{h})N_{0}^{\mu_{2}}N_{m-1}^{-\mu_{2}}.
\end{align}
As a consequence, by using the triangle inequality, \eqref{initial smallness condition in sh norm} and choosing $\varepsilon_{0}$ small enough we deduce
\begin{align*}
					\| V_{m}\|_{q}^{\gamma ,\mathcal{O}} & \leqslant  \| V_{m-1}\|_{q}^{\gamma ,\mathcal{O}}+\gamma\delta_{0}(s_{h})N_{0}^{\mu_{2}}N_{m-1}^{-\mu_{2}}\\
					& \leqslant  \displaystyle\| V_{0}\|_{q}^{\gamma ,\mathcal{O}}+\gamma\delta_{0}(s_h)N_0^{\mu_{2}}\left(\sum_{k=0}^{m-1}N_{k}^{-\mu_{2}}\right)\\
					& \leqslant  \displaystyle\| V_{0}\|_{q}^{\gamma ,\mathcal{O}}+\sum_{k=0}^{\infty}N_{k}^{-\mu_{2}}.
				\end{align*}
				Now, remark that \eqref{param-trans} implies in particular
				$$\mu_{2}\geqslant\tau_1 q+\tau_1+2.$$
				Hence, by the induction hypothesis \eqref{hypothesis of induction deltam}, \eqref{initial smallness condition in sh norm}, \eqref{mu2-1} and \eqref{Est-NO}, we have
				\begin{align}\label{useful estimate in the induction}
					\delta_{m}(s_{0})N_{m}^{\tau_{1} q+\tau_{1}+1}&\leqslant\delta_{0}(s_{h})N_{0}^{\mu_{2}}N_{m}^{\tau_{1} q+\tau_{1}+1-\mu_{2}}\nonumber\\
					&\leqslant\varepsilon_{0}N_{0}^{-1}\nonumber\\
					&\leqslant\varepsilon_{0}.
				\end{align}
				Using \eqref{boundedness of V0} and the previous estimate, we deduce that
				\begin{equation}\label{assumptions KAM iterations checked}
					\sup_{m\in\mathbb{N}}\|V_{m}\|_{q}^{\gamma,\mathcal{O}}\leqslant C\quad\textnormal{and}\quad \delta_{m}(s_{0})N_{m}^{\tau_{1}q+\tau_{1}+1}\leqslant\varepsilon_{0}.
				\end{equation}
				Thus, the KAM step applies and, in particular, the estimate \eqref{estimate KAM step transport} becomes
				\begin{equation}\label{recurrence estimate deltam}
					\delta_{m+1}(s)\leqslant N_{m}^{s-\overline{s}}\delta_{m}(\overline{s})+CN_{m}^{\tau_{1}q+\tau_{1}+1}\delta_{m}(s)\delta_{m}(s_{0}).
				\end{equation}
				If we apply \eqref{recurrence estimate deltam} with $s=s_{l}$ and $\overline{s}=s_{h},$ we obtain 
				$$\delta_{m+1}(s_{l})\leqslant N_{m}^{s_{l}-s_{h}}\delta_{m}(s_{h})+CN_{m}^{\tau_{1} q+\tau_{1}+1}\delta_{m}(s_{l})\delta_{m}(s_{0}).$$
				Using the induction assumption \eqref{hypothesis of induction deltam} and the fact that $s_{l}\geqslant s_{0}$ yields
				\begin{align*}
					\delta_{m+1}(s_{l})&\leqslant N_{m}^{s_{l}-s_{h}}\delta_{m}(s_{h})+CN_{m}^{\tau_{1} q+\tau_{1}+1}(\delta_{m}(s_{l}))^{2}\\
					&\leqslant\left(2-\tfrac{1}{m+1}\right)N_{m}^{s_{l}-s_{h}}\delta_{0}(s_{h})+CN_{0}^{2\mu_{2}}N_{m}^{\tau_{1}q+\tau_{1}+1-2\mu_{2}}(\delta_{0}(s_{h}))^{2}\\
					&\leqslant 2N_{m}^{s_{l}-s_{h}}\delta_{0}(s_{h})+CN_{0}^{2\mu_{2}}N_{m}^{\tau_{1}q+\tau_{1}+1-2\mu_{2}}(\delta_{0}(s_{h}))^{2}.
				\end{align*}
				The conditions \eqref{param-trans} imply
					$$s_h\geqslant \tfrac{3}{2}\mu_{2}+s_l+1,\quad\textnormal{and}\quad  \mu_{2}\geqslant2(\tau_{1}q+\tau_1+1)+1.$$
				Also, using the fact that $N_{0}\geqslant 2$ and choosing $\varepsilon_{0}$ small enough, we get in view of \eqref{initial smallness condition in sh norm},
			$$4N_{0}^{-\mu_{2}}\leqslant 1\quad\textnormal{and}\quad 2C\delta_{0}(s_{h})N_{0}^{\mu_{2}}\leqslant 1.$$
			As a consequence, one has
				\begin{equation}\label{conv-t1}
					N_{m}^{s_{l}-s_{h}}\leqslant\tfrac{1}{4}N_{0}^{\mu_{2}}N_{m+1}^{-\mu_{2}}\quad\textnormal{and}\quad CN_{0}^{2\mu_{2}}N_{m}^{\tau_{1}q+\tau_{1}+1-2\mu_{2}}\delta_{0}(s_{h})\leqslant\tfrac{1}{2}N_{0}^{\mu_{2}}N_{m+1}^{-\mu_{2}},
				\end{equation}
				which implies in turn
				$$\delta_{m+1}(s_{l})\leqslant \delta_{0}(s_{h})N_{0}^{\mu_{2}}N_{m+1}^{-\mu_{2}}.$$
				This proves the first statement of the induction in \eqref{hypothesis of induction deltam} and we now turn to the proof of the second statement. Applying \eqref{recurrence estimate deltam} with $s=\overline{s}=s_h$ and and using the induction \eqref{hypothesis of induction deltam}, we get
				\begin{align*}
					\delta_{m+1}(s_{h})&\leqslant\delta_{m}(s_{h})\left(1+CN_{m}^{\tau_{1} q+\tau_{1}+1}\delta_{m}(s_{0})\right)\\
					&\leqslant\left(2-\tfrac{1}{m+1}\right)\delta_{0}(s_{h})\left(1+CN_{0}^{\mu_{2}}N_{m}^{\tau_{1} q+\tau_{1}+1-\mu_{2}}\delta_{0}(s_{h})\right).
				\end{align*}
				Notice that if the condition
				\begin{equation}\label{conv-t2}
					\left(2-\tfrac{1}{m+1}\right)\left(1+CN_{0}^{\mu_{2}}N_{m}^{\tau_{1} q+\tau_{1}+1-\mu_{2}}\delta_{0}(s_{h})\right)\leqslant 2-\tfrac{1}{m+2}
				\end{equation}
				holds true, then
				$$\delta_{m+1}(s_{h})\leqslant\left(2-\tfrac{1}{m+2}\right)\delta_{0}(s_{h}),
				$$
				which achieves the induction argument of \eqref{hypothesis of induction deltam}. Notice that \eqref{conv-t2} is equivalent to
				$$\left(2-\tfrac{1}{m+1}\right)CN_{0}^{\mu_{2}}N_{m}^{\tau_{1} q+\tau_{1}+1-\mu_{2}}\delta_{0}(s_{h})\leqslant\tfrac{1}{(m+1)(m+2)}\cdot$$
				Using \eqref{mu2-1}, the preceding condition holds true if
				\begin{equation}\label{conv-t3}
					CN_{0}^{\mu_{2}}N_{m}^{-1}\delta_{0}(s_{h})\leqslant\tfrac{1}{(m+1)(m+2)}\cdot
				\end{equation}
				Since $N_0\geqslant 2$, then in view of \eqref{definition of Nm} there exists a small enough constant $c_{0}>0$ such that
				$$\forall m\in\mathbb{N},\quad c_{0}N_{m}^{-1}\leqslant\tfrac{1}{(m+1)(m+2)}\cdot
				$$
				Consequently, \eqref{conv-t3} is ensured provided that
				\begin{equation}\label{conv-t4}
					CN_{0}^{\mu_{2}}\delta_{0}(s_{h})\leqslant c_{0}.
				\end{equation}
				Choosing $\varepsilon_{0}$ small enough and making use of \eqref{initial smallness condition in sh norm}, we obtain 
				\begin{align*}
					CN_{0}^{\mu_{2}}\delta_{0}(s_{h})&\leqslant C\varepsilon_{0}\\
					&\leqslant c_{0}.
				\end{align*}
				Hence, the condition \eqref{conv-t4} is satisfied and  the proof of \eqref{hypothesis of induction deltam} is now achieved.\\
				\ding{226} \textit{Persistence of the regularity.} Putting together \eqref{recurrence estimate deltam}, applied with $\overline{s}=s\in[s_{0},S]$, \eqref{hypothesis of induction deltam} and \eqref{mu2-1}, we infer
				\begin{align*}
					\delta_{m+1}(s) & \leqslant\delta_{m}(s)\left(1+CN_{m}^{\tau_{1} q+\tau_{1}+1}\delta_{m}(s_{0})\right)\\
					& \leqslant\delta_{m}(s)\left(1+C\delta_{0}(s_{h})N_{0}^{\mu_{2}}N_{m}^{\tau_{1} q+\tau_{1}+1-\mu_{2}}\right)\\
					&\leqslant\delta_{m}(s)\left(1+CN_{m}^{-1}\right).
				\end{align*}
				Gathering this estimate with \eqref{estimate delta0 and I0}, implies, up to a trivial induction,
				\begin{align}\label{uniform estimate of deltams}
					\delta_{m}(s)&\leqslant\delta_{0}(s)\prod_{k=0}^{\infty}\left(1+CN_{k}^{-1}\right)\nonumber\\
					&\leqslant C\delta_{0}(s)\\
					&\leqslant C\varepsilon\gamma^{-1}\left(1+\|\mathfrak{I}_{0}\|_{q,s+1}^{\gamma,\mathcal{O}}\right).\nonumber
				\end{align}
				Then, \eqref{control of g by f}, interpolation inequality in Lemma \ref{Lem-lawprod} and \eqref{hypothesis of induction deltam} give
				\begin{align*}
					\|g_{m}\|_{q,s}^{\gamma,\mathcal{O}}&\leqslant C\delta_{m}(s+\tau_{1}q+\tau_{1})\\
					&\leqslant C\left(\delta_{m}(s_{0})\right)^{\overline{\theta}(s)}\big(\delta_{m}(s+\tau_{1}q+\tau_{1}+1)\big)^{1-\overline{\theta}(s)}\\
					&\leqslant C\delta_{0}^{\overline{\theta}(s)}(s_{h})\delta_{0}^{1-\overline{\theta}(s)}(s+\tau_{1}q+\tau_{1}+1)N_{0}^{\overline{\theta}(s)\mu_{2}}N_{m}^{-\overline{\theta}(s)\mu_{2}},
				\end{align*}
				with $\overline{\theta}(s):=\frac{1}{s+\tau_{1}q+\tau_{1}+1-s_{0}}.$
				From \eqref{uniform estimate of deltams}, \eqref{smallness condition transport} and \eqref{estimate delta0 and I0}, we deduce
				\begin{align}\label{bound gm}
					\|g_{m}\|_{q,s}^{\gamma,\mathcal{O}}&\leqslant C\varepsilon\gamma^{-1}\left(1+\|\mathfrak{I}_{0}\|_{q,s_{h}+1}^{\gamma,\mathcal{O}}\right)\left(1+\|\mathfrak{I}_{0}\|_{q,s+\tau_{1}q+\tau_{1}+2}^{\gamma,\mathcal{O}}\right)N_{0}^{\overline{\theta}(s)\mu_{2}}N_{m}^{-\overline{\theta}(s)\mu_{2}}\nonumber\\
					&\leqslant C\varepsilon\gamma^{-1}\left(1+\|\mathfrak{I}_{0}\|_{q,s+\tau_{1}q+\tau_{1}+2}^{\gamma,\mathcal{O}}\right)N_{0}^{\overline{\theta}(s)\mu_{2}}N_{m}^{-\overline{\theta}(s)\mu_{2}}.
				\end{align}
				Using \eqref{definition betam} and \eqref{tame comp}, we get for all $s\in[s_{0},S]$
				\begin{equation}\label{link betam and betam-1}
					\|\beta_{m}\|_{q,s}^{\gamma ,\mathcal{O}}\leqslant\|\beta_{m-1}\|_{q,s}^{\gamma ,\mathcal{O}}\left(1+C\| g_{m}\|_{q,s_{0}}^{\gamma ,\mathcal{O}}\right)+C\left(1+\| \beta_{m-1}\|_{q,s_{0}}^{\gamma,\mathcal{O}}\right)\| g_{m}\|_{q,s}^{\gamma ,\mathcal{O}}.
				\end{equation}
				If we apply this estimate with $s=s_{0}$ and use Sobolev embeddings, we deduce
				$$\|\beta_{m}\|_{q,s_{0}}^{\gamma ,\mathcal{O}}\leqslant\|\beta_{m-1}\|_{q,s_{0}}^{\gamma ,\mathcal{O}}\left(1+C\| g_{m}\|_{q,s_{0}}^{\gamma ,\mathcal{O}}\right)+C\| g_{m}\|_{q,s_{0}}^{\gamma ,\mathcal{O}}.$$
				The previous two expressions make appear recurrent relation for the weighted norms of the sequence $(\beta_m)_m.$ To get good  estimate for $\beta_{m}$, we shall make use the following result which is quite easy to prove by induction : Given three positive sequences $(a_n)_{n\in\mathbb{N}},(b_n)_{n\in\mathbb{N}}$ and $(c_n)_{n\in\mathbb{N}}$ satisfying
				$$
				\forall\, n\in\mathbb{N},\quad a_{n+1}\leqslant b_n a_n+c_n,
				$$
				we have
				\begin{align}\label{Ind-res}
					\nonumber \forall\, n\geqslant 2,\quad a_{n}&\leqslant a_0\prod_{i=0}^{n-1}b_i+\sum_{k=0}^{n-2}c_k\prod_{i=k+1}^{n-1}b_i+c_{n-1}\\
					&\leqslant \Big(a_0+\sum_{k=0}^{n-1} c_k\Big)\prod_{i=0}^{n-1}b_i.
				\end{align}
				In particular, if $\displaystyle \prod_{n=0}^{\infty}b_n$ and $\displaystyle \sum_{n=0}^{\infty} c_n$ converge then
				\begin{align}\label{Ind-res1}
					\sup_{n\in\mathbb{N}}a_{n}&\leqslant \Big(a_0+\sum_{n=0}^{\infty} c_n\Big)\prod_{n=0}^{\infty}b_i.
				\end{align}
				Since the conditions \eqref{param-trans} and \eqref{sigma1} imply
				\begin{equation}\label{p-tmu}
					 s_{0}+\tau_{1}q+\tau_{1}+2\leqslant s_{h}+\sigma_{1}\quad\textnormal{and}\quad \overline{\theta}(s_{0})\mu_{2}\geqslant 1,
				\end{equation}
				then, from \eqref{bound gm} and \eqref{smallness condition transport}, we deduce
				\begin{align*}
					\|g_{m}\|_{q,s_{0}}^{\gamma,\mathcal{O}}&\leqslant C\varepsilon\gamma^{-1}N_{0}^{\mu_{2}}\left(1+\|\mathfrak{I}_{0}\|_{q,s_{0}+\tau_{1}q+\tau_{1}+2}^{\gamma,\mathcal{O}}\right)N_{m}^{-\overline{\theta}(s_{0})\mu_{2}}\\
					&\leqslant C\varepsilon_{0}N_{m}^{-1}.
				\end{align*}
				Choosing $\varepsilon_{0}$ small enough to ensure $C\varepsilon_{0}\leqslant1$, $N_{0}$ sufficiently large to ensure $\displaystyle\sum_{m=0}^{\infty} N_{m}^{-1}<\infty$ and we can apply \eqref{Ind-res1} together with the fact that $\beta_{0}=g_{0}$ to obtain
				\begin{align}\label{unif betam s0}
					\sup_{m\in\mathbb{N}}\|\beta_{m}\|_{q,s_{0}}^{\gamma,\mathcal{O}}&\leqslant\left(\|\beta_{0}\|_{q,s_{0}}^{\gamma,\mathcal{O}}+C\sum_{k=0}^{\infty}\|g_{k}\|_{q,s_{0}}^{\gamma,\mathcal{O}}\right)\prod_{k=0}^{\infty}\left(1+C\|g_{k}\|_{q,s_{0}}^{\gamma,\mathcal{O}}\right)\nonumber\\
					&\leqslant \left(1+C\sum_{k=0}^{\infty}N_{k}^{-1}\right)\prod_{k=0}^{\infty}\left(1+N_{k}^{-1}\right)\nonumber\\
					&\leqslant C.
				\end{align}
				Hence the sequence $\left(\| \beta_{m}\|_{q,s_{0}}^{\gamma,\mathcal{O}}\right)_{m\in\mathbb{N}}$ is bounded and inserting this information in \eqref{link betam and betam-1} gives for all $s\in[s_{0},S]$
				$$\|\beta_{m}\|_{q,s}^{\gamma,\mathcal{O}}\leqslant\|\beta_{m-1}\|_{q,s}^{\gamma,\mathcal{O}}\left(1+C\|g_{m}\|_{q,s_{0}}^{\gamma,\mathcal{O}}\right)+C\|g_{m}\|_{q,s}^{\gamma,\mathcal{O}}.$$
				Similarly to what preceeds, if we apply \eqref{Ind-res1} and \eqref{bound gm}, we infer 
				\begin{align*}
					\sup_{m\in\mathbb{N}}\|\beta_{m}\|_{q,s}^{\gamma,\mathcal{O}}&\leqslant\left(\|\beta_{0}\|_{q,s}^{\gamma,\mathcal{O}}+C\sum_{k=0}^{\infty}\|g_{k}\|_{q,s}^{\gamma,\mathcal{O}}\right)\prod_{k=0}^{\infty}\left(1+C\|g_{k}\|_{q,s_{0}}^{\gamma,\mathcal{O}}\right)\\
					&\leqslant C\varepsilon\gamma^{-1}\left(1+\|\mathfrak{I}_{0}\|_{q,s+\tau_{1}q+\tau_{1}+2}^{\gamma,\mathcal{O}}\right)\left(1+N_{0}^{\overline{\theta}(s)\mu_{2}}\sum_{k=0}^{\infty}N_{k}^{-\overline{\theta}(s)\mu_{2}}\right).
				\end{align*}
				From Lemma \ref{lemma sum Nn} we get
				$$\forall s\in[s_{0},S],\quad N_{0}^{\overline{\theta}(s)\mu_{2}}\sum_{k=0}^{\infty}N_{k}^{-\overline{\theta}(s)\mu_{2}}\lesssim 1$$
				which implies in turn
				\begin{equation}\label{uniform estimate betam}
					\forall s\in[s_{0},S],\quad\sup_{m\in\mathbb{N}}\| \beta_{m}\|_{q,s}^{\gamma,\mathcal{O}}\leqslant C\varepsilon\gamma^{-1}\left(1+\|\mathfrak{I}_{0}\|_{q,s+\tau_{1}q+\tau_{1}+2}^{\gamma,\mathcal{O}}\right).
				\end{equation}
			From the condition \eqref{param} we have
			$s_{l}= s_{0}+\tau_{1}q+\tau_{1}+2,$
			and consequently we deduce from  \eqref{tame comp}, \eqref{unif betam s0}, \eqref{control of g by f} and \eqref{hypothesis of induction deltam}, 
				\begin{align}\label{decay telescopic beta}
					\|\beta_{m}-\beta_{m-1}\|_{q,s_{0}+2}^{\gamma,\mathcal{O}}&\leqslant C\|g_{m}\|_{q,s_{0}+2}^{\gamma,\mathcal{O}}\left(1+\|\beta_{m-1}\|_{q,s_{0}+2}^{\gamma,\mathcal{O}}\right)\nonumber\\
					&\leqslant C\|g_{m}\|_{q,s_{0}+2}^{\gamma,\mathcal{O}}\leqslant C\delta_{m}(s_{l})\nonumber\\
					&\leqslant CN_{0}^{\mu_{2}}N_{m}^{-\mu_{2}}\delta_{0}(s_{h}).
				\end{align}
				Applying once again Lemma \ref{lemma sum Nn}, we deduce that 
				$$\sum_{m=0}^{\infty}\|\beta_{m}-\beta_{m-1}\|_{q,s_{0}+2}^{\gamma,\mathcal{O}}\leqslant C\delta_{0}(s_{h}).$$
				Hence there exists $\beta\in W^{q,\infty,\gamma}(\mathcal{O},H^{s_{0}+2})$ such that
				$$\beta_{m}\underset{m\rightarrow\infty}{\longrightarrow}\beta\quad\textnormal{(strongly) in}\quad W^{q,\infty,\gamma}(\mathcal{O},H^{s_{0}+2}).$$
				%By \eqref{difference composition CVAR}, we have :
				%$$\| h_{m}-h_{m-1}\|_{q,s}^{\gamma ,\mathcal{O}}\lesssim_{s}\| \beta_{m}\|_{q,s+s_{0}}^{\gamma ,\mathcal{O}}\left(1+\| h_{m-1}\|_{q,s_{0}}^{\gamma ,\mathcal{O}}\right)+\| \beta_{m}\|_{q,s_{0}}^{\gamma ,\mathcal{O}}\| h_{m-1}\|_{q,s+1}^{\gamma ,\mathcal{O}}.$$
				%So by \eqref{estimate on gm}, we have :
				%\begin{equation}\label{Cauchy hm}
				%\| h_{m}-h_{m-1}\|_{q,s}^{\gamma ,\mathcal{O}}\lesssim_{s} (\texttt{C}(s+1))^{(1-t)m}N_{0}^{t\mu_{2}}N_{m}^{-t\mu_{2}}M(s+s_{0}+\tau_{1}q+\tau_{1}+2).
				%\end{equation}
				By \eqref{uniform estimate betam} the sequence $\left(\beta_{m}\right)_{m\in\mathbb{N}}$ is bounded in $W^{q,\infty,\gamma}(\mathcal{O},H^{s})$, then by a weak-compactness argument we find that $\beta\in W^{q,\infty,\gamma}(\mathcal{O},H^{s})$.
			Using  \eqref{uniform estimate betam}, we obtain 
				\begin{align}\label{estimate beta}
					\forall s\in[s_{0},S],\quad\|\beta\|_{q,s}^{\gamma,\mathcal{O}}&\leqslant\liminf_{m\rightarrow\infty}\|\beta_{m}\|_{q,s}^{\gamma,\mathcal{O}}\nonumber\\
					&\lesssim\varepsilon\gamma^{-1}\left(1+\|\mathfrak{I}_{0}\|_{q,s+\tau_{1}q+\tau_{1}+2}^{\gamma,\mathcal{O}}\right).
				\end{align}
				We then can consider the quasi-periodic symplectic change of variables $\mathscr{B}$ associated with $\beta$ and defined by
				\begin{align*}
					\mathscr{B}\rho(\lambda,\omega,\varphi,\theta)&=\left(1+\partial_{\theta}\beta(\lambda,\omega,\varphi,\theta)\right)\mathcal{B}\rho(\lambda,\omega,\varphi,\theta)\\
					&=\left(1+\partial_{\theta}\beta(\lambda,\omega,\varphi,\theta)\right)\rho(\lambda,\omega,\varphi,\theta+\beta(\lambda,\omega,\varphi,\theta)).
				\end{align*}
				By \eqref{estimate beta}, \eqref{p-tmu} and \eqref{smallness condition transport}, we have 
				\begin{equation}\label{beta s0 bound}
					\|\beta\|_{q,s_{0}}^{\gamma,\mathcal{O}}\lesssim\varepsilon\gamma^{-1}\left(1+\|\mathfrak{I}_{0}\|_{q,s_{0}+\tau_{1}q+\tau_{1}+2}^{\gamma,\mathcal{O}}\right)\lesssim\varepsilon_{0}.
				\end{equation}
			Proceeding as for \eqref{e-g-2s_0}, using interpolation \eqref{estimate beta}, \eqref{beta s0 bound}, \eqref{smallness condition transport} and the fact that $2s_0+\tau_{1}q+\tau_{1}+3\leqslant s_h+\sigma_1$, one obtains
			$$\|\beta\|_{q,2s_0}^{\gamma,\mathcal{O}}\lesssim\varepsilon_0.$$
				Therefore, choosing $\varepsilon_{0}$ small enough, we deduce in view of Lemma \ref{Compos1-lemm} that $\mathscr{B}$ is an invertible operator. Moreover, by \eqref{tame comp symp} and \eqref{estimate beta}, we get 
				\begin{equation}\label{est Bpm1}
					\|\mathscr{B}^{\pm 1}\rho\|_{q,s}^{\gamma,\mathcal{O}}\lesssim\|\rho\|_{q,s}^{\gamma,\mathcal{O}}+\varepsilon\gamma^{-1}\|\mathfrak{I}_{0}\|_{q,s+\tau_{1}q+\tau_{1}+3}^{\gamma,\mathcal{O}}\|\rho\|_{q,s_{0}}^{\gamma,\mathcal{O}}.
				\end{equation}
				In addition, by \eqref{symmetry for h_{m}}, and Sobolev embeddings (to get pointwise convergence), we find 
				$$\beta\in \bigcap_{s\in[s_{0},S]}W^{q,\infty,\gamma}(\mathcal{O},H_{\mbox{\tiny{odd}}}^{s}).
				$$
				We also have an estimate of the rate of  convergence for the sequence $(\beta_{m})_m$ towards $\beta$,
				\begin{align}\label{difference beta and betam}
					\|\beta-\beta_{m}\|_{q,s_{0}+2}^{\gamma,\mathcal{O}}&\leqslant\sum_{k=m}^{\infty}\|\beta_{k+1}-\beta_{k}\|_{q,s_{0}+2}^{\gamma,\mathcal{O}}\nonumber\\
					&\lesssim\gamma\delta_{0}(s_{h})N_{0}^{\mu_{2}}\sum_{k=m+1}^{\infty}N_{k}^{-\mu_{2}}.
				\end{align}
				From Lemma \ref{lemma sum Nn}, one obtains
				\begin{equation}\label{big O series Nm}
					\sum_{k=m}^{\infty}N_{k}^{-\mu_{2}}\underset{m\rightarrow\infty}{=}O\left(N_{m}^{-\mu_{2}}\right).
				\end{equation}
				Gathering \eqref{big O series Nm}, \eqref{difference beta and betam} and \eqref{estimate delta0 and I0}, we get
				\begin{align}\label{rate on convergence betam to beta}
					\nonumber \|\beta-\beta_{m}\|_{q,s_{0}+2}^{\gamma,\mathcal{O}}&\lesssim\gamma\delta_{0}(s_{h})N_{0}^{\mu_{2}}N_{m+1}^{-\mu_{2}}\\
					&\lesssim\varepsilon N_{0}^{\mu_{2}}N_{m+1}^{-\mu_{2}}.
				\end{align}
				\noindent $\blacktriangleright$ \textbf{KAM conclusion}\\
				By \eqref{Cauchy Vm}, we have 
				\begin{align*}
					\sum_{m=0}^{\infty}\| V_{m+1}-V_{m}\|_{q}^{\gamma ,\mathcal{O}}&\leqslant\gamma\delta_{0}(s_{h})N_{0}^{\mu_{2}}\sum_{m=0}^{\infty}N_{m}^{-\mu_{2}}\\
					&\lesssim \gamma\delta_{0}(s_{h}).
				\end{align*}
				We deduce that the sequence $(V_{m})_{m\in\mathbb{N}}$ is convergent in $W^{q,\infty,\gamma }(\mathcal{O},\mathbb{C})$ and let us   denote by  $c_{i_{0}}$ its limit. Moreover, we have by \eqref{big O series Nm}, \eqref{estimate delta0 and I0} and \eqref{smallness condition transport}
				\begin{align*}
					\| c_{i_{0}}-V_{0}\|_{q}^{\gamma ,\mathcal{O}}&\leqslant\sum_{m=0}^{\infty}\| V_{m+1}-V_{m}\|_{q}^{\gamma ,\mathcal{O}}\\
					&\lesssim\gamma \delta_{0}(s_{h})
					\lesssim \varepsilon\left(1+\|\mathfrak{I}_{0}\|_{q,s_{h}+1}^{\gamma,\mathcal{O}}\right)\\
					&\lesssim\varepsilon.
				\end{align*}
				Now, we introduce the truncated Cantor set
				$$\mathcal{O}_{\infty,n}^{\gamma,\tau_{1}}(i_{0})=\bigcap_{(l,j)\in\mathbb{Z}^{d}\times\mathbb{Z}\setminus\{(0,0\}\atop|l|\leqslant N_{n}}\left\lbrace\mu:=(\lambda,\omega)\in\mathcal{O}\quad\textnormal{s.t.}\quad\big|\omega\cdot l+jc_{i_{0}}(\mu)\big|>\tfrac{4\gamma^{\upsilon}\langle j\rangle}{\langle l\rangle^{\tau_{1}}}\right\rbrace.$$
				In what follows, we shall prove that the Cantor set $\mathcal{O}_{\infty,n}^{\gamma,\tau_{1}}(i_{0})$ satisfies the inclusion
				$$\mathcal{O}_{\infty,n}^{\gamma,\tau_{1}}(i_{0})\subset\bigcap_{m=0}^{n+1}\mathcal{O}_{m}^{\gamma}=\mathcal{O}_{n+1}^{\gamma},
				$$
				where the intermediate Cantor sets are defined in \eqref{definition mathcal Om+1}. For this aim, we shall argue by induction. We first remark that by construction $\mathcal{O}_{\infty,n}^{\gamma,\tau_{1}}(i_{0})\subset\mathcal{O}=:\mathcal{O}_{0}^{\gamma}.$ Now assume that $\mathcal{O}_{\infty,n}^{\gamma,\tau_{1}}(i_{0})\subset\mathcal{O}_{m}^{\gamma}$ for $m\leqslant n$ and let us check that
				\begin{equation}\label{inclusion Oinftyn in Om+1}
					\mathcal{O}_{\infty,n}^{\gamma,\tau_{1}}(i_{0})\subset\mathcal{O}_{m+1}^{\gamma}.
				\end{equation}
				Putting together  \eqref{Cauchy Vm} and \eqref{big O series Nm} we infer
				\begin{align}\label{estimate Vm-ci0}
				\| V_{m}-c_{i_{0}}\|_{q}^{\gamma ,\mathcal{O}}
					&\leqslant\sum_{l=m}^{\infty}\| V_{l+1}-V_{l}\|_{q}^{\gamma ,\mathcal{O}}\nonumber\\
					&\leqslant\gamma\delta_{0}(s_{h})N_{0}^{\mu_{2}}\sum_{l=m}^{\infty}N_{l}^{-\mu_{2}}\nonumber\\
					&\lesssim\gamma\delta_{0}(s_{h})N_{0}^{\mu_{2}}N_{m}^{-\mu_{2}}.
				\end{align}
				Given $\mu\in\mathcal{O}_{\infty,n}^{\gamma,\tau_{1}}(i_{0})$ and  $(l,j)\in\mathbb{Z}^{d}\times\mathbb{Z}\setminus\{(0,0)\}$ such that $0\leqslant|l|\leqslant N_{m},$ we have then  $|l|\leqslant N_{n}$ and by triangle inequality,
				\begin{align*}
					|\omega\cdot l+jV_{m}(\mu)| & \geqslant |\omega\cdot l+jc_{i_{0}}(\mu)|-|j||V_{m}(\mu)-c_{i_{0}}(\mu)|\\
					&\geqslant\displaystyle\tfrac{4\gamma^{\upsilon}\langle j\rangle}{\langle l\rangle^{\tau_{1}}}-C\langle j\rangle\gamma\delta_{0}(s_{h})N_{0}^{\mu_{2}}N_{m}^{-\mu_{2}}\\
					&\geqslant\displaystyle\tfrac{4\gamma^{\upsilon}\langle j\rangle}{\langle l\rangle^{\tau_{1}}}-C\langle j\rangle\gamma^{\upsilon}\varepsilon_{0}\langle l\rangle^{-\mu_{2}}.
				\end{align*}
				Since  \eqref{param-trans} implies $\mu_{2}\geqslant\tau_{1}$, then taking $\varepsilon_{0}\leqslant\frac{1}{C}$, we deduce from the previous estimate
				$$|\omega\cdot l+jV_{m}(\mu)|> \displaystyle\tfrac{\gamma^{\upsilon}\langle j\rangle}{\langle l\rangle^{\tau_{1}}}.$$
				Consequently,  $\mu\in\mathcal{O}_{m+1}^{\gamma}$ and the inclusion \eqref{inclusion Oinftyn in Om+1} holds.\\
				\textbf{(iii)} We can write for all $n\in\mathbb{N},$
				\begin{align*}
					\mathscr{B}^{-1}\Big(\omega\cdot\partial_{\varphi}+\partial_{\theta}\left(V_{0}+f_{0}\right)\Big)\mathscr{B}&=\left(\mathscr{B}^{-1}-\mathscr{B}_{n}^{-1}\right)\Big(\omega\cdot\partial_{\varphi}+\partial_{\theta}\left(V_{0}+f_{0}\right)\Big)\mathscr{B}\\
					&+\mathscr{B}_{n}^{-1}\Big(\omega\cdot\partial_{\varphi}+\partial_{\theta}\left(V_{0}+f_{0}\right)\Big)\left(\mathscr{B}-\mathscr{B}_{n}\right)\\
					&+\mathscr{B}_{n}^{-1}\Big(\omega\cdot\partial_{\varphi}+\partial_{\theta}\left(V_{0}+f_{0}\right)\Big)\mathscr{B}_{n}.
				\end{align*}
				In view of \eqref{inclusion Oinftyn in Om+1} and the definition of $\mathscr{B}_n$, we have in the Cantor set $\mathcal{O}_{\infty,n}^{\gamma,\tau_{1}}(i_{0})$ 
				$$\mathscr{B}_{n}^{-1}\Big(\omega\cdot\partial_{\varphi}+\partial_{\theta}\left(V_{0}+f_{0}\right)\Big)\mathscr{B}_{n}=\omega\cdot\partial_{\varphi}+\partial_{\theta}\left(V_{n+1}+f_{n+1}\right).$$
				Therefore, in the Cantor set $\mathcal{O}_{\infty,n}^{\gamma,\tau_{1}}(i_{0}),$ the following decomposition holds 
				$$\mathscr{B}^{-1}\Big(\omega\cdot\partial_{\varphi}+\partial_{\theta}\left(V_{0}+f_{0}\right)\Big)\mathscr{B}=\omega\cdot\partial_{\varphi}+c_{i_{0}}\partial_{\theta}+\mathtt{E}_{n}^{0}(i_{0}),$$
				where
				\begin{align*}
					\mathtt{E}_{n}^{0}(i_{0})&:=\left(V_{n+1}-c_{i_{0}}\right)\partial_{\theta}+\partial_{\theta}\left(f_{n+1}\cdot\right)+\left(\mathscr{B}^{-1}-\mathscr{B}_{n}^{-1}\right)\Big(\omega\cdot\partial_{\varphi}+\partial_{\theta}\left(V_{0}+f_{0}\right)\Big)\mathscr{B}\\
					&\quad+\mathscr{B}^{-1}\Big(\omega\cdot\partial_{\varphi}+\partial_{\theta}\left(V_{0}+f_{0}\right)\Big)\left(\mathscr{B}-\mathscr{B}_{n}\right)\\
					&:=\mathtt{E}_{n,1}^{0}+\mathtt{E}_{n,2}^{0}+\mathtt{E}_{n,3}^{0}+\mathtt{E}_{n,4}^{0}.
				\end{align*}
				By the law products in Lemma \ref{Lem-lawprod},  \eqref{estimate Vm-ci0} and \eqref{estimate delta0 and I0} we have 
								\begin{align}\label{En11 s}
					\|\mathtt{E}_{n,1}^{0}\rho\|_{q,s_0}^{\gamma,\mathcal{O}}&\lesssim\|V_{n+1}-c_{i_{0}}\|_{q}^{\gamma,\mathcal{O}}\|\rho\|_{q,s_0+1}^{\gamma,\mathcal{O}}\nonumber\\
					&\lesssim\gamma\delta_{0}(s_{h})N_{0}^{\mu_{2}}N_{n+1}^{-\mu_{2}}\|\rho\|_{q,s_0+1}^{\gamma,\mathcal{O}}\nonumber\\
					&\lesssim\varepsilon N_{0}^{\mu_{2}}N_{n+1}^{-\mu_{2}}\|\rho\|_{q,s_0+1}^{\gamma,\mathcal{O}}.
				\end{align}
%				and by \eqref{hypothesis of induction deltam} and \eqref{uniform estimate of deltams}
%				\begin{align}\label{En12 s}
%					\|\mathtt{E}_{n,2}^{0}\rho\|_{q,s}^{\gamma,\mathcal{O}}&=\|\partial_{\theta}\left(f_{n+1}\rho\right)\|_{q,s}^{\gamma,\mathcal{O}}\nonumber\\
%					&\lesssim\|f_{n+1}\|_{q,s+1}^{\gamma,\mathcal{O}}\|\rho\|_{q,s_{0}+1}^{\gamma,\mathcal{O}}+\|f_{n+1}\|_{q,s_{0}+1}^{\gamma,\mathcal{O}}\|\rho\|_{q,s+1}^{\gamma,\mathcal{O}}\nonumber\\
%					&\lesssim\|\rho\|_{q,s+1}^{\gamma,\mathcal{O}}+\varepsilon\gamma^{-1}\|\mathfrak{I}_{0}\|_{q,s+2}^{\gamma,\mathcal{O}}\|\rho\|_{q,s_{0}+1}^{\gamma,\mathcal{O}}.
%				\end{align}
				From \eqref{hypothesis of induction deltam} and since \eqref{param} implies in particular $s_{l}\geqslant s_{0}+1,$ we obtain
				\begin{align}\label{En12 s0}
					\|\mathtt{E}_{n,2}^{0}\rho\|_{q,s_{0}}^{\gamma,\mathcal{O}}&\lesssim\gamma\delta_{n+1}(s_{0}+1)\|\rho\|_{q,s_{0}+1}^{\gamma,\mathcal{O}}\nonumber\\
					&\lesssim\gamma\delta_{0}(s_{h})N_{0}^{\mu_{2}}N_{n+1}^{-\mu_{2}}\|\rho\|_{q,s_{0}+1}^{\gamma,\mathcal{O}}\nonumber\\
					&\lesssim\varepsilon  N_{0}^{\mu_{2}}N_{n+1}^{-\mu_{2}}\|\rho\|_{q,s_0+1}^{\gamma,\mathcal{O}}.
				\end{align}
				We now turn to the estimate of  $\mathtt{E}_{n,4}^{0}.$ First remark that by the law products in Lemma \ref{Lem-lawprod}, we have
				\begin{align*}
					\|\omega\cdot\partial_{\varphi}\rho+\partial_{\theta}\left(V_{\varepsilon r}\rho\right)\|_{q,s_0}^{\gamma,\mathcal{O}}&\leqslant\|\omega\cdot\partial_{\varphi}\rho\|_{q,s_0}^{\gamma,\mathcal{O}}+\|\partial_{\theta}\left(V_{\varepsilon r}\rho\right)\|_{q,s_0}^{\gamma,\mathcal{O}}\\
					&\lesssim\|\rho\|_{q,s_0+1}^{\gamma,\mathcal{O}}\left(1+\|V_{\varepsilon r}\|_{q,s_{0}+1}^{\gamma,\mathcal{O}}\right).
				\end{align*}
				But combining \eqref{V=V0+f0}, \eqref{boundedness of V0},  \eqref{uniform estimate of deltams} and \eqref{smallness condition transport}, we obtain
				\begin{align*}
					\|V_{\varepsilon r}\|_{q,s_0}^{\gamma,\mathcal{O}}&\leqslant\|V_{0}\|_{q}^{\gamma,\mathcal{O}}+\|f_{0}\|_{q,s_0}^{\gamma,\mathcal{O}}\\
					&\leqslant C+C\varepsilon\gamma^{-1}\|\mathfrak{I}_{0}\|_{q,s_0+1}^{\gamma,\mathcal{O}}\\
					&\leqslant C.
				\end{align*}
				Therefore, we get
				\begin{equation}\label{estimate initial transport operator}
					\|\omega\cdot\partial_{\varphi}\rho+\partial_{\theta}\left(V_{\varepsilon r}\rho\right)\|_{q,s_0}^{\gamma,\mathcal{O}}\lesssim\|\rho\|_{q,s_0+1}^{\gamma,\mathcal{O}}.				
					\end{equation}
					Putting together \eqref{estimate initial transport operator}, \eqref{est Bpm1} and \eqref{smallness condition transport}, gives
				\begin{align}\label{En14 sM0}
					\|\mathtt{E}_{n,4}^{0}\rho\|_{q,s_{0}}^{\gamma,\mathcal{O}}&\lesssim\|(\mathscr{B}-\mathscr{B}_{n})\rho\|_{q,s_{0}+1}^{\gamma,\mathcal{O}}.
									\end{align}
				Applying Taylor Formula, we may write				\begin{align*}
					(\mathscr{B}-\mathscr{B}_{n})\rho(\theta)&=(1+\partial_{\theta}\beta(\theta))\rho(\theta+\beta(\theta))-(1+\partial_{\theta}\beta_{n}(\theta))\rho(\theta+h_{n}(\theta))\\
					&=(1+\partial_{\theta}\beta(\theta))\left[\rho(\theta+\beta(\theta))-\rho(\theta+\beta_{n}(\theta))\right]+\partial_{\theta}\left(\beta-\beta_{n}\right)(\theta)\rho(\theta+\beta_{n}(\theta))\\
					&:=(1+\partial_\theta\beta(\theta))(\beta-\beta_{n})(\theta)\mathscr{I}_{n}(\theta)+\partial_{\theta}(\beta-\beta_{n})(\theta)\mathcal{B}_{n}\rho(\theta),
				\end{align*}
				where
				$$\mathscr{I}_{n}\rho(\theta):=\int_{0}^{1}(\partial_{\theta}\rho)\left(\theta+\beta_{n}(\theta)+t(\beta(\theta)-\beta_{n}(\theta))\right)dt.$$
%				By \ref{tame comp}, we have 
%				\begin{align*}
%					\|\partial_{\theta}\left(\beta-\beta_{n}\right)\mathcal{B}_{n}\rho\|_{q,s}^{\gamma,\mathcal{O}}&\lesssim\|\beta-\beta_{n}\|_{q,s+1}^{\gamma,\mathcal{O}}\|\rho\|_{q,s_{0}}^{\gamma,\mathcal{O}}\left(1+\|\beta_{n}\|_{q,s_{0}}^{\gamma,\mathcal{O}}\right)\\
%					&+\|\beta-\beta_{n}\|_{q,s_{0}+1}^{\gamma,\mathcal{O}}\|\rho\|_{q,s}^{\gamma,\mathcal{O}}\left(1+\|\beta_{n}\|_{q,s_{0}}^{\gamma,\mathcal{O}}\right)\\
%					&+\|\beta-\beta_{n}\|_{q,s_{0}+1}^{\gamma,\mathcal{O}}\|\rho\|_{q,s_{0}}^{\gamma,\mathcal{O}}\|\beta_{n}\|_{q,s}^{\gamma,\mathcal{O}}.
%				\end{align*}
%				Notice that from \eqref{param-trans}, one has $s_{l}\geqslant s_{0}+\tau_{1}q+\tau_{1}+1.$ Then one obtains the estimate \eqref{decay telescopic beta} with $s_0+1$ in place of $s_0.$ As a consequence, the estimates \eqref{rate on convergence betam to beta} and \eqref{beta s0 bound} also hold with $s_{0}+1$ in the place of $s_{0}.$ 
Hence, we get by  the law products, \ref{tame comp} and \eqref{rate on convergence betam to beta}
				\begin{align*}
					\|\partial_{\theta}\left(\beta-\beta_{n}\right)\mathcal{B}_{n}\rho\|_{q,s_{0}+1}^{\gamma,\mathcal{O}}&\lesssim\|\beta-\beta_{n}\|_{q,s_{0}+2}^{\gamma,\mathcal{O}}\|\rho\|_{q,s_{0}+1}\left(1+\|\beta_{n}\|_{q,s_{0}+1}^{\gamma,\mathcal{O}}\right)\\
					&\lesssim\varepsilon N_{0}^{\mu_{2}}N_{n+1}^{-\mu_{2}}\|\rho\|_{q,s_{0}+1}^{\gamma,\mathcal{O}}.
				\end{align*}
				Using the law products together with  \eqref{rate on convergence betam to beta}, \eqref{unif betam s0} and \eqref{beta s0 bound} we find
				\begin{align*}
					\Big\|(1+\partial_{\theta}\beta)\left(\beta-\beta_{n}\right)\mathscr{I}_{n}\rho\Big\|_{q,s_{0}+1}^{\gamma,\mathcal{O}}&\lesssim\left(1+\|\beta\|_{q,s_{0}+2}^{\gamma,\mathcal{O}}\right)\|\beta-\beta_{n}\|_{q,s_{0}+1}^{\gamma,\mathcal{O}}\|\rho\|_{q,s_{0}+2}^{\gamma,\mathcal{O}}\\
					&\qquad \times \left(1+\|\beta_{n}\|_{q,s_{0}+1}^{\gamma,\mathcal{O}}+\|\beta-\beta_{n}\|_{q,s_{0}+1}^{\gamma,\mathcal{O}}\right)\\
					&\lesssim\varepsilon N_{0}^{\mu_{2}}N_{n+1}^{-\mu_{2}}\|\rho\|_{q,s_{0}+2}^{\gamma,\mathcal{O}}.
				\end{align*}
%				For a general $s\in[s_{0},S],$ we deduce  from \eqref{estimate beta}, \eqref{uniform estimate betam} and \eqref{smallness condition transport}
%				\begin{align*}
%					\Big\|(1+\partial_{\theta}\beta)\left(\beta-\beta_{n}\right)\mathscr{I}_{n}\rho\Big\|_{q,s}^{\gamma,\mathcal{O}}\lesssim\|\rho\|_{q,s+1}^{\gamma,\mathcal{O}}+\varepsilon\gamma^{-1}\|\mathfrak{I}_{0}\|_{q,s+\tau_{1}q+\tau_{1}+3}^{\gamma,\mathcal{O}}\|\rho\|_{q,s_{0}+1}^{\gamma,\mathcal{O}}.
%				\end{align*}
				Gathering the foregoing estimates leads to
%				\begin{align}\label{est diff B and Bn}
%					\|(\mathscr{B}-\mathscr{B}_{n})\rho\|_{q,s}^{\gamma,\mathcal{O}}&\lesssim\|\rho\|_{q,s+1}^{\gamma,\mathcal{O}}+\varepsilon\gamma^{-1}\|\mathfrak{I}_{0}\|_{q,s+\tau_{1}q+\tau_{1}+3}^{\gamma,\mathcal{O}}\|\rho\|_{q,s_{0}+1}^{\gamma,\mathcal{O}}
%				\end{align}
%				and
				\begin{align}\label{dbns0}
					\|(\mathscr{B}-\mathscr{B}_{n})\rho\|_{q,s_{0}+1}^{\gamma,\mathcal{O}}&\lesssim\varepsilon N_{0}^{\mu_{2}}N_{n+1}^{-\mu_{2}}\|\rho\|_{q,s_{0}+2}^{\gamma,\mathcal{O}}.
				\end{align}
				Plugging \eqref{dbns0} into \eqref{En14 sM0}  gives
				\begin{align}\label{En14 s0}
					\|\mathtt{E}_{n,4}^{0}\rho\|_{q,s_{0}}^{\gamma,\mathcal{O}}&\lesssim \varepsilon N_{0}^{\mu_{2}}N_{n+1}^{-\mu_{2}}\|\rho\|_{q,s_{0}+2}^{\gamma,\mathcal{O}}.
									\end{align}
%				and by Lemma \ref{Compos1-lemm}, \eqref{uniform estimate betam} and \eqref{estimate beta}, we infer 
%				\begin{equation}\label{En14 s}
%					\|\mathtt{E}_{n,4}^{0}\rho\|_{q,s}^{\gamma,\mathcal{O}}\lesssim\|\rho\|_{q,s+2}^{\gamma,\mathcal{O}}+\varepsilon\gamma^{-1}\|\mathfrak{I}_{0}\|_{q,s+\tau_{1}q+\tau_{1}+4}^{\gamma,\mathcal{O}}\|\rho\|_{q,s_{0}+2}^{\gamma,\mathcal{O}}.
%				\end{equation}
				Proceeding in a similar way as before using in particular the identity \eqref{mathscrB1} and \eqref{link betah and beta} we find
				\begin{align}\label{En13 s0}
					\|\mathtt{E}_{n,3}^{0}\rho\|_{q,s_{0}}^{\gamma,\mathcal{O}}&\lesssim \varepsilon N_{0}^{\mu_{2}}N_{n+1}^{-\mu_{2}}\|\rho\|_{q,s_{0}+2}^{\gamma,\mathcal{O}}.
									\end{align}
				Putting together  \eqref{En11 s}, \eqref{En12 s0}, \eqref{En14 s0}, \eqref{En13 s0} allows to get
				$$\|\mathtt{E}_{n}^{0}\rho\|_{q,s_{0}}^{\gamma,\mathcal{O}}\lesssim\varepsilon N_{0}^{\mu_{2}}N_{n+1}^{-\mu_{2}}\|\rho\|_{q,s_{0}+2}^{\gamma,\mathcal{O}}.
				$$
				\textbf{(iv)} $\blacktriangleright$ \textbf{Estimate of $\Delta_{12}\beta$.}
				First notice that, since $\beta_{-1}=0$, then
				\begin{equation}\label{series Delta12beta}
					\Delta_{12}\beta=\sum_{m=0}^{\infty}\Delta_{12}(\beta_{m}-\beta_{m-1}).
				\end{equation}
				The triangle inequality allows us to write
				\begin{equation}\label{estimate series Delta12 beta}
					\|\Delta_{12}\beta\|_{q,\overline{s}_{h}+\mathtt{p}}^{\gamma,\mathcal{O}}\leqslant\sum_{m=0}^{\infty}\|\Delta_{12}(\beta_{m}-\beta_{m-1})\|_{q,\overline{s}_{h}+\mathtt{p}}^{\gamma,\mathcal{O}}.
				\end{equation}
				According to Taylor Formula and \eqref{definition betam}, we infer		\begin{align*}
					\Delta_{12}\beta_{m}(\theta)&=\Delta_{12}\beta_{m-1}(\theta)+(\mathcal{B}_{m-1})_{r_{1}}(\Delta_{12}g_{m})(\theta)\\
					&\quad+\Delta_{12}\beta_{m-1}(\theta)\int_{0}^{1}(\partial_{\theta}(g_{m})_{r_{2}})\big(\theta+(\beta_{m-1})_{r_{2}}(\theta)+t\Delta_{12}\beta_{m-1}(\theta)\big)dt.
				\end{align*}
				Thus,
				\begin{align*}
					\Delta_{12}(\beta_{m}-\beta_{m-1})(\theta)&=(\mathcal{B}_{m-1})_{r_{1}}(\Delta_{12}g_{m})(\theta)\\
					&\quad+\Delta_{12}\beta_{m-1}(\theta)\int_{0}^{1}(\partial_{\theta}(g_{m})_{r_{2}})\big(\theta+(\beta_{m-1})_{r_{2}}(\theta)+t\Delta_{12}\beta_{m-1}(\theta)\big)dt.
				\end{align*}
				Consequently, using the law product in Lemma \ref{Lem-lawprod}, Lemma \ref{Compos1-lemm} and Sobolev embeddings we obtain
				\begin{align*}
					\|\Delta_{12}(\beta_{m}-\beta_{m-1})\|_{q,\overline{s}_{h}+\mathtt{p}}^{\gamma,\mathcal{O}}&\leqslant\|\Delta_{12}g_{m}\|_{q,\overline{s}_{h}+\mathtt{p}}^{\gamma,\mathcal{O}}\left(1+C\|(\beta_{m-1})_{r_{1}}\|_{q,s_{0}}^{\gamma,\mathcal{O}}\right)+\|\Delta_{12}g_{m}\|_{q,s_{0}}^{\gamma,\mathcal{O}}\|(\beta_{m-1})_{r_{1}}\|_{q,\overline{s}_{h}+\mathtt{p}}^{\gamma,\mathcal{O}}\\
					&+C\|\Delta_{12}\beta_{m-1}\|_{q,s_{0}}^{\gamma,\mathcal{O}}\|(g_{m})_{r_{2}}\|_{q,\overline{s}_{h}+\mathtt{p}+1}^{\gamma,\mathcal{O}}\left(1+\|(\beta_{m-1})_{r_{2}}\|_{q,s_{0}}^{\gamma,\mathcal{O}}+\|\Delta_{12}\beta_{m-1}\|_{q,s_{0}}^{\gamma,\mathcal{O}}\right)\\
					&+C\|\Delta_{12}\beta_{m-1}\|_{q,s_{0}}^{\gamma,\mathcal{O}}\|(g_{m})_{r_{2}}\|_{q,s_{0}+1}^{\gamma,\mathcal{O}}\left(\|(\beta_{m-1})_{r_{2}}\|_{q,\overline{s}_{h}+\mathtt{p}}^{\gamma,\mathcal{O}}+\|\Delta_{12}\beta_{m-1}\|_{q,\overline{s}_{h}+\mathtt{p}}^{\gamma,\mathcal{O}}\right)\\
					&+C\|\Delta_{12}\beta_{m-1}\|_{q,\overline{s}_{h}+\mathtt{p}}^{\gamma,\mathcal{O}}\|(g_{m})_{r_{2}}\|_{q,s_{0}+1}^{\gamma,\mathcal{O}}\left(1+\|(\beta_{m-1})_{r_{2}}\|_{q,s_{0}}^{\gamma,\mathcal{O}}+\|\Delta_{12}\beta_{m-1}\|_{q,s_{0}}^{\gamma,\mathcal{O}}\right)
				\end{align*}
				and for all $s\in[s_{0},\overline{s}_{h}+\mathtt{p}]$
				\begin{align*}
					\|\Delta_{12}\beta_{m}\|_{q,s}^{\gamma,\mathcal{O}}&\leqslant\|\Delta_{12}g_{m}\|_{q,s}^{\gamma,\mathcal{O}}\left(1+C\|(\beta_{m-1})_{r_{1}}\|_{q,s_{0}}^{\gamma,\mathcal{O}}\right)+\|\Delta_{12}g_{m}\|_{q,s_{0}}^{\gamma,\mathcal{O}}\|(\beta_{m-1})_{r_{1}}\|_{q,s}^{\gamma,\mathcal{O}}\\
					&\quad+C\|\Delta_{12}\beta_{m-1}\|_{q,s_{0}}^{\gamma,\mathcal{O}}\|(g_{m})_{r_{2}}\|_{q,s+1}^{\gamma,\mathcal{O}}\left(1+\|(\beta_{m-1})_{r_{2}}\|_{q,s_{0}}^{\gamma,\mathcal{O}}+\|\Delta_{12}\beta_{m-1}\|_{q,s_{0}}^{\gamma,\mathcal{O}}\right)\\
					&\quad+C\|\Delta_{12}\beta_{m-1}\|_{q,s_{0}}^{\gamma,\mathcal{O}}\|(g_{m})_{r_{2}}\|_{q,s_{0}+1}^{\gamma,\mathcal{O}}\left(\|(\beta_{m-1})_{r_{2}}\|_{q,s}^{\gamma,\mathcal{O}}+\|\Delta_{12}\beta_{m-1}\|_{q,s}^{\gamma,\mathcal{O}}\right)\\
					&\quad+\|\Delta_{12}\beta_{m-1}\|_{q,s}^{\gamma,\mathcal{O}}\left(1+C\|(g_{m})_{r_{2}}\|_{q,s_{0}+1}^{\gamma,\mathcal{O}}\left(1+\|(\beta_{m-1})_{r_{2}}\|_{q,s_{0}}^{\gamma,\mathcal{O}}+\|\Delta_{12}\beta_{m-1}\|_{q,s_{0}}^{\gamma,\mathcal{O}}\right)\right).
				\end{align*}
				 Notice that \eqref{param-trans} implies in particular $\overline{s}_{h}+\mathtt{p}+\tau_{1}q+\tau_{1}+3\leqslant s_{h}+\sigma_{1}.$ Therefore, using \eqref{bound gm} and \eqref{smallness condition transport}, we get
				\begin{align}\label{maj gmrk}
					\sup_{m\in\mathbb{N}}\max_{k\in\{1,2\}}\|(g_{m})_{r_{k}}\|_{q,\overline{s}_{h}+\mathtt{p}+1}^{\gamma,\mathcal{O}}&\leqslant C \varepsilon\gamma^{-1}\left(1+\max_{k\in\{1,2\}}\|\mathfrak{I}_{k}\|_{q,\overline{s}_{h}+\mathtt{p}+\tau_{1}q+\tau_{1}+3}^{\gamma,\mathcal{O}}\right)\nonumber\\
					& \leqslant C.
				\end{align}
				Notice that the previous estimate is  sufficient to easily get rid of most of  terms in the estimates of $\Delta_{12}\beta_{m}$ and $\Delta_{12}(\beta_{m}-\beta_{m-1})$, but not enough to make the series \eqref{series Delta12beta} convergent. For this purpose, we shall refine the estimates. By \eqref{bound gm}, \eqref{param-trans} and \eqref{smallness condition transport}, we have 
				\begin{align}\label{deacr gm}
					\max_{k\in\{1,2\}}\|(g_{m})_{r_{k}}\|_{q,\overline{s}_{h}+\mathtt{p}+1}&\leqslant C\varepsilon\gamma^{-1}\left(1+\|\mathfrak{I}_{k}\|_{q,\overline{s}_{h}+\mathtt{p}+\tau_{1}q+\tau_{1}+3}^{\gamma,\mathcal{O}}\right)N_{0}^{\overline{\theta}\left(\overline{s}_{h}+\mathtt{p}+1\right)\mu_{2}}N_{m}^{-\overline{\theta}\left(\overline{s}_{h}+\mathtt{p}+1\right)\mu_{2}}\nonumber\\
					&\leqslant C\varepsilon\gamma^{-1}N_{0}^{\overline{\theta}\left(\overline{s}_{h}+\mathtt{p}+1\right)\mu_{2}}N_{m}^{-\overline{\theta}\left(\overline{s}_{h}+\mathtt{p}+1\right)\mu_{2}}.
				\end{align}
				Combining \eqref{uniform estimate betam} and \eqref{smallness condition transport}
				\begin{align}\label{maj betamrk}
					\sup_{m\in\mathbb{N}}\max_{k\in\{1,2\}}\|(\beta_{m})_{r_{k}}\|_{q,\overline{s}_{h}+\mathtt{p}}&\leqslant C\varepsilon\gamma^{-1}\left(1+\max_{k\in\{1,2\}}\|\mathfrak{I}_{k}\|_{q,\overline{s}_{h}+\mathtt{p}+\tau_{1}q+\tau_{1}+3}^{\gamma,\mathcal{O}}\right)\nonumber\\
					&\leqslant C.
				\end{align}
				Hence, using \eqref{maj gmrk}, \eqref{maj betamrk} and Sobolev embeddings, the previous two estimates can be reduced to
				\begin{align}\label{itmdt dlt12 diff bt}
					\|\Delta_{12}(\beta_{m}-\beta_{m-1})\|_{q,\overline{s}_{h}+\mathtt{p}}^{\gamma,\mathcal{O}}&\leqslant C\left(\|\Delta_{12}g_{m}\|_{q,\overline{s}_{h}+\mathtt{p}}^{\gamma,\mathcal{O}}+\|\Delta_{12}\beta_{m-1}\|_{q,\overline{s}_{h}+\mathtt{p}}^{\gamma,\mathcal{O}}\|(g_{m})_{r_{2}}\|_{q,\overline{s}_{h}+\mathtt{p}+1}^{\gamma,\mathcal{O}}\right),
				\end{align}
			
				\begin{align}\label{dlt12 btm s0}
					\|\Delta_{12}\beta_{m}\|_{q,s_{0}}^{\gamma,\mathcal{O}}&\leqslant C\|\Delta_{12}g_{m}\|_{q,s_{0}}^{\gamma,\mathcal{O}}+\|\Delta_{12}\beta_{m-1}\|_{q,s_{0}}^{\gamma,\mathcal{O}}\left(1+C\|(g_{m})_{r_{2}}\|_{q,s_{0}+1}^{\gamma,\mathcal{O}}\right)
				\end{align}
				and
				\begin{align}\label{dlt12 btm shb+1}
					\|\Delta_{12}\beta_{m}\|_{q,\overline{s}_{h}+\mathtt{p}}^{\gamma,\mathcal{O}}&\leqslant C\left(\|\Delta_{12}g_{m}\|_{q,\overline{s}_{h}+\mathtt{p}}^{\gamma,\mathcal{O}}+\|\Delta_{12}\beta_{m-1}\|_{q,s_{0}}^{\gamma,\mathcal{O}}\|(g_{m})_{r_{2}}\|_{q,\overline{s}_{h}+\mathtt{p}+1}\right)\nonumber\\
					&\quad+\|\Delta_{12}\beta_{m-1}\|_{q,\overline{s}_{h}+\mathtt{p}}^{\gamma,\mathcal{O}}\left(1+C\|(g_{m})_{r_{2}}\|_{q,s_{0}+1}^{\gamma,\mathcal{O}}\right).
				\end{align}
				From \eqref{dlt12 btm s0}, using \eqref{Ind-res1} and the fact that $\beta_{0}=g_{0}$, we deduce that
				$$\sup_{m\in\mathbb{N}}\|\Delta_{12}\beta_{m}\|_{q,s_{0}}^{\gamma,\mathcal{O}}\leqslant\left(\|\Delta_{12}g_{0}\|_{q,s_{0}}^{\gamma,\mathcal{O}}+C\sum_{k=0}^{\infty}\|\Delta_{12}g_{k}\|_{q,s_{0}}^{\gamma,\mathcal{O}}\right)\prod_{k=0}^{\infty}\left(1+\|(g_{k})_{r_{2}}\|_{q,s_{0}+1}^{\gamma,\mathcal{O}}\right).$$
				Adding \eqref{deacr gm}, we obtain
				$$\sup_{m\in\mathbb{N}}\|\Delta_{12}\beta_{m}\|_{q,s_{0}}^{\gamma,\mathcal{O}}\leqslant C\sum_{k=0}^{\infty}\|\Delta_{12}g_{k}\|_{q,s_{0}}^{\gamma,\mathcal{O}}.$$
				Similarly, \eqref{dlt12 btm shb+1}, \eqref{Ind-res1}, \eqref{deacr gm} and the previous estimate allow to get
				$$\sup_{m\in\mathbb{N}}\|\Delta_{12}\beta_{m}\|_{q,\overline{s}_{h}+\mathtt{p}}^{\gamma,\mathcal{O}}\leqslant C\sum_{k=0}^{\infty}\|\Delta_{12}g_{k}\|_{q,\overline{s}_{h}+\mathtt{p}}^{\gamma,\mathcal{O}}.$$ 
				Putting together the previous bounds, \eqref{itmdt dlt12 diff bt} and \eqref{deacr gm} gives
				\begin{equation}\label{telscopic diff betam}
					\|\Delta_{12}(\beta_{m}-\beta_{m-1})\|_{q,\overline{s}_{h}+\mathtt{p}}^{\gamma,\mathcal{O}}\lesssim \|\Delta_{12}g_{m}\|_{q,\overline{s}_{h}+\mathtt{p}}^{\gamma,\mathcal{O}}+\varepsilon\gamma^{-1}N_{0}^{\overline{\theta}(\overline{s}_{h}+\mathtt{p}+1)\mu_{2}}N_{m}^{-\overline{\theta}(\overline{s}_{h}+\mathtt{p}+1)\mu_{2}}\sum_{k=0}^{\infty}\|\Delta_{12}g_{k}\|_{q,\overline{s}_{h}+\mathtt{p}}^{\gamma,\mathcal{O}}.
				\end{equation}
				Thus, the main delicate point is   to estimate $\Delta_{12}g_{m}.$ First remark that according to \eqref{set gtlj} and \eqref{def gm}, we can make the splitting
				\begin{align*}
					g_{m}(\mu,\varphi,\theta)&=\ii\sum_{(l,j)\in\mathbb{Z}^{d+1}\setminus\{0\}\atop\langle l,j\rangle\leqslant N_{m}}a_{l,j}\widehat{\chi}\big(a_{l,j}(A_{l,j})_{r_{2}}(\mu)\big)(\Delta_{12}f_{m})_{l,j}(\mu)\mathbf{e}_{l,j}\\
					&\textnormal{\hspace{0.3cm}}+\ii\sum_{(l,j)\in\mathbb{Z}^{d+1}\setminus\{0\}\atop\langle l,j\rangle\leqslant N_{m}}a_{l,j}\Delta_{12}\widehat{\chi}\big(a_{l,j}A_{l,j}(\mu)\big)((f_{m})_{r_{1}})_{l,j}(\mu)\mathbf{e}_{l,j}\\
					&:=\mathbf{I}_{1}+\mathbf{I}_{2}.
				\end{align*}
			Similarly to \eqref{control of g by f}, one obtains
			\begin{equation}\label{I1}
				\|\mathbf{I}_{1}\|_{q,s}^{\gamma,\mathcal{O}}\lesssim\gamma^{-1}\|\Pi_{N_{m}}\Delta_{12}f_{m}\|_{q,s+\tau_{1}q+\tau_{1}}^{\gamma,\mathcal{O}}.
			\end{equation}
			We shall now estimate the second term. Applying Taylor Formula, we get
			\begin{align}\label{Taylor-I2}
				\mathbf{I}_{2}&=\ii\sum_{(l,j)\in\mathbb{Z}^{d+1}\setminus\{0\}\atop\langle l,j\rangle\leqslant N_{m}}a_{l,j}^{2}(\Delta_{12}A_{l,j})\int_{0}^{1}\widehat{\chi}'\Big(a_{l,j}\Big[\tau(A_{l,j})_{r_{1}}(\mu)+(1-\tau)(A_{l,j})_{r_{2}}(\mu)\Big]\Big)d\tau((f_{m})_{r_{1}})_{l,j}\mathbf{e}_{l,j}\nonumber\\
				&:=\sum_{(l,j)\in\mathbb{Z}^{d+1}\setminus\{0\}\atop\langle l,j\rangle\leqslant N_{m}}\widetilde{h}_{l,j}(\mu)((f_{m})_{r_{1}})_{l,j}(\mu)\mathbf{e}_{l,j}.
			\end{align}
		Remark that direct computations yield
		\begin{equation}\label{est dAlj}
			\forall q'\in\llbracket0,q\rrbracket,\quad\|\Delta_{12}A_{l,j}\|_{q'}^{\gamma,\mathcal{O}}\lesssim\langle l,j\rangle\|\Delta_{12}V_{m}\|_{q'}^{\gamma,\mathcal{O}}.
		\end{equation}
		Since that $\widehat{\chi}'\in C^{\infty}$ with $\widehat{\chi}'(0)=0,$ then applying Lemma \ref{Lem-lawprod}-(iv)-(vi) together with \eqref{est Alj} and \eqref{est dAlj}, we get
		\begin{align*}
			\forall q'\in\llbracket0,q\rrbracket,\quad \|\widetilde{h_{l,j}}\|_{q'}^{\gamma,\mathcal{O}}&\lesssim a_{l,j}^{3}\|\Delta_{12}A_{l,j}\|_{q'}^{\gamma,\mathcal{O}}\left(\|(A_{l,j})_{r_{1}}\|_{q'}^{\gamma,\mathcal{O}}+\|(A_{l,j})_{r_{2}}\|_{q'}^{\gamma,\mathcal{O}}\right)\\
			&\textnormal{\hspace{0.3cm}}\times\left(1+a_{l,j}^{q'-1}\left(\|(A_{l,j})_{r_{1}}\|_{L^{\infty}(\mathcal{O})}+\|(A_{l,j})_{r_{2}}\|_{L^{\infty}(\mathcal{O})}\right)^{q'-1}\right)\\
			&\lesssim\gamma^{-\upsilon(q'+2)}\langle l,j\rangle^{\tau_{1}q'+2\tau_{1}+q'+1}\|\Delta_{12}V_{m}\|_{q'}^{\gamma,\mathcal{O}}.
		\end{align*}
	By assumption in Proposition \ref{reduction of the transport part}, we have
	\begin{equation}\label{second choice of upsilon}
		\upsilon\leqslant\tfrac{1}{q+2}
	\end{equation}
	and using Leibniz rule, we deduce that
	\begin{equation}\label{I2}
		\|\mathbf{I}_{2}\|_{q}^{\gamma,\mathcal{O}}\lesssim\gamma^{-1}\|\Delta_{12}V_{m}\|_{q}^{\gamma,\mathcal{O}}\|\Pi_{N_{m}}(f_{m})_{r_{1}}\|_{q,s+\tau_{1}q+2\tau_{1}+1}^{\gamma,\mathcal{O}}.
	\end{equation}
	Putting together \eqref{I1} and \eqref{I2}, we obtain for all $s\geqslant s_{0}$
	\begin{align}\label{estimate delta12 gm}
		\|\Delta_{12}g_{m}\|_{q,s}^{\gamma,\mathcal{O}}&\lesssim\gamma^{-1}\|\Pi_{N_{m}}\Delta_{12}f_{m}\|_{q,s+\tau_{1}q+\tau_{1}}^{\gamma,\mathcal{O}}+\gamma^{-1}\|\Delta_{12}V_{m}\|_{q}^{\gamma,\mathcal{O}}\|\Pi_{N_{m}}(f_{m})_{r_{1}}\|_{q,s+\tau_{1}q+2\tau_{1}+1}^{\gamma,\mathcal{O}}\\
		&\lesssim\gamma^{-1}N_{m}^{\tau_{1}q+\tau_{1}}\|\Delta_{12}f_{m}\|_{q,s}^{\gamma,\mathcal{O}}+\gamma^{-2}N_{m}^{\tau_{1}q+2\tau_{1}+1}\|\Delta_{12}V_{m}\|_{q}^{\gamma,\mathcal{O}}\|(f_{m})_{r_{1}}\|_{q,s}^{\gamma,\mathcal{O}}\nonumber.
	\end{align}
	Therefore, estimating $\Delta_{12}g_{m}$ can be done through the estimate of $\Delta_{12}f_{m}.$ To do so, we shall argue by induction. For that purpose, we shall consider a parameter $\widetilde{\mathtt{p}}$ (which can depend on the parameter $\mathtt{p}$,  see for instance \eqref{link pt p}) satisfying the following constraint
	\begin{equation}\label{cond-pt}
		\overline{s}_{h}+\widetilde{\mathtt{p}}+3\leqslant s_{h}+\sigma_{1}.
	\end{equation}
	We denote 
	$$u_{m}:=\Pi_{N_{m}}^{\perp}f_{m}+f_{m}\partial_{\theta}g_{m}.$$
	Then, we can write
	$$\Delta_{12}f_{m+1}=(\mathcal{G}_{m}^{-1})_{r_{1}}\Delta_{12}u_{m}+\left(\Delta_{12}\mathcal{G}_{m}^{-1}\right)(u_{m})_{r_{2}},$$
	with
	$$\Delta_{12}u_{m}=\Pi_{N_{m}}^{\perp}\Delta_{12}f_{m}+\Delta_{12}f_{m}\partial_{\theta}(g_{m})_{r_{1}}+(f_{m})_{r_{2}}\partial_{\theta}\Delta_{12}g_{m}.
	$$
	By the triangle inequality, we have for all $s\geqslant s_{0}$
	\begin{align}\label{IT diff fm+1}
		\|\Delta_{12}f_{m+1}\|_{q,s}^{\gamma,\mathcal{O}}&\leqslant\|(\mathcal{G}_{m}^{-1})_{r_{1}}\Delta_{12}u_{m}\|_{q,s}^{\gamma,\mathcal{O}}+\|(\Delta_{12}\mathcal{G}_{m}^{-1})(u_{m})_{r_{2}}\|_{q,s}^{\gamma,\mathcal{O}}.
	\end{align}
				Therefore, combining  \eqref{tame comp}, \eqref{link betah and beta}, \eqref{control of g by f} and Lemma \ref{Lem-lawprod}-(ii), we get for all $s\geqslant s_{0}$
				\begin{align*}
					\|(\mathcal{G}_{m}^{-1})_{r_{1}}\Delta_{12}u_{m}\|_{q,s}^{\gamma,\mathcal{O}}&\leqslant\|\Delta_{12}u_{m}\|_{q,s}^{\gamma,\mathcal{O}}\left(1+C\|(\widehat{g_{m}})_{r_{1}}\|_{q,s_{0}}^{\gamma,\mathcal{O}}\right)+C\|(\widehat{g_{m}})_{r_{1}}\|_{q,s}^{\gamma,\mathcal{O}}\|\Delta_{12}u_{m}\|_{q,s_{0}}^{\gamma,\mathcal{O}}\\
					&\leqslant\|\Delta_{12}u_{m}\|_{q,s}^{\gamma,\mathcal{O}}\left(1+C\|(g_{m})_{r_{1}}\|_{q,s_{0}}^{\gamma,\mathcal{O}}\right)+C\|(g_{m})_{r_{1}}\|_{q,s}^{\gamma,\mathcal{O}}\|\Delta_{12}u_{m}\|_{q,s_{0}}^{\gamma,\mathcal{O}}\\
					&\leqslant\|\Delta_{12}u_{m}\|_{q,s}^{\gamma,\mathcal{O}}\left(1+C\gamma^{-1}N_{m}^{\tau_{1}q+\tau_{1}}\max_{k\in\{1,2\}}\|(f_{m})_{r_{k}}\|_{q,s_{0}}^{\gamma,\mathcal{O}}\right)\\
					&\quad+C\gamma^{-1}N_{m}^{\tau_{1}q+\tau_{1}}\max_{k\in\{1,2\}}\|(f_{m})_{r_{k}}\|_{q,s}^{\gamma,\mathcal{O}}\|\Delta_{12}u_{m}\|_{q,s_{0}}^{\gamma,\mathcal{O}}.
				\end{align*}
				Using \eqref{uniform estimate of deltams},\eqref{cond-pt} and \eqref{smallness condition transport}, one  gets
				\begin{align}\label{unif fm shb+1}
					\gamma^{-1}\sup_{m\in\mathbb{N}}\max_{k\in\{1,2\}}\|(f_{m})_{r_{k}}\|_{q,\overline{s}_{h}+\widetilde{\mathtt{p}}+1}^{\gamma,\mathcal{O}}&\leqslant C\varepsilon\gamma^{-1}\left(1+\max_{k\in\{1,2\}}\|\mathfrak{I}_{k}\|_{q,\overline{s}_{h}+\widetilde{\mathtt{p}}+2}^{\gamma,\mathcal{O}}\right)\nonumber\\
					&\leqslant C.
				\end{align}
				Therefore, from \eqref{hypothesis of induction deltam} and \eqref{unif fm shb+1}, we get for all $s\in[s_{0},\overline{s}_{h}]$
				$$\|(\mathcal{G}_{m}^{-1})_{r_{1}}\Delta_{12}u_{m}\|_{q,s}^{\gamma,\mathcal{O}}\leqslant\|\Delta_{12}u_{m}\|_{q,s}^{\gamma,\mathcal{O}}\left(1+CN_{0}^{\overline{\mu}_{2}}N_{m}^{\tau_{1}q+\tau_{1}-\overline{\mu}_{2}}\right)+CN_{m}^{\tau_{1}q+\tau_{1}}\|\Delta_{12}u_{m}\|_{q,s_{0}}^{\gamma,\mathcal{O}}.$$
				At this level we need to give a suitable estimate for $\Delta_{12}u_{m}.$ For this aim, we apply the  law products in Lemma \ref{Lem-lawprod}, ensuring that for all $s\geqslant s_{0}$
				\begin{align*}
					\|\Delta_{12}u_{m}\|_{q,s}^{\gamma,\mathcal{O}}&\leqslant \|\Pi_{N_{m}}^{\perp}\Delta_{12}f_{m}\|_{q,s}^{\gamma,\mathcal{O}}+C\|\Delta_{12}f_{m}\|_{q,s}^{\gamma,\mathcal{O}}\|\partial_{\theta}(g_{m})_{r_{1}}\|_{q,s_{0}}^{\gamma,\mathcal{O}}+C\|\Delta_{12}f_{m}\|_{q,s_{0}}^{\gamma,\mathcal{O}}\|\partial_{\theta}(g_{m})_{r_{1}}\|_{q,s}^{\gamma,\mathcal{O}}\\
					&\quad+C\|(f_{m})_{r_{2}}\|_{q,s}^{\gamma,\mathcal{O}}\|\Delta_{12}g_{m}\|_{q,s_{0}}^{\gamma,\mathcal{O}}+C\|(f_{m})_{r_{2}}\|_{q,s_{0}}^{\gamma,\mathcal{O}}\|\Delta_{12}g_{m}\|_{q,s}^{\gamma,\mathcal{O}}.
				\end{align*}
				Hence we deduce by  \eqref{control of g by f} and Lemma \ref{Lem-lawprod}-(ii), 				\begin{align*}
					\|\Delta_{12}u_{m}\|_{q,s}^{\gamma,\mathcal{O}}&\leqslant \|\Pi_{N_{m}}^{\perp}\Delta_{12}f_{m}\|_{q,s}^{\gamma,\mathcal{O}}+C\gamma^{-1}N_{m}^{\tau_{1}q+\tau_{1}+1}\|\Delta_{12}f_{m}\|_{q,s}^{\gamma,\mathcal{O}}\max_{k\in\{1,2\}}\|(f_{m})_{r_{k}}\|_{q,s_{0}}^{\gamma,\mathcal{O}}\\
					&\quad+C\gamma^{-1}N_{m}^{\tau_{1}q+\tau_{1}+1}\|\Delta_{12}f_{m}\|_{q,s_{0}}^{\gamma,\mathcal{O}}\max_{k\in\{1,2\}}\|(f_{m})_{r_{k}}\|_{q,s}^{\gamma,\mathcal{O}}\\
					&\quad+C\max_{k\in\{1,2\}}\|(f_{m})_{r_{k}}\|_{q,s}^{\gamma,\mathcal{O}}\|\Delta_{12}g_{m}\|_{q,s_{0}}^{\gamma,\mathcal{O}}+C\max_{k\in\{1,2\}}\|(f_{m})_{r_{k}}\|_{q,s_{0}}^{\gamma,\mathcal{O}}\|\Delta_{12}g_{m}\|_{q,s}^{\gamma,\mathcal{O}}.
				\end{align*}
				Added to \eqref{estimate delta12 gm}, we finally obtain for all $s\geqslant s_{0}$
				\begin{align*}
					\|\Delta_{12}u_{m}\|_{q,s}^{\gamma,\mathcal{O}}&\leqslant \|\Pi_{N_{m}}^{\perp}\Delta_{12}f_{m}\|_{q,s}^{\gamma,\mathcal{O}}+C\gamma^{-1}N_{m}^{\tau_{1}q+\tau_{1}+1}\|\Delta_{12}f_{m}\|_{q,s}^{\gamma,\mathcal{O}}\max_{k\in\{1,2\}}\|(f_{m})_{r_{k}}\|_{q,s_{0}}^{\gamma,\mathcal{O}}\\
					&\quad+C\gamma^{-1}N_{m}^{\tau_{1}q+\tau_{1}+1}\|\Delta_{12}f_{m}\|_{q,s_{0}}^{\gamma,\mathcal{O}}\max_{k\in\{1,2\}}\|(f_{m})_{r_{k}}\|_{q,s}^{\gamma,\mathcal{O}}\\
					&\quad+C\gamma^{-2}N_{m}^{\tau_{1}q+2\tau_{1}+1}\max_{k\in\{1,2\}}\|(f_{m})_{r_{k}}\|_{q,s}^{\gamma,\mathcal{O}}\max_{k\in\{1,2\}}\|(f_{m})_{r_{k}}\|_{q,s_{0}}^{\gamma,\mathcal{O}}\|\Delta_{12}V_{m}\|_{q}^{\gamma,\mathcal{O}}.
				\end{align*}
				Consequently, we find from  \eqref{hypothesis of induction deltam}, Lemma \ref{Lem-lawprod}-(ii) and \eqref{unif fm shb+1},
				\begin{align*}
					\|\Delta_{12}u_{m}\|_{q,s_{0}}^{\gamma,\mathcal{O}}&\leqslant N_{m}^{s_{0}-\overline{s}_{h}-\widetilde{\mathtt{p}}}\|\Delta_{12}f_{m}\|_{q,\overline{s}_{h}+\widetilde{\mathtt{p}}}^{\gamma,\mathcal{O}}+CN_{0}^{\overline{\mu}_{2}}N_{m}^{\tau_{1}q+\tau_{1}+1-\overline{\mu}_{2}}\delta_{0}^{1,2}(\overline{s}_{h})\|\Delta_{12}f_{m}\|_{q,s_{0}}^{\gamma,\mathcal{O}}\\
					&\quad+CN_{0}^{2\overline{\mu}_{2}}N_{m}^{\tau_{1}q+2\tau_{1}+1-2\overline{\mu}_{2}}\delta_{0}^{1,2}(\overline{s}_{h})\|\Delta_{12}V_{m}\|_{q}^{\gamma,\mathcal{O}}
				\end{align*}
				and
				\begin{align*}
					\|\Delta_{12}u_{m}\|_{q,\overline{s}_{h}+\widetilde{\mathtt{p}}}^{\gamma,\mathcal{O}}&\leqslant\|\Delta_{12}f_{m}\|_{q,\overline{s}_{h}+\widetilde{\mathtt{p}}}^{\gamma,\mathcal{O}}\left(1+CN_{0}^{\overline{\mu}_{2}}N_{m}^{\tau_{1}q+\tau_{1}+1-\overline{\mu}_{2}}\delta_{0}^{1,2}(\overline{s}_{h})\right)+CN_{m}^{\tau_{1}q+\tau_{1}+1}\delta_{0}^{1,2}(\overline{s}_{h})\|\Delta_{12}f_{m}\|_{q,s_{0}}^{\gamma,\mathcal{O}}\\
					&\quad+CN_{0}^{\overline{\mu}_{2}}N_{m}^{\tau_{1}q+2\tau_{1}+1-\overline{\mu}_{2}}\delta_{0}^{1,2}(\overline{s}_{h})\|\Delta_{12}V_{m}\|_{q}^{\gamma,\mathcal{O}},
				\end{align*}
				where we use the notation
				$$\delta_{0}^{1,2}(s):=\gamma^{-1}\max_{k\in\{1,2\}}\|(f_{0})_{r_{k}}\|_{q,s}^{\gamma,\mathcal{O}}.$$
				It follows from the preceding estimates that,
				\begin{align}\label{Gm-1 diff s0}
					\|(\mathcal{G}_{m}^{-1})_{r_{1}}\Delta_{12}u_{m}\|_{q,s_{0}}^{\gamma,\mathcal{O}}&\leqslant CN_m^{\tau_{1}q+\tau_1}\|\Delta_{12}u_m\|_{q,s_0}^{\gamma,\mathcal{O}}\nonumber\\
					&\leqslant C N_{m}^{s_{0}+\tau_{1}q+\tau_{1}-\overline{s}_{h}-\widetilde{\mathtt{p}}}\|\Delta_{12}f_{m}\|_{q,\overline{s}_{h}+\widetilde{\mathtt{p}}}^{\gamma,\mathcal{O}}+CN_{0}^{\overline{\mu}_{2}}N_{m}^{2(\tau_{1}q+\tau_{1})+1-\overline{\mu}_{2}}\delta_{0}^{1,2}(\overline{s}_{h})\|\Delta_{12}f_{m}\|_{q,s_{0}}^{\gamma,\mathcal{O}}\nonumber\\
					&\quad+CN_{0}^{2\overline{\mu}_{2}}N_{m}^{2\tau_{1}q+3\tau_{1}+1-2\overline{\mu}_{2}}\delta_{0}^{1,2}(\overline{s}_{h})\|\Delta_{12}V_{m}\|_{q}^{\gamma,\mathcal{O}}.
				\end{align}
				In a similar way, direct computations yield
				\begin{align}\label{Gm-1 diff shb}
					\|(\mathcal{G}_{m}^{-1})_{r_{1}}\Delta_{12}u_{m}\|_{q,\overline{s}_{h}+\widetilde{\mathtt{p}}}^{\gamma,\mathcal{O}}&\leqslant\|\Delta_{12}f_{m}\|_{q,\overline{s}_{h}+\widetilde{\mathtt{p}}}^{\gamma,\mathcal{O}}\left(1+N_m^{s_0+\tau_{1}q+\tau_{1}-\overline{s}_{h}-\widetilde{\mathtt{p}}}+C N_{0}^{\overline{\mu}_{2}}N_{m}^{\tau_{1}q+\tau_{1}+1-\overline{\mu}_{2}}\delta_{0}^{1,2}(\overline{s}_{h})\right)\nonumber\\
					&\quad+C\left(N_0^{\overline{\mu}_{2}}N_m^{2(\tau_{1}q+\tau_{1})+1-\overline{\mu}_{2}}+N_{m}^{\tau_{1}q+\tau_{1}+1}\right)\delta_{0}^{1,2}(\overline{s}_{h})\|\Delta_{12}f_{m}\|_{q,s_{0}}^{\gamma,\mathcal{O}}\nonumber\\
					&\quad+CN_{0}^{\overline{\mu}_{2}}N_{m}^{2\tau_{1}q+3\tau_{1}+1-\overline{\mu}_{2}}\delta_{0}^{1,2}(\overline{s}_{h})\|\Delta_{12}V_{m}\|_{q}^{\gamma,\mathcal{O}}.
				\end{align} 
				By a new use of Taylor Formula, we can write
				$$(\Delta_{12}\mathcal{G}_{m}^{-1})(u_{m})_{r_{2}}(\theta)=\Delta_{12}\widehat{g_{m}}(\theta)\int_{0}^{1}\partial_{\theta}(u_{m})_{r_{2}}\Big(\theta+(\widehat{g_{m}})_{r_{2}}(\theta)+t\Delta_{12}\widehat{g_{m}}(\theta)\Big)dt.$$
				Applying  Lemma \ref{Lem-lawprod} and \eqref{control of g by f}, we deduce  for all $s\geqslant s_{0}$
				\begin{align}\label{link um-fm}
					\|u_{m}\|_{q,s}^{\gamma,\mathcal{O}}&\leqslant\|\Pi_{N_{m}}^{\perp}f_{m}\|_{q,s}^{\gamma,\mathcal{O}}+C\|f_{m}\|_{q,s}^{\gamma,\mathcal{O}}\|\partial_{\theta}g_{m}\|_{q,s_{0}}^{\gamma,\mathcal{O}}+C\|f_{m}\|_{q,s_{0}}^{\gamma,\mathcal{O}}\|\partial_{\theta}g_{m}\|_{q,s}^{\gamma,\mathcal{O}}\nonumber\\
					&\leqslant \|f_{m}\|_{q,s}\left(1+CN_{m}^{\tau_{1}q+\tau_{1}+1}\|f_{m}\|_{q,s_{0}}^{\gamma,\mathcal{O}}\right)\nonumber\\
					&\leqslant C\|f_{m}\|_{q,s}^{\gamma,\mathcal{O}}.
				\end{align}
				Using once again the law products in Lemma \ref{Lem-lawprod} combined with  \eqref{tame comp} yield for all $s\geqslant s_{0}$
				\begin{align*}
					\|(\Delta_{12}\mathcal{G}_{m}^{-1})(u_{m})_{r_{2}}\|_{q,s}^{\gamma,\mathcal{O}}&\leqslant C\|\Delta_{12}\widehat{g_{m}}\|_{q,s}^{\gamma,\mathcal{O}}\|(u_{m})_{r_{2}}\|_{q,s_{0}+1}\left(1+\|(\widehat{g_{m}})_{r_{2}}\|_{q,s_{0}}^{\gamma,\mathcal{O}}+\|\Delta_{12}\widehat{g_{m}}\|_{q,s_{0}}^{\gamma,\mathcal{O}}\right)\\
					&\quad+C\|\Delta_{12}\widehat{g_{m}}\|_{q,s_{0}}^{\gamma,\mathcal{O}}\|(u_{m})_{r_{2}}\|_{q,s+1}^{\gamma,\mathcal{O}}\left(1+\|(\widehat{g_{m}})_{r_{2}}\|_{q,s_{0}}^{\gamma,\mathcal{O}}+\|\Delta_{12}\widehat{g_{m}}\|_{q,s_{0}}^{\gamma,\mathcal{O}}\right)\\
					&\quad+C\|\Delta_{12}\widehat{g_{m}}\|_{q,s_{0}}\|(u_{m})_{r_{2}}\|_{q,s_{0}+1}^{\gamma,\mathcal{O}}\left(\|(\widehat{g_{m}})_{r_{2}}\|_{q,s}^{\gamma,\mathcal{O}}+\|\Delta_{12}\widehat{g_{m}}\|_{q,s}^{\gamma,\mathcal{O}}\right).
				\end{align*}
				In view of \eqref{link diff beta hat and diff beta}, \eqref{deacr gm} and Sobolev embeddings, one gets for all $s\in[s_{0},\overline{s}_{h}+\mathtt{p}]$
				\begin{align*}
					\|\Delta_{12}\widehat{g_{m}}\|_{q,s}^{\gamma,\mathcal{O}}&\leqslant C\left(\|\Delta_{12}g_{m}\|_{q,s}^{\gamma,\mathcal{O}}+\|\Delta_{12}g_{m}\|_{q,s_{0}}^{\gamma,\mathcal{O}}\max_{k\in\{1,2\}}\|(g_{m})_{r_{k}}\|_{q,s+1}^{\gamma,\mathcal{O}}\right)\\
					&\leqslant C\|\Delta_{12}g_{m}\|_{q,s}^{\gamma,\mathcal{O}}.
				\end{align*}
				Putting together the previous estimates, \eqref{maj gmrk} and \eqref{link betah and beta} gives for all $s\in[s_{0},\overline{s}_{h}]$
				\begin{align*}
					\|(\Delta_{12}\mathcal{G}_{m}^{-1})(u_{m})_{r_{2}}\|_{q,s}^{\gamma,\mathcal{O}}&\leqslant C\|\Delta_{12}g_{m}\|_{q,s}^{\gamma,\mathcal{O}}\|(u_{m})_{r_{2}}\|_{q,s_{0}+1}+C\|\Delta_{12}g_{m}\|_{q,s_{0}}^{\gamma,\mathcal{O}}\|(u_{m})_{r_{2}}\|_{q,s+1}^{\gamma,\mathcal{O}}.
				\end{align*}
				Thus,  by virtue of  \eqref{estimate delta12 gm}, \eqref{link um-fm}, we get for all $s\in[s_{0},\overline{s}_{h}]$
				\begin{align*}
					\|(\Delta_{12}\mathcal{G}_{m}^{-1})(u_{m})_{r_{2}}\|_{q,s}^{\gamma,\mathcal{O}}&\leqslant C\gamma^{-1}N_{m}^{\tau_{1}q+\tau_{1}}\|\Delta_{12}f_{m}\|_{q,s}^{\gamma,\mathcal{O}}\max_{k\in\{1,2\}}\|(f_{m})_{r_{k}}\|_{q,s_{0}+1}^{\gamma,\mathcal{O}}\\
					&\quad+C\gamma^{-1}N_{m}^{\tau_{1}q+\tau_{1}}\|\Delta_{12}f_{m}\|_{q,s_{0}}^{\gamma,\mathcal{O}}\max_{k\in\{1,2\}}\|(f_{m})_{r_{k}}\|_{q,s+1}^{\gamma,\mathcal{O}}\\
					&\quad+C\gamma^{-2}N_{m}^{\tau_{1}q+2\tau_{1}+1}\|\Delta_{12}V_{m}\|_{q}^{\gamma,\mathcal{O}}\max_{k\in\{1,2\}}\|(f_{m})_{r_{k}}\|_{q,s_{0}+1}^{\gamma,\mathcal{O}}\max_{k\in\{1,2\}}\|(f_{m})_{r_{k}}\|_{q,s+1}^{\gamma,\mathcal{O}}.
				\end{align*}
				Hence,  \eqref{hypothesis of induction deltam}, \eqref{uniform estimate of deltams} and \eqref{unif fm shb+1} allow to get  (since $s_l\geqslant s_0+1$)
				\begin{align}\label{diff Gm-1 s0}
					\|(\Delta_{12}\mathcal{G}_{m}^{-1})(u_{m})_{r_{2}}\|_{q,s_{0}}^{\gamma,\mathcal{O}}&\leqslant CN_{0}^{\overline{\mu}_{2}}N_{m}^{\tau_{1}q+\tau_{1}-\overline{\mu}_{2}}\delta_{0}^{1,2}(\overline{s}_{h})\|\Delta_{12}f_{m}\|_{q,s_{0}}^{\gamma,\mathcal{O}}\nonumber\\
					&\quad+CN_{0}^{2\overline{\mu}_{2}}N_{m}^{\tau_{1}q+2\tau_{1}+1-2\overline{\mu}_{2}}\delta_{0}^{1,2}(\overline{s}_{h})\|\Delta_{12}V_{m}\|_{q}^{\gamma,\mathcal{O}}
				\end{align}
				and
				\begin{align}\label{diff Gm-1 shb}
					\|(\Delta_{12}\mathcal{G}_{m}^{-1})(u_{m})_{r_{2}}\|_{q,\overline{s}_{h}+\widetilde{\mathtt{p}}}^{\gamma,\mathcal{O}}&\leqslant CN_{0}^{\overline{\mu}_{2}}N_{m}^{\tau_{1}q+\tau_{1}-\overline{\mu}_{2}}\delta_{0}^{1,2}(\overline{s}_{h})\|\Delta_{12}f_{m}\|_{q,\overline{s}_{h}+\widetilde{\mathtt{p}}}^{\gamma,\mathcal{O}}\nonumber\\
					&\quad+CN_{m}^{\tau_{1}q+\tau_{1}}\delta_{0}^{1,2}(\overline{s}_{h}+\widetilde{\mathtt{p}}+1)\|\Delta_{12}f_{m}\|_{q,s_{0}}^{\gamma,\mathcal{O}}\nonumber\\
					&\quad+CN_{0}^{\overline{\mu}_{2}}N_{m}^{\tau_{1}q+2\tau_{1}+1-\overline{\mu}_{2}}\delta_{0}^{1,2}(\overline{s}_{h}+\widetilde{\mathtt{p}}+1)\|\Delta_{12}V_{m}\|_{q}^{\gamma,\mathcal{O}}.
				\end{align}
				Gathering \eqref{IT diff fm+1}, \eqref{Gm-1 diff s0} and \eqref{diff Gm-1 s0} implies (since $N_m^{-\widetilde{\mathtt{p}}}\leqslant 1$)
				\begin{align}\label{link diff fm+1 and fm s0}
					\|\Delta_{12}f_{m+1}\|_{q,s_{0}}^{\gamma,\mathcal{O}}&\leqslant N_{m}^{s_{0}+\tau_{1}q+\tau_{1}-\overline{s}_{h}}\|\Delta_{12}f_{m}\|_{q,\overline{s}_h+\widetilde{\mathtt{p}}}^{\gamma,\mathcal{O}}\nonumber\\
					&\quad+CN_{0}^{\overline{\mu}_{2}}N_{m}^{2(\tau_{1}q+\tau_{1})+1-\overline{\mu}_{2}}\delta_{0}^{1,2}(\overline{s}_{h})\|\Delta_{12}f_{m}\|_{q,s_{0}}^{\gamma,\mathcal{O}}\nonumber\\
					&\quad+CN_{0}^{2\overline{\mu}_{2}}N_{m}^{2\tau_{1}q+3\tau_{1}+1-2\overline{\mu}_{2}}\delta_{0}^{1,2}(\overline{s}_{h})\|\Delta_{12}V_{m}\|_{q}^{\gamma,\mathcal{O}}.
				\end{align}
				In a similar war, we get in view of  \eqref{IT diff fm+1}, \eqref{Gm-1 diff shb} and \eqref{diff Gm-1 shb} 
				\begin{align}\label{link diff fm+1 and fm shb}
					\|\Delta_{12}f_{m+1}\|_{q,\overline{s}_{h}+\widetilde{\mathtt{p}}}^{\gamma,\mathcal{O}}&\leqslant\|\Delta_{12}f_{m}\|_{q,\overline{s}_{h}+\widetilde{\mathtt{p}}}^{\gamma,\mathcal{O}}\left(1+N_m^{s_0+\tau_{1}q+\tau_{1}-\overline{s}_{h}-\widetilde{\mathtt{p}}}+CN_{0}^{\overline{\mu}_{2}}N_{m}^{\tau_{1}q+\tau_{1}+1-\overline{\mu}_{2}}\delta_{0}^{1,2}(\overline{s}_{h})\right)\nonumber\\
					&\quad+C\left(N_0^{\overline{\mu}_2}N_m^{2(\tau_{1}q+\tau_{1})+1-\overline{\mu}_2}+N_{m}^{\tau_{1}q+\tau_{1}+1}\right)\delta_{0}^{1,2}(\overline{s}_{h}+\widetilde{\mathtt{p}}+1)\|\Delta_{12}f_{m}\|_{q,s_{0}}^{\gamma,\mathcal{O}}\nonumber\\
					&\quad+CN_{0}^{\overline{\mu}_{2}}N_{m}^{2\tau_{1}q+3\tau_{1}+1-\overline{\mu}_{2}}\delta_{0}^{1,2}(\overline{s}_{h}+\widetilde{\mathtt{p}}+1)\|\Delta_{12}V_{m}\|_{q}^{\gamma,\mathcal{O}}.
				\end{align}
				In the sequel, we shall use the following notations
				$$\overline{\delta}_{m}(s)=\gamma^{-1}\|\Delta_{12}f_{m}\|_{q,s}^{\gamma,\mathcal{O}}\quad\textnormal{and}\quad\varkappa_{m}=\gamma^{-1}\|\Delta_{12}V_{m}\|_{q}^{\gamma,\mathcal{O}}.$$
				Notice that
				$$\Delta_{12}V_{m+1}=\Delta_{12}V_{m}+\langle\Delta_{12}f_{m}\rangle_{\varphi,\theta}\quad\textnormal{and}\Delta_{12}V_0=0.$$
				Then, by using  Sobolev embeddings, we obtain
				\begin{equation}\label{link varkappa and overline delta}
					\varkappa_{m}\leqslant\sum_{k=0}^{m-1}\overline{\delta}_{k}(s_{0}).
				\end{equation}
				We shall now prove by induction that, for all $\widetilde{\mathtt{p}}$ satisfying the condition \eqref{cond-pt}, we have
				\begin{equation}\label{hypothesis of induction overline delta}
					\forall k\leqslant m,\quad\overline{\delta}_{k}(s_{0})\leqslant N_{0}^{\overline{\mu}_{2}}N_{k}^{-\overline{\mu}_{2}}\nu(\overline{s}_{h}+\widetilde{\mathtt{p}})\quad\textnormal{and}\quad\overline{\delta}_{k}(\overline{s}_{h}+\widetilde{\mathtt{p}})\leqslant\left(2-\tfrac{1}{k+1}\right)\nu(\overline{s}_{h}+\widetilde{\mathtt{p}}),
				\end{equation}
				with
				$$\nu(s):=\overline{\delta}_{0}(s)+\varepsilon\gamma^{-1}\|\Delta_{12}i\|_{s_{0}+2}.
				$$
				First remark that the property \eqref{hypothesis of induction overline delta} is trivially satisfied for $m=0$ according to Sobolev embeddings. We now assume that \eqref{hypothesis of induction overline delta} is true at the order $m$ and let us check it at the next order. By the induction  assumption \eqref{hypothesis of induction overline delta} and \eqref{link varkappa and overline delta}, one obtains the following estimate
				\begin{equation}\label{est-vchim}
			\sup_{m\in\mathbb{N}}\varkappa_{m}\leqslant C\nu(\overline{s}_{h}+\widetilde{\mathtt{p}}).
		\end{equation}
				Using  \eqref{link diff fm+1 and fm s0}, \eqref{est-vchim} and hypothesis of induction \eqref{hypothesis of induction overline delta}, we find
				\begin{align*}
					\overline{\delta}_{m+1}(s_{0})&\leqslant N_{m}^{s_{0}+\tau_{1}q+\tau_{1}-\overline{s}_{h}}\overline{\delta}_{m}(\overline{s}_{h}+\widetilde{\mathtt{p}})+CN_{0}^{\overline{\mu}_{2}}N_{m}^{2(\tau_{1}q+\tau_{1})+1-\overline{\mu}_{2}}\delta_{0}^{1,2}(\overline{s}_{h})\overline{\delta}_{m}(s_{0})\\
					&\quad+CN_{0}^{2\overline{\mu}_{2}}N_{m}^{2\tau_{1}q+3\tau_{1}+1-2\overline{\mu}_{2}}\delta_{0}^{1,2}(\overline{s}_{h})\varkappa_{m}\\
					&\leqslant \left[2N_{m}^{s_{0}+\tau_{1}q+\tau_{1}-\overline{s}_{h}}+CN_{0}^{2\overline{\mu}_{2}}N_{m}^{2\tau_{1}q+3\tau_{1}+1-2\overline{\mu}_{2}}\delta_{0}^{1,2}(\overline{s}_{h})\right]\nu(\overline{s}_{h}+\widetilde{\mathtt{p}}).
				\end{align*}
			Then, in view of \eqref{param-trans}, we infer
			\begin{align*}
				2N_{m}^{s_{0}+\tau_{1}q+\tau_{1}-\overline{s}_{h}}&=2N_{m}^{-\frac{3}{2}\overline{\mu}_{2}-3}=2N_{m}^{-3}N_{m+1}^{-\overline{\mu}_{2}}\\
				&\leqslant 2N_{0}^{-3}N_{m+1}^{-\overline{\mu}_{2}}\\
				&\leqslant \frac{1}{2}N_{0}^{\overline{\mu}_{2}}N_{m+1}^{-\overline{\mu}_{2}}.
			\end{align*}
		To prove the last inequality, we remark that since $N_{0}\geqslant 2$ and $\overline{\mu}_{2}\geqslant 0$ (according to \eqref{param}), then 
		$$4\leqslant N_{0}^{\overline{\mu}_{2}+3}.$$
		Similarly,  from the expression of $\overline{\mu}_{2}$ in  \eqref{param-trans}
		and  using \eqref{definition of Nm} one obtains
		\begin{align*}
			CN_{0}^{2\overline{\mu}_{2}}N_{m}^{2\tau_{1}q+3\tau_{1}+1-2\overline{\mu}_{2}}\delta_{0}^{1,2}(\overline{s}_{h})&\leqslant C\varepsilon\gamma^{-1}N_{0}^{\overline{\mu}_{2}}N_{m}^{2\tau_{1}q+3\tau_{1}+1-\frac{1}{2}\overline{\mu}_{2}}N_{0}^{\overline{\mu}_{2}}N_{m+1}^{-\overline{\mu}_{2}}\\
			&\leqslant C\varepsilon\gamma^{-1}N_{0}^{2\tau_{1}q+3\tau_{1}+1+\frac{1}{2}\overline{\mu}_{2}}N_{0}^{\overline{\mu}_{2}}N_{m+1}^{-\overline{\mu}_{2}}\\
			&\leqslant C \varepsilon\gamma^{-1}N_{0}^{\overline{\mu}_{2}}N_{0}^{\overline{\mu}_{2}}N_{m+1}^{-\overline{\mu}_{2}}.
		\end{align*}
	Hence, choosing $\varepsilon_{0}$ small enough and using \eqref{smallness condition transport} we deduce that
	$$CN_{0}^{2\overline{\mu}_{2}}N_{m}^{2\tau_{1}q+3\tau_{1}+1-2\overline{\mu}_{2}}\delta_{0}^{1,2}(\overline{s}_{h})\leqslant\frac{1}{2}N_{0}^{\overline{\mu}_{2}}N_{m+1}^{-\overline{\mu}_{2}}.$$
	Gathering the preceding estimates gives
	$$\overline{\delta}_{m+1}(s_{0})\leqslant N_{0}^{\overline{\mu}_{2}}N_{m+1}^{-\overline{\mu}_{2}}\nu(\overline{s}_{h}+\widetilde{\mathtt{p}}).$$
				This ends the proof of the first statement in \eqref{hypothesis of induction overline delta}. As to the second one, we shall first write in view of \eqref{link diff fm+1 and fm shb}, 
				\begin{align*}
					\overline{\delta}_{m+1}(\overline{s}_{h}+\widetilde{\mathtt{p}})&\leqslant\overline{\delta}_{m}(\overline{s}_{h}+\widetilde{\mathtt{p}})\left(1+N_m^{s_0+\tau_{1}q+\tau_{1}-\overline{s}_h}+CN_{0}^{\overline{\mu}_{2}}N_{m}^{\tau_{1}q+\tau_{1}+1-\overline{\mu}_{2}}\delta_{0}^{1,2}(\overline{s}_{h})\right)\\
					&\quad+C\left(N_m^{\tau_{1}q+\tau_{1}+1}+N_0^{\overline{\mu}_2}N_{m}^{2(\tau_{1}q+\tau_{1})+1-\overline{\mu}_2}\right)\delta_{0}^{1,2}(\overline{s}_{h}+\widetilde{\mathtt{p}}+1)\overline{\delta}_{m}(s_{0})\\
					&\quad+CN_{0}^{\overline{\mu}_{2}}N_{m}^{2\tau_{1}q+3\tau_{1}+1-\overline{\mu}_{2}}\delta_{0}^{1,2}(\overline{s}_{h}+\widetilde{\mathtt{p}}+1)\varkappa_{m}.
				\end{align*}
			Notice that since $\overline{s}_h+\widetilde{\mathtt{p}}+2\leqslant s_h+\sigma_{1}$, then by \eqref{estimate delta0 and I0} and \eqref{smallness condition transport}, one has
			\begin{align*}
				\delta_0^{1,2}(\overline{s}_h+\widetilde{\mathtt{p}}+1)&\lesssim\varepsilon\gamma^{-1}\left(1+\max_{k\in\{1,2\}}\|\mathfrak{I}_k\|_{q,\overline{s}_h+\widetilde{\mathtt{p}}+2}^{\gamma,\mathcal{O}}\right)\\
			&\lesssim\varepsilon\gamma^{-1}.
		\end{align*}
				It follows from  \eqref{hypothesis of induction overline delta} and \eqref{est-vchim}, 
				\begin{align*}
					\overline{\delta}_{m+1}(\overline{s}_{h}+\widetilde{\mathtt{p}})&\leqslant\left(2-\tfrac{1}{m+1}\right)\left(1+N_m^{s_0+\tau_{1}q+\tau_{1}-\overline{s}_h}+CN_{0}^{\overline{\mu}_{2}}N_{m}^{\tau_{1}q+\tau_{1}+1-\overline{\mu}_{2}}\right)\nu(\overline{s}_{h}+\widetilde{\mathtt{p}})\\
					&\quad+C\left(N_m^{\tau_{1}q+\tau_{1}+1}+N_0^{\overline{\mu}_2}N_m^{2\tau_{1}q+3\tau_{1}+1-\overline{\mu}_2}\right)N_0^{\overline{\mu}_2}N_m^{-\overline{\mu}_2}\varepsilon\gamma^{-1}\nu(\overline{s}_h+\widetilde{\mathtt{p}}).
				\end{align*}
				Proceeding as for \eqref{conv-t2}, taking $\varepsilon_{0}$ small enough and thanks to   \eqref{param-trans}, we obtain
				\begin{align*}
					&\left(2-\tfrac{1}{m+1}\right)\left(1+N_m^{s_0+\tau_{1}q+\tau_{1}-\overline{s}_h}+CN_{0}^{\overline{\mu}_{2}}N_{m}^{\tau_{1}q+\tau_{1}+1-\overline{\mu}_{2}}\right)\\
					&+C\left(N_m^{\tau_{1}q+\tau_{1}+1}+N_0^{\overline{\mu}_2}N_m^{2\tau_{1}q+3\tau_{1}+1-\overline{\mu}_2}\right)N_0^{\overline{\mu}_2}N_m^{-\overline{\mu}_2}\varepsilon\gamma^{-1}\\
					&\leqslant 2-\tfrac{1}{m+2},
				\end{align*}
				so that
				$$\overline{\delta}_{m+1}(\overline{s}_{h}+\widetilde{\mathtt{p}})\leqslant\left(2-\tfrac{1}{m+2}\right)\nu(\overline{s}_{h}+\widetilde{\mathtt{p}}).$$
				This completes the proof of the second statement in \eqref{hypothesis of induction overline delta}.\\
				\ding{226} \textit{Conclusion.}
				From \eqref{estimate delta12 gm}, we get for $s=s_{0}.$
				\begin{align*}
					\|\Delta_{12}g_{m}\|_{q,s_{0}}^{\gamma,\mathcal{O}}&\lesssim \overline{\delta}_{m}(s_{0}+\tau_{1}q+\tau_{1})+\varkappa_{m}\delta_{m}(s_{0}+\tau_{1}q+2\tau_{1}+1).
				\end{align*}
				By interpolation inequality in Lemma \ref{Lem-lawprod}, \eqref{hypothesis of induction deltam} applied with $\mu_{2}=\overline{\mu}_{2}$, \eqref{hypothesis of induction overline delta} applied with $\widetilde{\mathtt{p}}=0$ and Sobolev embeddings, we have for some $\overline{\theta}\in(0,1)$
				\begin{align*}
					\overline{\delta}_{m}(s_{0}+\tau_{1}q+\tau_{1})&\leqslant  \overline{\delta}_{m}(s_{0}+\tau_{1}q+2\tau_{1}+1)\\
					&\lesssim\overline{\delta}_{m}(s_{0})^{\overline{\theta}}\overline{\delta}_{m}(\overline{s}_{h})^{1-\overline{\theta}}\\
					&\lesssim N_{0}^{\overline{\theta}\overline{\mu}_{2}}N_{m}^{-\overline{\theta}\overline{\mu}_{2}}\nu(\overline{s}_{h})
				\end{align*}
				and
				\begin{align*}
					\delta_{m}(s_{0}+\tau_{1}q+2\tau_{1}+1)&\lesssim \delta_{m}(s_{0})^{\overline{\theta}}\delta_{m}(\overline{s}_{h})^{1-\overline{\theta}}\\
					&\lesssim N_{0}^{\overline{\theta}\overline{\mu}_{2}}N_{m}^{-\overline{\theta}\overline{\mu}_{2}}\delta_{0}(\overline{s}_{h})\\
					&\lesssim N_{0}^{\overline{\theta}\overline{\mu}_{2}}N_{m}^{-\overline{\theta}\overline{\mu}_{2}}.
				\end{align*}
				Therefore
				$$\|\Delta_{12}g_{m}\|_{q,s_{0}}^{\gamma,\mathcal{O}}\lesssim N_{0}^{\overline{\theta}\overline{\mu}_{2}}N_{m}^{-\overline{\theta}\overline{\mu}_{2}}\nu(\overline{s}_{h}).$$
				Now from \eqref{estimate delta12 gm}, we have 
				\begin{align*}
					\|\Delta_{12}g_{m}\|_{q,\overline{s}_{h}+\mathtt{p}+1}^{\gamma,\mathcal{O}}&\lesssim \overline{\delta}_{m}(\overline{s}_{h}+\mathtt{p}+\tau_{1}q+\tau_{1}+1)+\varkappa_{m}\delta_{m}(\overline{s}_{h}+\mathtt{p}+\tau_{1}q+2\tau_{1}+2).
				\end{align*}
			Applying \eqref{hypothesis of induction overline delta} with \begin{equation}\label{link pt p}
				\widetilde{\mathtt{p}}=\mathtt{p}+\tau_{1}q+\tau_{1}+1,
			\end{equation}
		 which is possible since from \eqref{param}, \eqref{param-trans} and \eqref{sigma1}, one has\quad 
			$\overline{s}_{h}+\mathtt{p}+\tau_{1}q+\tau_{1}+4\leqslant s_{h}+\sigma_{1},$
			we find
			\begin{align*}
				\overline{\delta}_{m}(\overline{s}_{h}+\mathtt{p}+\tau_{1}q+\tau_{1}+1)&\leqslant 2\nu(\overline{s}_{h}+\mathtt{p}+\tau_{1}q+\tau_{1}+1)\\
				&\leqslant 2\overline{\delta}_{0}(\overline{s}_{h}+\mathtt{p}+\tau_{1}q+\tau_{1}+1)+2\varepsilon\gamma^{-1}\|\Delta_{12}i\|_{q,s_0+2}^{\gamma,\mathcal{O}}.
			\end{align*}
		Implementing a similar proof to    \eqref{Estim-Vr-October} based on  the kernel decomposition \eqref{fundamental decomposition of K0(Ar)},    the  composition laws \mbox{and \eqref{control difference r},} we find
				\begin{align*}
			\forall s\geqslant s_0,\quad \overline{\delta}_{0}(s)&=\gamma^{-1}\|\Delta_{12}V_{\varepsilon r}\|_{q,s}^{\gamma,\mathcal{O}}\\
			&\lesssim\varepsilon\gamma^{-1}\Big(\|\Delta_{12}i\|_{q,s+1}^{\gamma,\mathcal{O}}+\|\Delta_{12}i\|_{q,s_0+1}^{\gamma,\mathcal{O}}\max_{\ell=1,2}\|r_\ell\|_{q,s+1}^{\gamma,\mathcal{O}}\Big).
		\end{align*}
%	Consequently, we get in view of \eqref{param-trans} and \eqref{smallness condition transport}
%	\begin{align*}
%		\overline{\delta}_{m}(\overline{s}_{h}+\mathtt{p}+\tau_{1}q+\tau_{1}+1)\lesssim \varepsilon\gamma^{-1}\|\Delta_{12}i\|_{q,\overline{s}_{h}+\mathtt{p}+\tau_{1}q+\tau_{1}+3}^{\gamma,\mathcal{O}}.
%	\end{align*}
On the other hand, since
$$\overline{s}_{h}+\mathtt{p}+\tau_{1}q+2\tau_{1}+3\leqslant s_{h}+\sigma_{1},$$
one may obtain through combining \eqref{uniform estimate of deltams} and \eqref{smallness condition transport}
\begin{align*}
	\delta_{m}(\overline{s}_{h}+\mathtt{p}+\tau_{1}q+2\tau_{1}+2)&\leqslant C\varepsilon\gamma^{-1}\left(1+\|\mathfrak{I}_{0}\|_{q,\overline{s}_{h}+\mathtt{p}+\tau_{1}q+2\tau_{1}+3}^{\gamma,\mathcal{O}}\right)\\
	&\leqslant C\varepsilon\gamma^{-1}.
\end{align*}
				Thus, by interpolation inequality in Lemma \ref{Lem-lawprod}, we finally obtain for some  $\overline{\theta}\in(0,1)$
				\begin{equation}\label{final dlt12 gm shb+1}
					\|\Delta_{12}g_{m}\|_{q,\overline{s}_{h}+\mathtt{p}}^{\gamma,\mathcal{O}}\lesssim N_{0}^{\overline{\theta}\overline{\mu}_{2}}N_{m}^{-\overline{\theta}\overline{\mu}_{2}}\nu(\overline{s}_{h}+\mathtt{p}+\tau_{1}q+\tau_{1}+1).
				\end{equation}
				Choosing $N_{0}$ sufficiently large, then the composition law in Lemma \ref{Lem-lawprod} allows to get
				\begin{align}\label{final sum dlt12 gm shb+1}
					\sum_{k=0}^{\infty}\|\Delta_{12}g_{k}\|_{q,\overline{s}_{h}+\mathtt{p}}^{\gamma,\mathcal{O}}&\lesssim \nu(\overline{s}_{h}+\mathtt{p}+\tau_{1}q+\tau_{1}+1)N_{0}^{\overline{\theta}\overline{\mu}_{2}}\sum_{k=0}^{\infty}N_{m}^{-\overline{\theta}\overline{\mu}_{2}}\nonumber\\
					&\lesssim\varepsilon\gamma^{-1}\|\Delta_{12}i\|_{q,\overline{s}_{h}+\mathtt{p}+\tau_{1}q+\tau_{1}+2}^{\gamma,\mathcal{O}}.
				\end{align}
				Finally, gathering \eqref{estimate series Delta12 beta}, \eqref{telscopic diff betam}, \eqref{final dlt12 gm shb+1} and \eqref{final sum dlt12 gm shb+1}, we get
				$$\|\Delta_{12}\beta\|_{q,\overline{s}_{h}+\mathtt{p}}^{\gamma,\mathcal{O}}\lesssim\varepsilon\gamma^{-1}\|\Delta_{12}i\|_{q,\overline{s}_{h}+\mathtt{p}+\tau_{1}q+\tau_{1}+2}.$$
				Putting together this estimate,  \eqref{link diff beta hat and diff beta} and \eqref{maj betamrk} yields 
				\begin{align*}
					\|\Delta_{12}\widehat{\beta}\|_{q,\overline{s}_{h}+\mathtt{p}}^{\gamma,\mathcal{O}}&\lesssim\|\Delta_{12}\beta\|_{q,\overline{s}_{h}+\mathtt{p}}^{\gamma,\mathcal{O}}\\
					&\lesssim\varepsilon\gamma^{-1}\|\Delta_{12}i\|_{q,\overline{s}_{h}+\mathtt{p}+\tau_{1}q+\tau_{1}+2}.
				\end{align*}
				$\blacktriangleright$ \textbf{Estimate on $\Delta_{12}c_{i}.$} Since $V_{0}=\Omega+I_{1}K_{1}$ is independent of $r$, then
				$$
				\Delta_{12}c_{i}=\sum_{m=0}^{\infty}\Delta_{12}(V_{m+1}-V_{m}).$$
				Therefore we obtain in view of  \eqref{definition Vm+1 and fm+1}, Sobolev embeddings and \eqref{hypothesis of induction overline delta} applied with $\widetilde{\mathtt{p}}=0$, 
				\begin{align*}
					\|\Delta_{12}(V_{m+1}-V_{m})\|_{q}^{\gamma,\mathcal{O}}&=\|\langle \Delta_{12}f_{m}\rangle_{\varphi,\theta}\|_{q}^{\gamma,\mathcal{O}}\\
					&\leqslant C\gamma\overline{\delta}_{m}(s_{0})\\
					&\leqslant C\gamma N_{0}^{\overline{\mu}_{2}}N_{m}^{-\overline{\mu}_{2}}\nu(\overline{s}_{h}).
				\end{align*}
				Hence by the composition law in Lemma \ref{Lem-lawprod}, Lemma \ref{lemma sum Nn} and \eqref{control difference r} one may find
				\begin{align*}
					\|\Delta_{12}c_{i}\|_{q}^{\gamma,\mathcal{O}}&\leqslant \sum_{m=0}^{\infty}\|\Delta_{12}(V_{m+1}-V_{m})\|_{q}^{\gamma,\mathcal{O}}\\
					& \leqslant  C\gamma \nu(\overline{s}_{h})N_{0}^{\overline{\mu}_{2}}\sum_{m=0}^{\infty}N_{m}^{-\overline{\mu}_{2}}\\
					& \leqslant C\varepsilon\|\Delta_{12}i\|_{q,\overline{s}_{h}+2}^{\gamma,\mathcal{O}}.
				\end{align*}
				This achieves the proof of  Proposition \ref{reduction of the transport part}. 
			\end{proof}
					\subsubsection{Action on the non-local term}
			In this section, we shall analyze   the conjugation action  by $\mathscr{B}$ on the nonlocal term appearing in the linearized operator $\mathcal{L}_{\varepsilon r}$ described in Proposition \ref{lemma setting for Lomega}. The main result reads as follows.
			\begin{prop}\label{Action on the non local part}
				Let $(\gamma,q,d,\tau_{1},s_{0},\overline{s}_{h}, \sigma_{1},S)$ satisfy \eqref{initial parameter condition}, \eqref{setting tau1 and tau2}, \eqref{param} and \eqref{sigma1}. We set 
				\begin{equation}\label{sigma2}
					\sigma_2=s_0+\sigma_1+3.
				\end{equation}
				For any $(\mu_{2},\mathtt{p},s_h)$ satisfying the condition \eqref{param-trans}, there exists $\varepsilon_0>0$ such that if 
				\begin{equation}\label{small-2}
					\varepsilon\gamma^{-1}N_0^{\mu_{2}}\leqslant\varepsilon_0\quad\textnormal{and}\quad\|\mathfrak{I}_0\|_{q,s_h+\sigma_{2}}^{\gamma,\mathcal{O}}\leqslant1,
				\end{equation}
				then in the Cantor set $\mathcal{O}_{\infty,n}^{\gamma,\tau_{1}}(i_{0})$, we have 
				$$\mathfrak{L}_{\varepsilon r}:=\mathscr{B}^{-1}\mathcal{L}_{\varepsilon r}\mathscr{B}=\omega\cdot\partial_{\varphi}+c_{i_{0}}\partial_{\theta}-\partial_{\theta}\mathcal{K}_{\lambda}\ast\cdot+\partial_{\theta}\mathfrak{R}_{\varepsilon r}+\mathtt{E}_{n}^{0},$$
				where $\mathcal{K}_{\lambda}$ is defined in \eqref{definition of mathcalKlambda}, $\mathtt{E}_{n}^{0}$ is introduced  in Proposition $\ref{reduction of the transport part}$ and $\mathfrak{R}_{\varepsilon r}$ is a real and reversibility preserving self-adjoint integral operator satisfying
				\begin{equation}\label{estim partial thetamathfrakRr}
					\forall s\in[s_{0},S],\quad \max_{k\in\{0,1,2\}}\| \partial_{\theta}^k\mathfrak{R}_{\varepsilon r}\|_{\textnormal{\tiny{O-d}},q,s}^{\gamma,\mathcal{O}}\lesssim \varepsilon\gamma^{-1}\left(1+\|\mathfrak{I}_{0}\|_{q,s+\sigma_{2}}^{\gamma,\mathcal{O}}\right).
				\end{equation}
				In addition, if $i_{1}$ and $i_{2}$ are two tori satisfying the smallness property \eqref{small-2}, then 
				\begin{equation}\label{differences partial thetamathfrakRr}
					\max_{k\in\{0,1\}}\|\Delta_{12}\partial_{\theta}^{k}\mathfrak{R}_{\varepsilon r}\|_{\textnormal{\tiny{O-d}},q,\overline{s}_{h}+\mathtt{p}}^{\gamma,\mathcal{O}}\lesssim\varepsilon\gamma^{-1}\|\Delta_{12}i\|_{q,\overline{s}_{h}+\mathtt{p}+\sigma_{2}}^{\gamma,\mathcal{O}}.
				\end{equation}
			\end{prop}
			\begin{proof}
				We recall from Proposition \ref{lemma setting for Lomega} and Lemma \ref{lemma general form of the linearized operator}. that 
				$$\mathcal{L}_{\varepsilon r}=\omega\cdot\partial_{\varphi}+\partial_{\theta}\left(V_{\varepsilon r}\cdot\right)-\partial_{\theta}\mathbf{L}_{\varepsilon r},$$
				where $\mathbf{L}_{\varepsilon r}$ is a non-local operator defined by 
				$$\mathbf{L}_{\varepsilon r}(\rho)(\varphi,\theta)=\int_{\mathbb{T}}\rho(\varphi,\eta)K_{0}(\lambda A_{\varepsilon r}(\varphi,\theta,\eta))d\eta,$$
				with 
				$$A_{\varepsilon r}(\varphi,\theta,\eta)=\left(\big(R(\varphi,\eta)-R(\varphi,\theta)\big)^{2}+4R(\varphi,\eta)R(\varphi,\theta)\sin^{2}\left(\tfrac{\eta-\theta}{2}\right)\right)^{\frac{1}{2}}$$
				and
				$$R(\varphi,\theta)=\big(1+2\varepsilon r(\varphi,\theta)\big)^{\frac{1}{2}}.$$
				Notice that we have removed the dependance in $(\lambda,\omega)$ from the  functions in order to alleviate the notation. Hence by Proposition \ref{reduction of the transport part}, Lemma \ref{algeb1}-(i) and \eqref{link non local operators at state r and at equilibrium}, we have in the Cantor set $\mathcal{O}_{\infty,n}^{\gamma,\tau_{1}}(i_{0})$
				\begin{align}\label{Ler}
						\mathfrak{L}_{\varepsilon r} & :=  \mathscr{B}^{-1}\mathcal{L}_{\varepsilon r}\mathscr{B}\nonumber\\
						& =  \mathscr{B}^{-1}\big(\omega\cdot\partial_{\varphi}+\partial_{\theta}\left(V_{\varepsilon r}\cdot\right)\big)\mathscr{B}-\mathscr{B}^{-1}\partial_{\theta}\mathbf{L}_{\varepsilon r}\mathscr{B}\nonumber\\
						& =  \omega\cdot\partial_{\varphi}+c_{i_{0}}\partial_{\theta}-\partial_{\theta}\mathcal{B}^{-1}\mathbf{L}_{\varepsilon r}\mathscr{B}+\mathtt{E}_{n}^{0}\nonumber\\
						& =  \omega\cdot\partial_{\varphi}+c_{i_{0}}\partial_{\theta}-\partial_{\theta}\Big(\mathcal{B}^{-1}\left(\mathcal{K}_{\lambda}\ast\cdot\right)\mathscr{B}+\mathcal{B}^{-1}\mathbf{L}_{\varepsilon r,1}\mathscr{B}\Big)+\mathtt{E}_{n}^{0}.
				\end{align}
				From a direct computation using \eqref{definition of mathcalKlambda} combined with  \eqref{mathscrB1} and \eqref{definition symplectic change of variables},  we find
				$$\mathcal{B}^{-1}\big(\mathcal{K}_{\lambda}\ast\mathscr{B}\rho\big)(\varphi,\theta)=\int_{\mathbb{T}}\rho(\varphi,\eta)K_{0}\big(\lambda\mathscr{A}_{\widehat\beta}(\varphi,\theta,\eta)\big)d\eta,$$
				where 
				$$\mathscr{A}_{\widehat\beta}(\varphi,\theta,\eta):=2\left|\sin\left(\tfrac{\eta-\theta}{2}+\widehat{h}(\varphi,\theta,\eta)\right)\right|,$$
				with
				$$\widehat{h}(\varphi,\theta,\eta):=\tfrac{\widehat{\beta}(\varphi,\eta)-\widehat{\beta}(\varphi,\theta)}{2}\cdot
				$$
				Using elementary trigonometric identities, we can write
				\begin{align}\label{TU-1}\mathscr{A}_{\widehat\beta}(
					\varphi,\theta,\eta)=2\left|\sin\left(\tfrac{\eta-\theta}{2}\right)\right|v_{\widehat\beta,2}(\varphi,\theta,\eta),
				\end{align}
				with
				$$v_{\widehat\beta,2}(\varphi,\theta,\eta):=\cos\left(\widehat{h}(\varphi,\theta,\eta)\right)+\tfrac{\sin\left(\widehat{h}(\varphi,\theta,\eta)\right)}{\tan\left(\frac{\eta-\theta}{2}\right)}\cdot
				$$
				Notice that $v_{0,2}=1$ and one may write
				$$v_{\widehat\beta,2}(\theta,\eta)=1+\big(\cos\big(\widehat{h}(\theta,\eta)\big)-1\big)+\tfrac{\widehat{h}(\theta,\eta)}{\tan\left(\frac{\eta-\theta}{2}\right)}+\left(\tfrac{\sin(\widehat{h}(\theta,\eta))}{\widehat{h}(\theta,\eta)}-1\right)\tfrac{\widehat{h}(\theta,\eta)}{\tan\left(\frac{\eta-\theta}{2}\right)}$$
				and then using Lemma \ref{Lem-lawprod}-(iv)-(v), Lemma \ref{cheater lemma} and  \eqref{estimate beta and r},  we obtain
				\begin{align}\label{estimate vr2}
					\nonumber\sup_{\eta\in\mathbb{T}}\big\| v_{\widehat\beta,2}(\ast,\cdot,\centerdot,\eta+\centerdot)-1\big\|_{q,s}^{\gamma,\mathcal{O}}&\lesssim\|\widehat{\beta}\|_{q,s+1}^{\gamma,\mathcal{O}}\lesssim\varepsilon\gamma^{-1}\left(1+\|\mathfrak{I}_{0}\|_{q,s+\sigma_{1}+1}^{\gamma,\mathcal{O}}\right),\\
					\forall k\in\mathbb{N}^{*},\,\sup_{\eta\in\mathbb{T}}\big\| (\partial_{\theta}^{k}v_{\widehat\beta,2})(\ast,\cdot,\centerdot,\eta+\centerdot)\big\|_{q,s}^{\gamma,\mathcal{O}}&\lesssim\|\widehat{\beta}\|_{q,s+k+1}^{\gamma,\mathcal{O}}\lesssim\varepsilon\gamma^{-1}\left(1+\|\mathfrak{I}_{0}\|_{q,s+\sigma_{1}+1+k}^{\gamma,\mathcal{O}}\right).
				\end{align}
				Proceeding as for \eqref{fundamental decomposition of K0(Ar)}, one obtains the decomposition
				$$K_{0}(\lambda\mathscr{A}_{\widehat\beta}(\lambda, \omega,\varphi,\theta,\eta))=K_{0}\left(2\lambda\left|\sin\left(\tfrac{\eta-\theta}{2}\right)\right|\right)+\mathscr{K}(\eta-\theta)\mathscr{K}_{\widehat\beta,2}^{1}(\varphi,\theta,\eta)+\mathscr{K}_{\widehat\beta,2}^{2}(\varphi,\theta,\eta)$$
				with  similar estimates to \eqref{estimate kernel Kr11} and \eqref{estimate kernel Kr12}, that is,  for all $k\in{\mathbb{N}},$
				\begin{align}\label{estimates kernels Kr21 and Kr22}
					\nonumber\sup_{\eta\in\mathbb{T}}\left(\big\|(\partial_{\theta}^{k}\mathscr{K}_{\widehat\beta,2}^{1})(\ast,\cdot,\centerdot,\eta+\centerdot)\big\|_{q,s}^{\gamma,\mathcal{O}}+\big\|(\partial_{\theta}^{k}\mathscr{K}_{\widehat\beta,2}^{2})(\ast,\cdot,\centerdot,\eta+\centerdot)\big\|_{q,s}^{\gamma,\mathcal{O}}\right)&\lesssim\| \widehat\beta\|_{q,s+1+k}^{\gamma,\mathcal{O}}\\
					&\lesssim \varepsilon\gamma^{-1}\left(1+\|\mathfrak{I}_{0}\|_{q,s+\sigma_{1}+1+k}^{\gamma,\mathcal{O}}\right),
				\end{align}
				where the symbols $\ast,\cdot,\centerdot$ stand for $(\lambda,\omega),\varphi,\theta$, respectively.
				Now we shall denote by  $\mathbf{L}_{\varepsilon r,2}$  the integral operator with the  kernel $\mathbb{K}_{\varepsilon r,2}$ defined by 
				\begin{equation}\label{expression kernel K2}
					\mathbb{K}_{\varepsilon r,2}(\varphi,\theta,\eta):=\mathscr{K}(\eta-\theta)\mathscr{K}_{\widehat\beta,2}^{1}(\varphi,\theta,\eta)+\mathscr{K}_{\widehat\beta,2}^{2}(\varphi,\theta,\eta).
				\end{equation}
				Then  we find the decomposition
				$$\mathcal{B}^{-1}\left(\mathcal{K}_{\lambda}\ast\cdot\right)\mathscr{B}=\mathcal{K}_{\lambda}\ast\cdot+\mathbf{L}_{\varepsilon r,2}.
				$$
				Inserting this identity into \eqref{Ler} allows to get  
				$$\mathfrak{L}_{\varepsilon r}=\mathscr{B}^{-1}\mathcal{L}_{\varepsilon r}\mathscr{B}=\omega\cdot\partial_{\varphi}+c_{i_{0}}\partial_{\theta}-\partial_{\theta}\mathcal{K}_{\lambda}\ast\cdot+\partial_{\theta}\mathfrak{R}_{\varepsilon r}+\mathtt{E}_{n}^{0},$$
				with
				\begin{equation}\label{definition of mathfrakRr}
					\mathfrak{R}_{\varepsilon r}=-\mathbf{L}_{\varepsilon r,2}-\mathcal{B}^{-1}\mathbf{L}_{\varepsilon r,1}\mathscr{B}.
				\end{equation}
				Observe that by  \eqref{symmetry for r} and \eqref{symmetry for beta} we can easily check that  the kernel $\mathbb{K}_{\varepsilon r,2}$ satisfies the following symmetry property
				\begin{equation}\label{symmetry kernel K2}
					\mathbb{K}_{\varepsilon r,2}(-\varphi,-\theta,-\eta)=\mathbb{K}_{\varepsilon r,2}(\varphi,\theta,\eta)\in\mathbb{R},
				\end{equation}
				which implies in turn, according to Lemma \ref{lemma symmetry and reversibility}, that $\mathbf{L}_{\varepsilon r,2}$ is a real and reversibility preserving operator. Moreover, one obtains from \eqref{estimates kernels Kr21 and Kr22}
				\begin{align}\label{est-Yu1}
					\max_{k\in\{0,1,2\}}\,\big\|(\partial_{\theta}^{k}\mathbb{K}_{\varepsilon r,2})(\ast,\cdot,\centerdot,\eta+\centerdot)\big\|_{q,s}^{\gamma,\mathcal{O}}&\lesssim \varepsilon\gamma^{-1}\left(1+\|\mathfrak{I}_{0}\|_{q,s+\sigma_{1}+3}^{\gamma,\mathcal{O}}\right)\,\left(1-\log\left|\sin\left(\tfrac{\eta}{2} \right)\right|\right).
				\end{align}
				Our next purpose is  to highlight  some properties of the operator $\mathcal{B}^{-1}\mathbf{L}_{\varepsilon r,1}\mathscr{B}$ which takes the integral form
				\begin{align}\label{ZY1}\big(\mathcal{B}^{-1}\mathbf{L}_{\varepsilon r,1}\mathscr{B}\big)\rho(\varphi,\theta)=\int_{\mathbb{T}}\rho(\varphi,\eta)\widehat{\mathbb{K}}_{\varepsilon r,1}(\varphi,\theta,\eta)d\eta,
				\end{align}
				where the kernel $\widehat{\mathbb{K}}_{\varepsilon r,1}$  is related to the kernel $\mathbb{K}_{\varepsilon r,1}$ defined in  \eqref{definition of mathbbK1} through the formula,
				\begin{align}\label{KKLLPP1}\widehat{\mathbb{K}}_{\varepsilon r,1}(\varphi,\theta,\eta)=\mathbb{K}_{\varepsilon r,1}\big(\varphi,\theta+\widehat{\beta}(\varphi,\theta),\eta+\widehat{\beta}(\varphi,\eta)\big).
				\end{align}
				It is quite easy to check from  \eqref{symmetry kernel K1} and \eqref{symmetry for beta}, that
				\begin{equation}\label{symmetry kernel K1'}
					\widehat{\mathbb{K}}_{\varepsilon r,1}(-\varphi,-\theta,-\eta)=\widehat{\mathbb{K}}_{\varepsilon r,1}(\varphi,\theta,\eta)\in\mathbb{R}.
				\end{equation}
				According to 
				\eqref{definition of mathbbK1}, one gets the decomposition
				\begin{equation}\label{Sam-d1}
					\widehat{\mathbb{K}}_{\varepsilon r,1}(\varphi,\theta,\eta)=\widehat{\mathscr{K}}(\varphi,\theta,\eta)\widehat{\mathscr{K}}_{\varepsilon r,1}^{\,\,1}(\varphi,\theta,\eta)+\widehat{\mathscr{K}}_{\varepsilon r,1}^{\,\,2}(\varphi,\theta,\eta),
				\end{equation}
			with
			\begin{align*}
				\widehat{\mathscr{K}}(\varphi,\theta,\eta)&:=\mathscr{K}\big(\eta-\theta+\widehat{\beta}(\varphi,\eta)-\widehat{\beta}(\varphi,\theta)\big),\\
				\widehat{\mathscr{K}}_{\varepsilon r,1}^{\,\,1}(\varphi,\theta,\eta)&:=\mathscr{K}_{\varepsilon r,1}^{1}\big(\varphi,\theta+\widehat{\beta}(\varphi,\theta),\eta+\widehat{\beta}(\varphi,\eta)\big),\\
				\widehat{\mathscr{K}}_{\varepsilon r,1}^{\,\,2}(\varphi,\theta,\eta)&:=\mathscr{K}_{\varepsilon r,1}^{2}\big(\varphi,\theta+\widehat{\beta}(\varphi,\theta),\eta+\widehat{\beta}(\varphi,\eta)\big).
			\end{align*}
				Coming back to 
				\eqref{definition of mathscrK} and using the morphism property of the logarithm, combined with \eqref{TU-1} we deduce that
				\begin{align*}
					\widehat{\mathscr{K}}(\varphi,\theta,\eta)&=\sin^2\left(\tfrac{\eta-\theta}{2} \right) v_{\widehat{\beta},2}^2(\varphi,\theta,\eta)\left(\log\Big|\sin\Big(\tfrac{\eta-\theta}{2} \Big) \Big|+\log \big|v_{\widehat{\beta},2}(\varphi,\theta,\eta)\big|\right)\\
					&=\mathscr{K}\big(\eta-\theta\big)+\mathscr{K}\big(\eta-\theta\big)\Big( v_{\widehat{\beta},2}^2(\varphi,\theta,\eta)-1 \Big)\\
					&\quad  +\sin^2\left(\tfrac{\eta-\theta}{2} \right) v_{\widehat{\beta},2}^2(\varphi,\theta,\eta)\log \big|v_{\widehat{\beta},2}(\varphi,\theta,\eta)\big|.
				\end{align*}
				Combining Lemma \ref{Lem-lawprod}-(iv)-(v), \eqref{estimate vr2}, \eqref{estimate beta and r} gives for any $\eta\in\T$
				\begin{align}\label{Es-U-D0}
					\nonumber \max_{k\in\{0,1,2\}} \big\|(\partial_\theta^k\widehat{\mathscr{K}})\big(\ast,\cdot,\centerdot,&\eta+\centerdot\big)\big\|_{q,s}^{\gamma,\mathcal{O}}\lesssim\,\big\|\widehat\beta\|_{q,s+3}^{\gamma,\mathcal{O}}\, \left(1-\log\left|\sin\left(\frac{\eta}{2} \right)\right|\right){-\log\left|\sin\left(\tfrac{\eta}{2}\right)\right|}+1\\
					&\lesssim \varepsilon\gamma ^{-1}\left(1+\| \mathfrak{I}_{0}\|_{q,s+\sigma_{2}}^{\gamma ,\mathcal{O}}\right)\left(1-\log\left|\sin\left(\frac{\eta}{2} \right)\right|\right){-\log\left|\sin\left(\tfrac{\eta}{2}\right)\right|}+1.
				\end{align}			
				The next goal is to prove that 
				\begin{align}\label{Es-U-D}
					\max_{k\in\{0,1,2\}} \sup_{\eta\in\T}\big\|(\partial_\theta^k\widehat{\mathscr{K}}_{\varepsilon r,1}^{\,\,1})\big(\ast,\cdot,\centerdot,\eta+\centerdot\big)\big\|_{q,s}^{\gamma,\mathcal{O}}&\lesssim\, \varepsilon\gamma^{-1}\left(1+\| \mathfrak{I}_{0}\|_{q,s+\sigma_{2}}^{\gamma ,\mathcal{O}}\right). 
					\end{align}
				For this aim we first write from \eqref{ExpressionK} and \eqref{definition of modified Bessel function of first kind} 
				\begin{align}\label{ExpressionKML}
					\nonumber \mathscr{K}_{\varepsilon r,1}^{1}(\varphi,\theta,\eta)&=4\lambda^2\big(1-v_{\varepsilon r,1}(\varphi,\theta,\eta)\big)\tilde{I}_\lambda(\eta-\theta)\\
					\nonumber&\quad-4\lambda^{2}(v_{\varepsilon r,1}(\varphi,\theta)-1)^{2}\int_{0}^{1}(1-t)I_{0}''\left(2\lambda\sin\left(\tfrac{\eta-\theta}{2}\right)(1-t+tv_{\varepsilon r,1}(\varphi,\theta,\eta))\right)dt\\
					&:=4\lambda^2\big(1-v_{\varepsilon r,1}(\varphi,\theta,\eta)\big)\tilde{I}_\lambda(\eta-\theta)+G(\varphi,\theta,\eta),
				\end{align}
				with 
				\begin{align*}
					\tilde{I}_\lambda(\eta)&:=\frac{I_{0}'\left(2\lambda\left|\sin\left(\frac{\eta}{2}\right)\right|\right)}{2\lambda\left|\sin\left(\frac{\eta}{2}\right)\right|}\\
					&=
					\frac12\sum_{m=0}^{\infty}\frac{\lambda^{2m} \sin^{2m}\left(\frac{\eta}{2}\right)}{m!(m+1)!}\cdot
				\end{align*}
				Then we get the decomposition
				$$\widehat{\mathscr{K}}_{\varepsilon r,1}^{\,\,1}(\varphi,\theta,\eta)=4\lambda^2\big[1-\widehat{v}_{\varepsilon r,1}(\varphi,\theta,\eta)\big]\widehat{I}_{\lambda}(\varphi,\theta,\eta)+\widehat{G}(\varphi,\theta,\eta),
				$$
			with 
			\begin{align*}
				\widehat{v}_{\varepsilon r,1}(\varphi,\theta,\eta)&:=v_{\varepsilon r,1}\big(\varphi,\theta+\widehat{\beta}(\varphi,\theta),\eta+\widehat{\beta}(\varphi,\eta)\big),\\
				\widehat{I}_{\lambda}(\varphi,\theta,\eta)&:=\tilde{I}_\lambda\big(\eta+\widehat{\beta}(\varphi,\eta)-\theta-\widehat{\beta}(\varphi,\theta)\big),\\
				\widehat{G}(\varphi,\theta,\eta)&:=G\big(\varphi,\theta+\widehat{\beta}(\varphi,\theta),\eta+\widehat{\beta}(\varphi,\eta)\big).
			\end{align*}
				It follows that
				\begin{align}\label{ExpressionKMM}
					\nonumber \widehat{\mathscr{K}}_{\varepsilon r,1}^{\,\,1}(\varphi,\theta,\theta+\eta)=&4\lambda^2\big[1-v_{\varepsilon r,1}\big(\varphi,\theta+\widehat{\beta}(\varphi,\theta),\theta+\eta+\widehat{\beta}(\varphi,\theta+\eta)\big)\big]\tilde{I}_\lambda\big(\eta+\widehat{\beta}(\varphi,\theta+\eta)-\widehat{\beta}(\varphi,\theta)\big)\\
					&+G\big(\varphi,\theta+\widehat{\beta}(\varphi,\theta),\eta+\theta+\widehat{\beta}(\varphi,\eta+\theta)\big).
				\end{align}
				Notice that $(\lambda,\eta)\mapsto \tilde{I}_\lambda(\eta)$ is $\mathcal{C}^\infty$, then using Lemma \ref{Lem-lawprod}-(v) and \eqref{estimate beta and r} yields for any $k\in\mathbb{N}$
				\begin{align*}
					\sup_{\eta\in\T}\|(\partial_{\theta}^{k}\widehat{I}_{\lambda})(\ast,\cdot,\centerdot,\eta+\centerdot)\|_{q,s}^{\gamma,\mathcal{O}}&\lesssim 1+\|\widehat{\beta}\|_{q,s+k}^{\gamma,\mathcal{O}}\\
					&\lesssim 1+\varepsilon\gamma^{-1}\| \mathfrak{I}_{0}\|_{q,s+\sigma_{1}+k}^{\gamma ,\mathcal{O}}.
				\end{align*}
				Now using \eqref{expression of Ar with vr1}, Lemma \ref{Lem-lawprod}-(v), \eqref{estimate on the first reduction operator and its inverse}, \eqref{estimate r and mathfrakI0}, \eqref{estimate vr2} and proceeding as in \eqref{estimate vr1} we obtain
				\begin{align*}%\label{estimate vr1M}
					\sup_{\eta\in\mathbb{T}}\|\widehat{v}_{\varepsilon r,1}\big(\ast,\cdot,\centerdot,\eta+\centerdot\big)-1\|_{q,s}^{\gamma,\mathcal{O}}&\lesssim \varepsilon\| r\|_{q,s+1}^{\gamma,\mathcal{O}}+\varepsilon^{2}\gamma^{-1}\|\mathfrak{I}_0\|_{q,s+\sigma_{1}+1}^{\gamma,\mathcal{O}}\|r\|_{q,s_0+1}^{\gamma,\mathcal{O}}\\
					&\lesssim\varepsilon\gamma^{-1}\left(1+\|\mathfrak{I}_{0}\|_{q,s+\sigma_{1}+1}^{\gamma,\mathcal{O}}\right)
				\end{align*}
			and 
			\begin{align*}%\label{estimate vr1M}
				\max_{k\in\{1,2\}}\sup_{\eta\in\mathbb{T}}\|(\partial_{\theta}^{k}\widehat{v}_{\varepsilon r,1})\big(\ast,\cdot,\centerdot,\eta+\centerdot\big)\|_{q,s}^{\gamma,\mathcal{O}}&\lesssim \varepsilon\| r\|_{q,s+3}^{\gamma,\mathcal{O}}+\varepsilon^{2}\gamma^{-1}\|\mathfrak{I}_0\|_{q,s+\sigma_{1}+3}^{\gamma,\mathcal{O}}\|r\|_{q,s_0+3}^{\gamma,\mathcal{O}}\\
				&\lesssim\varepsilon\gamma^{-1}\left(1+\|\mathfrak{I}_{0}\|_{q,s+\sigma_{1}+3}^{\gamma,\mathcal{O}}\right).
			\end{align*}
				Arguing as above using the structure of $G$ detailed in  \eqref{ExpressionKML} allows to get
				\begin{align*}
					\sup_{\eta\in\mathbb{T}}\|\widehat{G}\big(\ast,\cdot,\centerdot,\eta+\centerdot\big)-1\|_{q,s}^{\gamma,\mathcal{O}}&\lesssim\varepsilon\gamma^{-1}\left(1+\|\mathfrak{I}_{0}\|_{q,s+\sigma_{1}+1}^{\gamma,\mathcal{O}}\right)
				\end{align*}
				and 
				\begin{align*}
					\max_{k\in\{1,2\}}\sup_{\eta\in\mathbb{T}}\|(\partial_{\theta}^{k}\widehat{G}\big(\ast,\cdot,\centerdot,\eta+\centerdot\big)\|_{q,s}^{\gamma,\mathcal{O}}&\lesssim\varepsilon\gamma^{-1}\left(1+\|\mathfrak{I}_{0}\|_{q,s+\sigma_{1}+3}^{\gamma,\mathcal{O}}\right).
				\end{align*} 
%				\begin{equation*}%\label{estimate vr1M}
%					\sup_{\eta\in\mathbb{T}}\|G\big(\cdot,\centerdot+\widehat{\beta}(\cdot,\centerdot),\centerdot+\eta+\widehat{\beta}(\cdot,\centerdot+\eta)\big)-1\|_{q,s}^{\gamma,\mathcal{O}}\lesssim \varepsilon \| r\|_{q,s+1}^{\gamma,\mathcal{O}}.
%				\end{equation*}
				Thus applying  the law products in Lemma \ref{Lem-lawprod} and using  the preceding estimates combined with \eqref{ExpressionKMM} imply
				\begin{align}\label{estimate vr1MD}
					 \max_{k\in\{0,1,2\}}\sup_{\eta\in\mathbb{T}}\|\big(\partial_\theta^k\widehat{\mathscr{K}}_{\varepsilon r,1}^{\,\,1}\big)(\ast,\cdot,\centerdot,\eta+\centerdot)\|_{q,s}^{\gamma,\mathcal{O}}
					&\lesssim \varepsilon\gamma^{-1}\left(1+ \| \mathfrak{I}_{0}\|_{q,s+\sigma_{1}+3}^{\gamma ,\mathcal{O}}\right),
				\end{align}
				which gives in particular \eqref{Es-U-D}. 
				The estimate of  the last  term $\mathscr{K}_{\varepsilon r,1}^{2}$  in \eqref{Sam-d1}, which is connected to \eqref{definition of mathscrK12}, can be treated  in a similar way to the estimate \eqref{estimate vr1MD} and one finds
				\begin{align}\label{Es-U-D1}
					\max_{k\in\{0,1,2\}} \sup_{\eta\in\T}\big\|(\partial_\theta^k\widehat{\mathscr{K}}_{\varepsilon r,1}^{\,\,2})\big(\ast,\cdot,\centerdot,\eta+\centerdot\big)\big\|_{q,s}^{\gamma,\mathcal{O}}&\lesssim\, \varepsilon\gamma^{-1}\left(1+\| \mathfrak{I}_{0}\|_{q,s+\sigma_{1}+3}^{\gamma ,\mathcal{O}}\right). \end{align}
				Consequently, putting together \eqref{Sam-d1},\eqref{Es-U-D0}, \eqref{Es-U-D} and \eqref{Es-U-D1} yields
				\begin{align}\label{Es-U-D3}
					\max_{k\in\{0,1,2\}} \big\|(\partial_\theta^k\widehat{\mathbb{K}}_{\varepsilon r,1})\big(\ast,\cdot,\centerdot,\eta+\centerdot\big)\big\|_{q,s}^{\gamma,\mathcal{O}}&\lesssim \varepsilon\gamma ^{-1}\left(1+\| \mathfrak{I}_{0}\|_{q,s+\sigma_{1}+3}^{\gamma ,\mathcal{O}}\right)\left(1-\log\left|\sin\left(\tfrac{\eta}{2} \right)\right|\right).
				\end{align}
				By \eqref{definition of mathfrakRr} we infer  that $\mathfrak{R}_{\varepsilon r}$ is an integral operator of kernel $\mathbb{K}_{\varepsilon r}$ given  by
				$$\mathbb{K}_{\varepsilon r}=-\widehat{\mathbb{K}}_{\varepsilon r,1}-\mathbb{K}_{\varepsilon r,2}.
				$$
				Therefore, by virtue of  Lemma \ref{lemma symmetry and reversibility} combined with \eqref{est-Yu1} and \eqref{Es-U-D3} we find, taking $\sigma_{2}= s_0+\sigma_{1}+3,$
				\begin{align*}
					\max_{k\in\{0,1,2\}}\| \partial_\theta^k \mathfrak{R}_{\varepsilon r}\|_{\textnormal{\tiny{O-d}},q,s}^{\gamma,\mathcal{O}}&\lesssim \max_{k\in\{0,1,2\}}\bigintssss_{\T}\left(\| (\partial_\theta^k\widehat{\mathbb{K}}_{\varepsilon r,1})(\ast,\cdot,\centerdot,\eta+\centerdot)\|_{q,s+s_{0}}^{\gamma,\mathcal{O}}+\| (\partial_\theta^k\mathbb{K}_{\varepsilon r,2})(\ast,\cdot,\centerdot,\eta+\centerdot)\|_{q,s+s_{0}}^{\gamma,\mathcal{O}}\right) d\eta\\
					&\lesssim \varepsilon\gamma ^{-1}\left(1+\| \mathfrak{I}_{0}\|_{q,s+s_0+\sigma_{1}+3}^{\gamma,\mathcal{O}}\right)\int_{\T}\left(1-\log\left|\sin\left(\tfrac{\eta}{2} \right)\right|\right)d\eta\\
					&\lesssim \varepsilon\gamma ^{-1}\left(1+\| \mathfrak{I}_{0}\|_{q,s+\sigma_{2}}^{\gamma,\mathcal{O}}\right).
				\end{align*}
				Notice that by \eqref{symmetry kernel K1'}, \eqref{symmetry kernel K2}, the kernel $\mathbb{K}_{\varepsilon r}$ satisfies the following symmetry property 
				\begin{equation}\label{symmertry mathbbK}
					\mathbb{K}_{\varepsilon r}(-\varphi,-\theta,-\eta)=\mathbb{K}_{\varepsilon r}(\varphi,\theta,\eta)\in\mathbb{R},
				\end{equation}
				which implies in view of Lemma \ref{lemma symmetry and reversibility} that $\mathfrak{R}_{\varepsilon r}$ is a real and reversibility preserving Toeplitz in time integral operator. It remains to estimate the quantity $\displaystyle\max_{k\in\{0,1\}}\|\Delta_{12}\partial_{\theta}^{k}\mathfrak{R}_{\varepsilon r}\|_{\textnormal{\tiny{O-d}},q,\overline{s}_{h}+\mathtt{p}}^{\gamma,\mathcal{O}}$. This is will be implemented as before and we shall here sketch the main ideas. First we observe that for $k\in\{0,1\}$ the kernel of $\Delta_{12}\partial_{\theta}^{k}\mathfrak{R}_{\varepsilon r}$ is given by 
				$$
				\Delta_{12}\partial_\theta^{k}\mathbb{K}_{\varepsilon r}=-\Delta_{12}\partial_\theta^{k}\widehat{\mathbb{K}}_{\varepsilon r,1}-\Delta_{12}\partial_\theta^{k}\mathbb{K}_{\varepsilon r,2}.
				$$
				To estimate  $\Delta_{12}\partial_\theta^{k}\mathbb{K}_{\varepsilon r,2}$ we shall use \eqref{expression kernel K2} leading to
				\begin{align}\label{GL-0}
					\Delta_{12}\mathbb{K}_{\varepsilon r,2}(\varphi,\theta,\eta)=\mathscr{K}(\theta-\eta)\Delta_{12}\mathscr{K}_{\widehat\beta,2}^{1}(\varphi,\theta,\eta)+\Delta_{12}\mathscr{K}_{\widehat\beta,2}^{2}(\varphi,\theta,\eta)
				\end{align}
			and
				\begin{align}\label{GL-1}
					\nonumber\Delta_{12}\partial_\theta\mathbb{K}_{\varepsilon r,2}(\varphi,\theta,\eta)&=\mathscr{K}(\theta-\eta)\Delta_{12}\partial_\theta\mathscr{K}_{\widehat\beta,2}^{1}(\varphi,\theta,\eta)+\mathscr{K}^\prime(\theta-\eta)\Delta_{12}\mathscr{K}_{\widehat\beta,2}^{1}(\varphi,\theta,\eta)\\
					&\qquad\qquad +\Delta_{12}\partial_\theta\mathscr{K}_{\widehat\beta,2}^{2}(\varphi,\theta,\eta).
				\end{align}
				Observe from  \eqref{TU-1} that the preceding  kernels can be expressed with respect to $\widehat\beta$. Then proceeding in a similar way to \eqref{estimate vr1PM} we obtain
				\begin{equation}\label{PMZ}
					\forall \, i\in\{1,2\},\quad \max_{k\in\{0,1\}}\sup_{\eta\in\mathbb{T}}\| d_{\widehat\beta}\partial_\theta^k\mathscr{K}_{\widehat\beta,2}^{i}[\rho](\ast,\cdot,\centerdot,\eta+\centerdot)\|_{q,s}^{\gamma,\mathcal{O}}\lesssim \| \rho\|_{q,s+2}^{\gamma,\mathcal{O}}+\| \rho\|_{q,s_0+1}^{\gamma,\mathcal{O}}\| \widehat\beta\|_{q,s+2}^{\gamma,\mathcal{O}}.
				\end{equation}
				Applying Taylor Formula yields for all $i\in\{1,2\}$ and for all $k\in\{0,1\},$
				\begin{align*}
					\Delta_{12}\partial_\theta^k\mathscr{K}_{\widehat\beta,2}^{i}(\varphi,\theta,\theta+\eta)=\bigintssss_0^1d_{\widehat\beta}\partial_\theta^k\mathscr{K}_{(1-\tau)\widehat\beta_2+\tau\widehat\beta_1,2}^{i}[\widehat\beta_1-\widehat\beta_2](\varphi,\theta,\theta+\eta)\, d\tau.				\end{align*}
It follows from \eqref{PMZ} that for all $i\in\{1,2\}$ and for all $k\in\{0,1\}$
				\begin{align*}
					\big\|\Delta_{12}\partial_\theta^k\mathscr{K}_{\widehat\beta,2}^{i}(\ast,\cdot,\centerdot,\eta+\centerdot)\big\|_{q,s}^{\gamma,\mathcal{O}}\lesssim \|\widehat\beta_2-\widehat\beta_1\|_{q,s+2}^{\gamma,\mathcal{O}}+\|\widehat\beta_2-\widehat\beta_1\|_{q,s_0+1}^{\gamma,\mathcal{O}}\bigintssss_0^1\|(1-\tau)\widehat\beta_2+\tau\widehat\beta_1\|_{q,s+2}^{\gamma,\mathcal{O}}\, d\tau.
				\end{align*}
				Therefore,  by our previous choice of $\sigma_{2}$, we obtain in view of \eqref{estimate beta and r}, \eqref{difference beta} (applied with $\mathtt{p}$ replaced by $\mathtt{p}+s_{0}$) and the smallness condition \eqref{small-2},
				\begin{align*}
					\forall i\in\{1,2\},\quad \max_{k\in\{0,1\}}\big\|\Delta_{12}\partial_\theta^k\mathscr{K}_{\widehat\beta,2}^{i}(\ast,\cdot,\centerdot,\eta+\centerdot)\big\|_{q,\overline{s}_{h}+\mathtt{p}+s_{0}}^{\gamma,\mathcal{O}}&\lesssim\varepsilon\gamma^{-1}\|\Delta_{12}i\|_{q,\overline{s}_{h}+\mathtt{p}+\sigma_{2}}^{\gamma,\mathcal{O}}\Big(1+\varepsilon\gamma ^{-1}\big(1+\| \mathfrak{I}_{0}\|_{q,\overline{s}_{h}+\mathtt{p}+\sigma_{2}}^{\gamma ,\mathcal{O}}\big)\Big)\\
					&\lesssim\varepsilon\gamma^{-1}\|\Delta_{12}i\|_{q,\overline{s}_{h}+\mathtt{p}+\sigma_{2}}^{\gamma,\mathcal{O}}.
				\end{align*}
				Inserting this estimate into \eqref{GL-0} and \eqref{GL-1} yields
				\begin{align}\label{GL-2}
					\max_{k\in\{0,1\}}\sup_{\eta\in\T}\big\|\Delta_{12}\partial_\theta^{k}\mathbb{K}_{\varepsilon r,2}(\ast,\cdot,\centerdot,\eta+\centerdot)\big\|_{q,\overline{s}_{h}+\mathtt{p}+s_{0}}^{\gamma,\mathcal{O}}
					&\lesssim \varepsilon\gamma^{-1}\|\Delta_{12}i\|_{q,\overline{s}_{h}+\mathtt{p}+\sigma_{2}}^{\gamma,\mathcal{O}}.
				\end{align}
				Using similar techniques based on Taylor Formula, one can estimate $\Delta_{12}\partial_\theta\widehat{\mathbb{K}}_{\varepsilon r,1}$. We use in particular the identity  \eqref{Sam-d1} combined with \eqref{estimate vr1PM},  \eqref{estimate beta and r}, \eqref{difference beta} and the smallness condition \eqref{small-2} allowing to get
				\begin{equation}\label{GL-3}
					\max_{k\in\{0,1\}}\sup_{\eta\in\T}\big\|\Delta_{12}\partial_{\theta}^{k}\widehat{\mathbb{K}}_{\varepsilon r,1}(\ast,\cdot,\centerdot,\eta+\centerdot)\big\|_{q,\overline{s}_{h}+\mathtt{p}+s_{0}}^{\gamma,\mathcal{O}}
					\lesssim \varepsilon\gamma^{-1}\|\Delta_{12}i\|_{q,\overline{s}_{h}+\mathtt{p}+\sigma_{2}}^{\gamma,\mathcal{O}}.
				\end{equation}
				Putting together \eqref{GL-2} and \eqref{GL-3}  gives
				\begin{equation*}
					\max_{k\in\{0,1\}}\sup_{\eta\in\T}\big\|\Delta_{12}\partial_{\theta}^{k}{\mathbb{K}}_{\varepsilon r}(\ast,\cdot,\centerdot,\eta+\centerdot)\big\|_{q,\overline{s}_{h}+\mathtt{p}+s_{0}}^{\gamma,\mathcal{O}}
					\lesssim \varepsilon\gamma^{-1}\|\Delta_{12}i\|_{q,\overline{s}_{h}+\mathtt{p}+\sigma_{2}}^{\gamma,\mathcal{O}}.
				\end{equation*}
				Comibining this estimate with  Lemma \ref{lemma symmetry and reversibility}  yields
				\begin{align*}
					\max_{k\in\{0,1\}}\|\Delta_{12}\partial_{\theta}^{k}\mathfrak{R}_{\varepsilon r}\|_{\textnormal{\tiny{O-d}},q,\overline{s}_{h}+\mathtt{p}}^{\gamma,\mathcal{O}}&\lesssim\max_{k\in\{0,1\}}\int_{\mathbb{T}}\|\Delta_{12}\partial_{\theta}^{k}\mathbb{K}_{\varepsilon r}(\ast,\cdot,\centerdot,\eta+\centerdot)\|_{q,\overline{s}_{h}+\mathtt{p}+s_{0}}^{\gamma,\mathcal{O}}d\eta\\
					&\lesssim\varepsilon\gamma^{-1}\|\Delta_{12}i\|_{q,\overline{s}_{h}+\mathtt{p}+\sigma_{2}}^{\gamma,\mathcal{O}}.
				\end{align*}
				This completes  the proof of the Proposition \ref{Action on the non local part}.
			\end{proof}
			\subsection{Diagonalization up to small errors}
			The main goal of this section is to diagonalize, up to small errors, the operator $\widehat{\mathcal{L}}_{\omega}$ discussed in Proposition \ref{lemma setting for Lomega} 
			and given by
			$$\widehat{\mathcal{L}}_{\omega}=\Pi_{\mathbb{S}_{0}}^{\perp}(\mathcal{L}_{\varepsilon r}-\varepsilon\partial_{\theta}\mathcal{R})\Pi_{\mathbb{S}_{0}}^{\perp}.$$
			This will be performed in two main steps. First, we shall explore the effect of the frequency localization in the normal direction  on  the transport  reduction discussed in Section \ref{Reduction of order 1-1}. We essentially get the same structure up to a small perturbation of finite-dimensional rank. Then, in the second step  we shall implement a KAM reducibility scheme in order  to reduce the remainder to a diagonal one modulo  small fast decaying  operators. This will be performed through the use of  a suitable strong  topology on continuous  operators given by  \eqref{Top-NormX}. With this topology one has  tame estimates and   the Toeplitz structure of the remainder  is very important in this part. The reduction will be conducted by assuming non resonance conditions stemming from the second order Melnikov conditions needed in the resolution of adequate homological equations during the scheme.
			\subsubsection{Projection in the normal directions}\label{Diagon-normal}
			In this section, we study the effects of the reduction of  the transport part  when the linearized operator is  localized in the normal directions. Notice that the change of coordinates does not stabilize the normal subspace and as we shall see the defect of the commutation can be modeled by projectors of finite ranks. Let us define 
			$$\mathscr{B}_{\perp}:=\Pi_{\mathbb{S}_{0}}^{\perp}\mathscr{B}\Pi_{\mathbb{S}_{0}}^{\perp},$$
			where the transformation $\mathscr{B}$ is introduced in \eqref{definition symplectic change of variables} and constructed in Proposition \ref{reduction of the transport part}. Recall that the projection $\Pi_{\mathbb{S}_0}^{\perp}$ and the Soboloev space $H_{\perp}^{s}$ were respectively defined in \eqref{def proj} and \eqref{decoacca}. We also recall the following notations
			$$\mathbf{e}_{l,j}(\varphi,\theta)=e^{\ii (l\cdot\varphi+j\theta)}\quad\textnormal{and}\quad e_m(\theta)=e^{\ii m\theta}.$$
			The first main result of this section reads as follows. 
			\begin{lem}\label{Lema-decomp}
			Let $\mathscr{B}$ the transformation constructed in Proposition $\ref{reduction of the transport part},$ then under the condition \eqref{small-2} and \eqref{param-trans}, the following assertions hold.
			\begin{enumerate}[label=(\roman*)]
						\item For all $s\in[s_{0},S]$, the operator $\mathscr{B}_{\perp}:W^{q,\infty,\gamma}\big(\mathcal{O}, H_{\perp}^{s}\big)\to W^{q,\infty,\gamma}\big(\mathcal{O}, H_{\perp}^{s}\big)$ is continuous and invertible, with
						\begin{equation}\label{estimate for mathscrBperp and its inverse}
							\|\mathscr{B}_{\perp}^{\pm 1}\rho\|_{q,s}^{\gamma ,\mathcal{O}}\lesssim\|\rho\|_{q,s}^{\gamma ,\mathcal{O}}+\varepsilon\gamma ^{-1}\| \mathfrak{I}_{0}\|_{q,s+\sigma_{3}}^{\gamma ,\mathcal{O}}\|\rho\|_{q,s_{0}}^{\gamma ,\mathcal{O}}.
						\end{equation}
						In addition, we have the representations
						$$
						\mathscr{B}_{\perp}\rho=\mathscr{B}\rho-\sum_{m\in\mathbb{S}_0}\big\langle\rho,  \big(\mathcal{B}^{-1}-\textnormal{Id}\big)e_m  \big\rangle_{L^2_\theta(\T)}e_m
						$$
						and
						\begin{align*}
							\mathscr{B}_{\perp}^{-1}\rho=\mathscr{B}^{-1}\rho-\sum_{m\in\mathbb{S}_0}\big\langle \rho, \big(\mathcal{B}-\textnormal{Id}\big)g_m  \big\rangle_{L^2_\theta(\T)}\mathscr{B}^{-1}e_m,
						\end{align*}
						where 
						$$
						\quad \mathbf{A}(\varphi)=\Big(\big\langle e_m,\mathcal{B}e_k\big\rangle_{L^2_\theta(\T)}\Big)_{m\in\mathbb{S}_0\atop k\in\mathbb{S}_0},\quad   \mathbf{A}^{-1}(\varphi)=\Big(\alpha_{k,m}\Big)_{m\in\mathbb{S}_0\atop k\in\mathbb{S}_0},\quad g_m(\varphi,\theta)=\sum_{
							k\in\mathbb{S}_0}\alpha_{k,m}(\varphi)e_k(\theta),
						$$
						with the estimate
						\begin{align*}
							\sup_{k,m\in\mathbb{S}_{0}}\|\alpha_{k,m}-\delta_{km}\|_{q,s}^{\gamma,\mathcal{O}}\lesssim \varepsilon\gamma ^{-1}\left(1+\| \mathfrak{I}_{0}\|_{q,s+\sigma_{1}+1}^{\gamma ,\mathcal{O}}\right) .
						\end{align*}
						\item Given two tori $i_{1}$ and $i_{2}$ satisfying the smallness condition \eqref{small-2}, one has
						\begin{equation}\label{diff gm}
							\max_{m\in\mathbb{S}_{0}}\|\Delta_{12}g_{m}\|_{q,\overline{s}_{h}+\mathtt{p}}^{\gamma,\mathcal{O}}\lesssim\varepsilon\gamma^{-1}\|\Delta_{12}i\|_{q,\overline{s}_{h}+\mathtt{p}+\sigma_{1}+1}^{\gamma,\mathcal{O}}.
						\end{equation}
					\end{enumerate}

			\end{lem}
			
			\begin{proof}

				{\bf{(i)}} 
				The first estimate concerning $\mathscr{B}_{\perp}$ follows easily from the continuity of the orthogonal projector $\Pi_{\mathbb{S}_{0}}^{\perp}$ on Sobolev spaces  $H^{s}_{\perp},$ combined with  \eqref{estimate on the first reduction operator and its inverse}. 
				 For the representation of $\mathscr{B}_{\perp},$ take $\rho\in W^{q,\infty,\gamma}(\mathcal{O}, H_{\perp}^{s})$ and set
				$$
				\mathscr{B}_{\perp}\rho=\Pi_{\mathbb{S}_{0}}^{\perp}\mathscr{B}\Pi_{\mathbb{S}_{0}}^{\perp}\rho=\Pi_{\mathbb{S}_{0}}^{\perp}\mathscr{B}\rho:=g.
				$$ 
				Next, we write the following splitting 
				\begin{equation}\label{Split29}
					\mathscr{B}\rho=g+h\quad\hbox{with}\quad \Pi_{\mathbb{S}_{0}}h=h.
				\end{equation}
			Notice that the projector $\Pi_{\mathbb{S}_0}$ is defined by 
			$$\Pi_{\mathbb{S}_0}\rho=\sum_{j\in\mathbb{S}_0}\rho_j e_{j}=\Pi_{\overline{\mathbb{S}}}\rho+\langle\rho\rangle_{\theta},$$
			where $\Pi_{\overline{\mathbb{S}}}$ is defined in \eqref{def proj} and $\langle\cdot\rangle_{\theta}$ denotes the avarage in the variable $\theta.$
				Therefore 
				$$\displaystyle{h(\varphi,\theta)=\sum_{m\in\mathbb{S}_{0}}h_m(\varphi) e_m(\theta)},$$
				supplemented with the orthogonal conditions
				$$
				\forall \, k\in \mathbb{S}_{0},\quad \langle  \mathscr{B}\rho-h, e_k\rangle_{L^2_\theta(\T)}=0.
				$$
				This implies 
				$$
				h(\varphi,\theta)=\sum_{m\in\mathbb{S}_{0}}\langle \mathscr{B}\rho, e_m\rangle_{L^2_\theta(\T)}e_m(\theta).
				$$
				Using Lemma \ref{algeb1}-(iii) leads to
				$$
				h(\varphi,\theta)=\sum_{m\in\mathbb{S}_{0}}\langle\rho,  \mathcal{B}^{-1}e_m\rangle_{L^2_\theta(\T)}e_m(\theta).
				$$
				Inserting this identity into \eqref{Split29} yields
				$$
				\mathscr{B}_{\perp}\rho=g=\mathscr{B}\rho-\sum_{m\in\mathbb{S}_{0}}\big\langle\rho,  \mathcal{B}^{-1}e_m\big\rangle_{L^2_\theta(\T)}e_m.
				$$
				Since $\forall m\in \mathbb{S}_{0},\,\big\langle\rho, e_m\big\rangle_{L^2_\theta(\T)}=0,$ then
				$$
				\mathscr{B}_{\perp}\rho=g=\mathscr{B}\rho-\sum_{m\in\mathbb{S}_{0}}\big\langle\rho,  \big(\mathcal{B}^{-1}-\textnormal{Id}\big)e_m\big\rangle_{L^2_\theta(\T)}e_m.
				$$
				This ensures the desired representation of $  \mathscr{B}_{\perp}$.\\
				Next, we intend to establish similar  representation for $\mathscr{B}_{\perp}^{-1}$. Let $g\in W^{q,\infty,\gamma}(\mathcal{O}, H_{\perp}^{s})$ and we need to solve the equation
				$$
				f\in W^{q,\infty,\gamma}(\mathcal{O}, H_{\perp}^{s}),\quad  \mathscr{B}_{\perp}f=\Pi_{\mathbb{S}_{0}}^{\perp}\mathscr{B}f=g.
				$$ 
				This is equivalent to
				$$
				\mathscr{B}f=g+h, \quad \hbox{with}\quad \Pi_{\mathbb{S}_{0}}h=h\quad\hbox{and}\quad  \Pi_{\mathbb{S}_{0}}f=0.
				$$
				Then we get 
				\begin{align}\label{f-Id}
					f= \mathscr{B}^{-1}\big(g+h\big), \quad \hbox{with}\quad \Pi_{\mathbb{S}_{0}}h=h\quad\textnormal{and}\quad\Pi_{\mathbb{S}_{0}}f=0.
				\end{align}
				The condition  $\Pi_{\mathbb{S}_{0}}f=0$ is equivalent to,
				$$
				\forall k\in \mathbb{S}_{0},\quad \big\langle\mathscr{B}^{-1}\big(g+h\big), e_k   \big\rangle_{L^2_\theta(\T)}=0.
				$$
				Therefore  using Lemma \ref{algeb1}-(iii) the latter equation reads
				$$
				\forall k\in \mathbb{S}_{0},\quad\big\langle g+h, \widehat{e}_k   \big\rangle_{L^2_\theta(\T)}=0\quad \hbox{with}\quad \widehat{e}_k(\varphi,\theta):=\mathcal{B}e_k(\varphi,\theta)=e^{\ii k(\theta+\beta(\varphi,\theta))},
				$$
				which will fix $h$. Indeed, by expanding  $\displaystyle h(\varphi,\theta)=\sum_{m\in \mathbb{S}_{0}} a_m(\varphi)e_m(\theta)$, we can transform the preceding system  into
				\begin{align}\label{matrix-inv}
					\forall k\in \mathbb{S}_{0},\quad \sum_{m\in \mathbb{S}_{0}} a_m(\varphi)\big\langle e_m, \widehat{e}_k \big\rangle_{L^2_\theta(\T)}=-\big\langle g, \widehat{e}_k   \big\rangle_{L^2_\theta(\T)}.
				\end{align}
				Define the matrix  
				\begin{equation}\label{def matrix A}
					\mathbf{A}(\varphi)=\big(c_{m,k}(\varphi)\big)_{(m,k)\in\mathbb{S}_{0}^2},\quad c_{m,k}(\varphi)=\langle e_m, \widehat{e}_k\rangle_{L^2_\theta(\T)}=\int_{\T}e^{\ii((m-k)\theta-k\beta(\varphi,\theta))}d\theta.
				\end{equation}
			Notice that according to \eqref{symmetry for beta} and the change of variables $\theta\mapsto-\theta$, one obtains
			\begin{equation}\label{symmetry cmk}
				\forall (m,k)\in\mathbb{S}_0^2,\quad \forall\varphi\in\mathbb{T}^d,\quad c_{m,k}(-\varphi)=c_{-m,-k}(\varphi)=\overline{c_{m,k}(\varphi,\theta)}.
		\end{equation}
				One can check by slight adaptation of   the composition law in Lemma \ref{Lem-lawprod} and using the smallness condition \eqref{small-2} and \eqref{estimate beta and r}
				\begin{align}\label{refm}
					\nonumber\|c_{m,m}-1\|_{q,s}^{\gamma,\mathcal{O}}&\leqslant  \bigintssss_{\T}\|e^{-\ii m\beta(\cdot,\theta))}-1\|_{q,H_\varphi^s}^{\gamma,\mathcal{O}}d\theta\\
					\nonumber &\lesssim 
					\|\beta\|_{q,s}^{\gamma,\mathcal{O}}\\
					&\lesssim \varepsilon\gamma ^{-1}\left(1+\| \mathfrak{I}_{0}\|_{q,s+\sigma_{1}}^{\gamma ,\mathcal{O}}\right).
				\end{align}
				For $k\neq m\in  \mathbb{S}_{0}$ we use  integration by parts, 
				$$
				c_{m,k}(\varphi)=\tfrac{k}{\ii (m-k)}\int_{\T}e^{\ii((m-k)\theta-k\beta(\varphi,\theta))}\partial_\theta\beta(\varphi,\theta)d\theta.
				$$
				Then using law products and composition laws in  Lemma \ref{Lem-lawprod} combined with \eqref{estimate beta and r} yield
				\begin{align*}
					\sup_{(m,k)\in \mathbb{S}_{0}^{2}\atop m\neq k} \|c_{m,k}\|_{q,s}^{\gamma,\mathcal{O}}&\lesssim \|\beta\|_{q,s+1}^{\gamma,\mathcal{O}}\\
					&\lesssim  \varepsilon\gamma ^{-1}\left(1+\| \mathfrak{I}_{0}\|_{q,s+\sigma_{1}+1}^{\gamma ,\mathcal{O}}\right) .
				\end{align*}
				Finally, we get that
				\begin{equation}\label{vdo1}
					\mathbf{A}(\varphi)=\textnormal{Id}+\mathbf{R}(\varphi)\quad\hbox{with} \quad \|\mathbf{R}\|_{q,s}^{\gamma,\mathcal{O}}\lesssim \|\beta\|_{q,s+1}^{\gamma,\mathcal{O}}.
				\end{equation}
				Hence under the smallness condition $\|\beta\|_{q,s_0}^{\gamma,\mathcal{O}}\ll1$ following from \eqref{small-2}, combined with the law products in  Lemma \ref{Lem-lawprod}   we get that $\mathbf{A}$ is invertible  with
				\begin{align}\label{H-D-P0}
					\nonumber\|\mathbf{A}^{-1}-\textnormal{Id}\|_{q,s}^{\gamma,\mathcal{O}}&\lesssim \|\beta\|_{q,s+1}^{\gamma,\mathcal{O}}
					\\
					&\lesssim  \varepsilon\gamma ^{-1}\left(1+\| \mathfrak{I}_{0}\|_{q,s+\sigma_{1}+1}^{\gamma ,\mathcal{O}}\right) .
				\end{align}
				Therefore the system \eqref{matrix-inv} is invertible and one gets a unique solution given by
				\begin{align}\label{matrix-inv0}
					a_m(\varphi)=-\sum_{k\in\mathbb{S}_0}\alpha_{m,k}(\varphi)\big\langle g, \widehat{e}_k   \big\rangle_{L^2_\theta(\T)}\quad\hbox{with}\quad \mathbf{A}^{-1}(\varphi):=\Big(\alpha_{m,k}(\varphi)\Big)_{(m,k)\in\mathbb{S}_0^2}.
					\end{align}
	We claim that the coefficients of $\mathbf{A}^{-1}$  admit the same symmetry conditions as  \eqref{symmetry cmk}, that is
				\begin{equation}\label{symmetry alphamk}
					\forall(m,k)\in\mathbb{S}_0^2,\quad\forall\varphi\in\mathbb{T}^d,\quad \alpha_{m,k}(-\varphi)=\alpha_{-m,-k}(\varphi)=\overline{\alpha_{m,k}(\varphi)}.
				\end{equation}
				This can be done through the series expansion $\displaystyle{A^{-1}=\sum_{n\in\mathbb{N}}(-1)^n(A-\hbox{Id})^n}$ together with the fact that the entries of the monomials $(A-\hbox{Id})^n$ satisfy in turn $ \eqref{symmetry cmk}$.
				Next, using the law products yields
				\begin{equation}\label{H-D-P}
					\sup_{m\in\mathbb{S}_{0}}\|a_{m}\|_{q,s}^{\gamma,\mathcal{O}}\lesssim \sup_{k\in\mathbb{S}_{0}} \Big(\|\mathbf{A}^{-1}\|_{q,s}^{\gamma,\mathcal{O}} \|\big\langle g, \widehat{e}_k   \big\rangle_{L^2_\theta(\T)}\|_{q,H^{s_0}_\varphi}^{\gamma,\mathcal{O}}+ \|\mathbf{A}^{-1}\|_{q,s_0}^{\gamma,\mathcal{O}} \|\big\langle g, \widehat{e}_k   \big\rangle_{L^2_\theta(\T)}\|_{q,H^{s}_\varphi}^{\gamma,\mathcal{O}}\Big).
				\end{equation}
				Notice that one gets from \eqref{H-D-P0}
				\begin{align*}
					\sup_{(m,k)\in\mathbb{S}_{0}^{2}}\|\alpha_{k,m}-\delta_{km}\|_{q,s}^{\gamma,\mathcal{O}}\lesssim  \varepsilon\gamma ^{-1}\left(1+\| \mathfrak{I}_{0}\|_{q,s+\sigma_{1}+1}^{\gamma ,\mathcal{O}}\right) ,
				\end{align*}
				where $\delta_{km}$ denotes the Kronecker symbol. 
				Let us now move to the estimate of the partial scalar product  containing  $g$ in \eqref{H-D-P}. Using the law products in  Lemma \ref{Lem-lawprod} with  Cauchy-Schwarz inequality gives
				\begin{align*}
					\|\big\langle g, \widehat{e}_k   \big\rangle_{L^2_\theta(\T)}\|_{q,H^{s}_\varphi}^{\gamma,\mathcal{O}}&\lesssim \bigintssss_{\T}\Big(\|g(\cdot,\theta)\|_{q,H^{s}_\varphi}^{\gamma,\mathcal{O}}\|e^{i\beta(\cdot,\theta)}\|_{q,H^{s_0}_\varphi}^{\gamma,\mathcal{O}}+\|g(\cdot,\theta)\|_{q,H^{s_0}_\varphi}^{\gamma,\mathcal{O}}\|e^{i\beta(\cdot,\theta)}\|_{q,H^{s}_\varphi}^{\gamma,\mathcal{O}}\Big)d\theta\\
					&\lesssim \|g\|_{q,L^2_\theta H^{s}_\varphi}^{\gamma,\mathcal{O}}\|e^{i\beta}\|_{q,L^2_\theta H^{s_0}_\varphi}^{\gamma,\mathcal{O}}+\|g\|_{q,L^2_\theta H^{s_0}_\varphi}^{\gamma,\mathcal{O}}\|e^{i\beta}\|_{q,L^2_\theta H^{s}_\varphi}^{\gamma,\mathcal{O}}
					\\
					&\lesssim \|g\|_{q,s}^{\gamma,\mathcal{O}}\|e^{i\beta}\|_{q,s_0}^{\gamma,\mathcal{O}}+\|g\|_{q,s_0}^{\gamma,\mathcal{O}}\|e^{i\beta}\|_{q,s}^{\gamma,\mathcal{O}}.
				\end{align*}
				Then applying  the composition law as in \eqref{refm} combined with  with \eqref{estimate beta and r} and the smallness condition \eqref{small-2} gives
				\begin{align*}
					\|\big\langle g, \widehat{e}_k   \big\rangle_{L^2_\theta(\T)}\|_{q,H^{s}_\varphi}^{\gamma,\mathcal{O}}&\lesssim \|g\|_{q,s}^{\gamma,\mathcal{O}}+\|g\|_{q,s_0}^{\gamma,\mathcal{O}}\|\beta\|_{q,s}^{\gamma,\mathcal{O}}\\
					&\lesssim \|g\|_{q,s}^{\gamma,\mathcal{O}}+ \varepsilon\gamma ^{-1}\left(1+\| \mathfrak{I}_{0}\|_{q,s+\sigma_{1}}^{\gamma ,\mathcal{O}}\right)\|g\|_{q,s_0}^{\gamma,\mathcal{O}}.
				\end{align*}
				Plugging this estimate into \eqref{H-D-P} and using \eqref{estimate beta and r},  \eqref {H-D-P0} combined with the smallness condition \eqref{small-2} and Sobolev embeddings implies 
				\begin{align*}
					\sup_{m\in\mathbb{S}_{0}}\|a_{m}\|_{q,s}^{\gamma,\mathcal{O}}&\lesssim \|g\|_{q,s}^{\gamma,\mathcal{O}}+ \varepsilon\gamma ^{-1}\left(1+\| \mathfrak{I}_{0}\|_{q,s+\sigma_{1}+1}^{\gamma ,\mathcal{O}}\right)\|g\|_{q,s_0}^{\gamma,\mathcal{O}}\\
					&\lesssim \|g\|_{q,s}^{\gamma,\mathcal{O}}+ \varepsilon\gamma ^{-1}\| \mathfrak{I}_{0}\|_{q,s+\sigma_{1}+1}^{\gamma ,\mathcal{O}}\|g\|_{q,s_0}^{\gamma,\mathcal{O}}.
				\end{align*}
				Therefore we obtain
				\begin{align*}
					\|h\|_{q,s}^{\gamma,\mathcal{O}}&\lesssim \sum_{m\in \mathbb{S}_{0}} \|a_m\|_{q,H^s_\varphi}^{\gamma,\mathcal{O}}\\
					&\lesssim \|g\|_{q,s}^{\gamma,\mathcal{O}}+ \varepsilon\gamma ^{-1}\| \mathfrak{I}_{0}\|_{q,s+\sigma_{1}+1}^{\gamma ,\mathcal{O}}\|g\|_{q,s_0}^{\gamma,\mathcal{O}}.
				\end{align*}
				Coming  back to \eqref{f-Id} and using \eqref{estimate on the first reduction operator and its inverse}, we get
				\begin{align*}
					\|f\|_{q,s}^{\gamma ,\mathcal{O}}&\lesssim\|g+h\|_{q,s}^{\gamma ,\mathcal{O}}+\varepsilon\gamma ^{-1}\| \mathfrak{I}_{0}\|_{q,s+\sigma_{1}+1}^{\gamma ,\mathcal{O}}\|g+h\|_{q,s_{0}}^{\gamma ,\mathcal{O}}\\
					&\lesssim \|g\|_{q,s}^{\gamma,\mathcal{O}}+ \varepsilon\gamma ^{-1}\| \mathfrak{I}_{0}\|_{q,s+\sigma_{1}+1}^{\gamma ,\mathcal{O}}\|g\|_{q,s_0}^{\gamma,\mathcal{O}}.
				\end{align*}
				It follows that
				\begin{align*}
					\|\mathscr{B}_{\perp}^{-1}g\|_{q,s}^{\gamma ,\mathcal{O}}&\lesssim \|g\|_{q,s}^{\gamma,\mathcal{O}}+ \varepsilon\gamma ^{-1}\| \mathfrak{I}_{0}\|_{q,s+\sigma_{1}+1}^{\gamma ,\mathcal{O}}\|g\|_{q,s_0}^{\gamma,\mathcal{O}}.
				\end{align*}
				In addition from \eqref{matrix-inv0} and \eqref{f-Id} we deduce the formula
				\begin{align}\label{Mer-11}
					\nonumber \mathscr{B}_{\perp}^{-1}g(\varphi,\theta)&=\mathscr{B}^{-1}g(\varphi,\theta)-\sum_{m\in\mathbb{S}_0\atop
						k\in\mathbb{S}_0}\alpha_{m,k}(\varphi)\big\langle g, \mathcal{B}{e}_k   \big\rangle_{L^2_\theta(\T)}\mathscr{B}^{-1}e_m(\theta)\\
					&=\mathscr{B}^{-1}g(\varphi,\theta)-\sum_{m\in\mathbb{S}_0}\big\langle g,\mathcal{B}g_m \big\rangle_{L^2_\theta(\T)}\big(\mathscr{B}^{-1}e_m\big)(\varphi,\theta),
				\end{align}
				with
				\begin{equation}\label{link gm alphamk}
				g_m(\varphi,\theta):=\sum_{k\in\mathbb{S}_0}\alpha_{m,k}(\varphi){e}_k (\theta). 
			\end{equation}
	From  \eqref{symmetry alphamk} and the symmetry of $\mathbb{S}_0$, we infer
		\begin{equation}\label{sym gm}
			\forall m\in\mathbb{S}_0,\quad\forall(\varphi,\theta)\in\mathbb{T}^{d+1},\quad g_m(-\varphi,-\theta)=g_{-m}(\varphi,\theta)=\overline{g_{m}(\varphi,\theta)}.
	\end{equation}
				Since $\Pi_{\mathbb{S}_{0}}^{\perp}g=g $ and  $\Pi_{\mathbb{S}_{0}}^{\perp}g_m=0 $ then $\big\langle g,g_m \big\rangle_{L^2_\theta(\T)}=0$ and therefore
				$$
				\big\langle g,\mathcal{B}g_m \big\rangle_{L^2_\theta(\T)}=\big\langle g,\big(\mathcal{B}-\textnormal{Id}\big)g_m \big\rangle_{L^2_\theta(\T)}.
				$$
				Plugging this identity into \eqref{Mer-11} yields
				\begin{align*}
					\mathscr{B}_{\perp}^{-1}g&=\mathscr{B}^{-1}g-\sum_{m\in\mathbb{S}_0}\big\langle g,\big(\mathcal{B}-\textnormal{Id}\big)g_m \big\rangle_{L^2_\theta(\T)}\mathscr{B}^{-1}e_m.
					\end{align*}
				\textbf{(ii)} Coming back to the definition of $c_{m,k}$ in \eqref{def matrix A}, one can write
				$$\forall(m,k)\in\mathbb{S}_{0}^{2},\quad\Delta_{12}c_{m,k}=\langle e_{m},(\Delta_{12}\mathcal{B})e_{k}\rangle_{L^{2}_{\theta}(\mathbb{T})}.$$
				Hence, using Taylor Formula and \eqref{difference beta}, we have
				$$\max_{(m,k)\in\mathbb{S}_{0}^{2}}\|\Delta_{12}c_{m,k}\|_{q,\overline{s}_{h}+\mathtt{p}}^{\gamma,\mathcal{O}}\lesssim\varepsilon\gamma^{-1}\|\Delta_{12}i\|_{q,\overline{s}_{h}+\mathtt{p}+\sigma_{1}}^{\gamma,\mathcal{O}}.$$
				From \eqref{link gm alphamk}, one has
				$$\Delta_{12}g_{m}=\sum_{k\in\mathbb{S}_{0}}\Delta_{12}\alpha_{m,k}\,e_{k}.$$
				Thus
				\begin{equation}\label{control diff gm by diff alphamk}
					\max_{m\in\mathbb{S}_{0}}\|\Delta_{12}g_{m}\|_{q,\overline{s}_{h}+\mathtt{p}}^{\gamma,\mathcal{O}}\lesssim\max_{(m,k)\in\mathbb{S}_{0}^{2}}\|\Delta_{12}\alpha_{m,k}\|_{q,\overline{s}_{h}+\mathtt{p}}^{\gamma,\mathcal{O}}.
				\end{equation}
				Using Neumann series, we can write 
				$$\mathbf{A}^{-1}(\varphi)=\textnormal{Id}+\sum_{n=1}^{\infty}(-1)^n\mathbf{R}^n(\varphi).$$
				Therefore, the law products in Lemma \ref{Lem-lawprod} combined with \eqref{vdo1} and the smallness condition \eqref{small-2} lead to
				\begin{align*}
					\|\Delta_{12}\mathbf{A}^{-1}\|_{q,\overline{s}_{h}+\mathtt{p}}^{\gamma,\mathcal{O}}&\lesssim\sum_{n=1}^{\infty}\|\Delta_{12}\mathbf{R}^{n}\|_{q,\overline{s}_{h}+\mathtt{p}}^{\gamma,\mathcal{O}}\\
					&\lesssim \|\Delta_{12}\mathbf{R}\|_{q,\overline{s}_{h}+\mathtt{p}}^{\gamma,\mathcal{O}}\\
					&\lesssim\varepsilon\gamma^{-1}\|\Delta_{12}i\|_{q,\overline{s}_{h}+\mathtt{p}+\sigma_{1}+1}^{\gamma,\mathcal{O}}.
				\end{align*}
				As a consequence,
				\begin{equation}\label{est diff alphamk}
					\max_{(m,k)\in\mathbb{S}_{0}^{2}}\|\Delta_{12}\alpha_{m,k}\|_{q,\overline{s}_{h}+\mathtt{p}}^{\gamma,\mathcal{O}}\lesssim\varepsilon\gamma^{-1}\|\Delta_{12}i\|_{q,\overline{s}_{h}+\mathtt{p}+\sigma_{1}+1}^{\gamma,\mathcal{O}}.
				\end{equation}
				Gathering \eqref{est diff alphamk} and \eqref{control diff gm by diff alphamk} finally gives
				$$\max_{m\in\mathbb{S}_{0}}\|\Delta_{12}g_{m}\|_{q,\overline{s}_{h}+\mathtt{p}}^{\gamma,\mathcal{O}}\lesssim\varepsilon\gamma^{-1}\|\Delta_{12}i\|_{q,\overline{s}_{h}+\mathtt{p}+\sigma_{1}+1}^{\gamma,\mathcal{O}}.$$
				This achieves the proof of Lemma \ref{Lema-decomp}.
				\end{proof}
 In Lemma \ref{Lema-decomp}, the parameter $\mathtt{p}$ is subject to the constraint \eqref{param-trans} and from now on, we shall fix it to the value
 \begin{equation}\label{choice ttp-1}
 	\mathtt{p}=4\tau_2q+4\tau_2.
 \end{equation}
 This particular choice is determined through some constraints  in the proof of the remainder reduction. More precisely, it appears in \eqref{choise ttp-2}.
 Next we shall establish the second main result of this section.
						\begin{prop}\label{projection in the normal directions}
				Let $(\gamma,q,d,\tau_{1},s_{0},s_{h},\overline{s}_{h},\mu_2,\mathtt{p},\sigma_2,S)$ satisfy the assumptions \eqref{initial parameter condition}, \eqref{setting tau1 and tau2}, \eqref{param}, \eqref{param-trans}, \eqref{sigma2} and \eqref{choice ttp-1}. Consider  the operator  $\widehat{\mathcal{L}}_{\omega}$ defined in Proposition $\ref{lemma setting for Lomega}.$\\
				There exists $\varepsilon_0>0$ and $\sigma_{3}=\sigma_{3}(\tau_{1},q,d,s_{0})\geqslant\sigma_{2}$ such that if 
				\begin{equation}\label{small-3}
					\varepsilon\gamma^{-1}N_0^{\mu_2}\leqslant\varepsilon_0\quad\textnormal{and}\quad \|\mathfrak{I}_0\|_{q,s_h+\sigma_3}^{\gamma,\mathcal{O}}\leqslant 1,
				\end{equation}
			then the following assertions hold true.
				\begin{enumerate}[label=(\roman*)]
					 		
						\item For any $n\in\mathbb{N}^{*},$ in the Cantor set $\mathcal{O}_{\infty,n}^{\gamma,\tau_{1}}(i_{0})$ introduced in Proposition $\ref{reduction of the transport part},$ we have 
					\begin{align*}
						\mathscr{B}_{\perp}^{-1}\widehat{\mathcal{L}}_{\omega}\mathscr{B}_{\perp}&=\big(\omega\cdot\partial_{\varphi}+c_{i_{0}}\partial_{\theta}-\partial_{\theta}\mathcal{K}_{\lambda}\ast\cdot\big)\Pi_{\mathbb{S}_0}^{\perp}+\mathscr{R}_{0}+\mathtt{E}_{n}^{1}\\
						&:=\big(\omega\cdot\partial_{\varphi}+\mathscr{D}_{0}\big)\Pi_{\mathbb{S}_0}^{\perp}+\mathscr{R}_{0}+\mathtt{E}_{n}^{1}\\
						&:=\mathscr{L}_{0}+\mathtt{E}_{n}^{1},
					\end{align*}
					where $\mathscr{D}_{0}$ is a reversible  Fourier multiplier    given by 
					$$
					\forall (l,j)\in \mathbb{Z}^{d}\times\mathbb{S}_{0}^{c},\quad  \mathscr{D}_{0}\mathbf{e}_{l,j}=\ii\mu_{j}^{0}\,\mathbf{e}_{l,j},$$
					with 
					$$\mu_{j}^{0}(\lambda,\omega,i_{0})=\Omega_{j}(\lambda)+jr^{1}(\lambda,\omega,i_{0})\quad\mbox{ and }\quad r^{1}(\lambda,\omega,i_{0})=c_{i_{0}}(\lambda,\omega)-V_0(\lambda)$$
					and such that
					\begin{equation}\label{differences mu0}
						\|r^1\|_{q}^{\gamma,\mathcal{O}}\lesssim \varepsilon \quad\textnormal{and}\quad \|\Delta_{12}r^{1}\|_{q}^{\gamma,\mathcal{O}}\lesssim \varepsilon \| \Delta_{12}i\|_{q,\overline{s}_{h}+2}^{\gamma,\mathcal{O}}.
					\end{equation}
				\item The operator $\mathtt{E}_{n}^{1}$ satisfies the following estimate
%				\begin{equation}\label{En1-s}
%					\forall s\in[s_{0},S],\quad\|\mathtt{E}_{n}^{1}\rho\|_{q,s}^{\gamma,\mathcal{O}}\lesssim \|\rho\|_{q,s+2}^{\gamma,\mathcal{O}}+\varepsilon\gamma^{-1}\|\mathfrak{I}_{0}\|_{q,s+\sigma_{3}}^{\gamma,\mathcal{O}}\|\rho\|_{q,s_{0}+2}^{\gamma,\mathcal{O}}
%				\end{equation}
%			and 
			\begin{equation}\label{En1-s0}
				\|\mathtt{E}_{n}^{1}\rho\|_{q,s_{0}}^{\gamma,\mathcal{O}}\lesssim \varepsilon N_{0}^{\mu_{2}}N_{n+1}^{-\mu_{2}}\|\rho\|_{q,s_{0}+2}^{\gamma,\mathcal{O}}.
			\end{equation}
		\item  $\mathscr{R}_{0}$   is a real and reversible Toeplitz in time operator satisfying $\mathscr{R}_{0}=\Pi_{\mathbb{S}_0}^\perp \mathscr{R}_{0}\Pi_{\mathbb{S}_0}^\perp$ with
					\begin{equation}\label{estimate mathscr R in off diagonal norm}
						\forall s\in [s_{0},S],\quad \max_{k\in\{0,1\}}\|\partial_{\theta}^{k}\mathscr{R}_{0}\|_{\textnormal{\tiny{O-d}},q,s}^{\gamma,\mathcal{O}}\lesssim\varepsilon\gamma^{-1}\left(1+\| \mathfrak{I}_{0}\|_{q,s+\sigma_{3}}^{\gamma,\mathcal{O}}\right)
					\end{equation}
					and
					\begin{equation}\label{estimate differences mathscr R in off diagonal norm}
						\|\Delta_{12}\mathscr{R}_{0}\|_{\textnormal{\tiny{O-d}},q,\overline{s}_{h}+\mathtt{p}}^{\gamma,\mathcal{O}}\lesssim\varepsilon\gamma^{-1}\| \Delta_{12}i\|_{q,\overline{s}_{h}+\mathtt{p}+\sigma_{3}}^{\gamma,\mathcal{O}}.
					\end{equation}
					\item The operator $\mathscr{L}_{0}$ satisfies
					\begin{equation}\label{Taptap1}
						\forall s\in[s_{0},S],\quad\|\mathscr{L}_{0}\rho\|_{q,s}^{\gamma ,\mathcal{O}}\lesssim\|\rho\|_{q,s+1}^{\gamma ,\mathcal{O}}+\varepsilon\gamma ^{-1}\| \mathfrak{I}_{0}\|_{q,s+\sigma_{3}}^{\gamma ,\mathcal{O}}\|\rho\|_{q,s_{0}}^{\gamma ,\mathcal{O}}.
					\end{equation}
				\end{enumerate}
			\end{prop}
			\begin{proof}
			
							\textbf{ (i)} We shall first start with finding a suitable expansion  for $\mathscr{B}_{\perp}^{-1}\widehat{\mathcal{L}}_{\omega}\mathscr{B}_{\perp}$. Using the expression of $\widehat{\mathcal{L}}_{\omega}$ given in Proposition \ref{lemma setting for Lomega} and the decomposition $\textnormal{Id}=\Pi_{\mathbb{S}_{0}}+\Pi_{\mathbb{S}_{0}}^{\perp}$  we write
				\begin{align*}
					\mathscr{B}_{\perp}^{-1}\widehat{\mathcal{L}}_{\omega}\mathscr{B}_{\perp}&=\mathscr{B}_{\perp}^{-1}\Pi_{\mathbb{S}_{0}}^{\perp}(\mathcal{L}_{\varepsilon r}-\varepsilon\partial_{\theta}\mathcal{R})\mathscr{B}_{\perp}
					\\
					&=\mathscr{B}_{\perp}^{-1}\Pi_{\mathbb{S}_{0}}^{\perp}\mathcal{L}_{\varepsilon r}\mathscr{B}\Pi_{\mathbb{S}_{0}}^{\perp}-\mathscr{B}_{\perp}^{-1}\Pi_{\mathbb{S}_{0}}^{\perp}\mathcal{L}_{\varepsilon r}\Pi_{\mathbb{S}_{0}}\mathscr{B}\Pi_{\mathbb{S}_{0}}^{\perp}-\varepsilon\mathscr{B}_{\perp}^{-1}\Pi_{\mathbb{S}_{0}}^{\perp}\partial_{\theta}\mathcal{R}\mathscr{B}_{\perp}.
				\end{align*}
				According to the definitions of $\mathfrak{L}_{\varepsilon r}$ and $\mathcal{L}_{\varepsilon r}$ seen  in Proposition \ref{Action on the non local part} and in Lemma \ref{lemma general form of the linearized operator} and using \eqref{link non local operators at state r and at equilibrium}, one has in the Cantor set $\mathcal{O}_{\infty,n}^{\gamma,\tau_{1}}(i_{0})$
				$$
				\mathcal{L}_{\varepsilon r}\mathscr{B}=\mathscr{B}\mathfrak{L}_{\varepsilon r}\quad\hbox{and}\quad
				\mathcal{L}_{\varepsilon r}=\omega\cdot\partial_{\varphi}+\partial_{\theta}\left(V_{\varepsilon r}\cdot\right)-\partial_{\theta}\mathbf{L}_{\varepsilon r,1}-\partial_{\theta}\mathcal{K}_{\lambda}\ast\cdot$$
				and therefore
				\begin{align*}
					\mathscr{B}_{\perp}^{-1}\widehat{\mathcal{L}}_{\omega}\mathscr{B}_{\perp}=&\mathscr{B}_{\perp}^{-1}\Pi_{\mathbb{S}_{0}}^{\perp}\mathscr{B}\mathfrak{L}_{\varepsilon r}\Pi_{\mathbb{S}_{0}}^{\perp}-\mathscr{B}_{\perp}^{-1}\Pi_{\mathbb{S}_{0}}^{\perp}\left(\partial_{\theta}\left(V_{\varepsilon r}\cdot\right)-\partial_{\theta}\mathbf{L}_{\varepsilon r,1}\right)\Pi_{\mathbb{S}_{0}}\mathscr{B}\Pi_{\mathbb{S}_{0}}^{\perp}-\varepsilon\mathscr{B}_{\perp}^{-1}\partial_{\theta}\mathcal{R}\mathscr{B}_{\perp},
				\end{align*}
				where we have used the identities
				$$
				\mathscr{B}_{\perp}^{-1}\Pi_{\mathbb{S}_{0}}^{\perp}=\mathscr{B}_{\perp}^{-1}\quad\hbox{and}\quad [\Pi_{\mathbb{S}_{0}}^{\perp},T]=0=[\Pi_{\mathbb{S}_{0}},T], 
				$$
				for any  Fourier multiplier $T$. The structure of $\mathfrak{L}_{\varepsilon r}$ is detailed  in Proposition \ref{Action on the non local part}, and from this we deduce that
				\begin{align*}\Pi_{\mathbb{S}_{0}}^{\perp}\mathscr{B}\mathfrak{L}_{\varepsilon r}\Pi_{\mathbb{S}_{0}}^{\perp}&=\Pi_{\mathbb{S}_{0}}^{\perp}\mathscr{B}\big(\omega\cdot\partial_{\varphi}+c_{i_{0}}\partial_{\theta}-\partial_{\theta}\mathcal{K}_{\lambda}\ast\cdot+\partial_{\theta}\mathfrak{R}_{\varepsilon r}+\mathtt{E}_{n}^{0}\big)\Pi_{\mathbb{S}_{0}}^{\perp}\\
					&= \Pi_{\mathbb{S}_{0}}^{\perp}\mathscr{B}\Pi_{\mathbb{S}_{0}}^{\perp}\big(\omega\cdot\partial_{\varphi}+c_{i_{0}}\partial_{\theta}-\partial_{\theta}\mathcal{K}_{\lambda}\ast\cdot\big)+\Pi_{\mathbb{S}_{0}}^{\perp}\mathscr{B}\partial_{\theta}\mathfrak{R}_{\varepsilon r}\Pi_{\mathbb{S}_{0}}^{\perp}+\Pi_{\mathbb{S}_{0}}^{\perp}\mathscr{B}\mathtt{E}_{n}^{0}\Pi_{\mathbb{S}_{0}}^{\perp}\\
					&=\mathscr{B}_\perp\big(\omega\cdot\partial_{\varphi}+c_{i_{0}}\partial_{\theta}-\partial_{\theta}\mathcal{K}_{\lambda}\ast\cdot\big)+\Pi_{\mathbb{S}_{0}}^{\perp}\mathscr{B}\partial_{\theta}\mathfrak{R}_{\varepsilon r}\Pi_{\mathbb{S}_{0}}^{\perp}+\Pi_{\mathbb{S}_{0}}^{\perp}\mathscr{B}\mathtt{E}_{n}^{0}\Pi_{\mathbb{S}_{0}}^{\perp}.
				\end{align*}
				It follows that
				\begin{align*}
					\mathscr{B}_{\perp}^{-1}\Pi_{\mathbb{S}_{0}}^{\perp}\mathscr{B}\mathfrak{L}_{\varepsilon r}\Pi_{\mathbb{S}_{0}}^{\perp}
					&=\big(\omega\cdot\partial_{\varphi}+c_{i_{0}}\partial_{\theta}-\partial_{\theta}\mathcal{K}_{\lambda}\ast\cdot\big)\Pi_{\mathbb{S}_0}^{\perp}+\mathscr{B}_{\perp}^{-1}\Pi_{\mathbb{S}_{0}}^{\perp}\mathscr{B}\partial_{\theta}\mathfrak{R}_{\varepsilon r}\Pi_{\mathbb{S}_{0}}^{\perp}+\mathscr{B}_{\perp}^{-1}\Pi_{\mathbb{S}_{0}}^{\perp}\mathscr{B}\mathtt{E}_{n}^{0}\Pi_{\mathbb{S}_{0}}^{\perp}\\
					&=\big(\omega\cdot\partial_{\varphi}+c_{i_{0}}\partial_{\theta}-\partial_{\theta}\mathcal{K}_{\lambda}\ast\cdot\big)\Pi_{\mathbb{S}_0}^{\perp}+\Pi_{\mathbb{S}_{0}}^{\perp}\partial_{\theta}\mathfrak{R}_{\varepsilon r}\Pi_{\mathbb{S}_{0}}^{\perp}+\mathscr{B}_{\perp}^{-1}\mathscr{B}\Pi_{\mathbb{S}_{0}}\partial_{\theta}\mathfrak{R}_{\varepsilon r}\Pi_{\mathbb{S}_{0}}^{\perp}\\
					&\quad+\mathscr{B}_{\perp}^{-1}\Pi_{\mathbb{S}_{0}}^{\perp}\mathscr{B}\mathtt{E}_{n}^{0}\Pi_{\mathbb{S}_{0}}^{\perp}.
				\end{align*}
				Consequently, in the Cantor set $\mathcal{O}_{\infty,n}^{\gamma,\tau_{1}}(i_{0})$, one has the following reduction
				\begin{align}\label{Form1A}
					\nonumber \mathscr{B}_{\perp}^{-1}\widehat{\mathcal{L}}_{\omega}\mathscr{B}_{\perp}=&\big(\omega\cdot\partial_{\varphi}+c_{i_{0}}\partial_{\theta}-\partial_{\theta}\mathcal{K}_{\lambda}\ast\cdot\big)\Pi_{\mathbb{S}_0}^{\perp}+\Pi_{\mathbb{S}_{0}}^{\perp}\partial_{\theta}\mathfrak{R}_{\varepsilon r}\Pi_{\mathbb{S}_{0}}^{\perp}+\mathscr{B}_{\perp}^{-1}\mathscr{B}\Pi_{\mathbb{S}_{0}}\partial_{\theta}\mathfrak{R}_{\varepsilon r}\Pi_{\mathbb{S}_{0}}^{\perp}\\
					\nonumber &-\mathscr{B}_{\perp}^{-1}\Pi_{\mathbb{S}_{0}}^{\perp}\left(\partial_{\theta}\left(V_{\varepsilon r}\cdot\right)-\partial_{\theta}\mathbf{L}_{\varepsilon r,1}\right)\Pi_{\mathbb{S}_{0}}\mathscr{B}\Pi_{\mathbb{S}_{0}}^{\perp}-\varepsilon\mathscr{B}_{\perp}^{-1}\partial_{\theta}\mathcal{R}\mathscr{B}_{\perp}+\mathscr{B}_{\perp}^{-1}\Pi_{\mathbb{S}_{0}}^{\perp}\mathscr{B}\mathtt{E}_{n}^{0}\Pi_{\mathbb{S}_{0}}^{\perp}\\
					:=&\big(\omega\cdot\partial_{\varphi}+c_{i_{0}}\partial_{\theta}-\partial_{\theta}\mathcal{K}_{\lambda}\ast\cdot\big)\Pi_{\mathbb{S}_0}^{\perp}+\mathscr{R}_0+\mathtt{E}_{n}^{1},
				\end{align}
				where we set
				$$\mathtt{E}_{n}^{1}:=\mathscr{B}_{\perp}^{-1}\Pi_{\mathbb{S}_{0}}^{\perp}\mathscr{B}\mathtt{E}_{n}^{0}\Pi_{\mathbb{S}_{0}}^{\perp}.
				$$
				Notice that the estimates \eqref{differences mu0} are simple reformulations of \eqref{estimate r1} and \eqref{difference ci} since $\Delta_{12}r^1=\Delta_{12}c_i$.\\
				\textbf{(ii)} By using \eqref{estimate for mathscrBperp and its inverse}, \eqref{estimate on the first reduction operator and its inverse}, the continuity of the projectors,  \eqref{estim En0 s0} and \eqref{small-3}, one obtains
%				\begin{align*}
%					\|\mathtt{E}_{n}^{1}\rho\|_{q,s}^{\gamma,\mathcal{O}}&\lesssim\|\mathscr{B}\mathtt{E}_{n}^{0}\Pi_{\mathbb{S}_0}^{\perp}\rho\|_{q,s}^{\gamma,\mathcal{O}}+\varepsilon\gamma^{-1}\|\mathfrak{I}_{0}\|_{q,s+\sigma_{3}}^{\gamma,\mathcal{O}}\|\mathscr{B}\mathtt{E}_{n}^{0}\Pi_{\mathbb{S}_0}^{\perp}\rho\|_{q,s_{0}}^{\gamma,\mathcal{O}}\\
%					&\lesssim\|\mathtt{E}_{n}^{0}\Pi_{\mathbb{S}_0}^{\perp}\rho\|_{q,s}^{\gamma,\mathcal{O}}+\varepsilon\gamma^{-1}\|\mathfrak{I}_{0}\|_{q,s+\sigma_{3}}^{\gamma,\mathcal{O}}\|\mathtt{E}_{n}^{0}\Pi_{\mathbb{S}_0}^{\perp}\rho\|_{q,s_{0}}^{\gamma,\mathcal{O}}\\
%					&\lesssim\|\rho\|_{q,s+2}^{\gamma,\mathcal{O}}+\varepsilon\gamma^{-1}\|\mathfrak{I}_{0}\|_{q,s+\sigma_{3}}^{\gamma,\mathcal{O}}\|\rho\|_{q,s_{0}+2}^{\gamma,\mathcal{O}}
%				\end{align*}
%				and 
				\begin{align*}
					\|\mathtt{E}_{n}^{1}\rho\|_{q,s_{0}}^{\gamma,\mathcal{O}}&\lesssim\|\mathscr{B}\mathtt{E}_{n}^{0}\Pi_{\mathbb{S}_0}^{\perp}\rho\|_{q,s_{0}}^{\gamma,\mathcal{O}}\\
					&\lesssim\|\mathtt{E}_{n}^{0}\Pi_{\mathbb{S}_0}^{\perp}\rho\|_{q,s_{0}}^{\gamma,\mathcal{O}}\\
					&\lesssim\varepsilon N_{0}^{\mu_{2}}N_{n+1}^{-\mu_{2}}\|\rho\|_{q,s_{0}+2}^{\gamma,\mathcal{O}}.
				\end{align*}
				\textbf{(iii)} Now, we shall prove the following estimates,
				\begin{equation}\label{F00}
					\max_{k\in\{0,1\}}\|\partial_{\theta}^{k}\mathscr{R}_{0}\|_{\textnormal{\tiny{O-d}},q,s}^{\gamma,\mathcal{O}}\lesssim\varepsilon\gamma^{-1}\left(1+\| \mathfrak{I}_{0}\|_{q,s+\sigma_{3}}^{\gamma,\mathcal{O}}\right)
				\end{equation}
			and
			\begin{equation}\label{F00-1}
				\|\Delta_{12}\mathscr{R}_{0}\|_{\textnormal{\tiny{O-d}},q,\overline{s}_{h}+\mathtt{p}}^{\gamma,\mathcal{O}}\lesssim\varepsilon\gamma^{-1}\|\Delta_{12}i\|_{q,\overline{s}_{h}+\mathtt{p}+\sigma_{3}}^{\gamma,\mathcal{O}}.
			\end{equation}
		To do that, we shall study separately the different terms appearing in \eqref{Form1A} in the definition of $\mathscr{R}_0.$\\
		Notice that in the various estimates below, we use the notation $\sigma_3$ to denote some loss of regularity. This index depends only on $\tau_{1}, q,d,s_0$ and may change increasingly from one line to another and it  is always taken greater than the $\sigma_2$ introduced in Proposition \ref{Action on the non local part}.\\
				$\blacktriangleright$ {\it Study of the term} $\Pi_{\mathbb{S}_{0}}^{\perp}\partial_{\theta}\mathfrak{R}_{\varepsilon r}\Pi_{\mathbb{S}_{0}}^{\perp}$. One gets easily according to \eqref{estim partial thetamathfrakRr} and \eqref{differences partial thetamathfrakRr}
				\begin{align}\label{F01}
					\max_{k\in\{0,1\}}\|\partial_{\theta}^{k}\Pi_{\mathbb{S}_{0}}^{\perp}\partial_{\theta}\mathfrak{R}_{\varepsilon r}\Pi_{\mathbb{S}_{0}}^{\perp}\|_{\textnormal{\tiny{O-d}},q,s}^{\gamma,\mathcal{O}}&\lesssim\max_{k\in\{0,1,2\}}\|\partial_{\theta}^{k}\mathfrak{R}_{\varepsilon r}\|_{\textnormal{\tiny{O-d}},q,s}^{\gamma,\mathcal{O}}\nonumber\\
					&\lesssim\varepsilon\gamma^{-1}\left(1+\| \mathfrak{I}_{0}\|_{q,s+\sigma_{3}}^{\gamma,\mathcal{O}}\right)
				\end{align}
			and
			\begin{align}\label{F01-1}
				\|\Delta_{12}\big(\Pi_{\mathbb{S}_{0}}^{\perp}\partial_{\theta}\mathfrak{R}_{\varepsilon r}\Pi_{\mathbb{S}_{0}}^{\perp}\big)\|_{\textnormal{\tiny{O-d}},q,\overline{s}_{h}+\mathtt{p}}^{\gamma,\mathcal{O}}&\lesssim\|\Delta_{12}\partial_{\theta}\mathfrak{R}_{\varepsilon r}\|_{\textnormal{\tiny{O-d}},q,\overline{s}_{h}+\mathtt{p}}^{\gamma,\mathcal{O}}\nonumber\\
				&\lesssim\varepsilon\gamma^{-1}\|\Delta_{12}i\|_{q,\overline{s}_{h}+\mathtt{p}+\sigma_{3}}^{\gamma,\mathcal{O}}.
			\end{align}
				$\blacktriangleright$ {\it Study of the term} $\mathscr{B}_{\perp}^{-1}\mathscr{B}\Pi_{\mathbb{S}_{0}}\partial_{\theta}\mathfrak{R}_{\varepsilon r}\Pi_{\mathbb{S}_{0}}^{\perp}$. Using the first point of Proposition \ref{projection in the normal directions} yields 
				\begin{align}\label{DT-U}
					\mathscr{B}_{\perp}^{-1}\mathscr{B}\Pi_{\mathbb{S}_{0}}\partial_{\theta}\mathfrak{R}_{\varepsilon r}\Pi_{\mathbb{S}_{0}}^{\perp}&=\Pi_{\mathbb{S}_{0}}\partial_{\theta}\mathfrak{R}_{\varepsilon r}\Pi_{\mathbb{S}_{0}}^{\perp}-\mathcal{T}_0 \mathcal{S}_1,
				\end{align}
				where
				\begin{align}\label{def op T0}
					\mathcal{T}_0\rho=\sum_{m\in\mathbb{S}_0}\big\langle \rho, \big(\mathcal{B}-\textnormal{Id}\big)g_m  \big\rangle_{L^2_\theta(\T)}\mathscr{B}^{-1}e_m\quad\hbox{and}\quad \mathcal{S}_1:=\mathscr{B}\Pi_{\mathbb{S}_{0}}\partial_{\theta}\mathfrak{R}_{\varepsilon r}\Pi_{\mathbb{S}_{0}}^{\perp}.
				\end{align}
				To estimate the first term,  we use Proposition \ref{Action on the non local part} 
				\begin{align}\label{Z00}
					\max_{k\in\{0,1\}}\|\partial_{\theta}^{k}\Pi_{\mathbb{S}_{0}}\partial_{\theta}\mathfrak{R}_{\varepsilon r}\Pi_{\mathbb{S}_{0}}^{\perp}\|_{\textnormal{\tiny{O-d}},q,s}^{\gamma,\mathcal{O}}&\lesssim\max_{k\in\{0,1,2\}}\|\partial_{\theta}^{k}\mathfrak{R}_{\varepsilon r}\|_{\textnormal{\tiny{O-d}},q,s}^{\gamma,\mathcal{O}}\nonumber\\
					&\lesssim\varepsilon\gamma^{-1}\left(1+\| \mathfrak{I}_{0}\|_{q,s+\sigma_{3}}^{\gamma,\mathcal{O}}\right).
				\end{align}
				As to the second term, we write 
				\begin{align*}
					\mathcal{T}_0\mathcal{S}_1\rho&=\sum_{m\in\mathbb{S}_0}\big\langle \mathcal{S}_1\rho, \big(\mathcal{B}-\textnormal{Id}\big)g_m  \big\rangle_{L^2_\theta(\T)}\mathscr{B}^{-1}e_m\\\
					&=\sum_{m\in\mathbb{S}_0}\big\langle \rho, \mathcal{S}_1^\star\big(\mathcal{B}-\textnormal{Id}\big)g_m  \big\rangle_{L^2_\theta(\T)}\mathscr{B}^{-1}e_m,
				\end{align*}
				where $\mathcal{S}_1^\star$ is the $L^2_\theta(\mathbb{T})$-adjoint of $\mathcal{S}_1.$ This is an integral operator taking the form 
				\begin{align*}
					\big(\mathcal{T}_0\mathcal{S}_1\rho\big)(\varphi,\theta)&=\int_{\T} \mathcal{K}_1(\varphi,\theta,\eta)\rho(\varphi,\eta) d\eta,\\
					\mathcal{K}_1(\varphi,\theta,\eta)&:=\sum_{m\in\mathbb{S}_0} \big(\mathcal{S}_1^\star\big(\mathcal{B}-\textnormal{Id}\big)g_m\big)(\varphi,\eta)  \big(\mathscr{B}^{-1}e_m\big)(\varphi,\theta).
				\end{align*}
			Recall from Proposition \ref{Action on the non local part} that   $\mathfrak{R}_{\varepsilon r}$ is self-adjoint and using Lemma \ref{algeb1} we have the identities 
			$\mathscr{B}^\star=\mathcal{B}^{-1}$ and $\mathcal{B}^\star=\mathscr{B}^{-1},$ then 
			\begin{equation}\label{S1star}
			\mathcal{S}_1^\star=-\Pi_{\mathbb{S}_{0}}^{\perp}\mathfrak{R}_{\varepsilon r}\partial_{\theta}\Pi_{\mathbb{S}_{0}}\mathcal{B}^{-1}.
			\end{equation}
		Therefore, combining \eqref{sym gm}, \eqref{symmetry for beta} and \eqref{symmertry mathbbK} imply
		\begin{equation}\label{sym K1}
			\mathcal{K}_1(-\varphi,-\theta,-\eta)=-\mathcal{K}_1(\varphi,\theta,\eta)\in\mathbb{R}.
		\end{equation}
				Applying Lemma \ref{lemma symmetry and reversibility} combined with the law products yield for any $k\in\mathbb{N}$
				\begin{align}\label{Z01}
					\|\partial_\theta^k \mathcal{T}_0\mathcal{S}_1\|_{\textnormal{\tiny{O-d}},q,s}^{\gamma,\mathcal{O}}&\lesssim \bigintssss_{\T}\|(\partial_\theta^k\mathcal{K}_1)(\ast,\cdot,\centerdot,\eta+\centerdot)\|_{q,s+s_{0}}^{\gamma,\mathcal{O}} d\eta\\
					\nonumber&\lesssim \sum_{m\in\mathbb{S}_0}\Big(\big\|\mathcal{S}_1^\star\big(\mathcal{B}-\textnormal{Id}\big)g_m\big\|_{q,s+s_0}^{\gamma,\mathcal{O}}\big\|\mathscr{B}^{-1}e_m\big\|_{q,s_0+k}^{\gamma,\mathcal{O}}+\big\|\mathcal{S}_1^\star\big(\mathcal{B}-\textnormal{Id}\big)g_m\big\|_{q,s_0}^{\gamma,\mathcal{O}}\big\|\mathscr{B}^{-1}e_m\big\|_{q,s+s_0+k}^{\gamma,\mathcal{O}} \Big).
				\end{align}
				Remark that \eqref{S1star} implies
				$$
				\mathcal{S}_1^\star\big(\mathcal{B}-\textnormal{Id}\big)g_m=-\Pi_{\mathbb{S}_{0}}^{\perp}\mathfrak{R}_{\varepsilon r}\partial_{\theta}\Pi_{\mathbb{S}_{0}}\big(\textnormal{Id}-\mathcal{B}^{-1}\big)g_m.
				$$
				Hence  according to  Lemma \ref{properties of Toeplitz in time operators} combined with  Proposition \ref{Action on the non local part}  we find
				\begin{align}\label{D-L-O}
					\nonumber \|\mathcal{S}_1^\star\big(\mathcal{B}-\textnormal{Id}\big)g_m\|_{q,s}^{\gamma,\mathcal{O}}&\lesssim \|\mathfrak{R}_{\varepsilon r}\|_{\textnormal{\tiny{O-d}},q,s}^{\gamma,\mathcal{O}}\|\partial_{\theta}\Pi_{\mathbb{S}_{0}}\big(\textnormal{Id}-\mathcal{B}^{-1}\big)g_m\|_{q,s_0}^{\gamma,\mathcal{O}}+\|\mathfrak{R}_{\varepsilon r}\|_{\textnormal{\tiny{O-d}},q,s_0}^{\gamma,\mathcal{O}}\|\partial_{\theta}\Pi_{\mathbb{S}_{0}}\big(\textnormal{Id}-\mathcal{B}^{-1}\big)g_m\|_{q,s}^{\gamma,\mathcal{O}}\\
					\nonumber &\lesssim \varepsilon\gamma^{-1}\left(1+\|\mathfrak{I}_{0}\|_{q,s+\sigma_{3}}^{\gamma,\mathcal{O}}\right)\|\big(\textnormal{Id}-\mathcal{B}^{-1}\big)g_m\|_{q,s_0+1}^{\gamma,\mathcal{O}}\\
					&\quad+ \varepsilon\gamma^{-1}\left(1+\|\mathfrak{I}_{0}\|_{q,s_0+\sigma_{3}}^{\gamma,\mathcal{O}}\right)\|\big(\textnormal{Id}-\mathcal{B}^{-1}\big)g_m\|_{q,s+1}^{\gamma,\mathcal{O}}.
				\end{align}
				Using   \eqref{estimate on the first reduction operator and its inverse}  together  with Lemma \ref{Lema-decomp} and  the smallness condition \eqref{small-3} leads to 
				\begin{align}\label{id-B-1 gm}
					\big\|\big(\textnormal{Id}-\mathcal{B}^{-1}\big)g_m\big\|_{q,s}^{\gamma,\mathcal{O}}&\lesssim\|g_m\|_{q,s}^{\gamma ,\mathcal{O}}+{\varepsilon\gamma ^{-1}}\| \mathfrak{I}_{0}\|_{q,s+\sigma_{3}}^{\gamma ,\mathcal{O}}\|g_m\|_{q,s_{0}}^{\gamma ,\mathcal{O}}\nonumber\\
					&\lesssim\sup_{k,m\in\mathbb{S}_0}\|\alpha_{k,m}\|_{q,H^s_\varphi}^{\gamma ,\mathcal{O}}+{\varepsilon\gamma ^{-1}}\| \mathfrak{I}_{0}\|_{q,s+\sigma_{3}}^{\gamma ,\mathcal{O}}\sup_{k,m\in\mathbb{S}_0}\|\alpha_{k,m}\|_{q,H^{s_0}_\varphi}^{\gamma ,\mathcal{O}}\nonumber\\
					&\lesssim 1+{\varepsilon\gamma ^{-1}}\| \mathfrak{I}_{0}\|_{q,s+\sigma_{3}}^{\gamma ,\mathcal{O}}.
				\end{align}
				Therefore, inserting this estimate into \eqref{D-L-O} and using  \eqref{small-3} allow to get
				\begin{align*}
					\|\mathcal{S}_1^\star\big(\mathcal{B}-\textnormal{Id}\big)g_m\|_{q,s}^{\gamma,\mathcal{O}}&\lesssim \varepsilon\gamma^{-1}\left(1+\|\mathfrak{I}_{0}\|_{q,s+\sigma_{3}}^{\gamma,\mathcal{O}}\right).
				\end{align*}
				Plugging this estimate into \eqref{Z01} and using  \eqref{estimate on the first reduction operator and its inverse}  ensure
				\begin{align}\label{Z02}
					\max_{k\in\{0,1\}}\|\partial_\theta^k \mathcal{T}_0\mathcal{S}_1\|_{\textnormal{\tiny{O-d}},q,s}^{\gamma,\mathcal{O}}&\lesssim \varepsilon\gamma^{-1}\left(1+\|\mathfrak{I}_{0}\|_{q,s+\sigma_{3}}^{\gamma,\mathcal{O}}\right).
				\end{align}
			Consequently, by  combining \eqref{DT-U}, \eqref{Z00} and \eqref{Z02}, we find
			\begin{equation}\label{F02}
				\max_{k\in\{0,1\}}\|\partial_\theta^k \mathscr{B}_{\perp}^{-1}\mathscr{B}\Pi_{\mathbb{S}_0}\partial_{\theta}\mathfrak{R}_{\varepsilon r}\Pi_{\mathbb{S}_0}^{\perp}\|_{\textnormal{\tiny{O-d}},q,s}^{\gamma,\mathcal{O}}\lesssim \varepsilon\gamma^{-1}\left(1+\|\mathfrak{I}_{0}\|_{q,s+\sigma_{3}}^{\gamma,\mathcal{O}}\right).
			\end{equation}
			We now turn to the difference estimate. From  \eqref{DT-U}, it is obvious that
			\begin{equation}\label{DT-U1}
				\Delta_{12}\big(\mathscr{B}_{\perp}^{-1}\mathscr{B}\Pi_{\mathbb{S}_{0}}\partial_{\theta}\mathfrak{R}_{\varepsilon r}\Pi_{\mathbb{S}_{0}}^{\perp}\big)=\Pi_{\mathbb{S}_{0}}\Delta_{12}\partial_{\theta}\mathfrak{R}_{\varepsilon r}\Pi_{\mathbb{S}_{0}}^{\perp}-\Delta_{12}(\mathcal{T}_{0}\mathcal{S}_{1}).
			\end{equation}
			To estimate the first term, we use \eqref{differences partial thetamathfrakRr}
			\begin{align}\label{D-00}
				\|\Pi_{\mathbb{S}_{0}}\Delta_{12}\partial_{\theta}\mathfrak{R}_{\varepsilon r}\Pi_{\mathbb{S}_{0}}^{\perp}\|_{\textnormal{\tiny{O-d}},q,\overline{s}_{h}+\mathtt{p}}^{\gamma,\mathcal{O}}&\lesssim \|\Delta_{12}\partial_{\theta}\mathfrak{R}_{\varepsilon r}\|_{\textnormal{\tiny{O-d}},q,\overline{s}_{h}+\mathtt{p}}^{\gamma,\mathcal{O}}\nonumber\\
				&\lesssim \varepsilon\gamma^{-1}\|\Delta_{12}i\|_{q,\overline{s}_{h}+\mathtt{p}+\sigma_{3}}^{\gamma,\mathcal{O}}.
			\end{align}
		As to the second term, we notice that $\Delta_{12}\big(\mathcal{T}_{0}\mathcal{S}_{1}\big)$ is an integral operator whose  kernel $\Delta_{12}\mathcal{K}_{1}$ is
		\begin{align*}
			\Delta_{12}\mathcal{K}_{1}(\varphi,\theta,\eta)&=\sum_{m\in\mathbb{S}_{0}}\Delta_{12}\big(\mathcal{S}_{1}^{\star}\big(\mathcal{B}-\textnormal{Id}\big)g_{m}\big)(\varphi,\eta)\big(\mathscr{B}_{r_{1}}e_{m}\big)(\varphi,\theta)\\
			&\quad+\big(\mathcal{S}_{1}^{\star}\big(\mathcal{B}-\textnormal{Id}\big)g_{m}\big)_{r_{2}}(\varphi,\eta)\big(\Delta_{12}\mathscr{B}e_{m}\big)(\varphi,\theta).
		\end{align*}
		Hence, using  Lemma  \ref{lemma symmetry and reversibility}-(ii) together with the law products  we deduce that
		\begin{align*}
			\|\Delta_{12}\mathcal{T}_{0}\mathcal{S}_{1}\|_{\textnormal{\tiny{O-d}},q,\overline{s}_{h}+\mathtt{p}}^{\gamma,\mathcal{O}}&\lesssim\int_{\mathbb{T}}\|\Delta_{12}\mathcal{K}_{1}(\ast,\cdot,\centerdot,\eta+\centerdot)\|_{q,\overline{s}_{h}+\mathtt{p}+s_{0}}^{\gamma,\mathcal{O}}d\eta\\
			&\lesssim \sum_{m\in\mathbb{S}_{0}}\|\Delta_{12}\big(\mathcal{S}_{1}^{\star}\big(\mathcal{B}-\textnormal{Id}\big)g_{m}\big)\|_{q,\overline{s}_{h}+\mathtt{p}+s_{0}}^{\gamma,\mathcal{O}}\|\mathscr{B}_{r_{1}}e_{m}\|_{q,\overline{s}_{h}+\mathtt{p}+s_{0}}^{\gamma,\mathcal{O}}\\
			&\quad +\|\big(\mathcal{S}_{1}^{\star}\big(\mathcal{B}-\textnormal{Id}\big)g_{m}\big)_{r_{2}}\|_{q,\overline{s}_{h}+\mathtt{p}+s_{0}}^{\gamma,\mathcal{O}}\|\Delta_{12}\mathscr{B}e_{m}\|_{q,\overline{s}_{h}+\mathtt{p}+s_{0}}^{\gamma,\mathcal{O}}.
		\end{align*}
	Notice that by Taylor Formula and \eqref{difference beta} (applied with $\mathtt{p}$ replaced by $\mathtt{p}+s_{0}$), one has
	\begin{equation}\label{D-01}
		\sup_{m\in\mathbb{S}_{0}}\|\Delta_{12}\mathscr{B}^{-1}e_{m}\|_{q,\overline{s}_{h}+\mathtt{p}+s_{0}}^{\gamma,\mathcal{O}}\lesssim\varepsilon\gamma^{-1}\|\Delta_{12}i\|_{q,\overline{s}_{h}+\mathtt{p}+\sigma_{3}}.
	\end{equation}
	On the other hand, we have
	$$\Delta_{12}\mathcal{S}_{1}^{\star}=-\Pi_{\mathbb{S}_{0}}^{\perp}\Delta_{12}\mathfrak{R}_{\varepsilon r}\partial_{\theta}\Pi_{\mathbb{S}_{0}}\mathcal{B}_{r_{1}}^{-1}-\Pi_{\mathbb{S}_{0}}^{\perp}\mathfrak{R}_{\varepsilon r_{2}}\partial_{\theta}\Pi_{\mathbb{S}_{0}}\Delta_{12}\mathcal{B}^{-1},$$
	leading to 
	\begin{align*}
		\Delta_{12}\big(\mathcal{S}_{1}^{\star}\big(\mathcal{B}-\textnormal{Id}\big)g_{m}\big)=&-\Pi_{\mathbb{S}_{0}}^{\perp}\Delta_{12}\mathfrak{R}_{\varepsilon r}\partial_{\theta}\Pi_{\mathbb{S}_{0}}\big(\textnormal{Id}-\mathcal{B}_{r_{1}}^{-1}\big)g_{m,r_{1}}-\Pi_{\mathbb{S}_{0}}^{\perp}\mathfrak{R}_{\varepsilon r_{2}}\partial_{\theta}\Pi_{\mathbb{S}_{0}}\Delta_{12}\mathcal{B}^{-1}\big(\mathcal{B}_{r_{1}}-\textnormal{Id}\big)g_{m,r_{1}}\\
		&+\mathcal{S}_{1,r_{2}}^{\star}\big(\Delta_{12}\mathcal{B}\big)g_{m,r_{1}}+\mathcal{S}_{1,r_{2}}^{\star}\big(\mathcal{B}_{r_{2}}-\textnormal{Id}\big)\Delta_{12}g_{m}.
	\end{align*}
AAccording to Lemma \ref{properties of Toeplitz in time operators}, we obtain
\begin{align*}
	\|\Pi_{\mathbb{S}_{0}}^{\perp}\Delta_{12}\mathfrak{R}_{\varepsilon r}\partial_{\theta}\Pi_{\mathbb{S}_{0}}\big(\textnormal{Id}-\mathcal{B}_{r_{1}}^{-1}\big)g_{m,r_{1}}\|_{q,\overline{s}_{h}+\mathtt{p}+s_{0}}^{\gamma,\mathcal{O}}&\lesssim\|\Delta_{12}\mathfrak{R}_{\varepsilon r}\|_{\textnormal{\tiny{O-d}},q,\overline{s}_{h}+\mathtt{p}+s_{0}}^{\gamma,\mathcal{O}}\|\big(\textnormal{Id}-\mathcal{B}_{r_{1}}^{-1}\big)g_{m,r_{1}}\|_{q,\overline{s}_{h}+\mathtt{p}+s_{0}+1}^{\gamma,\mathcal{O}}.
\end{align*}
From \eqref{id-B-1 gm}, one has
\begin{align*}
	\|\big(\textnormal{Id}-\mathcal{B}_{r_{1}}^{-1}\big)g_{m,r_{1}}\|_{q,s}^{\gamma,\mathcal{O}}\lesssim1+\varepsilon\gamma^{-1}\|\mathfrak{I}_{1}\|_{q,s+\sigma_3}^{\gamma,\mathcal{O}}.
\end{align*}
Thus, from \eqref{differences partial thetamathfrakRr} and\eqref{small-3}, we infer
\begin{equation}\label{D-02}
	\|\Pi_{\mathbb{S}_{0}}^{\perp}\Delta_{12}\mathfrak{R}_{\varepsilon r}\partial_{\theta}\Pi_{\mathbb{S}_{0}}\big(\textnormal{Id}-\mathcal{B}_{r_{1}}^{-1}\big)g_{m,r_{1}}\|_{q,\overline{s}_{h}+\mathtt{p}+s_{0}}^{\gamma,\mathcal{O}}\lesssim\varepsilon\gamma^{-1}\|\Delta_{12}i\|_{q,\overline{s}_{h}+\mathtt{p}+\sigma_3}^{\gamma,\mathcal{O}}.
\end{equation}
Applying   Lemma \ref{properties of Toeplitz in time operators}, \eqref{estim partial thetamathfrakRr} and \eqref{small-3} we deduce that
\begin{align*}
	\|\mathcal{S}_{1,r_{2}}^{\star}\big(\Delta_{12}\mathcal{B}\big)g_{m,r_{1}}\|_{q,\overline{s}_{h}+\mathtt{p}+s_{0}}^{\gamma,\mathcal{O}}&\lesssim\|\mathfrak{R}_{\varepsilon r_{2}}\|_{\textnormal{\tiny{O-d}},q,\overline{s}_{h}+s_{0}}\|\mathcal{B}_{r_{2}}\big(\Delta_{12}\mathcal{B}\big)g_{m,r_{1}}\|_{q,\overline{s}_{h}+\mathtt{p}+s_{0}+1}^{\gamma,\mathcal{O}}\\
	&\lesssim\|\mathcal{B}_{r_{2}}\big(\Delta_{12}\mathcal{B}\big)g_{m,r_{1}}\|_{q,\overline{s}_{h}+\mathtt{p}+s_{0}+1}^{\gamma,\mathcal{O}}.
\end{align*}
To estimate the right hand side member, it suffices to use \eqref{estimate on the first reduction operator and its inverse} and  \eqref{small-3}, leading to
$$\|\mathcal{B}_{r_{2}}\big(\Delta_{12}\mathcal{B}\big)g_{m,r_{1}}\|_{q,\overline{s}_{h}+\mathtt{p}+s_{0}+1}^{\gamma,\mathcal{O}}\lesssim\|\big(\Delta_{12}\mathcal{B}\big)g_{m,r_{1}}\|_{q,\overline{s}_{h}+\mathtt{p}+s_{0}+1}^{\gamma,\mathcal{O}}.$$
By Taylor Formula, we may write
$$\Delta_{12}\mathcal{B}\rho(\theta)=\Delta_{12}\beta(\theta)\int_{0}^{1}\partial_{\theta}\rho\big(\theta+\beta_{2}(\theta)+t\Delta_{12}\beta(\theta)\big)dt.$$
It follows from the  law products in Lemma \ref{Lem-lawprod}, \eqref{difference beta} and \eqref{small-3} that
\begin{align*}
	\|\big(\Delta_{12}\mathcal{B}\big)g_{m,r_{1}}\|_{q,\overline{s}_{h}+\mathtt{p}+1}^{\gamma,\mathcal{O}}&\lesssim\|\Delta_{12}\beta\|_{q,\overline{s}_{h}+\mathtt{p}+s_{0}+1}^{\gamma,\mathcal{O}}\|g_{m,r_{1}}\|_{q,\overline{s}_{h}+\mathtt{p}+s_{0}+2}^{\gamma,\mathcal{O}}\\
	&\lesssim\varepsilon\gamma^{-1}\|\Delta_{12}i\|_{q,\overline{s}_{h}+\mathtt{p}+\sigma_{3}}^{\gamma,\mathcal{O}}.
\end{align*}
Thus
\begin{equation}\label{D-03}
	\|\mathcal{S}_{1,r_{2}}^{\star}\big(\Delta_{12}\mathcal{B}\big)g_{m,r_{1}}\|_{q,\overline{s}_{h}+\mathtt{p}+s_{0}}^{\gamma,\mathcal{O}}\lesssim\varepsilon\gamma^{-1}\|\Delta_{12}i\|_{q,\overline{s}_{h}+\mathtt{p}+\sigma_{3}}^{\gamma,\mathcal{O}}.
\end{equation}
In the same way, using Taylor Formula together with \eqref{difference beta}, we get
\begin{equation}\label{D-04}
	\|\Pi_{\mathbb{S}_{0}}^{\perp}\mathfrak{R}_{\varepsilon r_{2}}\partial_{\theta}\Pi_{\mathbb{S}_{0}}\Delta_{12}\mathcal{B}^{-1}\big(\mathcal{B}_{r_{1}}-\textnormal{Id}\big)g_{m,r_{1}}\|_{q,\overline{s}_{h}+\mathtt{p}+s_{0}}^{\gamma,\mathcal{O}}\lesssim\varepsilon\gamma^{-1}\|\Delta_{12}i\|_{q,\overline{s}_{h}+\mathtt{p}+\sigma_{3}}^{\gamma,\mathcal{O}}.
\end{equation}
By Lemma \ref{properties of Toeplitz in time operators},\eqref{estim partial thetamathfrakRr} and \eqref{small-3}, one finds
\begin{align*}
	\|\mathcal{S}_{1,r_{2}}^{\star}\big(\mathcal{B}_{r_{2}}-\textnormal{Id}\big)\Delta_{12}g_{m}\|_{q,\overline{s}_{h}+\mathtt{p}+s_{0}}^{\gamma,\mathcal{O}}&\lesssim\|\mathfrak{R}_{\varepsilon r_{2}}\|_{\textnormal{\tiny{O-d}},q,\overline{s}_{h}+\mathtt{p}+s_{0}}^{\gamma,\mathcal{O}}\|\big(\mathcal{B}_{r_{2}}-\textnormal{Id}\big)\Delta_{12}g_{m}\|_{q,\overline{s}_{h}+\mathtt{p}+s_{0}+1}^{\gamma,\mathcal{O}}\\
	&\lesssim\|\big(\mathcal{B}_{r_{2}}-\textnormal{Id}\big)\Delta_{12}g_{m}\|_{q,\overline{s}_{h}+\mathtt{p}+s_{0}+1}^{\gamma,\mathcal{O}}.
\end{align*}
Applying  \eqref{estimate on the first reduction operator and its inverse} and \eqref{small-3}, we obtain
$$\|\big(\mathcal{B}_{r_{2}}-\textnormal{Id}\big)\Delta_{12}g_{m}\|_{q,\overline{s}_{h}+\mathtt{p}+s_{0}+1}^{\gamma,\mathcal{O}}\lesssim\|\Delta_{12}g_{m}\|_{q,\overline{s}_{h}+\mathtt{p}+s_{0}+1}^{\gamma,\mathcal{O}}.$$
Using \eqref{diff gm} (applied with $\mathtt{p}=s_{0}+1$), we finally get
\begin{equation}\label{D-05}
	\|\mathcal{S}_{1,r_{2}}^{\star}\big(\mathcal{B}_{r_{2}}-\textnormal{Id}\big)\Delta_{12}g_{m}\|_{q,\overline{s}_{h}+\mathtt{p}+s_{0}}^{\gamma,\mathcal{O}}\lesssim\varepsilon\gamma^{-1}\|\Delta_{12}i\|_{q,\overline{s}_{h}+\mathtt{p}+\sigma_{3}}^{\gamma,\mathcal{O}}.
\end{equation}
Gathering \eqref{D-01}, \eqref{D-02}, \eqref{D-03}, \eqref{D-04} and \eqref{D-05} implies
\begin{equation}\label{D-06}
	\|\Delta_{12}\big(\mathcal{T}_{0}\mathcal{S}_{1}\big)\|_{\textnormal{\tiny{O-d}},q,\overline{s}_{h}+\mathtt{p}}^{\gamma,\mathcal{O}}\lesssim\varepsilon\gamma^{-1}\|\Delta_{12}i\|_{q,\overline{s}_{h}+\mathtt{p}+\sigma_{3}}^{\gamma,\mathcal{O}}.
\end{equation}
Putting together \eqref{DT-U1}, \eqref{D-00} and \eqref{D-06}, one obtains
\begin{equation}\label{F02-1}
	\|\Delta_{12}\big(\mathscr{B}_{\perp}^{-1}\mathscr{B}\Pi_{\mathbb{S}_{0}}\partial_{\theta}\mathfrak{R}_{\varepsilon r}\Pi_{\mathbb{S}_{0}}^{\perp}\big)\|_{\textnormal{\tiny{O-d}},q,\overline{s}_{h}+\mathtt{p}}^{\gamma,\mathcal{O}}\lesssim\varepsilon\gamma^{-1}\|\Delta_{12}i\|_{q,\overline{s}_{h}+\mathtt{p}+\sigma_{3}}^{\gamma,\mathcal{O}}.
\end{equation}
				$\blacktriangleright$ {\it  Study of the term} $ \mathscr{B}_{\perp}^{-1}\Pi_{\mathbb{S}_{0}}^{\perp}\left(\partial_{\theta}\left(V_{\varepsilon r}\cdot\right)-\partial_{\theta}\mathbf{L}_{\varepsilon r,1}\right)\Pi_{\mathbb{S}_{0}}\mathscr{B}\Pi_{\mathbb{S}_{0}}^{\perp}.$
				We first write,
				$$
				\mathscr{B}_{\perp}^{-1}\Pi_{\mathbb{S}_{0}}^{\perp}\left(\partial_{\theta}\left(V_{\varepsilon r}\cdot\right)-\partial_{\theta}\mathbf{L}_{\varepsilon r,1}\right)\Pi_{\mathbb{S}_{0}}\mathscr{B}\Pi_{\mathbb{S}_{0}}^{\perp}:=\mathscr{B}_{\perp}^{-1}\partial_{\theta}\mathcal{S}_2\mathscr{B}\Pi_{\mathbb{S}_{0}}^{\perp},
				$$
				with
				$$
				\mathcal{S}_2=\left(\big(V_{\varepsilon r}-c_{i_0}\big)\cdot -\mathbf{L}_{\varepsilon r,1}\right)\Pi_{\mathbb{S}_{0}}.
				$$
				Notice that to get the above identity we have used   the identity
				$$
				\Pi_{\mathbb{S}_{0}}^{\perp}\partial_{\theta}(c_{i_0}\cdot)\Pi_{\mathbb{S}_{0}}=0.
				$$
				Recall from \eqref{link non local operators at state r and at equilibrium} and \eqref{definition of mathbbK1} that 
				$$
				\mathbf{L}_{\varepsilon r,1}\rho(\varphi,\theta)=\int_{\T}\mathbb{K}_{\varepsilon r,1}(\varphi,\theta,\eta) \rho(\varphi,\eta) d\eta.
				$$
				Then from elementary computations we find
				$$
				\mathcal{S}_2\rho(\varphi,\theta)=\int_{\T}\mathcal{K}_2(\varphi,\theta,\eta) \rho(\varphi,\eta) d\eta,
				$$
				with 
				\begin{align*}
					\mathcal{K}_2(\varphi,\theta,\eta) &:=\big(V_{\varepsilon r}(\varphi,\theta)-c_{i_0}\big)D_{\mathbb{S}_0}(\theta-\eta)-\int_{\T}\mathbb{K}_{\varepsilon r,1}(\varphi,\theta,\eta^\prime)D_{\mathbb{S}_0}(\eta^\prime-\eta)d\eta^\prime,\\
					D_{\mathbb{S}_0}(\theta)&:=\sum_{n\in\mathbb{S}_0} e^{i n\theta}.
				\end{align*}
			Combining \eqref{symmetry kernel K1}, \eqref{symmetry for r}, \eqref{symmetry for Vr} and the change of variables $\eta'\mapsto\eta'$, one gets
				\begin{equation}\label{symmetry mathcal K2}
					\mathcal{K}_2(-\varphi,-\theta,-\eta)=\mathcal{K}_2(\varphi,\theta,\eta)\in\mathbb{R}.
			\end{equation}
				Proceeding as in \eqref{DT-U} we obtain
				\begin{align}\label{mod-T-V}
					\mathscr{B}_{\perp}^{-1}\partial_{\theta}\mathcal{S}_2\mathscr{B}\Pi_{\mathbb{S}_{0}}^{\perp}&=\mathscr{B}^{-1}\partial_{\theta}\mathcal{S}_2\mathscr{B}\Pi_{\mathbb{S}_{0}}^{\perp}-\mathcal{T}_0\partial_{\theta}\mathcal{S}_2\mathscr{B}\Pi_{\mathbb{S}_{0}}^{\perp}.
				\end{align}
				It follows that 
				\begin{align}\label{WL1}
					\big\|\partial_\theta^k\mathscr{B}_{\perp}^{-1}\partial_{\theta}\mathcal{S}_2\mathscr{B}\Pi_{\mathbb{S}_{0}}^{\perp}\big\|_{\textnormal{\tiny{O-d}},q,s}^{\gamma,\mathcal{O}}&\lesssim \big\|\partial_\theta^{k+1}\mathcal{B}^{-1}\mathcal{S}_2\mathscr{B}\big\|_{\textnormal{\tiny{O-d}},q,s}^{\gamma,\mathcal{O}}+\big\|\partial_\theta^k\mathcal{T}_0\partial_{\theta}\mathcal{S}_2\mathscr{B}\big\|_{\textnormal{\tiny{O-d}},q,s}^{\gamma,\mathcal{O}}.
				\end{align}
				The expression  of the first term is similar to that of \eqref{ZY1}, namely, one has
				\begin{align*}\big(\mathcal{B}^{-1}\mathcal{S}_2\mathscr{B}\big)\rho(\varphi,\theta)=\int_{\mathbb{T}}\rho(\varphi,\eta)\widehat{\mathcal{K}}_{2}(\varphi,\theta,\eta)d\eta,
				\end{align*}
				with
				$$\widehat{\mathcal{K}}_{2}(\varphi,\theta,\eta):=\mathcal{K}_{2}\big(\varphi,\theta+\widehat{\beta}(\varphi,\theta),\eta+\widehat{\beta}(\varphi,\eta)\big).$$
				Combining \eqref{symmetry mathcal K2} and \eqref{symmetry for beta}, one gets
					\begin{equation}\label{symmetry mathcal K2hat}
						\widehat{\mathcal{K}}_2(-\varphi,-\theta,-\eta)=\widehat{\mathcal{K}}_2(\varphi,\theta,\eta)\in\mathbb{R}.
				\end{equation}
				Then coming bacl to \eqref{KKLLPP1} and arguing as for \eqref{Es-U-D3}, we find
				\begin{align}\label{Es-U-D4}
					\sup_{k\in\{0,1,2\}} \big\|(\partial_\theta^k\widehat{\mathcal{K}}_{2})\big(\ast,\cdot,\centerdot,\eta+\centerdot\big)\big\|_{q,s}^{\gamma,\mathcal{O}}&\lesssim \varepsilon\gamma ^{-1}\left(1+\| \mathfrak{I}_{0}\|_{q,s+\sigma_{3}}^{\gamma ,\mathcal{O}}\right)\Big(1-\log\big|\sin\big({\eta}/{2} \big)\big|\Big).
				\end{align}
				By virtue of Lemma \ref{lemma symmetry and reversibility} and \eqref{Es-U-D4} we obtain
				\begin{align}\label{Z04}
					\nonumber \sup_{k\in\{0,1\}}\big\|\partial_\theta^{k+1}\mathcal{B}^{-1}\mathcal{S}_2\mathscr{B}\big\|_{\textnormal{\tiny{O-d}},q,s}^{\gamma,\mathcal{O}}&\lesssim \sup_{k\in\{0,1\}}\bigintssss_{\T}\big\|(\partial_\theta^{k+1}\widehat{\mathcal{K}}_{2})\big(\ast,\cdot,\centerdot,\eta+\centerdot\big)\big\|_{q,s+s_0}^{\gamma,\mathcal{O}}d\eta\\
					&\lesssim \varepsilon\gamma ^{-1}\left(1+\| \mathfrak{I}_{0}\|_{q,s+\sigma_{3}}^{\gamma ,\mathcal{O}}\right).
				\end{align}
				Notice that from \eqref{def op T0}, we can write
				\begin{align*}
					\mathcal{T}_{0}\partial_{\theta}\mathcal{S}_{2}\mathscr{B}\rho&=\sum_{m\in\mathbb{S}_{0}}\big\langle \partial_{\theta}\mathcal{S}_{2}\mathscr{B}\rho,\big(\mathcal{B}-\textnormal{Id}\big)g_{m}\big\rangle_{L_{\theta}^{2}(\mathbb{T})}\mathscr{B}^{-1}e_{m}\\
					&=-\sum_{m\in\mathbb{S}_{0}}\big\langle \rho,\mathcal{B}^{-1}\mathcal{S}_{2}^{\star}\partial_{\theta}\big(\mathcal{B}-\textnormal{Id}\big)g_{m}\big\rangle_{L_{\theta}^{2}(\mathbb{T})}\mathscr{B}^{-1}e_{m},
				\end{align*}
				where $\mathcal{S}_{2}^{\star}$ is the adjoint of $\mathcal{S}_{2}$ and is given by
				\begin{equation}\label{def S2star}
					\mathcal{S}_{2}^{\star}=\Pi_{\mathbb{S}_{0}}\big(\big(V_{\varepsilon r}-c_{i_{0}}\big)\cdot-\mathbf{L}_{\varepsilon r,1}\big).
				\end{equation}
				This is an integral operator taking the form
				\begin{align*}
					\big(\mathcal{T}_{0}\partial_{\theta}\mathcal{S}_{2}\mathscr{B}\rho\big)(\varphi,\theta)&=\int_{\mathbb{T}}\mathcal{K}_{3}(\varphi,\theta,\eta)\rho(\varphi,\eta)d\eta,\\
					\mathcal{K}_{3}(\varphi,\theta,\eta)&:=\sum_{m\in\mathbb{S}_{0}}\big(\mathcal{B}^{-1}\mathcal{S}_{2}^{\star}\partial_{\theta}\big(\mathcal{B}-\textnormal{Id}\big)g_{m}\big)(\varphi,\eta)\big(\mathscr{B}^{-1}e_{m}\big)(\varphi,\theta).
				\end{align*}
			According to \eqref{sym gm}, \eqref{symmetry for beta}, \eqref{def S2star}, \eqref{symmetry for r}, \eqref{symmetry for Vr} and \eqref{symmetry kernel K1}, one gets
			\begin{equation}\label{sym K3}
				\mathcal{K}_3(-\varphi,-\theta,-\eta)=-\mathcal{K}_{3}(\varphi,\theta,\eta)\in\mathbb{R}.
			\end{equation}
				On the other hand, applying Lemma \ref{lemma symmetry and reversibility} combined with the law products yield for any $k\in\mathbb{N}$
				\begin{align*}
					&\|\partial_{\theta}^{k}\mathcal{T}_{0}\partial_{\theta}\mathcal{S}_{2}\mathscr{B}\|_{\textnormal{\tiny{O-d}},q,s}^{\gamma,\mathcal{O}}\lesssim\int_{\mathbb{T}}\|(\partial_{\theta}^{k}\mathcal{K}_{3})(\ast,\cdot,\centerdot,\eta+\centerdot)\|_{q,s+s_{0}}^{\gamma,\mathcal{O}}d\eta\\
					&\lesssim \sum_{m\in\mathbb{S}_{0}}\Big(\|\mathcal{B}^{-1}\mathcal{S}_{2}^{\star}\partial_{\theta}\big(\mathcal{B}-\textnormal{Id}\big)g_{m}\|_{q,s+s_{0}}^{\gamma,\mathcal{O}}\|\mathscr{B}^{-1}e_{m}\|_{q,s_{0}+k}^{\gamma,\mathcal{O}}+\|\mathcal{B}^{-1}\mathcal{S}_{2}^{\star}\partial_{\theta}\big(\mathcal{B}-\textnormal{Id}\big)g_{m}\|_{q,s_{0}}^{\gamma,\mathcal{O}}\|\mathscr{B}^{-1}e_{m}\|_{q,s+s_{0}+k}^{\gamma,\mathcal{O}}\Big).
				\end{align*}
				Applying \eqref{estimate on the first reduction operator and its inverse} we find
				$$\|\mathcal{B}^{-1}\mathcal{S}_{2}^{\star}\partial_{\theta}\big(\mathcal{B}-\textnormal{Id}\big)g_{m}\|_{q,s}^{\gamma,\mathcal{O}}\lesssim \|\mathcal{S}_{2}^{\star}\partial_{\theta}\big(\mathcal{B}-\textnormal{Id}\big)g_{m}\|_{q,s}^{\gamma,\mathcal{O}}+\varepsilon\gamma^{-1}\|\mathfrak{I}_{0}\|_{q,s+\sigma}^{\gamma,\mathcal{O}}\|\mathcal{S}_{2}^{\star}\partial_{\theta}\big(\mathcal{B}-\textnormal{Id}\big)g_{m}\|_{q,s_{0}}^{\gamma,\mathcal{O}}.$$
				Now, from \eqref{def S2star}, the law products and Lemma \ref{properties of Toeplitz in time operators}, we find
				\begin{align*}
					\|\mathcal{S}_{2}^{\star}\rho\|_{q,s}^{\gamma,\mathcal{O}}\lesssim&\|V_{\varepsilon r}-c_{i_{0}}\|_{q,s_{0}}^{\gamma,\mathcal{O}}\|\rho\|_{q,s}^{\gamma,\mathcal{O}}+\|V_{\varepsilon r}-c_{i_{0}}\|_{q,s}^{\gamma,\mathcal{O}}\|\rho\|_{q,s_{0}}^{\gamma,\mathcal{O}}\\
					&+\|\mathbf{L}_{\varepsilon r,1}\|_{\textnormal{\tiny{O-d}},q,s_{0}}^{\gamma,\mathcal{O}}\|\rho\|_{q,s}^{\gamma,\mathcal{O}}+\|\mathbf{L}_{\varepsilon r,1}\|_{\textnormal{\tiny{O-d}},q,s}^{\gamma,\mathcal{O}}\|\rho\|_{q,s_{0}}^{\gamma,\mathcal{O}}.
				\end{align*}
			From the composition law and \eqref{estimate r1}, one has 
			\begin{align*}
				\|V_{\varepsilon r}-c_{i_{0}}\|_{q,s}^{\gamma,\mathcal{O}}&\leqslant \|V_{\varepsilon r}-V_{0}\|_{q,s}^{\gamma,\mathcal{O}}+\|V_{0}-c_{i_{0}}\|_{q}^{\gamma,\mathcal{O}}\\
				&\lesssim\varepsilon\left(1+\|\mathfrak{I}_{0}\|_{q,s+\sigma_{3}}^{\gamma,\mathcal{O}}\right).
			\end{align*}
		According to Lemma \ref{lemma symmetry and reversibility} and \eqref{estimate kernel mathbbK1}, we deduce that 
		\begin{align*}
			\|\mathbf{L}_{\varepsilon r,1}\|_{\textnormal{\tiny{O-d}},q,s}^{\gamma,\mathcal{O}}&\lesssim\int_{\mathbb{T}}\|\mathbb{K}_{\varepsilon r,1}(\ast,\cdot,\centerdot,\eta+\centerdot)\|_{q,s+s_{0}}^{\gamma,\mathcal{O}}d\eta\\
			&\lesssim\varepsilon\gamma^{-1}\left(1+\|\mathfrak{I}_{0}\|_{q,s+\sigma_{3}}^{\gamma,\mathcal{O}}\right).
		\end{align*}
	Using \eqref{small-3}, one gets
	$$\|\mathcal{S}_{2}^{\star}\rho\|_{q,s}^{\gamma,\mathcal{O}}\lesssim\|\rho\|_{q,s}^{\gamma,\mathcal{O}}+\varepsilon\gamma^{-1}\|\mathfrak{I}_{0}\|_{q,s+\sigma_{3}}^{\gamma,\mathcal{O}}\|\rho\|_{q,s_{0}}^{\gamma,\mathcal{O}}.$$
			Combining this with \eqref{estimate on the first reduction operator and its inverse} allow to get
				\begin{align*}
					\|\mathcal{S}_{2}^{\star}\partial_{\theta}\big(\mathcal{B}-\textnormal{Id}\big)g_{m}\|_{q,s}^{\gamma,\mathcal{O}}&\lesssim\|g_{m}\|_{q,s+1}^{\gamma,\mathcal{O}}+\varepsilon\gamma^{-1}\|\mathfrak{I}_{0}\|_{q,s+\sigma_{3}}^{\gamma,\mathcal{O}}\|g_{m}\|_{q,s_{0}}^{\gamma,\mathcal{O}}\\
					&\lesssim 1+\varepsilon\gamma^{-1}\|\mathfrak{I}_{0}\|_{q,s+\sigma_{3}}^{\gamma,\mathcal{O}}.
				\end{align*}
				Therefore,
			\begin{equation}\label{Z05}
				\max_{k\in\{0,1\}}\big\|\partial_{\theta}^{k}\mathcal{T}_0\partial_{\theta}\mathcal{S}_2\mathscr{B}\big\|_{\textnormal{\tiny{O-d}},q,s}^{\gamma,\mathcal{O}}\lesssim \varepsilon\gamma^{-1}\left(1+\|\mathfrak{I}_{0}\|_{q,s+\sigma_{3}}^{\gamma,\mathcal{O}}\right).
			\end{equation}
				Plugging the estimates \eqref{Z04} and \eqref{Z05} into \eqref{WL1} we find
				\begin{align}\label{F03}
					\nonumber \max_{k\in\{0,1\}}\big\|\partial_\theta^k\mathscr{B}_{\perp}^{-1}\partial_{\theta}\mathcal{S}_2\mathscr{B}\Pi_{\mathbb{S}_{0}}^{\perp}\big\|_{\textnormal{\tiny{O-d}},q,s}^{\gamma,\mathcal{O}}&\lesssim \max_{k\in\{0,1\}}\big\|\partial_\theta^{k+1}\mathcal{B}^{-1}\mathcal{S}_2\mathscr{B}\big\|_{\textnormal{\tiny{O-d}},q,s}^{\gamma,\mathcal{O}}+\max_{k\in\{0,1\}}\big\|\partial_\theta^k\mathcal{T}_0\partial_{\theta}\mathcal{S}_2\mathscr{B}\big\|_{\textnormal{\tiny{O-d}},q,s}^{\gamma,\mathcal{O}}\\
					&\lesssim \varepsilon\gamma^{-1}\left(1+\|\mathfrak{I}_{0}\|_{q,s+\sigma_{3}}^{\gamma,\mathcal{O}}\right).
				\end{align}
			We now turn to the estimate of the difference. Coming back to \eqref{mod-T-V}, one can write
			$$\Delta_{12}\big(\mathscr{B}_{\perp}^{-1}\partial_{\theta}\mathcal{S}_2\mathscr{B}\Pi_{\mathbb{S}_{0}}^{\perp}\big)=\Delta_{12}\big(\mathscr{B}^{-1}\partial_{\theta}\mathcal{S}_2\mathscr{B}\Pi_{\mathbb{S}_{0}}^{\perp}\big)-\Delta_{12}\big(\mathcal{T}_0\partial_{\theta}\mathcal{S}_2\mathscr{B}\Pi_{\mathbb{S}_{0}}^{\perp}\big).$$
		It follows that 
		\begin{equation}\label{WL2}
			\big\|\Delta_{12}\big(\mathscr{B}_{\perp}^{-1}\partial_{\theta}\mathcal{S}_2\mathscr{B}\Pi_{\mathbb{S}_{0}}^{\perp}\big)\big\|_{\textnormal{\tiny{O-d}},q,\overline{s}_h+\mathtt{p}}^{\gamma,\mathcal{O}}\lesssim \big\|\Delta_{12}\big(\partial_\theta\mathcal{B}^{-1}\mathcal{S}_2\mathscr{B}\big)\big\|_{\textnormal{\tiny{O-d}},q,\overline{s}_h+\mathtt{p}}^{\gamma,\mathcal{O}}+\big\|\Delta_{12}\big(\mathcal{T}_0\partial_{\theta}\mathcal{S}_2\mathscr{B}\big)\big\|_{\textnormal{\tiny{O-d}},q,\overline{s}_h+\mathtt{p}}^{\gamma,\mathcal{O}}.
		\end{equation}
	Arguing as for \eqref{GL-3}, one obtains
	$$\|\Delta_{12}(\partial_{\theta}\widehat{\mathcal{K}}_{2})(\ast,\cdot,\centerdot,\eta+\centerdot)\|_{q,\overline{s}_h+\mathtt{p}+s_0}\lesssim\varepsilon\gamma^{-1}\|\Delta_{12}i\|_{q,\overline{s}_h+\mathtt{p}+\sigma_{3}}^{\gamma,\mathcal{O}}\left(1-\log\left|\sin\left(\tfrac{\eta}{2}\right)\right|\right).
	$$
	Then, using Lemma \ref{lemma symmetry and reversibility} implies
	\begin{align}\label{d-BS2B}
		\big\|\Delta_{12}\big(\partial_\theta\mathcal{B}^{-1}\mathcal{S}_2\mathscr{B}\big)\big\|_{\textnormal{\tiny{O-d}},q,\overline{s}_h+\mathtt{p}}^{\gamma,\mathcal{O}}&\lesssim\int_{\mathbb{T}}\|\Delta_{12}(\partial_{\theta}\widehat{\mathcal{K}}_{2})(\ast,\cdot,\centerdot,\eta+\centerdot)\|_{q,\overline{s}_h+\mathtt{p}+s_0}d\eta\nonumber\\
		&\lesssim\varepsilon\gamma^{-1}\|\Delta_{12}i\|_{q,\overline{s}_h+\mathtt{p}+\sigma_{3}}^{\gamma,\mathcal{O}}.
	\end{align}
On the other hand, proceeding as for \eqref{D-06}, and using in particular \eqref{difference ci}, 
\begin{equation}\label{d-T0S2}
	\big\|\Delta_{12}\big(\mathcal{T}_0\partial_{\theta}\mathcal{S}_2\mathscr{B}\big)\big\|_{\textnormal{\tiny{O-d}},q,\overline{s}_h+\mathtt{p}}^{\gamma,\mathcal{O}}\lesssim\varepsilon\gamma^{-1}\|\Delta_{12}i\|_{q,\overline{s}_h+\mathtt{p}+\sigma_{3}}^{\gamma,\mathcal{O}}.
\end{equation}
Putting together \eqref{d-BS2B}, \eqref{d-T0S2} and \eqref{WL2}, ensures that
\begin{equation}\label{F03-1}
	\big\|\Delta_{12}\big(\mathscr{B}_{\perp}^{-1}\partial_{\theta}\mathcal{S}_2\mathscr{B}\Pi_{\mathbb{S}_{0}}^{\perp}\big)\big\|_{\textnormal{\tiny{O-d}},q,\overline{s}_h+\mathtt{p}}^{\gamma,\mathcal{O}}\lesssim\varepsilon\gamma^{-1}\|\Delta_{12}i\|_{q,\overline{s}_h+\mathtt{p}+\sigma_{3}}^{\gamma,\mathcal{O}}.
\end{equation}
				$\blacktriangleright$ {\it  Study of the term} $\varepsilon\mathscr{B}_{\perp}^{-1}\partial_{\theta}\mathcal{R}\mathscr{B}_{\perp}.$ Using the relation $\textnormal{Id}=\Pi_{\mathbb{S}_{0}}+\Pi_{\mathbb{S}_{0}}^{\perp},$ we can write
			\begin{align}\label{dec-R}
				\partial_{\theta}^{k}\mathscr{B}_{\perp}^{-1}&\partial_{\theta}\mathcal{R}\mathscr{B}_{\perp}=\partial_{\theta}^{k+1}\mathcal{B}^{-1}\mathcal{R}\mathscr{B}_{\perp}-\partial_{\theta}^{k}\mathcal{T}_{0}\partial_{\theta}\mathcal{R}\mathscr{B}_{\perp}\nonumber\\
				&=\partial_{\theta}^{k+1}\mathcal{B}^{-1}\mathcal{R}\Pi_{\mathbb{S}_{0}}^{\perp}\mathscr{B}\Pi_{\mathbb{S}_{0}}^{\perp}-\partial_{\theta}^{k}\mathcal{T}_{0}\partial_{\theta}\mathcal{R}\Pi_{\mathbb{S}_{0}}^{\perp}\mathscr{B}\Pi_{\mathbb{S}_{0}}^{\perp}\nonumber\\
				&=\partial_{\theta}^{k+1}\mathcal{B}^{-1}\mathcal{R}\mathscr{B}\Pi_{\mathbb{S}_{0}}^{\perp}-\partial_{\theta}^{k+1}\mathcal{B}^{-1}\mathcal{R}\Pi_{\mathbb{S}_{0}}\mathscr{B}\Pi_{\mathbb{S}_{0}}^{\perp}-\partial_{\theta}^{k}\mathcal{T}_{0}\partial_{\theta}\mathcal{R}\mathscr{B}\Pi_{\mathbb{S}_{0}}^{\perp}+\partial_{\theta}^{k}\mathcal{T}_{0}\partial_{\theta}\mathcal{R}\Pi_{\mathbb{S}_{0}}\mathscr{B}\Pi_{\mathbb{S}_{0}}^{\perp}.
			\end{align}
		Hence
		\begin{align}\label{dec-BRB}
			\|\partial_{\theta}^{k}\mathscr{B}_{\perp}^{-1}\partial_{\theta}\mathcal{R}\mathscr{B}_{\perp}\|_{\textnormal{\tiny{O-d}},q,s}^{\gamma,\mathcal{O}}\leqslant&\|\partial_{\theta}^{k+1}\mathcal{B}^{-1}\mathcal{R}\mathscr{B}\|_{\textnormal{\tiny{O-d}},q,s}^{\gamma,\mathcal{O}}+\|\partial_{\theta}^{k+1}\mathcal{B}^{-1}\mathcal{R}\Pi_{\mathbb{S}_{0}}\mathscr{B}\|_{\textnormal{\tiny{O-d}},q,s}^{\gamma,\mathcal{O}}\nonumber\\
			&+\|\partial_{\theta}^{k}\mathcal{T}_{0}\partial_{\theta}\mathcal{R}\mathscr{B}\|_{\textnormal{\tiny{O-d}},q,s}^{\gamma,\mathcal{O}}+\|\partial_{\theta}^{k}\mathcal{T}_{0}\partial_{\theta}\mathcal{R}\Pi_{\mathbb{S}_{0}}\mathscr{B}\|_{\textnormal{\tiny{O-d}},q,s}^{\gamma,\mathcal{O}}.
		\end{align}
	Recall that from Proposition \ref{lemma setting for Lomega} that $\mathcal{R}$ is an integral operator of kernel $J$ and therefore direct computations give
	\begin{equation}\label{conj-R}
		\big(\mathcal{B}^{-1}\mathcal{R}\mathscr{B}\rho\big)(\varphi,\theta)=\int_{\mathbb{T}}\rho(\varphi,\eta)\widehat{J}(\varphi,\theta,\eta)d\eta,
	\end{equation}
	with
	\begin{equation}\label{Jhat}
		\widehat{J}(\varphi,\theta,\eta):=J\big(\varphi,\theta+\widehat{\beta}(\varphi,\theta),\eta+\widehat{\beta}(\varphi,\eta)\big).
	\end{equation}
Combining \eqref{symmetry for beta} and \eqref{symmetry for the kernel J}, one gets
\begin{equation}\label{symmetry Jhat}
	\widehat{J}(-\varphi,-\theta,-\eta)=\widehat{J}(\varphi,\theta,\eta)\in\mathbb{R}.
\end{equation}
	Using the composition law and \eqref{estimate J}, we obtain
	$$\max_{k\in\{0,1,2\}}\sup_{\eta\in\mathbb{T}}\|(\partial_{\theta}^{k}\widehat{J})\big(\ast,\cdot,\centerdot,\eta+\centerdot\big)\|_{q,s}^{\gamma,\mathcal{O}}\lesssim 1+\|\mathfrak{I}_{0}\|_{q,s+\sigma_{3}}^{\gamma,\mathcal{O}}.$$
	Thus, applying  Lemma \ref{lemma symmetry and reversibility}-(ii) implies
	\begin{align}\label{e-BRB}
		\max_{k\in\{0,1\}}\|\partial_{\theta}^{k+1}\mathcal{B}^{-1}\mathcal{R}\mathscr{B}\|_{\textnormal{\tiny{O-d}},q,s}^{\gamma,\mathcal{O}}&\lesssim\max_{k\in\{0,1,2\}}\int_{\mathbb{T}}\|(\partial_{\theta}^{k}\widehat{J})\big(\ast,\cdot,\centerdot,\eta+\centerdot\big)\|_{q,s+s_{0}}^{\gamma,\mathcal{O}}d\eta\nonumber\\
		&\lesssim 1+\|\mathfrak{I}_{0}\|_{q,s+\sigma_{3}}^{\gamma,\mathcal{O}}.
	\end{align}
On the other hand we notice  from \eqref{conj-R} that we get the structure
\begin{equation}\label{conj-R-pi}
	\big(\mathcal{B}^{-1}\mathcal{R}\Pi_{\mathbb{S}_{0}}\mathscr{B}\rho\big)(\varphi,\theta)=\int_{\mathbb{T}}\rho(\varphi,\eta)\widetilde{J}(\varphi,\theta,\eta)d\eta,
\end{equation}
with
\begin{equation}\label{Jtilde}
	\widetilde{J}(\varphi,\theta,\eta):=\int_{\mathbb{T}}J\big(\varphi,\theta+\widehat{\beta}(\varphi,\theta),\eta'\big)D_{\mathbb{S}_{0}}(\eta'-\eta)d\eta'.
\end{equation}
Combining \eqref{symmetry for beta}, \eqref{symmetry for the kernel J} and the change of variables $\eta'\mapsto-\eta'$, one finds
	\begin{equation}\label{symmetry Jtilde}
		\widetilde{J}(-\varphi,-\theta,-\eta)=\widetilde{J}(\varphi,\theta,\eta)\in\mathbb{R}.
\end{equation}
Using the change of variables $\eta'\mapsto\eta'+\theta$ yields
$$\widetilde{J}(\varphi,\theta,\eta+\theta):=\int_{\mathbb{T}}J\big(\varphi,\theta+\widehat{\beta}(\varphi,\theta),\eta'+\theta\big)D_{\mathbb{S}_{0}}(\eta'-\eta)d\eta'.$$
Then by the composition law, we infer
$$\max_{k\in\{0,1,2\}}\sup_{\eta\in\mathbb{T}}\|(\partial_{\theta}^{k}\widetilde{J})(\ast,\cdot,\centerdot,\eta+\centerdot)\|_{q,s}^{\gamma,\mathcal{O}}\lesssim 1+\|\mathfrak{I}_{0}\|_{q,s+\sigma_{3}}^{\gamma,\mathcal{O}}.$$
Consequently, we find in view of Lemma \ref{lemma symmetry and reversibility}  
\begin{align}\label{e-BRpiB}
	\max_{k\in\{0,1\}}\|\partial_{\theta}^{k}\mathcal{B}^{-1}\mathcal{R}\Pi_{\mathbb{S}_{0}}\mathscr{B}\|_{\textnormal{\tiny{O-d}},q,s}^{\gamma,\mathcal{O}}&\lesssim\max_{k\in\{0,1,2\}}\int_{\mathbb{T}}\|(\partial_{\theta}^{k}\widetilde{J})\big(\ast,\cdot,\centerdot,\eta+\centerdot\big)\|_{q,s+s_{0}}^{\gamma,\mathcal{O}}d\eta\nonumber\\
	&\lesssim 1+\|\mathfrak{I}_{0}\|_{q,s+\sigma_{3}}^{\gamma,\mathcal{O}}.
\end{align}
If wet set
$$\mathcal{S}_{3}=\partial_{\theta}\mathcal{R}\mathscr{B}\quad \textnormal{or}\quad\partial_{\theta}\mathcal{R}\Pi_{\mathbb{S}_{0}}\mathscr{B},$$
then using \eqref{def op T0}, we deduce that 
\begin{align*}
	\mathcal{T}_{0}\mathcal{S}_{3}\rho&=\sum_{m\in\mathbb{S}_{0}}\big\langle\mathcal{S}_{3}\rho,\big(\mathcal{B}-\textnormal{Id}\big)g_{m}\big\rangle_{L_{\theta}^{2}(\mathbb{T})}\mathscr{B}^{-1}e_{m}\\
	&=\sum_{m\in\mathbb{S}_{0}}\big\langle\rho,\mathcal{S}_{3}^{\star}\big(\mathcal{B}-\textnormal{Id}\big)g_{m}\big\rangle_{L_{\theta}^{2}(\mathbb{T})}\mathscr{B}^{-1}e_{m},
\end{align*}
with  $\mathcal{S}_{3}^{\star}$ is the adjoint of $\mathcal{S}_{3}$ given by
\begin{equation}\label{S3star}
	\mathcal{S}_{3}^{\star}=-\mathcal{B}^{-1}\mathcal{R}^{\star}\partial_{\theta}\quad\textnormal{or}\quad-\mathcal{B}^{-1}\Pi_{\mathbb{S}_{0}}\mathcal{R}^{\star}\partial_{\theta},
	\end{equation}
and  $\mathcal{R}^{\star}$ the adjoint of $\mathcal{R}$ which is an integral operator with kernel
\begin{equation}\label{Jstar}
	J^{\star}(\varphi,\theta,\eta):=\sum_{k'=1}^{3}\sum_{k=1}^{d}g_{k,k'}(\varphi,\theta)\chi_{k,k'}(\varphi,\eta),
\end{equation}
where we use the notations of  the proof of Proposition \ref{lemma setting for Lomega}. Notice that similarly to \eqref{estimate J} and \eqref{symmetry for the kernel J}, the kernel $J^{\star}$ satisfies
\begin{equation}\label{estimate Jstar}
	\max_{k\in\{0,1,2\}}\sup_{\eta\in\mathbb{T}}\|(\partial_{\theta}^{k}J^{\star})(\ast,\cdot,\centerdot,\eta+\centerdot)\|_{q,s}^{\gamma,\mathcal{O}}\lesssim 1+\|\mathfrak{I}_{0}\|_{q,s+\sigma_{3}}^{\gamma,\mathcal{O}}
\end{equation}
and
\begin{equation}\label{symmetry Jstar}
	J^\star(-\varphi,-\theta,-\eta)=J^\star(\varphi,\theta,\eta)\in\mathbb{R}.
\end{equation}
Now, we have the integral representation
\begin{align*}
	\big(\mathcal{T}_{0}\mathcal{S}_{3}\rho\big)(\varphi,\theta)&=\int_{\mathbb{T}}\mathcal{K}_{4}(\varphi,\theta,\eta)\rho(\varphi,\eta)d\eta,\\
	\mathcal{K}_{4}(\varphi,\theta,\eta)&:=\sum_{m\in\mathbb{S}_{0}}\big(\mathcal{S}_{3}^{\star}\big(\mathcal{B}-\textnormal{Id}\big)g_{m}\big)(\varphi,\eta)\big(\mathscr{B}^{-1}e_{m}\big)(\varphi,\theta).
\end{align*}
Then by virtue of  \eqref{sym gm}, \eqref{symmetry for beta}, \eqref{def S2star} and \eqref{symmetry Jstar} we obtain
	\begin{equation}\label{sym K4}
		\mathcal{K}_4(-\varphi,-\theta,-\eta)=-\mathcal{K}_{4}(\varphi,\theta,\eta)\in\mathbb{R}.
\end{equation}
Applying Lemma \ref{lemma symmetry and reversibility} combined with the law products, we get for all $k\in\{0,1\}$
\begin{align*}
	\|\partial_{\theta}^{k}\mathcal{T}_{0}\mathcal{S}_{3}\|_{\textnormal{\tiny{O-d}},q,s}^{\gamma,\mathcal{O}}&\lesssim\int_{\mathbb{T}}\|(\partial_{\theta}^{k}\mathcal{K}_{4})(\ast,\cdot,\centerdot,\eta+\centerdot)\|_{q,s+s_{0}}^{\gamma,\mathcal{O}}d\eta\\
	&\lesssim\sum_{m\in\mathbb{S}_{0}}\Big(\|\mathcal{S}_{3}^{\star}\big(\mathcal{B}-\textnormal{Id}\big)g_{m}\|_{q,s+s_{0}}^{\gamma,\mathcal{O}}\|\mathscr{B}^{-1}e_{m}\|_{q,s_{0}+k}^{\gamma,\mathcal{O}}+\|\mathcal{S}_{3}^{\star}\big(\mathcal{B}-\textnormal{Id}\big)g_{m}\|_{q,s_{0}}^{\gamma,\mathcal{O}}\|\mathscr{B}^{-1}e_{m}\|_{q,s+s_{0}+k}^{\gamma,\mathcal{O}}\Big).
\end{align*}
Consequently, using  Lemma \ref{lemma symmetry and reversibility} and \eqref{estimate Jstar}, we get
\begin{align*}
	\|\mathcal{R}^{\star}\|_{\textnormal{\tiny{O-d}},q,s}^{\gamma,\mathcal{O}}&\lesssim\int_{\mathbb{T}}\|J^{\star}(\ast,\cdot,\centerdot,\eta+\centerdot)\|_{q,s+s_{0}}^{\gamma,\mathcal{O}}\\
	&\lesssim1+\|\mathfrak{I}_{0}\|_{q,s+\sigma_{3}}^{\gamma,\mathcal{O}}.
\end{align*}
Applying  \eqref{estimate on the first reduction operator and its inverse}, Lemma \ref{properties of Toeplitz in time operators} and the previous estimate implies
\begin{align*}
	\|\mathcal{S}_{3}^{\star}\rho\|_{q,s}^{\gamma,\mathcal{O}}&\lesssim\|\mathcal{R}^{\star}\partial_{\theta}\rho\|_{q,s}^{\gamma,\mathcal{O}}+\varepsilon\gamma^{-1}\|\mathfrak{I}_{0}\|_{q,s+\sigma_{3}}^{\gamma,\mathcal{O}}\|\mathcal{R}^{\star}\partial_{\theta}\rho\|_{q,s_{0}}^{\gamma,\mathcal{O}}\\
	&\lesssim\left(\varepsilon\gamma^{-1}\left(1+\|\mathfrak{I}_{0}\|_{q,s+\sigma_{3}}^{\gamma,\mathcal{O}}\right)+\|\mathcal{R}^{\star}\|_{\textnormal{\tiny{O-d}},q,s}^{\gamma,\mathcal{O}}\right)\|\rho\|_{q,s_{0}+1}^{\gamma,\mathcal{O}}+\|\mathcal{R}^{\star}\|_{\textnormal{\tiny{O-d}},q,s_{0}}^{\gamma,\mathcal{O}}\|\rho\|_{q,s+1}^{\gamma,\mathcal{O}}\\
	&\lesssim\|\rho\|_{q,s+1}^{\gamma,\mathcal{O}}+\|\mathfrak{I}_{0}\|_{q,s+\sigma_{3}}^{\gamma,\mathcal{O}}\|\rho\|_{q,s_{0}+1}^{\gamma,\mathcal{O}}.
\end{align*}
Thus
\begin{align*}
	\|\mathcal{S}_{3}^{\star}\big(\mathcal{B}-\textnormal{Id}\big)g_{m}\|_{q,s}^{\gamma,\mathcal{O}}&\lesssim \|g_{m}\|_{q,s+1}^{\gamma,\mathcal{O}}+\|\mathfrak{I}_{0}\|_{q,s+\sigma_{3}}^{\gamma,\mathcal{O}}\|g_{m}\|_{q,s_{0}+1}^{\gamma,\mathcal{O}}\\
	&\lesssim 1+\|\mathfrak{I}_{0}\|_{q,s+\sigma_{3}}^{\gamma,\mathcal{O}}.
\end{align*}
Hence
\begin{equation}\label{e-TOS3}
	\max_{k\in\{0,1\}}\|\partial_{\theta}^{k}\mathcal{T}_{0}\mathcal{S}_{3}\|_{\textnormal{\tiny{O-d}},q,s}^{\gamma,\mathcal{O}}\lesssim 1+\|\mathfrak{I}_{0}\|_{q,s+\sigma_{3}}^{\gamma,\mathcal{O}}.
\end{equation}
Putting together \eqref{dec-BRB}, \eqref{e-BRB}, \eqref{e-BRpiB} and \eqref{e-TOS3} allows to get
\begin{align}\label{F04}
	\max_{k\in\{0,1\}}\varepsilon\big\|\partial_{\theta}^k\mathscr{B}_{\perp}^{-1}\partial_{\theta}\mathcal{R}\mathscr{B}_{\perp}\big\|_{\textnormal{\tiny{O-d}},q,s}^{\gamma,\mathcal{O}}
	&\lesssim \varepsilon\gamma^{-1}\left(1+\|\mathfrak{I}_{0}\|_{q,s+\sigma_{3}}^{\gamma,\mathcal{O}}\right).
\end{align}
We now  move to the estimate of the difference. From \eqref{dec-R}, one has
\begin{align}\label{WL4}
	\|\Delta_{12}\big(\mathscr{B}_{\perp}^{-1}\partial_{\theta}\mathcal{R}\mathscr{B}_{\perp}\big)&\|_{\textnormal{\tiny{O-d}},q,\overline{s}_h+\mathtt{p}}^{\gamma,\mathcal{O}}\leqslant\|\Delta_{12}\big(\partial_{\theta}\mathcal{B}^{-1}\mathcal{R}\mathscr{B}\big)\|_{\textnormal{\tiny{O-d}},q,\overline{s}_h+\mathtt{p}}^{\gamma,\mathcal{O}}+\|\Delta_{12}\big(\partial_{\theta}\mathcal{B}^{-1}\mathcal{R}\Pi_{\mathbb{S}_{0}}\mathscr{B}\big)\|_{\textnormal{\tiny{O-d}},q,\overline{s}_h+\mathtt{p}}^{\gamma,\mathcal{O}}\nonumber\\
	&+\|\Delta_{12}\big(\mathcal{T}_{0}\partial_{\theta}\mathcal{R}\mathscr{B}\big)\|_{\textnormal{\tiny{O-d}},q,\overline{s}_h+\mathtt{p}}^{\gamma,\mathcal{O}}+\|\Delta_{12}\big(\mathcal{T}_{0}\partial_{\theta}\mathcal{R}\Pi_{\mathbb{S}_{0}}\mathscr{B}\big)\|_{\textnormal{\tiny{O-d}},q,\overline{s}_h+\mathtt{p}}^{\gamma,\mathcal{O}}.
\end{align}
Combining Lemma \ref{lemma symmetry and reversibility} with Taylor Formula, \eqref{conj-R}, \eqref{Jhat}, \eqref{differences J} and \eqref{difference beta} one obtains 
\begin{align}\label{d-BRB}
	\|\Delta_{12}\big(\partial_{\theta}\mathcal{B}^{-1}\mathcal{R}\mathscr{B}\big)\|_{\textnormal{\tiny{O-d}},q,\overline{s}_h+\mathtt{p}}^{\gamma,\mathcal{O}}&\lesssim\int_{\mathbb{T}}\|\Delta_{12}(\partial_{\theta}\widehat{J})(\ast,\cdot,\centerdot,\eta+\centerdot)\|_{q,\overline{s}_h+\mathtt{p}+s_0}^{\gamma,\mathcal{O}}d\eta\nonumber\\
	&\lesssim\|\Delta_{12}i\|_{q,\overline{s}_h+\mathtt{p}+\sigma_{3}}^{\gamma,\mathcal{O}}.
\end{align}
In the same spirit,  \eqref{conj-R-pi} and \eqref{Jtilde} give
\begin{align}\label{d-BRpiB}
	\|\Delta_{12}\big(\partial_{\theta}\mathcal{B}^{-1}\mathcal{R}\Pi_{\mathbb{S}_0}\mathscr{B}\big)\|_{\textnormal{\tiny{O-d}},q,\overline{s}_h+\mathtt{p}}^{\gamma,\mathcal{O}}&\lesssim\int_{\mathbb{T}}\|\Delta_{12}(\partial_{\theta}\widetilde{J})(\ast,\cdot,\centerdot,\eta+\centerdot)\|_{q,\overline{s}_h+\mathtt{p}+s_0}^{\gamma,\mathcal{O}}d\eta\nonumber\\
	&\lesssim\|\Delta_{12}i\|_{q,\overline{s}_h+\mathtt{p}+\sigma_{3}}^{\gamma,\mathcal{O}}.
\end{align}
According to the structure of $J^{\star}$ detailed in \eqref{Jstar} one can check that $J^{\star}$ satisfies similar  estimates as  \eqref{differences J}. Then using \eqref{difference beta}, one finds in a similar way to \eqref{D-06}, 
\begin{align}\label{d-TOS3}
	\|\Delta_{12}(\mathcal{T}_0\mathcal{S}_3)\|_{\textnormal{\tiny{O-d}},q,\overline{s}_h+\mathtt{p}}^{\gamma,\mathcal{O}}&\lesssim\int_{\mathbb{T}}\|\Delta_{12}\mathcal{K}_4(\ast,\cdot,\centerdot,\eta+\centerdot)\|_{q,\overline{s}_h+\mathtt{p}+s_0}^{\gamma,\mathcal{O}}d\eta\nonumber\\
	&\lesssim \|\Delta_{12}i\|_{q,\overline{s}_h+\mathtt{p}+\sigma_{3}}^{\gamma,\mathcal{O}}.
\end{align}
Hence, putting together, \eqref{d-BRB}, \eqref{d-BRpiB}, \eqref{d-TOS3} and \eqref{WL4} gives
\begin{equation}\label{F04-1}
	\varepsilon\|\Delta_{12}\big(\mathscr{B}_{\perp}^{-1}\partial_{\theta}\mathcal{R}\mathscr{B}_{\perp}\big)\|_{\textnormal{\tiny{O-d}},q,\overline{s}_h+\mathtt{p}}^{\gamma,\mathcal{O}}\lesssim\varepsilon\gamma^{-1}\|\Delta_{12}i\|_{q,\overline{s}_h+\mathtt{p}+\sigma_{3}}.
\end{equation}
On the other hand, gathering  \eqref{Form1A}, \eqref{symmertry mathbbK}, \eqref{sym K1}, \eqref{symmetry mathcal K2hat}, \eqref{sym K3}, \eqref{sym K4}, \eqref{symmetry Jhat} and \eqref{symmetry Jtilde} together with Lemma \ref{lemma symmetry and reversibility}, we find that $\mathscr{R}_0$ is a real and reversible Toeplitz in time integral operator.\\
In addition, \eqref{Form1A}, \eqref{F01}, \eqref{F02}, \eqref{F03} and \eqref{F04} give \eqref{F00}.\\  
Furthermore, \eqref{Form1A}, \eqref{F01-1}, \eqref{F02-1}, \eqref{F03-1} and \eqref{F04-1} imply \eqref{F00-1}.\\
\textbf{(iv)} Using Lemma \ref{properties of Toeplitz in time operators} together with \eqref{estimate mathscr R in off diagonal norm}, \eqref{estimate kernel equilibrium}, \eqref{estimate r1} and \eqref{small-3}, one obtains for all $s\in[s_{0},S]$
\begin{align*}
	\|\mathscr{L}_{0}\rho\|_{q,s}^{\gamma,\mathcal{O}}&\leqslant\|\big(\omega\cdot\partial_{\varphi}+c_{i_{0}}\partial_{\theta}+\partial_{\theta}\mathcal{K}_{\lambda}\ast\cdot)\big)\rho\|_{q,s}^{\gamma,\mathcal{O}}+\|\mathscr{R}_{0}\rho\|_{q,s}^{\gamma,\mathcal{O}}\\
	&\lesssim \|\rho\|_{q,s+1}^{\gamma,\mathcal{O}}+\|\mathscr{R}_{0}\|_{\textnormal{\tiny{O-d}},q,s}^{\gamma,\mathcal{O}}\|\rho\|_{q,s_{0}}^{\gamma,\mathcal{O}}+\|\mathscr{R}_{0}\|_{\textnormal{\tiny{O-d}},q,s_{0}}^{\gamma,\mathcal{O}}\|\rho\|_{q,s}^{\gamma,\mathcal{O}}\\
	&\lesssim\|\rho\|_{q,s+1}^{\gamma,\mathcal{O}}+\varepsilon\gamma^{-1}\|\mathfrak{I}_{0}\|_{q,s+\sigma_{3}}^{\gamma,\mathcal{O}}\|\rho\|_{q,s_{0}}^{\gamma,\mathcal{O}}.
\end{align*}
This ends the proof of Proposition \ref{projection in the normal directions}.
			\end{proof}
			\subsubsection{KAM reduction of the remainder term}\label{KAM-red-Feb}
			The goal of this section is to conjugate $\mathscr{L}_{0}$ defined in Proposition \ref{projection in the normal directions} to a diagonal operator, up to a fast decaying small  remainder. This will be achieved through  a standard  KAM reducibility techniques in the spirit of Proposition \ref{reduction of the transport part} but well-adapted to the operators setting.  This will be implemented by taking advantage of  the exterior parameters which are restricted to  a suitable Cantor set that prevents the resonances in the second Melnikov assumption. Notice that one gets from this study  some estimates on the distribution of   the eigenvalues and their stability with respect to the torus parametrization. This is considered as the key step not only to  get an approximate inverse but also to achieve Nash-Moser scheme with a  final massive Cantor set.			  The main result of this section reads as follows.
			\begin{prop}\label{reduction of the remainder term}
				Let $(\gamma,q,d,\tau_{1},\tau_2,s_{0},s_l,\overline{s}_l,\overline{\mu}_{2},S)$ satisfy \eqref{initial parameter condition}, \eqref{setting tau1 and tau2} and \eqref{param}. For any $(\mu_2,s_h)$ satisfying  
				\begin{align}\label{Conv-T2}
					\mu_2\geqslant \overline{\mu}_2+2\tau_2q+2\tau_2\quad\textnormal{and}\quad  s_h\geqslant \frac{3}{2}\mu_{2}+\overline{s}_{l}+1,
				\end{align}
				there exist $\varepsilon_{0}\in(0,1)$ and $\sigma_{4}=\sigma_{4}(\tau_1,\tau_2,q,d)\geqslant\sigma_{3}$, with $\sigma_{3}$ defined in Proposition $\ref{projection in the normal directions},$ such that if 
				\begin{equation}\label{hypothesis KAM reduction of the remainder term}
					\varepsilon\gamma^{-2-q}N_{0}^{\mu_{2}}\leqslant \varepsilon_{0}\quad\textnormal{and}\quad\|\mathfrak{I}_{0}\|_{q,s_{h}+\sigma_{4}}^{\gamma,\mathcal{O}}\leqslant 1,
				\end{equation}
				then the following assertions hold true.
				\begin{enumerate}[label=(\roman*)]
					\item There exists a  family of invertible linear operator $\Phi_{\infty}:\mathcal{O}\to \mathcal{L}\big(H_{\perp}^{s}\big)$ satisfying the estimates
					\begin{equation}\label{estimate on Phiinfty and its inverse}
						\forall s\in[s_{0},S],\quad \mbox{ }\|\Phi_{\infty}^{\pm 1}\rho\|_{q,s}^{\gamma ,\mathcal{O}}\lesssim \|\rho\|_{q,s}^{\gamma,\mathcal{O}}+\varepsilon\gamma^{-2}\| \mathfrak{I}_{0}\|_{q,s+\sigma_{4}}^{\gamma,\mathcal{O}}\|\rho\|_{q,s_{0}}^{\gamma,\mathcal{O}}.
					\end{equation}
					There exists a diagonal operator $\mathscr{L}_\infty=\mathscr{L}_{\infty}(\lambda,\omega,i_{0})$ taking the form
					\begin{align*}\mathscr{L}_{\infty}&=\omega\cdot\partial_{\varphi}\Pi_{\mathbb{S}_0}^{\perp}+\mathscr{D}_{\infty}
					\end{align*}
					where $\mathscr{D}_{\infty}=\mathscr{D}_{\infty}(\lambda,\omega,i_{0})$ is a reversible Fourier multiplier operator given by,
					$$
					\forall (l,j)\in \mathbb{Z}^{d}\times\mathbb{S}_{0}^{c},\quad  \mathscr{D}_{\infty}\mathbf{e}_{l,j}=\ii\mu_{j}^{\infty}\,\mathbf{e}_{l,j},$$
					with
					\begin{equation}\label{estim mujinfty}
						\forall j\in\mathbb{S}_{0}^{c},\quad\mu_{j}^{\infty}(\lambda,\omega,i_{0})=\mu_{j}^{0}(\lambda,\omega,i_{0})+r_{j}^{\infty}(\lambda,\omega,i_{0}),\quad \|r_{j}^{\infty}\|_{q}^{\gamma,\mathcal{O}}\lesssim\varepsilon\gamma^{-1}
					\end{equation}
					and 
					\begin{align}\label{estimate rjinfty}
						\sup_{j\in\mathbb{S}_{0}^{c}}|j|\| r_{j}^{\infty}\|_{q}^{\gamma ,\mathcal{O}}\lesssim\varepsilon\gamma^{-1},
					\end{align}
					such that in the Cantor set
					\begin{align*}
						\mathscr{O}_{\infty,n}^{\gamma,\tau_1,\tau_{2}}(i_{0}):=\bigcap_{\underset{|l|\leqslant N_{n}}{(l,j,j_{0})\in\mathbb{Z}^{d}\times(\mathbb{S}_0^{c})^{2}}\atop(l,j)\neq(0,j_{0})}\Big\{&(\lambda,\omega)\in\mathcal{O}_{\infty,n}^{\gamma,\tau_{1}}(i_{0}),\big|\omega\cdot l+\mu_{j}^{\infty}(\lambda,\omega,i_{0})-\mu_{j_{0}}^{\infty}(\lambda,\omega,i_{0})\big|>\tfrac{2\gamma \langle j-j_{0}\rangle}{\langle l\rangle^{\tau_{2}}}\Big\}
					\end{align*}
					we have 
					\begin{align*}
						\Phi_{\infty}^{-1}\mathscr{L}_{0}\Phi_{\infty}&=\mathscr{L}_{\infty}+\mathtt{E}^2_n,
					\end{align*}
					and the linear operator  $\mathtt{E}_n^{2}$ satisfies the estimate
					%						\begin{equation}\label{Error-Est-1D}
					%							\forall s\in[s_{0},S],\quad \mbox{ }\|\mathtt{E}^2_n\rho\|_{q,s}^{\gamma ,\mathcal{O}}\lesssim  \|\rho\|_{q,s+1}^{\gamma,\mathcal{O}}+\varepsilon\gamma^{-2}\| \mathfrak{I}_{0}\|_{q,s+\sigma_{4}}^{\gamma,\mathcal{O}}\|\rho\|_{q,s_{0}+1}^{\gamma,\mathcal{O}}
					%						\end{equation}
					%						and
					\begin{equation}\label{Error-Est-2D}
						\|\mathtt{E}^2_n\rho\|_{q,s_0}^{\gamma ,\mathcal{O}}\lesssim \varepsilon\gamma^{-2}N_{0}^{{\mu}_{2}}N_{n+1}^{-\mu_{2}} \|\rho\|_{q,s_0+1}^{\gamma,\mathcal{O}}.
					\end{equation}
					Notice that the Cantor set $\mathcal{O}_{\infty,n}^{\gamma,\tau_{1}}(i_{0})$ was introduced in Proposition $\ref{reduction of the transport part},$ the operator $\mathscr{L}_{0}$ and the frequencies $\big(\mu_{j}^{0}(\lambda,\omega,i_{0})\big)_{j\in\mathbb{S}_0^c}$ were stated in Proposition  $\ref{projection in the normal directions}.$
					\item Given two tori $i_{1}$ and $i_{2}$ both satisfying \eqref{hypothesis KAM reduction of the remainder term}, then
					\begin{equation}\label{diffenrence rjinfty}
						\forall j\in\mathbb{S}_{0}^{c},\quad\|\Delta_{12}r_{j}^{\infty}\|_{q}^{\gamma,\mathcal{O}}\lesssim\varepsilon\gamma^{-1}\|\Delta_{12}i\|_{q,\overline{s}_{h}+\sigma_{4}}^{\gamma,\mathcal{O}}
					\end{equation}
					and
					\begin{equation}\label{estimate differences mujinfty}
						\forall j\in\mathbb{S}_{0}^{c},\quad\|\Delta_{12}\mu_{j}^{\infty}\|_{q}^{\gamma,\mathcal{O}}\lesssim\varepsilon\gamma^{-1}|j|\| \Delta_{12}i\|_{q,\overline{s}_{h}+\sigma_{4}}^{\gamma,\mathcal{O}}.
					\end{equation}
				\end{enumerate}
				
			\end{prop}
			%\noindent\textbf{Remarks :}
			%\begin{enumerate}[label=\textbullet]
			%\item The constant $2$ appearing in the Cantor set $\mathscr{O}_{\infty}^{\gamma}(i_{0})$ is here to ensure the inclusion of this Cantor set in all the Cantor sets built in the proof.
			%\item The constant $\tau_{2}>\tau_{1}$ (see \eqref{setting tau1 and tau2}) is here to ensure the convergence of the series in the estimates of the Cantor set obtained at the end of the Nash-Moser scheme.
			%\end{enumerate}
			\begin{proof}
				{\bf{(i)}} 
				We shall introduce the quantity
				$$
				\delta_{0}(s):=\gamma ^{-1}\|\mathscr{R}_{0}\|_{\textnormal{\tiny{O-d}},q,s}^{\gamma ,\mathcal{O}},
				$$
				where $\mathscr{R}_0$ is the remainder seen in  Proposition  \ref{projection in the normal directions}. By applying   \eqref{estimate mathscr R in off diagonal norm}, we deduce that
				\begin{align}\label{whab1}\delta_{0}(s)\leqslant C\varepsilon\gamma^{-2}\left(1+\|\mathfrak{I}_{0}\|_{q,s+\sigma_{3}}^{\gamma,\mathcal{O}}\right).
				\end{align}
				Therefore  with the notation of \eqref{Conv-T2},  \eqref{hypothesis KAM reduction of the remainder term} and the fact that $\sigma_4\geqslant \sigma_3$  we obtain 
				\begin{align}\label{Conv-P3}
					\nonumber N_{0}^{\mu_{2}}\delta_{0}(s_{h}) &\leqslant  C N_{0}^{\mu_{2}}\varepsilon\gamma^{-2}\\
					&\leqslant C\varepsilon_{0}.
				\end{align}
				$\blacktriangleright$ \textbf{KAM step.} Recall from  Proposition \ref{projection in the normal directions} that in the Cantor set $\mathcal{O}_{\infty,n}^{\gamma,\tau_1}(i_0)$ one has
				$$\mathscr{B}_{\perp}^{-1}\widehat{\mathcal{L}}_{\omega}\mathscr{B}_{\perp}=\mathscr{L}_{0}+\mathtt{E}_{n}^{1},$$
				where the operator $\mathscr{L}_{0}$ has the following structure 
				\begin{equation}\label{mouka1}
					\mathscr{L}_0=\big(\omega\cdot\partial_{\varphi}+\mathscr{D}_0\big)\Pi_{\mathbb{S}_0}^{\perp}+\mathscr{R}_0,
				\end{equation}
				with $\mathscr{D}_0$ a diagonal operator of pure imaginary spectrum and $\mathscr{R}_0$ a real and reversible Toeplitz in time operator of zero order satisfying $\Pi_{\mathbb{S}_0}^\perp \mathscr{R}_0\Pi_{\mathbb{S}_0}^\perp =\mathscr{R}_0.$   Similarly to the reduction of the transport part, we shall first expose a typical step of the iteration process of the KAM scheme whose goal is  to reduce to a diagonal part $\mathscr{R}_0$. Notice that the scheme is flexible and has been used in the literature to deal with various equations. Assume that we have a linear operator $\mathscr{L}$ taking the following form in restriction to some Cantor set  $\mathscr{O}$ one has 
				$$
				\mathscr{L}=\big(\omega\cdot\partial_{\varphi}+\mathscr{D}\big)\Pi_{\mathbb{S}_0}^{\perp}+\mathscr{R},
				$$
				where $\mathscr{D}$ is real and reversible diagonal Toeplitz in time  operator, that is,
				\begin{align}\label{TUma1}
					\mathscr{D}\mathbf{e}_{l,j}=i\mu_j(\lambda,\omega) \,\mathbf{e}_{l,j}\quad\hbox{and}\quad \mu_{-j}(\lambda,\omega)=-\mu_{j}(\lambda,\omega).
				\end{align}
				The operator  $\mathscr{R}$ is assumed to be a real and reversible Toeplitz in time operator  of zero order satisfying $\Pi_{\mathbb{S}_0}^\perp \mathscr{R}\Pi_{\mathbb{S}_0}^\perp =\mathscr{R}.$
				Consider   a linear invertible transformation close to the identity $$\Phi=\Pi_{\mathbb{S}_0}^{\perp}+\Psi:\mathcal{O}\rightarrow\mathcal{L}(H_{\perp}^{s}),$$
				where $\Psi$ is small and depends on $\mathscr{R}$.
				Then  straightforward calculus show that in $\mathscr{O}$
				\begin{align*}
					\Phi^{-1}\mathscr{L}\Phi & =  \Phi^{-1}\Big(\Phi\left(\omega\cdot\partial_{\varphi}+\mathscr{D}\right)\Pi_{\mathbb{S}_0}^{\perp}+\left[\omega\cdot\partial_{\varphi}\Pi_{\mathbb{S}_0}^{\perp}+\mathscr{D},\Psi\right]+\mathscr{R}+\mathscr{R}\Psi\Big)\\
					& =  \big(\omega\cdot\partial_{\varphi}+\mathscr{D}\big)\Pi_{\mathbb{S}_0}^{\perp}+\Phi^{-1}\Big(\big[\big(\omega\cdot\partial_{\varphi}+\mathscr{D}\big)\Pi_{\mathbb{S}_0}^{\perp},\Psi\big]+P_{N}\mathscr{R}+P_{N}^{\perp}\mathscr{R}+\mathscr{R}\Psi\Big),
				\end{align*}
				where the projector $P_{N}$ was defined in \eqref{definition of projections for operators}.
				The main idea consists in replacing the remainder $\mathscr{R}$ with another quadratic one  up to  a diagonal part and provided that the parameters $(\lambda,\omega)$ belongs to a Cantor set connected to non-resonance conditions associated to the {\it homological equation}. Iterating this scheme  will generate new remainders  which  become smaller and smaller  up to  new contributions on  the diagonal part and with  more extraction on the parameters. Then by passing   to the limit we expect to diagonalize completely the operators provided that the parameters belong to a limit Cantor set. Notice that the  Cantor set should be truncated in the time mode in order to get a stability form required later in Nash-Moser scheme and during the measure of the final Cantor set. This will induce   a diagonalization  up to small  fast decaying remainders modeled by the operators $\mathtt{E}^2_n$ in \mbox{Proposition \ref{reduction of the remainder term}.}      Now the first step is to impose the following  homological equation,
				\begin{equation}\label{equation Psi}
					\big[\big(\omega\cdot\partial_{\varphi}+\mathscr{D}\big)\Pi_{\mathbb{S}_0}^{\perp},\Psi\big]+P_{N}\mathscr{R}=\lfloor P_{N}\mathscr{R}\rfloor,
				\end{equation}
				where $\lfloor P_{N}\mathscr{R}\rfloor$ is the diagonal part of the operator $\mathscr{R}$. We emphasize that the notation   $\lfloor  \mathcal{R}\rfloor$ with a general  operator $ \mathcal{R}$ is defined as follows, {for all $(l_{0},j_{0})\in\mathbb{Z}^{d}\times\mathbb{S}_{0}^{c}$},
				\begin{equation}\label{Dp1X}
					\mathcal{R}{\bf e}_{l_0,j_0}=\sum_{(l,j)\in\mathbb{Z}^{d}\times\mathbb{S}_{0}^{c}}\mathcal{R}_{l_0,j_0}^{l,j}{\bf e}_{l,j}\Longrightarrow  \lfloor\mathcal{R}\rfloor {\bf e}_{l_0,j_0}=\mathcal{R}_{l_0,j_0}^{l_0,j_0}{\bf e}_{l_0,j_0}=\big \langle \mathcal{R}{\bf e}_{l_0,j_0}, {\bf e}_{l_0,j_0}\big\rangle_{L^2(\T^{d+1})}\,{\bf e}_{l_0,j_0}.
				\end{equation}
				Remind the notation  ${\bf e}_{l_0,j_0}(\varphi,\theta)=e^{\ii(l_0\cdot\varphi+j_0\theta)}$. 
				The Fourier coefficients  of $\Psi $ are defined through
				$$
				\Psi {\bf e}_{l_0,j_0}=\sum_{(l,j)\in\mathbb{Z}^{d}\times\mathbb{S}_{0}^{c}}\Psi_{l_0,j_0}^{l,j}{\bf e}_{l,j},\quad \Psi_{l_0,j_0}^{l,j}\in\mathbb{C}.
				$$
				From direct  computations based on the above Fourier decomposition, we infer
				$$
				\big[\omega\cdot\partial_{\varphi}\Pi_{\mathbb{S}_0}^{\perp},\Psi\big]\mathbf{e}_{l_{0},j_{0}}=\ii\sum_{(l,j)\in\mathbb{Z}^{d }\times\mathbb{S}_{0}^{c}}\Psi_{l_{0},j_{0}}^{l,j}\,\,\omega\cdot(l-l_{0})\,\,\mathbf{e}_{l,j}
				$$
				and using the diagonal structure of $\mathscr{D}$,
				$$
				[\mathscr{D}_{0},\Psi]\mathbf{e}_{l_{0},j_{0}}=\ii\sum_{(l,j)\in\mathbb{Z}^{d }\times\mathbb{S}_{0}^{c}}\Psi_{l_{0},j_{0}}^{l,j}\left(\mu_{j}(\lambda,\omega)-\mu_{j_0}(\lambda,\omega)\right)\mathbf{e}_{l,j}.
				$$
				By hypothesis, $\mathscr{R}$ is a real and reversible Toeplitz in time operator. Hence  its Fourier coefficients  write in view of  Proposition \ref{characterization of real operator by its Fourier coefficients},
				\begin{equation}\label{coefficients of the remainder operator R}
					\mathscr{R}_{l_{0},j_{0}}^{l,j}:=\ii\,r_{j_{0}}^{j}(\lambda,\omega,l_{0}-l)\in \ii\,\mathbb{R}\quad\mbox{ and }\quad \mathscr{R}_{-l_{0},-j_{0}}^{-l,-j}=-\mathscr{R}_{l_{0},j_{0}}^{l,j}.
				\end{equation}
				Consequently  $\Psi$ is a solution of \eqref{equation Psi} if and only if 
				$$
				\Psi {\bf e}_{l_0,j_0}=\sum_{|l-l_0|\leqslant N\atop |j-j_{0}|\leqslant N}\Psi_{l_0,j_0}^{l,j}{\bf e}_{l,j}
				$$
				and
				$$\Psi_{l_{0},j_{0}}^{l,j}\Big(\omega\cdot (l-l_{0})+\mu_{j}(\lambda,\omega)-\mu_{j_0}(\lambda,\omega)\Big)=\left\lbrace\begin{array}{ll}
					-r_{j_{0}}^{j}(\lambda,\omega,l_{0}-l) & \mbox{if }(l,j)\neq(l_{0},j_{0})\\
					0 & \mbox{if }(l,j)=(l_{0},j_{0}).
				\end{array}\right.$$
				In particular, we get  that $\Psi$ is a Toeplitz in time operator with $\Psi_{j_{0}}^{j}(l_{0}-l):=\Psi_{l_{0},j_{0}}^{l,j}$. Moreover, for  $(l,j,j_{0})\in\mathbb{Z}^{d }\times(\mathbb{S}_{0}^{c})^{2}$ with $|l|,|j-j_{0}|\leqslant N,$ one obtains
				\begin{equation}\label{Psi-gh}
					\Psi_{j_{0}}^{j}(\lambda,\omega,l)=\left\lbrace\begin{array}{ll}
						\frac{-r_{j_{0}}^{j}(\lambda,\omega,l)}{\omega\cdot l+\mu_{j}(\lambda,\omega)-\mu_{j_{0}}(\lambda,\omega)} & \mbox{if }(l,j)\neq(0,j_{0})\\
						0 & \mbox{if }(l,j)=(0,j_{0}),
					\end{array}\right.
				\end{equation}
				provided that the denominator is non zero. In addition, from $\Pi_{\mathbb{S}_0}^\perp \mathscr{R}\Pi_{\mathbb{S}_0}^\perp =\mathscr{R},$  one easily gets
				$$
				\forall \,l\,\in\mathbb{Z}^d,\,\forall\, j\,\,\hbox{or}\,\,\,j_0\in\mathbb{S}_0,\quad r_{j_{0}}^{j}(\lambda,\omega,l)=0.
				$$
				Therfore, we should impose the compatibility condition
				$$
				\forall \,l\,\in\mathbb{Z}^d,\,\forall\, j\,\,\hbox{or}\,\,\,j_0\in\mathbb{S}_0,\quad \Psi_{j_{0}}^{j}(\lambda,\omega,l)=0.
				$$
				This implies that $\Pi_{\mathbb{S}_0}^\perp \Psi\Pi_{\mathbb{S}_0}^\perp=\Psi.$
				To justify the formula  given by \eqref{Psi-gh} we need to avoid resonances and   restrict the parameters to the following open set according to the so-called  second Melnikov  condition,
				$$
				\mathscr{O}_{+}^{\gamma }=\bigcap_{\underset{|l|\leqslant N}{(l,j,j_{0})\in\mathbb{Z}^{d }\times({\mathbb{S}}_{0}^{c})^{2}}\atop (l,j)\neq(0,j_0)}\left\lbrace(\lambda,\omega)\in\mathscr{O}
				\quad\textnormal{s.t.}\quad |\omega\cdot l+\mu_{j}(\lambda,\omega)-\mu_{j_{0}}(\lambda,\omega)|>\tfrac{\gamma\langle j-j_{0}\rangle }{\langle l\rangle^{\tau_{2}}}\right\rbrace.
				$$ 
				%$$
				%\textcolor{red}{\mathscr{O}_{+}^{\gamma }=\bigcap_{\underset{|l||\leqslant N}{(l,j,j_{0})\in\mathbb{Z}^{d }\times({\mathbb{S}}_{0}^{c})^{2}}\atop (l,j)\neq(0,j_0)}\left\lbrace(\lambda,\omega)\in\mathcal{O}_{\infty,n}^{\gamma,\tau_1}(i_0)
				%, |\omega\cdot l+\mu_{j}(\lambda,\omega)-\mu_{j_{0}}(\lambda,\omega)|>\frac{\gamma|j-j_{0}|}{\langle l\rangle^{\tau_{2}}}\right\rbrace}
				%.
				%$$
				In view of this restriction, the identity  \eqref{Psi-gh} is well defined and to extend  $\Psi$ to the whole set $\mathcal{O}$ we shall   use the cut-off function $\chi$ of  \eqref{properties cut-off function first reduction}.  We set 
				\begin{equation}\label{Ext-psi-op}
					\Psi_{j_{0}}^{j}(\lambda,\omega,l)=\left\lbrace\begin{array}{ll}
						-\varrho_{j_{0}}^{j}(\lambda,\omega,l)\,\, r_{j_{0}}^{j}(\lambda,\omega,l),& \mbox{if }\quad (l,j)\neq(0,j_{0})\\
						0, & \mbox{if }\quad (l,j)=(0,j_{0}),
					\end{array}\right.
				\end{equation}
				with
				\begin{align}\label{varr-d}
					\varrho_{j_{0}}^{j}(\lambda,\omega,l):=\frac{\chi\left((\omega\cdot l+\mu_{j}^{0}(\lambda,\omega)-\mu_{j_{0}}^{0}(\lambda,\omega))(\gamma\langle j-j_{0}\rangle)^{-1}\langle l\rangle^{\tau_{2}}\right)}{\omega\cdot l+\mu_{j}(\lambda,\omega)-\mu_{j_{0}}(\lambda,\omega)}\cdot
				\end{align}
				To simplify the notation, in the sequel we shall still write $\Psi$ to denote this extension. Note that the extension \eqref{Ext-psi-op} is smooth and in restriction to the Cantor set $\mathscr{O}_{+}^{\gamma }$ coincides with $\Psi.$ 
				On the other hand,  \eqref{coefficients of the remainder operator R} and \eqref{Ext-psi-op} imply that $\Psi_{j_{0}}^{j}(l)\in\mathbb{R}.$ In addition,  \eqref{varr-d} combined with  \eqref{TUma1} give	$$\Psi_{-j_{0}}^{-j}(-l)=\Psi_{j_{0}}^{j}(l).$$
				Consequently, in view of  Proposition \ref{characterization of real operator by its Fourier coefficients}, we deduce that $\Psi$ is a real and reversibility preserving operator.
				Now consider,
				\begin{equation}\label{PU-RT}
					\mathscr{D}_{+}=\mathscr{D}+\lfloor P_{N}\mathscr{R}\rfloor,\quad \quad \mathscr{R}_{+}=\Phi^{-1}\big(-\Psi\,\lfloor P_{N}\mathscr{R}\rfloor +P_{N}^{\perp}\mathscr{R}+\mathscr{R}\Psi\big)
				\end{equation}
				and
				$$
				\mathscr{L}_{+}:=\big(\omega\cdot\partial_{\varphi}+\mathscr{D}_{+}+\mathscr{R}_{+}\big)\Pi_{\mathbb{S}_0}^{\perp}.$$
				Therefore, in restriction to the Cantor set $\mathscr{O}_{+}^{\gamma }$, we can write
				$$
				\mathscr{L}_{+}=\Phi^{-1}\mathscr{L}\Phi.
				$$
				Our next task is to  estimate $\varrho_{j_{0}}^{j}$  defined by \eqref{varr-d}. Notice that this quantity can be written in the following form 
				\begin{align}\label{varr-dX1}
					&\varrho_{j_{0}}^{j}(\lambda,\omega,l)= a_{l,j,j_{0}} \widehat \chi\big(a_{l,j,j_{0}} A_{l,j,j_{0}}(\lambda,\omega)\big), \quad \widehat{\chi}(x)=\frac{\chi(x)}{x},\nonumber\\
					&A_{l,j,j_{0}}(\lambda,\omega)=\omega\cdot l+\mu_{j}(\lambda,\omega)-\mu_{j_{0}}(\lambda,\omega), \quad a_{l,j,j_{0}}= (\gamma\langle j-j_{0}\rangle )^{-1}\langle l\rangle^{\tau_{2}},
				\end{align}
				where   $\widehat{\chi}(x):=\frac{\chi(x)}{x}$  is $\mathcal{C}^\infty$ with bounded derivatives.
				Assume now the following estimate 
				\begin{align}\label{reg-G-10}
					\forall\,(j,j_{0})\in(\mathbb{S}_0^c)^{2},\quad\max_{|\alpha|\in\llbracket 0, q\rrbracket}\sup_{(\lambda,\omega)\in\mathcal{O}}\left|\partial_{\lambda,\omega}^{\alpha}\left(\mu_{j}(\lambda,\omega)-\mu_{j_{0}}(\lambda,\omega)\right)\right|\leqslant C\,|j-j_{0}|.
				\end{align}
				%Now, according to  \mbox{Lemma \ref{lemma control of unperturbed linear frequencies},} 
				%\begin{align}\label{reg-G-1}\exists C>0,\,\, \forall(j,j_{0})\in\mathbb{Z}^{2},\quad\max_{0\leqslant\alpha \leqslant q}\sup_{\lambda\in(\lambda_{0},\lambda_{1})}\left|\partial_{\lambda}^{\alpha}\left(\Omega_{j}(\lambda)-\Omega_{j_{0}}(\lambda)\right)\right|\leqslant C\,|j-j_{0}|.
				%\end{align}
				%From the structure of $\mu_{j}^{0}$   in Proposition \ref{projection in the normal directions} combined with \eqref{differences mu0} and \eqref{reg-G-1}, \eqref{hypothesis KAM reduction of the remainder term} we find 
				%\begin{align}\label{reg-G-10}\exists C>0,\,\, \forall\,j,j_{0}\in\mathbb{S}_0^c,\quad\max_{0\leqslant|\alpha| \leqslant q}\sup_{\lambda\in(\lambda_{0},\lambda_{1})}\left|\partial_{\lambda,\omega}^{\alpha}\left(\mu_{j}^{0}(\lambda,\omega)-\mu_{j_{0}}^{0}(\lambda,\omega)\right)\right|\leqslant C\,|j-j_{0}|.
				%\end{align}
				Then, we find
				\begin{align}\label{reg-G-11}
					\forall\, (l,j,j_{0})\in\mathbb{Z}^{d}\times(\mathbb{S}_0^c)^2,\quad\max_{q'\in\llbracket 0, q\rrbracket}\|A_{l,j,j_{0}}\|_{q'}^{\gamma,\mathcal{O}}\leqslant C\,\langle l,j-j_{0}\rangle.
				\end{align}
				In a similar way to \eqref{gljtqprm}, using Lemma \ref{Lem-lawprod}-(vi) and \eqref{reg-G-11}, we obtain
				\begin{align}\label{eg-k-1}
					\forall q'\in\llbracket 0,q\rrbracket,\quad \|\varrho_{j_{0}}^{j}(\ast,l)\|_{q'}^{\gamma,\mathcal{O}}\leqslant C\gamma^{-(q'+1)} \langle l,j-j_{0}\rangle^{\tau_{2}q'+\tau_2+q'}.
				\end{align}
				Similarly  to \eqref{control of g by f}, using Leibniz rule, we get 
				\begin{align}\label{link Psi and R}
					\|\Psi\|_{\textnormal{\tiny{O-d}},q,s}^{\gamma ,\mathcal{O}}&
					\leqslant C \gamma ^{-1}\| P_{N}\mathscr{R}\|_{\textnormal{\tiny{O-d}},q,s+\tau_{2} q+\tau_{2}}^{\gamma ,\mathcal{O}}.
				\end{align}
				We also assume that the following smallness condition holds
				\begin{align}\label{assum-Zi1}
					\gamma ^{-1}\| \mathscr{R}\|_{\textnormal{\tiny{O-d}},q,s_0+\tau_{2} q+\tau_{2}}^{\gamma ,\mathcal{O}}\leqslant  C\varepsilon_{0}.
				\end{align}
				Hence, by virtue of  \eqref{link Psi and R}, we get 
				\begin{align}\label{ZZ-KL}
					\nonumber \|\Psi\|_{\textnormal{\tiny{O-d}},q,s_0}^{\gamma ,\mathcal{O}}&\leqslant C \gamma ^{-1}\| \mathscr{R}\|_{\textnormal{\tiny{O-d}},q,s_0+\tau_{2} q+\tau_{2}}^{\gamma ,\mathcal{O}}\\
					&\leqslant C\varepsilon_{0}.
				\end{align}
				As a consequence, up to take $\varepsilon_0$ small enough, the operator  $\Phi$ is invertible and 
					$$\Phi^{-1}=\displaystyle\sum_{n=0}^{\infty}(-1)^{n}\Psi^{n}:=\textnormal{Id}+\Sigma.
					$$
				According to the law products in Lemma \ref{Lem-lawprod}, Lemma \ref{properties of Toeplitz in time operators}, \eqref{link Psi and R} and \eqref{ZZ-KL} one gets
				\begin{align}\label{QWX0}
					\displaystyle\|\Sigma\|_{\textnormal{\tiny{O-d}},q,s}^{\gamma ,\mathcal{O}} & \leqslant  \displaystyle\|\Psi\|_{\textnormal{\tiny{O-d}},q,s}^{\gamma ,\mathcal{O}}\left(1+\sum_{n=1}^{\infty}\left(C\|\Psi\|_{\textnormal{\tiny{O-d}},q,s_{0}}^{\gamma ,\mathcal{O}}\right)^{n}\right)\nonumber\\
					&\leqslant \displaystyle C\,\gamma ^{-1}N^{\tau_{2} q+\tau_{2}}\|\mathscr{R}\|_{\textnormal{\tiny{O-d}},q,s}^{\gamma ,\mathcal{O}}.
				\end{align}
				Therefore, we conclude with the assumption  \eqref{assum-Zi1} that  $\Phi^{-1}$ is satisfies the following estimate
				\begin{equation}\label{QWX1}
					\|\Phi^{-1}-\hbox{Id}\|_{\textnormal{\tiny{O-d}},q,s}^{\gamma ,\mathcal{O}}  {\leqslant}  C\gamma ^{-1}N^{\tau_{2} q+\tau_{2}}\|\mathscr{R}\|_{\textnormal{\tiny{O-d}},q,s}^{\gamma ,\mathcal{O}}.
				\end{equation}
				From \eqref{PU-RT}, we can write
				$$
				\mathscr{R}_{+}=P_{N}^{\perp}\mathscr{R}+\Phi^{-1}\mathscr{R}\Psi-\Psi\,\lfloor P_{N}\mathscr{R}\rfloor +\Sigma\big(P_{N}^{\perp}\mathscr{R}-\Psi\,\lfloor P_{N}\mathscr{R}\rfloor \big).
				$$
				Thus, by virtue of Lemma  \ref{properties of Toeplitz in time operators} and \eqref{QWX1}, we infer
				\begin{align}\label{k-12}
					\nonumber \|\mathscr{R}_{+}\|_{\textnormal{\tiny{O-d}},q,s}^{\gamma ,\mathcal{O}} & \leqslant  \| P_{N}^{\perp}\mathscr{R}\|_{\textnormal{\tiny{O-d}},q,s}^{\gamma ,\mathcal{O}}+C\|\Sigma\|_{\textnormal{\tiny{O-d}},q,s}^{\gamma ,\mathcal{O}}\left(\| P_{N}^{\perp}\mathscr{R}\|_{\textnormal{\tiny{O-d}},q,s_{0}}^{\gamma ,\mathcal{O}}+\|\Psi\|_{\textnormal{\tiny{O-d}},q,s_{0}}^{\gamma ,\mathcal{O}}\|\mathscr{R}\|_{\textnormal{\tiny{O-d}},q,s_{0}}^{\gamma ,\mathcal{O}}\right)\\
					&+C\left(1+\|\Sigma\|_{\textnormal{\tiny{O-d}},q,s_{0}}^{\gamma ,\mathcal{O}}\right)
					\left(\|\Psi\|_{\textnormal{\tiny{O-d}},q,s}^{\gamma ,\mathcal{O}}\|\mathscr{R}\|_{\textnormal{\tiny{O-d}},q,s_{0}}^{\gamma ,\mathcal{O}}+\|\Psi\|_{\textnormal{\tiny{O-d}},q,s_{0}}^{\gamma ,\mathcal{O}}\|\mathscr{R}\|_{\textnormal{\tiny{O-d}},q,s}^{\gamma ,\mathcal{O}}\right).
				\end{align}
				By Lemma  \ref{properties of Toeplitz in time operators} , \eqref{link Psi and R},\eqref{ZZ-KL} and  \eqref{QWX1}, we get for all $S\geqslant \overline{s}\geqslant s\geqslant s_{0}$,
				\begin{equation}\label{KAM step remainder term}
					\|\mathscr{R}_{+}\|_{\textnormal{\tiny{O-d}},q,s}^{\gamma ,\mathcal{O}}\leqslant N^{s-\overline{s}}\|\mathscr{R}\|_{\textnormal{\tiny{O-d}},q,\overline{s}}^{\gamma ,\mathcal{O}}+C\gamma ^{-1}N^{\tau_{2} q+\tau_{2}}\|\mathscr{R}\|_{\textnormal{\tiny{O-d}},q,s_{0}}^{\gamma ,\mathcal{O}}\|\mathscr{R}\|_{\textnormal{\tiny{O-d}},q,s}^{\gamma ,\mathcal{O}}.
				\end{equation}
				%
				%$\blacktriangleright$ \textbf{Initialization}
				%
				$\blacktriangleright$ \textbf{Initialization} We shall verify that the assumptions \eqref{reg-G-10} and \eqref{assum-Zi1} required along the KAM step to get  the final form \eqref{KAM step remainder term} are satisfied for $\mathscr{L}=\mathscr{L}_0$ in \eqref{mouka1}. Indeed, \eqref{reg-G-10} is an immediate consequence of Lemma \ref{lemma properties linear frequencies}-(vi), that is
				\begin{align}\label{reg-G-1}\exists C>0,\,\, \forall(j,j_{0})\in\mathbb{Z}^{2},\quad\max_{|\alpha|\in\llbracket 0, q\rrbracket}\sup_{\lambda\in(\lambda_{0},\lambda_{1})}\left|\partial_{\lambda}^{\alpha}\left(\Omega_{j}(\lambda)-\Omega_{j_{0}}(\lambda)\right)\right|\leqslant C\,|j-j_{0}|.
				\end{align}
				Thus, applying \eqref{differences mu0} we obtain 
				\begin{align*}%\label{reg-G-10}
					\exists C>0,\,\, \forall(j,j_{0})\in\mathbb{Z}^{2},\quad\max_{|\alpha|\in\llbracket 0, q\rrbracket}\sup_{(\lambda,\omega)\in\mathcal{O}}\left|\partial_{\lambda,\omega}^{\alpha}\left(\mu_{j}^{0}(\lambda,\omega)-\mu_{j_{0}}^{0}(\lambda,\omega)\right)\right|\leqslant C\,|j-j_{0}|.
				\end{align*}
				Concerning the second assumption \eqref{assum-Zi1}, we may combine  \eqref{estimate mathscr R in off diagonal norm} and \eqref{hypothesis KAM reduction of the remainder term} to find 
				\begin{align*}%\label{ZZ-KLPP}
					\nonumber\gamma ^{-1}\| \mathscr{R}_{0}\|_{\textnormal{\tiny{O-d}},q,s_0+\tau_{2} q+\tau_{2}}^{\gamma ,\mathcal{O}}&\leqslant C \varepsilon \gamma ^{-2}\left(1+\|\mathfrak{I}_{0}\|_{q,s_h+\sigma_{4}}^{\gamma,\mathcal{O}}\right)\\
					&\leqslant C\varepsilon_{0}.
				\end{align*}
				%
				%$\blacktriangleright$ \textbf{KAM iteration}
				%
				$\blacktriangleright$ \textbf{KAM iteration}. Let $m\in\mathbb{N}$ and consider a linear  operator 
				\begin{align}\label{Op-Lm}
					\mathscr{L}_{m}:=\big(\omega\cdot\partial_{\varphi}+\mathscr{D}_{m}+\mathscr{R}_{m}\big)\Pi_{\mathbb{S}_0}^{\perp}
				\end{align}
				with $\mathscr{D}_m$ a diagonal real reversible  operator and $\mathscr{R}_m$  a  real and reversible Toeplitz in time operator  of zero order satisfying $\Pi_{\mathbb{S}_0}^\perp \mathscr{R}_m\Pi_{\mathbb{S}_0}^\perp =\mathscr{R}_m.$ We assume that both assumptions \eqref{reg-G-10} and \eqref{assum-Zi1} are satisfied for $\mathscr{D}_m$ and $\mathscr{R}_m$. Remark  that for $m=0$ we take the operator $\mathscr{L}_0$ defined in \eqref{mouka1}. Let  $\Phi_{m}=\hbox{Id}+\Psi_{m}$ be a linear invertible operator  such that 
				\begin{align}\label{Op-Lm1}
					\Phi_{m}^{-1}\mathscr{L}_{m}\Phi_{m}:=\big(\omega\cdot\partial_{\varphi}+\mathscr{D}_{m+1}+\mathscr{R}_{m+1}\big)\Pi_{\mathbb{S}_0}^{\perp},
				\end{align}
				with $\Psi_{m}$ satisfying the homological equation
				$$
				\big[\big(\omega\cdot\partial_{\varphi}+\mathscr{D}_{m}\big)\Pi_{\mathbb{S}_0}^{\perp},\Psi_{m}\big]+P_{N_{m}}\mathscr{R}_{m}=\lfloor P_{N_{m}}\mathscr{R}_{m}\rfloor.
				$$
				Recall that $N_{m}$ was defined in \eqref{definition of Nm}. The diagonal parts   $(\mathscr{D}_{m})_{m\in\mathbb{N}}$ and the remainders $(\mathscr{R}_{m})_{m\in\mathbb{N}}$ are defined similarly to \eqref{PU-RT}  by the recursive formulas,
				\begin{align}\label{Deco-T1}
					\mathscr{D}_{m+1}=\mathscr{D}_{m}+\lfloor P_{N_{m}}\mathscr{R}_{m}\rfloor\quad \mbox{ and }\quad \mathscr{R}_{m+1}=\Phi_{m}^{-1}\left(-\Psi_{m}\,\lfloor P_{N_{m}}\mathscr{R}_{m}\rfloor +P_{N_{m}}^{\perp}\mathscr{R}_{m}+\mathscr{R}_{m}\Psi_{m}\right).
				\end{align}
				Remark  that $\mathscr{D}_{m}$ and $\lfloor P_{N_{m}}\mathscr{R}_{m}\rfloor$ are Fourier multiplier Toeplitz operators that can be identified to their spectra  $(\ii \mu_{j}^{m})_{j\in\mathbb{S}_{0}^{c}}$ and $(\ii r_{j}^{m})_{j\in\mathbb{S}_{0}^{c}}$, namely 
				\begin{align}\label{Spect-T1}
					\forall (l,j)\in\mathbb{Z}^d\times\mathbb{S}_0^c,\quad \mathscr{D}_{m} {\bf e}_{l,j}=\ii\mu_{j}^{m}\,{\bf e}_{l,j}\quad\hbox{and}\quad \lfloor P_{N_{m}}\mathscr{R}_{m}\rfloor \mathbf{e}_{l,j}=\ii r_{j}^{m}\,{\bf e}_{l,j}.
				\end{align}
				By construction, we find
				\begin{align}\label{Spect-T2}
					\mu_{j}^{m+1}=\mu_{j}^{m}+r_{j}^{m}.
				\end{align}
				In a similar way to \eqref{Psi-gh} we obtain 
				\begin{equation}\label{Psi-gh1}
					(\Psi_{m})_{j_{0}}^{j}(\lambda,\omega,l)=\left\lbrace\begin{array}{ll}
						\frac{-r_{j_{0},m}^{j}(\lambda,\omega,l)}{\omega\cdot l+\mu_{j}^{0}(\lambda,\omega)-\mu_{j_{0}}^{0}(\lambda,\omega)} & \mbox{if }(l,j)\neq(0,j_{0})\\
						0 & \mbox{if }(l,j)=(0,j_{0}),
					\end{array}\right.
				\end{equation}
				where  the collection $\{r_{j_{0},m}^{j}(\lambda,\omega,l)\}$ describes the Fourier coefficients of $\mathscr{R}_{m}$, that is,
				$$
				\mathscr{R}_{m} {\bf e}_{l_0,j_0}=\ii\sum_{(l,j)\in\mathbb{Z}^{d+1}}{r_{j_{0},m}^{j}(\lambda,\omega,l_0-l){\bf e}_{l,j}}.
				$$
				Now we  shall define  the open Cantor set where the preceding formula is meaningful,
				\begin{equation}\label{Cantor-SX}
					\mathscr{O}_{m+1}^{\gamma }=\bigcap_{\underset{(l,j)\neq(0,j_0)}{\underset{|l|\leqslant N_{m}}{(l,j,j_{0})\in\mathbb{Z}^{d }\times(\mathbb{S}_{0}^{c})^{2}}}}\left\lbrace(\lambda,\omega)\in\mathscr{O}_{m}^{\gamma }\quad\textnormal{s.t.}\quad |\omega\cdot l+\mu_{j}^{m}(\lambda,\omega)-\mu_{j_{0}}^{m}(\lambda,\omega)|>\tfrac{\gamma\langle j-j_{0}\rangle }{\langle l\rangle^{\tau_{2}}}\right\rbrace.
				\end{equation}
				Similarly to  \eqref{Ext-psi-op} and \eqref{varr-d} we can  extend \eqref{Psi-gh1} as follows 
				\begin{equation}\label{coeff-psim} 
					(\Psi_{m})_{j_{0}}^{j}(\lambda,\omega,l)=\left\lbrace\begin{array}{ll}
						-\frac{\chi\left((\omega\cdot l+\mu_{j}^{m}(\lambda,\omega)-\mu_{j_{0}}^{m}(\lambda,\omega))(\gamma|j-j_{0}|)^{-1}\langle l\rangle^{\tau_{2}}\right){r}_{j_{0},m}^{j}(\lambda,\omega,l)}{\omega\cdot l+\mu_{j}^{m}(\lambda,\omega)-\mu_{j_{0}}^{m}(\lambda,\omega)} & \mbox{if }(l,j)\neq(0,j_{0})\\
						0 & \mbox{if }(l,j)=(0,j_{0}).
					\end{array}\right.
				\end{equation}
				We point out  that working with this extension for $\Psi_m$ allows to extend  both $\mathscr{D}_{m+1}$ and   the remainder $\mathscr{R}_{m+1}$ provided that the  operators $\mathscr{D}_{m}$ and $\mathscr{R}_{m}$ are defined in the whole range of parameters. Thus the  operator defined by   the right-hand side in \eqref{Op-Lm1} can be  extended to  the whole set $\mathcal{O}$ and we denote this extension by $\mathscr{L}_{m+1}$, that is,
				\begin{align}\label{Op-Lm2}
					\big(\omega\cdot\partial_{\varphi}+\mathscr{D}_{m+1}+\mathscr{R}_{m+1}\big)\Pi_{\mathbb{S}_0}^{\perp}:=\mathscr{L}_{m+1}.
				\end{align}
				This enables to construct by induction the sequence of operators $\left(\mathscr{L}_{m+1}\right)$ in the full set  $\mathcal{O}$. Similarly the operator $\Phi_{m}^{-1}\mathscr{L}_{m}\Phi_{m}$ admits an extension in $\mathcal{O}$ induced by the extension of $\Phi_m^{\pm1}$ . However, by   construction the identity $\mathscr{L}_{m+1}=\Phi_{m}^{-1}\mathscr{L}_{m}\Phi_{m}$ in \eqref{Op-Lm1} occurs in the  Cantor set $\mathscr{O}_{m+1}^{\gamma }$ and may fail outside this set.
				We define
				\begin{equation}\label{Def-Taou}
					\delta_{m}(s):=\gamma ^{-1}\|\mathscr{R}_{m}\|_{\textnormal{\tiny{O-d}},q,s}^{\gamma ,\mathcal{O}}
				\end{equation}
				and we want to   prove by induction in $m\in\mathbb{N}$ that 
				\begin{equation}\label{hypothesis of induction for deltamprime}
					\forall\, m\in\mathbb{N},\quad \forall s\in[s_{0},\overline{s}_{l}],\,\delta_{m}(s)\leqslant \delta_{0}(s_{h})N_{0}^{\mu_{2}}N_{m}^{-\mu_{2}}\quad \mbox{ and }\quad \delta_{m}(s_{h})\leqslant\left(2-\tfrac{1}{m+1}\right)\delta_{0}(s_{h}),
				\end{equation}
				with $\overline{s}_{l}$ and $s_h$ fixed by \eqref{Conv-T2}. Moreover, we should check the validity of the  assumptions \eqref{reg-G-10} and \eqref{assum-Zi1}  for $\mathscr{D}_{m+1}$ and $\mathscr{R}_{m+1}$.
				Notice that by Sobolev embeddings, it is sufficient to prove the first inequality with $s=\overline{s}_{l}.$ The  property is obvious  for $m=0$. Now, assume that  the property \eqref{hypothesis of induction for deltamprime} is true for $m\in\mathbb{N}$ and let us check it at the next order.
				We write 
				\begin{align}\label{phi-inverse1}
					\Phi_{m}^{-1}=\hbox{Id}+\Sigma_{m}\quad\hbox{with}\quad \Sigma_{m}=\displaystyle\sum_{n=1}^{\infty}(-1)^{n}\Psi_{m}^{n}.
				\end{align} 
				Thus similarly to  \eqref{QWX0}, using in particular \eqref{link Psi and R} and \eqref{ZZ-KL} we deduce   successively
				\begin{align*}
					\|\Sigma_{m}\|_{\textnormal{\tiny{O-d}},q,s_0}^{\gamma,\mathcal{O}} & \leqslant \displaystyle\|\Psi_{m}\|_{\textnormal{\tiny{O-d}},q,s_0}^{\gamma,\mathcal{O}}\left(1+\sum_{n=0}^{\infty}\big(C\|\Psi_{m}\|_{\textnormal{\tiny{O-d}},q,s_{0}}^{\gamma,\mathcal{O}}\big)^{n}\right)\\
					& \leqslant  \displaystyle \delta_{m}(s_0+\tau_{2} q+\tau_{2})\left(1+\sum_{n=0}^{\infty}\big(C\delta_{m}(s_{0}+\tau_{2} q+\tau_{2})\big)^{n}\right)
				\end{align*}
				and for any $s\in[s_0,S],$
				\begin{align*}
					\|\Sigma_{m}\|_{\textnormal{\tiny{O-d}},q,s}^{\gamma,\mathcal{O}} & \leqslant  \displaystyle\|\Psi_{m}\|_{\textnormal{\tiny{O-d}},q,s}^{\gamma,\mathcal{O}}\left(1+\sum_{n=0}^{\infty}\big(C\|\Psi_{m}\|_{\textnormal{\tiny{O-d}},q,s_{0}}^{\gamma,\mathcal{O}}\big)^{n}\right)\\
					& \leqslant  \displaystyle N_{m}^{\tau_{2} q+\tau_{2}}\delta_{m}(s)\left(1+\sum_{n=0}^{\infty}\big(C\delta_{m}(s_{0}+\tau_{2} q+\tau_{2})\big)^{n}\right).
				\end{align*}
				Hence, from the induction assumption, the fact that $N_{m}\geqslant N_{0}$ and since \eqref{param} implies in particular $s_{0}+\tau_{2} q+\tau_{2}\leqslant \overline{s}_{l}$, we obtain
				\begin{align*}
					\|\Sigma_{m}\|_{\textnormal{\tiny{O-d}},q,s_0}^{\gamma,\mathcal{O}} & \leqslant  \displaystyle CN_{0}^{\mu_{2}}N_{m}^{-\mu_{2}}\delta_{0}(s_{h})\left(1+\sum_{n=0}^{\infty}\big(CN_{0}^{\mu_{2}}N_{m}^{-\mu_{2}}\delta_{0}(s_{h})\big)^{n}\right)\\
					& \leqslant \displaystyle CN_{0}^{\mu_{2}}N_{m}^{-\mu_{2}}\delta_{0}(s_{h})\left(1+\sum_{n=0}^{\infty}\big(C\delta_{0}(s_{h})\big)^{n}\right)
				\end{align*}
				and for any $s\in[s_0,S],$
				\begin{align*}
					\|\Sigma_{m}\|_{\textnormal{\tiny{O-d}},q,s}^{\gamma,\mathcal{O}}  &\leqslant  \displaystyle N_{m}^{\tau_{2} q+\tau_{2}}\delta_{m}(s)\left(1+\sum_{n=0}^{\infty}\big(CN_{0}^{\mu_{2}}N_{m}^{-\mu_{2}}\delta_{0}(s_{h})\big)^{n}\right)\\
					&\leqslant  \displaystyle N_{m}^{\tau_{2} q+\tau_{2}}\delta_{m}(s)\left(1+\sum_{n=0}^{\infty}\big(C\delta_{0}(s_{h})\big)^{n}\right).
				\end{align*}
				It follows from the condition \eqref{Conv-P3} that 
				\begin{align}\label{AP-W34}
					\|\Sigma_{m}\|_{\textnormal{\tiny{O-d}},q,s_0}^{\gamma,\mathcal{O}}\leqslant CN_{0}^{\mu_{2}}N_{m}^{-\mu_{2}}\delta_{0}(s_{h})\quad\hbox{and}\quad \|\Sigma_{m}\|_{\textnormal{\tiny{O-d}},q,s}^{\gamma,\mathcal{O}} & \leqslant  CN_{m}^{\tau_{2} q+\tau_{2}}\delta_{m}(s).
				\end{align}
				One also gets
				\begin{align}\label{Bir-11}
					\|\Sigma_{m}\|_{\textnormal{\tiny{O-d}},q,s}^{\gamma,\mathcal{O}} & \leqslant  C\delta_{m}(s+\tau_2q+\tau_2).
				\end{align}
				From   KAM step \eqref{KAM step remainder term} and  Sobolev embeddings,  we infer
				$$\delta_{m+1}(\overline{s}_{l})\leqslant  N_{m}^{\overline{s}_{l}-s_{h}}\delta_{m}(s_{h})+CN_{m}^{\tau_{2}q+\tau_{2}}\left(\delta_{m}(\overline{s}_{l})\right)^{2}.
				$$
				Using the induction assumption \eqref{hypothesis of induction for deltamprime} yields
				\begin{align*}
					\delta_{m+1}(\overline{s}_{l}) & \leqslant N_{m}^{\overline{s}_{l}-s_{h}}\left(2-\tfrac{1}{m+1}\right)\delta_{0}(s_{h})+CN_{m}^{\tau_{2}q+\tau_{2}}\delta_{0}^2(s_{h})N_{0}^{2\mu_{2}}N_{m}^{-2\mu_{2}}\\
					& \leqslant 2 N_{m}^{\overline{s}_{l}-s_{h}}\delta_{0}(s_{h})+CN_{m}^{\tau_{2}q+\tau_{2}}\delta_{0}^2(s_{h})N_{0}^{2\mu_{2}}N_{m}^{-2\mu_{2}}.
				\end{align*}
				At this level we need to select  the parameters $\overline{s}_{l},s_h$ and $\mu_2$ in such a way 
				\begin{equation}\label{Conv-T1}
					N_{m}^{\overline{s}_{l}-s_{h}} \leqslant  \frac{1}{4}N_{0}^{\mu_{2}}N_{m+1}^{-\mu_{2}}\qquad\hbox{and}\qquad 
					CN_{m}^{\tau_{2}q+\tau_{2}}\delta_{0}(s_{h})N_{0}^{2\mu_{2}}N_{m}^{-2\mu_{2}}  \leqslant  \frac{1}{2}N_{0}^{\mu_{2}}N_{m+1}^{-\mu_{2}}
				\end{equation}
				leading to 
				$$
				\delta_{m+1}(\overline{s}_{l}) \leqslant \delta_{0}(s_{h})N_{0}^{\mu_{2}}N_{m+1}^{-\mu_{2}}.
				$$
				The conditions \eqref{Conv-T2} imply in particular 
				$$s_h\geqslant \frac{3}{2}\mu_{2}+\overline{s}_{l}+1\quad\textnormal{and}\quad \mu_{2}\geqslant2(\tau_2q+\tau_2)+1.$$
				Then, using \eqref{definition of Nm}, we conclude that the assumptions of \eqref{Conv-T1} hold true provided that 
				\begin{equation}\label{Cond1}
					4N_{0}^{-\mu_2}\leqslant 1\qquad\hbox{and}\qquad 
					2C N_{0}^{\mu_{2}}\delta_{0}(s_{h}) \leqslant  1,
				\end{equation}
				which follow from  \eqref{Conv-P3}, since the first condition $4N_{0}^{-\mu_2}\leqslant 1$ is automatically satisfied because $N_0\geqslant 2$ and $\mu_2\geqslant 2,$ according to \eqref{Conv-T2}.
				Therefore, under the assumptions \eqref{Conv-T2} we get  the first statement of the induction in \eqref{hypothesis of induction for deltamprime}. The next goal is  to establish the second  estimate  in \eqref{hypothesis of induction for deltamprime}. By  KAM step \eqref{KAM step remainder term} combined with  the induction assumptions \eqref{hypothesis of induction for deltamprime} we deduce that
				\begin{align*}
					\delta_{m+1}(s_{h}) & \leqslant  \delta_{m}(s_{h})+C N_{m}^{\tau_{2}q+\tau_{2}}\delta_{m}(s_{0}) \delta_{m}(s_{h})\\
					&\leqslant  \Big(2-\tfrac{1}{m+1}\Big)\delta_{0}(s_{h})\Big(1+CN_{0}^{\mu_{2}}N_{m}^{\tau_{2}q+\tau_{2}-\mu_{2}}\delta_{0}(s_{h})\Big).
				\end{align*}
				Thus if one has
				\begin{equation}\label{Conv-4}
					\Big(2-\tfrac{1}{m+1}\Big)\Big(1+CN_{0}^{\mu_{2}}N_{m}^{\tau_{2}q+\tau_{2}-\mu_{2}}\delta_{0}(s_{h})\Big)\leqslant 2-\tfrac{1}{m+2},
				\end{equation}
				then we get
				$$
				\delta_{m+1}(s_{h})  \leqslant  \left(2-\tfrac{1}{m+2}\right)\delta_{0}(s_{h}),
				$$
				which ends the induction argument of \eqref{hypothesis of induction for deltamprime}. Remark  that with the choice $\mu_2\geqslant 2(\tau_2q+\tau_2)$ fixed in \eqref{Conv-T2}, the condition  \eqref{Conv-4} is staisfied  if  
				\begin{equation}\label{Conv-5}
					CN_{0}^{\mu_{2}}N_{m}^{-\tau_{2}q-\tau_{2}}\delta_{0}(s_{h})\leqslant\tfrac{1}{(2m+1)(m+2)}\cdot
				\end{equation}
				Since $N_{0}\geqslant 2$ we may find a constant $c_0>0$ small enough  such that
				$$
				\forall\, m\in\mathbb{N},\quad c_0N_m^{-1}\leq \tfrac{1}{(2m+1)(m+2)}\cdot
				$$
				Consequently,  \eqref{Conv-5} is satisfied  provided that
				\begin{equation}\label{Conv-6}
					CN_{0}^{\mu_{2}}N_{m}^{-\tau_{2}q-\tau_{2}+1}\delta_{0}(s_{h})\leqslant c_0.
				\end{equation}
				By virtue of  the assumption   \eqref{setting tau1 and tau2} we get in particular
				\begin{equation}\label{Conv-7}
					\tau_{2}q+\tau_{2}-1\geqslant0.
				\end{equation}
				Thus  \eqref{Conv-6} is satisfied in view of  \eqref{Conv-P3}. To conclude  the induction   proof of \eqref{hypothesis of induction for deltamprime} it remains to check that the  assumptions \eqref{reg-G-10} and \eqref{assum-Zi1} are satisfied for $\mathscr{D}_{m+1}$ and $\mathscr{R}_{m+1}$. First, the assumption  \eqref{assum-Zi1} is a consequence of  the first inequality of  \eqref{hypothesis of induction for deltamprime}  applied at the order $m+1$  with the regularity index $s=s_0+\tau_2q+\tau_2\leqslant\overline{s}_{l}$ supplemented with \eqref{Conv-P3}. Concerning the validity of  \eqref{reg-G-10} for $\mathscr{D}_{m+1}$, we combine \eqref{Spect-T1} and \eqref{Spect-T2} and  \eqref{Dp1X}, in order to find
				\begin{align*}
					\|\mu_{j}^{m+1}-\mu_{j}^{m}\|_{q}^{\gamma ,\mathcal{O}}&=\big\|\langle P_{N_{m}}\mathscr{R}_{m}\mathbf{e}_{l,j},\mathbf{e}_{l,j}\rangle_{L^{2}(\mathbb{T}^{d +1})}\big\|_{q}^{\gamma ,\mathcal{O}}.
				\end{align*}
				From the Topeplitz structure of  $\mathscr{R}_{m}$ we may write 
				\begin{align*}
					\|\mu_{j}^{m+1}-\mu_{j}^{m}\|_{q}^{\gamma ,\mathcal{O}}&=\big\|\langle P_{N_{m}}\mathscr{R}_{m}\mathbf{e}_{0,j},\mathbf{e}_{0,j}\rangle_{L^{2}(\mathbb{T}^{d +1})}\big\|_{q}^{\gamma ,\mathcal{O}}.
				\end{align*}
				By a  duality argument combined with  Lemma \ref{properties of Toeplitz in time operators} and \eqref{Def-Taou} we infer
				\begin{align}\label{dual-1X}
					\nonumber \|\mu_{j}^{m+1}-\mu_{j}^{m}\|_{q}^{\gamma ,\mathcal{O}}&\lesssim \big\|\mathscr{R}_{m}\mathbf{e}_{0,j}\|_{q,s_0}^{\gamma,\mathcal{O}}\,\,\langle j\rangle^{-s_0}\\
					\nonumber &\lesssim\big\|\mathscr{R}_{m}\big\|_{\textnormal{\tiny{O-d}},q,s_0}^{\gamma,\mathcal{O}}\|\mathbf{e}_{0,j}\|_{H^{s_0}}\,\,\langle j\rangle^{-s_0}\\
					&\lesssim \big\|\mathscr{R}_{m}\big\|_{\textnormal{\tiny{O-d}},q,s_0}^{\gamma,\mathcal{O}}=\gamma \delta_{m}(s_{0}).
				\end{align}
				Hence  we deduce from \eqref{hypothesis of induction for deltamprime}, \eqref{whab1} and \eqref{hypothesis KAM reduction of the remainder term} 
				\begin{align}\label{Mahma1-R}
					\|\mu_{j}^{m+1}-\mu_{j}^{m}\|_{q}^{\gamma ,\mathcal{O}}
					&\leqslant C \gamma\,\delta_{0}(s_{h})N_{0}^{\mu_{2}}N_{m}^{-\mu_{2}}\nonumber\\
					&\leqslant C \varepsilon\gamma^{-1}N_0^{\mu_{2}}N_{m}^{-\mu_{2}}.
				\end{align}
				As the assumption \eqref{reg-G-10} is satisfied  with  $\mathscr{D}_{m}$, that is, 
				\begin{align}\label{e-mujm} \forall\,(j,j_{0})\in(\mathbb{S}_0^c)^{2},\quad\max_{|\alpha| \in\llbracket 0,q\rrbracket}\sup_{(\lambda,\omega)\in\mathcal{O}}\left|\partial_{\lambda,\omega}^{\alpha}\left(\mu_{j}^m(\lambda,\omega)-\mu_{j_{0}}^m(\lambda,\omega)\right)\right|\leqslant C\,|j-j_{0}|,
				\end{align}
				then we obtain  by  \eqref{Mahma1-R}
				\begin{align*} \forall\,(j,j_{0})\in(\mathbb{S}_0^c)^{2},\quad\max_{|\alpha| \in\llbracket 0,q\rrbracket}\sup_{(\lambda,\omega)\in\mathcal{O}}\left|\partial_{\lambda,\omega}^{\alpha}\left(\mu_{j}^{m+1}(\lambda,\omega)-\mu_{j_{0}}^{m+1}(\lambda,\omega)\right)\right|\leqslant C\big(1+\varepsilon\gamma^{-1-q}N_0^{\mu_2}N_{m}^{-\mu_{2}}\big)\,|j-j_{0}|.
				\end{align*}
				Consequently, the convergence of the series $\sum N_m^{-\mu_2}$  gives  the required  assumption with the  same constant $C$ independently of $m$. This completes the induction principle. In what follows,  we shall provide some estimates for $\Psi_m$ that will be used later to study the string convergence.  
				Using \eqref{link Psi and R} combined with Lemma  \ref{properties of Toeplitz in time operators} and $s_0+\tau_2q+\tau_2+1\leqslant\overline{s}_l$ we find
				\begin{align}\label{link-ZD3}
					\nonumber  \|\Psi_m\|_{\textnormal{\tiny{O-d}},q,s_0+1}^{\gamma ,\mathcal{O}}&\leqslant C \gamma ^{-1}\| P_{N_{m}}\mathscr{R}_{m}\|_{\textnormal{\tiny{O-d}},q,{s_0+\tau_2q+\tau_2+1}}^{\gamma ,\mathcal{O}}\\
					&\leqslant\, C\, \delta_m(\overline{s}_l).
				\end{align}
				Thus  \eqref{hypothesis of induction for deltamprime} and \eqref{Conv-P3} yield
				\begin{align}\label{link-Z3}
					\nonumber  \|\Psi_m\|_{\textnormal{\tiny{O-d}},q,s_0+1}^{\gamma ,\mathcal{O}}&\leqslant\, C\, \delta_{0}(s_{h})N_{0}^{\mu_{2}}N_{m}^{-\mu_{2}}\\
					&\leqslant\, C\,\varepsilon \gamma^{-2}N_0^{\mu_2} \,N_{m}^{-\mu_{2}}.
				\end{align}
				Next, we  discuss the persistence of higher regularity. Let  $s\in[s_0,S]$, then from \eqref{KAM step remainder term}, \eqref{hypothesis of induction for deltamprime} and   \eqref{Conv-P3}  and \eqref{Conv-7}
				\begin{align*}
					\delta_{m+1}(s) & \leqslant  \delta_{m}(s)\Big(1+CN_{m}^{\tau_{2} q+\tau_{2}}\delta_{m}(s_{0})\Big)\\
					& \leqslant  \delta_{m}(s)\Big(1+CN_{0}^{{\mu}_{2}}N_{m}^{\tau_{2} q+\tau_{2}-\mu_{2}}\delta_{0}(s_{h})\Big)\\
					& \leqslant  \delta_{m}(s)\big(1+CN_{m}^{-1}\big).
				\end{align*}
				Combining  this estimate with  \eqref{definition of Nm} and \eqref{whab1} yields
				\begin{align}\label{uniform estimate of deltamsprime}
					\nonumber \forall\, s\geqslant  s_0,\, \forall\, m\in\mathbb{N},\quad \delta_{m}(s)&\leqslant\delta_{0}(s)\prod_{n=0}^{\infty}\big(1+CN_{n}^{-1}\big)\\
					\nonumber&\leqslant C\delta_{0}(s)\\
					&\leqslant C\varepsilon\gamma^{-2}\left(1+\|\mathfrak{I}_{0}\|_{q,s+\sigma_{4}}^{\gamma,\mathcal{O}}\right).
				\end{align}
				Using \eqref{link Psi and R} combined with Lemma  \ref{properties of Toeplitz in time operators}, applying in particular interpolation inequalities,   leads to
				\begin{align}\label{link-Z1}
					\nonumber \|\Psi_m\|_{\textnormal{\tiny{O-d}},q,s}^{\gamma ,\mathcal{O}}&\leqslant C \gamma ^{-1}\| P_{N_{m}}\mathscr{R}_{m}\|_{\textnormal{\tiny{O-d}},q,s+\tau_{2} q+\tau_{2}}^{\gamma ,\mathcal{O}}\\
					\nonumber&\leqslant\, C\, \delta_m(s+\tau_{2} q+\tau_{2})\\
					&\leqslant\, C\, \delta_m^{\overline\theta}(s_0)\delta_m^{1-\overline\theta}(s+\tau_{2} q+\tau_{2}+1).
				\end{align}
				with $\overline\theta=\frac{1}{s-s_0+\tau_2q+\tau_2+1}.$ Inserting  \eqref{hypothesis of induction for deltamprime} and  \eqref{uniform estimate of deltamsprime}  into \eqref{link-Z1} and using    \eqref{Conv-P3} give
				\begin{align}\label{link-Z2}
					\nonumber \|\Psi_m\|_{\textnormal{\tiny{O-d}},q,s}^{\gamma ,\mathcal{O}}&\leqslant\, C\, \delta_0^{\overline\theta}(s_2)\delta_0^{1-\overline\theta}(s+\tau_{2} q+\tau_{2}+1)N_0^{\mu_2\overline\theta} N_m^{-\mu_2\overline\theta}\\
					&\leqslant\, C\,{\varepsilon_{0}}^{\overline\theta} \delta_0^{1-\overline\theta}(s+\tau_{2} q+\tau_{2}+1) N_m^{-\mu_2\overline\theta}.
				\end{align}
				We point out  that one also finds from \eqref{Bir-11},  the second inequality of  \eqref{link-Z1} and \eqref{uniform estimate of deltamsprime} that
				\begin{align}\label{link-ZP2}
					\forall\, s\in[s_0,S],\quad \sup_{m\in\mathbb{N}}\left(\|\Sigma_m\|_{\textnormal{\tiny{O-d}},q,s}^{\gamma ,\mathcal{O}}+\|\Psi_m\|_{\textnormal{\tiny{O-d}},q,s}^{\gamma ,\mathcal{O}}\right)&\leqslant\,C\varepsilon\gamma^{-2}\left(1+\|\mathfrak{I}_{0}\|_{q,s+\sigma_{4}}^{\gamma,\mathcal{O}}\right).
				\end{align}
				$\blacktriangleright$ \textbf{KAM conclusion}. Let us  examine  the sequence of operators 
				$\left(\widehat{\Phi}_m\right)_{m\in\mathbb{N}}$ defined by 
				\begin{equation}\label{Def-Phi}\widehat\Phi_0:=\Phi_0\quad\hbox{and}\quad \quad  \forall m\geqslant1,\,\, \widehat\Phi_m:=\Phi_0\circ\Phi_1\circ...\circ\Phi_m.
				\end{equation} 
				It is obvious  from the identity ${\Phi}_{m}=\hbox{Id}+\Psi_m$ that $\widehat\Phi_{m+1}=\widehat\Phi_{m}+\widehat\Phi_{m}\Psi_{m+1}$.  Applying   the law products yields
				$$\|\widehat{\Phi}_{m+1}\|_{\textnormal{\tiny{O-d}},q,s_{0}+1}^{\gamma,\mathcal{O}}\leqslant \|\widehat{\Phi}_{m}\|_{\textnormal{\tiny{O-d}},q,s_{0}+1}^{\gamma,\mathcal{O}}\Big(1+C\|\Psi_{m+1}\|_{\textnormal{\tiny{O-d}},q,s_{0}+1}^{\gamma,\mathcal{O}}\Big).$$
				By  iterating this inequality and using \eqref{link-Z3} we infer
				\begin{align*}
					\|\widehat{\Phi}_{m+1}\|_{\textnormal{\tiny{O-d}},q,s_{0}+1}^{\gamma,\mathcal{O}} 
					& \leqslant  \displaystyle\|\Phi_{0}\|_{\textnormal{\tiny{O-d}},q,s_{0}+1}^{\gamma,\mathcal{O}}\prod_{n=1}^{m+1}\Big(1+C\|\Psi_{n}\|_{\textnormal{\tiny{O-d}},q,s_{0}+1}^{\gamma,\mathcal{O}}\Big)\\
					& \leqslant  \displaystyle\prod_{n=0}^{\infty}\Big(1+C\,\varepsilon_{0}N_{n}^{-\mu_{2}}\Big).
				\end{align*}
				Using the  first condition of \eqref{Conv-P3} and  \eqref{definition of Nm} imply
				$$\|\widehat{\Phi}_{m+1}\|_{\textnormal{\tiny{O-d}},q,s_{0}+1}^{\gamma,\mathcal{O}}\leqslant\displaystyle\prod_{n=0}^{\infty}\Big(1+C\,\varepsilon_{0}4^{-(\frac32)^n}\Big)$$
				and since  the infinite product converges, we obtain for  $\varepsilon_{0}$ small enough  
				\begin{align}\label{phim-T}
					\sup_{m\in\mathbb{N}}\|\widehat{\Phi}_{m}\|_{\textnormal{\tiny{O-d}},q,s_{0}+1}^{\gamma,\mathcal{O}}\leqslant 2.
				\end{align}
				Now we shall estimate the difference $\widehat{\Phi}_{m+1}-\widehat{\Phi}_{m}$ and for this aim we use the law products combined with \eqref{link-Z3} and \eqref{phim-T}
				\begin{align}\label{Trank-P1}
					\|\widehat{\Phi}_{m+1}-\widehat{\Phi}_{m}\|_{\textnormal{\tiny{O-d}},q,s_{0}+1}^{\gamma ,\mathcal{O}}
					\nonumber& \leqslant   C \|\widehat{\Phi}_{m}\|_{\textnormal{\tiny{O-d}},q,s_{0}+1}^{\gamma ,\mathcal{O}}\|\Psi_{m+1}\|_{\textnormal{\tiny{O-d}},q,s_{0}+1}^{\gamma ,\mathcal{O}}\\
					& \leqslant C\,\delta_{0}(s_{h})N_{0}^{\mu_{2}} \,N_{m+1}^{-\mu_2}.
				\end{align}
				Applying Lemma \ref{lemma sum Nn} gives 
				\begin{align}\label{Conv-Z1}
					\sum_{m=0}^\infty\|\widehat{\Phi}_{m+1}-\widehat{\Phi}_{m}\|_{\textnormal{\tiny{O-d}},q,s_{0}+1}^{\gamma ,\mathcal{O}}\leqslant C\,\delta_{0}(s_{h}).
				\end{align}
				Therefore, by a completeness argument we deduce that the series $\displaystyle \sum_{m\in\mathbb{N}}(\widehat\Phi_{m+1}-\widehat\Phi_{m})$ converges to an element $\Phi_\infty$. In addition, we get in view of \eqref{Trank-P1} and Lemma \ref{lemma sum Nn}
				\begin{align}\label{Conv-op-od}
					\nonumber\|\widehat\Phi_{m}-\Phi_\infty\|_{\textnormal{\tiny{O-d}},q,s_{0}+1}^{\gamma ,\mathcal{O}}&\leqslant \sum_{j=m}^{\infty}\|\widehat\Phi_{j+1}-\widehat\Phi_{j}\|_{\textnormal{\tiny{O-d}},q,s_{0}+1}^{\gamma ,\mathcal{O}}\\
					\nonumber&\leqslant  C\,\delta_{0}(s_{h})N_{0}^{\mu_{2}} \, \sum_{j=m}^{\infty}N_{j+1}^{-\mu_2}\\
					&\leqslant  C\,\delta_{0}(s_{h})N_{0}^{\mu_{2}}N_{m+1}^{-\mu_2}.
				\end{align} 
				Remark  that one also  finds from \eqref{phim-T}
				\begin{equation}\label{Conv-op-ES}
					\|\Phi_\infty\|_{\textnormal{\tiny{O-d}},q,s_{0}+1}^{\gamma ,\mathcal{O}}\leqslant 2.
				\end{equation} 
				Using \eqref{Conv-Z1} combined with  \eqref{link-Z1} for $m=0$ and \eqref{Conv-T2}
				\begin{align}\label{Conv-UL1}
					\nonumber \|\Phi_\infty-\textnormal{Id}\|_{\textnormal{\tiny{O-d}},q,s_0+1}^{\gamma ,\mathcal{O}}&\leqslant \sum_{m=0}^\infty\|\widehat{\Phi}_{m+1}-\widehat{\Phi}_{m}\|_{\textnormal{\tiny{O-d}},q,s_0+1}^{\gamma ,\mathcal{O}}+\|{\Psi}_{0}\|_{\textnormal{\tiny{O-d}},q,s_0+1}^{\gamma ,\mathcal{O}}\nonumber\\ &\leqslant C\,\delta_{0}(s_{h}).
				\end{align} 
				Let us now check   the convergence with  higher  order norms. Take $s\in[s_0,S]$, then using the  law products, \eqref{link-Z3},  \eqref{link-Z2} and \eqref{phim-T}  we infer
				\begin{align}\label{Rec-YP1}
					\|\widehat{\Phi}_{m+1}\|_{\textnormal{\tiny{O-d}},q,s}^{\gamma ,\mathcal{O}} & \leqslant  \|\widehat{\Phi}_{m}\|_{\textnormal{\tiny{O-d}},q,s}^{\gamma ,\mathcal{O}}\left(1+C\|\Psi_{m+1}\|_{\textnormal{\tiny{O-d}},q,s_{0}}^{\gamma ,\mathcal{O}}\right)+C\|\widehat{\Phi}_{m}\|_{\textnormal{\tiny{O-d}},q,s_{0}}^{\gamma ,\mathcal{O}}\|\Psi_{m+1}\|_{\textnormal{\tiny{O-d}},q,s}^{\gamma ,\mathcal{O}}\\
					\nonumber& \leqslant  \|\widehat{\Phi}_{m}\|_{\textnormal{\tiny{O-d}},q,s}^{\gamma ,\mathcal{O}}\left(1+C\,\varepsilon_{0} \,N_{m+1}^{-\mu_{2}}\right)+C\,\delta_{0}^{\overline\theta}(s_{h})N_{0}^{\mu_{2} \overline\theta} \delta_0^{1-\overline\theta}(s+\tau_{2} q+\tau_{2}+1) N_m^{-\mu_2\overline\theta}.
				\end{align}
				According to the first condition of \eqref{Conv-P3} and  \eqref{definition of Nm}  one finds
				\begin{align*}
					\prod_{n=0}^{\infty}\left(1+C\,\varepsilon_{0} \,N_{n}^{-\mu_{2}}\right)
					& \leqslant  \prod_{n=0}^{\infty}\Big(1+C\,\varepsilon_{0}4^{-(\frac32)^n}\Big)\\
					& \leqslant 2,
				\end{align*}
				where the last inequality holds if $\varepsilon_{0}$ is chosen small enough. 
				Applying \eqref{Ind-res1} together with  \eqref{Rec-YP1} and Lemma \ref{lemma sum Nn} and using \eqref{link Psi and R} yield
				\begin{align*}
					\nonumber \sup_{m\in\mathbb{N}}\|\widehat{\Phi}_{m}\|_{\textnormal{\tiny{O-d}},q,s}^{\gamma ,\mathcal{O}} 
					& \leqslant  C\Big( \|{\Phi}_{0}\|_{\textnormal{\tiny{O-d}},q,s}^{\gamma ,\mathcal{O}}+\delta_{0}^{\overline\theta}(s_{h})N_{0}^{\mu_{2} \overline\theta} \delta_0^{1-\overline\theta}(s+\tau_{2} q+\tau_{2}+1) \Big)\\
					& \leqslant  C\Big(1+ \delta_0(s+\tau_{2} q+\tau_{2})+ \delta_{0}^{\overline\theta}(s_{h})N_{0}^{\mu_{2} \overline\theta}\delta_0^{1-\overline\theta}\big(s+\tau_{2} q+\tau_{2}+1\big) \Big).
				\end{align*}
				Interpolation inequalities and \eqref{Conv-P3} allow to get
				\begin{align}\label{Rec-YP2}
					\sup_{m\in\mathbb{N}}\|\widehat{\Phi}_{m}\|_{\textnormal{\tiny{O-d}},q,s}^{\gamma ,\mathcal{O}} 
					& \leqslant  C\Big(  1+\delta_0\big(s+\tau_{2} q+\tau_{2}+1\big) \Big).
				\end{align}
				%Since $\overline\theta\in(0,1)$ then  we infer
				%\begin{align}\label{Rec-YP2}
				% \sup_{m\in\mathbb{N}}\|\widehat{\Phi}_{m}\|_{\textnormal{\tiny{O-d}},q,s}^{\gamma ,\mathcal{O}} 
				%& \leqslant  C\Big(1+ \delta_0(s+\tau_{2} q+\tau_{2}+1) \Big).
				%\end{align}
				The next task  is to estimate the difference $\|\widehat{\Phi}_{m+1}-\widehat{\Phi}_{m}\|_{q,s}^{\gamma ,\mathcal{O}} $. By the law products combined with the first inequality in \eqref{link-Z3}, \eqref{link-Z2}, \eqref{phim-T} and \eqref{Rec-YP2} we obtain
				\begin{align*}
					\|\widehat{\Phi}_{m+1}-\widehat{\Phi}_{m}\|_{\textnormal{\tiny{O-d}},q,s}^{\gamma ,\mathcal{O}}
					&\leqslant C\Big(\|\widehat{\Phi}_{m}\|_{\textnormal{\tiny{O-d}},q,s}^{\gamma ,\mathcal{O}}\|\Psi_{m+1}\|_{\textnormal{\tiny{O-d}},q,s_{0}}^{\gamma ,\mathcal{O}}+\|\widehat{\Phi}_{m}\|_{\textnormal{\tiny{O-d}},q,s_{0}}^{\gamma ,\mathcal{O}}\|\Psi_{m+1}\|_{\textnormal{\tiny{O-d}},q,s}^{\gamma ,\mathcal{O}}\Big)\\
					&\leqslant C\,\delta_{0}(s_{h})N_{0}^{\mu_{2}} \,N_{m+1}^{-\mu_2}\Big(1+\delta_0\big(s+\tau_{2} q+\tau_{2}+1\big) \Big)\\
					&\quad+C\,\delta^{\overline\theta}_{0}(s_{h})N_{0}^{\mu_{2} \overline\theta}\,\delta_0^{1-\overline\theta}\big(s+\tau_{2} q+\tau_{2}+1\big) N_{m+1}^{-\mu_2\overline\theta}.
				\end{align*}
				Thus, we obtain in view of Lemma \ref{lemma sum Nn} 
				\begin{align*}
					\sum_{m=0}^\infty\|\widehat{\Phi}_{m+1}-\widehat{\Phi}_{m}\|_{\textnormal{\tiny{O-d}},q,s}^{\gamma ,\mathcal{O}}
					&\leqslant C\,\delta_{0}(s_{h}) \Big(1+\delta_0\big(s+\tau_{2} q+\tau_{2}+1\big) \Big)\\
					&\quad+C\,\delta^{\overline\theta}_{0}(s_{h}) \delta_0^{1-\overline\theta}\big(s+\tau_{2} q+\tau_{2}+1\big).
				\end{align*}
				Combining the  interpolation inequalities with  the second condition in \eqref{Conv-P3} gives
				\begin{align}\label{P-R-S1}
					\sum_{m=0}^\infty\|\widehat{\Phi}_{m+1}-\widehat{\Phi}_{m}\|_{\textnormal{\tiny{O-d}},q,s}^{\gamma ,\mathcal{O}}
					&\leqslant C\Big(\delta_{0}(s_{h})+\delta_0\big(s+\tau_{2} q+\tau_{2}+1\big) \Big).
				\end{align}
				From this latter inequality  combined with \eqref{Conv-P3} and \eqref{Rec-YP2} we infer
				\begin{align}\label{Conv-opV1}
					\nonumber \|\Phi_\infty\|_{\textnormal{\tiny{O-d}},q,s}^{\gamma ,\mathcal{O}}&\leqslant \sum_{m=0}^\infty\|\widehat{\Phi}_{m+1}-\widehat{\Phi}_{m}\|_{\textnormal{\tiny{O-d}},q,s}^{\gamma ,\mathcal{O}}+\|\widehat{\Phi}_{0}\|_{\textnormal{\tiny{O-d}},q,s}^{\gamma ,\mathcal{O}}\\
					&\leqslant C\Big(1+\delta_0\big(s+\tau_{2} q+\tau_{2}+1\big) \Big).
				\end{align} 
				On the other hand, using \eqref{P-R-S1} and the second inequality in  \eqref{link-Z1} with $m=0,$  one can  check that
				\begin{align}\label{Conv-U1}
					\nonumber \|\Phi_\infty-\textnormal{Id}\|_{\textnormal{\tiny{O-d}},q,s}^{\gamma ,\mathcal{O}}&\leqslant \sum_{m=0}^\infty\|\widehat{\Phi}_{m+1}-\widehat{\Phi}_{m}\|_{\textnormal{\tiny{O-d}},q,s}^{\gamma ,\mathcal{O}}+\|{\Psi}_{0}\|_{\textnormal{\tiny{O-d}},q,s}^{\gamma ,\mathcal{O}}\\
					&\leqslant C\Big(\delta_{0}(s_{h})+\delta_0\big(s+\tau_{2} q+\tau_{2}+1\big) \Big).
				\end{align} 
				Therefore,  Lemma \ref{properties of Toeplitz in time operators} together with \eqref{Conv-op-ES}, \eqref{Conv-opV1} and Sobolev embeddings give 
				\begin{align}\label{Ti-Ti}
					\nonumber \|\Phi_{\infty}\rho\|_{q,s}^{\gamma ,\mathcal{O}}&\lesssim\| \Phi_{\infty}\|_{\textnormal{\tiny{O-d}},q,s_{0}}^{\gamma,\mathcal{O}}\|\rho\|_{q,s}^{\gamma,\mathcal{O}}+\| \Phi_{\infty}\|_{\textnormal{\tiny{O-d}},q,s}^{\gamma,\mathcal{O}}\|\rho\|_{q,s_{0}}^{\gamma,\mathcal{O}}\\
					\nonumber &\lesssim \|\rho\|_{q,s}^{\gamma,\mathcal{O}}+\Big(1+\delta_0\big(s+\tau_{2} q+\tau_{2}+1\big) \Big)\|\rho\|_{q,s_{0}}^{\gamma,\mathcal{O}}\\
					&\lesssim \|\rho\|_{q,s}^{\gamma,\mathcal{O}}+\delta_0\big(s+\tau_{2} q+\tau_{2}+1\big) \|\rho\|_{q,s_{0}}^{\gamma,\mathcal{O}}.
				\end{align}
				Applying   \eqref{estimate mathscr R in off diagonal norm} and \eqref{Def-Taou} we obtain				\begin{align}\label{fgh1} 
					\nonumber \delta_0\big(s+\tau_{2} q+\tau_{2}+1\big)&=\gamma^{-1}
					\|\mathscr{R}_{0}\|_{\textnormal{\tiny{O-d}},q,s}^{\gamma,\mathcal{O}}\\
					&\lesssim\varepsilon\gamma^{-2}\left(1+\| \mathfrak{I}_{0}\|_{q,s+\sigma_{4}}^{\gamma,\mathcal{O}}\right).
				\end{align}
				Plugging  \eqref{fgh1} into \eqref{Ti-Ti} and using \eqref{hypothesis KAM reduction of the remainder term} combined with Sobolev embeddings and \eqref{Conv-T2}  yield
				\begin{align}\label{Ti-TiX}
					\nonumber \|\Phi_{\infty}\rho\|_{q,s}^{\gamma ,\mathcal{O}}&\lesssim \|\rho\|_{q,s}^{\gamma,\mathcal{O}}+\varepsilon\gamma^{-2}\left(1+\| \mathfrak{I}_{0}\|_{q,s+\sigma_{4}}^{\gamma,\mathcal{O}}\right) \|\rho\|_{q,s_{0}}^{\gamma,\mathcal{O}}\\
					&\lesssim \|\rho\|_{q,s}^{\gamma,\mathcal{O}}+\varepsilon\gamma^{-2}\| \mathfrak{I}_{0}\|_{q,s+\sigma_{4}}^{\gamma,\mathcal{O}} \|\rho\|_{q,s_{0}}^{\gamma,\mathcal{O}}.
				\end{align}
				In a similar way to \eqref{Ti-Ti} we get  by  Lemma \ref{properties of Toeplitz in time operators}  combined with \eqref{Conv-U1} and \eqref{Conv-UL1}
				\begin{align*}
					\nonumber \big\|\big(\Phi_{\infty}-\textnormal{Id}\big)\rho\big\|_{q,s}^{\gamma ,\mathcal{O}}&\lesssim\| \Phi_{\infty}-\textnormal{Id}\|_{\textnormal{\tiny{O-d}},q,s_{0}}^{\gamma,\mathcal{O}}\|\rho\|_{q,s}^{\gamma,\mathcal{O}}+\| \Phi_{\infty}-\textnormal{Id}\|_{\textnormal{\tiny{O-d}},q,s}^{\gamma,\mathcal{O}}\|\rho\|_{q,s_{0}}^{\gamma,\mathcal{O}}\\
					\nonumber &\lesssim \delta_0(s_h)\|\rho\|_{q,s}^{\gamma,\mathcal{O}}+\Big(\delta_{0}(s_{h})+\delta_0\big(s+\tau_{2} q+\tau_{2}+1\big) \Big)\|\rho\|_{q,s_{0}}^{\gamma,\mathcal{O}}\\
					&\lesssim \delta_0(s_h)\|\rho\|_{q,s}^{\gamma,\mathcal{O}}+\delta_0\big(s+\tau_{2} q+\tau_{2}+1\big) \|\rho\|_{q,s_{0}}^{\gamma,\mathcal{O}}.
				\end{align*}
				Hence we find from  \eqref{fgh1} and \eqref{hypothesis KAM reduction of the remainder term} combined with Sobolev embeddings and \eqref{whab1}
				\begin{align}\label{Ti-TiP}
					\nonumber \big\|\big(\Phi_{\infty}-\textnormal{Id}\big)\rho\big\|_{q,s}^{\gamma ,\mathcal{O}}
					&\lesssim \big(\varepsilon\gamma^{-2}+\delta_0(s_h)\big)\|\rho\|_{q,s}^{\gamma,\mathcal{O}}+\varepsilon\gamma^{-2}\| \mathfrak{I}_{0}\|_{q,s+\sigma_{4}}^{\gamma,\mathcal{O}}  \|\rho\|_{q,s_{0}}^{\gamma,\mathcal{O}}\\
					&\lesssim \varepsilon\gamma^{-2}\|\rho\|_{q,s}^{\gamma,\mathcal{O}}+\varepsilon\gamma^{-2}\| \mathfrak{I}_{0}\|_{q,s+\sigma_{4}}^{\gamma,\mathcal{O}}  \|\rho\|_{q,s_{0}}^{\gamma,\mathcal{O}}.
				\end{align}
				The estimates $\Phi_{\infty}^{-1}$ and $ \Phi_{\infty}^{-1}-\widehat{\Phi}_{n}^{-1}$  follow from the same type pf arguments.\\
				\ding{226} In what follows  we plan to  study the asymptotic of the eigenvalues. Summing up  in $m$  the estimates  \eqref{Mahma1-R} and using Lemma \ref{lemma sum Nn}, we find
				\begin{align}\label{sum-R1}
					\nonumber \sum_{m=0}^\infty\|\mu_{j}^{m+1}-\mu_{j}^{m}\|_{q}^{\gamma ,\mathcal{O}}
					&\leqslant C\gamma\,\delta_{0}(s_{h})N_{0}^{\mu_{2}}\sum_{m=0}^\infty N_{m}^{-\mu_{2}}\\
					& \leqslant C\gamma\delta_{0}(s_{h}).
				\end{align}
				Thus  for each $j\in\mathbb{S}_0^c$  the sequence $(\mu_{j}^{m})_{m\in\mathbb{N}}$ converges in the space $W^{q,\infty,\gamma }(\mathcal{O},\mathbb{C})$ to an element denoted by  $\mu_{j}^{\infty}\in W^{q,\infty,\gamma }(\mathcal{O},\mathbb{C})$. 
				Moreover, for any $m\in\mathbb{N},$ we find in view of \eqref{Mahma1-R}
				\begin{align*}
					\|\mu_{j}^{\infty}-\mu_{j}^{m}\|_{q}^{\gamma ,\mathcal{O}}&\leqslant  \sum_{n=m}^\infty\|\mu_{j}^{n+1}-\mu_{j}^{n}\|_{q}^{\gamma ,\mathcal{O}}\\
					&\leqslant C\gamma\,\delta_{0}(s_{h})N_{0}^{\mu_{2}}\sum_{n=m}^\infty N_{n}^{-\mu_{2}}.
				\end{align*}
				Applying  Lemma \ref{lemma sum Nn}
				\begin{align}\label{Spect-TD1}
					\sup_{j\in \mathbb{S}_0^c}\|\mu_{j}^{\infty}-\mu_{j}^{m}\|_{q}^{\gamma ,\mathcal{O}}
					&\leqslant C \gamma\delta_{0}(s_{h})N_{0}^{\mu_{2}} N_{m}^{-\mu_{2}}.
				\end{align}
				Therefore, we deduce 
				\begin{align}\label{dekomp}
					\nonumber \mu_{j}^{\infty}&=\mu_{j}^{0}+\sum_{m=0}^\infty\big(\mu_{j}^{m+1}-\mu_{j}^{m}\big)\\
					&:=\mu_{j}^{0}+ r_{j}^{\infty},
				\end{align}
				where $(\mu_{j}^{0})$ is  described in Proposition \ref{projection in the normal directions} and takes the form  
				$$\mu_{j}^{0}(\lambda,\omega,i_{0})=\Omega_{j}(\lambda)+j\big(c_{i_{0}}(\lambda,\omega)-I_{1}(\lambda)K_{1}(\lambda)\big).
				$$
				Hence   \eqref{sum-R1}, \eqref{whab1} and \eqref{hypothesis KAM reduction of the remainder term} yield
				\begin{align*}
					\|r_{j}^{\infty}\|_{q}^{\gamma ,\mathcal{O}}
					&\leqslant C\, \gamma\,\delta_{0}(s_{h})\\
					&\leqslant C\, \varepsilon \gamma^{-1}
				\end{align*}
				and this gives the first  result in \eqref{estimate rjinfty}. Define the diagonal  operator $\mathscr{D}_{\infty}$ defined on the normal modes  by
				\begin{align}\label{Dinfty-op}
					\forall (l,j)\in\mathbb{Z}^d\times\mathbb{S}_0^c,\quad \mathscr{D}_{\infty} {\bf e}_{l,j}=\ii\mu_{j}^{\infty}{\bf e}_{l,j}.
				\end{align}
				By the norm definition we obtain
				\[
				\|\mathscr{D}_{m}-\mathscr{D}_{\infty}\|_{\textnormal{\tiny{O-d}},q,s_{0}}^{\gamma ,\mathcal{O}}=\displaystyle\sup_{j\in\mathbb{S}_{0}^{c}}\|\mu_{j}^{m}-\mu_{j}^{\infty}\|_{q}^{\gamma ,\mathcal{O}},
				\]
				which gives by virtue  of \eqref{Spect-TD1}
				\begin{align}\label{Conv-Dinf}
					\|\mathscr{D}_{m}-\mathscr{D}_{\infty}\|_{\textnormal{\tiny{O-d}},q,s_{0}}^{\gamma ,\mathcal{O}}\leqslant C\, \gamma\,\delta_{0}(s_{h})N_{0}^{\mu_{2}} N_{m}^{-\mu_{2}}.
				\end{align}
				% Similarly, using \eqref{Spect-TDMM1} yields
				% \begin{align}\label{Conv-DinPPf}
				% \|\mathscr{D}_{m}-\mathscr{D}_{\infty}\|_{\textnormal{\tiny{O-d}},q,s_{0}}^{\gamma ,\mathcal{O}}\leqslant C\, \gamma\,\delta_{0}(\overline{s}_{2})N_{0}^{\overline{\mu}_{2}} \textcolor{red}{N_{m}^{-\overline{\mu}_{2}}}.
				% \end{align}
				%
				%
				\ding{226}  The next goal is to  prove that the Cantor set $\mathscr{O}_{\infty,n}^{\gamma,\tau_1,\tau_{2}}(i_{0})$ defined in Proposition \ref{reduction of the remainder term} satisfies
				$$\mathscr{O}_{\infty,n}^{\gamma,\tau_1,\tau_{2}}(i_{0})\subset\bigcap_{m=0}^{n+1}\mathscr{O}_{m}^{\gamma}=\mathscr{O}_{n+1}^{\gamma}.$$
				where  the intermediate Cantor sets are defined in \eqref{Cantor-SX}. For this aim we shall proceed by finite induction on $m$ with $n$ fixed. First, we get by construction 
				$\mathscr{O}_{\infty,n}^{\gamma,\tau_1,\tau_{2}}(i_{0})\subset\mathcal{O}=:\mathscr{O}_{0}^{\gamma }.$ 
				Now  assume  that $\mathscr{O}_{\infty,n}^{\gamma,\tau_{2}}(i_{0})\subset\mathscr{O}_{m}^{\gamma }$ for $m\leqslant n$ and let us check that 
				\begin{align}\label{Inc-H}
					\mathscr{O}_{\infty,n}^{\gamma,\tau_1,\tau_{2}}(i_{0})\subset\mathscr{O}_{m+1}^{\gamma}.
				\end{align}
				Let  $(\lambda,\omega)\in\mathscr{O}_{\infty,n}^{\gamma,\tau_1,\tau_{2}}(i_{0})$ and $(l,j,j_{0})\in\mathbb{Z}^{d }\times(\mathbb{S}_{0}^{c})^{2}$ such that $0\leqslant |l|\leqslant N_{m}$ and $(l,j)\neq(0,j_0).$ Then, the triangle inequality, \eqref{Spect-TD1}, \eqref{Conv-T2} and \eqref{Conv-P3} imply
				\begin{align*}
					|\omega\cdot l+\mu_{j}^{m}(\lambda,\omega)-\mu_{j_{0}}^{m}(\lambda,\omega)| & \geqslant  \displaystyle|\omega\cdot l+\mu_{j}^{\infty}(\lambda,\omega)-\mu_{j_{0}}^{\infty}(\lambda,\omega)|-2\sup_{j\in\mathbb{S}_{0}^{c}}\|\mu_{j}^{m}-\mu_{j}^{\infty}\|_{q}^{\gamma ,\mathcal{O}}\\
					& \geqslant  \displaystyle\tfrac{2\gamma\langle j-j_{0}\rangle}{\langle l\rangle^{\tau_{2}}}-2\gamma\delta_{0}(s_{h})N_{0}^{\mu_{2}} N_{m}^{-\mu_{2}}\\
					& \geqslant \displaystyle \tfrac{2\gamma\langle j-j_{0}\rangle}{\langle l\rangle^{\tau_{2}}}-2\gamma\varepsilon_{0} \langle l\rangle^{-\mu_2}\langle j-j_{0}\rangle.
				\end{align*}
				Thus  for $\varepsilon_0$ small enough and by  \eqref{Conv-T2}(implying that $\mu_2\geqslant \tau_2$) we get
				$$
				\big|\omega\cdot l+\mu_{j}^{m}(\lambda,\omega)-\mu_{j_{0}}^{m}(\lambda,\omega)\big| > \tfrac{\gamma\langle j-j_{0}\rangle}{\langle l\rangle^{\tau_{2}}}
				$$
				which shows that $(\lambda,\omega)\in\mathscr{O}_{m+1}^{\gamma}$ and therefore the inclusion \eqref{Inc-H} is satisfied.\\
				\ding{226} Next we shall discuss  the convergence of the sequence $\left(\mathscr{L}_{m}\right)_{m\in\mathbb{N}}$ introduced in \eqref{Op-Lm} towards  the diagonal operator  $\mathscr{L}_{\infty}:=\omega\cdot\partial_{\varphi}\Pi_{\mathbb{S}_0}^{\perp}+\mathscr{D}_{\infty},$ where $\mathscr{D}_\infty$ is detailed in \eqref{Dinfty-op}. Applying  \eqref{Conv-Dinf} and \eqref{hypothesis of induction for deltamprime}
				\begin{align}\label{FD-TU-1}
					\nonumber \|\mathscr{L}_{m}-\mathscr{L}_{\infty}\|_{\textnormal{\tiny{O-d}},q,s_{0}}^{\gamma ,\mathcal{O}}&\leqslant\|\mathscr{D}_{m}-\mathscr{D}_{\infty}\|_{\textnormal{\tiny{O-d}},q,s_{0}}^{\gamma ,\mathcal{O}}+\|\mathscr{R}_{m}\|_{\textnormal{\tiny{O-d}},q,s_{0}}^{\gamma ,\mathcal{O}}\\
					&\leqslant C\, \gamma\,\delta_{0}(s_{h})N_{0}^{\mu_{2}} N_{m}^{-\mu_{2}},
				\end{align}
				which gives in particular that
				\begin{align}\label{strong-Cv}
					\lim_{m\rightarrow\infty} \|\mathscr{L}_{m}-\mathscr{L}_{\infty}\|_{\textnormal{\tiny{O-d}},q,s_{0}}^{\gamma ,\mathcal{O}}=0.
				\end{align}
				{By virtue of  \eqref{Def-Phi} and \eqref{Op-Lm1} one gets 
					\begin{align*}
						\forall (\lambda,\omega)\in \mathscr{O}_{n+1}^{\gamma },\quad\widehat{\Phi}_{n}^{-1}\mathscr{L}_{0}\widehat{\Phi}_{n}&=\big(\omega\cdot\partial_{\varphi}+\mathscr{D}_{n+1}+\mathscr{R}_{n+1}\big)\Pi_{\mathbb{S}_0}^{\perp}\\
						&=\mathscr{L}_{\infty}+\big(\mathscr{D}_{n+1}-\mathscr{D}_{\infty}+\mathscr{R}_{n+1}\big)\Pi_{\mathbb{S}_0}^{\perp},
					\end{align*} 
					It follows that any $(\lambda,\omega)\in \mathscr{O}_{n+1}^{\gamma }$ 					\begin{align*}
						\Phi_{\infty}^{-1}\mathscr{L}_{0}\Phi_{\infty}&=\mathscr{L}_{\infty}+\big(\mathscr{D}_{n+1}-\mathscr{D}_{\infty}+\mathscr{R}_{n+1}\big)\Pi_{\mathbb{S}_0}^{\perp}\\
						&+\Phi_{\infty}^{-1}\mathscr{L}_{0}\left(\Phi_{\infty}-\widehat{\Phi}_{n}\right)+\left(\Phi_{\infty}^{-1}-\widehat{\Phi}_{n}^{-1}\right)\mathscr{L}_{0}\widehat{\Phi}_{n}\\
						&:=\mathscr{L}_{\infty}+\mathtt{E}_{n,1}^2+\mathtt{E}_{n,2}^2+\mathtt{E}_{n,3}^2:=\mathscr{L}_{\infty}+\mathtt{E}_{n}^2.
					\end{align*}
					For the    estimate $\mathtt{E}_{n,1}^2$ we use \eqref{Conv-Dinf} combined with \eqref{Def-Taou},
					\eqref{hypothesis of induction for deltamprime}, \eqref{whab1} and  \eqref{hypothesis KAM reduction of the remainder term}
					\begin{align}\label{i-En21}
						\|\mathtt{E}_{n,1}^2\|_{\textnormal{\tiny{O-d}},q,s_0}^{\gamma,\mathcal{O}}
						&\leqslant C\, \gamma\,\delta_{0}(s_{h})N_{0}^{\mu_{2}} N_{n+1}^{-\mu_{2}}\nonumber\\
						&\leqslant C\varepsilon\gamma^{-1}N_{0}^{{\mu}_{2}}N_{n+1}^{-\mu_{2}}.
					\end{align}
					%					Applying \eqref{Conv-Dinf}   we find for any $s\in[s_0,S],$
					%					\begin{align*}
					%						\|\mathtt{E}_{n,1}^2\|_{\textnormal{\tiny{O-d}},q,s}^{\gamma,\mathcal{O}}&\leqslant \|\mathscr{D}_{n+1}-\mathscr{D}_{\infty}\|_{\textnormal{\tiny{O-d}},q,s}^{\gamma,\mathcal{O}}+\gamma \delta_{n+1}(s)\\
					%						&\leqslant C\, \gamma\,\delta_{0}({s}_{h})N_{0}^{{\mu}_{2}} N_{n+1}^{-\mu_{2}}+\gamma \delta_{n+1}(s)\\
					%						&\leqslant C\, \varepsilon\gamma^{-1}N_{0}^{{\mu}_{2}} N_{n+1}^{-\mu_{2}}+\gamma \delta_{n+1}(s).
					%					\end{align*}
					%					Then using \eqref{uniform estimate of deltamsprime}  and  \eqref{whab1} yields
					%					\begin{align*}
					%						\|\mathtt{E}_{n,1}^2\|_{\textnormal{\tiny{O-d}},q,s}^{\gamma,\mathcal{O}}&\leqslant C\varepsilon\gamma^{-1}\left(1+\|\mathfrak{I}_{0}\|_{q,{s}+\sigma_{4}}^{\gamma,\mathcal{O}}\right)
					%					\end{align*}
					According to Lemma \ref{definition of projections for operators} with \eqref{i-En21} we obtain
					\begin{align*}
						\|\mathtt{E}_{n,1}^2\rho\|_{q,s_0}^{\gamma,\mathcal{O}}
						&\leqslant C\varepsilon\gamma^{-1}N_{0}^{{\mu}_{2}}N_{n+1}^{-\mu_{2}} \|\rho\|_{q,s_0}^{\gamma,\mathcal{O}}.
					\end{align*}
					%					Similarly, we get for any $s\in[s_0,S],$
					%					\begin{align*}
					%						\|\mathtt{E}_{n,1}^2\rho\|_{q,s}^{\gamma,\mathcal{O}}
					%						&\leqslant C\varepsilon\gamma^{-1}\left(1+\|\mathfrak{I}_{0}\|_{q,{s}+\sigma_{4}}^{\gamma,\mathcal{O}}\right) \|\rho\|_{q,s_0}^{\gamma,\mathcal{O}}+C\varepsilon\gamma^{-1}N_{0}^{{\mu}_{2}}N_{n+1}^{-\mu_{2}}\|\rho\|_{q,s}^{\gamma,\mathcal{O}}\\
					%						&\leqslant C\|\rho\|_{q,s}^{\gamma,\mathcal{O}}+ C\varepsilon\gamma^{-1}\|\mathfrak{I}_{0}\|_{q,{s}+\sigma_{4}}^{\gamma,\mathcal{O}}\|\rho\|_{q,s_0}^{\gamma,\mathcal{O}}.
					%					\end{align*}
					Now let us move to the estimates of  $\mathtt{E}_{n,2}^2$ and $ \mathtt{E}_{n,3}^2.$ They can be treated in a similar way. Therefore we shall restrict the discussion to the term $\mathtt{E}_{n,2}^2$. Using \eqref{estimate on Phiinfty and its inverse} yields
					\begin{align}\label{Taptap2}
						\|\mathtt{E}_{n,2}^2\rho\|_{q,s_0}^{\gamma ,\mathcal{O}}\lesssim\big\|\mathscr{L}_{0}\big(\Phi_{\infty}-\widehat{\Phi}_{n}\big)\rho\big\|_{q,s_0}^{\gamma,\mathcal{O}}+\varepsilon\gamma^{-2}\| \mathfrak{I}_{0}\|_{q,s_0+\sigma_{4}}^{\gamma,\mathcal{O}}\big\|\mathscr{L}_{0}\big(\Phi_{\infty}-\widehat{\Phi}_{n}\big)\rho\big\|_{q,s_{0}}^{\gamma,\mathcal{O}}.
					\end{align}
					Therefore we get from     \eqref{Taptap1} combined with \eqref{hypothesis KAM reduction of the remainder term} 
					\begin{align*}%\label{Taptap3}
						\|\mathtt{E}_{n,2}^2\rho\|_{q,s_0}^{\gamma ,\mathcal{O}}&\lesssim\big\|\mathscr{L}_{0}\big(\Phi_{\infty}-\widehat{\Phi}_{n}\big)\rho\big\|_{q,s_0}^{\gamma,\mathcal{O}}\\
						&\lesssim \big\|\big(\Phi_{\infty}-\widehat{\Phi}_{n}\big)\rho\big\|_{q,s_0+1}^{\gamma,\mathcal{O}}.
					\end{align*}
					Applying  \eqref{Conv-op-od} with   Lemma \ref{definition of projections for operators},  \eqref{hypothesis KAM reduction of the remainder term} and \eqref{whab1} allow to get
					\begin{align*}%\label{Taptap3}
						\|\mathtt{E}_{n,2}^2\rho\|_{q,s_0}^{\gamma ,\mathcal{O}}&\lesssim\big\|\Phi_{\infty}-\widehat{\Phi}_{n}\big\|_{\textnormal{\tiny{O-d}},q,s_0+1}^{\gamma,\mathcal{O}}\|\rho\|_{q,s_0+1}^{\gamma,\mathcal{O}}\\
						&\leqslant  C\,\delta_{0}(s_{h})N_{0}^{\mu_{2}}N_{m+1}^{-\mu_2}\|\rho\|_{q,s_0+1}^{\gamma,\mathcal{O}}\\
						&\leqslant  C\,\varepsilon\gamma^{-2}N_0^{\mu_{2}}N_{m+1}^{-\mu_2}\|\rho\|_{q,s_0+1}^{\gamma,\mathcal{O}}.
					\end{align*}
					Notice that for $\mathtt{E}_{n,3}^2$ we get the same estimate as the preceding one. Consequently,  putting together the foregoing estimates yields \eqref{Error-Est-2D}.
					%
					%\begin{equation}\label{Def-Taou}
					%\delta_{m}(s):=\gamma ^{-1}\|\mathscr{R}_{m}\|_{\textnormal{\tiny{O-d}},q,s}^{\gamma ,\mathcal{O}}.
					%\end{equation}
					%In what follows, we  shall prove by induction in $m\in\mathbb{N}$ that 
					%\begin{equation}\label{hypothesis of induction for deltamprime}
					%\forall\, m\in\mathbb{N},\, \forall s\in[s_{0},s_{1}],\,\delta_{m}(s)\leqslant \delta_{0}(s_{2})N_{0}^{\mu_{2}}N_{m}^{-\mu_{2}}\quad \mbox{ and }\quad \delta_{m}(s_{2})\leqslant\left(2-\frac{1}{m+1}\right)\delta_{0}(s_{2})
					%\end{equation}
					%Let $\rho\in H^{s_0+1}(\T^{d+1})$ then from the preceding identity combined with   the law products and \eqref{Conv-op-od} we find
					%$$
					%\forall\,\, (\lambda,\omega)\in\mathscr{O}_{\infty}^{\gamma,\tau_{2}}(i_{0}),\quad \lim_{m\to\infty} \big\|\big(\mathscr{L}_{m}-\Phi_{\infty}^{-1}\mathscr{L}_{0}\Phi_{\infty}\big)\rho\big\|_{H^{s_0}(\T^{d+1})}=0.
					%$$
					%This ensures that $\left(\mathscr{L}_{m}  \right)_{m\in\mathbb{N}}$ converges {\it weakly} to $\Phi_{\infty}^{-1}\mathscr{L}_{0}\Phi_{\infty}$. Hence from \eqref{strong-Cv} and the uniqueness of the limit we get
					%$$
					%\forall \, (\lambda,\omega)\in \mathscr{O}_{\infty}^{\gamma,\tau_{2}}(i_{0}),\quad \Phi_{\infty}^{-1}\mathscr{L}_{0}\Phi_{\infty}=\mathscr{L}_{\infty}.
					%$$
				}\\
				\ding{226} The goal now is to prove  \eqref{estimate rjinfty}. We set
				$$
				\widehat{\delta}_{m}(s):=\max\left(\gamma ^{-1}\|\partial_{\theta}\mathscr{R}_{m}\|_{\textnormal{\tiny{O-d}},q,s}^{\gamma ,\mathcal{O}},\gamma^{-1}\|\mathscr{R}_{m}\|_{\textnormal{\tiny{O-d}},q,s}^{\gamma,\mathcal{O}}\right).
				$$
				Then we shall prove by induction on $m\in\mathbb{N}$ that
				\begin{align}\label{Ind-Ty1}
					\widehat{\delta}_{m}(s_{0})\leqslant\widehat{\delta}_{0}(s_{h})N_{0}^{\mu_{2}}N_{m}^{-\mu_{2}}\quad \mbox{ and }\quad \widehat{\delta}_{m}(s_{h})\leqslant\left(2-\tfrac{1}{m+1}\right)\widehat{\delta}_{0}(s_{h}).
				\end{align}
				According to Sobolev embeddings, the property is trivially satisfied  for $m=0.$ Notice that from  \eqref{estimate mathscr R in off diagonal norm} and \eqref{hypothesis KAM reduction of the remainder term} one  gets
				\begin{align}\label{tip-op1}
					\nonumber\widehat{\delta}_{0}(s_{h}) 
					&\lesssim\varepsilon\gamma^{-2}\left(1+\|\mathfrak{I}_{0}\|_{q,s_{h}+\sigma_{4}}^{\gamma,\mathcal{O}}\right)\\
					& \lesssim\varepsilon\gamma^{-2}.
				\end{align}
				We ssume that \eqref{Ind-Ty1} is satisfied at the order $m$ and let us check it at the order $m+1$. Applying $\partial_\theta$ to the second identity in \eqref{Deco-T1} and using \eqref{phi-inverse1} we obtain the expression%{AP-W34}yields
				\begin{align*}
					\nonumber\partial_{\theta}\mathscr{R}_{m+1}  &=  \Phi_{m}^{-1}\Big(P_{N_{m}}^{\perp}\partial_{\theta}\mathscr{R}_{m}+\partial_{\theta}\mathscr{R}_{m}\Psi_{m}-\Psi_{m}\partial_{\theta}\big\lfloor P_{N_{m}}\mathscr{R}_{m}\big\rfloor -\big[\partial_{\theta},\Psi_{m}\big]\big\lfloor P_{N_{m}}\mathscr{R}_{m}\big\rfloor \Big)\\
					\nonumber & \quad+\big[\partial_{\theta},\Sigma_{m}\big]\Big(P_{N_{m}}^{\perp}\mathscr{R}_{m}+\mathscr{R}_{m}\Psi_{m}-\Psi_{m}\big\lfloor P_{N_{m}}\mathscr{R}_{m}\big\rfloor \Big).\\
					&:=\mathcal{U}_m^1+\mathcal{U}_m^2
				\end{align*}
				with
				$$
				\mathcal{U}_m^2=\big[\partial_{\theta},\Sigma_{m}\big]\Big(P_{N_{m}}^{\perp}\mathscr{R}_{m}+\mathscr{R}_{m}\Psi_{m}-\Psi_{m}\big\lfloor P_{N_{m}}\mathscr{R}_{m}\big\rfloor \Big).
				$$
				It is easy to check that  for any Toeplitz in time operator $T(\lambda,\omega)$, we have
				$$
				\big[\partial_{\theta},T(\lambda,\omega)\big]\mathbf{e}_{l_{0},j_{0}}=\ii\sum_{(l,j)\in\mathbb{Z}^{d+1 }}(j-j_{0})T_{j_{0}}^{j}(\lambda,\omega,l-l_{0})\mathbf{e}_{l,j},
				$$
				which implies using the norm definition  
				\begin{equation}\label{estimate commutator with partialtheta in off diagonal norm}
					\big\|[\partial_{\theta},T]\big\|_{\textnormal{\tiny{O-d}},q,s}^{\gamma ,\mathcal{O}}\leqslant \| T\|_{\textnormal{\tiny{O-d}},q,s+1}^{\gamma ,\mathcal{O}}.
				\end{equation}
				{Since $\Phi_{m}^{-1}=\hbox{Id}+\Sigma_{m}$, then applying  Lemma \ref{properties of Toeplitz in time operators},  we obtain successively for $S\geqslant \overline{s}\geqslant s\geqslant s_{0}$ 
					\begin{align}\label{RPP-1}
						\nonumber\|\mathcal{U}_m^1\|_{\textnormal{\tiny{O-d}},q,s}^{\gamma,\mathcal{O}}&\leqslant C\|\Sigma_{m}\|_{\textnormal{\tiny{O-d}},q,s}^{\gamma,\mathcal{O}}\left[\|\partial_{\theta}\mathscr{R}_{m}\|_{\textnormal{\tiny{O-d}},q,s_{0}}^{\gamma,\mathcal{O}}\left(1+\|\Psi_{m}\|_{\textnormal{\tiny{O-d}},q,s_{0}}^{\gamma,\mathcal{O}}\right)+\|\left[\partial_{\theta},\Psi_{m}\right]\|_{\textnormal{\tiny{O-d}},q,s_{0}}^{\gamma,\mathcal{O}}\|\mathscr{R}_{m}\|_{\textnormal{\tiny{O-d}},q,s_{0}}^{\gamma,\mathcal{O}}\right]\\
						&+C\|\Sigma_{m}\|_{\textnormal{\tiny{O-d}},q,s_{0}}^{\gamma,\mathcal{O}}\left[\|\partial_{\theta}\mathscr{R}_{m}\|_{\textnormal{\tiny{O-d}},q,s}^{\gamma,\mathcal{O}}\left(1+\|\Psi_{m}\|_{\textnormal{\tiny{O-d}},q,s_{0}}^{\gamma,\mathcal{O}}\right)+\|\partial_{\theta}\mathscr{R}_{m}\|_{\textnormal{\tiny{O-d}},q,s_{0}}^{\gamma,\mathcal{O}}\|\Psi_{m}\|_{\textnormal{\tiny{O-d}},q,s}^{\gamma,\mathcal{O}}\right.\\
						\nonumber&\left.+\|\left[\partial_{\theta},\Psi_{m}\right]\|_{\textnormal{\tiny{O-d}},q,s}^{\gamma,\mathcal{O}}\|\mathscr{R}_{m}\|_{\textnormal{\tiny{O-d}},q,s_{0}}^{\gamma,\mathcal{O}}+\|\left[\partial_{\theta},\Psi_{m}\right]\|_{\textnormal{\tiny{O-d}},q,s_{0}}^{\gamma,\mathcal{O}}\|\mathscr{R}_{m}\|_{\textnormal{\tiny{O-d}},q,s}^{\gamma,\mathcal{O}}\right]+\|P_{N_{m}}^{\perp}\partial_{\theta}\mathscr{R}_{m}\|_{\textnormal{\tiny{O-d}},q,s}^{\gamma,\mathcal{O}}
					\end{align}
					and 
					\begin{align}\label{RPP-2}
						\nonumber \|\mathcal{U}_m^2\|_{\textnormal{\tiny{O-d}},q,s}^{\gamma,\mathcal{O}}&\lesssim \|\left[\partial_{\theta},\Sigma_{m}\right]\|_{\textnormal{\tiny{O-d}},q,s_{0}}^{\gamma,\mathcal{O}}\left(\|\mathscr{R}_{m}\|_{\textnormal{\tiny{O-d}},q,s_{0}}^{\gamma,\mathcal{O}}\|\Psi_{m}\|_{\textnormal{\tiny{O-d}},q,s}^{\gamma,\mathcal{O}}+\|\mathscr{R}_{m}\|_{\textnormal{\tiny{O-d}},q,s}^{\gamma,\mathcal{O}}\left(1+\|\Psi_{m}\|_{\textnormal{\tiny{O-d}},q,s_{0}}^{\gamma,\mathcal{O}}\right)\right)\\
						&\quad+\|\left[\partial_{\theta},\Sigma_{m}\right]\|_{\textnormal{\tiny{O-d}},q,s}^{\gamma,\mathcal{O}}\|\mathscr{R}_{m}\|_{\textnormal{\tiny{O-d}},q,s_{0}}^{\gamma,\mathcal{O}}\left(1+\|\Psi_{m}\|_{\textnormal{\tiny{O-d}},q,s_{0}}^{\gamma,\mathcal{O}}\right).
					\end{align}
					By using \eqref{estimate commutator with partialtheta in off diagonal norm}, \eqref{link Psi and R} and Lemma \ref{properties of Toeplitz in time operators}, we obtain
					\begin{align*}
						\|\left[\partial_{\theta},\Psi_{m}\right]\|_{\textnormal{\tiny{O-d}},q,s}^{\gamma,\mathcal{O}}&\leqslant\|\Psi_{m}\|_{\textnormal{\tiny{O-d}},q,s+1}^{\gamma,\mathcal{O}}\\
						&\leqslant C\gamma^{-1}\|P_{N_{m}}\mathscr{R}_{m}\|_{\textnormal{\tiny{O-d}},q,s+\tau_{2}q+\tau_{2}+1}^{\gamma,\mathcal{O}}\\
						&\leqslant CN_{m}^{\tau_{2}q+\tau_{2}+1}\delta_{m}(s).
					\end{align*}
					Coming back to  \eqref{AP-W34}, we obtain
					\begin{align*}
						\|\left[\partial_{\theta},\Sigma_{m}\right]\|_{\textnormal{\tiny{O-d}},q,s}^{\gamma,\mathcal{O}}&\leqslant\|\Sigma_{m}\|_{\textnormal{\tiny{O-d}},q,s+1}^{\gamma,\mathcal{O}}\\
						&\leqslant CN_{m}^{\tau_{2}q+\tau_{2}+1}\delta_{m}(s).
					\end{align*}
					Then inserting the preceding estimates and  \eqref{link-Z3} into   \eqref{RPP-1}   we deduce  that
					\begin{align}\label{DidaM1}
						\forall\, S\geqslant\overline{s}\geqslant s\geqslant s_{0},\quad \widehat{\delta}_{m+1}(s)\leqslant N_{m}^{s-\overline{s}}\,\widehat{\delta}_{m}(\overline{s})+CN_{m}^{\tau_{2}q+\tau_{2}+1}\widehat{\delta}_{m}(s)\widehat{\delta}_{m}(s_{0}).
					\end{align}
					In particular,  for $s=s_{0}$ we get by the induction assumption \eqref{Ind-Ty1}, %\eqref{hypothesis of induction for deltamprime} and  :
					\begin{align*}
						\widehat{\delta}_{m+1}(s_{0})&\leqslant  N_{m}^{s_{0}-s_{h}}\widehat{\delta}_{m}(s_{h})+CN_{m}^{\tau_{2}q+\tau_{2}+1}\left(\widehat{\delta}_{m}(s_{0})\right)^{2}\\
						& \leqslant \left(2-\tfrac{1}{m+1}\right)\widehat{\delta}_{0}(s_{h})N_{m}^{s_{0}-s_{h}}+CN_{0}^{2\mu_{2}}N_{m}^{\tau_{2}q+\tau_{2}+1-2\mu_{2}}\left(\widehat{\delta}_{0}(s_{h})\right)^{2}\\
						& \leqslant \widehat{\delta}_{0}(s_{h})\left(2N_{m}^{s_{0}-s_{h}}+CN_{0}^{2\mu_{2}}N_{m}^{\tau_{2}q+\tau_{2}+1-2\mu_{2}}\widehat{\delta}_{0}(s_{h})\right).
					\end{align*}
					If we fix $s_{2}$ and $\mu_{2}$ such that
					\begin{equation}\label{Conv-T3}
						N_{m}^{s_{0}-s_{h}}\leqslant\tfrac{1}{4}N_{0}^{\mu_{2}}N_{m+1}^{-\mu_{2}}\quad\textnormal{and}\quad CN_{0}^{2\mu_{2}}N_{m}^{\tau_{2}q+\tau_{2}+1-2\mu_{2}}\widehat{\delta}_{0}(s_{h})\leqslant\tfrac{1}{2}N_{0}^{\mu_{2}}N_{m+1}^{-\mu_{2}},
					\end{equation}
					then we find
					$$\widehat{\delta}_{m+1}(s_{0})\leqslant\widehat{\delta}_{0}(s_{h})N_{0}^{\mu_{2}}N_{m+1}^{-\mu_{2}}.$$
					Notice that \eqref{Conv-T2} implies in particular
					$$s_h\geqslant \tfrac{3}{2}\mu_{2}+s_0+1\quad\textnormal{and}\quad\mu_{2}\geqslant 2(\tau_{2}q+\tau_{2}+1)+1.$$
					Hence,  using \eqref{definition of Nm}, we see that the assumptions  of \eqref{Conv-T3} hold true provided that 
					\begin{equation*}
						4N_{0}^{-\mu_2}\leqslant 1\qquad\hbox{and}\qquad 
						2C\widehat{\delta}_{0}(s_{h})\leqslant  N_{0}^{-\mu_{2}}.
					\end{equation*}
					Remark  that these conditions are satisfied thanks to \eqref{Cond1}, \eqref{tip-op1} and \eqref{hypothesis KAM reduction of the remainder term}. Now, we turn to the proof of the second estimate in \eqref{Ind-Ty1}. By \eqref{DidaM1} and \eqref{Ind-Ty1}
					\begin{align*}
						\widehat{\delta}_{m+1}(s_{h})&\leqslant\widehat{\delta}_{m}(s_{h})+CN_{m}^{\tau_{2}q+\tau_{2}+1}\widehat{\delta}_{m}(s_{h})\widehat{\delta}_{m}(s_{0})\\
						&\leqslant\left(2-\tfrac{1}{m+1}\right)\widehat{\delta}_{0}(s_{h})\left(1+CN_{0}^{\mu_{2}}N_{m}^{\tau_{2}q+\tau_{2}+1-\mu_{2}}\widehat{\delta}_{0}(s_{h})\right).
					\end{align*}
					Taking  the parameters $s_{2}$ and $\mu_{2}$ such that
					\begin{equation}\label{Conv-9}
						\left(2-\tfrac{1}{m+1}\right)\left(1+CN_{0}^{\mu_{2}}N_{m}^{\tau_{2}q+\tau_{2}+1-\mu_{2}}\widehat{\delta}_{0}(s_{h})\right)\leqslant 2-\tfrac{1}{m+2},
					\end{equation}
					then we obtain
					$$\widehat{\delta}_{m+1}(s_{h})\leqslant\left(2-\tfrac{1}{m+2}\right)\widehat{\delta}_{0}(s_{h}),$$
					which achieves the induction argument in  \eqref{Ind-Ty1}. Now observe that \eqref{Conv-9}  is quite similar to \eqref{Conv-4} using in particular $\mu_{2}\geqslant 2(\tau_2q+\tau_2)+1$ and one may proceed following the same lines. Next let us see how to get the estimate  \eqref{estimate rjinfty}.
					Recall that
					$$r_{j}^{\infty}=\sum_{m=0}^{\infty}r_{j}^{m}\quad\quad \textnormal{with}\quad\quad r_{j}^{m}=\big\langle P_{N_{m}}\mathscr{R}_{m}\mathbf{e}_{0,j},\mathbf{e}_{0,j}\big\rangle_{L^{2}(\mathbb{T}^{d +1})}.$$
					Then it is clear that
					$$
					\big\langle P_{N_{m}}\mathscr{R}_{m}\mathbf{e}_{0,j},\mathbf{e}_{0,j}\big\rangle_{L^{2}(\mathbb{T}^{d +1},\mathbb{C})}=\tfrac{\ii}{j}\big\langle P_{N_{m}}\mathscr{R}_{m}\mathbf{e}_{0,j},\partial_{\theta}\mathbf{e}_{0,j}\big\rangle_{L^{2}(\mathbb{T}^{d +1})}.
					$$
					Therefore   integration by parts leads to
					$$\big\langle P_{N_{m}}\mathscr{R}_{m}\mathbf{e}_{0,j},\partial_{\theta}\mathbf{e}_{0,j}\big\rangle_{L^{2}(\mathbb{T}^{d +1},\mathbb{C})}=-\big\langle P_{N_{m}}\partial_{\theta}\mathscr{R}_{m}\mathbf{e}_{0,j},\mathbf{e}_{0,j}\big\rangle_{L^{2}(\mathbb{T}^{d +1})}.$$
					Using a duality argument $H^{s_{0}}-H^{-s_{0}}$ combined with   Lemma \ref{properties of Toeplitz in time operators} and  \eqref{Ind-Ty1}, we obtain
					\begin{align*}
						\|\langle P_{N_{m}}\partial_{\theta}\mathscr{R}_{m}\mathbf{e}_{0,j},\mathbf{e}_{0,j}\rangle_{L^{2}(\mathbb{T}^{d +1})}\|_{q}^{\gamma,\mathcal{O}}&\leqslant C\gamma\,\widehat{\delta}_{m}(s_{0})\\
						&\leqslant C\gamma\,\widehat{\delta}_{0}(s_{h})N_{0}^{\mu_{2}}N_{m}^{-\mu_{2}}.
					\end{align*}
					Putting together the preceding estimates  with \eqref{tip-op1} and Lemma \ref{lemma sum Nn} yields
					\begin{align*}
						\| r_{j}^{\infty}\|_{q}^{\gamma ,\mathcal{O}}&\lesssim\gamma\,\displaystyle|j|^{-1}\widehat{\delta}_{0}(s_{h})N_{0}^{\mu_{2}}\sum_{m=0}^{\infty}N_{m}^{-\mu_{2}}\\
						& \lesssim\displaystyle|j|^{-1}\varepsilon\gamma^{-1}.
					\end{align*}
					This achieves  the proof of \eqref{estimate rjinfty}.\\
					%{\bf{2.}}]
					%
					{\bf{(ii)}} We shall now work with   fixed  values (minimal) of $\mu_{2}$ and $s_h$ denoted respectively by $\mu_{c}$ and $s_c$, namely
					\begin{equation}\label{musc}
						\mu_{c}:=\overline{\mu_{2}}+2\tau_{2}q+2\tau_{2}\quad\textnormal{and}\quad s_c:=\frac{3}{2}\mu_c+\overline{s}_l+1=\overline{s}_h+4\tau_{2}q+4\tau_{2}.
					\end{equation}
					From \eqref{Deco-T1} and \eqref{phi-inverse1}, we can write
					$$\mathscr{R}_{m+1}=\left(\textnormal{Id}+\Sigma_{m}\right)U_m,$$
					where
					\begin{align}\label{bkl1}
						U_{m}:=P_{N_{m}}^{\perp}\mathscr{R}_{m}+\mathscr{R}_{m}\Psi_{m}-\Psi_{m}\big\lfloor P_{N_{m}}\mathscr{R}_{m}\big\rfloor.
					\end{align}
					After straightforward computations, we get
					\begin{align}\label{Gall-01}
						\nonumber\Delta_{12}U_{m}=&P_{N_{m}}^{\perp}\Delta_{12}\mathscr{R}_{m}+(\Delta_{12}\mathscr{R}_{m})(\Psi_{m})_{r_{1}}+(\mathscr{R}_{m})_{r_{2}}(\Delta_{12}\Psi_{m})\\
						&-(\Delta_{12}\Psi_{m})\big\lfloor P_{N_{m}}(\mathscr{R}_{m})_{r_1}\big\rfloor-(\Psi_{m})_{r_{2}}\big\lfloor P_{N_{m}}\Delta_{12}\mathscr{R}_{m}\big\rfloor
					\end{align}
					and
					\begin{align}\label{gui1}
						\Delta_{12}\mathscr{R}_{m+1}=\Delta_{12}U_{m}+(\Delta_{12}\Sigma_{m})(U_{m})_{r_{1}}+(\Sigma_{m})_{r_{2}}\Delta_{12}U_{m}.
					\end{align}
					We have used the notation $(f)_r=f(r).$ 
					Elementary manipulations  based on \eqref{phi-inverse1} give
					$$
					\Delta_{12}\Sigma_{m}=\Delta_{12}\Phi_{m}^{-1}=-(\Phi_{m}^{-1})_{r_{2}}(\Delta_{12}\Psi_{m})(\Phi_{m}^{-1})_{r_{1}}.$$
					The law products of Lemma \ref{properties of Toeplitz in time operators} together with  \eqref{link-ZP2} and  \eqref{hypothesis KAM reduction of 
						the remainder term} imply
					\begin{align}\label{link-ZP3}
						\forall \, s\in[s_0,s_c],\quad   \|\Delta_{12}\Sigma_{m}\|_{\textnormal{\tiny{O-d}},q,s}^{\gamma ,\mathcal{O}}&\lesssim 
						\|\Delta_{12}\Psi_{m}\|_{\textnormal{\tiny{O-d}},q,s}^{\gamma ,\mathcal{O}}.
					\end{align}
					Using once again  the law products of Lemma \ref{properties of Toeplitz in time operators}, \eqref{link-ZP3} and  \eqref{gui1} we obtain
					\begin{align}\label{gui2}
						\nonumber \|\Delta_{12}\mathscr{R}_{m+1}\|_{\textnormal{\tiny{O-d}},q,s_0}^{\gamma ,\mathcal{O}}&\leqslant  \|\Delta_{12}U_{m}\|_{\textnormal{\tiny{O-d}},q,s_0}^{\gamma ,\mathcal{O}}+\|(\Delta_{12}\Psi_{m})\|_{\textnormal{\tiny{O-d}},q,s_0}^{\gamma ,\mathcal{O}}\|(U_{m})_{r_{1}}\|_{\textnormal{\tiny{O-d}},q,s_0}^{\gamma ,\mathcal{O}}\\
						&\quad+\|(\Sigma_{m})_{r_{2}}\|_{\textnormal{\tiny{O-d}},q,s_0}^{\gamma ,\mathcal{O}}\|\Delta_{12}U_{m}\|_{\textnormal{\tiny{O-d}},q,s_0}^{\gamma ,\mathcal{O}}
					\end{align}
					and
					\begin{align}\label{gui3}
						\nonumber \|\Delta_{12}\mathscr{R}_{m+1}\|_{\textnormal{\tiny{O-d}},q,s_c}^{\gamma ,\mathcal{O}}&\leqslant  \|\Delta_{12}U_{m}\|_{\textnormal{\tiny{O-d}},q,s_c}^{\gamma ,\mathcal{O}}+\|(\Delta_{12}\Psi_{m})\|_{\textnormal{\tiny{O-d}},q,s_0}^{\gamma ,\mathcal{O}}\|(U_{m})_{r_{1}}\|_{\textnormal{\tiny{O-d}},q,s_c}^{\gamma ,\mathcal{O}}\\
						\nonumber &\quad+\|(\Delta_{12}\Psi_{m})\|_{\textnormal{\tiny{O-d}},q,s_c}^{\gamma ,\mathcal{O}}\|(U_{m})_{r_{1}}\|_{\textnormal{\tiny{O-d}},q,s_0}^{\gamma ,\mathcal{O}}+\|(\Sigma_{m})_{r_{2}}\|_{\textnormal{\tiny{O-d}},q,s_0}^{\gamma ,\mathcal{O}}\|\Delta_{12}U_{m}\|_{\textnormal{\tiny{O-d}},q,s_c}^{\gamma ,\mathcal{O}}\\
						&\quad+ \|(\Sigma_{m})_{r_{2}}\|_{\textnormal{\tiny{O-d}},q,s_c}^{\gamma ,\mathcal{O}}\|\Delta_{12}U_{m}\|_{\textnormal{\tiny{O-d}},q,s_0}^{\gamma ,\mathcal{O}}.
					\end{align}
					For the  estimate $(U_m)_{r_1}$ (to alleviate the notation we shall remove in this part  remove the subscript $r_1$) described by \eqref{bkl1} we use the law products leading to
					\begin{align}\label{bkl3}
						\|U_{m}\|_{\textnormal{\tiny{O-d}},q,s_0}^{\gamma ,\mathcal{O}}&\leqslant \|\mathscr{R}_{m}\|_{\textnormal{\tiny{O-d}},q,s_0}^{\gamma ,\mathcal{O}}+\|\mathscr{R}_{m}\|_{\textnormal{\tiny{O-d}},q,s_0}^{\gamma ,\mathcal{O}} \|\Psi_m\|_{\textnormal{\tiny{O-d}},q,s_0}^{\gamma ,\mathcal{O}}
					\end{align}
					and
					\begin{align}\label{bkl4}
						\nonumber \|U_{m}\|_{\textnormal{\tiny{O-d}},q,s_c}^{\gamma ,\mathcal{O}}&\leqslant \|\mathscr{R}_{m}\|_{\textnormal{\tiny{O-d}},q,s_c}^{\gamma ,\mathcal{O}}+\|\mathscr{R}_{m}\|_{\textnormal{\tiny{O-d}},q,s_0}^{\gamma ,\mathcal{O}} \|\Psi_m\|_{\textnormal{\tiny{O-d}},q,s_c}^{\gamma ,\mathcal{O}}\\
						&+\|\mathscr{R}_{m}\|_{\textnormal{\tiny{O-d}},q,s_c}^{\gamma ,\mathcal{O}} \|\Psi_m\|_{\textnormal{\tiny{O-d}},q,s_0}^{\gamma ,\mathcal{O}}.
					\end{align}
					By \eqref{hypothesis of induction for deltamprime},\eqref{whab1}  and \eqref{link-Z3} together with \eqref{bkl3} we infer
					\begin{align}\label{bkl03}
						\|U_{m}\|_{\textnormal{\tiny{O-d}},q,s_0}^{\gamma ,\mathcal{O}}&\leqslant C\varepsilon \gamma^{-1} N_0^{\mu_c} N_m^{-\mu_c}.\end{align}
					Putting together 
					the first estimate of \eqref{link-Z1}, \eqref{hypothesis of induction for deltamprime} and \eqref{whab1} we deduce that
					\begin{align}\label{bahya-11}
						\nonumber \max_{j=1,2}\|(\Psi_{m})_{r_{j}}\|_{\textnormal{\tiny{O-d}},q,s_c}^{\gamma ,\mathcal{O}}&\lesssim N_m^{\tau_2q+\tau_{2}}\delta_m(s_c)\\
						&\lesssim \varepsilon \gamma^{-2}N_m^{\tau_2q+\tau_{2}}.
					\end{align}
					Hence we get in view of   \eqref{bkl4},  \eqref{hypothesis of induction for deltamprime} and \eqref{hypothesis KAM reduction of the remainder term}
					\begin{align}\label{bkl04}
						\|U_{m}\|_{\textnormal{\tiny{O-d}},q,s_c}^{\gamma ,\mathcal{O}}&\leqslant C\varepsilon \gamma^{-1}.
					\end{align}
					Plugging  \eqref{bkl03} and \eqref{bkl04} into \eqref{gui3} implies
					\begin{align}\label{guiM3}
						\nonumber \|\Delta_{12}\mathscr{R}_{m+1}\|_{\textnormal{\tiny{O-d}},q,s_c}^{\gamma ,\mathcal{O}}&\leqslant  \|\Delta_{12}U_{m}\|_{\textnormal{\tiny{O-d}},q,s_c}^{\gamma ,\mathcal{O}}+C\varepsilon \gamma^{-1}\|\Delta_{12}\Psi_{m}\|_{\textnormal{\tiny{O-d}},q,s_0}^{\gamma ,\mathcal{O}}\\
						\nonumber &+C\varepsilon \gamma^{-1} N_0^{\mu_c} N_m^{-\mu_c}\|\Delta_{12}\Psi_{m}\|_{\textnormal{\tiny{O-d}},q,s_c}^{\gamma ,\mathcal{O}}+\|(\Sigma_{m})_{r_{2}}\|_{\textnormal{\tiny{O-d}},q,s_0}^{\gamma ,\mathcal{O}}\|\Delta_{12}U_{m}\|_{\textnormal{\tiny{O-d}},q,s_c}^{\gamma ,\mathcal{O}}\\
						&+ \|(\Sigma_{m})_{r_{2}}\|_{\textnormal{\tiny{O-d}},q,s_c}^{\gamma ,\mathcal{O}}\|\Delta_{12}U_{m}\|_{\textnormal{\tiny{O-d}},q,s_0}^{\gamma ,\mathcal{O}}.
					\end{align}
					Applying  \eqref{AP-W34} and \eqref{whab1} gives
					\begin{align}\label{tahya-b1}
						\|(\Sigma_{m})_{r_{2}}\|_{\textnormal{\tiny{O-d}},q,s_0}^{\gamma ,\mathcal{O}}\leqslant C\varepsilon \gamma^{-2}N_0^{\mu_c} N_m^{-\mu_c}\quad\hbox{and}\quad \|(\Sigma_{m})_{r_{2}}\|_{\textnormal{\tiny{O-d}},q,s_c}^{\gamma ,\mathcal{O}}\leqslant C\varepsilon \gamma^{-2}N_m^{\tau_2 q+\tau_{2}}.
					\end{align}
					Inserting \eqref{tahya-b1} into \eqref{guiM3} allows to get
					\begin{align}\label{guiM4}
						\nonumber \|\Delta_{12}\mathscr{R}_{m+1}\|_{\textnormal{\tiny{O-d}},q,s_c}^{\gamma ,\mathcal{O}}&\leqslant \big(1+C\varepsilon \gamma^{-2}N_0^{\mu_c} N_m^{-\mu_c}\big) \|\Delta_{12}U_{m}\|_{\textnormal{\tiny{O-d}},q,s_c}^{\gamma ,\mathcal{O}}+C\varepsilon \gamma^{-2}N_m^{\tau_2 q+\tau_{2}}\|\Delta_{12}U_{m}\|_{\textnormal{\tiny{O-d}},q,s_0}^{\gamma ,\mathcal{O}}\\
						&+C\varepsilon \gamma^{-1} N_0^{\mu_c} N_m^{-\mu_c}\|\Delta_{12}\Psi_{m}\|_{\textnormal{\tiny{O-d}},q,s_c}^{\gamma ,\mathcal{O}}+C\varepsilon \gamma^{-1}\|\Delta_{12}\Psi_{m}\|_{\textnormal{\tiny{O-d}},q,s_0}^{\gamma ,\mathcal{O}}.
					\end{align}
					In a similar way, by combining \eqref{bkl03}, \eqref{tahya-b1} with  \eqref{gui2} we find
					\begin{align}\label{guiP2}
						\nonumber \|\Delta_{12}\mathscr{R}_{m+1}\|_{\textnormal{\tiny{O-d}},q,s_0}^{\gamma ,\mathcal{O}}&\leqslant \big(1+C\varepsilon \gamma^{-2}N_0^{\mu_c} N_m^{-\mu_c}\big) \|\Delta_{12}U_{m}\|_{\textnormal{\tiny{O-d}},q,s_0}^{\gamma ,\mathcal{O}}\\
						&+C\varepsilon \gamma^{-1} N_0^{\mu_c} N_m^{-\mu_c}\|\Delta_{12}\Psi_{m}\|_{\textnormal{\tiny{O-d}},q,s_0}^{\gamma ,\mathcal{O}}.
					\end{align}
					From   \eqref{Gall-01}  and   the law products of Lemma \ref{properties of Toeplitz in time operators} we obtain $ \forall \, s\in[s_0,s_c],$
					\begin{align*}
						\|\Delta_{12}U_{m}\|_{\textnormal{\tiny{O-d}},q,s}^{\gamma ,\mathcal{O}}\leqslant&N_{m}^{s-s_c}   
						\|\Delta_{12}\mathscr{R}_{m}\|_{\textnormal{\tiny{O-d}},q,s_c}^{\gamma ,\mathcal{O}}
						+C\|\Delta_{12}\mathscr{R}_{m}\|_{\textnormal{\tiny{O-d}},q,s}^{\gamma ,\mathcal{O}}\max_{j=1,2}\|(\Psi_{m})_{r_{j}}\|_{\textnormal{\tiny{O-d}},q,s_0}^{\gamma ,\mathcal{O}}\\
						&+C\|\Delta_{12}\mathscr{R}_{m}\|_{\textnormal{\tiny{O-d}},q,s_0}^{\gamma ,\mathcal{O}}\max_{j=1,2}\|(\Psi_{m})_{r_{j}}\|_{\textnormal{\tiny{O-d}},q,s}^{\gamma ,\mathcal{O}}+C\|\Delta_{12}\Psi_{m}\|_{\textnormal{\tiny{O-d}},q,s}^{\gamma ,\mathcal{O}}\max_{j=1,2}\|(\mathscr{R}_m)_{r_{j}}\|_{\textnormal{\tiny{O-d}},q,s_0}^{\gamma ,\mathcal{O}}\\
						&+C\|\Delta_{12}\Psi_{m}\|_{\textnormal{\tiny{O-d}},q,s_0}^{\gamma ,\mathcal{O}}\max_{j=1,2}\|(\mathscr{R}_m)_{r_{j}}\|_{\textnormal{\tiny{O-d}},q,s}^{\gamma ,\mathcal{O}}.
					\end{align*}
					Combining the foregoing  estimate with \eqref{bahya-11},  \eqref{hypothesis of induction for deltamprime} and \eqref{link-Z3} yields
					\begin{align}\label{jima1}
						\nonumber \|\Delta_{12}U_{m}\|_{\textnormal{\tiny{O-d}},q,s_0}^{\gamma ,\mathcal{O}}\leqslant&N_{m}^{s_0-s_c}   
						\|\Delta_{12}\mathscr{R}_{m}\|_{\textnormal{\tiny{O-d}},q,s_c}^{\gamma ,\mathcal{O}}
						+C\varepsilon \gamma^{-2} N_0^{\mu_c} N_m^{-\mu_c}\|\Delta_{12}\mathscr{R}_{m}\|_{\textnormal{\tiny{O-d}},q,s_0}^{\gamma ,\mathcal{O}}\\
						&+C\varepsilon \gamma^{-1} N_0^{\mu_c} N_m^{-\mu_c}\|\Delta_{12}\Psi_{m}\|_{\textnormal{\tiny{O-d}},q,s_0}^{\gamma ,\mathcal{O}}\end{align}
					and
					\begin{align*}
						\|\Delta_{12}U_{m}\|_{\textnormal{\tiny{O-d}},q,s_c}^{\gamma ,\mathcal{O}}\leqslant &  
						\big(1+C\varepsilon \gamma^{-2} N_0^{\mu_c} N_m^{-\mu_c}\big)\|\Delta_{12}\mathscr{R}_{m}\|_{\textnormal{\tiny{O-d}},q,s_c}^{\gamma ,\mathcal{O}}+CN_m^{\tau_2q+\tau_{2}}\varepsilon \gamma^{-2}\|\Delta_{12}\mathscr{R}_{m}\|_{\textnormal{\tiny{O-d}},q,s_0}^{\gamma ,\mathcal{O}}\\
						&+C\varepsilon \gamma^{-1} N_0^{\mu_c} N_m^{-\mu_c}\|\Delta_{12}\Psi_{m}\|_{\textnormal{\tiny{O-d}},q,s_c}^{\gamma ,\mathcal{O}}+C\varepsilon \gamma^{-1}\|\Delta_{12}\Psi_{m}\|_{\textnormal{\tiny{O-d}},q,s_0}^{\gamma ,\mathcal{O}}.\end{align*}
					Putting together the preceding  estimate with \eqref{guiM4}, \eqref{guiP2}, \eqref{Conv-T2} and  \eqref{hypothesis KAM reduction of the remainder term} we deduce that
					\begin{align}\label{guiPMS}
						\nonumber \|\Delta_{12}\mathscr{R}_{m+1}\|_{\textnormal{\tiny{O-d}},q,s_c}^{\gamma ,\mathcal{O}}&\leqslant \Big(1+C\varepsilon 
						\gamma^{-2} N_0^{\mu_c} N_m^{-\mu_c}+C\varepsilon \gamma^{-2}N_{m}^{s_0-s_c+\tau_2q+\tau_2}\Big)\|\Delta_{12}\mathscr{R}_{m}\|_{\textnormal{\tiny{O-d}},q,s_c}^{\gamma ,\mathcal{O}}\\
						\nonumber &\quad+CN_m^{\tau_2q+\tau_{2}}\varepsilon \gamma^{-2}\|\Delta_{12}\mathscr{R}_{m}\|_{\textnormal{\tiny{O-d}},q,s_0}^{\gamma ,\mathcal{O}}+C\varepsilon \gamma^{-1} N_0^{\mu_c} N_m^{-\mu_c}\|\Delta_{12}\Psi_{m}\|_{\textnormal{\tiny{O-d}},q,s_c}^{\gamma ,\mathcal{O}}\\
						&\quad+C\varepsilon \gamma^{-1}\|\Delta_{12}\Psi_{m}\|_{\textnormal{\tiny{O-d}},q,s_0}^{\gamma ,\mathcal{O}}.
					\end{align}
					In a similar way, by making appeal to  \eqref{guiP2}, \eqref{jima1} and \eqref{hypothesis KAM reduction of the remainder term} we find
					\begin{align}\label{guiPM2}
						\nonumber \|\Delta_{12}\mathscr{R}_{m+1}\|_{\textnormal{\tiny{O-d}},q,s_0}^{\gamma ,\mathcal{O}}&\leqslant N_{m}^{s_0-s_c}   
						\|\Delta_{12}\mathscr{R}_{m}\|_{\textnormal{\tiny{O-d}},q,s_c}^{\gamma ,\mathcal{O}}
						+C\varepsilon \gamma^{-2} N_0^{\mu_c} N_m^{-\mu_c}\|\Delta_{12}\mathscr{R}_{m}\|_{\textnormal{\tiny{O-d}},q,s_0}^{\gamma ,\mathcal{O}}\\
						&+C\varepsilon \gamma^{-1} N_0^{\mu_c} N_m^{-\mu_c}\|\Delta_{12}\Psi_{m}\|_{\textnormal{\tiny{O-d}},q,s_0}^{\gamma ,\mathcal{O}}.
					\end{align}
					We shall now estimate $\Delta_{12}\Psi_{m}.$ Remark that
					$$\|\Delta_{12}\Psi_m\|_{\textnormal{\tiny{O-d}},q,s}^{\gamma,\mathcal{O}}=\sum_{\alpha\in\mathbb{N}^{d+1}\atop|\alpha|\leqslant q}\gamma^{\alpha}\sup_{(\lambda,\omega)\in\mathcal{O}}\left(\sum_{(l,k)\in\mathbb{Z}^{d+1}\atop|l|,|k|\leqslant N_m}\langle l,k\rangle^{2(s-|\alpha|)}\sup_{j\in\mathbb{Z}}\left|\partial_{\lambda,\omega}^{\alpha}\Delta_{12}(\Psi_m)_{j+k}^{j}(\lambda,\omega,l)\right|^{2}\right)^{\frac{1}{2}}.$$
					By virtue of  \eqref{coeff-psim}, we get
					$$(\Psi_m)_{j_0}^{j}(\lambda,\omega,l)=\left\lbrace\begin{array}{ll}
						-(\varrho_m)_{j_0}^{j}(\lambda,\omega,l)r_{j_0,m}^{j}(\lambda,\omega,l) & \textnormal{if }(l,j)\neq (0,j_0)\\
						0 & \textnormal{if }(l,j)= (0,j_0),
					\end{array}\right.$$
					where
					$$(\varrho_m)_{j_0}^j(\lambda,\omega,l):=\frac{\chi\big((\omega\cdot l+\mu_{j}^{m}(\lambda,\omega)-\mu_{j_0}^{m}(\lambda,\omega))(\gamma\langle j-j_0\rangle)^{-1}\langle l\rangle^{\tau_2}\big)}{\omega\cdot l+\mu_{j}^{m}(\lambda,\omega)-\mu_{j_0}^{m}(\lambda,\omega)}\cdot$$
					Recall from \eqref{coefficients of the remainder operator R}, that $\{\ii r_{j_0,m}^{j}(\lambda,\omega,l)\}$ are the Fourier coefficients of $P_{N_m}\mathscr{R}_m$, that is
					\begin{equation}\label{coef-R-2}
						\ii\,r_{j_0,m}^{j}(\lambda,\omega,l)=\langle P_{N_m}\mathscr{R}_m\mathbf{e}_{0,j_0},\mathbf{e}_{l,j}\rangle_{L^{2}(\mathbb{T}^{d+1})}.
					\end{equation}
					We can write for non-zero coefficients
					\begin{align*}\Delta_{12}(\Psi_m)_{j+k}^{j}(\lambda,\omega,l)&=\Delta_{12}(\varrho_m)_{j+k}^{j}(\lambda,\omega,l)\big(r_{j+k,m}^{j}\big)_{r_1}(\lambda,\omega,l)\\
						&+\big((\varrho_m)_{j+k}^{j}\big)_{r_2}(\lambda,\omega,l)\Delta_{12}r_{j+k,m}^{j}(\lambda,\omega,l).
					\end{align*}
					Hence, using Lemma \ref{Lem-lawprod}-(iv)
					\begin{align}\label{leib-dPsi}
						\forall q'\in\llbracket 0,q\rrbracket,\quad \|\Delta_{12}(\Psi_m)_{j+k}^{j}(\ast,l)\|_{q'}^{\gamma,\mathcal{O}}&\lesssim\|\Delta_{12}(\varrho_m)_{j+k}^{j}(\ast,l)\|_{q'}^{\gamma,\mathcal{O}}\max_{i\in\{1,2\}}\|\big(r_{j+k,m}^{j}\big)_{r_i}(\ast,l)\|_{q'}^{\gamma,\mathcal{O}}\nonumber\\
						&\quad+\max_{i\in\{1,2\}}\|\big((\varrho_m)_{j+k}^{j}\big)_{r_i}(\ast,l)\|_{q'}^{\gamma,\mathcal{O}}\|\Delta_{12}r_{j+k,m}^{j}(\ast,l)\|_{q'}^{\gamma,\mathcal{O}}.
					\end{align}
					From \eqref{coef-R-2}, we deduce
					$$\ii\Delta_{12}r_{j_0,m}^{j}(\lambda,\omega,l)=\langle P_{N_m}\Delta_{12}\mathscr{R}_m\mathbf{e}_{0,j_0},\mathbf{e}_{l,j}\rangle_{L^{2}(\mathbb{T}^{d+1})}.$$ 
					One can write
					$$(\varrho_m)_{j_0}^{j}(\lambda,\omega,l)=b_{l,j,j_{0},m}\widehat{\chi}\Big(b_{l,j,j_0,m}B_{l,j,j_0,m}(\lambda,\omega)\Big),$$
					with
					$$b_{l,j,j_0,m}:=(\gamma\langle j-j_0\rangle)^{-1}\langle l\rangle^{\tau_2},\quad B_{l,j,j_0,m}(\lambda,\omega):=\omega\cdot l+\mu_{j}^{m}(\lambda,\omega)-\mu_{j_0}^{m}(\lambda,\omega),\quad\widehat{\chi}(x)=\frac{\chi(x)}{x}\cdot$$
					Notice that from \eqref{e-mujm}, one obtains
					\begin{equation}\label{e-Bljj0m}
						\forall q'\in\llbracket 0,q\rrbracket,\quad\|B_{l,j,j_0,m}\|_{q'}^{\gamma,\mathcal{O}}\lesssim\langle l,j-j_0\rangle.
					\end{equation}
					In a similar way to \eqref{eg-k-1}, one gets from \eqref{e-mujm} 
					\begin{equation}\label{evrhom}
						\forall q'\in\llbracket 0,q\rrbracket,\quad \|(\varrho_m)_{j_0}^{j}(\ast,l)\|_{q'}^{\gamma,\mathcal{O}}\lesssim\gamma^{-(q'+1)}\langle l,j-j_0\rangle^{\tau_{2}q'+\tau_{2}+q'}.
					\end{equation}
					Using Taylor formula in a similar way to \eqref{Taylor-I2}, we find (to simplify the notation we remove the dependence in $(\lambda,\omega)$)
					$$\Delta_{12}(\varrho_m)_{j_0}^{j}(l)=b_{l,j,j_0,m}^{2}(\Delta_{12}B_{l,j,j_0,m})\int_{0}^{1}\widehat{\chi}'\Big(b_{l,j,j_0,m}\Big[(1-\tau)(B_{l,j,_0,m})_{r_1}+\tau(B_{l,j,_0,m})_{r_2}\Big]\Big)d\tau.$$
					We shall estimate $\Delta_{12}B_{l,j,j_0,m}.$ For that purpose, we use \eqref{Spect-T2} to write
					$$\mu_{j}^{m}=\mu_{j}^{0}+\sum_{n=0}^{m-1}\langle P_{N_n}\mathscr{R}_n\mathbf{e}_{0,j},\mathbf{e}_{0,j}\rangle_{L^{2}(\mathbb{T}^{d+1})}.$$
					We recall from Proposition \ref{projection in the normal directions} that
					$$\mu_{j}^{0}(\lambda,\omega,i_0)=\Omega_{j}(\lambda)+jr^{1}(\lambda,\omega,i_0),\quad r^{1}(\lambda,\omega,i_0)=c_{i_0}(\lambda,\omega)-V_0(\lambda).$$
					Therefore
					$$\Delta_{12}\mu_{j}^{m}=\Delta_{12}\mu_{j}^{0}+\sum_{n=0}^{m-1}\langle \Delta_{12}P_{N_n}\mathscr{R}_n\mathbf{e}_{0,j},\mathbf{e}_{0,j}\rangle_{L^{2}(\mathbb{T}^{d+1})}$$
					and
					\begin{align*}
						\Delta_{12}B_{l,j,j_0,m}&=\Delta_{12}\big(\mu_{j}^{m}-\mu_{j_0}^{m}\big)\\
						&=(j-j_0)\Delta_{12}c_{i}+\sum_{n=0}^{m-1}\langle \Delta_{12}P_{N_n}\mathscr{R}_n\mathbf{e}_{0,j},\mathbf{e}_{0,j}\rangle_{L^{2}(\mathbb{T}^{d+1})}\\
						&\quad-\sum_{n=0}^{m-1}\langle \Delta_{12}P_{N_n}\mathscr{R}_n\mathbf{e}_{0,j_{0}},\mathbf{e}_{0,j_{0}}\rangle_{L^{2}(\mathbb{T}^{d+1})}.
					\end{align*}
					Hence, using \eqref{difference ci}, one gets
					\begin{equation}\label{ed-Bljj0m}
						\forall q'\in\llbracket 0,q\rrbracket,\quad \|\Delta_{12}B_{l,j,j_0,m}\|_{q'}^{\gamma,\mathcal{O}}\lesssim\varepsilon|j-j_0|\|\Delta_{12}i\|_{q',\overline{s}_h+2}^{\gamma,\mathcal{O}}+\sum_{n=0}^{m-1}\|P_{N_n}\Delta_{12}\mathscr{R}_n\|_{\textnormal{\tiny{O-d}},q',s_0}^{\gamma,\mathcal{O}}.
					\end{equation}
					Then, one obtains from Lemma \ref{Lem-lawprod}-(vi), \eqref{e-Bljj0m} and \eqref{ed-Bljj0m} 
					\begin{align}\label{edvrom}
						\forall q'\in\llbracket 0,q\rrbracket, \quad \|\Delta_{12}(\varrho_m)_{j_0}^{j}(\ast,l)\|_{q'}^{\gamma,\mathcal{O}}&\lesssim\varepsilon\gamma^{-2-q'}\langle l,j-j_0\rangle^{\tau_2q'+2\tau_2+q'+1}\|\Delta_{12}i\|_{q',\overline{s}_h+2}^{\gamma,\mathcal{O}}\\
						&\quad+\gamma^{-2-q'}\langle l,j-j_0\rangle^{\tau_2q'+2\tau_2+q'}\sum_{n=0}^{m-1}\|P_{N_n}\Delta_{12}\mathscr{R}_n\|_{\textnormal{\tiny{O-d}},q',s_0}^{\gamma,\mathcal{O}}.\nonumber
					\end{align}
					Gathering \eqref{leib-dPsi}, \eqref{evrhom} and \eqref{edvrom} gives for all $q'\in\llbracket 0,q\rrbracket$,
					\begin{align*}
						\|\Delta_{12}(\Psi_m)_{j+k}^{j}(\ast,l)\|_{q'}^{\gamma,\mathcal{O}}&\lesssim\varepsilon\gamma^{-2-q'}\langle l,k\rangle^{\tau_2q'+2\tau_2+q'+1}\|\Delta_{12}i\|_{q',\overline{s}_h+2}^{\gamma,\mathcal{O}}\max_{i\in\{1,2\}}\|\big(r_{j+k,m}^{j}\big)_{r_i}(\ast,l)\|_{q'}^{\gamma,\mathcal{O}}\\
						&+\gamma^{-2-q'}\langle l,k\rangle^{\tau_2q'+2\tau_2+q'}\max_{i\in\{1,2\}}\|\big(r_{j+k,m}^{j}\big)_{r_i}(\ast,l)\|_{q'}^{\gamma,\mathcal{O}}\sum_{n=0}^{m-1}\|P_{N_n}\Delta_{12}\mathscr{R}_n\|_{\textnormal{\tiny{O-d}},q',s_0}^{\gamma,\mathcal{O}}\\
						&+\gamma^{-1-q'}\langle l,k\rangle^{\tau_{2}q'+\tau_{2}+q'}\|\Delta_{12}r_{j+k,m}^{j}(\ast,l)\|_{q'}^{\gamma,\mathcal{O}}.
					\end{align*}
					We deduce that for all $s\in[s_0,S],$
					\begin{align}\label{e-dPsim-od}
						\|\Delta_{12}\Psi_m\|_{\textnormal{\tiny{O-d}},q,s}^{\gamma,\mathcal{O}}&\lesssim\varepsilon\gamma^{-2-q}\|\Delta_{12}i\|_{q,\overline{s}_h+2}^{\gamma,\mathcal{O}}\|P_{N_m}\mathscr{R}_m\|_{\textnormal{\tiny{O-d}},q,s+\tau_{2}q+2\tau_{2}+1}^{\gamma,\mathcal{O}}\nonumber\\
						&\quad+\gamma^{-2-q}\|P_{N_m}\mathscr{R}_m\|_{\textnormal{\tiny{O-d}},q,s+\tau_{2}q+2\tau_{2}}^{\gamma,\mathcal{O}}\sum_{n=0}^{m-1}\|P_{N_n}\Delta_{12}\mathscr{R}_n\|_{\textnormal{\tiny{O-d}},q,s_0}^{\gamma,\mathcal{O}}\nonumber\\
						&\quad+\gamma^{-1-q}\|P_{N_m}\Delta_{12}\mathscr{R}_m\|_{\textnormal{\tiny{O-d}},q,s+\tau_{2}q+\tau_{2}}^{\gamma,\mathcal{O}}.
					\end{align}
					We set 
					\begin{equation}\label{rec-dp-1}
						\overline{\delta}_m(s)=\gamma^{-1}\|\Delta_{12}\mathscr{R}_m\|_{\textnormal{\tiny{O-d}},q,s}^{\gamma,\mathcal{O}}\quad\textnormal{and}\quad \varkappa_{m}(s):=\sum_{n=0}^{m-1}\overline{\delta}_n(s).
					\end{equation}
					Then, using \eqref{e-dPsim-od}, \eqref{Def-Taou} and \eqref{param}, we get
					\begin{align}\label{galeria-1}
						\|\Delta_{12}\Psi_m\|_{\textnormal{\tiny{O-d}},q,s_{0}}^{\gamma,\mathcal{O}}&\lesssim\varepsilon\gamma^{-1-q}N_m^{\tau_2}\|\Delta_{12}i\|_{q,\overline{s}_h+2}^{\gamma,\mathcal{O}}\delta_{m}(\overline{s}_l)+\gamma^{-q}N_m^{\tau_{2}}\delta_{m}(\overline{s}_l)\varkappa_{m}(s_0)\nonumber\\
						&\quad+\gamma^{-q}N_m^{\tau_{2}q+\tau_{2}}\overline{\delta}_m(s_0)
					\end{align}
					and
					\begin{align}\label{Galeria-1}
						\|\Delta_{12}\Psi_m\|_{\textnormal{\tiny{O-d}},q,s_{c}}^{\gamma,\mathcal{O}}&\lesssim\varepsilon\gamma^{-1-q}N_m^{\tau_2q+2\tau_{2}+1}\|\Delta_{12}i\|_{q,\overline{s}_h+2}^{\gamma,\mathcal{O}}\delta_{m}(s_c)+\gamma^{-q}N_m^{\tau_{2}q+2\tau_{2}}\delta_{m}(s_c)\varkappa_{m}(s_0)\nonumber\\
						&\quad +\gamma^{-q}N_m^{\tau_{2}q+\tau_{2}}\overline{\delta}_m(s_c).
					\end{align}
					According to \eqref{tip-op1}, one has 
					\begin{equation}\label{exp-hrrrt}
						\delta_{m}(\overline{s}_l)\lesssim\varepsilon\gamma^{-2}N_0^{\mu_{c}}N_m^{-\mu_{c}}\quad\textnormal{and}\quad\sup_{m\in\mathbb{N}}\delta_{m}(s_c)\lesssim\varepsilon\gamma^{-2}\lesssim 1.
					\end{equation}
					Putting together \eqref{exp-hrrrt} and \eqref{galeria-1} and using \eqref{hypothesis KAM reduction of the remainder term} yields
					\begin{align}\label{galeria-2}
						\|\Delta_{12}\Psi_m\|_{\textnormal{\tiny{O-d}},q,s_0}^{\gamma,\mathcal{O}}&\lesssim \varepsilon\gamma^{-1}N_m^{\tau_{2}-\mu_c}\|\Delta_{12}i\|_{q,\overline{s}_h+2}^{\gamma,\mathcal{O}}+N_m^{\tau_{2}-\mu_{c}}\varkappa_{m}(s_0)+\gamma^{-q}N_m^{\tau_{2}q+\tau_2}\overline{\delta}_m(s_0).
					\end{align}
					In a similar way, on gets by \eqref{exp-hrrrt}, \eqref{Galeria-1} and \eqref{hypothesis KAM reduction of the remainder term}
					\begin{align}\label{Galeria-2}
						\nonumber \|\Delta_{12}\Psi_m\|_{\textnormal{\tiny{O-d}},q,s_c}^{\gamma,\mathcal{O}}&\lesssim\varepsilon\gamma^{-1}N_{m}^{\tau_{2}q+2\tau_{2}+1}\|\Delta_{12}i\|_{q,\overline{s}_h+2}^{\gamma,\mathcal{O}}+N_m^{\tau_{2}q+2\tau_{2}}\varkappa_{m}(s_0)\\
						&\quad+\gamma^{-q}N_m^{\tau_{2}q+\tau_{2}}\overline{\delta}_m(s_c).
					\end{align}
					Plugging \eqref{galeria-2} into \eqref{guiPM2} yields by virtue of  \eqref{hypothesis of induction for deltamprime} and \eqref{hypothesis KAM reduction of the remainder term}
					\begin{align}\label{guiPM3}
						\nonumber \overline{\delta}_{m+1}(s_0)&\leqslant N_{m}^{s_0-s_c}   
						\overline{\delta}_{m}(s_c)
						+CN_m^{\tau_2q+\tau_{2}-\mu_c}\overline{\delta}_{m}(s_0)+CN_m^{\tau_2-2\mu_c}\varkappa_{m}(s_0)\\
						&\quad+C\varepsilon\gamma^{-1}N_m^{\tau_2-2\mu_c}\|\Delta_{12}i\|_{q,\overline{s}_h+2}^{\gamma,\mathcal{O}}.
					\end{align}
					Therefore, inserting \eqref{galeria-2}, \eqref{Galeria-2} into \eqref{guiPMS} and using \eqref{hypothesis KAM reduction of the remainder term} implies
					\begin{align}\label{guiPMSY}
						\nonumber \overline{\delta}_{m+1}(s_c)&\leqslant \Big(1+C N_m^{\tau_2q+\tau_{2}-\mu_c}+CN_{m}^{s_0-s_c+\tau_2q+\tau_{2}}\Big)\overline{\delta}_{m}(s_c)+C\varepsilon\gamma^{-2-q}N_m^{\tau_2q+\tau_{2}}\overline{\delta}_{m}(s_0)\\
						&\quad+C N_m^{\tau_2q+\tau_{2}-\mu_c}\varkappa_{m}(s_0)+C\varepsilon\gamma^{-1}N_m^{\tau_2 q+2\tau_{2}+1-\mu_c}\|\Delta_{12}i\|_{q,\overline{s}_h+2}^{\gamma,\mathcal{O}}.
					\end{align}
					Next, we intend to  prove by induction in $m\in\mathbb{N}$ that
					\begin{equation}\label{hypothesis of induction differences of mathscr Rm}
						\forall\, k\leqslant m,\quad \overline{\delta}_{k}(s_{0})\leqslant N_{0}^{\mu_{c}}N_{k}^{-\mu_{c}}\nu(s_c)\quad \mbox{ and }\quad \overline{\delta}_{k}(s_{c})\leqslant\left(2-\tfrac{1}{k+1}\right)\nu(s_c),
					\end{equation}
					with
					$$
					\nu(s):=\overline\delta_0(s)+\varepsilon\gamma^{-1}\|\Delta_{12}i\|_{q,\overline{s}_h+2}^{\gamma,\mathcal{O}}.
					$$
					The estimate  \eqref{hypothesis of induction differences of mathscr Rm} is obvious  for $m=0$ by Sobolev embeddings. Now let us 
					assume that the preceding property holds true at the order $m$ and let us check it at the order $m+1$. Thus by applying \eqref{rec-dp-1} and Lemma \ref{lemma sum Nn}, we get
					$$
					\sup_{m\in\mathbb{N}}\varkappa_{m}(s_0)\leqslant C\nu(s_c).
					$$
					Putting together  this estimate with the induction assumption,  \eqref{guiPM3}, \eqref{guiPMSY},\eqref{musc}  and \eqref{hypothesis KAM reduction of the remainder term} yields
					\begin{align}%\label{guiPM6}
						\nonumber\overline{\delta}_{m+1}(s_0)&\leqslant \left(2N_{m}^{s_0-s_c}+CN_0^{\mu_c} N_m^{\tau_2q+\tau_{2}-2\mu_c}\right)\nu(s_c)
					\end{align}
					and
					\begin{align*}
						\overline{\delta}_{m+1}(s_c)&\leqslant \Big(1+C N_m^{\tau_2q+\tau_2-\mu_c}+CN_{m}^{s_0-s_c+\tau_2q+\tau_{2}}\Big)\left(2-\frac{1}{m+1}\right)\nu(s_{c})\\
						&\quad+C\left(N_m^{\tau_2q+\tau_2-\mu_c}+N_m^{\tau_2 q+2\tau_{2}+1-\mu_c}\right)\nu(s_c).
					\end{align*}
					Since \eqref{musc} implies in particular
					$$\mu_{c}\geqslant2\tau_2q+2\tau_{2}+1\quad\textnormal{and}\quad s_c\geqslant\tfrac{3}{2}\mu_{c}+s_0+\tau_{2}q+\tau_{2}+1,$$
					then proceeding similarly to the proof of \eqref{hypothesis of induction for deltamprime}, we conclude that
					\begin{align*}
						\overline{\delta}_{m+1}(s_0)&\leqslant N_0^{\mu_c}N_{m+1}^{-\mu_c}\nu(s_c)\quad\hbox{and}\quad \overline{\delta}_{m+1}(s_c)\leqslant \left(2-\tfrac{1}{m+2}\right)\nu(s_c),
					\end{align*}
					which achieves the induction. The next target is to estimate $\Delta_{12}r_{j}^{\infty}$. Then  similarly to  \eqref{dual-1X} we obtain  through  
					a duality argument, Lemma \ref{properties of Toeplitz in time operators},  \eqref{hypothesis of induction differences of mathscr Rm} and Lemma \ref{lemma sum Nn} 
					\begin{align*}
						\|\Delta_{12}r_{j}^{\infty}\|_{q}^{\gamma,\mathcal{O}}&\leqslant\sum_{m=0}^{\infty}\Big\|\langle P_{N_{m}}\Delta_{12}\mathscr{R}_{m}\mathbf{e}_{0,j},\mathbf{e}_{0,j}\rangle_{L^{2}(\mathbb{T}^{d+1})}\Big\|_{q}^{\gamma,\mathcal{O}}\\
						&\lesssim \sum_{m=0}^{\infty}\big\|\Delta_{12}\mathscr{R}_{m}\big\|_{\textnormal{\tiny{O-d}},q,s_0}^{\gamma,\mathcal{O}}\\
						&\lesssim \gamma \nu(s_c)\sum_{m=0}^{\infty}N_0^{\mu_c}N_{m}^{-\mu_c}\\
						&\leqslant C \gamma\,\nu(s_c).
					\end{align*}
					From the particular value of $\mathtt{p}$ in \eqref{choice ttp-1}, we infer
					\begin{equation}\label{choise ttp-2}
						s_c=\overline{s}_h+4\tau_{2}q+4\tau_{2}=\overline{s}_h+\mathtt{p}.
					\end{equation} 
					Then, applying \eqref{estimate differences mathscr R in off diagonal norm} we obtain
					\begin{align*}
						\|\Delta_{12}r_{j}^{\infty}\|_{q}^{\gamma,\mathcal{O}}
						&\leqslant C \gamma\,\nu(\overline{s}_h+4\tau_{2}q+4\tau_{2})
						\\
						&\leqslant C\varepsilon  \gamma^{-1}\,\|\Delta_{12} i\|_{q, \overline{s}_h+\sigma_{4}}^{\gamma,\mathcal{O}}.
					\end{align*}
					Finally, combining the previous estimate with   \eqref{dekomp} and \eqref{differences mu0} we deduce 
					\begin{align*}
						\forall j\in \mathbb{S}_0^c,\quad \|\Delta_{12}\mu_{j}^{\infty}\|_{q}^{\gamma,\mathcal{O}}&\lesssim\|\Delta_{12}\mu_{j}^{0}\|_{q}^{\gamma,\mathcal{O}}+\|\Delta_{12}r_{j}^{\infty}\|_{q}^{\gamma,\mathcal{O}}\\
						&\lesssim\varepsilon\gamma^{-1}|j|\|\Delta_{12}i\|_{q,\overline{s}_h+\sigma_{4}}^{\gamma,\mathcal{O}}.
					\end{align*}
					The  achieves the proof of Proposition \ref{reduction of the remainder term}.}
			\end{proof}
			%
			%\subsubsection{Conclusion : Estimates on the inverse operator in the normal directions}
			%
			\subsection{Approximate inverse  in the normal directions}
			In this section we plan  to construct an approximate  right inverse in the normal directions for the linearized operator $\widehat{\mathcal{L}}_{\omega}$ defined in \eqref{Lomega def} when the parameters are restricted in a Cantor like set.  Our  main result is the following.
			\begin{prop}\label{inversion of the linearized operator in the normal directions}
				Let $(\gamma,q,d,\tau_{1},s_{0},S)$ satisfying \eqref{initial parameter condition} and \eqref{setting tau1 and tau2}. There exists $\sigma:=\sigma(\tau_1,\tau_2,q,d)\geqslant\sigma_{4}$ such that if 
				\begin{equation}\label{small-1}
\varepsilon\gamma^{-2-q}N_0^{\mu_2}\leqslant\varepsilon_0\quad\textnormal{and}\quad\|\mathfrak{I}_0\|_{q,s_h+\sigma}^{\gamma,\mathcal{O}}\leqslant 1,
\end{equation}
then the following assertions hold true.
\begin{enumerate}[label=(\roman*)]
\item Consider the operator  $\mathscr{L}_{\infty}$ defined in Proposition $\ref{reduction of the remainder term},$ then there exists a family of  linear operators $\big(\mathtt{T}_n\big)_{n\in\mathbb{N}}$  defined in $\mathcal{O}$ satisfying the estimate
					$$
					\forall s\in[s_0,S],\quad \sup_{n\in\mathbb{N}}\|\mathtt{T}_{n}\rho\|_{q,s}^{\gamma ,\mathcal{O}}\lesssim \gamma ^{-1}\|\rho\|_{s+\tau_{1}q+\tau_{1}}^{\gamma ,\mathcal{O}}$$
					and such that for any $n\in\mathbb{N}$, in the Cantor set
					$$\Lambda_{\infty,n}^{\gamma,\tau_{1}}(i_{0})=\bigcap_{(l,j)\in\mathbb{Z}^{d }\times\mathbb{S}_{0}^{c}\atop |l|\leqslant N_{n}}\left\lbrace(\lambda,\omega)\in\mathcal{O}\quad\textnormal{s.t.}\quad\left|\omega\cdot l+\mu_{j}^{\infty}(\lambda,\omega,i_{0})\right|>\tfrac{\gamma \langle j\rangle }{\langle l\rangle^{\tau_{1}}}\right\rbrace,
					$$
					we have
					$$
					\mathscr{L}_{\infty}\mathtt{T}_n=\textnormal{Id}+\mathtt{E}^3_n,
					$$
					with 
					$$
					\forall s_{0}\leqslant s\leqslant\overline{s}\leqslant S, \quad \|\mathtt{E}^3_{n}\rho\|_{q,s}^{\gamma ,\mathcal{O}} \lesssim 
					N_n^{s-\overline{s}}\gamma^{-1}\|\rho\|_{\overline{s}+1+\tau_{1}q+\tau_{1}}^{\gamma ,\mathcal{O}}.
					$$
					\item 
					There exists a family of linear  operators $\big(\mathtt{T}_{\omega,n}\big)_{n\in\mathbb{N}}$ satisfying
					\begin{equation}\label{estimate mathcalTomega}
						\forall \, s\in\,[ s_0, S],\quad\sup_{n\in\mathbb{N}}\|\mathtt{T}_{\omega,n}\rho\|_{q,s}^{\gamma ,\mathcal{O}}\lesssim\gamma^{-1}\left(\|\rho\|_{q,s+\sigma}^{\gamma ,\mathcal{O}}+\| \mathfrak{I}_{0}\|_{q,s+\sigma}^{\gamma ,\mathcal{O}}\|\rho\|_{q,s_{0}+\sigma}^{\gamma,\mathcal{O}}\right)
					\end{equation}
					and such that in the Cantor set
					$$\mathcal{G}_n(\gamma,\tau_{1},\tau_{2},i_{0}):=\mathcal{O}_{\infty,n}^{\gamma,\tau_{1}}(i_{0})\cap\mathcal{O}_{\infty,n}^{\gamma,\tau_{1},\tau_{2}}(i_{0})\cap\Lambda_{\infty,n}^{\gamma,\tau_{1}}(i_{0}),$$
					we have
					$$
					\widehat{\mathcal{L}}_{\omega}\mathtt{T}_{\omega,n}=\textnormal{Id}+\mathtt{E}_n,
					$$
					where $\mathtt{E}_n$ satisfies the following estimate
					\begin{align*}
%						\forall\, s\in [s_0,S],\quad &\sup_{n\in\mathbb{N}}\|\mathtt{E}_n\rho\|_{q,s}^{\gamma ,\mathcal{O}} \lesssim \gamma^{-1}\|\rho\|_{{s}+\sigma}^{\gamma ,\mathcal{O}}+\varepsilon\gamma^{-3}\| \mathfrak{I}_{0}\|_{q,s+\sigma}^{\gamma,\mathcal{O}}\|\rho\|_{q,s_{0}+\sigma}^{\gamma,\mathcal{O}},\\
						\forall\, s\in [s_0,S],\quad  &\|\mathtt{E}_n\rho\|_{q,s_0}^{\gamma ,\mathcal{O}}
						\nonumber\lesssim N_n^{s_0-s}\gamma^{-1}\Big( \|\rho\|_{q,s+\sigma}^{\gamma,\mathcal{O}}+\varepsilon\gamma^{-2}\| \mathfrak{I}_{0}\|_{q,s+\sigma}^{\gamma,\mathcal{O}}\|\rho\|_{q,s_{0}+\sigma}^{\gamma,\mathcal{O}} \Big)\\
						&\qquad\qquad\qquad\quad+ \varepsilon\gamma^{-3}N_{0}^{{\mu}_{2}}N_{n+1}^{-\mu_{2}} \|\rho\|_{q,s_0+\sigma}^{\gamma,\mathcal{O}}.
					\end{align*}
					Recall that  $\widehat{\mathcal{L}}_{\omega},$ $  \mathcal{O}_{\infty,n}^{\gamma,\tau_{1}}(i_{0})$ and $\mathcal{O}_{\infty,n}^{\gamma,\tau_{1},\tau_{2}}(i_{0})$ are given in Propositions $\ref{lemma setting for Lomega}, \ref{reduction of the transport part}$ and $\ref{reduction of the remainder term},$ respectively.
					\item In the Cantor set $\mathcal{G}_{n}(\gamma,\tau_{1},\tau_{2},i_{0})$, we have the following splitting
					$$\widehat{\mathcal{L}}_{\omega}=\widehat{\mathtt{L}}_{\omega,n}+\widehat{\mathtt{R}}_{n}\quad\textnormal{with}\quad\widehat{\mathtt{L}}_{\omega,n}\mathtt{T}_{\omega,n}=\textnormal{Id}\quad\textnormal{and}\quad\widehat{\mathtt{R}}_{n}=\mathtt{E}_{n}\widehat{\mathtt{L}}_{\omega,n},$$
					where the operators $\widehat{\mathtt{L}}_{\omega,n}$ and $\widehat{\mathtt{R}}_{n}$ are defined in  $\mathcal{O}$ and satisfy the following estimates
					\begin{align*}
						\forall s\in[s_{0},S],\quad& \sup_{n\in\mathbb{N}}\|\widehat{\mathtt{L}}_{\omega,n}\rho\|_{q,s}^{\gamma,\mathcal{O}}\lesssim\|\rho\|_{q,s+1}^{\gamma,\mathcal{O}}+\varepsilon\gamma^{-2}\|\mathfrak{I}_{0}\|_{q,s+\sigma}^{\gamma,\mathcal{O}}\|\rho\|_{q,s_{0}+1}^{\gamma,\mathcal{O}},\\
						%\forall s\in[s_{0},S],\quad &\sup_{n\in\mathbb{N}}\|\widehat{\mathtt{R}}_{n}\rho\|_{q,s}^{\gamma,\mathcal{O}}\lesssim\gamma^{-1}\|\rho\|_{q,s+\sigma}^{\gamma,\mathcal{O}}+\varepsilon\gamma^{-3}\|\mathfrak{I}_{0}\|_{q,s+\sigma}^{\gamma,\mathcal{O}}\|\rho\|_{q,s_{0}}^{\gamma,\mathcal{O}},\\
						\forall s\in[s_{0},S],\quad &\|\widehat{\mathtt{R}}_{n}\rho\|_{q,s_{0}}^{\gamma,\mathcal{O}}\lesssim N_{n}^{s_{0}-s}\gamma^{-1}\left(\|\rho\|_{q,s+\sigma}^{\gamma,\mathcal{O}}+\varepsilon\gamma^{-2}\|\mathfrak{I}_{0}\|_{q,s+\sigma}^{\gamma,\mathcal{O}}\|\rho\|_{q,s_{0}+\sigma}^{\gamma,\mathcal{O}}\right)\\
						&\qquad\qquad\qquad\quad+\varepsilon\gamma^{-3}N_{0}^{\mu_{2}}N_{n+1}^{-\mu_{2}}\|\rho\|_{q,s_{0}+\sigma}^{\gamma,\mathcal{O}}.
					\end{align*}
				\end{enumerate}
			\end{prop}
			\begin{proof}
				{\bf{(i)}} From Proposition \ref{reduction of the remainder term} we recall that
				\begin{align*}
					\mathscr{L}_{\infty}&=\omega\cdot\partial_{\varphi}\Pi_{\mathbb{S}_0}^{\perp}+\mathscr{D}_{\infty}.
				\end{align*}
				Then we may split this operator as follows, using the projectors defined in \eqref{definition of projections for functions}
				\begin{align}\label{Tikl1}
					\nonumber\mathscr{L}_{\infty}&=\Pi_{N_n}\omega\cdot\partial_{\varphi}\Pi_{N_n}\Pi_{\mathbb{S}_0}^{\perp}+\mathscr{D}_{\infty}-\Pi_{N_n}^\perp\omega\cdot\partial_{\varphi}\Pi_{N_n}^\perp\Pi_{\mathbb{S}_0}^{\perp}\nonumber\\
					&:=\mathtt{L}_n-\mathtt{R}_n,
				\end{align}
				with $\mathtt{R}_{n}:=\Pi_{N_n}^\perp\omega\cdot\partial_{\varphi}\Pi_{N_n}^\perp\Pi_{\mathbb{S}_0}^{\perp}.$ From this definition and the structure of $\mathscr{D}_{\infty}$ in Proposition \ref{reduction of the remainder term} we deduce that
				$$\forall(l,j)\in\mathbb{Z}^{d}\times\mathbb{S}_{0}^{c},\quad {\bf e}_{-l,-j}\mathtt{L}_n{\bf e}_{l,j}=\left\lbrace\begin{array}{ll}
					\ii\big(\omega\cdot l+\mu_j^\infty\big) & \hbox{if }  |l|\leqslant N_n\\
					\ii\,\mu_j^\infty&\hbox{if } |l|> N_n.
				\end{array}\right.$$
				Define the diagonal  operator  $\mathtt{T}_n$ by 
				\begin{align*}
					\mathtt{T}_{n}\rho(\lambda,\omega,\varphi,\theta):=&
					-\ii\sum_{(l,j)\in\mathbb{Z}^{d}\times\mathbb{S}_{0}^{c}\atop |l|\leqslant N_n }\tfrac{\chi\left((\omega\cdot l+\mu_{j}^{\infty}(\lambda,\omega,i_{0}))\gamma ^{-1}\langle l\rangle^{\tau_{1}}\right)}{\omega\cdot l+\mu_{j}^{\infty}(\lambda,\omega,i_{0})}\rho_{l,j}(\lambda,\omega)\,e^{\ii(l\cdot\varphi+j\theta)}\\
					&-\ii\sum_{(l,j)\in\mathbb{Z}^{d}\times\mathbb{S}_{0}^{c}\atop|l|>N_n}\tfrac{\rho_{l,j}(\lambda,\omega)}{\mu_{j}^{\infty}(\lambda,\omega,i_{0})}\,e^{\ii(l\cdot\varphi+j\theta)},
				\end{align*}
				where $\chi$ is the cut-off function defined in \eqref{properties cut-off function first reduction} and $\left(\rho_{l,j}(\lambda,\omega)\right)_{l,j}$ are  the Fourier coefficients of $\rho$. We recall from Proposition \ref{reduction of the remainder term} that
				$$\mu_{j}^{\infty}(\lambda,\omega,i_{0})=\Omega_{j}(\lambda)+jr^{1}(\lambda,\omega)+r_{j}^{\infty}(\lambda,\omega)\quad \mbox{ with }\quad r^{1}(\lambda,\omega)=c_{i_{0}}(\lambda,\omega)-V_0(\lambda),$$
				with the estimates  
				$$
				\forall j\in\mathbb{S}_{0}^{c},\quad \|\mu^\infty_{j}\|_{q}^{\gamma ,\mathcal{O}}\lesssim |j|,
				$$
				where we use in part the estimate \eqref{reg-G-1}, \eqref{estimate rjinfty} and \eqref{differences mu0}.  According  to Lemma \ref{lemma properties linear frequencies}-(iii), \eqref{reg-G-1}, \eqref{estimate rjinfty} and the smallness condition \eqref{hypothesis KAM reduction of the remainder term} we infer
				$$|j|\lesssim \|\mu^\infty_{j}\|_{0}^{\gamma ,\mathcal{O}}\leqslant  \|\mu^\infty_{j}\|_{q}^{\gamma ,\mathcal{O}},$$
				Implementing the  same arguments as for  \eqref{link Psi and R} one gets  
				\begin{align}\label{Es-Mu1}
					\forall s\geqslant s_0,\quad  \|\mathtt{T}_{n}\rho\|_{q,s}^{\gamma ,\mathcal{O}}\lesssim \gamma ^{-1}\|\rho\|_{s+\tau_{1}q+\tau_{1}}^{\gamma ,\mathcal{O}}.
				\end{align}
				Moreover, by construction 
				\begin{equation}\label{Inv-Ty1}
					\mathtt{L}_{n}\mathtt{T}_{n}=\textnormal{Id}\quad\textnormal{in }\Lambda_{\infty,n}^{\gamma,\tau_{1}}(i_{0})
				\end{equation} 
				since $\chi(\cdot)=1$ in this set. It follows from \eqref{Tikl1} that 
				\begin{align}\label{LinftyTn}
					\forall\, (\lambda,\omega)\in \Lambda_{\infty,n}^{\gamma,\tau_{1}}(i_{0}),\quad \mathscr{L}_\infty \mathtt{T}_{n}&=\textnormal{Id}-\mathtt{R}_n\mathtt{T}_{n}\nonumber\\
					&:=\textnormal{Id}+\mathtt{E}_n^3.
				\end{align}
				Notice that by Lemma \ref{Lem-lawprod}-(ii),
				$$
				\forall \, s_0\leqslant s\leqslant \overline{s},\quad \|\mathtt{R}_{n}\rho\|_{q,s}^{\gamma ,\mathcal{O}}\lesssim N_n^{s-\overline{s}}\|\rho\|_{\overline{s}+1}^{\gamma ,\mathcal{O}}.
				$$
				Combining this estimate with \eqref{Es-Mu1} yields
				\begin{align}\label{Mila-1}
					\nonumber  \forall \, s_0\leqslant s\leqslant \overline{s},\quad \|\mathtt{E}^3_{n}\rho\|_{q,s}^{\gamma ,\mathcal{O}}&\lesssim N_n^{s-\overline{s}}\|\mathtt{T}_{n}\rho\|_{\overline{s}+1}^{\gamma ,\mathcal{O}}\\
					&\lesssim 
					N_n^{s-\overline{s}}\gamma^{-1}\|\rho\|_{\overline{s}+1+\tau_{1}q+\tau_{1}}^{\gamma ,\mathcal{O}}.
				\end{align}
				{{\bf{(ii)}} }Let us define  
				\begin{align}\label{tomega}
					\mathtt{T}_{\omega,n}:=\mathscr{B}_{\perp}\Phi_{\infty}\mathtt{T}_{n}\Phi_{\infty}^{-1}\mathscr{B}_{\perp}^{-1},
				\end{align}
				where the operators $ \mathscr{B}_{\perp}$ and $\Phi_{\infty}$ are defined in Propositions \ref{projection in the normal directions} and \ref{reduction of the remainder term} respectively. Notice that $\mathtt{T}_{\omega,n}$ is defined in the whole range of parameters $\mathcal{O}.$ Since the condition \eqref{small-1} is satisfied, then, both Propositions \ref{reduction of the transport part} and \ref{reduction of the remainder term} apply and from \eqref{estimate for mathscrBperp and its inverse}  we obtain
				$$\forall s\in[s_{0},S],\quad\|\mathtt{T}_{\omega,n}\rho\|_{q,s}^{\gamma,\mathcal{O}}\lesssim\|\Phi_{\infty}\mathtt{T}_{n}\Phi_{\infty}^{-1}\mathscr{B}_{\perp}^{-1}\rho\|_{q,s}^{\gamma,\mathcal{O}}+\|\mathfrak{I}_{0}\|_{q,s+\sigma}^{\gamma,\mathcal{O}}\|\Phi_{\infty}\mathtt{T}_{n}\Phi_{\infty}^{-1}\mathscr{B}_{\perp}^{-1}\rho\|_{q,s_{0}}^{\gamma,\mathcal{O}}.$$
				By using \eqref{estimate on Phiinfty and its inverse} and \eqref{small-1}, one gets
				$$\forall s\in[s_{0},S],\quad\|\Phi_{\infty}\mathtt{T}_{n}\Phi_{\infty}^{-1}\mathscr{B}_{\perp}^{-1}\rho\|_{q,s}^{\gamma,\mathcal{O}}\lesssim\|\mathtt{T}_{n}\Phi_{\infty}^{-1}\mathscr{B}_{\perp}^{-1}\rho\|_{q,s}^{\gamma,\mathcal{O}}+\|\mathfrak{I}_{0}\|_{q,s+\sigma}^{\gamma,\mathcal{O}}\|\mathtt{T}_{n}\Phi_{\infty}^{-1}\mathscr{B}_{\perp}^{-1}\rho\|_{q,s_{0}}^{\gamma,\mathcal{O}}.$$
				Thus  the point (i) of the current proposition implies
				$$\forall s\geqslant s_{0},\quad\|\mathtt{T}_{n}\Phi_{\infty}^{-1}\mathscr{B}_{\perp}^{-1}\rho\|_{q,s}^{\gamma,\mathcal{O}}\lesssim\gamma^{-1}\|\Phi_{\infty}^{-1}\mathscr{B}_{\perp}^{-1}\rho\|_{q,s+\tau_{1}q+\tau_{1}}^{\gamma,\mathcal{O}}.$$
				Applying  \eqref{estimate on Phiinfty and its inverse} and \eqref{estimate for mathscrBperp and its inverse} with \eqref{small-1} yields
				\begin{align*}
					\forall s\in[s_{0},S],\quad\|\Phi_{\infty}^{-1}\mathscr{B}_{\perp}^{-1}\rho\|_{q,s}^{\gamma,\mathcal{O}}&\lesssim\|\mathscr{B}_{\perp}^{-1}\rho\|_{q,s}^{\gamma,\mathcal{O}}+\|\mathfrak{I}_{0}\|_{q,s+\sigma}^{\gamma,\mathcal{O}}\|\mathscr{B}_{\perp}^{-1}\rho\|_{q,s_{0}}^{\gamma,\mathcal{O}}\\
					&\lesssim\|\rho\|_{q,s}^{\gamma,\mathcal{O}}+\|\mathfrak{I}_{0}\|_{q,s+\sigma}^{\gamma,\mathcal{O}}\|\rho\|_{q,s_{0}}^{\gamma,\mathcal{O}}.
				\end{align*}
				Putting  together the preceding three  estimates gives \eqref{estimate mathcalTomega}. Now combining Propositions \ref{projection in the normal directions} and \ref{reduction of the remainder term}, we find that in the Cantor set  $\mathcal{O}_{\infty,n}^{\gamma,\tau_{1}}(i_{0})\cap\mathcal{O}_{\infty,n}^{\gamma,\tau_{1},\tau_{2}}(i_{0})$
				the following decomposition holds
				\begin{align*}
					\Phi^{-1}_{\infty}\mathscr{B}_{\perp}^{-1}\widehat{\mathcal{L}}_{\omega}\mathscr{B}_{\perp}\Phi_{\infty}&=\Phi^{-1}_{\infty}\mathscr{L}_{0}\Phi_{\infty}+\Phi^{-1}_{\infty}\mathtt{E}_{n}^1\Phi_{\infty}\\
					&=\mathscr{L}_{\infty}+\mathtt{E}_n^2+\Phi^{-1}_{\infty}\mathtt{E}_{n}^1\Phi_{\infty}.
				\end{align*}
				It follows that in the Cantor set  $\mathcal{O}_{\infty,n}^{\gamma,\tau_{1}}(i_{0})\cap\mathscr{O}_{\infty,n}^{\gamma,\tau_{2}}(i_{0})\cap \Lambda_{\infty,n}^{\gamma,\tau_{1}}(i_{0})$ one has  by virtue of the identity \eqref{LinftyTn}
				\begin{align*}
					\Phi^{-1}_{\infty}\mathscr{B}_{\perp}^{-1}\widehat{\mathcal{L}}_{\omega}\mathscr{B}_{\perp}\Phi_{\infty}\mathtt{T}_{n}
					&=\textnormal{Id}+\mathtt{E}_n^3+\mathtt{E}_n^2\mathtt{T}_{n}+\Phi^{-1}_{\infty}\mathtt{E}_{n}^1\Phi_{\infty}\mathtt{T}_{n},
				\end{align*}
				which gives, using \eqref{tomega}, the following identity in $\mathcal{G}_n(\gamma,\tau_{1},\tau_{2},i_{0})$
				\begin{align}\label{Id-Moh1}
					\widehat{\mathcal{L}}_{\omega}\mathtt{T}_{\omega,n}&=\textnormal{Id}+\mathscr{B}_{\perp}\Phi_{\infty}\big(\mathtt{E}_n^3+\mathtt{E}_n^2\mathtt{T}_{n}+\Phi^{-1}_{\infty}\mathtt{E}_{n}^1\Phi_{\infty}\mathtt{T}_{n}\big)\Phi^{-1}_{\infty}\mathscr{B}_{\perp}^{-1}
					\nonumber\\
					&:=\textnormal{Id}+\mathscr{B}_{\perp}\Phi_{\infty}\mathtt{E}_n^4\Phi^{-1}_{\infty}\mathscr{B}_{\perp}^{-1}\nonumber\\
					&:=\textnormal{Id}+\mathtt{E}_n.
				\end{align}
				The  estimate of the first term of  $\mathtt{E}_n^4$ is given in \eqref{Mila-1}. For the second term of $\mathtt{E}_n^4$
%				we use the estimates \eqref{Error-Est-1D} and \eqref{Es-Mu1} leading to
%				\begin{align}\label{Error-Est-1LD}
%					\nonumber \forall s\geqslant s_{0},\quad \mbox{ }\|\mathtt{E}^2_n\mathtt{T}_{n}\rho\|_{q,s}^{\gamma ,\mathcal{O}}&\lesssim \|\mathtt{T}_{n}\rho\|_{q,s+1}^{\gamma,\mathcal{O}}+\varepsilon\gamma^{-2}\| \mathfrak{I}_{0}\|_{q,s+\sigma}^{\gamma,\mathcal{O}}\|\mathtt{T}_{n}\rho\|_{q,s_{0}+1}^{\gamma,\mathcal{O}}\\
%					&\lesssim \gamma^{-1}\|\rho\|_{q,s+1+\tau_1q+\tau_1}^{\gamma,\mathcal{O}}+\varepsilon\gamma^{-3}\| \mathfrak{I}_{0}\|_{q,s+\sigma}^{\gamma,\mathcal{O}}\|\rho\|_{q,s_{0}+1+\tau_1q+\tau_1}^{\gamma,\mathcal{O}}
%				\end{align}
				we use \eqref{Error-Est-2D} and \eqref{Es-Mu1} leading to
				\begin{align}\label{Error-Est-2LD}
					\nonumber \|\mathtt{E}^2_n\mathtt{T}_{n}\rho\|_{q,s_0}^{\gamma ,\mathcal{O}}&\lesssim\varepsilon\gamma^{-2}N_{0}^{{\mu}_{2}}N_{n+1}^{-\mu_{2}} \|\mathtt{T}_{n}\rho\|_{q,s_0+1}^{\gamma,\mathcal{O}}\\
					&\lesssim\varepsilon\gamma^{-3}N_{0}^{{\mu}_{2}}N_{n+1}^{-\mu_{2}} \|\rho\|_{q,s_0+1+\tau_1q+\tau_1}^{\gamma,\mathcal{O}}.
				\end{align}
%				To estimate $\Phi^{-1}_{\infty}\mathtt{E}_{n}^1\Phi_{\infty}\mathtt{T}_{n}$, we combine \eqref{estimate on Phiinfty and its inverse}, \eqref{En1-s}, \eqref{Es-Mu1} and \eqref{small-1}
%				\begin{align}\label{Error-Est-1LD-1}
%					\|\Phi^{-1}_{\infty}\mathtt{E}_{n}^1\Phi_{\infty}\mathtt{T}_{n}\rho\|_{q,s}^{\gamma,\mathcal{O}}&\lesssim\|\mathtt{E}_{n}^1\Phi_{\infty}\mathtt{T}_{n}\rho\|_{q,s}^{\gamma,\mathcal{O}}+\varepsilon\gamma^{-2}\|\mathfrak{I}_{0}\|_{q,s+\sigma}^{\gamma,\mathcal{O}}\|\mathtt{E}_{n}^1\Phi_{\infty}\mathtt{T}_{n}\rho\|_{q,s_0}^{\gamma,\mathcal{O}}\nonumber\\
%					&\lesssim\|\Phi_{\infty}\mathtt{T}_{n}\rho\|_{q,s+2}^{\gamma,\mathcal{O}}+\varepsilon\gamma^{-2}\|\mathfrak{I}_{0}\|_{q,s+\sigma}^{\gamma,\mathcal{O}}\|\Phi_{\infty}\mathtt{T}_{n}\rho\|_{q,s_{0}+2}^{\gamma,\mathcal{O}}\nonumber\\
%					&\lesssim\|\mathtt{T}_n\rho\|_{q,s+2}^{\gamma,\mathcal{O}}+\varepsilon\gamma^{-2}\|\mathfrak{I}_{0}\|_{q,s+\sigma}^{\gamma,\mathcal{O}}\|\mathtt{T}_n\rho\|_{q,s_{0}+2}^{\gamma,\mathcal{O}}\nonumber\\
%					&\lesssim\gamma^{-1}\|\rho\|_{q,s+2+\tau_{1}q+\tau_{1}}^{\gamma,\mathcal{O}}+\varepsilon\gamma^{-3}\|\mathfrak{I}_{0}\|_{q,s+\sigma}^{\gamma,\mathcal{O}}\|\rho\|_{q,s_{0}+2+\tau_{1}q+\tau_{1}}^{\gamma,\mathcal{O}}.
%				\end{align} 
				For the estimate of $\Phi^{-1}_{\infty}\mathtt{E}_{n}^1\Phi_{\infty}\mathtt{T}_{n}$, we combine \eqref{estimate on Phiinfty and its inverse}, \eqref{En1-s0}, \eqref{Es-Mu1} and \eqref{small-1} to get 
				\begin{align}\label{Error-Est-2LD-1}
					\|\Phi^{-1}_{\infty}\mathtt{E}_{n}^1\Phi_{\infty}\mathtt{T}_{n}\rho\|_{q,s_{0}}^{\gamma,\mathcal{O}}&\lesssim\|\mathtt{E}_{n}^1\Phi_{\infty}\mathtt{T}_{n}\rho\|_{q,s_{0}}^{\gamma,\mathcal{O}}\nonumber\\
					&\lesssim\varepsilon N_0^{\mu_{2}}N_{n+1}^{-\mu_2}\|\Phi_{\infty}\mathtt{T}_{n}\rho\|_{q,s_{0}+2}^{\gamma,\mathcal{O}}\nonumber\\
					%&\lesssim\varepsilon N_0^{\mu_{2}}N_{n+1}^{-\mu_2}\|\mathtt{T}_{n}\rho\|_{q,s_{0}+2}^{\gamma,\mathcal{O}}\nonumber\\
					&\lesssim\varepsilon\gamma^{-1}N_0^{\mu_{2}}N_{n+1}^{-\mu_2}\|\rho\|_{q,s_{0}+2+\tau_1q+\tau_1}^{\gamma,\mathcal{O}}.
				\end{align}
%				Putting together \eqref{Mila-1} and \eqref{Error-Est-1LD} and \eqref{Error-Est-1LD-1} we find
%				\begin{align}\label{MLLD}
%					\forall s\in[s_{0},S],\quad \mbox{ }\|\mathtt{E}^4_n\rho\|_{q,s}^{\gamma ,\mathcal{O}}&\lesssim \gamma^{-1}\|\rho\|_{s+2+\tau_{1}q+\tau_{1}}^{\gamma ,\mathcal{O}}+\varepsilon\gamma^{-3}\| \mathfrak{I}_{0}\|_{q,s+\sigma}^{\gamma,\mathcal{O}}\|\rho\|_{q,s_{0}+2+\tau_1q+\tau_1}^{\gamma,\mathcal{O}}.
%				\end{align}
				Putting together \eqref{Mila-1} and \eqref{Error-Est-2LD} and \eqref{Error-Est-2LD-1} we find
				\begin{align}\label{Fif-1}
					\|\mathtt{E}^4_n\rho\|_{q,s_0}^{\gamma ,\mathcal{O}}\lesssim N_n^{s_0-s}\gamma^{-1}\|\rho\|_{s+2+\tau_{1}q+\tau_{1}}^{\gamma ,\mathcal{O}}+ \varepsilon\gamma^{-3}N_{0}^{{\mu}_{2}}N_{n+1}^{-\mu_{2}} \|\rho\|_{q,s_0+2+\tau_1q+\tau_1}^{\gamma,\mathcal{O}}.
				\end{align}
				Set $\Psi=\mathscr{B}_{\perp}\Phi_{\infty}$ then from  \eqref{estimate on Phiinfty and its inverse}, \eqref{estimate for mathscrBperp and its inverse} and \eqref{small-1}  we deduce that
				\begin{align}\label{MLLD01}
					\forall s\in[s_{0},S],\quad \mbox{ }\|\Psi^{\pm1}\rho\|_{q,s}^{\gamma ,\mathcal{O}}\lesssim \|\rho\|_{q,s}^{\gamma,\mathcal{O}}+\varepsilon\gamma^{-2}\| \mathfrak{I}_{0}\|_{q,s+\sigma}^{\gamma,\mathcal{O}}\|\rho\|_{q,s_{0}}^{\gamma,\mathcal{O}}.
				\end{align}
				Straightforward computations based on  \eqref{Fif-1}, \eqref{MLLD01}  and  \eqref{small-1} yields 
				\begin{align*}
					\nonumber \|\Psi\mathtt{E}^4_n\Psi^{-1}\rho\|_{q,s_0}^{\gamma ,\mathcal{O}}&\lesssim \|\mathtt{E}^4_n\Psi^{-1}\rho\|_{q,s_0}^{\gamma,\mathcal{O}}\\
					\nonumber&\lesssim N_n^{s_0-s}\gamma^{-1}\|\Psi^{-1}\rho\|_{s+2+\tau_{1}q+\tau_{1}}^{\gamma ,\mathcal{O}}+ \varepsilon\gamma^{-3}N_{0}^{{\mu}_{2}}N_{n+1}^{-\mu_{2}} \|\Psi^{-1}\rho\|_{q,s_0+2+\tau_1q+\tau_1}^{\gamma,\mathcal{O}}\\
					\nonumber&\lesssim N_n^{s_0-s}\gamma^{-1}\big( \|\rho\|_{q,s+2+\tau_1q+\tau_1}^{\gamma,\mathcal{O}}+\varepsilon\gamma^{-2}\| \mathfrak{I}_{0}\|_{q,s+\sigma}^{\gamma,\mathcal{O}}\|\rho\|_{q,s_{0}}^{\gamma,\mathcal{O}} \big)\\
					&\quad+ \varepsilon\gamma^{-3}N_{0}^{{\mu}_{2}}N_{n+1}^{-\mu_{2}} \|\rho\|_{q,s_0+2+\tau_1q+\tau_1}^{\gamma,\mathcal{O}}.
				\end{align*}
				Consequently, taking $\sigma$ large enough, we get
				\begin{equation}\label{MLLD3}
					\sup_{n\in\mathbb{N}}\|\mathtt{E}_n\rho\|_{q,s_0}^{\gamma ,\mathcal{O}}
					\nonumber\lesssim N_n^{s_0-s}\gamma^{-1}\big( \|\rho\|_{q,s+\sigma}^{\gamma,\mathcal{O}}+\varepsilon\gamma^{-2}\| \mathfrak{I}_{0}\|_{q,s+\sigma}^{\gamma,\mathcal{O}}\|\rho\|_{q,s_{0}+\sigma}^{\gamma,\mathcal{O}} \big)+ \varepsilon\gamma^{-3}N_{0}^{{\mu}_{2}}N_{n+1}^{-\mu_{2}} \|\rho\|_{q,s_0+\sigma}^{\gamma,\mathcal{O}}.
				\end{equation}
				\textbf{(iii)} According to \eqref{Id-Moh1}, one can write in the Cantor set $\mathcal{G}_{n}(\gamma,\tau_{1},\tau_{2},i_{0})$
				\begin{equation}\label{Id-Moh11}
					\widehat{\mathcal{L}}_{\omega}=\mathtt{T}_{\omega,n}^{-1}+\mathtt{E}_{n}\mathtt{T}_{\omega,n}^{-1}.
				\end{equation}
				Gathering \eqref{tomega} and \eqref{Inv-Ty1}, one obtains in the Cantor set $\mathcal{G}_n(\gamma,\tau_1,\tau_2,i_0)$
				$$\widehat{\mathtt{L}}_{\omega,n}:=\mathtt{T}_{\omega,n}^{-1}=\mathscr{B}_{\perp}\Phi_{\infty}\mathtt{L}_{n}\Phi_{\infty}^{-1}\mathscr{B}_{\perp}^{-1}=\Psi\mathtt{L}_n\Psi^{-1}.$$
				Hence, \eqref{Id-Moh11} can be rewritten
				\begin{equation}\label{dec-Lom}
					\widehat{\mathcal{L}}_{\omega}=\widehat{\mathtt{L}}_{\omega,n}+\widehat{\mathtt{R}}_{n}\quad\textnormal{with}\quad\widehat{\mathtt{R}}_{n}:=\mathtt{E}_{n}\widehat{\mathtt{L}}_{\omega,n}.
				\end{equation}
				Putting together \eqref{Tikl1}, \eqref{MLLD01} and \eqref{small-1}, we obtain
				\begin{align}\label{est hatLomegan}
					\forall s\in[s_{0},S],\quad \|\widehat{\mathtt{L}}_{\omega,n}\rho\|_{q,s}^{\gamma,\mathcal{O}}&=\|\Psi\mathtt{L}_{n}\Psi^{-1}\rho\|_{q,s}^{\gamma,\mathcal{O}}\nonumber\\
					&\lesssim\|\mathtt{L}_n\Psi^{-1}\rho\|_{q,s}^{\gamma,\mathcal{O}}+\varepsilon\gamma^{-2}\|\mathfrak{I}_{0}\|_{q,s+\sigma}^{\gamma,\mathcal{O}}\|\mathtt{L}_n\Psi^{-1}\rho\|_{q,s_{0}}^{\gamma,\mathcal{O}}\nonumber\\
					&\lesssim\|\Psi^{-1}\rho\|_{q,s+1}^{\gamma,\mathcal{O}}+\varepsilon\gamma^{-2}\|\mathfrak{I}_{0}\|_{q,s+\sigma}^{\gamma,\mathcal{O}}\|\Psi^{-1}\rho\|_{q,s_{0}+1}^{\gamma,\mathcal{O}}\nonumber\\
					&\lesssim\|\rho\|_{q,s+1}^{\gamma,\mathcal{O}}+\varepsilon\gamma^{-2}\|\mathfrak{I}_{0}\|_{q,s+\sigma}^{\gamma,\mathcal{O}}\|\rho\|_{q,s_{0}+1}^{\gamma,\mathcal{O}}.
				\end{align}
				Hence combining this estimate with  \eqref{MLLD3} yields
				\begin{align*}
				\forall s\in[s_{0},S],\quad \sup_{n\in\mathbb{N}}\|\widehat{\mathtt{R}}_{n}\rho\|_{q,s_{0}}^{\gamma,\mathcal{O}}&\lesssim N_{n}^{s_{0}-s}\gamma^{-1}\left(\|\rho\|_{q,s+\sigma}^{\gamma,\mathcal{O}}+\varepsilon\gamma^{-2}\|\mathfrak{I}_{0}\|_{q,s+\sigma}^{\gamma,\mathcal{O}}\|\rho\|_{q,s_{0}+\sigma}^{\gamma,\mathcal{O}}\right)\\
				&+\varepsilon\gamma^{-3}N_{0}^{\mu_{2}}N_{n+1}^{-\mu_{2}}\|\rho\|_{q,s_{0}+\sigma}^{\gamma,\mathcal{O}}.
				\end{align*}
				This achieves the proof of the third point and the proof of the proposition is now complete.  
			\end{proof}
			\section{Proof of the main result}\label{Nash-Moser scheme and measure estimates}
			This section is devoted to the proof of Theorem \ref{main theorem}. For this aim we intend to  implement Nash-Moser scheme in order to construct 
			zeros for the nonlinear functional $\mathcal{F}\big(i, \alpha, \lambda,\omega, \varepsilon\big)$ defined in \eqref{main function}. We shall be able to capture the solutions when the  parameters $(\lambda, \omega)$ belong to a suitable   final Cantor \mbox{set $\mathcal{G}_{\infty}^{\gamma}$} obtained as the intersection of all the Cantor sets required during the steps of the scheme to invert the linearized operator. More precisely, we get a relatively smooth function
			$(\lambda,\omega)\in\mathcal{O}\mapsto U_{\infty}(\lambda,\omega)=\big(i_{\infty}(\lambda,\omega),\alpha_{\infty}(\lambda,\omega)\big)$ such that
			$$\forall(\lambda,\omega)\in\mathcal{G}_{\infty}^{\gamma},\quad\mathcal{F}\big(U_{\infty}(\lambda,\omega),\lambda,\omega,\varepsilon\big)=0.$$
			To generate  solutions to the initial Hamiltonian equation \eqref{perturbed hamiltonian} we should adjust the parameters so that $\alpha_{\infty}(\lambda,\omega)= -\omega_{\textnormal{Eq}}(\lambda)$, where $\omega_{\textnormal{Eq}}$ corresponds to the equilibrium frequency vector defined in \eqref{def linear frequency vector}. As a consequence, nontrivial solutions are constructed when the scalar parameter  $\lambda$ is  selected in  the final  Cantor set
			$$\mathcal{C}_{\infty}^{\varepsilon}=\Big\{\lambda\in(\lambda_{0},\lambda_{1})\quad\textnormal{s.t.}\quad(\lambda,\omega(\lambda,\varepsilon))\in\mathcal{G}_{\infty}^{\gamma}\quad\textnormal{with}\quad \alpha_{\infty}(\lambda,\omega(\lambda,\varepsilon))=-\omega_{\textnormal{Eq}}(\lambda)\Big\}.$$
			The measure of this  set will be discussed in Section \ref{Section measure of the final Cantor set}.
			\subsection{Nash-Moser scheme}\label{Sec Nash-Moser}
			In this section we implement the Nash-Moser scheme, which is  a modified Newton method implemented with a suitable Banach scales and through a frequency cut-off.  Basically, it  consists in a recursive construction of approximate solutions to the equation $\mathcal{F}\big(i,\alpha,\lambda,\omega,\varepsilon\big)=0$ where the functional  $\mathcal{F}$ is defined in \eqref{main function}. At each step of this scheme, we need to construct an approximate inverse of  the linearized operator at a state near the equilibrium by applying the reduction procedure developed in Section \ref{sec red lin op}. This enables  to get the result of Theorem \ref{proposition alomst approximate inverse}  with the suitable  tame estimates associated to  the final  loss of regularity $\overline{\sigma}$ that could be arranged to be large enough.  We point out that  $\overline{\sigma}$ depends only on the shape of the Cantor set through the  parameters $\tau_1,\tau_2, d$ and on the non degeneracy of the equilibrium frequency through $q=1+q_0$,  where $q_{0}$ be defined in Lemma \ref{lemma transversality}. However, $\overline{\sigma}$  is independent of  the regularity of the solutions that we want to construct.  Now, we shall fix the following parameters needed  to implement   Nash-Moser scheme and related to the geometry of the Cantor sets encoded in $\tau_1,\tau_2,d$ fixed by \eqref{setting tau1 and tau2} and to the parameter $q=q_0+1,$ 
			\begin{equation}\label{param NM}
				\left\lbrace\begin{array}{rcl}
					\overline{a} & = & \tau_{2}+2\\
					\mu_{1} & = & 3q(\tau_{2}+2)+6\overline{\sigma}+6\\
					a_{1} & = & 6q(\tau_{2}+2)+12\overline{\sigma}+15\\
					a_{2} & = & 3q(\tau_{2}+2)+6\overline{\sigma}+9\\
					\mu_{2} & = & 2q(\tau_{2}+2)+5\overline{\sigma}+7\\
					s_{h} & = & s_{0}+4q(\tau_{2}+2)+9\overline{\sigma}+11\\
					b_{1} & = & 2s_{h}-s_{0}.
				\end{array}\right.
			\end{equation}
			The numbers $a_1$ and $a_2$ will describe the rate of convergence for the regularity $s_0$ and $s_0+\overline{\sigma}$, respectively. They appear in the statements $(\mathcal{P}1)_n$ and $(\mathcal{P}2)_n$ in the Proposition \ref{Nash-Moser}. The parameter $\mu_1$ controls the norm inflation at the high regularity index $b_1$ and appears in $(\mathcal{P}3)_n.$ As to the parameter $\overline{a}$, it  is linked to the thickness  of a suitable  enlargement of the intermediate Cantor sets, needed   to construct classical extensions of our approximate solutions. Finally, the numbers $\mu_2$ and $s_h$ corresponds to those already encountered before in the reduction of the linearized operator and are now fixed to their minimal required values. In particular, we recall that $\mu_2$ corresponds to the rate of convergence of the error terms emerging in the almost reducibility of the linearized operator, for instance we refer  to Theorem \ref{proposition alomst approximate inverse}.
			We should emphasize that,  by taking  $\overline{\sigma}$ large enough, the choice for $\mu_2$ and $s_h$ done in \eqref{param NM}  enables to cover all the required assumptions in \eqref{param-trans} and \eqref{Conv-T2}. 
			Another assumption that we need to fix is related to $\gamma, N_{0}$ and  $\varepsilon$ 
			\begin{equation}\label{choice of gamma and N0 in the Nash-Moser}
				0<a<\tfrac{1}{\mu_{2}+q+2}, \quad \gamma:=\varepsilon^{a}, \quad N_{0}:=\gamma^{-1}.
			\end{equation}
			This constraint is required for the measuring the final Cantor set and to check that it is massive, for more details we refer to Proposition \ref{lem-meas-es1}. \\
			We shall start with  defining  the finite dimensional subspaces where  the approximate solutions are expected to live with controlled estimates. Consider the space,
			\begin{align*}
				E_{n}:=&\Big\{\mathfrak{I}=(\Theta,I,z)\quad\textnormal{s.t.}\quad\Theta=\Pi_{n}\Theta,\quad I=\Pi_{n}I\quad\textnormal{and}\quad z=\Pi_{n}z\Big\},
			\end{align*}
			where $\Pi_{n}$ is the projector defined by
			$$f(\varphi,\theta)=\sum_{(l,j)\in\mathbb{Z}^{d}\times\mathbb{Z}}f_{l,j}e^{\ii(l\cdot\varphi+j\theta)}\quad\Rightarrow\quad\Pi_{n}f(\varphi,\theta)=\sum_{\langle l,j\rangle\leqslant N_{n}}f_{l,j}e^{\ii(l\cdot\varphi+j\theta)},
			$$
			where the sequence $(N_n)$ is defined in \eqref{definition of Nm}. We observe that the same definition applies without ambiguity  when the functions depend only on $\varphi$ such as the action and the angles unknowns. The main result of this section is to prove the following induction statement. 		
			\begin{prop}[Nash-Moser]\label{Nash-Moser}
				{Let $(\tau_{1},\tau_{2},q,d,s_{0})$ satisfy \eqref{initial parameter condition} and \eqref{setting tau1 and tau2}. Consider the parameters fixed by \eqref{param NM} and \eqref{choice of gamma and N0 in the Nash-Moser}. There exist $C_{\ast}>0$ and $\varepsilon_{0}>0$ such that for any $\varepsilon\in[0,\varepsilon_{0}]$ we get for all $n\in\mathbb{N}$ the following properties,
					\begin{itemize}
						\item [$(\mathcal{P}1)_{n}$] There exists a $q$-times differentiable function 
						$$W_{n}:\begin{array}[t]{rcl}
							\mathcal{O} & \rightarrow &  E_{n-1}\times\mathbb{R}^{d}\times\mathbb{R}^{d+1}\\
							(\lambda,\omega) & \mapsto & \big(\mathfrak{I}_{n},\alpha_{n}-\omega,0\big)
						\end{array}$$
						satisfying 
						$$W_{0}=0\quad\mbox{ and }\quad\hbox{for}\quad n\geqslant 1,\quad\,\| W_{n}\|_{q,s_{0}+\overline{\sigma}}^{\gamma,\mathcal{O}}\leqslant C_{\ast}\varepsilon\gamma^{-1}N_{0}^{q\overline{a}}.$$
						By setting 
						$$U_{0}=\Big((\varphi,0,0),\omega,(\lambda,\omega)\Big)\quad\textnormal{and}\quad  \hbox{for}\quad n\in\mathbb{N}^*,\,U_{n}=U_{0}+W_{n}\quad \textnormal{and}\quad H_{n} =U_{n}-U_{n-1},$$
						then 
						$$\forall s\in[s_{0},S],\,\| H_{1}\|_{q,s}^{\gamma,\mathcal{O}}\leqslant \tfrac{1}{2}C_{\ast}\varepsilon\gamma^{-1}N_{0}^{q\overline{a}}\quad\mbox{ and }\quad\forall\, 2\leqslant k\leqslant n,\,\| H_{k}\|_{q,s_{0}+\overline{\sigma}}^{\gamma,\mathcal{O}}\leqslant C_{\ast}\varepsilon\gamma^{-1}N_{k-1}^{-a_{2}}.$$
						\item [$(\mathcal{P}2)_{n}$] Define 
						$$i_{n}=(\varphi,0,0)+\mathfrak{I}_{n},\quad \gamma_{n}=\gamma(1+2^{-n})\in[\gamma,2\gamma].$$
						The embedded torus $i_{n}$ statisfies the reversibility condition
						$$\mathfrak{S}i_n(\varphi)=i_n(-\varphi),$$
						where the involution $\mathfrak{S}$ is defined in \eqref{rev_aa}.
						Introduce
						$$\mathcal{A}_{0}^{\gamma}=\mathcal{O}\quad \mbox{ and }\quad \mathcal{A}_{n+1}^{\gamma}=\mathcal{A}_{n}^{\gamma}\cap\mathcal{G}_{n}(\gamma_{n+1},\tau_{1},\tau_{2},i_{n}),$$
						where $\mathcal{G}_{n}(\gamma_{n+1},\tau_{1},\tau_{2},i_{n})$ is described in Proposition $\ref{inversion of the linearized operator in the normal directions}$ and consider the open sets 
						$$
						\forall \mathtt{r}>0,\quad \mathrm{O}_{n}^\mathtt{r}:=\Big\{(\lambda,\omega)\in\mathcal{O}\quad\textnormal{s.t.}\quad {\mathtt{dist}}\big((\lambda,\omega),\mathcal{A}_{n}^{\gamma}\big)< \mathtt{r} N_{n}^{-\overline{a}}\Big\},
						$$
						where $\displaystyle {\mathtt{dist}}(x,A)=\inf_{y\in A}\|x-y\|$.  Then we have the following estimate 
						$$ \|\mathcal{F}(U_{n})\|_{q,s_{0}}^{\gamma,\mathrm{O}_{n}^{\gamma}}\leqslant C_{\ast}\varepsilon N_{n-1}^{-a_{1}}.
						$$
						\item[$(\mathcal{P}3)_{n}$] $\| W_{n}\|_{q,b_{1}+\overline{\sigma}}^{\gamma,\mathcal{O}}\leqslant C_{\ast}\varepsilon\gamma^{-1}N_{n-1}^{\mu_{1}}.$
				\end{itemize}}
			\end{prop}
			\begin{remark}
				Let $\mathrm{O}$ be an open subset of $\mathcal{O}.$ Since $\forall n\in\mathbb{N},\gamma_{n} \in[\gamma,2\gamma],$ then the norms $\|\cdot\|_{q,s}^{\gamma,\mathrm{O}}$ and  $\|\cdot\|_{q,s}^{\gamma_{n},\mathrm{O}}.$ are equivalent uniformly in $n.$
			\end{remark}
			\begin{proof}
				$\bullet$ \textbf{Initialization :} By construction, $U_{0}=\Big((\varphi,0,0),\omega,(\lambda,\omega)\Big)$ and  the flat torus $i_{\textnormal{\tiny{flat}}}(\varphi)=(\varphi,0,0)$ satisfies obviously the reversibility condition.
				%$$\mathfrak{S}i_{\textnormal{\tiny{flat}}}(\varphi)=i_{\textnormal{\tiny{flat}}}(-\varphi).$$
				By \eqref{main function}, we have
				$$\mathcal{F}(U_{0})=\varepsilon\left(\begin{array}{c}
					-\partial_{I}\mathcal{P}_{\varepsilon}((\varphi,0,0))\\
					\partial_{\vartheta}\mathcal{P}_{\varepsilon}((\varphi,0,0))\\
					-\partial_{\theta}\nabla_{z}\mathcal{P}_{\varepsilon}((\varphi,0,0))
				\end{array}\right).$$
				Using Lemma \ref{tame estimates for the vector field XmathcalPvarepsilon}, we get
				\begin{equation}\label{estimate mathcalF(U0)}
					\forall s\geqslant 0,\quad \|\mathcal{F}(U_{0})\|_{q,s}^{\gamma,\mathcal{O}}\leqslant C_{\ast}\varepsilon,
				\end{equation}
				up to take $C_{\ast}$ large enough. The properties $(\mathcal{P}1)_{0},$ $(\mathcal{P}2)_{0}$ and $(\mathcal{P}3)_{0}$ then follow immediately since $N_{-1}=1$ and $\mathrm{O}_{0}^{\gamma}=\mathcal{O}$ and by setting $W_0=0.$\\
				
				\noindent $\bullet$ \textbf{Induction step :} Given $n\in\mathbb{N},$ assume that $(\mathcal{P}1)_{k},$ $(\mathcal{P}2)_{k}$ and $(\mathcal{P}3)_{k}$ are true  for all $k\in\llbracket 0,n\rrbracket$ and let us check them at the next order $n+1$. Introduce the linearized operator of $\mathcal{F}$ at the state $(i_n,\alpha_n)$
				$$L_{n}:=L_{n}(\lambda,\omega):=d_{i,\alpha}\mathcal{F}(i_{n}\big(\lambda,\omega),\alpha_{n}(\lambda,\omega),(\lambda,\omega)\big).$$
				In order to construct the next approximation $U_{n+1},$ we need an approximate right inverse for $L_n.$ Its construction was performed along   the preceding sections and we refer to  Theorem \ref{proposition alomst approximate inverse} for a precise statement. To apply this result and get some bounds on $U_{n+1}$ we need to establish first some intermediate results connected to the smallness condition and to some Cantor set inclusions.\\
				
				\noindent $\blacktriangleright$ \textbf{Smallness/boundedness properties.} First of all, remark that the parameters conditions \eqref{param} are automatically satisfied by \eqref{param NM}. Then, provided that the smallness assumption \eqref{small-1} is satisfied, Proposition \ref{inversion of the linearized operator in the normal directions} applies.  It remains to check that \eqref{small-1} is satisfied. According to the first condition in \eqref{choice of gamma and N0 in the Nash-Moser} 			and  choosing $\varepsilon$ small enough, we can ensure
				\begin{equation}\label{small-NM}
					\varepsilon\gamma^{-2-q}N_0^{\mu_{2}}=\varepsilon^{1-a(\mu_2+q+2)}\leqslant\varepsilon_0
				\end{equation}
				for some a priori fixed $\varepsilon_0>0.$ Therefore the first assumption in \eqref{small-1}  holds. We now turn to the second assumption. Since from \eqref{param NM} $b_1=2s_h-s_0,$ then by interpolation inequality in Lemma \ref{Lem-lawprod}, we have
				\begin{equation}\label{interp NM}
					\|H_{n}\|_{q,s_{h}+\overline{\sigma}}^{\gamma,\mathcal{O}}\lesssim\left(\|H_{n}\|_{q,s_{0}+\overline{\sigma}}^{\gamma,\mathcal{O}}\right)^{\frac{1}{2}}\left(\|H_{n}\|_{q,b_{1}+\overline{\sigma}}^{\gamma,\mathcal{O}}\right)^{\frac{1}{2}}.		\end{equation}
				Besides, by using $(\mathcal{P}1)_{n}$, we find 
				\begin{equation}\label{e-H1-Hm}
					\forall s\in[s_{0},S],\quad \|H_{1}\|_{q,s}^{\gamma,\mathcal{O}}\leqslant \tfrac{1}{2}C_{\ast}\varepsilon\gamma^{-1}N_{0}^{q\overline{a}}\quad \textnormal{and}\quad \|H_{n}\|_{q,s_{0}+\overline{\sigma}}^{\gamma,\mathcal{O}}\leqslant C_{\ast}\varepsilon\gamma^{-1}N_{n-1}^{-a_{2}}.
				\end{equation}
				Now $(\mathcal{P}3)_{n}$ and $(\mathcal{P}3)_{n-1}$ imply
				\begin{align*}
					\|H_{n}\|_{q,b_{1}+\overline{\sigma}}^{\gamma,\mathcal{O}}&=\|U_{n}-U_{n-1}\|_{q,b_{1}+\overline{\sigma}}^{\gamma,\mathcal{O}}\\
					&=\|W_{n}-W_{n-1}\|_{q,b_{1}+\overline{\sigma}}^{\gamma,\mathcal{O}}\\
					&\leqslant\|W_{n}\|_{q,b_{1}+\overline{\sigma}}^{\gamma,\mathcal{O}}+\|W_{n-1}\|_{q,b_{1}+\overline{\sigma}}^{\gamma,\mathcal{O}}\\
					&\leqslant 2C_{\ast}\varepsilon\gamma^{-1}N_{n-1}^{\mu_{1}}.
				\end{align*}
				Putting together the foregoing estimates into \eqref{interp NM} gives for  $n\geqslant 2$,
				\begin{equation}\label{est Hmtilde sh+sigma}
					\|H_{n}\|_{q,s_{h}+\overline{\sigma}}^{\gamma,\mathcal{O}}\leqslant CC_{\ast}\varepsilon\gamma^{-1}N_{n-1}^{\frac{1}{2}(\mu_{1}-a_{2})}
				\end{equation}
				and for $n=1$,
				\begin{equation}\label{est H1tilde sh+sigma}
					\|H_{1}\|_{q,s_{h}+\overline{\sigma}}^{\gamma,\mathcal{O}}\leqslant \tfrac{1}{2}C_{\ast}\varepsilon\gamma^{-1}N_{0}^{q\overline{a}}.
				\end{equation}
				Now from \eqref{param NM} we infer
				\begin{equation}\label{cond_b1_1}
					a_{2}\geqslant \mu_{1}+2.
				\end{equation}
				Thus, by \eqref{choice of gamma and N0 in the Nash-Moser} and Lemma \ref{lemma sum Nn}, we get for small $\varepsilon$
				\begin{align*}
					\|W_{n}\|_{q,s_{h}+\overline{\sigma}}^{\gamma,\mathcal{O}}&\leqslant\|H_{1}\|_{q,s_{h}+\overline{\sigma}}^{\gamma,\mathcal{O}}+\sum_{k=2}^{n}\|H_{k}\|_{q,s_{h}+\overline{\sigma}}^{\gamma,\mathcal{O}}\\
					&\leqslant \tfrac{1}{2}C_{\ast}\varepsilon\gamma^{-1}N_{0}^{q\overline{a}}+ CC_{\ast}\varepsilon\gamma^{-1}\sum_{k=0}^{n}N_{k}^{-1}\\
					&\leqslant \tfrac{1}{2}C_{\ast}\varepsilon\gamma^{-1}N_{0}^{q\overline{a}}+ CN_{0}^{-1}C_{\ast}\varepsilon\gamma^{-1}\\
					&\leqslant C_{\ast}\varepsilon^{1-a(1+q\overline{a})}.
				\end{align*}
				One can check from  \eqref{param NM} and \eqref{choice of gamma and N0 in the Nash-Moser} that
				\begin{equation}\label{cond-lu1}
					a\leqslant\tfrac{1}{2(1+q\overline{a})}
				\end{equation}
				and therefore, by choosing $\varepsilon$ small enough and since $\overline{\sigma}\geqslant\sigma$, we get
				\begin{align*}
					\|\mathfrak{I}_{n}\|_{q,s_{h}+\sigma}^{\gamma,\mathcal{O}}\leqslant\|W_{n}\|_{q,s_{h}+\overline{\sigma}}^{\gamma,\mathcal{O}}&\leqslant C_{\ast}\varepsilon^{\frac{1}{2}}\\
					&\leqslant 1.
				\end{align*}
				As we have already mentioned, the parameter  $\overline{\sigma}$ is the final loss of regularity constructed in  \mbox{Theorem \ref{proposition alomst approximate inverse}} and depending only on the parameters $\tau_1,\tau_2,q$ and $d$ but it is independent of the state and the regularity.  Hence it can be selected large enough  such that $s_0+\overline{\sigma}\geqslant\overline{s}_h+\sigma_{4}$ where $\overline{s}_h$ and $\sigma_{4}$ are respectively defined in \eqref{param-trans} and Proposition \ref{reduction of the remainder term}. Then using \eqref{e-H1-Hm} and Sobolev embeddings, we obtain
				\begin{equation}\label{e-Hn-diff}\forall n\geqslant 2,\quad \|H_{n}\|_{q,\overline{s}_h+\sigma_{4}}^{\gamma,\mathcal{O}}\leqslant C_{\ast}\varepsilon\gamma^{-1}N_{n-1}^{-a_{2}}.
				\end{equation}\\
				$\blacktriangleright$ \textbf{Set inclusions.} From the previous point, Propositions \ref{reduction of the transport part},  \ref{reduction of the remainder term} and \ref{inversion of the linearized operator in the normal directions} apply and allow us to perform the reduction of the linearized operator in the normal directions at the current step. Therefore, the sets $\mathcal{A}_{k}^{\gamma}$ for all $k\leqslant n+1$ are well-defined. We shall now prove the following inclusions needed later to establish  suitable estimates for the extensions.
				\begin{equation}\label{hypothesis of induction set inclusions Nash-Moser}
					\mathcal{A}_{n+1}^{\gamma}\subset\mathrm{O}_{n+1}^{2\gamma}\subset\left(\mathcal{A}_{n+1}^{\frac{\gamma}{2}}\cap\mathrm{O}_{n}^{\gamma}\right).
				\end{equation}
				Notice that the first inclusion is obvious by construction since $\mathrm{O}_{n+1}^{2\gamma}$ is an enlargement of $\mathcal{A}_{n+1}^{\gamma}.$ It remains to prove the last inclusion.
				We have the inclusion
				\begin{equation}\label{inco}
					\forall k\in \llbracket 0,n\rrbracket,\quad  \mathrm{O}_{k+1}^{2\gamma}\subset \mathrm{O}_{k}^{\gamma}.
				\end{equation}
				Indeed, since by construction $\mathcal{A}_{k+1}^{\gamma}\subset\mathcal{A}_{k}^{\gamma}$ then taking $(\lambda,\omega)\in \mathrm{O}_{k+1}^{2\gamma}$ we have the following estimates
				\begin{align*}
					{\mathtt{dist}}\big((\lambda,\omega),\mathcal{A}_{k}^{\gamma}\big)&\leqslant  {\mathtt{dist}}\big((\lambda,\omega),\mathcal{A}_{k+1}^{\gamma}\big)\\
					&< 2\gamma N_{k+1}^{-\overline{a}}= 2\gamma N_{k}^{-\overline{a}}N_{0}^{-\frac{1}{2}\overline{a}}\\
					&< \gamma N_{k}^{-\overline{a}},
				\end{align*}
				provided that $2N_0^{-\frac12\overline a}<1,$ which is true up to take $N_0$ large enough, that is in view of \eqref{choice of gamma and N0 in the Nash-Moser} for $\varepsilon$ small enough. We shall now prove by induction in $k$ that
				\begin{equation}\label{hyprec O in A}
					\forall k\in \llbracket 0,n+1\rrbracket,\quad  \mathrm{O}_{k}^{2\gamma}\subset\mathcal{A}_{k}^{\frac{\gamma}{2}}.
				\end{equation}
				The case $k=0$ is trivial since $\mathrm{O}_{0}^{2\gamma}=\mathcal{O}=\mathcal{A}_{0}^{\frac{\gamma}{2}}.$ Let us now assume that \eqref{hyprec O in A} is true for the index $k\in \llbracket 0,n\rrbracket$ and let us check  it  at the next order. From \eqref{inco} and \eqref{hyprec O in A}, we obtain
				$$\mathrm{O}_{k+1}^{2\gamma}\subset\mathrm{O}_{k}^{\gamma}\subset\mathrm{O}_{k}^{2\gamma}\subset\mathcal{A}_{k}^{\frac{\gamma}{2}}.$$
				Therefore, we are left to check that 
				$$\mathrm{O}_{k+1}^{2\gamma}\subset\mathcal{G}_{k}\Big(\frac{\gamma_{k+1}}{2},\tau_{1},\tau_{2},i_{k}\Big).$$
				Let $(\lambda,\omega)\in\mathcal{O}_{k+1}^{2\gamma},$ then by construction, there exists $(\lambda',\omega')\in\mathcal{A}_{k+1}^{\gamma}$ such that 
				$$\mathtt{dist}\left((\lambda,\omega),(\lambda',\omega')\right)<2\gamma N_{k+1}^{-\overline{a}}.$$
				Hence, for all $(l,j)\in\mathbb{Z}^{d}\times\mathbb{S}_{0}^{c}$ with $|l|\leqslant N_{k},$ we have by left triangle and Cauchy-Schwarz inequalities together with $(\lambda',\omega')\in\Lambda_{\infty,k}^{\gamma_{k+1},\tau_{1}}(i_{k})$
				\begin{align*}
					\left|\omega\cdot l+\mu_{j}^{\infty}(\lambda,\omega,i_{k})\right|&\geqslant\left|\omega'\cdot l+\mu_{j}^{\infty}(\lambda',\omega',i_{k})\right|-|\omega-\omega'||l|-\left|\mu_{j}^{\infty}(\lambda,\omega,i_{k})-\mu_{j}^{\infty}(\lambda',\omega',i_{k})\right|\\
					&>\tfrac{\gamma_{k+1}\langle j\rangle}{\langle l\rangle^{\tau_{1}}}-2\gamma N_{k}N_{k+1}^{-\overline{a}}-\left|\mu_{j}^{\infty}(\lambda,\omega,i_{k})-\mu_{j}^{\infty}(\lambda',\omega',i_{k})\right|\\
					&>\tfrac{\gamma_{k+1}\langle j\rangle}{\langle l\rangle^{\tau_{1}}}-2\gamma N_{k+1}^{1-\overline{a}}-\left|\mu_{j}^{\infty}(\lambda,\omega,i_{k})-\mu_{j}^{\infty}(\lambda',\omega',i_{k})\right|.
				\end{align*}
				Using the Mean Value Theorem and the definition of $\mathrm{O}_{k+1}^{2\gamma}$ yields
				\begin{align*}
					\left|\mu_{j}^{\infty}(\lambda,\omega,i_{k})-\mu_{j}^{\infty}(\lambda',\omega',i_{k})\right|&\leqslant|(\lambda,\omega)-(\lambda',\omega')|\gamma^{-1}\|\mu_{j}^{\infty}(i_{k})\|_{q}^{\gamma,\mathcal{O}}\\
					&\leqslant 2N_{k+1}^{-\overline{a}}\|\mu_{j}^{\infty}(i_{k})\|_{q}^{\gamma,\mathcal{O}}.
				\end{align*}
				On the other hand,
				$$\forall j\in\mathbb{S}_{0}^{c},\quad \|\mu_{j}^{\infty}(i_{k})\|_{q}^{\gamma,\mathcal{O}}\leqslant\|\mu_{j}^{\infty}(i_{k})-\Omega_{j}\|_{q}^{\gamma,\mathcal{O}}+\|\Omega_{j}\|_{q}^{\gamma,\mathcal{O}}.$$
				Using the asymptotic \eqref{asymptotic expansion of the eigenvalues at the equilibrium} and the smoothness of $\lambda\mapsto I_{j}(\lambda)K_{j}(\lambda)$ for all $j\in\mathbb{N}^{*}$, one has
				$$\|\Omega_{j}\|_{q}^{\gamma,\mathcal{O}}\leqslant C|j|.$$
				Since \eqref{hypothesis KAM reduction of the remainder term} is satisfied by the previous point, we can apply \eqref{estim mujinfty} and obtain
				$$\forall j\in\mathbb{S}_{0}^{c},\quad \|\mu_{j}^{\infty}(i_{k})-\Omega_{j}\|_{q}^{\gamma,\mathcal{O}}\leqslant C|j|.$$
				Hence
				$$\forall j\in\mathbb{S}_{0}^{c},\quad \|\mu_{j}^{\infty}(i_{k})\|_{q}^{\gamma,\mathcal{O}}\leqslant C|j|.$$
				It follows that
				$$\left|\mu_{j}^{\infty}(\lambda,\omega,i_{k})-\mu_{j}^{\infty}(\lambda',\omega',i_{k})\right|\leqslant C\langle j\rangle N_{k+1}^{-\overline{a}}\leqslant C\gamma\langle j\rangle N_{k+1}^{1-\overline{a}}.$$
				Since $|l|\leqslant N_{k}\leqslant N_{k+1}$ and $\gamma_{k+1}\geqslant \gamma$, we obtain
				\begin{align*}
					\left|\omega\cdot l+\mu_{j}^{\infty}(\lambda,\omega,i_{k})\right|&\geqslant\tfrac{\gamma_{k+1}\langle j\rangle}{\langle l\rangle^{\tau_{1}}}-C\gamma\langle j\rangle N_{k+1}^{1-\overline{a}}\\
					&\geqslant\tfrac{\gamma_{k+1}\langle j\rangle}{\langle l\rangle^{\tau_{1}}}\left(1-CN_{k+1}^{\tau_{1}+1-\overline{a}}\right).
				\end{align*}
				From \eqref{param NM} and \eqref{setting tau1 and tau2} we infer
				\begin{equation}\label{cond-abarre-1}
					\overline{a}\geqslant  \tau_{2}+2\geqslant \tau_{1}+2
				\end{equation} and we can take $N_{0}$  sufficiently large to ensure
				$$CN_{k+1}^{\tau_{1}+1-\overline{a}}\leqslant CN_{0}^{-1}<\tfrac{1}{2},$$
				allowing to finally get
				$$\left|\omega\cdot l+\mu_{j}^{\infty}(\lambda,\omega,i_{k})\right|>\tfrac{\gamma_{k+1}\langle j\rangle}{2\langle l\rangle^{\tau_{1}}}\cdot
				$$
				This shows that, $(\lambda,\omega)\in\Lambda_{\infty,k}^{\frac{\gamma_{k+1}}{2},\tau_{1}}(i_{k}).$ Let us now check that $(\lambda,\omega)\in\mathcal{O}_{\infty,k}^{\frac{\gamma_{k+1}}{2},\tau_{1}}(i_{k}).$ For all $(l,j)\in\mathbb{Z}^{d}\times\mathbb{S}_{0}^{c}$ with $|l|\leqslant N_{k},$ we have by Cauchy-Schwarz inequality together with $(\lambda',\omega')\in\mathcal{O}_{\infty,k}^{\gamma_{k+1},\tau_{1}}(i_{k})$
				\begin{align*}
					\left|\omega\cdot l+jc_{i_{k}}(\lambda,\omega)\right|&\geqslant\left|\omega'\cdot l+jc_{i_{k}}(\lambda',\omega')\right|-|\omega-\omega'||l|-|j|\left|c_{i_{k}}(\lambda,\omega)-c_{i_{k}}(\lambda',\omega')\right|\\
					&>\tfrac{4\gamma_{k+1}^{\upsilon}\langle j\rangle}{\langle l\rangle^{\tau_{1}}}-2\gamma N_{k+1}^{1-\overline{a}}-\langle j\rangle\left|c_{i_{k}}(\lambda,\omega)-c_{i_{k}}(\lambda',\omega')\right|.
				\end{align*}
				Using the Mean Value Theorem and the definition of $\mathrm{O}_{k+1}^{2\gamma}$ yields
				$$\left|c_{i_{k}}(\lambda,\omega)-c_{i_{k}}(\lambda',\omega')\right|\leqslant CN_{k+1}^{-\overline{a}}\|c_{i_{k}}\|_{q}^{\gamma,\mathcal{O}}.$$
				Since \eqref{smallness condition transport} is satisfied by the previous point, we can apply \eqref{estimate r1} leading to 
				\begin{align*}
					\|c_{i_{k}}\|_{q}^{\gamma,\mathcal{O}}&\leqslant\|c_{i_{k}}-V_{0}\|_{q}^{\gamma,\mathcal{O}}+\|V_{0}\|_{q}^{\gamma,\mathcal{O}}\\
					&\leqslant C.
				\end{align*}
				Thus
				$$\left|c_{i_{k}}(\lambda,\omega)-c_{i_{k}}(\lambda',\omega')\right|\leqslant C\gamma\gamma^{-1}N_{k+1}^{-\overline{a}}\leqslant C\gamma N_{k+1}^{1-\overline{a}}.
				$$
				Therefore, we obtain from the definition of $\gamma_{k}$ and $\upsilon\in(0,1)$
				\begin{align*}
					\left|\omega\cdot l+c_{i_{k}}(\lambda,\omega)\right|&>\tfrac{4\gamma_{k+1}^{\upsilon}\langle j\rangle}{\langle l\rangle^{\tau_{1}}}-C\gamma\langle j\rangle N_{k+1}^{1-\overline{a}}\\
					&\geqslant\tfrac{4\gamma_{k+1}^{\upsilon}\langle j\rangle}{2^{\upsilon}\langle l\rangle^{\tau_{1}}}\left(2^{\upsilon}-C2^{\upsilon}N_{k+1}^{\tau_{1}+1-\overline{a}}\right).
				\end{align*}
				By the  choice of $\overline{a}$ made  in  \eqref{cond-abarre-1}, we can ensure, up to take $N_{0}$ sufficiently large,
				$$CN_{k+1}^{\tau_{1}+1-\overline{a}}\leqslant CN_{0}^{-1}<1-2^{-\upsilon},$$
				so that
				$$\left|\omega\cdot l+jc_{i_{k}}(\lambda,\omega)\right|>\tfrac{4\gamma_{k+1}^{\upsilon}\langle j\rangle}{2^{\upsilon}\langle l\rangle^{\tau_{1}}}\cdot
				$$
				As a consequence, $(\lambda,\omega)\in\mathcal{O}_{\infty,k}^{\frac{\gamma_{k+1}}{2},\tau_{1}}(i_{k}).$ Let us now check that $(\lambda,\omega)\in\mathscr{O}_{\infty,k}^{\frac{\gamma_{k+1}}{2},\tau_{1},\tau_{2}}(i_{k}).$ For all $(l,j,j_0)\in\mathbb{Z}^{d}\times(\mathbb{S}_{0}^{c})^2$ with $|l|\leqslant N_{k},$ we have by the  triangle and Cauchy-Schwarz inequalities together with $(\lambda',\omega')\in\mathscr{O}_{\infty,k}^{\gamma_{k+1},\tau_{1},\tau_{2}}(i_{k})$
				\begin{align*}
					\left|\omega\cdot l+\mu_{j}^{\infty}(\lambda,\omega,i_{k})-\mu_{j_{0}}^{\infty}(\lambda,\omega,i_{k})\right|&\geqslant\left|\omega'\cdot l+\mu_{j}^{\infty}(\lambda',\omega',i_{k})-\mu_{j_{0}}^{\infty}(\lambda',\omega',i_{k})\right|-|\omega-\omega'||l|\\
					&\quad-\left|\mu_{j}^{\infty}(\lambda,\omega,i_{k})-\mu_{j_{0}}^{\infty}(\lambda,\omega,i_{k})+\mu_{j_{0}}^{\infty}(\lambda',\omega',i_{k})-\mu_{j}^{\infty}(\lambda',\omega',i_{k})\right|\\
					&>\tfrac{2\gamma_{k+1}\langle j-j_{0}\rangle}{\langle l\rangle^{\tau_{2}}}-2\gamma N_{k+1}^{1-\overline{a}}\\
					&\quad-\left|\mu_{j}^{\infty}(\lambda,\omega,i_{k})-\mu_{j_{0}}^{\infty}(\lambda,\omega,i_{k})+\mu_{j_{0}}^{\infty}(\lambda',\omega',i_{k})-\mu_{j}^{\infty}(\lambda',\omega',i_{k})\right|.
				\end{align*}
				We recall by virtue of  Proposition \ref{reduction of the remainder term} that
				$$\mu_{j}^{\infty}(\lambda,\omega,i_{k})=\mu_{j}^{0}(\lambda,\omega,i_{k})+r_{j}^{\infty}(\lambda,\omega,i_{k}).$$
				Thus
				\begin{align*}
					&\left|\mu_{j}^{\infty}(\lambda,\omega,i_{k})-\mu_{j_{0}}^{\infty}(\lambda,\omega,i_{k})+\mu_{j_{0}}^{\infty}(\lambda',\omega',i_{k})-\mu_{j}^{\infty}(\lambda',\omega',i_{k})\right|\\
					&\leqslant\left|\mu_{j}^{0}(\lambda,\omega,i_{k})-\mu_{j_{0}}^{0}(\lambda,\omega,i_{k})+\mu_{j_{0}}^{0}(\lambda',\omega',i_{k})-\mu_{j}^{0}(\lambda',\omega',i_{k})\right|\\
					&\quad+\left|r_{j}^{\infty}(\lambda,\omega,i_{k})-r_{j}^{\infty}(\lambda',\omega',i_{k})\right|+\left|r_{j_{0}}^{\infty}(\lambda,\omega,i_{k})-r_{j_{0}}^{\infty}(\lambda',\omega',i_{k})\right|.
				\end{align*}
				According to  the  Mean Value Theorem, \eqref{reg-G-10} and the definition of $\mathrm{O}_{k+1}^{2\gamma}$  we find
				$$\left|\mu_{j}^{0}(\lambda,\omega,i_{k})-\mu_{j_{0}}^{0}(\lambda,\omega,i_{k})+\mu_{j_{0}}^{0}(\lambda',\omega',i_{k})-\mu_{j}^{0}(\lambda',\omega',i_{k})\right|\leqslant \gamma CN_{k+1}^{1-\overline{a}}\langle j-j_{0}\rangle.$$
				Applying once again the  Mean Value Theorem, \eqref{estimate rjinfty}, \eqref{small-NM} and the definition of $\mathrm{O}_{n+1}^{2\gamma}$ yields
				$$\left|r_{j}^{\infty}(\lambda,\omega,i_{k})-r_{j}^{\infty}(\lambda',\omega',i_{k})\right|\leqslant C\gamma N_{k+1}^{-\overline{a}}\varepsilon \gamma^{-2}\leqslant \gamma CN_{k+1}^{1-\overline{a}}\langle j-j_{0}\rangle.$$
				Putting together the foregoing estimates and the facts that $|l|\leqslant N_{k}$ and $\gamma_{k+1}\geqslant\gamma$ we infer
				\begin{align*}
					\left|\omega\cdot l+\mu_{j}^{\infty}(\lambda,\omega,i_{k})-\mu_{j_{0}}^{\infty}(\lambda,\omega,i_{k})\right|&\geqslant\tfrac{\gamma_{k+1}\langle j-j_{0}\rangle}{\langle l\rangle^{\tau_{2}}}\left(2-CN_{k+1}^{\tau_{2}+1-\overline{a}}\right).
				\end{align*}
				By virtue of \eqref{cond-abarre-1} and taking 
				$N_{0}$ sufficiently large we get
				$$CN_{n}^{\tau_{2}+1-\overline{a}}\leqslant CN_0^{-1}<1.
				$$
				This implies 
				$$\left|\omega\cdot l+\mu_{j}^{\infty}(\lambda,\omega,i_{k})-\mu_{j_{0}}^{\infty}(\lambda,\omega,i_{k})\right|>\tfrac{2\gamma_{k+1}\langle j-j_{0}\rangle}{\langle l\rangle^{\tau_{2}}}.$$
				As a consequence, we deduce that $(\lambda,\omega)\in\mathscr{O}_{\infty,k}^{\frac{\gamma_{k+1}}{2},\tau_{1},\tau_{2}}(i_{n}).$ Finally, $(\lambda,\omega)\in\mathcal{G}_{k}\big(\frac{\gamma_{k+1}}{2},\tau_{1},\tau_{2},i_{k}\big)$ and therefore $(\lambda,\omega)\in\mathcal{A}_{k+1}^{\frac{\gamma}{2}}.$ This achieves the induction proof of \eqref{hyprec O in A}.\\
				
				\noindent $\blacktriangleright$ \textbf{Construction of the next approximation.} We are now going to construct  the next  approximation $U_{n+1}$ by using  a modified Nash-Moser scheme. The assumption \eqref{small-1} being satisfied, we can apply Theorem \ref{proposition alomst approximate inverse} with $L_{n}$ and obtain the existence of an operator $\mathrm{T}_n:=\mathrm{T}_{n}(\lambda,\omega)$ well-defined in the whole set of parameters $\mathcal{O}$ and satisfying the following estimates
				\begin{equation}\label{estimate Tm}
					\forall s\in[s_{0},S],\quad\|\mathrm{T}_{n}\rho\|_{q,s}^{\gamma,\mathcal{O}}\lesssim\gamma^{-1}\left(\|\rho\|_{q,s+\overline{\sigma}}^{\gamma,\mathcal{O}}+\|\mathfrak{I}_{n}\|_{q,s+\overline{\sigma}}^{\gamma,\mathcal{O}}\|\rho\|_{q,s_{0}+\overline{\sigma}}^{\gamma,\mathcal{O}}\right)
				\end{equation}
				and
				\begin{equation}\label{estimate Tm in norm s0}
					\|\mathrm{T}_{n}\rho\|_{q,s_{0}}^{\gamma,\mathcal{O}}\lesssim\gamma^{-1}\|\rho\|_{q,s_{0}+\overline{\sigma}}^{\gamma,\mathcal{O}}.
				\end{equation}
				Moreover, when it is  restricted to the Cantor set $\mathcal{G}_{n}(\gamma_{n+1},\tau_{1},\tau_{2},i_{n}),$ $\mathrm{T}_{n}$ is an approximate right inverse of $L_{n}$ with suitable tame estimates needed later, see Theorem \ref{proposition alomst approximate inverse}. 
				Next we define,
				$$\widetilde{U}_{n+1}:=U_{n}+\widetilde{H}_{n+1}\quad\mbox{ with }\quad \widetilde{H}_{n+1}:=(\widehat{\mathfrak{I}}_{n+1},\widehat{\alpha}_{n+1},0):=-\mathbf{\Pi}_{n}\mathrm{T}_{n}\Pi_{n}\mathcal{F}(U_{n})\in E_{n}\times\mathbb{R}^{d}\times\mathbb{R}^{d+1},
				$$
				where $\mathbf{\Pi}_{n}$ is defined by
				\begin{equation}\label{grproj}
					\mathbf{\Pi}_{n}(\mathfrak{I},\alpha,0)=(\Pi_{n}\mathfrak{I},\alpha,0)\quad \mbox{ and }\quad\mathbf{\Pi}_{n}^{\perp}(\mathfrak{I},\alpha,0)=(\Pi_{n}\mathfrak{I},0,0).
				\end{equation}
				Notice that the projectors $\Pi_n$ are reversibility preserving due to the symmetry with respect to the Fourier modes. Then, using the reversibility of $\mathrm{T}_n$ together with \eqref{rev_Halpha} and Lemma \ref{lem rev H-X}, one deduces from $\mathfrak{S}i_n(\varphi)=i_n(-\varphi)$ that
				\begin{equation}\label{sym frih} \mathfrak{S}\widehat{\mathfrak{I}}_{n+1}(\varphi)=\widehat{\mathfrak{I}}_{n+1}(-\varphi).
				\end{equation}
				Note that $U_{n}$ is defined in the full set $\mathcal{O}$ and so does $\widetilde{U}_{n+1}$. Nevertheless,
				we will not be working  with this natural extension but rather with a suitable  localized version of it around the Cantor set $\mathcal{A}_{n+1}^{\gamma}.$ Doing so, we shall get a nice decay property allowing the scheme to converge. Now, introduce  the quadratic function
				\begin{align}\label{Def-Qm}
					Q_{n}=\mathcal{F}(U_{n}+\widetilde{H}_{n+1})-\mathcal{F}(U_{n})-L_{n}\widetilde{H}_{n+1},
				\end{align}
				then simple transformations give
				\begin{align}\label{Decom-RTT1}
					\nonumber \mathcal{F}(\widetilde{U}_{n+1})& =  \mathcal{F}(U_{n})-L_{n}\mathbf{\Pi}_{n}\mathrm{T}_{n}\Pi_{n}\mathcal{F}(U_{n})+Q_{n}\\
					\nonumber& =  \mathcal{F}(U_{n})-L_{n}\mathrm{T}_{n}\Pi_{n}\mathcal{F}(U_{n})+L_{n}\mathbf{\Pi}_{n}^{\perp}\mathrm{T}_{n}\Pi_{n}\mathcal{F}(U_{n})+Q_{n}\\
					\nonumber& =  \mathcal{F}(U_{n})-\Pi_{n}L_{n}\mathrm{T}_{n}\Pi_{n}\mathcal{F}(U_{n})+(L_{n}\mathbf{\Pi}_{n}^{\perp}-\Pi_{n}^{\perp}L_{n})\mathrm{T}_{n}\Pi_{n}\mathcal{F}(U_{n})+Q_{n}\\
					& =  \Pi_{n}^{\perp}\mathcal{F}(U_{n})-\Pi_{n}(L_{n}\mathrm{T}_{n}-\textnormal{Id})\Pi_{n}\mathcal{F}(U_{n})+(L_{n}\mathbf{\Pi}_{n}^{\perp}-\Pi_{n}^{\perp}L_{n})\mathrm{T}_{n}\Pi_{n}\mathcal{F}(U_{n})+Q_{n}.
				\end{align} In the sequel we shall prove 
				$$\|\mathcal{F}(U_{n+1})\|_{q,s_{0}}^{\gamma,\mathrm{O}_{n+1}^{\gamma}}\leqslant C_{\ast}\varepsilon N_{n}^{-a_{1}},$$
				with $U_{n+1}$ a suitable extension of $\widetilde{U}_{n+1}|_{\mathrm{O}_{n+1}^{\gamma}}.$\\
				
				\noindent $\blacktriangleright$ \textbf{Estimates of $\mathcal{F}(\widetilde{U}_{n+1})$.}
				We shall now estimate $\mathcal{F}(\widetilde{U}_{n+1})$ with the norm $\|\cdot\|_{q,s_{0}}^{\gamma,\mathrm{O}_{n+1}^{2\gamma}}$ by using \eqref{Decom-RTT1}. The localization in $\mathrm{O}_{n+1}^{2\gamma}$ is required for the classical extension in the next point, see \eqref{ext Hn+1}.
				\\
				\ding{226} \textit{Estimate of $\Pi_{n}^{\perp}\mathcal{F}(U_{n}).$} 
				We apply Taylor formula combined with \eqref{main function} and Lemma \ref{tame estimates for the vector field XmathcalPvarepsilon} together with \eqref{estimate mathcalF(U0)} and $(\mathcal{P}1)_{n}.$ Therefore, we obtain  
				\begin{align}\label{link mathcalF(Um) and Wm}
					\nonumber \forall s\geqslant s_{0},\quad\|\mathcal{F}(U_{n})\|_{q,s}^{\gamma,\mathrm{O}_{n}^{\gamma}}&\leqslant\|\mathcal{F}(U_{0})\|_{q,s}^{\gamma,\mathcal{O}}+\|\mathcal{F}(U_{n})-\mathcal{F}(U_{0})\|_{q,s}^{\gamma,\mathrm{O}_{n}^{\gamma}}\\
					&\lesssim\varepsilon+\| W_{n}\|_{q,s+\overline{\sigma}}^{\gamma,\mathcal{O}}.
				\end{align}
				As a consequence, \eqref{choice of gamma and N0 in the Nash-Moser} and $(\mathcal{P}1)_{n}$ imply
				\begin{equation}\label{estimate mathcalF(U0) in norm s0}
					\gamma^{-1}\|\mathcal{F}(U_{n})\|_{q,s_{0}}^{\gamma,\mathrm{O}_{n}^{\gamma}}\leqslant 1.
				\end{equation}
				From Lemma \ref{Lem-lawprod}-(ii) and \eqref{link mathcalF(Um) and Wm}, we get
				\begin{align}\label{pipFuns0}
					\|\Pi_{n}^{\perp}\mathcal{F}(U_{n})\|_{q,s_{0}}^{\gamma,\mathrm{O}_{n}^{\gamma}}&\leqslant N_{n}^{s_{0}-b_{1}}\|\mathcal{F}(U_{n})\|_{q,b_{1}}^{\gamma,\mathrm{O}_{n}^{\gamma}}\nonumber\\
					&\lesssim  N_{n}^{\overline{\sigma}-b_{1}}\left(\varepsilon+\|W_{n}\|_{q,b_{1}+\overline{\sigma}}^{\gamma,\mathcal{O}}\right).
				\end{align}
				Now, $(\mathcal{P}3)_{n}$ together \eqref{definition of Nm} and \eqref{choice of gamma and N0 in the Nash-Moser} yield 
				\begin{align}\label{Wm in high norm}
					\nonumber\varepsilon+\|W_{n}\|_{q,b_{1}+\overline{\sigma}}^{\gamma,\mathcal{O}}&\leqslant\varepsilon\left(1+C_{\ast}\gamma^{-1}N_{n-1}^{\mu_{1}}\right)\\
					&\leqslant 2C_{\ast}\varepsilon N_{n}^{\frac{2}{3}\mu_{1}+1}.
				\end{align}
				By putting together \eqref{Wm in high norm} and  \eqref{pipFuns0} and by making appeal to \eqref{inco}, we infer for any $n\in\mathbb{N}$,
				\begin{align}\label{final estimate PiperpF(Um)}
					\|\Pi_{n}^{\perp}\mathcal{F}(U_{n})\|_{q,s_{0}}^{\gamma,\mathrm{O}_{n+1}^{2\gamma}} & \leqslant  \|\Pi_{n}^{\perp}\mathcal{F}(U_{n})\|_{q,s_{0}}^{\gamma,\mathrm{O}_{n}^{\gamma}}\nonumber\\
					&\lesssim C_{\ast}\varepsilon N_{n}^{s_{0}+\frac{2}{3}\mu_{1}+1-b_{1}}.
				\end{align}
				Remark that one also obtains, combining \eqref{link mathcalF(Um) and Wm} and \eqref{Wm in high norm},
				\begin{equation}\label{HDP10}
					\|\mathcal{F}(U_{n})\|_{q,b_{1}+\overline{\sigma}}^{\gamma,\mathrm{O}_{n}^{\gamma}}  \leqslant  C_{\ast}\varepsilon N_{n}^{\overline{\sigma}+\frac{2}{3}\mu_{1}+1}.
				\end{equation}
				\\
				\ding{226} \textit{Estimate of $\Pi_{n}(L_{n}\mathrm{T}_{n}-\textnormal{Id})\Pi_{n}\mathcal{F}(U_{n})$.} In view of \eqref{hyprec O in A}, one has
				$$\mathrm{O}_{n+1}^{2\gamma}\subset\mathcal{A}_{n+1}^{\frac{\gamma}{2}}\subset\mathcal{G}_{n}\Big(\frac{\gamma_{n+1}}{2},\tau_{1},\tau_{2},i_{n}\Big).$$
				Then, applying Theorem  \ref{proposition alomst approximate inverse}, we can write
				$$\Pi_{n}(L_{n}\mathrm{T}_{n}-\textnormal{Id})\Pi_{n}\mathcal{F}(U_{n})=\mathscr{E}_{1,n}+\mathscr{E}_{2,n}+\mathscr{E}_{3,n},$$
				with 
				\begin{align*}
					\mathscr{E}_{1,n}:=\Pi_{n}{\mathcal{E}_{1}^{(n)}}\Pi_{n}\mathcal{F}(U_{n}),\\
					\mathscr{E}_{2,n}:=\Pi_{n}\mathcal{E}_{2}^{(n)}\Pi_{n}\mathcal{F}(U_{n}),\\
					\mathscr{E}_{3,n}:=\Pi_{n}\mathcal{E}_{3}^{(n)}\Pi_{n}\mathcal{F}(U_{n})
				\end{align*}
				where $\mathcal{E}_{1}^{(n)}$, $\mathcal{E}_{2}^{(n)}$ and $\mathcal{E}_{3}^{(n)}$ satisfy the estimates \eqref{eaai-E1}, \eqref{eaai-E2} and \eqref{eaai-E3} respectively. By \eqref{inco}, we get
				\begin{equation}\label{e-ai-NM}
					\|\Pi_{n}(L_{n}\mathrm{T}_{n}-\textnormal{Id})\Pi_{n}\mathcal{F}(U_{n})\|_{q,s_0}^{\gamma,\mathrm{O}_{n+1}^{2\gamma}}\leqslant\|\mathscr{E}_{1,n}\|_{q,s_0}^{\gamma,\mathrm{O}_{n}^{\gamma}}+\|\mathscr{E}_{2,n}\|_{q,s_0}^{\gamma,\mathrm{O}_{n}^{\gamma}}+\|\mathscr{E}_{3,n}\|_{q,s_0}^{\gamma,\mathrm{O}_{n}^{\gamma}}.
				\end{equation}
				We shall first focus  on $\mathscr{E}_{1,n}.$ We need the following interpolation-type inequality
				\begin{align}\label{int-hln}
					\|\mathcal{F}(U_{n})\|_{q,s_{0}+\overline{\sigma}}^{\gamma,\mathrm{O}_{n}^{\gamma}}&\leqslant\|\Pi_{n}\mathcal{F}(U_n)\|_{q,s_0+\overline{\sigma}}^{\gamma,\mathrm{O}_{n}^{\gamma}}+\|\Pi_{n}^{\perp}\mathcal{F}(U_n)\|_{q,s_0+\overline{\sigma}}^{\gamma,\mathrm{O}_{n}^{\gamma}}\nonumber\\
					&\leqslant N_{n}^{\overline{\sigma}}\|\mathcal{F}(U_n)\|_{q,s_0}^{\gamma,\mathrm{O}_{n}^{\gamma}}+N_{n}^{s_0-b_{1}}\|\mathcal{F}(U_n)\|_{q,b_{1}+\overline{\sigma}}^{\gamma,\mathrm{O}_{n}^{\gamma}}.
				\end{align}
				Combining \eqref{eaai-E1}, \eqref{int-hln}, $(\mathcal{P}_{1})_n$, \eqref{small-NM} and \eqref{HDP10}, we obtain
				\begin{align}\label{Esc1n}
					\|\mathscr{E}_{1,n}\|_{q,s_0}^{\gamma,\mathrm{O}_{n}^{\gamma}}&\lesssim\gamma^{-1}\|\mathcal{F}(U_n)\|_{q,s_0+\overline{\sigma}}^{\gamma,\mathrm{O}_{n}^{\gamma}}\|\Pi_{n}\mathcal{F}(U_n)\|_{q,s_0+\overline{\sigma}}^{\gamma,\mathrm{O}_{n}^{\gamma}}\left(1+\|\mathfrak{I}_n\|_{q,s_0+\overline{\sigma}}^{\gamma,\mathcal{O}}\right)\nonumber\\
					&\lesssim \gamma^{-1}N_n^{\overline{\sigma}}\left(N_n^{\overline{\sigma}}\|\mathcal{F}(U_n)\|_{q,s_0}^{\gamma,\mathrm{O}_{n}^{\gamma}}+N_{n}^{s_0-b_{1}}\|\mathcal{F}(U_n)\|_{q,b_{1}+\overline{\sigma}}^{\gamma,\mathrm{O}_{n}^{\gamma}}\right)\|\mathcal{F}(U_n)\|_{q,s_0}^{\gamma,\mathrm{O}_{n}^{\gamma}}\left(1+\|W_{n}\|_{q,s_0+\overline{\sigma}}^{\gamma,\mathcal{O}}\right)\nonumber\\
					&\lesssim C_{\ast}\varepsilon\left(N_{n}^{2\overline{\sigma}-\frac{4}{3}a_1}+N_n^{s_0+2\overline{\sigma}+\frac{2}{3}\mu_{1}+1-\frac{2}{3}a_1-b_1}\right).
				\end{align}
				We now turn to $\mathscr{E}_{2,n}$ and $\mathscr{E}_{3,n}.$ Applying \eqref{eaai-E2} with $b=b_{1}-s_0$ and using \eqref{small-NM}, $(\mathcal{P}_2)_n$ and $(\mathcal{P}_3)_n$, we get
				\begin{align}\label{Esc2n}
					\|\mathscr{E}_{2,n}\|_{q,s_0}^{\gamma,\mathrm{O}_{n}^{\gamma}}&\lesssim\gamma^{-1}N_n^{s_0-b_1}\left(\|\Pi_{n}\mathcal{F}(U_n)\|_{q,b_1+\overline{\sigma}}^{\gamma,\mathrm{O}_{n}^{\gamma}}+\varepsilon\|\mathfrak{I}_n\|_{q,b_1+\overline{\sigma}}^{\gamma,\mathcal{O}}\|\Pi_{n}\mathcal{F}(U_n)\|_{q,s_0+\overline{\sigma}}^{\gamma,\mathrm{O}_{n}^{\gamma}}\right)\nonumber\\
					&\lesssim \gamma^{-1}N_n^{s_0-b_1}\left(\|\mathcal{F}(U_n)\|_{q,b_1+\overline{\sigma}}^{\gamma,\mathrm{O}_{n}^{\gamma}}+\varepsilon N_n^{\overline{\sigma}}\|W_n\|_{q,b_1+\overline{\sigma}}^{\gamma,\mathcal{O}}\|\mathcal{F}(U_n)\|_{q,s_0}^{\gamma,\mathrm{O}_{n}^{\gamma}}\right)\nonumber\\
					&\lesssim C_{\ast}\varepsilon N_n^{s_0+\overline{\sigma}+\frac{2}{3}\mu_{1}+2-b_1}+C_{\ast}\varepsilon N_n^{s_0+\overline{\sigma}+\frac{2}{3}\mu_{1}+2-\frac{2}{3}a_{1}-b_1}\nonumber\\
					&\lesssim C_{\ast}\varepsilon N_n^{s_0+\overline{\sigma}+\frac{2}{3}\mu_{1}+2-b_1}.
				\end{align}
				Using the same techniques together with \eqref{eaai-E3}, \eqref{definition of Nm}, \eqref{choice of gamma and N0 in the Nash-Moser} and \eqref{small-NM}, we infer
				\begin{align}\label{Esc3n}
					\|\mathscr{E}_{3,n}\|_{q,s_0}^{\gamma,\mathrm{O}_{n}^{\gamma}}&\lesssim N_n^{s_0-b_1}\gamma^{-2}\left(\|\Pi_n\mathcal{F}(U_n)\|_{q,b_{1}+\overline{\sigma}}^{\gamma,\mathrm{O}_{n}^{\gamma}}+\varepsilon\gamma^{-2}\|\mathfrak{I}_n\|_{q,b_{1}+\overline{\sigma}}^{\gamma,\mathcal{O}}\|\Pi_n\mathcal{F}(U_n)\|_{q,s_0+\overline{\sigma}}^{\gamma,\mathrm{O}_{n}^{\gamma}}\right)\nonumber\\
					&\quad+\varepsilon\gamma^{-4}N_{0}^{\mu_{2}}N_n^{-\mu_{2}}\|\Pi_n\mathcal{F}(U_n)\|_{q,s_0+\overline{\sigma}}^{\gamma,\mathrm{O}_{n}^{\gamma}}\nonumber\\
					&\lesssim C_{\ast}\varepsilon\left(N_n^{s_0+\overline{\sigma}+\frac{2}{3}\mu_{1}+2-b_1}+N_n^{\overline{\sigma}+1-\mu_{2}-\frac{2}{3}a_1}\right).
				\end{align}
				Putting together \eqref{e-ai-NM}, \eqref{Esc1n}, \eqref{Esc2n} and \eqref{Esc2n}, we obtain
				\begin{equation}\label{est-aait}
					\|\Pi_{n}(L_{n}\mathrm{T}_{n}-\textnormal{Id})\Pi_{n}\mathcal{F}(U_{n})\|_{q,s_0}^{\gamma,\mathrm{O}_{n+1}^{2\gamma}}\leqslant CC_{\ast}\varepsilon\left(N_{n}^{2\overline{\sigma}-\frac{4}{3}a_1}+N_n^{s_0+2\overline{\sigma}+\frac{2}{3}\mu_{1}+1-b_1}+N_n^{\overline{\sigma}+1-\mu_{2}-\frac{2}{3}a_1}\right).
				\end{equation}
				For $n=0$, we deduce from \eqref{estimate mathcalF(U0)},\eqref{small-NM} and by slight modifications of the preceding  computations 
				\begin{align}\label{e-ai-0}
					\|\Pi_{0}(L_{0}\mathrm{T}_{0}-\textnormal{Id})\Pi_{0}\mathcal{F}(U_{0})\|_{q,s_{0}}^{\gamma,\mathrm{O}_{1}^{2\gamma}}&\leqslant\|\mathscr{E}_{1,0}\|_{q,s_0}^{\gamma,\mathrm{O}_{n}^{\gamma}}+\|\mathscr{E}_{2,0}\|_{q,s_0}^{\gamma,\mathrm{O}_{n}^{\gamma}}+\|\mathscr{E}_{3,0}\|_{q,s_0}^{\gamma,\mathrm{O}_{n}^{\gamma}}\nonumber\\
					& \lesssim \varepsilon^{2}\gamma^{-1}+\varepsilon\gamma^{-1}+\big(\varepsilon \gamma^{-2}N_{0}^{s_{0}-b_{1}}+\varepsilon^{2}\gamma^{-4}\big)\nonumber\\
					&\lesssim\varepsilon\gamma^{-2}.
				\end{align}
				\ding{226} \textit{Estimate of $\big(L_{n}\mathbf{\Pi}_{n}^{\perp}-\Pi_{n}^{\perp}L_{n}\big)\mathrm{T}_{n}\Pi_{n}\mathcal{F}(U_{n}).$} Combining \eqref{lin func} and \eqref{main function}, we get for $H=(\widehat{\mathfrak{I}},\widehat{\alpha})$ with $\widehat{\mathfrak{I}}=(\widehat{\Theta},\widehat{I},\widehat{z}),$
				\begin{equation}\label{exp-LmH}
					L_{n}H=\omega\cdot\partial_{\varphi}\widehat{\mathfrak{I}}-(0,0,\partial_{\theta}\mathrm{L}(\lambda)\widehat{z})-\varepsilon d_{i}X_{\mathcal{P}_{\varepsilon}}(i_{n})\widehat{\mathfrak{I}}-(\widehat{\alpha},0,0).
				\end{equation}
				Using \eqref{grproj} and the fact that $\omega\cdot\partial_{\varphi}$ and $\partial_{\theta}\mathrm{L}(\lambda)$ are diagonal leading to $[\Pi_{n}^{\perp},\omega\cdot\partial_{\varphi}]=[\Pi_{n}^{\perp},\partial_{\theta}\mathrm{L}(\lambda)]=0$, one has for $H=(\widehat{\mathfrak{I}},\widehat{\alpha}),$ 
				$$\left(L_{n}\mathbf{\Pi}_{n}^{\perp}-\Pi_{n}^{\perp}L_{n}\right)H=-\varepsilon[d_{i}X_{\mathcal{P}_{\varepsilon}}(i_{n}),\Pi_{n}^{\perp}]\widehat{\mathfrak{I}}.$$
				In view of  Lemma \ref{tame estimates for the vector field XmathcalPvarepsilon}-{(ii)}, Lemma \ref{properties of Toeplitz in time operators}, \eqref{inco} and $(\mathcal{P}1)_{n}$ we get
				$$
				\left\|\left(L_{n}\mathbf{\Pi}_{n}^{\perp}-\Pi_{n}^{\perp}L_{n}\right)H\right\|_{q,s_{0}}^{\gamma,\mathrm{O}_{n+1}^{2\gamma}}\lesssim\varepsilon N_{n}^{s_{0}-b_{1}}\left(\|\widehat{\mathfrak{I}}\|_{q,b_{1}+1}^{\gamma,\mathrm{O}_{n}^{\gamma}}+\|\mathfrak{I}_{n}\|_{q,b_{1}+\overline{\sigma}}^{\gamma,\mathcal{O}}\|\widehat{\mathfrak{I}}\|_{q,s_{0}+1}^{\gamma,\mathrm{O}_{n}^{\gamma}}\right).
				$$
				Consequently,
				\begin{align*}
					\mathtt{N}_{\textnormal{\tiny{com}}}(s_0):=\left\|\left(L_{n}\mathbf{\Pi}_{n}^{\perp}-\Pi_{n}^{\perp}L_{n}\right)\mathrm{T}_{n}\Pi_{n}\mathcal{F}(U_{n})\right\|_{q,s_{0}}^{\gamma,\mathrm{O}_{n+1}^{2\gamma}}&\lesssim\varepsilon N_{n}^{s_{0}-b_{1}}\|\mathrm{T}_{n}\Pi_{n}\mathcal{F}(U_{n})\|_{q,b_{1}+1}^{\gamma,\mathrm{O}_{n}^{\gamma}}\\
					&\quad +\varepsilon N_{n}^{s_{0}-b_{1}}\|\mathfrak{I}_{n}\|_{q,b_{1}+\overline{\sigma}}^{\gamma,\mathcal{O}}\|\mathrm{T}_{n}\Pi_{n}\mathcal{F}(U_{n})\|_{q,s_{0}+1}^{\gamma,\mathrm{O}_{n}^{\gamma}}.
				\end{align*}
				Hence, gathering \eqref{estimate Tm},  Lemma \ref{Lem-lawprod}, Sobolev embeddings, \eqref{small-NM}, \eqref{choice of gamma and N0 in the Nash-Moser} and $(\mathcal{P}1)_{n}$ yields
				\begin{align*}
					\mathtt{N}_{\textnormal{\tiny{com}}}(s_0)&\lesssim \varepsilon\gamma^{-1}N_{n}^{s_{0}-b_{1}}\|\Pi_{n}\mathcal{F}(U_{n})\|_{q,b_{1}+\overline{\sigma}+1}^{\gamma,\mathrm{O}_{n}^{\gamma}}+\|\mathfrak{I}_{n}\|_{q,b_{1}+\overline{\sigma}+1}^{\gamma,\mathcal{O}}\|\Pi_{n}\mathcal{F}(U_{n})\|_{q,s_{0}+\overline{\sigma}}^{\gamma,\mathrm{O}_{n}^{\gamma}}\\
					&\quad+\varepsilon\gamma^{-1}N_{n}^{s_{0}-b_{1}}\|\mathfrak{I}_{n}\|_{q,b_{1}+\overline{\sigma}}^{\gamma,\mathcal{O}}\left(\|\Pi_{n}\mathcal{F}(U_{n})\|_{q,s_{0}+\overline{\sigma}+1}^{\gamma,\mathrm{O}_{n}^{\gamma}}+\|\mathfrak{I}_{n}\|_{q,s_{0}+\overline{\sigma}+1}^{\gamma,\mathcal{O}}\|\Pi_{n}\mathcal{F}(U_{n})\|_{q,s_{0}+\overline{\sigma}}^{\gamma,\mathrm{O}_{n}^{\gamma}}\right)\\
					&\lesssim \varepsilon N_{n}^{s_{0}+2-b_{1}}\left(\|\mathcal{F}(U_{n})\|_{q,b_{1}+\overline{\sigma}}^{\gamma,\mathrm{O}_{n}^{\gamma}}+\| W_{n}\|_{q,b_{1}+\overline{\sigma}}^{\gamma,\mathcal{O}}\|\Pi_{n}\mathcal{F}(U_{n})\|_{q,s_{0}+\overline{\sigma}}^{\gamma,\mathrm{O}_{n}^{\gamma}}\right).
				\end{align*}
				Applying Lemma \ref{Lem-lawprod}-(ii), $(\mathcal{P}2)_{n}$ and \eqref{definition of Nm}, we infer
				\begin{align*}
					\|\Pi_{n}\mathcal{F}(U_{n})\|_{q,s_{0}+\overline{\sigma}}^{\gamma,\mathrm{O}_{n}^{\gamma}}&\leqslant N_{n}^{\overline{\sigma}}\|\mathcal{F}(U_{n})\|_{q,s_{0}}^{\gamma,\mathrm{O}_{n}^{\gamma}}\\
					&\leqslant C_{\ast}\varepsilon N_{n}^{\overline{\sigma}}N_{n-1}^{-a_{1}}\\
					&\leqslant C_{\ast}\varepsilon N_{n}^{\overline{\sigma}-\frac{2}{3}a_{1}}.
				\end{align*}
				Added to \eqref{param NM}, \eqref{HDP10} and $(\mathcal{P}3)_{n}$, we obtain for  $n\in\mathbb{N}$,
				\begin{equation}\label{final estimate commutator}
					\|(L_{n}\mathbf{\Pi}_{n}^{\perp}-\Pi_{n}^{\perp}L_{n})\mathrm{T}_{n}\Pi_{n}\mathcal{F}(U_{n})\|_{q,s_{0}}^{\gamma,\mathrm{O}_{n+1}^{2\gamma}}\leqslant CC_{\ast}\varepsilon N_{n}^{s_{0}+\overline{\sigma}+\frac{2}{3}\mu_{1}+3-b_{1}}.
				\end{equation}
				\ding{226} \textit{Estimate of $Q_{n}$.} We apply Taylor formula together with \eqref{Def-Qm} leading to
				$$Q_{n}=\int_{0}^{1}(1-t)d_{i,\alpha}^{2}\mathcal{F}(U_{n}+t\widetilde{H}_{n+1})[\widetilde{H}_{n+1},\widetilde{H}_{n+1}]dt.$$
				Thus,   \eqref{exp-LmH} and Lemma \ref{tame estimates for the vector field XmathcalPvarepsilon}-{(iii)} allow to get
				\begin{align}\label{mahma-YDa1}
					\| Q_{n}\|_{q,s_{0}}^{\gamma,\mathrm{O}_{n+1}^{2\gamma}}\lesssim\varepsilon\left(1+\|W_{n}\|_{q,s_{0}+2}^{\gamma,\mathcal{O}}+\| \widetilde{H}_{n+1}\|_{q,s_{0}+2}^{\gamma,\mathrm{O}_{n+1}^{2\gamma}}\right)\left(\| \widetilde{H}_{n+1}\|_{q,s_{0}+2}^{\gamma,\mathrm{O}_{n+1}^{2\gamma}}\right)^{2}.
				\end{align}
				Combining \eqref{hyprec O in A}, \eqref{estimate Tm}, \eqref{link mathcalF(Um) and Wm} and \eqref{estimate mathcalF(U0) in norm s0}, we find for all $s\in[s_{0},S]$ 
				\begin{align}\label{link Hm+1 and Wm}
					\| \widetilde{H}_{n+1}\|_{q,s}^{\gamma,\mathrm{O}_{n+1}^{2\gamma}} & =  \|\mathbf{\Pi}_{n}\mathrm{T}_{n}\Pi_{n}\mathcal{F}(U_{n})\|_{q,s}^{\gamma,\mathrm{O}_{n+1}^{2\gamma}}\nonumber\\
					& \lesssim  \gamma^{-1}\left(\|\Pi_{n}\mathcal{F}(U_{n})\|_{q,s+\overline{\sigma}}^{\gamma,\mathrm{O}_{n}^{\gamma}}+\|\mathfrak{I}_{n}\|_{q,s+\overline{\sigma}}^{\gamma,\mathcal{O}}\|\Pi_{n}\mathcal{F}(U_{n})\|_{q,s_{0}+\overline{\sigma}}^{\gamma,\mathrm{O}_{n}^{\gamma}}\right)\nonumber\\
					& \lesssim  \gamma^{-1}\left(N_{n}^{\overline{\sigma}}\|\mathcal{F}(U_{n})\|_{q,s}^{\gamma,\mathrm{O}_{n}^{\gamma}}+N_{n}^{2\overline{\sigma}}\|\mathfrak{I}_{n}\|_{q,s}^{\gamma,\mathcal{O}}\|\mathcal{F}(U_{n})\|_{q,s_{0}}^{\gamma,\mathrm{O}_{n}^{\gamma}}\right)\nonumber\\
					& \lesssim  \gamma^{-1}N_{n}^{2\overline{\sigma}}\left(\varepsilon+\| W_{n}\|_{q,s}^{\gamma,\mathcal{O}}\right).
				\end{align}
				In the same way, according to \eqref{estimate Tm in norm s0}, $(\mathcal{P}1)_{n}$ and $(\mathcal{P}2)_{n}$, we infer
				\begin{align}\label{link Hm+1 and mathcalF(Um) in norm s0}
					\nonumber \| \widetilde{H}_{n+1}\|_{q,s_{0}}^{\gamma,\mathrm{O}_{n+1}^{2\gamma}}&\lesssim\gamma^{-1}N_{n}^{\overline{\sigma}}\|\mathcal{F}(U_{n})\|_{q,s_{0}}^{\gamma,\mathrm{O}_{n}^{2\gamma}}\\
					&\lesssim C_{\ast}\varepsilon\gamma^{-1} N_{n}^{\overline{\sigma}}N_{n-1}^{-a_{1}}.
				\end{align}
				Choosing $\varepsilon$ small enough and using $(\mathcal{P}1)_{n}$ and \eqref{link Hm+1 and mathcalF(Um) in norm s0},  we find
				\begin{align*}
					\|W_{n}\|_{q,s_{0}+2}^{\gamma,\mathcal{O}}+\| \widetilde{H}_{n+1}\|_{q,s_{0}+2}^{\gamma,\mathrm{O}_{n+1}^{2\gamma}} & \leqslant C_{\ast}\varepsilon\gamma^{-1}+N_{n}^{2}\| \widetilde{H}_{n+1}\|_{q,s_{0}}^{\gamma,\mathrm{O}_{n+1}^{2\gamma}}\\
					& \leqslant 1+C\varepsilon\gamma^{-1}N_{n}^{\overline{\sigma}+2}N_{n-1}^{-a_{1}}\\
					& \leqslant 1+C\varepsilon\gamma^{-1}N_{n-1}^{3+\frac{3}{2}\overline{\sigma}-a_{1}}.
				\end{align*}
				Now notice that \eqref{param NM} implies
				\begin{equation}\label{hyp-a1-1}
					a_{1}\geqslant 3+\tfrac{3}{2}\overline{\sigma}.
				\end{equation}
				Therefore, we obtain
				$$\|W_{n}\|_{q,s_{0}+2}^{\gamma,\mathcal{O}}+\| \widetilde{H}_{n+1}\|_{q,s_{0}+2}^{\gamma,\mathrm{O}_{n+1}^{2\gamma}}\leqslant 2.$$
				Hence,  plugging this estimate and \eqref{link Hm+1 and mathcalF(Um) in norm s0} into \eqref{mahma-YDa1}  and using \eqref{choice of gamma and N0 in the Nash-Moser} and \eqref{small-NM}, we find 
				\begin{align*}
					\| Q_{n}\|_{q,s_{0}}^{\gamma,\mathrm{O}_{n+1}^{2\gamma}} & \lesssim \varepsilon\left(\| \widetilde{H}_{n+1}\|_{q,s_{0}+2}^{\gamma,\mathrm{O}_{n+1}^{2\gamma}}\right)^{2}\\
					& \leqslant \varepsilon N_{n}^{4}\left(\| \widetilde{H}_{n+1}\|_{q,s_{0}}^{\gamma,\mathrm{O}_{n+1}^{2\gamma}}\right)^{2}\\
					& \lesssim \varepsilon C_{\ast}N_{n}^{2\overline{\sigma}+4}N_{n-1}^{-2a_1}.
				\end{align*}
				By using \eqref{definition of Nm}, we deduce when  $n\geqslant 1$,
				\begin{equation}\label{final estimate for Qm}
					\|Q_{n}\|_{q,s_{0}}^{\gamma,\mathcal{O}_{n+1}^{2\gamma}}\leqslant CC_{\ast}\varepsilon N_{n}^{2\overline{\sigma}+4-\frac{4}{3}a_{1}}.
				\end{equation}
				For $n=0$, we come back to \eqref{link Hm+1 and Wm} and \eqref{estimate mathcalF(U0)} to obtain for all $s\in[s_{0},S]$
				\begin{align}\label{H1}
					\|\widetilde{H}_{1}\|_{q,s}^{\gamma,\mathrm{O}_{1}^{2\gamma}}&\lesssim\gamma^{-1}\|\Pi_{0}\mathcal{F}(U_{0})\|_{q,s+\overline{\sigma}}^{\gamma,\mathrm{O}_{0}^{2\gamma}}\nonumber\\
					&\lesssim C_{\ast}\varepsilon\gamma^{-1}.
				\end{align}
				Finally, the inequality \eqref{final estimate for Qm} becomes for $n=0$,
				\begin{equation}\label{e-Q0}
					\|Q_{0}\|_{q,s_{0}}^{\gamma,\mathrm{O}_{0}^{2\gamma}}\lesssim C_{\ast}\varepsilon^{3}\gamma^{-2}.
				\end{equation}
				\ding{226} \textit{Conclusion.}
				Inserting  \eqref{final estimate PiperpF(Um)}, \eqref{est-aait}, \eqref{final estimate commutator} and \eqref{final estimate for Qm}, into \eqref{Decom-RTT1} implies for $n\in\mathbb{N}^{*}$,
				\begin{align*}
					\|\mathcal{F}(\widetilde{U}_{n+1})\|_{q,s_{0}}^{\gamma,\mathrm{O}_{n+1}^{2\gamma}}&\leqslant CC_{\ast}\varepsilon\left(N_{n}^{s_{0}+2\overline{\sigma}+\frac{2}{3}\mu_{1}+1-b_{1}}+N_{n}^{\overline{\sigma}+1-\mu_{2}-\frac{2}{3}a_{1}}+N_{n}^{2\overline{\sigma}+4-\frac{4}{3}a_{1}}\right).
				\end{align*} 
				The parameters conditions stated in \eqref{param NM} give
				\begin{equation}\label{Assump-DR1}\left\lbrace\begin{array}{rcl}
						s_{0}+2\overline{\sigma}+\frac{2}{3}\mu_{1}+2+a_{1} & \leqslant & b_{1}\\
						\overline{\sigma}+\frac{1}{3}a_{1}+2 & \leqslant & \mu_{2}\\
						
						2\overline{\sigma}+5& \leqslant & \frac{1}{3}a_{1}.
					\end{array}\right.
				\end{equation}
				Thus, by taking $N_{0}$ large enough, that is $\varepsilon$ small enough, we obtain for  $n\in\mathbb{N},$
				\begin{equation}\label{cond Nn}
					\left\lbrace\begin{array}{rcl}
						CN_{n}^{s_{0}+2\overline{\sigma}+\frac{2}{3}\mu_{1}+1-b_{1}} & \leqslant & \frac{1}{3}N_{n}^{-a_{1}}\\
						CN_{n}^{\overline{\sigma}+1-\mu_{2}-\frac{2}{3}a_{1}}& \leqslant & \frac{1}{3}N_{n}^{-a_{1}}\\
						CN_{n}^{2\overline{\sigma}+4-\frac{4}{3}a_{1}} & \leqslant & \frac{1}{3}N_{n}^{-a_{1}},
					\end{array}\right.
				\end{equation}
				which implies in turn
				that when  $n\in\mathbb{N}^{*}$,
				\begin{align}\label{Dabdoub1}
					\|\mathcal{F}(\widetilde{U}_{n+1})\|_{q,s_{0}}^{\gamma,\mathrm{O}_{n+1}^{2\gamma}}\leqslant C_{\ast}\varepsilon N_{n}^{-a_{1}}.
				\end{align}
				However, when  $n=0$, we plug \eqref{final estimate PiperpF(Um)}, \eqref{e-ai-0},  \eqref{final estimate commutator} and \eqref{e-Q0} into \eqref{Decom-RTT1}  in order to get
				$$\|\mathcal{F}(\widetilde{U}_{1})\|_{q,s_{0}}^{\gamma,\mathrm{O}_{1}^{2\gamma}}\leqslant CC_{\ast}\varepsilon\left(N_{0}^{s_{0}+2\overline{\sigma}+\frac{3}{2}\mu_{1}+1-b_{1}}+\varepsilon\gamma^{-2}+\varepsilon^{2}\gamma^{-2}\right).$$
				From \eqref{cond Nn}, one already has
				$$CN_{0}^{s_{0}+2\overline{\sigma}+\frac{3}{2}\mu_{1}+1-b_{1}}\leqslant\tfrac{1}{3}N_{0}^{-a_{1}}.$$
				Therefore, we need at this level to take $\varepsilon$ small enough to ensure
				$$C\left(\varepsilon\gamma^{-2}+\varepsilon^{2}\gamma^{-2}\right)\leqslant \tfrac{2}{3}N_{0}^{-a_{1}}.$$
				This occurs since \eqref{choice of gamma and N0 in the Nash-Moser} and \eqref{param NM} imply
				$$0<a<\tfrac{1}{2+a_{1}}.$$
				Hence
				$$\|\mathcal{F}(\widetilde{U}_{1})\|_{q,s_{0}}^{\gamma,\mathrm{O}_{1}^{2\gamma}}\leqslant C_{\ast}\varepsilon N_{0}^{-a_{1}}.$$ 
				This completes the proof of the estimates in $(\mathcal{P}2)_{n+1}.$ \\
				
				\noindent $\blacktriangleright$ \textbf{Extension and verification of $(\mathcal{P}1)_{n+1}-(\mathcal{P}3)_{n+1}.$} We shall now construct an extention of $\widetilde{H}_{n+1}$ living in the whole set of parameters and enjoying suitable decay properties. This is done by using the $C^{\infty}$ cut-off function $\chi_{n+1}:\mathcal{O}\rightarrow[0,1]$ defined by
				$$\chi_{n+1}(\lambda,\omega)=\left\lbrace\begin{array}{ll}
					1 & \textnormal{in }\mathrm{O}_{n+1}^{\gamma}\\
					0 & \textnormal{in }\mathcal{O}\setminus \mathrm{O}_{n+1}^{2\gamma}
				\end{array}\right.$$
				and satisfying the additional growth conditions
				\begin{equation}\label{dext-cutoff}
					\forall \alpha\in\mathbb{N}^{d},\quad |\alpha|\in\llbracket 0,q\rrbracket,\quad  \|\partial_{\lambda,\omega}^{\alpha}\chi_{n+1}\|_{L^{\infty}(\mathcal{O})}\lesssim\left(\gamma^{-1}N_{n}^{\overline{a}}\right)^{|\alpha|}.
				\end{equation}
				Next, we shall deal with the  extension $H_{n+1}$ of $\widetilde{H}_{n+1}$  defined by
				\begin{equation}\label{ext Hn+1}
					H_{n+1}(\lambda,\omega):=\left\lbrace\begin{array}{ll}
						\chi_{n+1}(\lambda,\omega)\widetilde{H}_{n+1}(\lambda,\omega) &\quad  \textnormal{in }\quad \mathrm{O}_{n+1}^{2\gamma}\\
						0 & \quad \textnormal{in }\quad\mathcal{O}\setminus\mathrm{O}_{n+1}^{2\gamma}
					\end{array}\right.
				\end{equation}
				and the extension $U_{n+1}$ of $\widetilde{U}_{n+1}$ by
				\begin{equation}\label{scheme}
					U_{n+1}:=U_{n}+H_{n+1}.
				\end{equation}
				We remark that
				$$H_{n+1}=\widetilde{H}_{n+1}\quad\textnormal{and}\quad\mathcal{F}(U_{n+1})=\mathcal{F}(\widetilde{U}_{n+1})\quad\textnormal{in}\quad \mathrm{O}_{n+1}^{\gamma}.
				$$
				Looking at the first component of \eqref{scheme}, one can write with obvious notations
				$$i_{n+1}=i_n+\mathfrak{I}_{n+1}.$$
				By the  induction assumption $(\mathcal{P}2)_n$,  \eqref{ext Hn+1} and \eqref{sym frih}, one has 
				$$\mathfrak{S}i_n(\varphi)=i_n(-\varphi)\quad\textnormal{and}\quad \mathfrak{S}\mathfrak{I}_{n+1}(\varphi)=\mathfrak{I}_{n+1}(-\varphi).$$
				Thus
				\begin{equation}\label{rev in+1}
					\mathfrak{S}i_{n+1}(\varphi)=i_{n+1}(-\varphi).
				\end{equation}
				Using Lemma \ref{Lem-lawprod}-(iv) together with \eqref{dext-cutoff} and the fact that  $H_{n+1}=0$ in $\mathcal{O}\setminus\mathrm{O}_{n+1}^{2\gamma}$, we obtain
				\begin{equation}\label{link Hm+1tilde Hm+1}
					\forall s\geqslant s_{0},\quad\|H_{n+1}\|_{q,s}^{\gamma,\mathcal{O}}\lesssim N_{n}^{q\overline{a}}\|\widetilde{H}_{n+1}\|_{q,s}^{\gamma,\mathrm{O}_{n+1}^{2\gamma}}.
				\end{equation}
				Applying \eqref{link Hm+1tilde Hm+1} and \eqref{link Hm+1 and mathcalF(Um) in norm s0} we deduce that for $n\in\mathbb{N}^{*},$
				\begin{align*}
					\| H_{n+1}\|_{q,s_{0}+\overline{\sigma}}^{\gamma,\mathcal{O}}&\leqslant C N_{n}^{q\overline{a}}\|\widetilde{H}_{n+1}\|_{q,s_{0}+\overline{\sigma}}^{\gamma,\mathrm{O}_{n+1}^{2\gamma}}\\
					&\leqslant C N_{n}^{q\overline{a}+\overline\sigma}\|\widetilde{H}_{n+1}\|_{q,s_{0}}^{\gamma,\mathrm{O}_{n+1}^{2\gamma}}\\
					&\leqslant C C_{\ast}\varepsilon\gamma^{-1} N_{n}^{q\overline{a}+2\overline{\sigma}-\frac{2}{3}a_{1}}.
				\end{align*}
				From \eqref{param NM}, we have 
				\begin{align}\label{ConDDR4}
					a_{2}=\tfrac{2}{3}a_{1}-q\overline{a}-2\overline{\sigma}-1\geqslant 1.
				\end{align}
				Therefore, choosing $\varepsilon$ small enough, we obtain 
				\begin{align}\label{estim Hm+1 tilde}
					\|H_{n+1}\|_{q,s_{0}+\overline{\sigma}}^{\gamma,\mathcal{O}}&\leqslant CN_0^{-1}C_{\ast}\varepsilon\gamma^{-1}N_{n}^{-a_{2}}\nonumber\\
					&\leqslant C_{\ast}\varepsilon\gamma^{-1}N_{n}^{-a_{2}}.
				\end{align}
				As to the case $n=0$, we combine \eqref{link Hm+1tilde Hm+1} and \eqref{H1} to obtain, up to take $C_{\ast}$ large enough,
				\begin{equation}\label{estim H1 tilde}
					\|H_{1}\|_{q,s}^{\gamma,\mathcal{O}}\leqslant \tfrac{1}{2}C_{\ast}\varepsilon\gamma^{-1}N_{0}^{q\overline{a}}.
				\end{equation}
				We now set
				\begin{equation}\label{Constru-1}
					W_{n+1}:=W_{n}+H_{n+1},
				\end{equation}
				then by construction, we infer
				$$U_{n+1}=U_0+W_{n+1}.$$
				Moreover, applying $(\mathcal{P}1)_{n}$, \eqref{estim H1 tilde} and \eqref{estim Hm+1 tilde} and Lemma \ref{lemma sum Nn}, we infer
				\begin{align*}
					\|W_{n+1}\|_{q,s_{0}+\overline{\sigma}}^{\gamma,\mathcal{O}}&\leqslant\|H_{1}\|_{q,s_{0}+\overline{\sigma}}^{\gamma,\mathcal{O}}+\sum_{k=2}^{n+1}\|H_{k}\|_{q,s_{0}+\overline{\sigma}}^{\gamma,\mathcal{O}}\\
					&\leqslant \tfrac{1}{2}C_{\ast}\varepsilon\gamma^{-1}N_{0}^{q\overline{a}}+C_{\ast}\varepsilon\gamma^{-1}\sum_{k=0}^{\infty}N_{k}^{-1}\\
					&\leqslant \tfrac{1}{2}C_{\ast}\varepsilon\gamma^{-1}N_{0}^{q\overline{a}}+CN_{0}^{-1}C_{\ast}\varepsilon\gamma^{-1}\\
					&\leqslant C_{\ast}\varepsilon\gamma^{-1}N_{0}^{q\overline{a}}.
				\end{align*}
				This completes the proof of $(\mathcal{P}1)_{n+1}.$ Now gathering \eqref{link Hm+1 and Wm}, \eqref{link Hm+1tilde Hm+1} and $(\mathcal{P}3)_{n}$ allows to write
				\begin{align*}
					\|W_{n+1}\|_{q,b_{1}+\overline{\sigma}}^{\gamma,\mathcal{O}} & \leqslant  \| W_{n}\|_{q,b_{1}+\overline{\sigma}}^{\gamma,\mathcal{O}}+CN_{n}^{q\overline{a}}\|H_{n+1}\|_{q,b_{1}+\overline{\sigma}}^{\gamma,\mathcal{O}}\\
					& \leqslant C_{\ast} \varepsilon\gamma^{-1}N_{n-1}^{\mu_{1}}+CC_{\ast}\gamma^{-1}N_{n}^{q\overline{a}+2\overline{\sigma}}\left(\varepsilon+\| W_{n}\|_{q,b_{1}+\overline{\sigma}}^{\gamma,\mathcal{O}_{n}^{\gamma}}\right)\\
					& \leqslant  CC_{\ast}\varepsilon\gamma^{-1}N_{n}^{q\overline{a}+2\overline{\sigma}+1+\frac{2}{3}\mu_{1}}.
				\end{align*}
				Fom \eqref{param NM}, we can ensure the condition
				\begin{align}\label{ConDDR1}
					q\overline{a}+2\overline{\sigma}+2=\tfrac{\mu_{1}}{3},
				\end{align}
				in order to get
				\begin{align*}
					\|W_{n+1}\|_{q,b_{1}+\overline{\sigma}}^{\gamma,\mathcal{O}_{n+1}^\gamma}&\leqslant CN_{0}^{-1}C_{\ast}\varepsilon\gamma^{-1}N_{n}^{\mu_{1}}\\
					&\leqslant C_{\ast}\varepsilon\gamma^{-1}N_{n}^{\mu_{1}}
				\end{align*}
				by taking $\varepsilon$ small enough and using \eqref{choice of gamma and N0 in the Nash-Moser}. This proves $(\mathcal{P}3)_{n+1}$ and the proof of Proposition \ref{Nash-Moser} is now complete.
			\end{proof}
			Once this sequence of approximate solutions is constructed, we may obtain a non-trivial solution by passing to the limit. This is possible due the decay properties given in Proposition \ref{Nash-Moser}. Actually, we obtain the following corollary.
			\begin{cor}\label{Corollary NM}
				There exists $\varepsilon_0>0$ such that for all $\varepsilon\in(0,\varepsilon_0),$ the following assertions hold true.
				We consider the Cantor set $\mathcal{G}_{\infty}^{\gamma}$, depending on $\varepsilon$ through $\gamma,$ and defined by
				$$\mathcal{G}_{\infty}^{\gamma}:=\bigcap_{n\in\mathbb{N}}\mathcal{A}_{n}^{\gamma}.$$
				There exists a function
				$$U_{\infty}:\begin{array}[t]{rcl}
					\mathcal{O} & \rightarrow & \left(\mathbb{T}^{d}\times\mathbb{R}^{d}\times H_{\mathbb{S}}^{\perp}\right)\times\mathbb{R}^{d}\times\mathbb{R}^{d+1}\\
					(\lambda,\omega) & \mapsto & \big(i_{\infty}(\lambda,\omega),\alpha_{\infty}(\lambda,\omega),(\lambda,\omega)\big)
				\end{array}$$
				such that 
				$$\forall(\lambda,\omega)\in\mathcal{G}_{\infty}^{\gamma},\quad\mathcal{F}(U_{\infty}(\lambda,\omega))=0.$$
				In addition, $i_{\infty}$ is reversible and $\alpha_{\infty}\in W^{q,\infty,\gamma}(\mathcal{O},\mathbb{R}^d)$ with
				\begin{equation}\label{alpha infty}
					\alpha_{\infty}(\lambda,\omega)=\omega+\mathrm{r}_{\varepsilon}(\lambda,\omega)\quad\mbox{ and }\quad\|\mathrm{r}_{\varepsilon}\|_{q}^{\gamma,\mathcal{O}}\lesssim\varepsilon\gamma^{-1}N_{0}^{q\overline{a}}.
				\end{equation}
				Moreover, there exists a $q$-times differentiable function $\lambda\in(\lambda_{0},\lambda_{1})\mapsto\omega(\lambda,\varepsilon)\in\mathbb{R}^d$ with
				\begin{equation}\label{alpha infty-1}
					\omega(\lambda,\varepsilon)=-\omega_{\textnormal{Eq}}(\lambda)+\bar{r}_{\varepsilon}(\lambda),\quad \|\bar{r}_{\varepsilon}\|_{q}^{\gamma,\mathcal{O}}\lesssim\varepsilon\gamma^{-1}N_{0}^{q\overline{a}}
				\end{equation}
				and 
				$$\forall\lambda\in \mathcal{C}_{\infty}^{\varepsilon},\quad \mathcal{F}\Big(U_{\infty}\big(\lambda,\omega(\lambda,\varepsilon)\big)\Big)=0\quad\textnormal{and}\quad\alpha_{\infty}\big(\lambda,\omega(\lambda,\varepsilon)\big)=-\omega_{\textnormal{Eq}}(\lambda),
				$$
				where the  Cantor set $\mathcal{C}_{\infty}^{\varepsilon}$ is defined by 
				\begin{equation}\label{def final Cantor set}
					\mathcal{C}_{\infty}^{\varepsilon}=\Big\{\lambda\in(\lambda_{0},\lambda_{1})\quad\textnormal{s.t.}\quad\big(\lambda,\omega(\lambda,\varepsilon)\big)\in\mathcal{G}_{\infty}^{\gamma}\Big\}.
				\end{equation}

			\end{cor}
			\begin{proof}
				Putting together \eqref{Constru-1} and \eqref{estim Hm+1 tilde}, we infer
				$$\|W_{n+1}-W_{n}\|_{q,s_{0}}^{\gamma,\mathcal{O}}=\|H_{n+1}\|_{q,s_{0}}^{\gamma,\mathcal{O}}\leqslant\|H_{n+1}\|_{q,s_{0}+\overline{\sigma}}^{\gamma,\mathcal{O}}\leqslant C_{\ast}\varepsilon\gamma^{-1}N_{n}^{-a_{2}}.$$
				Thus, the telescopic series associated with the sequence $\left(W_{n}\right)_{n\in\mathbb{N}}$ is convergent, so the sequence itself converges. We denote its limit
				$$W_{\infty}:=\lim_{n\rightarrow\infty}W_{n}:=(\mathfrak{I}_{\infty},\alpha_{\infty}-\omega,0,0)$$
				and
				$$U_{\infty}:=\big(i_{\infty},\alpha_{\infty},(\lambda,\omega)\big)=U_0+W_{\infty}.$$
				Passing to the limit in \eqref{rev in+1}, one obtains the reversibility property 
				$$\mathfrak{S}i_\infty(\varphi)=i_\infty(-\varphi).$$
				By the point $(\mathcal{P}2)_{n}$ of Proposition \ref{Nash-Moser}, we have for small $\varepsilon$ 
				\begin{equation}\label{sol before rigidity}
					\forall(\lambda,\omega)\in\mathcal{G}_{\infty}^{\gamma},\quad\mathcal{F}\Big(i_{\infty}(\lambda,\omega),\alpha_{\infty}(\lambda,\omega),(\lambda,\omega),\varepsilon\Big)=0,
				\end{equation}
				with $\mathcal{F}$ the functional  defined in \eqref{main function}. We highlight that the Cantor set $\mathcal{G}_{\infty}^{\gamma}$ depends on $\varepsilon$ through $\gamma$ and \eqref{choice of gamma and N0 in the Nash-Moser}. By the point $(\mathcal{P}1)_{n}$ of the Proposition \ref{Nash-Moser}, we have 
				$$\alpha_{\infty}(\lambda,\omega)=\omega+\mathrm{r}_{\varepsilon}(\lambda,\omega) \quad\mbox{ with }\quad\|\mathrm{r}_{\varepsilon}\|_{q}^{\gamma,\mathcal{O}}\lesssim\varepsilon\gamma^{-1}N_{0}^{q\overline{a}}.$$
				We now prove the second result and check the existence of solutions to the original Hamiltonian equation. First recall that the open set $\mathcal{O}$ is defined in \eqref{def initial parameters set} by
				$$\mathcal{O}=(\lambda_{0},\lambda_{1})\times\mathscr{U}\quad\textnormal{with}\quad\mathscr{U}=B(0,R_{0})\quad\textnormal{for some large }R_{0}>0,$$
				where the ball $\mathscr{U}$ is taken to contain the equilibrium frequency vector $\lambda\mapsto{\omega}_{\textnormal{Eq}}(\lambda).$ According to \eqref{alpha infty}, we deduce that for any $\lambda\in(\lambda_{0},\lambda_{1}),$ the mapping $\omega\mapsto\alpha_{\infty}(\lambda,\omega)$ is invertible from $\mathscr{U}$ into its image $\alpha_{\infty}(\lambda,\mathscr{U})$ and we have
				$$\widehat{\omega}=\alpha_{\infty}(\lambda,\omega)=\omega+\mathrm{r}_{\varepsilon}(\lambda,\omega)\Leftrightarrow\omega=\alpha_{\infty}^{-1}(\lambda,\widehat{\omega})=\widehat{\omega}+\widehat{\mathrm{r}}_{\varepsilon}(\lambda,\widehat{\omega}).$$
				This gives the identity
				$$\widehat{\mathrm{r}}_{\varepsilon}(\lambda,\widehat{\omega})=-\mathrm{r}_{\varepsilon}(\lambda,\omega),$$
				which implies in turn after using successive differentiation and \eqref{alpha infty} that $\widehat{\mathrm{r}}_{\varepsilon}$ satisfies the estimate
				\begin{equation}\label{estimate mathrm repsilon}
					\|\widehat{\mathrm{r}}_{\varepsilon}\|_{q}^{\gamma,\mathcal{O}}\lesssim\varepsilon\gamma^{-1}N_{0}^{q\overline{a}}.
				\end{equation}
				We now set 
				$${\omega}(\lambda,\varepsilon):=\alpha_{\infty}^{-1}(\lambda,-{\omega}_{\textnormal{Eq}}(\lambda))=-{\omega}_{\textnormal{Eq}}(\lambda)+\overline{\mathrm{r}}_{\varepsilon}(\lambda)\quad \mbox{ with }\quad \overline{\mathrm{r}}_{\varepsilon}(\lambda):=\widehat{\mathrm{r}}_{\varepsilon}\big(\lambda,-{\omega}_{\textnormal{Eq}}(\lambda)\big).
				$$
				As a consequence of \eqref{sol before rigidity}, if we denote
				$$\mathcal{C}_{\infty}^{\varepsilon}:=\Big\{\lambda\in(\lambda_{0},\lambda_{1})\quad\textnormal{s.t.}\quad\big(\lambda,\omega(\lambda,\varepsilon)\big)\in\mathcal{G}_{\infty}^{\gamma}\Big\},$$
				then we have
				$$\forall\lambda\in\mathcal{C}_{\infty}^{\varepsilon},\quad\mathcal{F}\Big(U_{\infty}\big(\lambda,\omega(\lambda,\varepsilon)\big)\Big)=0.$$
				This gives a nontrivial reversible solution for the original Hamiltonian equation provided that $\lambda\in\mathcal{C}_{\infty}^{\varepsilon}.$
				Since all the derivatives up to order $q$ of $\omega_{\textnormal{Eq}}$ are uniformly bounded on $[\lambda_{0},\lambda_{1}],$ see Lemma \ref{lemma properties linear frequencies}-(vi), then by chain rule and \eqref{estimate mathrm repsilon}, we obtain 
				\begin{equation}\label{estimate repsilon1}
					\|\overline{\mathrm{r}}_{\varepsilon}\|_{q}^{\gamma,\mathcal{O}}\lesssim\varepsilon\gamma^{-1}N_{0}^{q\overline{a}}\quad\textnormal{and}\quad\|\omega(\cdot,\varepsilon)\|_{q}^{\gamma,\mathcal{O}}\lesssim 1+\varepsilon\gamma^{-1}N_{0}^{q\overline{a}}\lesssim 1.
				\end{equation}
				This ends the proof of Corollary \ref{Corollary NM}.
			\end{proof}
			\subsection{Measure of the final Cantor set}\label{Section measure of the final Cantor set}
			The purpose of this final section is to give a lower bound of  the Lebesgue measure of the Cantor set $\mathcal{C}_{\infty}^{\varepsilon}$ constructed in Corollary \ref{Corollary NM} via \eqref{def final Cantor set}. We  show that this set is massive and asymptotically when  $\varepsilon\rightarrow 0$ it   tends to be of  full
			measure in $(\lambda_0,\lambda_1)$. Note that Corollary
			\ref{Corollary NM} allows us to write the Cantor set $\mathcal{C}_{\infty}^{\varepsilon}$  in the following form
			\begin{equation}\label{definition of the final Cantor set}
				\mathcal{C}_{\infty}^{\varepsilon}:=\bigcap_{n\in\mathbb{N}}\mathcal{C}_{n}^{\varepsilon}\quad \mbox{ where }\quad \mathcal{C}_{n}^{\varepsilon}:=\Big\{\lambda\in(\lambda_{0},\lambda_{1})\quad\hbox{s.t}\quad \big(\lambda,{\omega}(\lambda,\varepsilon)\big)\in\mathcal{A}_n^{\gamma}\Big\}.
			\end{equation}
			The sets $\mathcal{A}_n^{\gamma}$ and the perturbed frequency vector $\omega(\lambda,\varepsilon)$ are respectively defined in Proposition \ref{Nash-Moser} and in \eqref{alpha infty}. The main result of this section reads as follows.
			\begin{prop}\label{lem-meas-es1}
				Let $q_{0}$ be defined as in Lemma $\ref{lemma transversality}$ and assume that \eqref{param NM} and \eqref{choice of gamma and N0 in the Nash-Moser} hold with $q=q_0+1.$ Assume the additional conditions
				\begin{equation}\label{choice tau 1 tau2 upsilon} 
					\left\lbrace\begin{array}{l}
						\tau_{1}> dq_{0}\\
						\tau_{2}>\tau_{1}+dq_{0}\\
						\upsilon=\frac{1}{q_{0}+3}.
					\end{array}\right.
				\end{equation}
				Then, there exists $C>0$ such that 
				$$
				\big|\mathcal{C}_{\infty}^{\varepsilon}\big|\geqslant (\lambda_1-\lambda_0)-C \varepsilon^{\frac{a\upsilon}{q_{0}}}.
				$$
				In particular, 
				$$\displaystyle \lim_{\varepsilon\to0}\big|\mathcal{C}_{\infty}^{\varepsilon}\big|=\lambda_1-\lambda_0.$$
			\end{prop}
			%			\begin{remark}
			%				The condition \eqref{choice tau 1 tau2 upsilon} is imposed to solve some problem of convergence of series and smallness conditions. We refer to \eqref{assum2-d} and \eqref{assum-co1} to understand these needed conditions. In particular, one has $\tau_{2}>\tau_{1}>d$ and $\upsilon\leqslant\frac{1}{q_{0}+3}$ whose conditions were required from \eqref{setting tau1 and tau2} and Proposition \ref{reduction of the transport part}. Also notice that \eqref{choice of gamma and N0 in the Nash-Moser} and \eqref{param NM} implies that $\upsilon$ is indeed strictly positive.
			%			\end{remark}
			The remainder of this section is devoted to the proof of Proposition \ref{lem-meas-es1}. We shall begin by giving the proof using some a priori results. These results will be  proved later  in Lemmas \ref{lemm-dix1}, \ref{some cantor set are empty} and \ref{lemma Russeman condition for the perturbed frequencies}. We first give a short insight about   the strategy to prove Proposition \ref{lem-meas-es1}. The idea is to measure the complementary set of $\mathcal{C}_{\infty}^{\varepsilon}$ in $(\lambda_{0},\lambda_{1}).$ To proceed with, we write
			\begin{equation}\label{decomp Cantor}
				(\lambda_{0},\lambda_{1})\setminus\mathcal{C}_{\infty}^{\varepsilon}=\big((\lambda_{0},\lambda_{1})\setminus\mathcal{C}_{0}^{\varepsilon}\big)\sqcup\bigsqcup_{n=0}^{\infty}\big(\mathcal{C}_{n}^{\varepsilon}\setminus\mathcal{C}_{n+1}^{\varepsilon}\big).
			\end{equation}
			The measure of  each set which appears, we estimate it by using Lemma \ref{lemma useful for measure estimates}. We shall now give the proof of Proposition \ref{lem-meas-es1}.
			%Next we shall establish  a precise statement of Proposition \ref{lem-meas-es1}.
			%\begin{prop}\label{Prop-Pre-StX1}
			%Let $q_{0}$ as in Lemma \ref{lemma Russeman condition for the perturbed frequencies} and  $\gamma$ and $a$ as in \eqref{choice of gamma and N0 in the Nash-Moser} and fix $\upsilon=\frac{1}{2q_0}.$ Assume that $\tau_1$ and $\tau_2$ satisfy the conditions
			%$$
			%\textcolor{blue}{\tau_1>d \,q_0\quad\hbox{and}\quad \tau_2>\tau_1+ d\, q_0,}
			%$$
			%then the Lebesgue measure  of the set $\mathcal{C}_{\infty}^{\varepsilon}$ satisfies 
			%$$\big|[\lambda_{0},\lambda_{1}]\setminus\mathcal{C}_{\infty}^{\varepsilon}\big|\lesssim\varepsilon^{\frac{a}{2q_{0}}}.$$
			%
			%\end{prop}
			\begin{proof}
				By choosing $R_{0}$ large enough, one can ensure using \eqref{alpha infty-1} that
				$$\forall\lambda\in(\lambda_0,\lambda_1),\quad \omega(\lambda,\varepsilon)\in \mathscr{U}=B(0,R_{0}).$$
				Indeed,  $\mathscr{U}$ contains by construction the curve $\lambda\in(\lambda_0,\lambda_1)\mapsto\pm\omega_{\textnormal{Eq}}(\lambda)$ and by \eqref{alpha infty-1} and \eqref{choice of gamma and N0 in the Nash-Moser}, one has
				$$\sup_{\lambda\in(\lambda_0,\lambda_1)}\left|\omega(\lambda,\varepsilon)+\omega_{\textnormal{Eq}}(\lambda)\right|\leqslant\|\overline{\mathrm{r}}_{\varepsilon}\|_{q}^{\gamma,\mathcal{O}}\leqslant C\varepsilon\gamma^{-1}N_{0}^{q\overline{a}}=C\varepsilon^{1-a(1+q\overline{a})}.$$
				Now, the conditions \eqref{param NM} and \eqref{choice of gamma and N0 in the Nash-Moser} imply in particular
				$$0<a<\frac{1}{1+q\overline{a}}\cdot$$
				Hence, by taking $\varepsilon$ small enough, we find
				$$\sup_{\lambda\in(\lambda_0,\lambda_1)}\left|\omega(\lambda,\varepsilon)+\omega_{\textnormal{Eq}}(\lambda)\right|\leqslant\|\overline{\mathrm{r}}_{\varepsilon}\|_{q}^{\gamma,\mathcal{O}}\leqslant 1.$$
				As a consequence, 
				$$\mathcal{C}_{0}^{\varepsilon}=(\lambda_0,\lambda_1).$$
				By \eqref{decomp Cantor}, we can write 
				\begin{align}\label{Mir-1}
					\nonumber \Big|(\lambda_{0},\lambda_{1})\setminus\mathcal{C}_{\infty}^{\varepsilon}\Big|\leqslant&\sum_{n=0}^{\infty}\Big|\mathcal{C}_{n}^{\varepsilon}\setminus\mathcal{C}_{n+1}^{\varepsilon}\Big|\\
					&:=\sum_{n=0}^{\infty}\mathcal{S}_{n}.
				\end{align}
				According to the notation introduced in  Proposition \ref{reduction of the remainder term} and Proposition \ref{projection in the normal directions} one may write 
				\begin{align}\label{asy-z1}
					\nonumber \mu_{j}^{\infty,n}(\lambda,\varepsilon)&:=\mu_{j}^{\infty}\big(\lambda,{\omega}(\lambda,\varepsilon),i_{n}\big)\\
					&=\Omega_{j}(\lambda)+jr^{1,n}(\lambda,\varepsilon)+r_{j}^{\infty,n}(\lambda,\varepsilon),
				\end{align}
				with
				\begin{align*}
					r^{1,n}(\lambda,\varepsilon)&:=c_{n}(\lambda,\varepsilon)-\Omega-I_{1}(\lambda)K_{1}(\lambda),\\
					c_{n}(\lambda,\varepsilon)&:=c_{i_{n}}(\lambda,{\omega}(\lambda,\varepsilon)),\\
					r_{j}^{\infty,n}(\lambda,\varepsilon)&:=r_{j}^{\infty}\big(\lambda,{\omega}(\lambda,\varepsilon),i_{n}\big).
				\end{align*}
				Coming back to \eqref{definition of the final Cantor set} and using the Cantor sets introduced in Proposition \ref{reduction of the remainder term}, Proposition \ref{inversion of the linearized operator in the normal directions} and \mbox{Proposition \ref{reduction of the transport part}} one obtains by construction that   for any $n\in\mathbb{N}$, 
				\begin{equation}\label{set-U0}
					\mathcal{C}_{n}^{\varepsilon}\setminus\mathcal{C}_{n+1}^{\varepsilon}=\bigcup_{(l,j)\in\mathbb{Z}^{d}\times\mathbb{Z}\setminus\{(0,0)\}\atop |l|\leqslant N_{n}}\mathcal{R}_{l,j}^{(0)}(i_{n})\bigcup_{(l,j,j_{0})\in\mathbb{Z}^{d}\times(\mathbb{S}_{0}^{c})^{2}\atop |l|\leqslant N_{n}}\mathcal{R}_{l,j,j_{0}}(i_{n})\bigcup_{(l,j)\in\mathbb{Z}^{d}\times\mathbb{S}_{0}^{c}\atop |l|\leqslant N_{n}}\mathcal{R}_{l,j}^{(1)}(i_{n}),
				\end{equation}
				with
				\begin{align*}
					\mathcal{R}_{l,j}^{(0)}(i_{n})&:=\left\lbrace\lambda\in\mathcal{C}_{n}^{\varepsilon}\quad\textnormal{s.t.}\quad|{\omega}(\lambda,\varepsilon)\cdot l+jc_{n}(\lambda,\varepsilon)|\leqslant\tfrac{4\gamma_{n+1}^{\upsilon}\langle j\rangle}{\langle l\rangle^{\tau_{1}}}\right\rbrace,\\
					\mathcal{R}_{l,j,j_{0}}(i_{n})&:=\left\lbrace\lambda\in\mathcal{C}_{n}^{\varepsilon}\quad\textnormal{s.t.}\quad|{\omega}(\lambda,\varepsilon)\cdot l+\mu_{j}^{\infty,n}(\lambda,\varepsilon)-\mu_{j_{0}}^{\infty,n}(\lambda,\varepsilon)|\leqslant\tfrac{2\gamma_{n+1}\langle j-j_{0}\rangle}{\langle l\rangle^{\tau_{2}}}\right\rbrace,\\
					\mathcal{R}_{l,j}^{(1)}(i_{n})&:=\left\lbrace\lambda\in\mathcal{C}_{n}^{\varepsilon}\quad\textnormal{s.t.}\quad|{\omega}(\lambda,\varepsilon)\cdot l+\mu_{j}^{\infty,n}(\lambda,\varepsilon)|\leqslant\tfrac{\gamma_{n+1}\langle j\rangle}{\langle l\rangle^{\tau_{1}}}\right\rbrace.
				\end{align*}
				Notice that using the inclusion
				$$W^{q,\infty,\gamma}(\mathcal{O},\mathbb{C})\hookrightarrow C^{q-1}(\mathcal{O},\mathbb{C})$$
				and the fact that $q=q_0+1$, one gets that for all $n\in\mathbb{N}$ and $(l,j,j_0)\in\mathbb{Z}^{d}\times(\mathbb{S}_0^c)^{2},$ the curves
				\begin{align*}
					&\lambda\mapsto\omega(\lambda,\varepsilon)\cdot l+c_{n}(\lambda,\varepsilon),\\
					&\lambda\mapsto\omega(\lambda,\varepsilon)\cdot l+\mu_{j}^{\infty,n}(\lambda,\varepsilon)-\mu_{j_0}^{\infty,n}(\lambda,\varepsilon),\\
					&\lambda\mapsto\omega(\lambda,\varepsilon)\cdot l+\mu_{j}^{\infty,n}(\lambda,\varepsilon)	
				\end{align*}
				have a $C^{q_0}$ regularity. Then, applying Lemma \ref{lemma useful for measure estimates} combined with  Lemma \ref{lemma Russeman condition for the perturbed frequencies} gives for any  $n\in\mathbb{N}$, 
				\begin{align}\label{kio1}
					\Big|\mathcal{R}_{l,j}^{(0)}(i_{n})\Big|&\lesssim\gamma^{\frac{\upsilon}{q_{0}}}\langle j\rangle^{\frac{1}{q_{0}}}\langle l\rangle^{-1-\frac{\tau_{1}+1}{q_{0}}},\nonumber\\
					\Big|\mathcal{R}_{l,j}^{(1)}(i_{n})\Big|&\lesssim\gamma^{\frac{1}{q_{0}}}\langle j\rangle^{\frac{1}{q_{0}}}\langle l\rangle^{-1-\frac{\tau_{1}+1}{q_{0}}},\\
					\Big|\mathcal{R}_{l,j,j_{0}}(i_{n})\Big|&\lesssim\gamma^{\frac{1}{q_{0}}}\langle j-j_{0}\rangle^{\frac{1}{q_{0}}}\langle l\rangle^{-1-\frac{\tau_{2}+1}{q_{0}}}.\nonumber
				\end{align}
				Let us now move to the estimate of $\mathcal{S}_0$ and $\mathcal{S}_{1}$  defined in \eqref{Mir-1} that should  be treated differently from  the other terms. This is related to the discussion done at   the beginning of the proof of Lemma \ref{lemm-dix1} dealing  with the validity of the estimate \eqref{Bio-X1}. By using Lemma  \ref{some cantor set are empty}, we find for all $k\in\{0,1\}$,
				\begin{align}\label{set-U5}
					\mathcal{S}_{k}\lesssim  \sum_{\underset{|j|\leqslant C_{0}\langle l\rangle, |l|\leqslant N_{k}}{(l,j)\in\mathbb{Z}^{d}\times\mathbb{Z}\setminus\{(0,0)\}}}\Big|\mathcal{R}_{l,j}^{(0)}(i_{k})\Big|+\sum_{\underset{\underset{\min(|j|,|j_{0}|)\leqslant c_{2}\gamma_{1}^{-\upsilon}\langle l\rangle^{\tau_{1}}}{|j-j_{0}|\leqslant C_{0}\langle l\rangle, |l|\leqslant N_{k}}}{(l,j,j_{0})\in\mathbb{Z}^{d}\times(\mathbb{S}_{0}^{c})^{2}}}\Big|\mathcal{R}_{l,j,j_{0}}(i_{k})\Big|+\sum_{\underset{|j|\leqslant C_{0}\langle l\rangle, |l|\leqslant N_{k}}{(l,j)\in\mathbb{Z}^{d}\times\mathbb{S}_{0}^{c}}}\Big|\mathcal{R}_{l,j}^{(1)}(i_{k})\Big|.
				\end{align}
				Plugging \eqref{kio1} into \eqref{set-U5} yields for all $k\in\{0,1\}$,
				\begin{align*}
					\mathcal{S}_{k}&\lesssim  \gamma^{\frac{1}{q_{0}}}\Bigg(\sum_{|j|\leqslant C_{0}\langle l\rangle}|j|^{\frac{1}{q_{0}}}\langle l\rangle^{-1-\frac{\tau_{1}+1}{q_{0}}}+\sum_{\underset{\min(|j|,|j_{0}|)\leqslant c_{2}\gamma^{-\upsilon}\langle l\rangle^{\tau_{1}}}{|j-j_{0}|\leqslant C_{0}\langle l\rangle}} |j-j_{0}|^{\frac{1}{q_{0}}}\langle l\rangle^{-1-\frac{\tau_{2}+1}{q_{0}}}\Bigg)\\
					&\quad+\gamma^{\frac{\upsilon}{q_0}}\sum_{{|j|\leqslant C_{0}\langle l\rangle}}|j|^{\frac{1}{q_{0}}}\langle l\rangle^{-1-\frac{\tau_{1}+1}{q_{0}}}.
				\end{align*}
				Consequently, we obtain 
				\begin{align}\label{set-U6}
					\max_{k\in\{0,1\}}\mathcal{S}_{k}&\lesssim  \gamma^{\frac{1}{q_{0}}}\Bigg(\sum_{l\in\mathbb{Z}^d}\langle l\rangle^{-\frac{\tau_{1}}{q_{0}}}+\gamma^{-\upsilon}\sum_{l\in\mathbb{Z}^d}\langle l\rangle^{\tau_1-1-\frac{\tau_{2}}{q_{0}}}\Bigg)
					+\gamma^{\frac{\upsilon}{q_0}}\sum_{l\in\mathbb{Z}^d}\langle l\rangle^{-\frac{\tau_{1}}{q_{0}}}\\
					\nonumber &\lesssim  \gamma^{\min\left(\frac{\upsilon}{q_{0}},\frac{1}{q_{0}}-\upsilon\right)}.
				\end{align}
				Notice that the last estimate is obtained provided that we choose the parameters $\tau_1$ and $\tau_2$ in the following way in order to make the series convergent
				\begin{align}\label{assum2-d}
					\tau_1>d \,q_0\quad\hbox{and}\quad \tau_2>\tau_1+ d\, q_0.
				\end{align}
				This condition is exactly what we required in \eqref{choice tau 1 tau2 upsilon}. Concerning the estimate of $\mathcal{S}_{n}$ for $n\geqslant 2$ in \eqref{Mir-1} we  may use  Lemma \ref{lemm-dix1}  and Lemma \ref{some cantor set are empty}, in order to get
				\begin{align*}
					\mathcal{S}_{n}\leqslant \sum_{\underset{|j|\leqslant C_{0}\langle l\rangle,N_{n-1}<|l|\leqslant N_{n}}{(l,j)\in\mathbb{Z}^{d}\times\mathbb{Z}\setminus\{(0,0)\}}}\Big|\mathcal{R}_{l,j}^{(0)}(i_{n})\Big|+\sum_{\underset{\underset{\min(|j|,|j_{0}|)\leqslant c_{2}\gamma_{n+1}^{-\upsilon}\langle l\rangle^{\tau_{1}}}{|j-j_{0}|\leqslant C_{0}\langle l\rangle,N_{n-1}<|l|\leqslant N_{n}}}{(l,j,j_{0})\in\mathbb{Z}^{d}\times(\mathbb{S}_{0}^{c})^{2}}}\Big|\mathcal{R}_{l,j,j_{0}}(i_{n})\Big|+\sum_{\underset{|j|\leqslant C_{0}\langle l\rangle,N_{n-1}<|l|\leqslant N_{n}}{(l,j)\in\mathbb{Z}^{d}\times\mathbb{S}_{0}^{c}}}\Big|\mathcal{R}_{l,j}^{(1)}(i_{n})\Big|.
				\end{align*}
				Remark that if $|j-j_{0}|\leqslant C_{0}\langle l\rangle$ and $\min(|j|,|j_{0}|)\leqslant\gamma_{n+1}^{-\upsilon}\langle l\rangle^{\tau_{1}}$, then 
				$$\max(|j|,|j_{0}|)=\min(|j|,|j_{0}|)+|j-j_{0}|\leqslant\gamma_{n+1}^{-\upsilon}\langle l\rangle^{\tau_{1}}+C_{0}\langle l\rangle\lesssim\gamma^{-\upsilon}\langle l\rangle^{\tau_{1}}.$$
				Therefore, \eqref{kio1} implies
				\begin{align*}
					\mathcal{S}_{n}\lesssim 
					\gamma^{\frac{1}{q_{0}}}\Bigg(\sum_{|l|>N_{n-1}}\langle l\rangle^{-\frac{\tau_{1}}{q_{0}}}+\gamma^{-\upsilon}\sum_{|l|>N_{n-1}}\langle l\rangle^{\tau_1-1-\frac{\tau_{2}}{q_{0}}}\Bigg)
					+\gamma^{\frac{\upsilon}{q_0}}\sum_{|l|>N_{n-1}}\langle l\rangle^{-\frac{\tau_{1}}{q_{0}}}.
				\end{align*}
				Under the assumption, we obtain \eqref{assum2-d}
				\begin{align}\label{dtu-c1}
					\sum_{n=2}^\infty \mathcal{S}_{n}&\lesssim  \gamma^{\min\left(\frac{\upsilon}{q_{0}},\frac{1}{q_{0}}-\upsilon\right)}.
				\end{align}
				Plugging \eqref{dtu-c1} and \eqref{set-U6} into \eqref{Mir-1} gives
				$$\Big|(\lambda_{0},\lambda_{1})\setminus\mathcal{C}_{\infty}^{\varepsilon}\Big|  \lesssim  \gamma^{\min\left(\frac{\upsilon}{q_{0}},\frac{1}{q_{0}}-\upsilon\right)}$$
				provided that the condition \eqref{assum2-d} is satisfied. 
				The condition \eqref{choice tau 1 tau2 upsilon} implies that
				$$\min\left(\tfrac{\upsilon}{q_{0}},\tfrac{1}{q_{0}}-\upsilon\right)=\tfrac{\upsilon}{q_{0}}.$$
				We then find,  since $\gamma=\varepsilon^a$ according to \eqref{choice of gamma and N0 in the Nash-Moser},
				$$\Big|(\lambda_{0},\lambda_{1})\setminus\mathcal{C}_{\infty}^{\varepsilon}\Big|\lesssim\varepsilon^{\frac{a\upsilon}{q_{0}}}.$$
				This completes the proof of Proposition \ref{lem-meas-es1}. 
			\end{proof}
			Now we are left to prove Lemma \ref{lemm-dix1} and Lemma \ref{some cantor set are empty}  used in the proof of Proposition \ref{lem-meas-es1}.
			\begin{lem}\label{lemm-dix1}
				{Let $n\in\mathbb{N}\setminus\{0,1\}$ and  $l\in\mathbb{Z}^{d}$ such that $|l|\leqslant N_{n-1}.$ Then the following assertions hold true.
					\begin{enumerate}[label=(\roman*)]
						\item For $ j\in\mathbb{Z}$ with $(l,j)\neq(0,0)$, we get  $\,\,\mathcal{R}_{l,j}^{(0)}(i_{n})=\varnothing.$
						\item For  $ (j,j_{0})\in(\mathbb{S}_{0}^{c})^{2}$ with $(l,j)\neq(0,j_0),$ we get $\,\,\mathcal{R}_{l,j,j_{0}}(i_{n})=\varnothing.$
						\item For  $j\in\mathbb{S}_{0}^{c}$, we get $\,\,\mathcal{R}_{l,j}^{(1)}(i_{n})=\varnothing.$
						\item For any $n\in\mathbb{N}\setminus\{0,1\},$
						\begin{equation*}
							\mathcal{C}_{n}^{\varepsilon}\setminus\mathcal{C}_{n+1}^{\varepsilon}=\bigcup_{\underset{N_{n-1}<|l|\leqslant N_{n}}{(l,j)\in\mathbb{Z}^{d}\times\mathbb{Z}\setminus\{(0,0)\}}}\mathcal{R}_{l,j}^{(0)}(i_{n})\cup\bigcup_{\underset{N_{n-1}<|l|\leqslant N_{n}}{(l,j,j_{0})\in\mathbb{Z}^{d}\times(\mathbb{S}_{0}^{c})^{2}}}\mathcal{R}_{l,j,j_{0}}(i_{n})\cup\bigcup_{\underset{N_{n-1}<|l|\leqslant N_{n}}{(l,j)\in\mathbb{Z}^{d}\times\mathbb{S}_{0}^{c}}}\mathcal{R}_{l,j}^{(1)}(i_{n}).
						\end{equation*}
				\end{enumerate}}
			\end{lem}
			\begin{proof}
				In all the proof, we shall use the following estimate coming from \eqref{e-Hn-diff}, namely, for all $n\geqslant 2$, 
				\begin{align}\label{Bio-X1}
					\nonumber \| i_{n}-i_{n-1}\|_{q,\overline{s}_{h}+\sigma_{4}}^{\gamma,\mathcal{O}}&\leqslant\| U_{n}-U_{n-1}\|_{q,\overline{s}_{h}+\sigma_{4}}^{\gamma,\mathcal{O}}\\
					&\leqslant\| H_{n}\|_{q,s_{h}+\sigma_{4}}^{\gamma,\mathcal{O}}\nonumber\\
					&\leqslant C_{\ast}\varepsilon\gamma^{-1}N_{n-1}^{-a_{2}}.
				\end{align}
				The fact that the previous estimate is valid only for $n\geqslant 2$ is the reason why we had to treat the cases of $\mathcal{S}_{0}$ and $\mathcal{S}_{1}$ sparately in the proof of Proposition \ref{lem-meas-es1}.\\
				
				\noindent\textbf{(i)} We begin by proving that if $|l|\leqslant N_{n-1}$ and $(l,j)\neq (0,0),$  then $\mathcal{R}_{l,j}^{(0)}(i_{n}) \subset\mathcal{R}_{l,j}^{(0)}(i_{n-1}).$ Assume for a while this inclusion and let us check how this implies that $\mathcal{R}_{l,j}^{(0)}(i_{n})=\varnothing.$ In view of \eqref{set-U0} one obtains  
				$$
				\mathcal{R}_{l,j}^{(0)}(i_{n}) \subset\mathcal{R}_{l,j}^{(0)}(i_{n-1})\subset \mathcal{C}_{n-1}^{\varepsilon}\setminus\mathcal{C}_{n}^{\varepsilon}.
				$$
				Now \eqref{set-U0} implies in particular $\mathcal{R}_{l,j}^{(0)}(i_{n}) \subset\mathcal{C}_{n}^{\varepsilon}\setminus\mathcal{C}_{n+1}^{\varepsilon}$ and thus we conclude
				$$
				\mathcal{R}_{l,j}^{(0)}(i_{n}) \subset\big(\mathcal{C}_{n}^{\varepsilon}\setminus\mathcal{C}_{n+1}^{\varepsilon}\big)\cap \big(\mathcal{C}_{n-1}^{\varepsilon}\setminus\mathcal{C}_{n}^{\varepsilon}\big)=\varnothing.
				$$
				We now turn to the proof of the inclusion. Let us consider $\lambda\in\mathcal{R}_{l,j}^{(0)}(i_{n}).$ By construction, we get in particular that  $\lambda\in \mathcal{C}_{n}^{\varepsilon}\subset \mathcal{C}_{n-1}^{\varepsilon}.$ Moreover, by the triangle inequality, we obtain 
				\begin{align*}
					\big|{\omega}(\lambda,\varepsilon)\cdot l+jc_{n-1}(\lambda,\varepsilon)\big| & \leqslant \big|{\omega}(\lambda,\varepsilon)\cdot l+jc_{n}(\lambda,\varepsilon)\big|+|j|\big|c_{n}(\lambda,\varepsilon)-c_{n-1}(\lambda,\varepsilon)\big|\\
					&\leqslant \displaystyle\tfrac{4\gamma_{n+1}^{\upsilon}\langle j\rangle}{\langle l\rangle^{\tau_{1}}}+C|j|\|c_{ i_{n}}-c_{i_{n-1}}\|_{q}^{\gamma,\mathcal{O}}.
				\end{align*}
				Therefore, combining \eqref{difference ci}, \eqref{Bio-X1}, \eqref{choice of gamma and N0 in the Nash-Moser} and the fact tht $\sigma_{4}\geqslant 2$, we infer
				\begin{align*}
					\big|{\omega}(\lambda,\varepsilon)\cdot l+jc_{n-1}(\lambda,\varepsilon)\big| & \leqslant  \displaystyle\tfrac{4\gamma_{n+1}^{\upsilon}\langle j\rangle}{\langle l\rangle^{\tau_{1}}}+C\varepsilon\langle j\rangle\| i_{n}-i_{n-1}\|_{q,\overline{s}_{h}+2}^{\gamma,\mathcal{O}}\\
					& \leqslant \displaystyle\tfrac{4\gamma_{n+1}^{\upsilon}\langle j\rangle}{\langle l\rangle^{\tau_{1}}}+C\varepsilon^{2-a}\langle j\rangle N_{n-1}^{-a_{2}}.
				\end{align*}
				In view of the definition of $\gamma_n$ in Proposition \ref{Nash-Moser}-$(\mathcal{P}2)_n$ one gets
				$$
				\exists c_0>0,\quad \forall n\in \mathbb{N},\quad \gamma_{n+1}^{\upsilon}-\gamma_{n}^{\upsilon}\leqslant - c_0\,\gamma^{\upsilon} 2^{-n}.
				$$
				Now remark that \eqref{choice tau 1 tau2 upsilon},  \eqref{param NM} and \eqref{choice of gamma and N0 in the Nash-Moser} imply 
				\begin{align}\label{assum-co1}
					2-a-a \upsilon>1\quad\hbox{and}\quad a_2>\tau_1,
				\end{align}
				and therefore one gets $\displaystyle \sup_{n\in\mathbb{N}}2^{n}N_{n-1}^{-a_2+\tau_1}<\infty.$ It follows that, for $\varepsilon$ small enough  and $|l|\leqslant N_{n-1}$,  
				\begin{align*}
					\big|{\omega}(\lambda,\varepsilon)\cdot l+jc_{n-1}(\lambda,\varepsilon)\big| & \leqslant \displaystyle\tfrac{4\gamma_{n}^{\upsilon}\langle j\rangle}{\langle l\rangle^{\tau_{1}}}+C\tfrac{\langle j\rangle\gamma^{\upsilon}}{2^n\langle l\rangle^{\tau_{1}}}\Big(-4c_0+C\varepsilon2^{n}N_{n-1}^{-a_2+\tau_1}\Big)\\
					& \leqslant \displaystyle\tfrac{4\gamma_{n}^{\upsilon}\langle j\rangle}{\langle l\rangle^{\tau_{1}}}.
				\end{align*}
				Consequently $\lambda\in \mathcal{R}_{l,j}^{(0)}(i_{n-1})$ and this achieves the proof.\\
				\textbf{(ii)} Let $(j,j_{0})\in(\mathbb{S}_{0}^{c})^{2}$ and  $(l,j)\neq(0,j_0).$ If $j=j_0$ then by construction $\mathcal{R}_{l,j_0,j_{0}}(i_{n})=\mathcal{R}_{l,0}^{(0)}(i_{n})$ and then the result follows from the  point (i). Now let us discuss the case when $j\neq j_0.$ Similarly to the point (i), in order to get the result it is enough to  check that $\mathcal{R}_{l,j,j_{0}}(i_{n})\subset \mathcal{R}_{l,j,j_{0}}(i_{n-1}).$
				Let  $\lambda\in\mathcal{R}_{l,j,j_{0}}(i_{n})$ then from the definition of this set introduced  in \eqref{set-U0} we deduce that  $\lambda\in\mathcal{C}_{n}^{\varepsilon}\subset \mathcal{C}_{n-1}^{\varepsilon}$ and 
				\begin{equation}\label{poiH1}
					\big|{\omega}(\lambda,\varepsilon)\cdot l+\mu_{j}^{\infty,n-1}(\lambda,\varepsilon)-\mu_{j_{0}}^{\infty,n-1}(\lambda,\varepsilon)\big|  \leqslant  \displaystyle\tfrac{2\gamma_{n+1}\langle j-j_{0}\rangle}{\langle l\rangle^{\tau_{2}}}+\varrho^n_{j,j_0}(\lambda,\varepsilon),
				\end{equation}
				where we set 
				$$ \varrho^n_{j,j_0}(\lambda,\varepsilon):=\big|\mu_{j}^{\infty,n}(\lambda,\varepsilon)-\mu_{j_0}^{\infty,n}(\lambda,\varepsilon)-\mu_{j}^{\infty,n-1}(\lambda,\varepsilon)+\mu_{j_{0}}^{\infty,n-1}(\lambda,\varepsilon)\big|.
				$$ 
				Then coming back to \eqref{asy-z1}, one gets
				\begin{align}\label{asy-z2}
					\nonumber \varrho^n_{j,j_0}(\lambda,\varepsilon)\leqslant |j-j_0|\big|r^{1,n}(\lambda,\varepsilon)-r^{1,n-1}(\lambda,\varepsilon)\big|&+\big|r_{j}^{\infty,n}(\lambda,\varepsilon)-r_{j}^{\infty,n-1}(\lambda,\varepsilon)\big|\\
					&+\big|r_{j_0}^{\infty,n}(\lambda,\varepsilon)-r_{j_0}^{\infty,n-1}(\lambda,\varepsilon)\big|.
				\end{align}
				In view of \eqref{differences mu0}, \eqref{Bio-X1}, \eqref{choice of gamma and N0 in the Nash-Moser} and the fact that $\sigma_{4}\geqslant\sigma_{3},$ one obtains 
				\begin{align*}
					\big|r^{1,n}(\lambda,\varepsilon)-r^{1,n-1}(\lambda,\varepsilon)\big|&\lesssim \varepsilon\| i_{n}-i_{n-1}\|_{q,\overline{s}_{h}+\sigma_{3}}^{\gamma,\mathcal{O}}\\
					&\lesssim \varepsilon^{2}\gamma^{-1}N_{n-1}^{-a_2}\\
					&\lesssim \varepsilon^{2-a}\langle j-j_0\rangle N_{n-1}^{-a_2}.
				\end{align*}
				In a similar line,
				using  \eqref{diffenrence rjinfty}, \eqref{Bio-X1} and \eqref{choice of gamma and N0 in the Nash-Moser}  yields %\eqref{estimate differences mujinfty} and
				\begin{align*}
					\big|r_{j}^{\infty,n}(\lambda,\varepsilon)-r_{j}^{\infty,n-1}(\lambda,\varepsilon)\big|&\lesssim \varepsilon\gamma^{-1}\| i_{n}-i_{n-1}\|_{q,\overline{s}_{h}+\sigma_{4}}^{\gamma,\mathcal{O}}\\
					&\lesssim \varepsilon^{2}\gamma^{-2}N_{n-1}^{-a_2}\\
					&\lesssim \varepsilon^{2(1-a)}\langle j-j_0\rangle N_{n-1}^{-a_2}.
				\end{align*}
				Inserting the  preceding two  estimates into \eqref{asy-z2} gives
				\begin{align}\label{asy-z3}
					\varrho^n_{j,j_0}(\lambda,\varepsilon)\lesssim \varepsilon^{2(1-a)}\langle j-j_0\rangle N_{n-1}^{-a_{2}}.
				\end{align}
				Putting together \eqref{asy-z3} and \eqref{poiH1} and using  $\gamma_{n+1}=\gamma_{n}-\varepsilon^a 2^{-n-1},$ we deduce
				\begin{align*}
					\big|{\omega}(\lambda,\varepsilon)\cdot l+\mu_{j}^{\infty,n-1}(\lambda,\varepsilon)-\mu_{j_{0}}^{\infty,n-1}(\lambda,\varepsilon)\big| 
					&\leqslant \displaystyle\tfrac{2\gamma_{n}\langle j-j_0\rangle }{\langle l\rangle^{\tau_{2}}}-{\varepsilon^a \langle j-j_{0}\rangle}2^{-n}\langle l\rangle ^{-\tau_2}\\
					&\quad+C\varepsilon^{2(1-a)}\langle j-j_0\rangle N_{n-1}^{-a_{2}}.
				\end{align*} 
				Since $|l|\leqslant N_{n-1}$,we can write
				\begin{align*}
					-{\varepsilon^a }2^{-n}\langle l\rangle ^{-\tau_2}+C\varepsilon^{2(1-a)}N_{n-1}^{-a_{2}}\leqslant {\varepsilon^a }2^{-n}\langle l\rangle ^{-\tau_2}\Big(-1+C\varepsilon^{2-3a}2^{n}N_{n-1}^{-a_{2}+\tau_{2}}\Big).
				\end{align*}
				Now remark that \eqref{param NM} and \eqref{choice of gamma and N0 in the Nash-Moser} yield in particular
				\begin{align}\label{Condor-1}
					a_2>\tau_{2}\quad\hbox{and}\quad a<\tfrac{2}{3}.
				\end{align}
				Hence, we find for  $\varepsilon$ small enough
				\begin{align*}
					\forall\, n\in\mathbb{N},\quad -1+C\varepsilon^{2-3a}2^{n}N_{n-1}^{-a_{2}+\tau_{2}}\leqslant 0
				\end{align*}
				and therefore 
				$$\big|{\omega}(\lambda,\varepsilon)\cdot l+\mu_{j}^{\infty,n-1}(\lambda,\varepsilon)-\mu_{j_{0}}^{\infty,n-1}(\lambda,\varepsilon)\big| \leqslant  \displaystyle\tfrac{2\gamma_{n}\langle j-j_0\rangle }{\langle l\rangle^{\tau_{2}}}\cdot$$
				Consequently, $\lambda \in  \mathcal{R}_{l,j,j_{0}}(i_{n-1})$ and the proof of the second point is now achieved.
				\\
				\textbf{(iii)} Let $j\in\mathbb{S}_{0}^{c}.$ In particular, one has $(l,j)\neq (0,0).$ We shall first prove that if $|l|\leqslant N_{n-1}$ and then $\mathcal{R}_{l,j}^{(1)}(i_{n}) \subset\mathcal{R}_{l,j}^{(1)}(i_{n-1}).$ As in the point (i) this implies that  $\mathcal{R}_{l,j}^{(1)}(i_{n})=\varnothing.$ Remind that the set $\mathcal{R}_{l,j}^{(1)}(i_{n})$ is defined below \eqref{set-U0}. Consider    $\lambda\in\mathcal{R}_{l,j}^{(1)}(i_{n})$ then by construction $\lambda\in \mathcal{C}_{n}^{\varepsilon}\subset \mathcal{C}_{n-1}^{\varepsilon}.$ Now by the triangle inequality we may write  in view of \eqref{estimate differences mujinfty} and \eqref{Bio-X1} and the choice $\gamma=\varepsilon^a$
				\begin{align*}
					\big|{\omega}(\lambda,\varepsilon)\cdot l+\mu_{j}^{\infty,n-1}(\lambda,\varepsilon)\big| & \leqslant  \big|{\omega}(\lambda,\varepsilon)\cdot l+\mu_{j}^{\infty,n}(\lambda,\varepsilon)\big|+|\mu_{j}^{\infty,n}(\lambda,\varepsilon)-\mu_{j}^{\infty,n-1}(\lambda,\varepsilon)|\\
					& \leqslant \displaystyle\tfrac{\gamma_{n+1}\langle j\rangle }{\langle l\rangle^{\tau_{1}}}+C\varepsilon\gamma^{-1}|j|\| i_{n}-i_{n-1}\|_{q,\overline{s}_{h}+\sigma_{4}}^{\gamma,\mathcal{O}}\\
					& \leqslant  \displaystyle\tfrac{\gamma_{n+1}\langle j\rangle }{\langle l\rangle^{\tau_{1}}}+C\varepsilon^{2(1-a)}\langle j\rangle N_{n-1}^{-a_{2}}.
				\end{align*}
				Since   $\gamma_{n+1}=\gamma_{n}-\varepsilon^a 2^{-n-1}$ and $|l|\leqslant N_{n-1}$, then 
				$$\big|{\omega}(\lambda,\varepsilon)\cdot l+\mu_{j}^{\infty,n-1}(\lambda,\varepsilon)\big| \leqslant \displaystyle\tfrac{\gamma_{n}\langle j\rangle }{\langle l\rangle^{\tau_{1}}}+\tfrac{\langle j\rangle \varepsilon^a}{2^{n+1}\langle l\rangle^{\tau_1}}\Big(-1+\varepsilon^{2-3a} 2^{n+1}N_{n-1}^{-a_{2}+\tau_{1}}\Big).$$
				Notice that \eqref{Condor-1} implies in particular 
				\begin{align}\label{samedi-1}
					a_2>\tau_{1}\quad\hbox{and}\quad a<\tfrac{2}{3}
				\end{align}
				and taking  $\varepsilon$ small enough we find that 
				$$
				\forall\,n\in \mathbb{N},\quad -1+\varepsilon^{2-3a} 2^{n+1}N_{n-1}^{-a_{2}+\tau_{1}}\leqslant 0,
				$$
				which implies in turn that 
				$$
				\big|{\omega}(\lambda,\varepsilon)\cdot l+\mu_{j}^{\infty,n-1}(\lambda,\varepsilon)\big|  \leqslant \tfrac{\gamma_{n}\langle j\rangle }{\langle l\rangle^{\tau_{1}}}.
				$$
				Consequently, $\lambda\in \mathcal{R}_{l,j}^{(1)}(i_{n-1})$ and  this ends the proof of the third point. 
				\\
				\textbf{(iv)} It is an immediate consequence of \eqref{set-U0} and the points (i)-(ii) and (iii) of Lemma \ref{lemm-dix1}.
			\end{proof}
			The next result deals with necessary  conditions such that the sets in \eqref{set-U0} are nonempty.
			\begin{lem}\label{some cantor set are empty}
				There exists $\varepsilon_0$ such that for any $\varepsilon\in[0,\varepsilon_0]$ and $n\in\mathbb{N}$ the following assertions hold true. 
				\begin{enumerate}[label=(\roman*)]
					\item Let $(l,j)\in\mathbb{Z}^{d}\times\mathbb{Z}\setminus\{(0,0)\}.$ If $\,\displaystyle\mathcal{R}_{l,j}^{(0)}(i_{n})\neq\varnothing,$ then $|j|\leqslant C_{0}\langle l\rangle.$
					\item Let $(l,j,j_{0})\in\mathbb{Z}^{d}\times(\mathbb{S}_{0}^{c})^{2}.$ If $\,\displaystyle\mathcal{R}_{l,j,j_{0}}(i_{n})\neq\varnothing,$ then $|j-j_{0}|\leqslant C_{0}\langle l\rangle.$
					\item Let $(l,j)\in\mathbb{Z}^{d}\times\mathbb{S}_{0}^{c}.$ If $\,\displaystyle \mathcal{R}_{l,j}^{(1)}(i_{n})\neq\varnothing,$ then $|j|\leqslant C_{0}\langle l\rangle.$
					\item Let $(l,j,j_{0})\in\mathbb{Z}^{d}\times(\mathbb{S}_{0}^{c})^{2}.$ There exists $c_{2}>0$ such that if $\displaystyle \min(|j|,|j_{0}|)\geqslant c_{2}\gamma_{n+1}^{-\upsilon}\langle l\rangle^{\tau_{1}},$ then 
					$$\mathcal{R}_{l,j,j_{0}}(i_{n})\subset\mathcal{R}_{l,j-j_{0}}^{(0)}(i_{n}).$$
				\end{enumerate}
			\end{lem}
			\begin{proof}
				\textbf{(i)} Assume $\mathcal{R}_{l,j}^{(0)}(i_{n})\neq\varnothing$, then we can find   $\lambda\in(\lambda_{0},\lambda_{1})$ such that, using triangle and Cauchy-Schwarz inequalities, 
				\begin{align*}
					|c_{n}(\lambda,\varepsilon)||j|&\leqslant 4|j|\gamma_{n+1}^{\upsilon}\langle l\rangle^{-\tau_{1}}+|{\omega}(\lambda,\varepsilon)\cdot l|\\
					&\leqslant 4|j|\gamma_{n+1}^{\upsilon}+C\langle l\rangle\\
					&\leqslant 8\varepsilon^{a \upsilon}|j|+C\langle l\rangle,
				\end{align*}
				where we have used $\gamma=\varepsilon^a$ and  the fact that $(\lambda,\varepsilon)\mapsto \omega(\lambda, \varepsilon)$ is bounded. Notice that 
				$$c_{n}(\lambda,\varepsilon)=I_{1}(\lambda)K_{1}(\lambda)+r^{1,n}(\lambda,\varepsilon)\quad\hbox{and}\quad \displaystyle\inf_{\lambda\in(\lambda_{0},\lambda_{1})}I_{1}(\lambda)K_{1}(\lambda):=c_{1}>0.
				$$
				Then, from  \eqref{estimate r1}, \eqref{estimate rjinfty} and Proposition \ref{Nash-Moser} $(\mathcal{P}1)_{n}$, we obtain
				\begin{align}\label{uniform estimate r1}
					\forall k\in\llbracket 0,q\rrbracket,\quad\sup_{n\in\mathbb{N}}\sup_{\lambda\in(\lambda_{0},\lambda_{1})}|\partial_{\lambda}^{k}r^{1,n}(\lambda,\varepsilon)|&\leqslant\gamma^{-k}\sup_{n\in\mathbb{N}}\| r^{1,n}\|_{q}^{\gamma,\mathcal{O}}\nonumber\\
					&\lesssim\varepsilon\gamma^{-k}\nonumber\\
					&\lesssim\varepsilon^{1-ak}.
				\end{align}
				Thus, by choosing $\varepsilon$ small enough, we can ensure by \eqref{uniform estimate r1}
				$$\inf_{n\in\mathbb{N}}\inf_{\lambda\in(\lambda_{0},\lambda_{1})}|c_{n}(\lambda,\varepsilon)|\geqslant \tfrac{c_{1}}{2}.$$
				Hence, by taking $\varepsilon$ small enough we find that  $|j|\leqslant C_{0}\langle l\rangle$ for some $C_{0}>0.$\\
				\textbf{(ii)} In the  case $j=j_0$ we get by definition $\mathcal{R}_{l,j_0,j_{0}}(i_{n})=\mathcal{R}^{(0)}_{l,0}(i_{n}),$ and then we use the point $(i).$ In what follows we take  $j\neq j_0$ and we assume that  $\mathcal{R}_{l,j,j_{0}}(i_{n})\neq\varnothing$ then there exists $\lambda\in(\lambda_{0},\lambda_{1})$ such that 
				\begin{align*}
					|\mu_{j}^{\infty,n}(\lambda,\varepsilon)-\mu_{j_{0}}^{\infty,n}(\lambda,\varepsilon)|&\leqslant 2\gamma_{n+1}|j-j_{0}|\langle l\rangle^{-\tau_{2}}+|{\omega}(\lambda,\varepsilon)\cdot l|\\
					&\leqslant 2\gamma_{n+1}|j-j_{0}|+C\langle l\rangle\\
					&\leqslant 4\varepsilon^a|j-j_{0}|+C\langle l\rangle.
				\end{align*}
				Similarly to \eqref{uniform estimate r1}, we can prove
				\begin{align}\label{uniform estimate rjinfty}
					\forall k\in\llbracket 0,q\rrbracket,\quad\sup_{n\in\mathbb{N}}\sup_{j\in\mathbb{S}_{0}^{c}}\sup_{\lambda\in(\lambda_{0},\lambda_{1})}|j||\partial_{\lambda}^{k}r_{j}^{\infty,n}(\lambda,\varepsilon)|&\leqslant\gamma^{-k}\sup_{n\in\mathbb{N}}\sup_{j\in\mathbb{S}_{0}^{c}}|j|\| r_{j}^{\infty,n}\|_{q}^{\gamma,\mathcal{O}}\nonumber\\
					&\lesssim\varepsilon\gamma^{-1-k}\nonumber\\
					&\lesssim\varepsilon^{1-a(1+k)}.
				\end{align} 
				By using the triangle inequality, Lemma \ref{lemma properties linear frequencies}-(v), \eqref{uniform estimate r1} and \eqref{uniform estimate rjinfty} we get for $j\neq j_0$,
				\begin{align*}
					|\mu_{j}^{\infty,n}(\lambda,\varepsilon)-\mu_{j_{0}}^{\infty,n}(\lambda,\varepsilon)| & \geqslant  |\Omega_{j}(\lambda)-\Omega_{j_{0}}(\lambda)|-|r^{1,n}(\lambda,\varepsilon)||j-j_{0}|-|r_{j}^{\infty,n}(\lambda,\varepsilon)|-|r_{j_{0}}^{\infty,n}(\lambda,\varepsilon)|\\
					& \geqslant  \big(C_{0}-C\varepsilon^{1-a}\big)|j-j_{0}|\\
					& \geqslant \tfrac{C_{0}}{2}|j-j_{0}|
				\end{align*}
				provided that $\varepsilon$ is small enough. Putting together the previous inequalities yields  for $\varepsilon $ small enough $|j-j_{0}|\leqslant C_{0}\langle l\rangle,$ for some $C_{0}>0.$\\
				%\ding{226} Assume $j\neq j_{0}.$ Assume by contradiction that $\mathcal{R}_{0,j,j_{0}}(i_{m})\neq
				%\varnothing.$ There exists $\lambda\in]\lambda_{0},\lambda_{1}[$ such that
				%$$\frac{C_{0}}{2}|j-j_{0}|\leqslant|\mu_{j}^{\infty,m}(\lambda,\varepsilon)-\mu_{j_{0}}^{\infty,m}(\lambda,\varepsilon)|\leqslant 2\gamma_{m+1}|j-j_{0}|\leqslant 4\gamma|j-j_{0}|.$$
				%Hence $\frac{c_{0}}{2}\leqslant 4\gamma$ which is a contradiction with the fact that we took $\gamma\leqslant\frac{c_{0}}{16}.$\\
				\textbf{(iii)} First remark that the case  $j=0$ is trivial. Now  for $j\neq 0$ we assume that  $\mathcal{R}_{l,j}^{(1)}(i_{n})\neq\varnothing$ then there exists   $\lambda\in(\lambda_{0},\lambda_{1})$ such that 
				\begin{align*}
					|\mu_{j}^{\infty,n}(\lambda,\varepsilon)|&\leqslant \gamma_{n+1}|j|\langle l\rangle^{\tau_{1}}+|{\omega}(\lambda,\varepsilon)\cdot l|\\
					&\leqslant  2\varepsilon^a|j|+C\langle l\rangle.
				\end{align*}
				Using the definition \eqref{asy-z1} combined with 
				the triangle  inequality, Lemma \ref{lemma properties linear frequencies}-(iv), \eqref{uniform estimate r1} and \eqref{uniform estimate rjinfty}, we get
				\begin{align*}
					|\mu_{j}^{\infty,n}(\lambda,\varepsilon)|&\geqslant \Omega|j|-|j||r^{1,n}(\lambda,\varepsilon)|-|r_{j}^{\infty,n}(\lambda,\varepsilon)|\\
					&
					\geqslant \Omega|j|-C\varepsilon^{1-a}|j|.
				\end{align*}
				Combining the previous two inequalities and the second condition in \eqref{samedi-1} implies
				\begin{align*}
					\big( \Omega-C\varepsilon^{1-a}-2\varepsilon^a\big)|j|
					&\leqslant  C\langle l\rangle.
				\end{align*}
				Thus, by taking $\varepsilon$ small enough we obtain $|j|\leqslant C_0\langle  l\rangle,$ for some $C_{0}>0.$\\
				\textbf{(iv)} First notice that the case $j=j_0$ is trivial and follows from the definition \eqref{set-U0}. Let $j\neq j_0$ and $\lambda\in\mathcal{R}_{l,j,j_{0}}(i_{n})$, then by definition
				$$\big|{\omega}(\lambda,\varepsilon)\cdot l+\mu_{j}^{\infty,n}(\lambda,\varepsilon)-\mu_{j_{0}}^{\infty,n}(\lambda,\varepsilon)\big|\leqslant\tfrac{2\gamma_{n+1}\langle j-j_{0}\rangle}{\langle l\rangle^{\tau_{2}}}.
				$$
				Combining \eqref{asy-z1} and \eqref{def eigenvalues at equilibrium} with the triangle inequality we infer
				\begin{align*}
					\big|{\omega}(\lambda,\varepsilon)\cdot l+(j-j_{0})c_{n}(\lambda,\varepsilon)\big|  &\leqslant  \big|{\omega}(\lambda,\varepsilon)\cdot l+\mu_{j}^{\infty,n}(\lambda,\varepsilon)-\mu_{j_{0}}^{\infty,n}(\lambda,\varepsilon)\big|\\
					&\quad+|jI_{j}(\lambda)K_{j}(\lambda)-j_{0}I_{j_{0}}(\lambda)K_{j_{0}}(\lambda)|+\big|r_{j}^{\infty,n}(\lambda,\varepsilon)-r_{j_{0}}^{\infty,n}(\lambda,\varepsilon)\big|.
				\end{align*}
				Thus, we find
				\begin{align}\label{ZaraX1}
					\nonumber \big|{\omega}(\lambda,\varepsilon)\cdot l+(j-j_{0})c_{n}(\lambda,\varepsilon)\big|  \leqslant & \tfrac{2\gamma_{n+1}\langle j-j_{0}\rangle}{\langle l\rangle^{\tau_{2}}}+|jI_{j}(\lambda)K_{j}(\lambda)-j_{0}I_{j_{0}}(\lambda)K_{j_{0}}(\lambda)|\\
					&+\big|r_{j}^{\infty,n}(\lambda,\varepsilon)-r_{j_{0}}^{\infty,n}(\lambda,\varepsilon)\big|.
				\end{align}
				Without loss of generality, we can assume that $|j_0|\geqslant |j|$ and remind that $j\neq j_0$. Then, from \eqref{HUP1} and \eqref{decayk}, we easily  find
				\begin{align*}
					|jI_{j}(\lambda)K_{j}(\lambda)-j_{0}I_{j_{0}}(\lambda)K_{j_{0}}(\lambda)|&\leqslant|j|\big| I_{j}(\lambda)K_{j}(\lambda)-I_{j_{0}}(\lambda)K_{j_{0}}(\lambda)\big|+|j-j_0|\big|I_{j_{0}}(\lambda)K_{j_{0}}(\lambda)\big|\\
					& \leqslant \tfrac{\langle j-j_{0}\rangle}{\min(|j|,|j_{0}|)}\cdot
				\end{align*}
				Applying  \eqref{estimate rjinfty}, we find for $j\neq j_{0}\in\mathbb{S}_{0}^{c}$,
				\begin{align*}
					\big|r_{j}^{\infty,n}(\lambda,\varepsilon)-r_{j_{0}}^{\infty,n}(\lambda,\varepsilon)\big|  \leqslant &C \varepsilon^{1-a}\big(|j|^{-1}+|j_0|^{-1}\big)\\
					\leqslant & C \varepsilon^{1-a} \tfrac{\langle j-j_{0}\rangle}{\min(|j|,|j_{0}|)}\cdot
				\end{align*}
				Plugging the preceding estimates into \eqref{ZaraX1} yields 
				\begin{align*}
					\nonumber \big|{\omega}(\lambda,\varepsilon)\cdot l+(j-j_{0})c_{n}(\lambda,\varepsilon)\big|  \leqslant & \tfrac{2\gamma_{n+1}\langle j-j_{0}\rangle}{\langle l\rangle^{\tau_{2}}}+ C \tfrac{\langle j-j_{0}\rangle}{\min(|j|,|j_{0}|)}\cdot
				\end{align*}
				Therefore,  if we assume $\displaystyle \min(|j|,|j_{0}|)\geqslant \tfrac{1}{2} C\gamma_{n+1}^{-\upsilon}\langle l\rangle^{\tau_{1}}$  and $\tau_{2}>\tau_{1}$, then we deduce 
				$$\big|{\omega}(\lambda,\varepsilon)\cdot l+(j-j_{0})c_{n}(\lambda,\varepsilon)\big| \leqslant  \tfrac{4\gamma_{n+1}^{\upsilon}\langle j-j_{0}\rangle}{\langle l\rangle^{\tau_{1}}}\cdot$$
				This ends the proof of the lemma by taking $c_{2}=\frac{C}{2}.$
			\end{proof}
			We shall now establish that the perturbed frequencies $\omega(\lambda,\varepsilon)$ satisfy the Rüssemann conditions. This is done by a perturbation argument from the equilibrium linear frequencies $\omega_{\textnormal{Eq}}(\lambda)$ for which we already know by Lemma \ref{lemma transversality} that they satisfy the transversality conditions. 
			\begin{lem}\label{lemma Russeman condition for the perturbed frequencies}
				{Let $q_{0}$, $C_{0}$ and $\rho_{0}$ as in Lemma \ref{lemma transversality}. There exist $\varepsilon_{0}>0$ small enough such that for any   $\varepsilon\in[0,\varepsilon_{0}]$ the following assertions hold true.
					\begin{enumerate}[label=(\roman*)]
						\item For all $l\in\mathbb{Z}^{d}\setminus\{0\}, $ we have 
						$$\inf_{\lambda\in[\lambda_{0},\lambda_{1}]}\max_{k\in\llbracket 0,q_{0}\rrbracket}\big|\partial_{\lambda}^{k}\left({\omega}(\lambda,\varepsilon)\cdot l\right)\big|\geqslant\tfrac{\rho_{0}\langle l\rangle}{2}.$$
						\item For all $(l,j)\in\mathbb{Z}^{d+1}\setminus\{(0,0)\}$ such that $|j|\leqslant C_{0}\langle l\rangle,$ we have
						$$\forall n\in\mathbb{N},\quad\inf_{\lambda\in[\lambda_{0},\lambda_{1}]}\max_{k\in\llbracket 0,q_{0}\rrbracket}|\partial_{\lambda}^{k}\big({\omega}(\lambda,\varepsilon)\cdot l+jc_{n}(\lambda,\varepsilon)\big)|\geqslant\tfrac{\rho_{0}\langle l\rangle}{2}.$$
						\item For all $(l,j)\in\mathbb{Z}^{d}\times\mathbb{S}_{0}^{c}$ such that $|j|\leqslant C_{0}\langle l\rangle,$ we have
						$$\forall n\in\mathbb{N},\quad\inf_{\lambda\in[\lambda_{0},\lambda_{1}]}\max_{k\in\llbracket 0,q_{0}\rrbracket}\big|\partial_{\lambda}^{k}\big({\omega}(\lambda,\varepsilon)\cdot l+\mu_{j}^{\infty,n}(\lambda,\varepsilon)\big)\big|\geqslant\tfrac{\rho_{0}\langle l\rangle}{2}.$$
						\item For all $(l,j,j_{0})\in\mathbb{Z}^{d}\times(\mathbb{S}_{0}^{c})^{2}$ such that $|j-j_{0}|\leqslant C_{0}\langle l\rangle,$ we have
						$$\forall n\in\mathbb{N},\quad\inf_{\lambda\in[\lambda_{0},\lambda_{1}]}\max_{k\in\llbracket 0,q_{0}\rrbracket}\big|\partial_{\lambda}^{k}\big({\omega}(\lambda,\varepsilon)\cdot l+\mu_{j}^{\infty,n}(\lambda,\varepsilon)-\mu_{j_{0}}^{\infty,n}(\lambda,\varepsilon)\big)\big|\geqslant\tfrac{\rho_{0}\langle l\rangle}{2}.$$
				\end{enumerate}}
			\end{lem}
			\begin{proof}
				\textbf{(i)} From the triangle and Cauchy-Schwarz  inequalities together with  \eqref{estimate repsilon1}, \eqref{choice of gamma and N0 in the Nash-Moser} and Lemma \ref{lemma transversality}-(i), we deduce
				\begin{align*}
					\displaystyle\max_{k\in\llbracket 0,q_{0}\rrbracket}|\partial_{\lambda}^{k}\left({\omega}(\lambda,\varepsilon)\cdot l\right)| & \geqslant \displaystyle\max_{k\in\llbracket 0,q_{0}\rrbracket}|\partial_{\lambda}^{k}\left({\omega}_{\textnormal{Eq}}(\lambda)\cdot l\right)|-\max_{k\in\llbracket 0,q\rrbracket}|\partial_{\lambda}^{k}\left(\overline{\mathrm{r}}_{\varepsilon}(\lambda)\cdot l\right)|\\
					& \geqslant\displaystyle\rho_{0}\langle l\rangle-C\varepsilon\gamma^{-1-q}N_{0}^{q\overline{a}}\langle l\rangle\\
					& \geqslant \displaystyle\rho_{0}\langle l\rangle-C\varepsilon^{1-a(1+q+q\overline{a})}\langle l\rangle\\
					& \geqslant \displaystyle\tfrac{\rho_{0}\langle l\rangle}{2}
				\end{align*}
				provided that $\varepsilon$ is small enough and 
				\begin{equation}\label{condition param Russ}
					1-a(1+q+q\overline{a})>0.
				\end{equation}
				Notice that the condition \eqref{condition param Russ} is automatically satisfied by \eqref{choice of gamma and N0 in the Nash-Moser} and \eqref{param NM}.\\
				\textbf{(ii)} As before, using the triangle and Cauchy-Schwarz inequalities combined with  \eqref{estimate repsilon1}, \eqref{uniform estimate r1}, Lemma \ref{lemma transversality}-(ii) and the fact that $|j|\leqslant C_{0}\langle l\rangle$, we get
				\begin{align*}
					\displaystyle\max_{k\in\llbracket 0,q_{0}\rrbracket}|\partial_{\lambda}^{k}\left({\omega}(\lambda,\varepsilon)\cdot l+jc_{n}(\lambda,\varepsilon)\right)| & \geqslant \displaystyle\max_{k\in\llbracket 0,q_{0}\rrbracket}|\partial_{\lambda}^{k}\left({\omega}_{\textnormal{Eq}}(\lambda)\cdot l+jI_{1}(\lambda)K_{1}(\lambda)\right)|\\
					& \displaystyle\quad-\max_{k\in\llbracket 0,q\rrbracket}|\partial_{\lambda}^{k}\left(\overline{\mathrm{r}}_{\varepsilon}(\lambda)\cdot l+jr^{1,n}(\lambda,\varepsilon)\right)|\\
					& \geqslant \displaystyle\rho_{0}\langle l\rangle-C\varepsilon^{1-a(1+q+q\overline{a})}\langle l\rangle-C\varepsilon^{1-aq}|j|\\
					& \geqslant \displaystyle\tfrac{\rho_{0}\langle l\rangle}{2}
				\end{align*}
				for $\varepsilon$ small enough and with the condition \eqref{condition param Russ}.\\
				\textbf{(iii)} As before, performing the triangle  and Cauchy-Schwarz inequalities combined with \eqref{estimate repsilon1}, \eqref{uniform estimate r1}, \eqref{uniform estimate rjinfty}, Lemma \ref{lemma transversality}-(iii) and the fact that  $|j|\leqslant C_{0}\langle l\rangle$, we get
				\begin{align*}
					\displaystyle\max_{k\in\llbracket 0,q_{0}\rrbracket}\big|\partial_{\lambda}^{k}\big({\omega}(\lambda,\varepsilon)\cdot l+\mu_{j}^{\infty,n}(\lambda,\varepsilon)\big)\big| & \geqslant  \displaystyle\max_{k\in\llbracket 0,q_{0}\rrbracket}|\partial_{\lambda}^{k}\left({\omega}_{\textnormal{Eq}}(\lambda)\cdot l+\Omega_{j}(\lambda)\right)|\\
					& \displaystyle\quad-\max_{k\in\llbracket 0,q\rrbracket}\big|\partial_{\lambda}^{k}\big(\overline{\mathrm{r}}_{\varepsilon}(\lambda)\cdot l+jr^{1,n}(\lambda,\varepsilon)+r_{j}^{\infty,n}(\lambda,\varepsilon))\big)\big|\\
					& \geqslant  \displaystyle\rho_{0}\langle l\rangle-C\varepsilon^{1-a(1+q+q\overline{a})}\langle l\rangle-C\varepsilon^{1-a(1+q)}|j|\\
					& \geqslant  \displaystyle\tfrac{\rho_{0}\langle l\rangle}{2}
				\end{align*}
				for $\varepsilon$ small enough with the condition \eqref{condition param Russ}.\\
				\textbf{(iv)} Arguing as in the preceding point, using  \eqref{uniform estimate r1}, \eqref{uniform estimate rjinfty}, Lemma \ref{lemma transversality}-(iv)-(v) and the fact that $0<|j-j_{0}|\leqslant C_{0}\langle l\rangle$ (notice that the case $j=j_0$ is trivial), we have 
				\begin{align*}
					\max_{k\in\llbracket 0,q_{0}\rrbracket}\big|\partial_{\lambda}^{k}\big({\omega}(\lambda,\varepsilon)\cdot l&+\mu_{j}^{\infty,n}(\lambda,\varepsilon)-\mu_{j_{0}}^{\infty,n}(\lambda,\varepsilon)\big)\big|
					\geqslant\displaystyle\max_{k\in\llbracket 0,q_{0}\rrbracket}\big|\partial_{\lambda}^{k}\big({\omega}_{\textnormal{Eq}}(\lambda)\cdot l+\Omega_{j}(\lambda)-\Omega_{j_{0}}(\lambda)\big)\big|\displaystyle\\
					&-\max_{k\in\llbracket 0,q\rrbracket}\big|\partial_{\lambda}^{k}\big(\overline{\mathrm{r}}_{\varepsilon}(\lambda)\cdot l+(j-j_{0})r^{1,n}(\lambda,\varepsilon)+r_{j}^{\infty,n}(\lambda,\varepsilon)-r_{j_{0}}^{\infty,n}(\lambda,\varepsilon)\big)\big|\\
					&\geqslant\displaystyle\rho_{0}\langle l\rangle-C\varepsilon^{1-a(1+q+q\overline{a})}\langle l\rangle-C\varepsilon^{1-a(1+q)}|j-j_{0}|\\
					&\geqslant \displaystyle\tfrac{\rho_{0}\langle l\rangle}{2}
				\end{align*}
				for $\varepsilon$ small enough. This ends the proof of Lemma \ref{lemma Russeman condition for the perturbed frequencies}.
			\end{proof}
			
			\appendix
			\section{Modified Bessel functions}\label{appendix Bessel}
			In this short section we shall collect some properties about Bessel and modified Bessel functions that were used in the preceding sections. We refer to \cite{W95} for an almost exhaustive presentation of these special functions.\\
			We define first the Bessel functions of order $\nu\in\mathbb{C}$ by 
			$$J_{\nu}(z)=\sum_{m=0}^{\infty}\frac{(-1)^{m}\left(\frac{z}{2}\right)^{\nu +2m}}{m!\Gamma(\nu+m+1)},\quad|\mbox{arg}(z)|<\pi.$$
			Notice that when $\nu\in\mathbb{N}$ we have the following integral representation, see \cite[p. 115]{Leb65}.
			\begin{equation}\label{Besse-repr}
				J_\nu(x)=\frac{1}{\pi}\int_0^\pi\cos\big(x \sin \theta-\nu \theta\big) d\theta.
			\end{equation}
			We shall also introduce the  Bessel functions of imaginary argument also called modified Bessel functions of first and second kind
			\begin{equation}\label{definition of modified Bessel function of first kind}
				I_{\nu}(z)=\sum_{m=0}^{\infty}\frac{\left(\frac{z}{2}\right)^{\nu+2m}}{m!\Gamma(\nu+m+1)},\quad|\mbox{arg}(z)|<\pi
			\end{equation}
			and $$K_{\nu}(z)=\frac{\pi}{2}\frac{I_{-\nu}(z)-I_{\nu}(z)}{\sin(\nu\pi)},\quad\nu\in\mathbb{C}\setminus\mathbb{Z},\quad|\mbox{arg}(z)|<\pi.$$
			For $j\in\mathbb{Z},$ we define $K_{j}(z)=\displaystyle\lim_{\nu\rightarrow j}K_{\nu}(z).$
			We give now useful properties of modified Bessel functions.\\ 
			$\blacktriangleright$ \textit{Symmetry and positivity,} see \cite[p. 375]{AS64}
			\begin{equation}\label{symmetry Bessel}
				\forall j\in\mathbb{N},\quad\forall\lambda\in\mathbb{R}_{+}^{*},\quad I_{-j}(\lambda)=I_{j}(\lambda)\in\mathbb{R}_{+}^{*}\quad\mbox{and}\quad K_{-j}(\lambda)=K_{j}(\lambda)\in\mathbb{R}_{+}^{*}.
			\end{equation}
			$\blacktriangleright$ {\it Anti-derivative,} see \cite[p. 376]{AS64}. If we set $\mathcal{Z}_{\nu}(z)=I_{\nu}(z)$ or $e^{i\nu\pi}K_{\nu}(z)$, then for all $\nu\in\mathbb{R}$, we have 
			\begin{equation}\label{Bessel and anti-derivatives}
				\frac{d}{dz}\left(z^{\nu+1}\mathcal{Z}_{\nu+1}(z)\right)=z^{\nu+1}\mathcal{Z}_{\nu}(z).
			\end{equation}
			$\blacktriangleright$ {\it Power series expansion for $K_{j}$,} see \cite[p. 375]{AS64}.\\
			\begin{equation}\label{power series Kj}
				\begin{array}{l}
					\displaystyle K_{j}(z)=\tfrac{1}{2}\left(\tfrac{z}{2}\right)^{-j}\sum_{k=0}^{j-1}\tfrac{(j-k-1)!}{k!}\left(\tfrac{-z^{2
					}}{4}\right)^{k}+(-1)^{n+1}\log\left(\tfrac{z}{2}\right)I_{j}(z)\\
					\displaystyle\mbox{\hspace{4cm}}+\tfrac{1}{2}\left(\tfrac{-z}{2}\right)^{j}\sum_{k=0}^{\infty}\left(\psi(k+1)+\psi(j+k+1)\right)\tfrac{\left(\frac{z^{2}}{4}\right)^{k}}{k!(j+k)!},
				\end{array}
			\end{equation}
			where $\psi(1)=-\boldsymbol{\gamma}$ (Euler's constant) \quad and \quad $\forall m\in\mathbb{N}^{*},\,\psi(m+1)=\displaystyle\sum_{k=1}^{m}\tfrac{1}{k}-\boldsymbol{\gamma}.$\\
			In particular we have the expansion
			\begin{equation}\label{explicit form for K0}
				K_{0}(z)=-\log\left(\frac{z}{2}\right)I_{0}(z)+\sum_{m=0}^{\infty}\frac{\left(\frac{z}{2}\right)^{2m}}{(m!)^{2}}\psi(m+1).
			\end{equation}
			\\
			$\blacktriangleright$ {\it Integral representation for $K_{\nu}$,} see \cite[p. 133]{Leb65}
			For all $a,b>0$ for any $\nu,\mu\in\mathbb{C}$ satisfying $-1<\textnormal{Re}(\nu)<2\textnormal{Re}(\mu)+\frac{3}{2}$ one has
			\begin{equation}\label{integral representation for Knu}
				\int_{0}^{\infty}\frac{x^{\nu+1}J_{\nu}\left(bx\right)}{(x^{2}+a^{2})^{\mu+1}}dx=\frac{a^{\nu-\mu}b^{\mu}}{2^{\mu}\Gamma(\mu+1)}K_{\nu-\mu}(ab).
			\end{equation}
			$\blacktriangleright$ {\it Nicholson's integral representation,} see \cite[p. 441]{W95}. Let $ j\in\mathbb{N}$ then
			\begin{equation}\label{Bessel}
				(I_{j}K_{j})(z)=\tfrac{2(-1)^{j}}{\pi}\int_{0}^{\frac{\pi}{2}}K_{0}(2z\cos(\tau))\cos(2j\tau)d\tau.
			\end{equation}
			Another similar representation can be found in \cite[p. 140]{Leb65}
			\begin{align}\label{form-Lebe0}
				(I_{j}K_{j})(\lambda)=\frac12\int_{0}^\infty J_0\big(2\lambda\sinh(t/2)\big)e^{-jt}dt.
			\end{align}
			$\blacktriangleright$ {\it Holomorphic property of the product $I_{j}K_{j}$.} Let $j\in\mathbb{N}$ then the function $z\mapsto (I_{j}K_{j})(z)$ is holomorphic on the half plane $\mbox{Re}(z)>0.$\\
			$\blacktriangleright$ {\it Monotonicity  of  $I_{\nu}K_{\nu}$}, see for instance  \cite{B09,DHR19}.The map $(\lambda,\nu)\in(\mathbb{R}_{+}^{*})^{2}\mapsto I_{\nu}(\lambda)K_{\nu}(\lambda)$ is strictly decreasing in each variable. \\
			$\blacktriangleright$ {\it Asymptotic expansion of small argument,} see for instance \cite[p. 375]{AS64}.
			\begin{equation}\label{asymptotic expansion of small argument}
				\forall j\in\mathbb{N}^{*},\quad I_{j}(\lambda)\underset{\lambda\rightarrow 0}{\sim}\tfrac{\left(\frac{1}{2}\lambda\right)^{j}}{\Gamma(j+1)}\quad \mbox{ and }\quad K_{j}(\lambda)\underset{\lambda\rightarrow 0}{\sim}\tfrac{\Gamma(j)}{2\left(\frac{1}{2}\lambda\right)^{j}}.
			\end{equation} 
			$\blacktriangleright$ {\it Asymptotic expansion of large argument for the product $I_{j}K_{j}$}, see for instance \cite[p. 378]{AS64}. 
			\begin{equation}\label{asymptotic expansion of large argument}
				\forall N\in\mathbb{N}^{*},\quad I_{j}(\lambda)K_{j}(\lambda)\underset{\lambda\rightarrow\infty}{\sim}\frac{1}{2\lambda}\left(1+\sum_{m=1}^{N}\tfrac{\alpha_{j,m}}{(2\lambda)^{2m}}\right),
			\end{equation}
			with
			\begin{align}\label{Form-Pol}\alpha_{j,m}=(-1)^{m}\tfrac{(2m)!}{4^m\big(m!\big)^2}P_m(\mu_j), \quad P_m(X)=\prod_{\ell=1}^{m}\big(X-(2\ell-1)^2\big),\quad \mu_j=4{j^2}.
			\end{align}
			In particular, 
			\begin{equation}\label{limit at infinity of product Bessel}
				I_{j}(\lambda)K_{j}(\lambda)\underset{\lambda\rightarrow\infty}{\longrightarrow}0.
			\end{equation}
			$\blacktriangleright$ {\it  Asymptotic expansion of high order for the product $I_{j}K_{j}$}, for more details see \cite{HS11}.
			\begin{equation}\label{asymptotic expansion of high order}
				I_{j}(\lambda)K_{j}(\lambda)\underset{j\rightarrow\infty}{\sim}\tfrac{1}{2j}\left(\sum_{m=0}^{\infty}\tfrac{b_{m}(\lambda)}{j^{m}}\right)\left(\sum_{m=0}^{\infty}(-1)^{m}\tfrac{b_{m}(\lambda)}{j^{m}}\right),
			\end{equation}
			where for each $m\in\mathbb{N}$, $b_{m}(\lambda)$ is a polynomial of degree $m$ in $\lambda^{2}$ defined by 
			$$b_{0}(\lambda)=1\quad\mbox{ and }\quad\forall m\in\mathbb{N}^{*},\,b_{m}(\lambda)=\sum_{k=1}^{m}(-1)^{m-k}\tfrac{S(m,k)}{k!}\left(\tfrac{\lambda^{2}}{4}\right)^{k}$$
			and the $S(m,k)$ are Stirling numbers of second kind defined recursively by 
			$$\forall (m,k)\in\mathbb{N}^{*}\times\mathbb{N},\quad S(m,k)=S(m-1,k-1)+kS(m-1,k),$$
			with $$S(0,0)=1,\mbox{ }\forall m\in\mathbb{N}^{*},S(m,1)=1\mbox{ and }S(m,0)=0\mbox{ and if }m<k\mbox{ then }S(m,k)=0.$$
			
			We shall also prove the following result which has been  frequently  used before.
		\begin{lem}\label{lemma sum Nn}
			Let $N_{0}\geqslant 2.$ Consider the sequence $(N_{m})_{m\in\mathbb{N}}$ defined by \eqref{definition of Nm}. Then for all $\alpha>0$, we have
			$$\sum_{k=m}^{\infty}N_{k}^{-\alpha}\underset{m\rightarrow\infty}{\sim}N_{m}^{-\alpha}.$$
		\end{lem}
		\begin{proof}
			We consider the positive decaying function 
			$$t\in\mathbb{R}_{+}^{*}\mapsto N_{0}^{-\alpha\left(\frac{3}{2}\right)^{t}}=\exp\left(-\alpha\ln(N_{0})e^{t\ln\left(\frac{3}{2}\right)}\right),$$
			and apply to it a series-integral comparison, namely
			$$\sum_{k=m+1}^{\infty}N_{k}^{-\alpha}\leqslant\int_{m}^{\infty}\exp\left(-\alpha\ln(N_{0})e^{t\ln\left(\frac{3}{2}\right)}\right)dt=\int_{0}^{\infty}\exp\left(-\alpha\ln(N_{0})e^{u\ln\left(\frac{3}{2}\right)}e^{m\ln\left(\frac{3}{2}\right)}\right)du.$$
			Now remark that
			$$N_{m}^{\alpha}\exp\left(-\alpha\ln(N_{0})e^{u\ln\left(\frac{3}{2}\right)}e^{m\ln\left(\frac{3}{2}\right)}\right)=\exp\left(\alpha\ln(N_{0})\left(1-e^{u\ln\left(\frac{3}{2}\right)}\right)e^{m\ln\left(\frac{3}{2}\right)}\right).$$
			Since 
			$$\forall u\in\mathbb{R}_{+}^{*},\quad 1-e^{u\ln\left(\frac{3}{2}\right)}<0,$$
			then we deduce that
			$$\forall u\in\mathbb{R}_{+}^{*},\quad N_{m}^{\alpha}\exp\left(-\alpha\ln(N_{0})e^{u\ln\left(\frac{3}{2}\right)}e^{m\ln\left(\frac{3}{2}\right)}\right)\underset{m\rightarrow\infty}{\longrightarrow}0$$
			and
			$$\forall u\in\mathbb{R}_{+}^{*},\forall m\in\mathbb{N},\quad 0\leqslant N_{m}^{\alpha}\exp\left(-\alpha\ln(N_{0})e^{u\ln\left(\frac{3}{2}\right)}e^{m\ln\left(\frac{3}{2}\right)}\right)\leqslant N_{0}^{\alpha}\exp\left(-\alpha\ln(N_{0})e^{u\ln\left(\frac{3}{2}\right)}\right)\in L^{1}(\mathbb{R}_{+}).$$
			Applying dominated convergence theorem, we obtain
			$$\sum_{k=m+1}^{\infty}N_{k}^{-\alpha}\underset{m\rightarrow\infty}{=}o\left(N_{m}^{-\alpha}\right).$$
			As a consequence
			$$\sum_{k=m}^{\infty}N_{k}^{-\alpha}=N_{m}^{-\alpha}+\sum_{k=m+1}^{\infty}N_{k}^{-\alpha}\underset{m\rightarrow\infty}{\sim}N_{m}^{-\alpha}.$$	
		\end{proof}

			\newpage
			
			%\bibliography{tho}

\begin{thebibliography}{9999}
			 	\bibitem{AB11} T. Alazard, P. Baldi, {\it Gravity capillary standing water waves}, Arch. Ration. Mech. Anal. 217 (2015), no. 3, 741--830.
			 	
			 	%\bibitem{A-H} M. H. P. \ Ambaum and B. J. \ Harvey, {\it  Perturbed Rankine vortices in surface quasi-geostrophic dynamics}. Geophysical and Astrophysical Fluid Dynamics, 105  (2011), no. 4-5, 377--391. 
			 	
			 	%\bibitem{A83} H. Aref,  {\it Integrable, chaotic, and turbulent vortex motion in two-dimensional flows.}  Ann. Rev. Fluid Mech., 15 (1983),  345--389.
			 	
			 	\bibitem{AS64} M. Abramowitz, I. A. Stegun, \textit{Handbook of mathematical functions with formulas, graphs, and mathematical tables}, volume 55 of National Bureau of Standards Applied Mathematics Series, (1964).
			 	
			 	\bibitem{A63} V. I. Arnold, {\it Small denominators and problems of stability of motion in classical mechanics and celestial mechanics,} Uspekhi Mat. Nauk 18 (1963), 91--192.
			 	
			 	\bibitem{BBMH18} P. Baldi, M. Berti, E. Haus, R. Montalto, {\it Time quasi-periodic gravity water waves in finite depth,} Invent. Math. 214 (2018), no. 2, 739--911.
			 	
			 	\bibitem{BBM14} P. Baldi, M. Berti, R. Montalto, {\it KAM for quasi-linear and fully nonlinear forced perturbations of Airy equation,} Math. Ann. 359 (2014), no. 1-2, 471--536.
			 	
			 	\bibitem{BBM16} P. Baldi, M. Berti, R. Montalto, {\it KAM for autonomous quasi-linear perturbations of KdV}, 
			 	Ann. Inst. H. Poincar\'e Analyse Non. Lin. 33 (2016), no. 6, 1589--1638. 
			 	
			 	\bibitem{BM20} P. Baldi, R. Montalto, \textit{Quasi-periodic incompresible Euler flows in 3D}, Advances in Mathematics 384 (2021), 107730.
			 	
			 	\bibitem{BBM11} D. Bambusi, M. Berti, E. Magistrelli, {\it Degenerate KAM theory for partial differential equations}, Journal Diff. Equations, 250 (2011), no. 8, 3379--3397.
			 	\bibitem{B09} A. Baricz, \textit{On a product of modified Bessel functions}, Proc. Amer. Math. Soc. 137 (2009), no. 1, 189--193. 
			 	
			 	%\bibitem{BP13} A. Baricz, S. Ponnusamy, \textit{On Tur\'{a}n type inequalities for modified Bessel functions}, Proc. Amer. Math. Soc. 141 (2013), no. 2, 523--532.
			 	
			 	\bibitem{B19} M. Berti, \textit{KAM theory for partial differential equations}, Analysis in Theory and Applications 35 (2019), no. 3, 235--267.
			 	
			 	%\bibitem{BBP10} M. Berti, P. Bolle, M. Procesi, \textit{An abstract Nash-Moser theorem with parameters and applications to PDEs}, Ann. Inst. H. Poincar\'e Anal. Non Lin\'eaire 27 (2010), no. 1, 377--399.
			 	
			 	\bibitem{BB15} M. Berti, P. Bolle, {\it A Nash-Moser approach to KAM theory}, Fields Institute Communications, special volume ``Hamiltonian PDEs and Applications'' (2015), 255--284.
			 	
			 	\bibitem{BFM21} M. Berti, L. Franzoi, A. Maspero, \textit{Traveling quasi-periodic water waves with constant vorticity}, Archive for Rational Mechanics and Analysis 240 (2021), 99--202.
			 	
			 	\bibitem{BFM21-1} M. Berti, L. Franzoi, A. Maspero, \textit{Pure gravity traveling quasi-periodic water waves with constant vorticity}, arXiv:2101.12006.
			 	
			 	\bibitem{BM18} M. Berti , R. Montalto, {\it Quasi-periodic standing wave solutions of gravity-capillary water waves}, MEMO, Volume 263, 1273, Memoires AMS, ISSN 0065-9266, (2020).
			 	
			 	\bibitem{BHM} M. Berti, Z. Hassainia, N. Masmoudi.	{\it Time quasi-periodic vortex patches}. In progress.
			 	
			 	
			 	%\bibitem{BCP15} M. Berti, L. Corsi, M. Procesi  {\it An Abstract Nash Moser Theorem and Quasi-Periodic Solutions for NLW and NLS on Compact Lie Groups and Homogeneous Manifolds.} Comm. Math. Phys. 334 (2015), no. 3, 1413--1454.
			 	
			 	\bibitem{BC93} A. L. Bertozzi, P. Constantin, {\it  Global regularity for vortex patches,} Comm. Math. Phys. 152 (1993), no. 1, 9--28.
			 	
			 	\bibitem{BM02} A. L. Bertozzi, A. J. Majda, {\it Vorticity and Incompressible Flow,} Cambridge texts in applied Mathematics, Cambridge University Press, Cambridge, (2002).
			 	
			 	%\bibitem{BGLV16} A. L. Bertozzi, J. B. Garnett, T. Laurent, J. Verdera, \textit{The regularity of the boundary of a multidimensional aggregation patch}, SIAM J. Math. Anal. 48 (2016), no. 6, 3789--3819.
			 	
			 	%\bibitem{Bourgain}  J.  Bourgain, {\it Construction of quasi-periodic solutions for Hamiltonian perturbations of linear equations and applications to nonlinear PDE}. Internat. Math. Res. Notices (1994), no. 11, 475ff., approx. 21 pp
			 	
			 	\bibitem{B82} J.\ Burbea, {\it Motions of vortex patches,} Lett. Math. Phys. 6 (1982), no. 1,  1--16.
			 	
			 	%\bibitem{CMOV} J. C. Cantero, J. Mateu, J. Orobitg, J. Verdera, \textit{Regularity of the boundary of vortex patches for some non-linear transport equations}, arXiv:2103.05356.
			 	
			 	\bibitem{CCG16} A. Castro, D. C\'ordoba, J. G\'omez-Serrano, {\it   Uniformly rotating analytic global patch solutions for active scalars,} Ann. PDE, 2 (2016), no. 1, 1--34.
			 	
			 	%\bibitem{CCG16-2} A. Castro, D. C\'ordoba, J. G\'omez-Serrano, {\it  Global smooth solutions for the inviscid SQG equation}, arXiv:1603.03325, (2016).
			 	
			 	%\bibitem{CCG19} A. Castro, D. C\'ordoba, J. G\'omez-Serrano, {\it  Uniformly rotating smooth solutions for the incompressible 2D Euler equations}, Arch. Ration. Mech. Anal. 231 (2019), no. 2, 719--785.
			 	
			 	\bibitem{CCG16-3} A. Castro, D. C\'ordoba, J. G\'omez-Serrano, {\it  Existence and regularity of rotating global solutions for the generalized surface quasi-geostrophic equations}. Duke Math. J. 165(5) (2016), 935--984.
			 	
			 	\bibitem{C95} J. Y. Chemin, {\it Fluides parfaits incompressibles,} Ast\'{e}risque 230, Soci\'{e}t\'{e} Math\'{e}matique de France, (1995).
			 	
			 	%\bibitem{Cha2} D.~ Chae, {\it  Local existence and blow-up criterion for the Euler equations in the Besov spaces.} Asymptotic Analysis {\bf 38} (2004) 339--358.
			 	
			 	%\bibitem{C-C-C-G-W} D. Chae, P. Constantin, D. C\'ordoba, F. Gancedo, and J. Wu, {\it Generalized surface quasi-geostrophic equations with singular velocities.}   Comm. Pure Appl. Math., 65 (2012), no. 8, 1037-1066.
			 	
			 	%\bibitem{Chapl} S. A. Chaplygin, {\it On a pulsating cylindrical vortex.} Translated from the 1899 Russian original by G. Krichevets, edited by D. Blackmore and with comments by V. V. Meleshko. Regul. Chaotic Dyn. 12 (2007), no. 1, 101?116.
			 	
			 	
			 	%\bibitem{C-M-T} P. Constantin, A. J. Majda, and E. Tabak, {\it Formation of strong fronts in the 2-D quasigeostrophic thermal active scalar.}  Nonlinearity, 7 (1994), no. 6, 1495-1533.
			 	
			 	%\bibitem{C-F-M-R} D. C\'ordoba, M. A. Fontelos, A. M. Mancho and J. L. Rodrigo, {\it Evidence of singularities for a family of contour dynamics equations.}  Proc. Natl. Acad. Sci. USA 102 (2005),  5949--5952.
			 	
			 	%\bibitem{CR71} M. G. \ Crandall and P.H.\ Rabinowitz,  {\it Bifurcation from simple eigenvalues.} J. of Func. Analysis,  8 (1971),  321--340.
			 	
			 	\bibitem{DHR19} D. G. Dritschel, T. Hmidi, C. Renault, \textit{Imperfect bifurcation for the shallow-water quasi-geostrophic equations}, Arch. Ration. Mech. Anal. 231 (2019), no. 3, 1853--1915.
			 	
			 	\bibitem{DJ20} D. G. Dritschel, M. M. Jalali, \textit{Stability and evolution of two opposite-signed quasi-geostrophic shallow-water vortex patches}, Geophys. Astrophys. Fluid Dyn. 114 (2020), no. 4-5, 561--587.
			 	
			 	\bibitem{DP11} D. G. Dritschel, H. Plotka, \textit{Shallow-water vortex equilibria and their stability}, Journal of Physics: Conference Series, 318 (2011), no. 6.
			 	
			 	%				\bibitem{DZ16} R. Danchin, X. Zhang, \textit{Global persistence of geometrical structures for the boussinesq equation with no diffusion,} Communications in Partial Differential Equations, (2016).
			 	%				
			 	%				\bibitem{DZ16-2} R. Danchin, X. Zhang, \textit{On the persistence of hölder regular patches of density for the inhomogeneous navier-stokes equations}, Journal de l’Ecole polytechnique - Math\'{e}matiques 4, (2016).
			 	%				
			 	\bibitem{DZ78} G. S. Deem, N. J. Zabusky, {\it Vortex waves : Stationary "V-states", Interactions, Recurrence, and Breaking,}
			 	Phys. Rev. Lett.  40  (1978), no. 13, 859--862.
			 	
			 	
			 	%\bibitem{De} J.-M.~ Delort, {\it  Existence de nappes de tourbillon en dimension deux}, J. Amer. Math. Sot., Vol. {\bf 4} (1991)
			 	% 553--586. 
			 	
			 	%\bibitem{DM}  R.~ DiPerna and A.~ Madja, {\it Concentrations in regularization for 2D incompressible flow}. Comm. Pure Appl.
			 	%Math. {\bf 40} (1987), 301--345.
			 	
			 	%\bibitem{DR} D. G.  Dritschel, {\it The nonlinear evolution of rotating configurations of uniform vorticity}, J.Fluid Mech. 172 (1986), 157--182.
			 	
			 	%\bibitem{Du} P. L. Duren, {\it Univalent functions.} Grundlehren der mathematischen Wissenschaften 259. Springer-Verlag, New York, (1983).
			 	
			 	%\bibitem{Ell-Kuk } L.H. Eliasson and S.B. Kuksin.{\it  KAM for the nonlinear Schr\"odinger equation.} Ann. Math 172 (2010), 371-435.
			 	
			 	%\bibitem{Fas} C.D.  Fassnacht, C. R. Keeton, D.  Khavinson,  {\it Gravitational lensing by elliptical galaxies, and the Schwarz function.} Analysis and mathematical physics, 115?129, Trends Math., Birkh?user, Basel, 2009.
			 	
			 	\bibitem{FGMP19}
			 	R. Feola, F. Giuliani, R. Montalto, M. Procesi,
			 	\textit{Reducibility of first order linear operators on tori via Moser's
			 		theorem}, J. Funct. Anal. 276 (2019), no. 3, 932--970.
			 	
			 	
			 	%\bibitem{Flierl} G. R. Flierl, L. M. Polvani,{  \it Generalized Kirchhoff vortices.} Phys. Fluids 29 (1986), 2376--2379.
			 	
			 	%\bibitem{G} F. Gancedo,{  \it Existence for the $\alpha$-patch model and the QG sharp front in Sobolev spaces.}  Adv. Math., 217 (2008), no. 6 2569--2598.
			 	
			 	%\bibitem{Ga} A. E.\ Gatto, {\it On the boundedness on inhomogeneous Lipschitz spaces of fractional integrals,
			 	% singular integrals and hypersingular integrals associated to non-doubling measures on metric spaces,}
			 	% Collect. Math. {\bf 60} (2009), 101-114.
			 	
			 	%\bibitem{Har} B. J. Harvey, M. H. P.  Ambaum, {\it Perturbed Rankine vortices in surface quasi-geostrophic dynamics.} Geophys. Astrophys. Fluid Dyn. 105 (2011), no. 4-5, 377--391.
			 	
			 	%\bibitem{HAC} B. J. Harvey, M. H. P.  Ambaum, X. J. Carton {\it Instability of Shielded Surface Temperature Vortices}. J. Atmos. Sci., 68 (2010) 964--971.
			 	
			 	%\bibitem{GHS} C. Garc\`ia, T. Hmidi, J. Soler, {\it Non uniform rotating vortices and periodic orbits for the two-dimensional Euler Equations.} Arch. for Ration. Mech. and Anal.  238, 929--1085 (2020).
			 	
			 	\bibitem{G20} C. Garc\'{i}a, \textit{K\'{a}rm\'{a}n vortex street in incompressible fluid models}, Nonlinearity 33 (2020), no. 4, 1625--1676.
			 	
			 	\bibitem{G21} C. Garc\'{i}a, \textit{Vortex patches choreography for active scalar equations}, Journal of Nonlinear Science 31 (2021), no. 75, 1432--1467.
			 	
			 	\bibitem{G19}  J. G\'omez-Serrano, {\it On the existence of stationary patches,} Advances in Mathematics 343 (2019), 110--140.
			 	
			 	%\bibitem{G-S-Y}  J. G\'omez-Serrano,  J. Park, J. Shi, Y.  Yao {\it Symmetry in stationary and uniformly-rotating solutions of active scalar equations.} arXiv:1908.01722 
			 	
			 	%				\bibitem{HH15-1} Z. Hassainia, T. Hmidi, \textit{On the inviscid boussinesq system with rough initial data}, Journal of
			 	%				Mathematical Analysis and Applications 430 (2015), no. 2, 777--809.
			 	%				
			 	\bibitem{HH15}  Z. Hassainia, T. Hmidi, {\it On the V-States for the generalized quasi-geostrophic equations,} Comm. Math. Phys. 337
			 	(2015), no. 1, 321--377.
			 	
			 	\bibitem{HH21} Z. Hassainia, T. Hmidi, \textit{Steady asymmetric vortex pairs for Euler equations}, American Institut of Mathematical Science 41 (2021), no. 4, 1939--1969.
			 	
			 	\bibitem{HHH18} Z. Hassainia, T. Hmidi, F. de la Hoz, \textit{Doubly connected V-states for the generalized surface quasi-geostrophic equations}, Cambridge University Press 439 (2018), 90--117.
			 	
			 	\bibitem{HHHM15} Z. Hassainia, T. Hmidi, F. de la Hoz, J. Mateu, \textit{An analytical and numerical study of steady patches in the disc}, Analysis and PDE 9 (2015), no. 10.
			 	
			 	\bibitem{HMW20} Z. Hassainia, N. Masmoudi, M. H. Wheeler,  {\textit Global bifurcation of rotating vortex patches},  Comm. Pure Appl. Math. Vol. LXXIII (2020), 1933--1980 
			 	
			 	\bibitem{HHM21} Z. Hassainia, T. Hmidi, N. Masmoudi, {\textit KAM theory for active scalar equations}, arXiv:2110.08615
			 	
			 	\bibitem{HW21} Z. Hassainia, M. Wheeler, \textit{Multipole vortex patch equilibria for active scalar equations}, arXiv:2103.06839.
			 	
			 	%				\bibitem{H05} T. Hmidi, \textit{R\'{e}gularit\'{e} höld\'{e}rienne des poches de tourbillon visqueuses}, Journal de Math\'{e}matiques
			 	%				Pures et Appliqu\'{e}es 84 (2005), no. 11, 1455--1495.
			 	
			 	%\bibitem{H15}  T. Hmidi,  {\it On the trivial solutions for the rotating patch model,} J. Evol. Equ. 15 (2015), no. 4, 801--816.
			 	
			 	\bibitem{HHMV16} T. Hmidi,  F. de la Hoz, J.  Mateu, J.  Verdera, {\it  Doubly connected V-states for the planar Euler equations,}  SIAM J. Math. Anal. 48 (2016), no. 3, 1892--1928. 
			 	
			 	
			 	\bibitem{HM16} T. Hmidi, J. Mateu, {\it Bifurcation of rotating patches from Kirchhoff vortices,} Discrete Contin. Dyn. Syst. 36 (2016), no. 10, 5401--5422. 
			 	
			 	\bibitem{HM16-2} T. Hmidi, J. Mateu, {\it Degenerate bifurcation of the rotating patches,} Adv. Math. 302 (2016), 799--850. 
			 	
			 	\bibitem{HM17} T. Hmidi,  J.  Mateu, {\it  Existence of corotating and counter-rotating vortex pairs for active scalar equations,} Comm. Math. Phys. 350 (2017), no. 2, 699--747.   
			 	
			 	\bibitem{HMV13} T. Hmidi, J.  Mateu, J.  Verdera, {\it Boundary Regularity of Rotating Vortex Patches,} Arch. Ration. Mech. Anal. 209 (2013), no. 1, 171--208.
			 	
			 	\bibitem{HMV15} T. Hmidi, J.  Mateu, J.  Verdera, {\it On rotating doubly connected vortices,}  J. Differential Equations 258 (2015), no. 4, 1395--1429.
			 	
			 	%\bibitem{HR17}  T. Hmidi, C. Renault, {\it Existence of small loops in a bifurcation diagram near degenerate eigenvalues}. Nonlinearity 30 (2017), no. 10, 3821--3852.
			 	
			 	\bibitem{HS11} P. E. Hoggan, A. Sidi, \textit{Asymptotics of modified Bessel functions of high order}, Int. J. Pure Appl. Math. 71 (2011), no. 3, 481--498.
			 	
			 	\bibitem{IPT05} G. Iooss, P. Plotnikov, J. Toland, \textit{Standing waves on an infinitely deep perfect fluid under gravity}, Arch. Ration. Mech.
			 	Anal. 177 (2005), no. 3, 367--478.
			 	
			 	%\bibitem{Juk} M. Juckes, {\it Quasigeostrophic dynamics of the tropopause.} J. Armos. Sci. (1994) 2756--2768.
			 	
			 	%\bibitem{Kato} T. Kato and G. Ponce, {\it  Well-posedness of the Euler and Navier-Stokes equations in the Lebesgue spaces $L^p_s(\RR^2)$}. Rev. Mat. Iberoamericana {\bf 2} (1986), no. 1-2, 73--88.
			 	
			 	
			 	%\bibitem{HAC} B. J. Harvey, M. H. P.  Ambaum, X. J. Carton {\it Perturbed Rankine vortices in surface quasi-geostrophic dynamics.} Geophys. Astrophys. Fluid Dyn. (2010)1--15.
			 	
			 	%\bibitem{Kida} S. Kida, {\it Motion of an elliptical vortex in a uniform shear flow,} J. Phys. Soc. Japan. 50 (1981) \mbox{3517--3520.}
			 	
			 	%\bibitem{K11} H. Kielh\"ofer, {\it Bifurcation Theory: An Introduction With Applications to Partial Differential Equations,} Springer (2011)
			 	
			 	\bibitem{K74} G. Kirchhoff, \textit{Vorlesungen uber mathematische Physik}, Leipzig, (1874).
			 	
			 	\bibitem{K54} A. N. Kolmogorov, {\it On the persistence of conditionally periodic motions under a small change of the hamiltonian function,} Doklady Akad. Nauk SSSR 98 (1954),
			 	527--530.
			 	
			 	%\bibitem{Kuksin1} S.  B. Kuksin,{\it  Nearly integrable infinite-dimensional Hamiltonian systems,} Lecture Notes in Mathematics, vol. 1556, Springer-Verlag, Berlin, 1993
			 	
			 	%\bibitem{Kuksin-Poschel} S. Kuksin and J. P\"oschel, {\it  Invariant cantor manifolds of quasi-periodicoscillations for a nonlinear schr\"odinger equation,} Annals of Math 2 (1996), no. 143, 149--179.
			 	
			 	%\bibitem{L} H. \ Lamb, {\it Hydrodynamics}, Dover Publications, New York, (1945).
			 	
			 	%\bibitem{Lap} G. Lapeyre, P. Klein, {\it Dynamics of the upper oceanic layers in terms of surface quasigeostrophic theory.} J. Phys. Oceanogr. 36, (2006), 165--176.
			 	
			 	\bibitem{Leb65} N. N. Lebedev, \textit{Special Functions and their applications}, Prentice-Hall, (1965).
			 	
			 	%\bibitem{MZ20} O. Melkemi, M. Zerguine, \textit{Local persistence of geometric structures of the inviscid nonlinear boussinesq system}, arXiv:2005.11605.
			 	
			 	\bibitem{M62} J. Moser, \textit{On invariant curves of area-preserving mappings of an annulus}, Nachr. Akad. Wiss., Göttingen, Math. Phys. Kl. (1962), 1--20. 
			 	
			 	%\bibitem{M-O} W.\ Magnus, F.\ Oberhettinger {\it Formeln und satze fur die speziellen funktionen der mathematischen physik},  Berlin.Gottingen.Heidelberg : Springer , (1948) . 
			 	
			 	
			 	%\bibitem{MOV} J.\ Mateu, J.\ Orobitg and J.\ Verdera, {\it Extra cancellation of even Calder\'{o}n-Zygmund operators and quasiconformal mappings}, J. Math. Pures Appl. \textbf{91}(4)(2009), 402--431.
			 	
			 	\bibitem{N54} J. Nash, \textit{C1-isometric imbeddings}, Annals of Mathematics 60 (1954), no. 3, 383--396.
			 	
			 	%\bibitem{Neu} J. Neu, {\it  The dynamics of columnar vortex in an imposed strain?} Phys. Fluids 27 (1984) 12397--2402.
			 	
			 	%\bibitem{M-O} W.\ Magnus, F.\ Oberhettinger {\it Formeln und satze fur die speziellen funktionen der mathematischen physik},  Berlin.Gottingen.Heidelberg : Springer , 1948 . 
			 	
			 	%\bibitem{New} P. K. Newton, {\it The N-Vortex Problem. Analytical Techniques,} Springer, New York, (2001).
			 	
			 	%\bibitem{poschel1} J. P\"oschel, {\it A KAM-theorem for some nonlinear partial differential equations,} Annali della Scuola Normale Superiore di Pisa. Classe di Scienze. Serie IV 23 (1996), no. 1, 119-148.
			 	
			 	%\bibitem{Nov} P. S. Novikov, {\it On the uniqueness for the inverse problem of potential theory}, Dokl. Akad. Nauk SSSR ${\bf18}$ (1938), 165-168.
			 	
			 	%\bibitem{P} Ch. Pommerenke,  {\it Boundary behaviour of conformal maps.}  Springer-Verlag, Berlin, (1992).
			 	
			 	%\bibitem{Rain} E. D. Rainville,{ Special Functions} Macmillan, New York, (1960).
			 	
			 	\bibitem{PT01} P. Plotnikov, J. Toland, \textit{Nash-Moser theory for standing water waves}, Arch. Ration. Mech. Anal. 159 (2001), no. 1, 1--83.
			 	
			 	%\bibitem{Pyralti} Pyartli. A.S. {\it Approximations diophantiennes sur les sous-vari\'et\'es de l'espace euclidien} (en russe). Funkcional. Anal. i Prilozen. 3 (1969) 59--62 (trad. anglaise Functional Anal. Appl. 3 (1969), 303--306.
			 	
			 	%\bibitem{R} J. L. Rodrigo,  {\it On the evolution of sharp fronts for the quasi-geostrophic equation.}   Comm. Pure Appl. Math., 58, no.6  (2005), 821--866.
			 	
			 	\bibitem{R01} H. R\"ussmann, {\it Invariant tori in non-degenerate nearly integrable Hamiltonian systems,} Regul. Chaotic Dyn. 6 (2001), no. 2, 119--204.
			 	
			 	%\bibitem{sevryuk} M. B. Sevryuk, {\it Invariant tori in quasi-periodic non-autonomous dynamical systems via Herman's method,} Discrete Contin. Dyn. Syst. 18 (2007), no. 2--3, 569--595.
			 	
			 	%\bibitem{Uk}
			 	%M. R. Ukhovskii, V. I. Iudovich, {\it Axially symmetric flows of ideal and viscous fluids filling the whole space},
			 	%Prikl. Mat. Meh. {\bf 32} (1968), no. 1, 59-69.
			 	
			 	%\bibitem{Var} A. N.  Varchenko, P. I.  Etingof,{\it  Why the boundary of a round drop becomes a curve of order four.} University Lecture Series, 3. American Mathematical Society, Providence, RI, 1992.
			 	%\bibitem{V} J.\ Verdera, {\it $L^2$ boundedness of the Cauchy Integral and Menger curvature,} Contemporary Mathematics 277 (2001), 139--158.
			 	
			 	\bibitem{V17} G. K. Vallis, \textit{Atmospheric and Oceanic Fluid Dynamics: Fundamentals and Large-Scale Circulation}, Cambridge University Press, 2nd edition, (2017).
			 	
			 	%\bibitem{Vishik2} M. Vishik, {\it Hydrodynamics in Besov Spaces}, Arch. Rational Mech. Anal. {\bf{145}} (1998), 197--214.  
			 	
			 	%\bibitem{WS} S. E. Warschawski, {\it On the higher derivatives at the boundary in conformal mapping,} Trans. Amer. Math. Soc.  38 (1935), no. 2, 310--340.
			 	
			 	\bibitem{W95} G. N. Watson, {\it A Treatise on the Theory of Bessel Functions}, Cambridge University Press, (1922).
			 	
			 	%\bibitem{W} R.\ Wittmann, {\it Application of a Theorem of M.G. Krein to singular integrals,}
			 	%Trans. Amer. Math. Soc. {\bf 299}(2) (1987), 581--599.
			 	
			 	%\bibitem{Wolib} W. Wolibner,{\it  Un th\'eor\`eme sur lÍexistence du mouvement plan dÍun fluide parfait homog\`ene, incompressible, pendant un temps infiniment long,} Math. Z, Vol. 37, 1933, pp. 698--627.
			 	
			 	%\bibitem{W90} C. E. Wayne, {\it Periodic and quasi-periodic solutions of nonlinear wave equationsvia KAM theory,} Communications in Mathematical Physics 127 (1990), no. 3, 479--528.
			 	
			 	%\bibitem{WOZ84} H.M.\ Wu, E.A.\ Overman II and N.J. Zabusky {\it Steady-state solutions of the Euler equations in two dimensions : rotating and translating  V-states with limiting cases I. Algorithms ans results}, J. Comput. Phys. { 53} (1984), 42--71.
			 	
			 	\bibitem{Y63} Y. Yudovich, {\it Nonstationary flow of an ideal incompressible liquid,} Zh. Vych. Mat. 3 (1963), 1032--1066.
			 	
			 	%\bibitem{Y95} Y.~ Yudovich, {\it Uniqueness theorem for the basic nonstationary problem in the dynamics of an ideal incompressible fluid}. Math. Res. Lett., {\bf 2} (1995), 27--38.
			 	
			 	%\bibitem{HZ79} N. Zabusky, M. H. Hughes, K.V. Roberts, {\it Contour dynamics for the Euler equations in two dimensions.} J. Comput. Phys. 30 (1979), no. 1, 96--106.
			 \end{thebibliography}
			%\bibliographystyle{plain}
			 
		\end{document}